\documentclass[english]{article}
\usepackage[T1]{fontenc}
\usepackage[utf8]{inputenc}
\usepackage{amsmath,amsthm} 
\usepackage{graphicx}
\usepackage{amssymb}
\usepackage{enumitem}
\usepackage{natbib} 
\usepackage{bussproofs}
\usepackage[letterpaper]{geometry} 
\usepackage{appendix}
\newcommand\pfun{\mathrel{\ooalign{\hfil$\mapstochar\mkern5mu$\hfil\cr$\to$\cr}}}

\usepackage{thmtools,thm-restate}

\usepackage{chngcntr}
\counterwithin{figure}{section}

\usepackage{tikz} 
\usepackage{pgf} 
\usepackage{setspace}
\usetikzlibrary{patterns,automata,arrows,shapes,snakes,topaths,trees,backgrounds,positioning,through,calc}
\usepackage{multicol}
\usepackage{graphicx}
\usepackage{stmaryrd}  
\usepackage{hyperref}

\hypersetup{colorlinks=true, citecolor=blue, linkcolor=blue}  

\usepackage{comment}

\theoremstyle{definition}
\newtheorem{theorem}{Theorem}[section]
\newtheorem{observation}[theorem]{Observation}
\newtheorem{proposition}[theorem]{Proposition}
\newtheorem{corollary}[theorem]{Corollary}
\newtheorem{definition}[theorem]{Definition}

\newtheorem{lemma}[theorem]{Lemma}
\newtheorem{fact}[theorem]{Fact}
\newtheorem{example}[theorem]{Example}
\newtheorem{notation}[theorem]{Notation}
\newtheorem{remark}[theorem]{Remark}  
\newtheorem{problem}{Problem}

\makeatletter

\newcommand{\sig}{\Phi} 
\newcommand{\ind}{ I} 
\newcommand{\wo}{\mathrm} 

\newcommand{\inc}{\bot} 
\newcommand{\pow}{\wp} 
 
\newcommand{\comp}{\between}
\newcommand{\incomp}{\,\bot\,}

\newcommand{\cof}{\sqsubseteq_{s}}
\newcommand{\cofsub}{{s}}
\newcommand{\meet}{\wedge}
\newcommand{\bigmeet}{\bigwedge}
\newcommand{\join}{\vee}
\newcommand{\bigjoin}{\bigvee}

\newcommand{\Rlatbox}{\boxdot}
\newcommand{\latbox}{\boxdot}

\newcommand{\adm}{P} 
\newcommand{\under}{\mathsf{b}} 
\newcommand{\rela}{\mathsf{p}} 
\newcommand{\fullV}{\mathsf{u}} 
\newcommand{\ff}{\mathsf{f}} 
\newcommand{\gff}{\mathsf{g}} 
\newcommand{\found}{\sharp} 
\newcommand{\full}{\sharp} 
\newcommand{\tight}{\Box}

\newcommand{\upR}{\textbf{up}-$\boldsymbol{R}$}
\newcommand{\Rcomm}{$\boldsymbol{R}$-\textbf{com}}
\newcommand{\Rwin}{$\boldsymbol{R\mathord{\Rightarrow}}\textbf{\underline{win}}$}

\newcommand{\RWin}{$\boldsymbol{R\mathord{\Leftrightarrow}}\textbf{\underline{win}}$}
\newcommand{\Rwinweak}{$\boldsymbol{R\mathord{\Rightarrow}}\textbf{win}$}
\newcommand{\Rrule}{$\boldsymbol{R}$-\textbf{rule}}
\newcommand{\Fwinweak}{$\boldsymbol{f\mathord{\Rightarrow}}\textbf{win}$}
\newcommand{\Frule}{$\boldsymbol{f}$-\textbf{rule}}
\newcommand{\Fmon}{$\boldsymbol{f}$-\textbf{monotonicity}}
\newcommand{\Rdown}{$\boldsymbol{R}$-\textbf{down}}
\newcommand{\Rdense}{$\boldsymbol{R}$-\textbf{dense}}
\newcommand{\Rref}{$\boldsymbol{R}$-\textbf{refinability}}
\newcommand{\RrefPlus}{$\boldsymbol{R}$-\textbf{refinability}$^+$}
\newcommand{\RrefPlusPlus}{$\boldsymbol{R}$-\textbf{refinability}$^{++}$}
\newcommand{\Rmax}{$\boldsymbol{R}$-\textbf{max}}
\newcommand{\Rmaxe}{$\boldsymbol{R}$-\textbf{maxe}}
\newcommand{\Rprinc}{$\boldsymbol{R}$-\textbf{princ}}

\newcommand{\SqMatch}{\textit{$\sqsubseteq$-matching}}
\newcommand{\RMatch}{\textit{$R$-matching}}
\newcommand{\PullBack}{\textit{pull back}}

\newcommand{\SqForth}{\textit{$\sqsubseteq$-forth}}

\newcommand{\SqBack}{\textit{$\sqsubseteq$-back}}

\newcommand{\RForth}{\textit{$R$-forth}}

\newcommand{\SRBack}{\textit{$R$-back}}

\newcommand{\pRBack}{\textit{$p$-$R$-back}}
\newcommand{\pSqBack}{\textit{$p$-$\sqsubseteq$-back}}

\newcommand{\primvee}{\mbox{\,\textcircled{$\vee$}\,}}
\newcommand{\Bethvee}{\veebar}

\newcommand{\var}{\mathrm}

\makeatother

\usepackage{babel} 

\begin{document} 
     
\title{Possibility Frames and Forcing for Modal Logic\footnote{Previous versions of this article circulated online as the working papers \citealt{Holliday2015}, \citealt{Holliday2016}, and \citealt{Holliday2018}. The present version updates \citealt{Holliday2018} based on the review process for \textit{The Australasian Journal of Logic} and errata.}}
\author{Wesley H. Holliday\\{\small University of California, Berkeley}}

\date{{\small Published in \textit{The Australasian Journal of Logic}, Vol.~22, No.~2, 2025, 44--288.}}
         
\maketitle
 
\begin{abstract}
This paper develops the model theory of normal modal logics based on partial ``possibilities'' instead of total ``worlds,'' following Humberstone \citeyearpar{Humberstone1981} instead of Kripke \citeyearpar{Kripke1963}. Possibility semantics can be seen as extending to modal logic the semantics for classical logic used in weak forcing in set theory, or as semanticizing a negative translation of classical modal logic into intuitionistic modal logic. Thus, possibility frames are based on posets with accessibility relations, like intuitionistic modal frames, but with the constraint that the interpretation of every formula is a regular open set in the Alexandrov topology on the poset. The standard world frames for modal logic are the special case of possibility frames wherein the poset is discrete. The analogues of classical Kripke frames, i.e., \textit{full} world frames, are \textit{full} possibility frames,  in which propositional variables may be interpreted as any regular open sets. 

We develop the beginnings of duality theory, definability/correspondence theory, and completeness theory for possibility frames. The duality theory, relating possibility frames to Boolean algebras with operators (BAOs), shows the way in which full possibility frames are a generalization of Kripke frames.  Whereas Thomason \citeyearpar{Thomason1975} established a duality between the category of Kripke frames with p-morphisms and the category of \textit{complete} ($\mathcal{C}$), \textit{atomic} ($\mathcal{A}$), and \textit{completely additive} ($\mathcal{V}$) BAOs with complete BAO-homomorphisms (these categories being dually equivalent), we establish a duality between the category of full possibility frames with possibility morphisms and the category of $\mathcal{CV}$-BAOs, i.e., allowing \textit{non-atomic} BAOs, with complete BAO-homomorphisms (the latter category being dually equivalent to a reflective subcategory of the former). This parallels the connection between forcing posets and Boolean-valued models in set theory, but now with accessibility relations and modal operators involved. Generalizing further, we introduce a class of \textit{principal} possibility frames that capture the generality of $\mathcal{V}$-BAOs. If we do not require a full or principal frame, then every nontrivial BAO has an equivalent possibility frame, whose possibilities are proper filters in the BAO. With this filter representation, which does not require the ultrafilter axiom, we are led to a notion of \textit{filter-descriptive} possibility frames. Whereas Goldblatt \citeyearpar{Goldblatt1974} showed that the category of BAOs with BAO-homomorphisms is dually equivalent to the category of descriptive world frames with p-morphisms, we show that the category of BAOs with BAO-homomorphisms is dually equivalent to the category of filter-descriptive possibility frames with p-morphisms. Applying our duality theory to definability theory, we prove analogues for possibility semantics of theorems of Goldblatt \citeyearpar{Goldblatt1974} and Goldblatt and Thomason \citeyearpar{Goldblatt1975} characterizing modally definable classes of frames. In addition, we discuss analogues for possibility semantics of first-order correspondence results in the style of Lemmon and Scott \citeyearpar{Lemmon1977}, Sahlqvist \citeyearpar{Sahlqvist1975}, and van Benthem \citeyearpar{Benthem1976}. Finally, applying our duality theory to completeness theory, we show that there are continuum many normal modal logics that can be characterized by full possibility frames but not by Kripke frames, that all consistent Sahlqvist logics can be characterized by full possibility frames that contain no worlds, and that all consistent normal modal logics can be characterized by filter-descriptive possibility frames.\end{abstract} 

\hspace{.14in}\textbf{Keywords:} modal logic, possibility semantics, Boolean algebras with operators,

\hspace{.878in} duality theory, regular open algebras, Kripke-frame incompleteness\vspace{.1in}

\hspace{.14in}\textbf{MSC:} 03B45, 03G05 

\noindent 

\tableofcontents
  
\section{Introduction}\label{intro}
 
The model theory of modal logic has been developed extensively on the basis of \textit{possible world semantics}, as presented in its now standard form by Kripke \citeyearpar{Kripke1963}.\footnote{For a historical overview, see \citealt{Goldblatt2006}; for surveys of the theory, see \citealt{Blackburn2007}, \citealt{Goldblatt2006}, and \citealt{Goranko2007}; and for textbooks, see, e.g., \citealt{Blackburn2001} and \citealt{Chagrov1997}.} In this paper, we develop the beginnings of a more general theory of \textit{possibility semantics} for modal logic, building on the work of Humberstone \citeyearpar{Humberstone1981}. 

 As presented in \S~\ref{FromPartToPoss}, possibility semantics can be seen as extending to modal logic the semantics for classical logic used in weak forcing in set theory, or as semanticizing a negative translation of classical modal logic into intuitionistic modal logic. Thus, possibility frames are based on posets with accessibility relations, like intuitionistic modal frames, but with the constraint that the interpretation of every formula is a regular open set in the Alexandrov topology on the poset. Unlike intuitionistic Kripke semantics, possibility semantics allows a partial possibility to determine that a disjunction is true without determining which disjunct is true, leaving that to be determined by refinements of the possibility. Worlds are totally determinate possibilities---endpoints in the poset---which not all possibility frames contain; and the standard world frames for modal logic are the special case of possibility frames wherein all possibilities are worlds---so the poset is discrete. The analogues of classical Kripke frames, i.e., \textit{full} world frames, are \textit{full} possibility frames,  in which propositional variables may be interpreted as any regular open sets. At the end of \S~\ref{FromPartToPoss}, we give examples of full possibility frames whose logics are not validated by any Kripke frames.
 
An important motivation for the move from total worlds to partial possibilities is philosophical (\citealt{Humberstone1981}; \citealt[p. 564]{Edgington1985}; \citealt[\S~4]{Edgington2010}; \citealt[\S\S~6.1-6.2]{Rumfitt2015}), but in this paper we focus on the logical ramifications. After expanding the pr\'{e}cis of possibility semantics above in \S~\ref{FromPartToPoss}, we introduce a number of useful concepts for the study of possibility semantics in \S~\ref{morphisms} and \S~\ref{SpecialClasses}: in \S~\ref{morphisms}, we introduce notions of \textit{possibility morphisms} between possibility frames, which reduce to the standard p-morphisms (bounded morphisms) when the frames in question are Kripke frames; and in \S~\ref{SpecialClasses}, we provide a catalogue of classes of frames that are of special importance in possibility semantics. With this toolkit acquired, we proceed to develop for possibility semantics the beginnings of duality theory, definability/correspondence theory, and completeness theory, the ``three pillars of wisdom'' supporting modal logic \citep[p. 331]{Benthem1980}.
 
In \S~\ref{DualityTheory}, we prove our main results concerning dualities between possibility frames on the one hand and Boolean algebras with operators (BAOs) \citep{Jonsson1952a,Jonsson1952b} on the other (throughout, BAO-to-frame constructions and categories of BAOs are restricted to nontrivial BAOs as in Definition~\ref{BAOs}). From an algebraic perspective, the possibility semantics developed in this paper occupies a special place among generalizations of Kripke semantics.\footnote{It is not known (at least to this author) whether Humberstone's original frames for possibility semantics are more general than Kripke frames. We compare our possibility frames with Humberstone's frames at the end of \S~\ref{FullFrames} and in Appendix \S~\ref{Strengths}.}  Thomason \citeyearpar{Thomason1975} showed that Kripke frames correspond to BAOs that are \textit{complete} ($\mathcal{C}$), \textit{atomic} ($\mathcal{A}$), and \textit{completely additive} ($\mathcal{V}$).\footnote{The notation `$\mathcal{C}$', `$\mathcal{A}$', `$\mathcal{V}$', and `$\mathcal{T}$' (below) is from \citealt{Litak2005,Litak2005b,Litak2008}. See \S~\ref{DualityTheory} for the definitions of these classes of BAOs.} Semantically: any Kripke frame can be turned into a modally equivalent $\mathcal{CAV}$-BAO, and vice versa. Categorically: the category of Kripke frames with p-morphisms is dually equivalent to the category of $\mathcal{CAV}$-BAOs with complete BAO-homomorphisms. More general types of frames have been proposed that correspond to $\mathcal{CA}$-BAOs, viz. normal neighborhood frames as in \citealt{Dosen1989}, and $\mathcal{AV}$-BAOs, viz. discrete general frames as in \citealt{tenCate2007}. 

Possibility frames generalize in a different direction, by dropping atomicity. We show that \textit{full} possibility frames correspond to $\mathcal{CV}$-BAOs. Semantically: any full possibility frame can be turned into a modally equivalent $\mathcal{CV}$-BAO (\S~\ref{PossToBAO}), and vice versa (\S~\ref{VtoPossSection}). Categorically: the category \textbf{$\mathcal{CV}$-BAO} of $\mathcal{CV}$-BAOs with complete BAO-homomorphisms is dually equivalent to a \textit{reflective subcategory} of the category \textbf{FullPoss} of full possibility frames with what we call ``strict possibility morphisms'' (\S~\ref{DualEquiv}). This connection between $\mathcal{CV}$-BAOs and possibility frames parallels the connection between Boolean-valued models and forcing posets in set theory (see, e.g., \citealt{Takeuti1973}), but now with modal operators and accessibility relations involved. (As we briefly note at the end of \S~\ref{FullFrames}, if we consider normal \textit{neighborhood} versions of full possibility frames, then such frames are to $\mathcal{C}$-BAOs as our relational versions of full possibility frames are to $\mathcal{CV}$-BAOs.) 

One of our duality results for $\mathcal{CV}$-BAOs can be sketched as follows. On the algebraic side, given any complete Boolean algebra, we form a $\mathcal{CV}$-BAO by adding an operator $\blacklozenge$ on the algebra such that for the bottom element $\bot$ of the algebra, $\blacklozenge\bot=\bot$, and for any set $X$ of elements, $\blacklozenge\bigjoin X= \bigjoin \{\blacklozenge x\mid x\in X\}$. For morphisms between $\mathcal{CV}$-BAOs, we consider Boolean algebra homomorphisms that also preserve arbitrary joins and $\blacklozenge$ (complete BAO-homomorphisms). Then the question naturally arises: can we have an analogue for the category \textbf{$\mathcal{CV}$-BAO} of Thomason's theorem for \textbf{$\mathcal{CAV}$-BAO}? On the frame-theoretic side, instead of adding a modal operator $\blacklozenge$, we add an accessibility relation $R$. Given any complete Boolean algebra, we form what we call a \textit{rich} possibility frame (\S~\ref{RichFrames}) by deleting the bottom element of the algebra and then adding a binary relation $R$ on the underlying set that interacts with the order relation $\leq$ from the bottomless algebra in a nice first-order way (called \RWin{} and motivated game-theoretically in \S~\ref{FullFrames}): $xRy$ iff $\forall y'\leq y$ $\exists x'\leq x$ $\forall x''\leq x'$ $\exists y''\leq y'$ $x''Ry''$.\footnote{In \S~\ref{FromPartToPoss} and following, we use the symbol `$\sqsubseteq$' instead of `$\leq$' for the order relation in possibility frames. It is also important to note that our notation in this paper is flipped relative to most of the literature on intuitionistic Kripke semantics: we take $x\sqsubseteq y$ to mean that $x$ is a refinement/further specification/extension/etc.~of $y$ (in agreement with most of the literature on set-theoretic forcing, as well as \citealt{Dragalin1988}), whereas much of the intuitionistic literature would take it to mean that $y$ is an extension of $x$.} For morphisms between such frames, we consider maps $h$ that are p-morphisms with respect to both $\leq$ and $R$: $x\leq y$ implies $h(x)\leq' h(y)$; $xRy$ implies $h(x)R'h(y)$; if $y'\leq ' h(x)$, then $\exists y$: $y\leq x$ and $h(y)=y'$; and  if $h(x)R'y'$, then $\exists y$: $xRy$ and $h(y)=y'$. These maps are a special case of the strict possibility morphisms mentioned above (see \S~\ref{morphisms}). Together rich possibility frames with p-morphisms form the category \textbf{RichPoss}. In \S~\ref{DualEquiv}, we show that  \textbf{RichPoss} is dually equivalent to~\textbf{$\mathcal{CV}$-BAO}.

The rich possibility frames just sketched are a special case of the \textit{full} possibility frames described in \S~\ref{FromPartToPoss}, which have much looser constraints in general, e.g., the  order relation on possibilities can be an arbitrary partial order. Yet we show that for every full possibility frame, there is a strict possibility morphism from that frame to a rich possibility frame that validates exactly the same modal formulas. In addition to this semantic fact, we prove a categorical fact: \textbf{RichPoss} is a reflective subcategory of  \textbf{FullPoss}. 

Going beyond full possibility frames, we show that our \textit{principal} possibility frames (from \S~\ref{PrincFrames}) correspond semantically to $\mathcal{V}$-BAOs (\S\S~\ref{PossToBAO}-\ref{VtoPossSection}). We turn any $\mathcal{V}$-BAO into a semantically equivalent principal possibility frame by deleting the bottom element of the algebra and then defining a binary relation $R$ on the underlying set by: $xRy$ iff for all $y'\leq y$ such that $y'\not=\bot$, we have $x\meet\blacklozenge y'\not=\bot$. A key fact for proving the semantic equivalence of the resulting principal frame and the original $\mathcal{V}$-BAO is that complete additivity of a BAO implies the following condition: if $x\meet \blacklozenge y\not=\bot$, then there is a $y'\leq y$ such that $y'\not=\bot$ and $xRy'$ for the $R$ just defined. In fact, as observed in \citealt{Litak2015}, the ostensibly second-order condition of complete additivity is \textit{equivalent} to the first-order condition just given. In this way, thinking in terms of possibility frames has led to a new view of $\mathcal{V}$-BAOs. In addition, thinking in terms of possibility frames leads to a new characterization (in \S~\ref{MacNeille}) of the (lower) \textit{MacNeille completion} of a $\mathcal{V}$-BAO as in \citealt{Monk1970}. 

\begin{figure}
\begin{center}
\begin{tabular}{ccccc}

\textbf{$\mathcal{CV}$-BAO} & is dually equivalent to & \textbf{RichPoss} & which is a reflective subcategory of & \textbf{FullPoss} \\ &&&& \\
\textbf{BAO} & is dually equivalent to & \textbf{FiltPoss} & which is a reflective subcategory of & \textbf{Poss} \\

\end{tabular}
\end{center}
\caption{main categorical relationships.}\label{CatRel}
\end{figure}

In addition to this connection with complete additivity, we show that the special case of \textit{functional} principal possibility frames corresponds semantically to the class of what Litak~\citeyearpar{Litak2005} calls $\mathcal{T}$-BAOs (\S\S~\ref{PossToBAO}-\ref{VtoPossSection}). On this point, one of the interesting options that becomes available with the move from worlds to possibilities is a functional semantics for modalities according to which $\Box\varphi$ is true at a possibility $x$ iff $\varphi$ is true at the unique possibility $f(x)$ determined by an accessibility function $f$ (see \S~\ref{FuncFrames} and \citealt{Holliday2014}).  
  
If we do not require a full or principal possibility frame, then every BAO has an equivalent possibility frame, whose possibilities are \textit{proper filters} in the BAO (\S~\ref{GFPF}). In such a frame, a possibility $X$ refines a possibility $Y$ iff the proper filter $X$ is a superset of the proper filter $Y$. An accessibility relation $R$ on proper filters is given by a standard definition: $XRY$ iff for all elements $z$ in the algebra, $\blacksquare z\in X$ implies $z\in Y$, where $\blacksquare$ is the dual of the $\blacklozenge$ operator in the BAO. With this filter representation of BAOs, which does not require the ultrafilter axiom,\footnote{The ultrafilter axiom, stating that every nontrivial Boolean algebra contains an ultrafilter (equivalently, any filter disjoint from an ideal can be extended to an ultrafilter disjoint from the ideal), is weaker than the axiom of choice (\citealt{Halpern1971}) but still goes beyond ZF set theory (\citealt{Feferman1964}) as well as ZF with dependent choice (\citealt{Pincus1977a}).} we arrive at a notion of descriptive frames in the context of possibility semantics, which we call \textit{filter-descriptive} frames. Whereas Goldblatt \citeyearpar{Goldblatt1974} showed that the category \textbf{BAO} of BAOs with BAO-homomorphisms is dually equivalent to the category \textbf{DFm} of descriptive world frames with p-morphisms, we show that \textbf{BAO} is dually equivalent to the category \textbf{FiltPoss} of filter-descriptive possibility frames with p-morphisms (\S~\ref{Fdes}). In addition, parallel to Goldblatt's \citeyearpar{Goldblatt2006b} result that \textbf{DFm} is a reflective subcategory of the category of all general frames with what he calls \textit{modal maps}, we show that \textbf{FiltPoss} is a reflective subcategory of the category \textbf{Poss} of all possibility frames with what we call \textit{possibility morphisms} (see \S~\ref{morphisms}). Finally, our filter representation of a BAO leads to what can be considered a choice-free construction of the \textit{canonical extension} of a BAO (\S~\ref{MacNeille}).

\begin{figure}[h]
\begin{center}
\begin{tikzpicture}[->,>=stealth',shorten >=1pt,shorten <=1pt, auto,node
distance=2cm,thick,every loop/.style={<-,shorten <=1pt}]
\tikzstyle{every state}=[fill=gray!20,draw=none,text=black]

\node (CAV) at (0,-.5) {{$\mathcal{CAV}$-BAOs/}};
\node (Kripke) at (0,-1) {{Kripke frames$^\ast$}};
\node (CAVdot) at (0,-.25) {{}};

\node (CA) at (-4.5,2) {{$\mathcal{CA}$-BAOs/}};
\node (CAdot) at (-3.5,1.5) {{}};
\node (CAdotUp) at (-4.5,2.2) {{}};
\node (neigh1) at (-4.5,1.5) {{normal}};
\node (neigh2) at (-4.5,1) {{neighborhood}};
\node (neigh3) at (-4.5,.5) {{$\,$frames$^\ast$}};

\node (CV) at (0,2) {{{$\mathcal{CV}$-BAOs}/}};
\node (full) at (0,1.5) {{\textbf{full possibility frames}$^\dagger$}};
\node (CVdot) at (0,1.25) {{}};
\node (CVdotUp) at (0,2.2) {{}};

\node (AV) at (4.5,2) {{{$\mathcal{AV}$-BAOs/}}};
\node (dis1) at (4.5,1.5) {{discrete}};
\node (dis2) at (4.5,1) {{general}};
\node (dis3) at (4.5,.5) {{$\,$frames}$^\ast$};
\node (AVdot) at (3.5, 1.5) {{}};
\node (AVdotUp) at (4.5,2.2) {{}}; 

\node (BAOs) at (0,7) {{{BAOs}/}};
\node (BAOs) at (0,6.5) {{descriptive world frames$^\ast$}};
\node (BAOs) at (0,6) {{\textbf{filter-descriptive possibility frames}$^\ast$}};
\node (BAOdot) at (0,5.75) {{}};

\node (CBAOs) at (-3,4.5) {{{$\mathcal{C}$-BAOs}/}};
\node (CBAOs2) at (-3,4) {{\textbf{normal neighborhood version}}};
\node (CBAOs2) at (-3,3.5) {{\textbf{of full possibility frames}}};
\node (CBAOdot) at (-3,3.25) {{}};

\node (VBAOs) at (3,4.5) {{{$\mathcal{V}$-BAOs}/}};
\node (VBAOs2) at (3,4) {{\textbf{principal}}};
\node (VBAOs2) at (3,3.5) {{\textbf{possibility frames}}};
\node (VBAOdot) at (3,3.25) {{}};

\path (CAVdot) edge[->] node {{}} (CAdot);
\path (CAVdot) edge[->] node {{}} (CVdot);
\path (CAVdot) edge[->] node {{}} (AVdot);
\path (CAdotUp) edge[->] node {{}} (CBAOdot);
\path (CVdotUp) edge[->] node {{}} (VBAOdot);
\path (CVdotUp) edge[->] node {{}} (CBAOdot);
\path (VBAOs) edge[->] node {{}} (BAOdot);
\path (CBAOs) edge[->] node {{}} (BAOdot);
\path (AVdotUp) edge[->] node {{}} (VBAOdot);

\end{tikzpicture}
\end{center}
\caption{classes of BAOs and semantically equivalent frames---any BAO in the class before the / can be turned into a frame in the class after the / that validates the same formulas, and vice versa. A $\ast$ indicates there is also a dual equivalence between associated categories of BAOs and frames as described in the main text or references. See Thomason \citeyearpar{Thomason1975} on $\mathcal{CAV}$/Kripke, Do\v{s}en \citeyearpar{Dosen1989} on $\mathcal{CA}$/neighborhood, ten Cate \& Litak \citeyearpar{tenCate2007} on $\mathcal{AV}$/discrete, and Goldblatt \citeyearpar{Goldblatt1974} on BAO/descriptive world frames. The $\dagger$ indicates that we will prove a dual equivalence involving a reflective subcategory of the category of frames (see \S~\ref{DualEquiv}).}\label{BAOfig}
\end{figure}
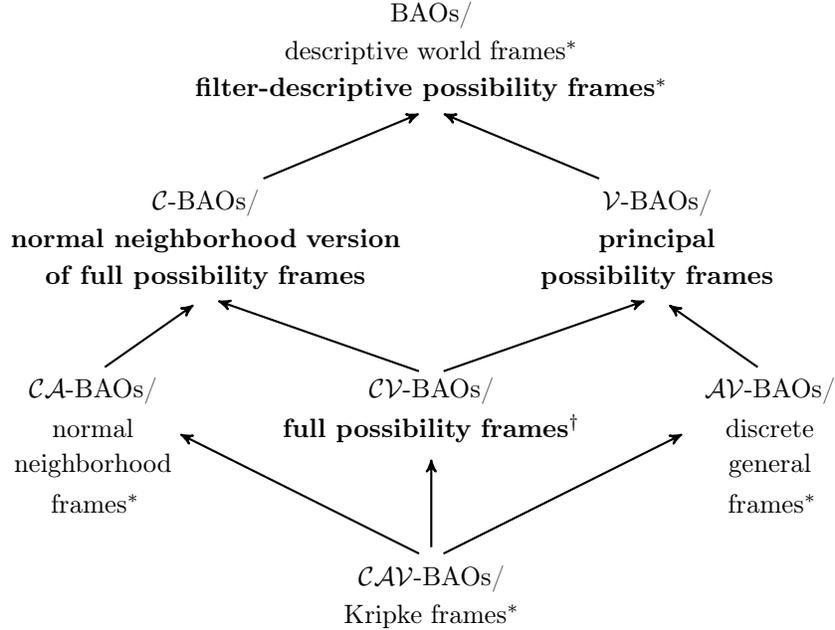

Putting together the links between BAOs and frames just described, we arrive at the pictures in Figures \ref{CatRel} and \ref{BAOfig}. Finally, in \S~\ref{OpFrameAlg} we compare possibility frame constructions that preserve the validity of modal formulas with algebraic constructions that preserve algebraic equations, in ways we exploit in \S\S~\ref{DefViaDual}-\ref{DefViaDual2}.

In \S~\ref{special}, we turn to definability and correspondence theory (cf. \citealt{Benthem1983,Benthem1980}). In \S\S~\ref{DefViaDual}-\ref{DefViaDual2}, we treat the question of which classes of possibility frames are definable by modal formulas. Using the duality theory of \S~\ref{DualityTheory}, we prove analogues for possibility semantics of Goldblatt's \citeyearpar{Goldblatt1974} characterization of modally definable classes of world frames and Goldblatt and Thomason's \citeyearpar{Goldblatt1975} characterization of modally definable classes of full world frames. In \S~\ref{LemmScottCorr}, we discuss the question of which modal formulas define classes of possibility frames that are first-order definable. As noted by Goldblatt \citeyearpar[p. 51]{Goldblatt2006}, ``a substantial reason for the great success'' of Kripke semantics is the way in which many natural modal axioms correspond to first-order properties of the accessibility relations in Kripke frames, such as seriality ($\Box\varphi\rightarrow\Diamond\varphi$), reflexivity ($\Box\varphi\rightarrow\varphi$), transitivity ($\Box\varphi\rightarrow\Box\Box\varphi$), symmetry ($\varphi\rightarrow\Box\Diamond\varphi$), and so on. Similarly, many natural modal axioms correspond to first-order properties of full possibility frames. As an illustration, we give an analogue for full possibility frames of one of the most elegant first-order correspondence results for Kripke frames, namely Lemmon and Scott's \citeyearpar[\S~4]{Lemmon1977} result for axioms of the form $\Diamond^k\Box^l p\rightarrow \Box^m\Diamond^np$. We also discuss the extent to which the ``minimal valuation'' heuristic for first-order correspondence in Kripke semantics applies in possibility semantics. More generally, using the connection between full possibility frames and $\mathcal{CV}$-BAOs (\S\S~\ref{PossToBAO}-\ref{DualEquiv}) and the methods of algebraic modal correspondence \citep{Conradie2014}, Yamamoto \citeyearpar{Yamamoto2016} establishes the possibility-semantic version of Sahlqvist correspondence (\citealt{Sahlqvist1975}, \citealt{Benthem1976}): every Sahlqvist modal formula corresponds to a first-order formula over full possibility frames. 

In \S~\ref{CompletenessTheory}, we draw out some consequences of \S~\ref{DualityTheory} for completeness theory. The discovery in the 1970s of \textit{Kripke-frame incompleteness}, the fact that not all normal modal logics are sound and complete with respect to a class of Kripke frames, has been considered one of the two forces that gave rise to the ``modern era'' of modal logic \citep[p. 44]{Blackburn2001}.\footnote{\label{IncompRefs}For a brief historical overview of Kripke-frame incompleteness results, see \citealt[\S~6.1]{Goldblatt2006}, and for textbook presentations, see, e.g., \citealt[Ch. 6]{Chagrov1997} or \citealt[Ch. 9]{Hughes1996}. Classic papers on Kripke-frame incompleteness include \citealt{Thomason1972,Thomason1974,Fine1974b,Benthem1978,Benthem1979}, and \citealt{Boolos1985}.}  By duality, this amounts to the fact that not every variety of BAOs is HSP-generated by the $\mathcal{CAV}$-BAOs that it contains. For a class $\mathcal{X}$ of BAOs or frames, let $\mathrm{ML}(\mathcal{X})$ be the set of modal logics \textbf{L} such that \textbf{L} is the logic of some subclass of $\mathcal{X}$. Let $\mathcal{ALG}$ be the class of all BAOs. Then $\mathrm{ML}(\mathcal{ALG})$ is the class of all normal modal logics (see \S~\ref{AlgSem}), and we have the following inequalities, the first three of which are from \citealt{Litak2005} and the last of which is from \citealt{Litak2015}:
\begin{equation}\mathrm{ML}(\mathcal{CAV})\subsetneq \mathrm{ML}(\mathcal{CV})\subsetneq \mathrm{ML}(\mathcal{T})\subsetneq \mathrm{ML}(\mathcal{V})\subsetneq\mathrm{ML}(\mathcal{ALG}).\label{Ineq}\end{equation}
Where \textsf{K} is the class of \textit{Kripke} frames, \textsf{FP} is the class of \textit{full} possibility frames, \textsf{PR} is the class of \textit{principal} possibility frames, \textsf{f-PR} is the class of \textit{functional} principal possibility frames, and \textsf{P} is the class of all possibility frames---or we could take just \textit{filter-descriptive} frames---it follows from (\ref{Ineq}) and the duality theory of \S~\ref{DualityTheory} that:
\[\begin{array}{ccccccccc}
\mathrm{ML}(\mathcal{CAV})&& \mathrm{ML}(\mathcal{CV})&& \mathrm{ML}(\mathcal{T})&& \mathrm{ML}(\mathcal{V})&&\mathrm{ML}(\mathcal{ALG})\\
 \rotatebox{90}{=} && \rotatebox{90}{=} && \rotatebox{90}{=} && \rotatebox{90}{=} &&\rotatebox{90}{=} \\
\mathrm{ML}(\mathsf{K})&\subsetneq& \mathrm{ML}(\mathsf{FP})&\subsetneq& \mathrm{ML}(\mbox{\textsf{f-PR}})&\subsetneq& \mathrm{ML}(\mathsf{PR})&\subsetneq&\mathrm{ML}(\mathsf{P}).\end{array}\]
Thus, every normal modal logic is sound and complete with respect to a class of filter-descriptive possibility frames, but the other types of possibility frames give rise to distinct notions of completeness. We stress the first and last inequalities above. The first means that there are \textit{Kripke-frame incomplete} normal modal logics that are sound and complete with respect to a class of \textit{full} possibility frames---indeed, we will show there are uncountably many such logics in \S~\ref{CompFull}.  The last inequality means that completeness with respect to \textit{principal} possibility frames is still an informative notion of completeness. To obtain better understandings of $\mathrm{ML}(\mathcal{CV})$ ($=\mathrm{ML}(\mathsf{FP})$) and $ \mathrm{ML}(\mathcal{V})$ ($=\mathrm{ML}(\mathsf{PR})$) is a major open problem in the theory of possibility semantics (see \S~\ref{OpenProb}). By contrast, $\mathrm{ML}(\mathcal{T})$ is relatively well understood. For instance, there is a known syntactic characterization of $\mathrm{ML}(\mathcal{T})$: a logic is in $\mathrm{ML}(\mathcal{T})$ iff its \textit{minimal tense extension} is a conservative extension. In \S~\ref{SyntacticProp}, we review other sufficient syntactic conditions for $\mathcal{T}$-completeness, as well as $\mathcal{V}$-completeness.

Another topic in completeness theory for possibility semantics, emphasized in \citealt{Humberstone1981} and \citealt{Holliday2014}, is completeness with respect to \textit{atomless} possibility frames, in which every possibility can be further refined. These frames have \textit{no worlds}. In \S~\ref{AtomlessFull}, we show that all consistent normal modal logics given by \textit{Sahlqvist} axioms are sound and complete with respect to an atomless, functional, full possibility frame. 

Finally, in \S~\ref{Canonical}, we discuss \textit{canonical} possibility frames for normal modal logics. Unlike in possible world semantics, where the domain of the canonical frame for a logic is the set of all maximally consistent sets of formulas (or ultrafilters in the Lindenbaum algebra for the logic), in possibility semantics we can take the domain of the canonical frame to be the set of all consistent and deductively closed sets of formulas (or proper filters in the Lindenbaum algebra). As an illustration, we prove that all normal modal logics given by Lemmon-Scott axioms are sound and complete with respect to their canonical full possibility frames.

In \S~\ref{Conclusion}, we conclude with a review of related work (\S~\ref{RelatedWork}) and a list of open problems (\S~\ref{OpenProb}).  Since this paper was submitted for publication, a number of papers on or related to possibility semantics, including papers focused on applications, have appeared. I have decided not to update \S~\ref{RelatedWork} with this more recent work, though I have added other references suggested by a referee. For a survey of classical possibility semantics including more recent work, see \citealt{Holliday2021}; for a topological, non-modal version of the duality in \S~\ref{Fdes}, see \citealt{BH2020}; and for non-classical generalizations of possibility semantics, see \citealt{Holliday2022,Holliday2023,Holliday2024}, \citealt{HM2024}, and \citealt{Massas2023}.

The appendices contain a review of Kripke semantics (\S~\ref{WRM}), general frame semantics (\S~\ref{GFS}), and algebraic semantics (\S~\ref{AlgSem}), followed by some deferred topics (\S~\ref{Deferred}), including a comparison of our definition of possibility frames with Humberstone's \citeyearpar{Humberstone1981} original definition and an intermediate alternative (\S~\ref{Strengths}).

Before embarking on the plan above, we fix our languages and logics in \S~\ref{LangLog} and our notation in \S~\ref{NotationSection}.

\subsection{Languages and Logics}\label{LangLog}

Possibility semantics naturally invites one to consider novel languages that cannot be treated by standard possible world semantics. However, in this study we restrict attention to the basic polymodal language, in order to see how possibility semantics and world semantics compare on a common playing field.

\begin{definition}[Modal Language]\label{language} Given a nonempty set $\sig$ of propositional variables and a set $\ind$ of modal operator indices, the language $\mathcal{L}(\sig,\ind)$ is generated by the following grammar:
\[\varphi::= p\mid \neg\varphi\mid (\varphi\wedge\varphi)\mid (\varphi\rightarrow\varphi)\mid \Box_i\varphi,\]
where $p\in \sig$ and $i\in \ind$. The language $\mathcal{L}(\sig,\emptyset)$ is the non-modal fragment of $\mathcal{L}(\sig,\ind)$.

As abbreviations, we define $(\varphi\vee\psi):=\neg (\neg\varphi\wedge\neg\psi)$, $(\varphi\leftrightarrow\psi):= ((\varphi\rightarrow\psi)\wedge (\psi\rightarrow\varphi))$, $\Diamond_i\varphi:=\neg\Box_i\neg\varphi$, $\bot:=(p\wedge\neg p)$ for some $p\in \sig $, and $\top:=\neg\bot$.

When there is no risk of ambiguity, we omit parentheses in formulas in standard fashion. \hfill $\triangleleft$
\end{definition}
 
As is common in the modal logic literature (e.g., \citealt{Chagrov1997,Blackburn2001}), we take a modal logic $\mathbf{L}$ for $\mathcal{L}(\sig,\ind)$ to be a subset of $\mathcal{L}(\sig,\ind)$ satisfying certain closure conditions. For familiar notation: if $\varphi\in\mathbf{L}$, then $\vdash_\mathbf{L}\varphi$; otherwise $\nvdash_\mathbf{L}\varphi$.

\begin{definition}[Classical Normal Modal Logic]\label{NML}  $\mathbf{L}\subseteq\mathcal{L}(\sig, \ind)$ is a \textit{classical normal modal logic} iff for all $\varphi,\psi\in\mathcal{L}(\sig, \ind)$:\footnote{The term `classical' has unfortunately also been used to mean that a logic $\mathbf{L}$ is \textit{congruential}: if $\vdash_\mathbf{L}\varphi\leftrightarrow\psi$, then $\vdash_\mathbf{L}\Box_i\varphi\leftrightarrow\Box_i\psi$ for all $i\in\ind$ \citep{Segerberg1971,Chellas1980}. We use `classical' to contrast with \textit{intuitionistic}, not non-congruential.}
\begin{enumerate}
\item Uniform Substitution -- if $\vdash_\mathbf{L}\varphi$ and $\psi$ is a substitution instance of $\varphi$, then $\vdash_\mathbf{L}\psi$;
\item\label{NMLTaut} Tautologies -- if $\varphi$ is a classical propositional tautology, then $\vdash_\textbf{L}\varphi$;
\item Modus Ponens -- if $\vdash_\textbf{L}\varphi$ and $\vdash_\textbf{L}\varphi\rightarrow\psi$, then $\vdash_\textbf{L}\psi$;
\item K axiom -- $\vdash_\textbf{L} \Box_i(\varphi\rightarrow\psi)\rightarrow (\Box_i\varphi\rightarrow\Box_i\psi)$ for all $i\in \ind$;
\item Necessitation -- if $\vdash_\textbf{L}\varphi$, then $\vdash_\textbf{L}\Box_i\varphi$ for all $i\in \ind$. \hfill $\triangleleft$
\end{enumerate}
\end{definition}
The smallest classical normal modal logic for a unimodal language, i.e., with $|\ind|=1$, is denoted by `\textbf{K}'. The smallest classical normal modal logic for a polymodal language with $|\ind|=n$ is usually denoted by `$\textbf{K}_n$'. Since the logic depends on both $\sig$ and $\ind$, the definition really defines a logic $\textbf{K}_{\sig,\ind}$, but to avoid proliferating subscripts, we will simply write `\textbf{K}', with $\sig$ and $\ind$ understood, and similarly for extensions of \textbf{K}. 

We will also mention \textit{intuitionistic} normal modal logic, for which we use the following language. 

\begin{definition}[Intuitionistic Modal Language]\label{IntLanguage} The language $\mathcal{L}'(\sig,\ind)$ is generated by the following grammar:
\[\varphi::= p\mid \neg\varphi\mid (\varphi\wedge\varphi)\mid (\varphi\rightarrow\varphi)\mid (\varphi\primvee\varphi)\mid  \Box_i\varphi,\]
where $p\in \sig$ and $i\in \ind$.  \hfill $\triangleleft$
\end{definition}
We think of $\primvee$ as the primitive intuitionistic disjunction, to be distinguished from the $\vee$ in Definition \ref{language} that we use to abbreviate $\neg (\neg\varphi\wedge\neg\psi)$. Also note that the intuitionistic diamond modality is often treated as a primitive, not definable in terms of $\neg$ and $\Box$, but our concern here is the $\Box$-only language $\mathcal{L}'(\sig,\ind)$. 

An intuitionistic normal modal logic for $\mathcal{L}'(\sig,\ind)$ is defined as in Definition \ref{NML}, except with $\mathcal{L}'(\sig,\ind)$ in place of $\mathcal{L}(\sig,\ind)$ and \textit{theorem of intuitionistic} (\textit{Heyting}) \textit{propositional calculus} in place of \textit{classical propositional tautology} in part \ref{NMLTaut}. We denote the smallest intuitionistic normal modal logic by `\textbf{HK}' \citep{Bozic1984} (called `\textbf{IntK}$_\Box$' in \citealt{Wolter1997}).

Whenever we use the term `normal modal logic' without specifying `classical' or `intuitionistic', the intended meaning is \textit{classical} normal modal logic.

\subsection{Notation}\label{NotationSection}

The following notation will be used throughout.

\begin{notation}[Posets and Relations]\label{notation} For a poset $\langle S,\sqsubseteq\rangle$ and $x,y\in S$:
\begin{enumerate}[label=\arabic*.,ref=\arabic*]
\item\label{notation1} $\mathord{\downarrow}x=\{x'\in S\mid x'\sqsubseteq x\}$;
\item $x\comp y$ iff $\exists z\in S$: $z\sqsubseteq x$ and $z\sqsubseteq y$ (``$x$ and $y$ are compatible'');
\item $x\incomp y$ iff \textit{not} $x\comp y$ (``$x$ and $y$ are incompatible'').
\end{enumerate}
For a binary relation $R$ on $S$, $X\subseteq S$, and $x\in S$:
\begin{enumerate}[label=\arabic*.,ref=\arabic*,resume]
\item $R[X]$ is the image of $X$ under $R$, i.e., $R[X]=\{y\in S\mid \exists x\in X\colon xRy\}$;
\item $R^{-1}[X]$ is the preimage of $X$ under $R$, i.e., $R^{-1}[X]=\{y\in S\mid \exists x \in X\colon yRx\}$; 
\item $R(x)=R[\{x\}]$.  \hfill $\triangleleft$
\end{enumerate}
\end{notation} 

\noindent Other new pieces of notation will be introduced as we go.

\section{From Partial-State Frames to Possibility Frames}\label{FromPartToPoss}

\subsection{Partial-State Frames and Semantics}\label{PSFramesSem}

We take the following structures as our starting point for the semantics of normal modal logics.

\begin{definition}[Partial-State Frames and Models]\label{PosetMod}
A \textit{partial-state frame} is a tuple $\mathcal{F}=\langle S, \sqsubseteq , \{R_i\}_{i\in \ind},\adm\rangle$ where: 
\begin{enumerate}[label=\arabic*.,ref=\arabic*]
\item $S$ is a nonempty set (the set of \textit{states});
\item $\sqsubseteq$ is a partial order on $S$ (the \textit{refinement} relation);
\item $R_i$ is a binary relation on $S$ (the \textit{$i$-accessibility} relation);
\item\label{PosetMod4} $\adm$ is a subset of $\wp(S)$ (the set of \textit{admissible propositions}) such that $\emptyset\in \adm$ and for all $X,Y\in \adm$:
\begin{enumerate}
\item $X\cap Y\in\adm$;
\item $X\supset Y=\{s\in S\mid \forall s'\sqsubseteq s\colon s'\in X\Rightarrow s'\in Y\}\in \adm$;
\item $\blacksquare_i Y=\{s\in S\mid R_i(s)\subseteq Y\}\in\adm$.
\end{enumerate}
\end{enumerate}
A partial-state \textit{model} $\mathcal{M}$ \textit{based on} $\mathcal{F}$ is a tuple $\mathcal{M}=\langle \mathcal{F},\pi\rangle$ where: 
\begin{enumerate}[label=\arabic*.,ref=\arabic*,resume]
\item $\pi\colon \sig\to\adm$ (an \textit{admissible valuation}). \hfill $\triangleleft$
\end{enumerate}
\end{definition}

We may abuse notation and write `$x\in\mathcal{M}$' to mean that $x\in S$ where $S$ is the set of states of $\mathcal{M}$.

\begin{remark}[Flipped Notation]\label{Flipped} For $x,y\in S$, we take `$x\sqsubseteq y$' to mean that the state $x$ is a \textit{refinement} or \textit{further specification} or \textit{extension} of the state $y$. 
 By contrast, Humberstone \citeyearpar{Humberstone1981} writes `$x\geqslant y$' for our $x\sqsubseteq y$.  Our notation agrees with common practice in set-theoretic forcing, where `$x\leq y$' usually means that $x$ is at least as strong of a forcing condition as $y$; but it is flipped relative to standard practice in intuitionistic semantics, which puts stronger states on the right of $\leq$ (an exception being \citealt{Dragalin1988}). It will prove convenient for us to go \textit{down} (left) rather than \textit{up} (right) for refinements, with more specific states below less specific ones, to match the standard practice of interpreting conjunction in a Boolean algebra of propositions as greatest \textit{lower} bound, so stronger propositions are below weaker ones. \hfill$\triangleleft$
\end{remark}

We relate the formal language $\mathcal{L}(\sig,\ind)$ (Definition \ref{language}) to partial-state frames and models as follows.

\begin{definition}[Partial-State Semantics]\label{pmtruth1} Given a partial-state model $\mathcal{M}=\langle S,\sqsubseteq,\{R_i\}_{i\in \ind},\pi\rangle$ with $x\in S$ and $\varphi\in\mathcal{L}(\sig,\ind)$,  define $\mathcal{M},x\Vdash\varphi$ (``$\varphi$ is forced at $x$ in $\mathcal{M}$'') as follows:
\begin{enumerate}[label=\arabic*.,ref=\arabic*]
\item\label{AtClause} $\mathcal{M},x\Vdash p$ iff $x\in \pi(p)$;
\item\label{NegClause} $\mathcal{M},x\Vdash\neg\varphi$ iff $\forall x'\sqsubseteq x$, $\mathcal{M},x'\nVdash\varphi$;
\item $\mathcal{M},x\Vdash \varphi\wedge\psi$ iff $\mathcal{M},x\Vdash\varphi$ and $\mathcal{M},x\Vdash\psi$;
\item\label{ImpClause} $\mathcal{M},x\Vdash \varphi\rightarrow\psi$ iff $\forall x'\sqsubseteq x$, if $\mathcal{M},x'\Vdash\varphi$ then $\mathcal{M},x'\Vdash\psi$;
\item $\mathcal{M},x\Vdash\Box_i\varphi$ iff $\forall y\in R_i(x)$, $\mathcal{M},y\Vdash\varphi$.
\end{enumerate}
The \textit{truth set} of $\varphi$ in $\mathcal{M}$ is the set $\llbracket \varphi\rrbracket^\mathcal{M}=\{x\in S\mid \mathcal{M},x\Vdash \varphi\}$. 

We also have the following derived notions, where $\mathcal{F}=\langle S,\sqsubseteq,\{R_i\}_{i\in\ind},\adm\rangle$ is a partial-state frame:
\begin{enumerate}[label=\arabic*.,ref=\arabic*,resume]
\item $\mathcal{F},x\Vdash\varphi$ iff  $\mathcal{M},x\Vdash\varphi$ for every model $\mathcal{M}$ based on $\mathcal{F}$;
\item $\mathcal{M}\Vdash\varphi$ ($\varphi$ is \textit{globally true in} $\mathcal{M}$) iff $\mathcal{M},x\Vdash\varphi$ for all $x\in S$;
\item $\mathcal{F}\Vdash\varphi$ ($\varphi$ is \textit{valid over} $\mathcal{F}$) iff $\mathcal{M}\Vdash\varphi$ for every model $\mathcal{M}$ based on $\mathcal{F}$ \\ (equivalently, $\mathcal{F},x\Vdash\varphi$ for all $x\in S$).
\end{enumerate}
Given a class $\mathsf{M}$ of partial-state models, $\varphi\in\mathcal{L}(\sig,\ind)$, and $\Sigma\subseteq\mathcal{L}(\sig,\ind)$: 
\begin{enumerate}[label=\arabic*.,ref=\arabic*,resume]
\item $\Vdash_\mathsf{M}\varphi$ ($\varphi$ is \textit{valid} over $\mathsf{M}$) iff for every $\mathcal{M}\in\mathsf{M}$, $\mathcal{M}\Vdash \varphi$;
\item $\Sigma$ is \textit{satisfiable} in \textsf{M} iff for some $\mathcal{M}\in\mathsf{M}$ and $x\in\mathcal{M}$, $\mathcal{M},x\Vdash \sigma$ for all $\sigma\in\Sigma$.
\end{enumerate}

 \noindent For a class $\mathsf{F}$ of partial-state \textit{frames}, validity ($\Vdash_\mathsf{F}$) and satisfiability with respect to $\mathsf{F}$ are defined as validity and satisfiability with respect to the class of all models based on frames in $\mathsf{F}$. 

Given a logic $\mathbf{L}\subseteq \mathcal{L}(\sig,\ind)$ and a class $\mathsf{C}$ of frames or of models:
\begin{enumerate}[label=\arabic*.,ref=\arabic*,resume]
\item $\mathbf{L}$ is \textit{sound} with respect to $\mathsf{C}$ iff for all $\varphi\in\mathcal{L}(\sig,\ind)$, if $\vdash_\mathbf{L}\varphi$, then $\Vdash_\mathsf{C} \varphi$;
\item $\mathbf{L}$ is \textit{complete} with respect to $\mathsf{C}$ iff for all $\varphi\in\mathcal{L}(\sig,\ind)$, if $\Vdash_\mathsf{C} \varphi$, then $\vdash_\mathbf{L}\varphi$. \hfill $\triangleleft$
\end{enumerate}
\end{definition}

The forcing clauses \ref{NegClause}-\ref{ImpClause} in Definition \ref{pmtruth1}, together with the following derived clause for $\vee$, are just the standard clauses for the connectives used in weak forcing in set theory (see, e.g., \citealt[\S~5.1.3]{Jech2008}).\footnote{This notion of forcing differs from Paul Cohen's original notion, which is sometimes called ``strong'' forcing. For the historical origins of weak forcing, owing to Solomon Feferman and Dana Scott, see \citealt[p.~160]{Moore1988}.} 

\begin{fact}[Forcing $\vee$]\label{ForcingDis} Given $\varphi\vee\psi:=\neg (\neg\varphi\wedge\neg\psi)$, $\mathcal{M},x\Vdash\varphi\vee\psi$ iff $\forall x'\sqsubseteq x$ $\exists x''\sqsubseteq x'$: $\mathcal{M},x''\Vdash\varphi$ or $\mathcal{M},x''\Vdash\psi$.
\end{fact}
\noindent The connection with set-theoretic forcing will reappear below (e.g., in Remark \ref{Persp2} and \S~\ref{DualityTheory}).

Together the closure conditions on the set $\adm$ of admissible propositions in Definition \ref{PosetMod}.\ref{PosetMod4} and the semantic clauses for the operators in Definition \ref{pmtruth1} guarantee that the truth set of every formula of $\mathcal{L}(\sig,\ind)$ is an admissible proposition, as in the following fact.

\begin{fact}[Truth Sets and Substitution]\label{TruthSub} For any partial-state frame $\mathcal{F}=\langle S, \sqsubseteq , \{R_i\}_{i\in \ind},\adm\rangle$ and model $\mathcal{M}=\langle \mathcal{F},\pi\rangle$ based on $\mathcal{F}$:
\begin{enumerate}
\item\label{TruthSub1} Where $\supset$ and $\blacksquare_i$ are the operations from Definition \ref{PosetMod}.\ref{PosetMod4}:
\begin{enumerate}
\item $\llbracket p\rrbracket^\mathcal{M}=\pi(p)$;
\item $\llbracket \neg\varphi\rrbracket^\mathcal{M}=\llbracket \varphi\rrbracket^\mathcal{M}\supset \emptyset$;
\item $\llbracket \varphi\wedge\psi\rrbracket^\mathcal{M}=\llbracket \varphi\rrbracket^\mathcal{M}\cap \llbracket \psi\rrbracket^\mathcal{M}$;
\item $\llbracket \varphi\rightarrow\psi\rrbracket^\mathcal{M}=\llbracket \varphi\rrbracket^\mathcal{M}\supset \llbracket\psi\rrbracket^\mathcal{M}$;
\item $\llbracket \Box_i\varphi\rrbracket^\mathcal{M}=\blacksquare_i\llbracket\varphi\rrbracket^\mathcal{M}$.
\end{enumerate}
\item\label{TruthSub2} for any $\varphi\in\mathcal{L}(\sig,\ind)$, $\llbracket \varphi\rrbracket^\mathcal{M}\in\adm$;
\item\label{TruthSub3} the set of formulas valid over $\mathcal{F}$ is closed under Uniform Substitution (Definition \ref{NML}).
\end{enumerate}
\end{fact}
\begin{proof} Part \ref{TruthSub1} is just a matter of checking the semantic clauses in Definition \ref{pmtruth1}. Part \ref{TruthSub2} is immediate from part \ref{TruthSub1} and Definition \ref{PosetMod}.\ref{PosetMod4}. Part \ref{TruthSub3} is a straightforward consequence of part \ref{TruthSub2}.
\end{proof}

We will now give several examples of partial-state frames. 

\begin{example}[World Frames]\label{KripkeExample} Classical \textit{Kripke frames} $\langle \wo{W},\{\wo{R}_i\}_{i\in\ind}\rangle$ (reviewed in Appendix \S~\ref{WRM}) may be regarded as partial-state frames $\langle \wo{W},\sqsubseteq, \{\wo{R}_i\}_{i\in\ind},P \rangle$ where:
\begin{enumerate}
\item $\sqsubseteq$ is the identity relation;
\item $P=\wp (\wo{W})$.\footnote{Alternatively, a Kripke frame may be regarded as the \textit{foundation}, as in Definition \ref{Foundation}, of such a partial-state frame.}
\end{enumerate}
For models based on Kripke frames, the forcing relation $\Vdash$ defined in Definition \ref{pmtruth1} is equivalent to the standard satisfaction relation, i.e., $\mathcal{M},x\Vdash \neg\varphi$ iff $\mathcal{M},x\nVdash \varphi$, and $\mathcal{M},x\Vdash \varphi\rightarrow\psi$ iff $\mathcal{M},x\nVdash \varphi$ or $\mathcal{M},x\Vdash\psi$. 

Classical \textit{general frames} $\langle \wo{W},\{\wo{R}_i\}_{i\in\ind},\wo{A}\rangle$ (reviewed in Appendix \S~\ref{GFS}) may also be regarded as partial-state frames $\langle \wo{W},\sqsubseteq, \{\wo{R}_i\}_{i\in\ind},\wo{A} \rangle$ in which $\sqsubseteq$ is the identity relation. When $\sqsubseteq$ is identity, the closure conditions on the set $\adm$ of admissible propositions in a partial-state frame coincide with the closure conditions on the set $\wo{A}$ of admissible propositions in a general frame (Definition \ref{GenSem}). As usual, we may view Kripke frames as \textit{full} general frames, i.e., general frames in which $\wo{A}=\wp(\wo{W})$.

Henceforth we will  call general frames \textit{world frames}. Kripke frames are then \textit{full} world frames. Models based on world frames will be called \textit{world models}. Note that we use the term `partial state' frame/model to cover world frame/models as well, just as `partial function' covers total functions. 
\hfill $\triangleleft$
\end{example}

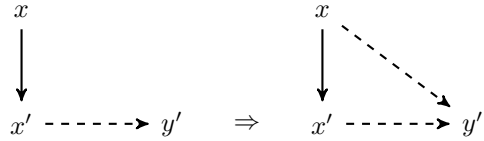
\begin{figure}[h]
\begin{center}
\begin{tikzpicture}[->,>=stealth',shorten >=1pt,shorten <=1pt, auto,node
distance=2cm,thick,every loop/.style={<-,shorten <=1pt}]
\tikzstyle{every state}=[fill=gray!20,draw=none,text=black]

\node (x-up) at (0,0) {{$x$}};
\node (x) at (0,-1.5) {{$x'$}};
\node (y) at (2,-1.5) {{$y'$}};

\node at (3,-1.5) {{\textit{$\Rightarrow$}}};

\path (x) edge[dashed,->] node {{}} (y);
\path (x-up) edge[->] node {{}} (x);

\node (x-up') at (4,0) {{$x$}};
\node (x') at (4,-1.5) {{$x'$}};
\node (y') at (6,-1.5) {{$y'$}};

\path (x') edge[dashed,->] node {{}} (y');
\path (x-up') edge[->] node {{}} (x');
\path (x-up') edge[dashed,->] node {{}} (y');

\end{tikzpicture}
\end{center}
\caption{the \upR{} condition from Example \ref{IntMod}. Given $x'R_iy'$, we may go \textit{up} in the first coordinate to any $x$ above $x'$ to obtain $xR_iy'$.  A solid arrow from $s$ to $t$ indicates that $t$ is a refinement of $s$ ($t\sqsubseteq s$), while a dashed arrow indicates that $t$ is accessible from $s$ ($sR_i t$).}\label{upRFig}
\end{figure}

\begin{figure}[h]
\begin{center}
\begin{tikzpicture}[->,>=stealth',shorten >=1pt,shorten <=1pt, auto,node
distance=2cm,thick,every loop/.style={<-,shorten <=1pt}]
\tikzstyle{every state}=[fill=gray!20,draw=none,text=black]

\node (x) at (0,-1.5) {{$x$}};
\node (y) at (2,-1.5) {{$y$}};
\node (y-down) at (2,-3) {{$y'$}};
\node at (3,-1.5) {{\textit{$\Rightarrow$}}};

\path (x) edge[dashed,->] node {{}} (y);

\path (y) edge[->] node {{}} (y-down);

\node (x') at (4,-1.5) {{$x$}};
\node (y') at (6,-1.5) {{$y$}};
\node (y-down') at (6,-3) {{$y'$}};

\path (x') edge[dashed,->] node {{}} (y');
\path (y') edge[->] node {{}} (y-down');
\path (x') edge[dashed,->] node {{}} (y-down');

\end{tikzpicture}
\end{center}
\caption{the \Rdown{} condition from Example \ref{IntMod}. Given $xR_iy$, we may go \textit{down} in the second coordinate to any $y'$ below $y$ to obtain $xR_iy'$.}\label{RdownFig}
\end{figure}
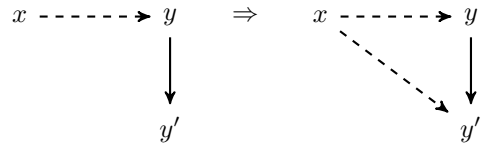

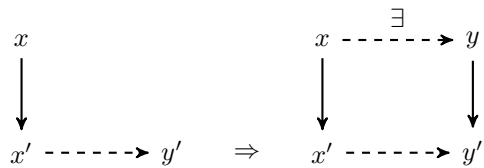
\begin{figure}[h]
\begin{center}
\begin{tikzpicture}[->,>=stealth',shorten >=1pt,shorten <=1pt, auto,node
distance=2cm,thick,every loop/.style={<-,shorten <=1pt}]
\tikzstyle{every state}=[fill=gray!20,draw=none,text=black]
\node (x-up) at (0,0) {{$x$}};
\node (x) at (0,-1.5) {{$x'$}};
\node (y) at (2,-1.5) {{$y'$}};

\node at (3,-1.5) {{\textit{$\Rightarrow$}}};

\path (x) edge[dashed,->] node {{}} (y);
\path (x-up) edge[->] node {{}} (x);

\node at (5,.3) {{$\exists$}};
\node (x-up') at (4,0) {{$x$}};
\node (x') at (4,-1.5) {{$x'$}};
\node (y') at (6,-1.5) {{$y'$}};
\node (y-up') at (6,0) {{$y$}};

\path (x') edge[dashed,->] node {{}} (y');
\path (x-up') edge[->] node {{}} (x');
\path (y') edge[<-] node {{}} (y-up');
\path (x-up') edge[dashed,->] node {{}} (y-up');

\end{tikzpicture}
\end{center}
\caption{the \Rcomm{} condition from Example \ref{IntMod}.}\label{RcommFig}
\end{figure}

\begin{example}[Intuitionistic Modal Frames]\label{IntMod} \textit{Full intuitionistic modal frames} as in \citealt{Wolter1997} (full $\Box$-IM frames) are partial-state frames $\mathcal{F}=\langle S,\sqsubseteq, \{R_i\}_{i\in \ind},P \rangle$ satisfying the following conditions:
\begin{enumerate}
\item \upR{} -- if $x'\sqsubseteq x$ and $x'R_iy'$, then $xR_iy'$ (see Figure \ref{upRFig}); 
\item \Rdown{} -- if $y'\sqsubseteq y$ and $xR_iy$, then $xR_iy'$ (see Figure \ref{RdownFig});
\item\label{Int2} $\adm$ is the set of \textit{all downsets} in $\langle S,\sqsubseteq\rangle$.\footnote{As usual, a downset in $\langle S,\sqsubseteq\rangle$ is an $X\subseteq S$ such that for all $x,x'\in S$, if $x\in X$ and $x'\sqsubseteq x$, then $x'\in X$.}
\end{enumerate}
The set of all downsets is closed under $\cap$ and $\supset$, as required by Definition \ref{PosetMod}.\ref{PosetMod4}; and the \upR{} condition guarantees that the set of all downsets is also closed under $\blacksquare_i$ from Definition \ref{PosetMod}.\ref{PosetMod4}. Thus, by Fact \ref{TruthSub}, the truth set of every formula will be a downset. This is the property of Persistence (or Heredity) that is necessary and sufficient for a partial-state frame to validate the intuitionistic principle $\varphi\rightarrow (\psi\rightarrow\varphi)$:
\begin{itemize}
\item Persistence: if $\mathcal{M},x\Vdash \varphi$ and $x'\sqsubseteq x$, then $\mathcal{M},x'\Vdash\varphi$.
\end{itemize}
For the $\Box_i$ case of the inductive proof of Persistence, if $\mathcal{M},x'\nVdash \Box_i\varphi$, so there is a $y'$ with $x'R_iy'$ and $\mathcal{M},y'\nVdash \varphi$, then by \upR{}, $x'\sqsubseteq x$ implies $xR_i y'$, so $\mathcal{M},x\nVdash \Box_i\varphi$. 

Having established Persistence using \upR{}, one can then motivate \Rdown{} as follows. The intended interpretation of $xR_iy'$ is that whatever is necessary/known/believed/henceforth true/etc. at $x$ is true at $y'$. If $\mathcal{M},x\Vdash \Box_i\varphi$ and $xR_iy$, so $\mathcal{M},y\Vdash\varphi$, then $y'\sqsubseteq y$ and Persistence together imply $\mathcal{M},y'\Vdash\varphi$. Thus, according to our interpretation of $R_i$, we should be able to have $xR_iy'$.

The condition \upR{} is not \textit{necessary} for Persistence. A weaker commutativity condition from \citealt{Bozic1984} (cf. \citealt{Celani1997}, \citealt[\S~11.2]{Restall2000}) is necessary and sufficient:
\begin{itemize}
\item \Rcomm{} -- if $x'\sqsubseteq x$ and $x'R_i y'$, then $\exists y$: $xR_iy$ and $y'\sqsubseteq y$ (see Figure \ref{RcommFig}).
\end{itemize}
Although \upR{} is not necessary, every partial-state frame $\mathcal{F}$ in which $\adm$ is a set of downsets may be turned into a semantically equivalent partial-state frame satisfying  \upR{} and \Rdown{} (see Proposition \ref{Representation}).

For the extended language $\mathcal{L}'(\sig,\ind)$ of Definition \ref{IntLanguage}, we extend the forcing relation $\Vdash$ from Definition \ref{pmtruth1} with the clause:
\begin{itemize}
\item $\mathcal{M},x\Vdash \varphi\primvee\psi$ iff $\mathcal{M},x\Vdash\varphi$ or $\mathcal{M},x\Vdash\psi$.
\end{itemize}
Since the set of all downsets is closed under unions, the Persistence property still holds for $\mathcal{L}'(\sig,\ind)$. If we go from \textit{full} intuitionistic modal frames to general intuitionistic modal frames for $\mathcal{L}'(\sig,\ind)$, then $\adm$ is any set of downsets containing $\emptyset$ and closed under $\cap$, $\supset$, $\blacksquare_i$, and $\cup$. These frames are clearly also partial-state frames.

(The above treatment of $\primvee$, as in intuitionistic Kripke semantics \citep{Kripke1965}, makes the forcing clauses for intuitionistic and classical disjunction different (recall Fact \ref{ForcingDis}). In \S~\ref{RelatedWork}, we discuss a more unified approach to intuitionistic and classical forcing based on Beth-style semantics from \citealt{Dragalin1988}.)
\hfill $\triangleleft$
\end{example} 

The logic of intuitionistic modal frames is the logic \textbf{HK} introduced in \S~\ref{LangLog}.

\begin{theorem}[\citealt{Bozic1984}]\label{Bozic}
 $\mathbf{HK}$ is sound with respect to the class of all intuitionistic modal frames and complete with respect to the class of full intuitionistic modal frames.
\end{theorem}

Our third example of partial-state frames will be quite important for our purposes. According to a world-based view of possibilities, a ``possibility'' is simply a \textit{set of worlds}; a possibility $X$ ``refines'' a possibility $Y$ iff $X\subseteq Y$; and however we define ``truth at a possibility,'' it should turn out that truth at a possibility is equivalent to truth at every world in that possibility (see \citealt{Cresswell2004} for this line of thinking, \S~\ref{RelatedWork} for other semantics that evaluate formulas at sets of worlds, and \citealt{Humberstone2019} for discussion of the ``flatness'' condition that truth at a set $X$ of worlds is equivalent to truth at all singleton subsets of $X$). Thus, where $\Vdash_p$ is a proposed forcing relation between possibilities and formulas, and $\vDash_w$ is the usual satisfaction relation between worlds and formulas (Definition \ref{pwmtruth1}), it should hold that for all possibilities $X$ and formulas~$\varphi$: 
\begin{equation}\mathcal{M},X\Vdash_p \varphi\mbox{ iff }\forall x\in X\colon \mathcal{M},\{x\}\Vdash_p \varphi\mbox{, where $\mathcal{M},\{x\}\Vdash_p\varphi$ iff $\mathcal{M},x\vDash_w \varphi$.}\label{WorldPoss}\end{equation}
Fact \ref{WtoP1} below records that the forcing relation $\Vdash$ from Definition \ref{pmtruth1} is indeed such a $\Vdash_p$ satisfying (\ref{WorldPoss}). First, let us make official the construction of possibilities out of arbitrary sets of worlds.\footnote{Compare the use of powerset (minus empty set) frames as intuitionistic frames for Medvedev's \citeyearpar{Medvedev1966} ``logic of finite problems''  and Skvortsov's \citeyearpar{Skvortsov1979} ``logic of infinite problems,'' as reviewed in, e.g., Definition 1 of  \citealt{Shatrov2008}.}
 
\begin{example}[Powerset Possibilization]\label{PowerPoss} Given a world frame $\mathfrak{F}=\langle \wo{W},\{\wo{R}_i\}_{i\in\ind},\wo{A}\rangle$ and a world model $\mathfrak{M}=\langle \mathfrak{F},\wo{V}\rangle$, the \textit{powerset possibilizations} of $\mathfrak{F}$ and $\mathfrak{M}$, $\mathfrak{F}^\pow=\langle S,\sqsubseteq,\{R_i\}_{i\in\ind},\adm\rangle$ and $\mathfrak{M}^\pow=\langle \mathfrak{F}^\pow,\pi\rangle$, are defined by:
\begin{enumerate}
\item $S=\wp(\wo{W})\setminus \{\emptyset\}$;
\item $X\sqsubseteq Y$ iff $X\subseteq Y$;
\item $XR_iY$ iff $Y\subseteq \mathrm{R}_i[X]$;
\item\label{PowerPoss4} $\adm = \{\mathord{\downarrow}X\mid X\in\wo{A}\}$;\footnote{Recall from Notation \ref{notation}.\ref{notation1} that $\mathord{\downarrow}X=\{X'\in S\mid X'\sqsubseteq X\}$.}
\item $\pi(p)=\{X\in S\mid X\subseteq\mathrm{V}(p)\}$.  
\end{enumerate}
One can check that since $\wo{A}$ satisfies the closure conditions required by a world frame (Definition \ref{GenSem}), the $\adm$ defined in part \ref{PowerPoss4} satisfies the closure conditions required by a partial-state frame. 

Note that if $\mathfrak{F}$ is a full world frame (Kripke frame), then $P=\{\mathord{\downarrow}X\mid X\in \wp (\wo{W})\}$.

A key fact about this construction is that the logic of the powerset possibilization $\mathfrak{F}^\pow$ coincides with the logic of the world frame $\mathfrak{F}$, as in Fact \ref{WtoP1}.\ref{WtoP1b}.\hfill $\triangleleft$
\end{example}

\begin{fact}[Powerset Possibilization]\label{WtoP1} For any world frame $\mathfrak{F}$ and world model $\mathfrak{M}=\langle \mathfrak{F},\wo{V}\rangle$:
\begin{enumerate} 
\item\label{WtoP1a} for any $X\in\mathfrak{M}^\pow$ and $\varphi\in\mathcal{L}(\sig,\ind)$, $\mbox{$\mathfrak{M}^\pow,X\Vdash\varphi$ iff  $\forall x\in X$: $\mathfrak{M},x\vDash\varphi$}$;
\item\label{WtoP1b} for any $\Sigma\subseteq\mathcal{L}(\sig,\ind)$, $\Sigma$ is satisfiable over $\mathfrak{F}^\pow$ iff $\Sigma$ is satisfiable over $\mathfrak{F}$.
\end{enumerate}
\end{fact}

\begin{proof} Part \ref{WtoP1a} is provable by an easy induction on $\varphi$. 

For part \ref{WtoP1b}, if $\Sigma$ is satisfied at a world $w$ in a world model $\mathfrak{M}$ based on $\mathfrak{F}$, then by part \ref{WtoP1a}, $\Sigma$ is satisfied at $\{w\}$ in $\mathfrak{M}^\pow$, which is easily seen to be an admissible model based on $\mathfrak{F}^\pow$. In the other direction, for any partial-state model $\mathcal{M}=\langle \mathfrak{F}^\pow,\pi\rangle$ based on $\mathfrak{F}^\pow$, the world model $\mathcal{M}^{-\pow}=\langle \mathfrak{F},\pi^{-\pow}\rangle$ defined by $w\in \pi^{-\pow}(p)$ iff $\{w\}\in\pi(p)$ is easily seen to be an admissible model based on $\mathfrak{F}$. Moreover,  $(\mathcal{M}^{-\pow})^\pow=\mathcal{M}$, so if $\Sigma$ is satisfied at a state $X$ in $\mathcal{M}$, then part \ref{WtoP1a} implies that there is an $x\in X$ such that $\Sigma$ is satisfied at $x$ in $\mathcal{M}^{-\pow}$.\end{proof}

From Fact \ref{WtoP1}.\ref{WtoP1b} and the fact that \textbf{K} is the logic of (full) world frames, we have the following.

\begin{corollary}[Logics of Powerset Possibilizations] \textbf{K} is sound with respect to the class of all powerset possibilizations of world frames and complete with respect to the class of powerset possibilizations of full world frames. Moreover, any normal modal logic that is sound and complete with respect to a class $\mathsf{F}$ of world frames, according to standard Kripke semantics (Definition \ref{RelWorld}), is also sound and complete with respect to the class of powerset possibilizations of frames from $\mathsf{F}$, according to partial-state semantics (Definition \ref{pmtruth1}). 
\end{corollary}

Whatever examples of partial-state frames we consider, the following properties are evident from the semantic clauses for $\rightarrow$ and $\Box_i$ and the definition of validity in Definition \ref{pmtruth1}.

\begin{fact}[Modus Ponens, K, and Necessitation]\label{Nec&K} For any class $\mathsf{F}$ of partial-state frames and $\varphi,\psi\in\mathcal{L}(\sig,\ind)$:
\begin{enumerate}
\item if $\Vdash_\mathsf{F}\varphi$ and $\Vdash_\mathsf{F}\varphi\rightarrow\psi$, then $\Vdash_\mathsf{F}\psi$;
\item if $\Vdash_\mathsf{F}\varphi$, then $\Vdash_\mathsf{F}\Box_i\varphi$.
\end{enumerate}
Moreover, if for each $\mathcal{F}=\langle S, \sqsubseteq , \{R_i\}_{i\in \ind},\adm\rangle\in \mathsf{F}$, each $X\in \adm$ is a downset in $\langle S,\sqsubseteq\rangle$, then:
\begin{enumerate}
\item[3.] $\Vdash_\mathsf{F} \Box_i(\varphi\rightarrow \psi)\rightarrow (\Box_i\varphi\rightarrow\Box_i\psi)$.
\end{enumerate}
\end{fact}
Classes of partial-state frames may differ with respect to their propositional logics and the extra modal principles they validate. The logic of all partial-state frames is a subintuitionistic modal logic (cf. \citealt{Restall1994,Celani2001}), but we will not go into the details. Our interest here is in partial-state semantics for \textit{classical} modal logic as in \S~\ref{ClassicalFrames}.
  
\subsection{Possibility Frames}\label{ClassicalFrames}

In this section, we present two ways of thinking about partial-state frames for classical modal logic, in Remarks \ref{Persp1} and \ref{Persp2}, both of which will converge in our definition of \textit{possibility frames} below.

\begin{remark}[Perspective 1 -- Persistence and Refinability]\label{Persp1} Over classical frames, the classical principle $\neg\neg\varphi\leftrightarrow\varphi$ must be valid, so $\neg\neg\varphi$ must be equivalent to $\varphi$.  Let us analyze both directions of this equivalence.

First, since $\varphi$ is to be a consequence of $\neg\neg\varphi$, any classical partial-state model is such that if $\mathcal{M},x\Vdash\neg\neg\varphi$, then $\mathcal{M},x\Vdash\varphi$, or equivalently, if $\mathcal{M},x\nVdash\varphi$, then $\mathcal{M},x\nVdash\neg\neg\varphi$, which is equivalent to:
\begin{itemize}
\item Refinability -- if $\mathcal{M},x\nVdash \varphi$, then $\exists x'\sqsubseteq x$: $\mathcal{M},x'\Vdash\neg\varphi$.
\end{itemize}
Note that if $\mathcal{M},x\nVdash\neg\varphi$, then $\exists x'\sqsubseteq x$: $\mathcal{M},x'\Vdash\varphi$, by the semantics of $\neg$. So the idea behind Refinability is this: if $\varphi$ is indeterminate at $x$, i.e., if $\mathcal{M},x\nVdash\varphi$ and $\mathcal{M},x\nVdash\neg\varphi$, then there is a refinement of $x$ that decides $\varphi$ negatively and a refinement of $x$ that decides $\varphi$ affirmatively. This gives us the \textit{classical view of possibilities}: indeterminacy with respect to $\varphi$ is equivalent to having refinements that determine $\varphi$ each way.

Second, since $\neg\neg\varphi$ is to be a consequence of $\varphi$, any classical partial-state model is such that if $\mathcal{M},x\Vdash\varphi$, then $\mathcal{M},x\Vdash \neg\neg\varphi$. Given Refinability, this condition is equivalent to the condition from Example \ref{IntMod}:
\begin{itemize}
\item Persistence: if $\mathcal{M},x\Vdash \varphi$ and $x'\sqsubseteq x$, then $\mathcal{M},x'\Vdash\varphi$.
\end{itemize}
 For if $\mathcal{M},x\Vdash\varphi$ but $\mathcal{M},x'\nVdash\varphi$, then by Refinability there is an $x''\sqsubseteq x'$ such that $\mathcal{M},x''\Vdash\neg\varphi$, which with $x''\sqsubseteq x$ from the transitivity of $\sqsubseteq$ implies $\mathcal{M},x\nVdash\neg\neg\varphi$. So for $\neg\neg\varphi$ to be a consequence of $\varphi$, Persistence is necessary. It is also easy to see that Persistence is sufficient for $\neg\neg\varphi$ to be a consequence of $\varphi$. 

Thus, in classical partial-state frames, every admissible proposition $X\in \adm$ will satisfy:
\begin{itemize}
\item \textit{persistence}: if $x\in X$ and $x'\sqsubseteq x$, then $x'\in X$;
\item \textit{refinability}: if $x\not\in X$, then $\exists x'\sqsubseteq x$ $\forall x''\sqsubseteq x'$, $x''\not\in X$. 
\end{itemize}
The first condition is just the condition that $\adm$ be a set of downsets, familiar from Example \ref{IntMod}. However, we cannot assume that $\adm$ contains \textit{all} downsets, but only those satisfying \textit{refinability}.\footnote{\label{CofinalFootnote}Our use of the term `refinability' comes from Humberstone \citeyearpar{Humberstone1981}. It is not to be confused with other uses of `refined' in modal logic, e.g., refined general frames as in \citealt[Def. 5.65]{Blackburn2001}. Van Benthem \citeyearpar{Benthem1981} uses `cofinality' for the same idea, but given our flipped perspective noted in Remark \ref{Flipped}, we would have to talk of `coinitiality'. Our \textit{refinability} is equivalent to: if $\mathord{\downarrow} x\cap X$ is a \textit{coinitial} subset of $\mathord{\downarrow}x$, or in the terms of set-theoretic forcing, a \textit{dense} subset of $\mathord{\downarrow}x$, then $x\in X$.}

The contrapositive form of \textit{refinability}, if $\forall x'\sqsubseteq x$ $\exists x''\sqsubseteq x'$: $x''\in X$, then $x\in X$, is often convenient to use. It is also useful to note that $X$ satisfying both \textit{persistence} and \textit{refinability} is equivalent to $X$ satisfying:
\begin{itemize}
\item $x\in X$ iff $\forall x'\sqsubseteq x$ $\exists x''\sqsubseteq x'$: $x''\in X$.
\end{itemize}
The conditions of \textit{persistence} and \textit{refinability} on admissible propositions are not only necessary but also sufficient for a partial-state frame to be classical. To see this, one can check that over frames satisfying those conditions, the axioms of some complete Hilbert-style calculus for classical propositional logic are valid. Instead, we will give a more illuminating proof in the style of Cohen \citeyearpar[p. 119]{Cohen1966} for the case where $\sig$ is countable, as in Lemma \ref{classicality} below. The idea is that every state $x$ belongs to a chain $x_0\sqsupseteq x_1\sqsupseteq x_2\dots$ that decides the truth value of every formula eventually, i.e., $\forall \varphi\in\mathcal{L}(\sig,\emptyset)$ $\exists k\in\mathbb{N}$: $\mathcal{M},x_k\Vdash\varphi$ or $\mathcal{M},x_k\Vdash\neg\varphi$, and we can read off a classical propositional valuation from such a chain. \hfill $\triangleleft$
\end{remark} 

\begin{lemma}[Persistence and Refinability imply Classicality]\label{classicality} Fix a language $\mathcal{L}(\sig,\ind)$ with $\sig$ countable. If $\mathcal{F}=\langle S,\sqsubseteq, \{R_i\}_{i\in\ind},\adm\rangle$ is a partial-state frame such that every $X\in\adm$ satisfies \textit{persistence} and \textit{refinability}, and if $\varphi\in \mathcal{L}(\sig,\emptyset)$ is a classical propositional tautology, then $\varphi$ is valid over $\mathcal{F}$; so by Fact \ref{TruthSub}.\ref{TruthSub3}, every substitution instance of $\varphi$ in $\mathcal{L}(\sig,\ind)$ is valid over $\mathcal{F}$.
\end{lemma}

\begin{proof} Suppose $\mathcal{F}$ is a frame as in the statement of the lemma, so any model based on $\mathcal{F}$ satisfies the properties of Persistence and Refinability from Remark \ref{Persp1}. First, we observe that over $\mathcal{F}$, $\varphi\rightarrow\psi$ is equivalent to $\neg (\varphi\wedge\neg\psi)$, so in the inductive proof below, we do not need a separate case for $\rightarrow$. That $\varphi\rightarrow\psi$ entails $\neg (\varphi\wedge\neg\psi)$ is obvious. In the other direction, suppose $\mathcal{M},x\nVdash \varphi\rightarrow\psi$, so there is some $y\sqsubseteq x$  such that $\mathcal{M},y\Vdash\varphi$ but $\mathcal{M},y\nVdash\psi$. Then by Refinability, there is some $z\sqsubseteq y$ such that $\mathcal{M},z\Vdash\neg\psi$. By Persistence, $\mathcal{M},y\Vdash\varphi$  implies $\mathcal{M},z\Vdash\varphi$, and by the transitivity of $\sqsubseteq$, $y\sqsubseteq x$ and $z\sqsubseteq y$ together imply $z\sqsubseteq x$. Then $z\sqsubseteq x$,  $\mathcal{M},z\Vdash\varphi$, and $\mathcal{M},z\Vdash\neg\psi$ imply $\mathcal{M},x\nVdash \neg (\varphi\wedge\neg\psi)$ by the semantic clauses for $\neg$ and $\wedge$.

Now suppose that $\varphi\in \mathcal{L}(\sig,\emptyset)$ is not valid over $\mathcal{F}$, so there is a model $\mathcal{M}$ based on $\mathcal{F}$ and an $x\in\mathcal{M}$ such that $\mathcal{M},x\nVdash\varphi$. It follows by Refinability that there is an $x'\sqsubseteq x$ such that $\mathcal{M},x'\Vdash \neg\varphi$. By the semantic clause for $\neg$, for any $y\in \mathcal{M}$ and $\psi\in\mathcal{L}(\sig,\ind)$ we can choose a $y^\psi\sqsubseteq y$ with $\mathcal{M},y^\psi\Vdash\psi$ or $\mathcal{M},y^\psi\Vdash\neg\psi$. Enumerating the formulas of $\mathcal{L}(\sig,\emptyset)$ as $\psi_1,\psi_2,\dots$, define a sequence $x_0,x_1,x_2,\dots$ such that $x_0 =x'$ and $x_{n+1}  = x_n^{\psi_{n+1}}$. Thus, $x_0\sqsupseteq x_1\sqsupseteq x_2\dots$ is a chain that decides every formula eventually. Define a propositional valuation $v\colon \sig\rightarrow\{0,1\}$ such that $v(p)=1$ if for some $k\in\mathbb{N}$, $\mathcal{M},x_k\Vdash p$, and $v(p)=0$ otherwise. Where $\overline{v}\colon\mathcal{L}(\sig,\emptyset)\rightarrow \{0,1\}$ is the usual classical extension of $v$, we claim that for all $\alpha\in\mathcal{L}(\sig,\emptyset)$:
\begin{equation}\overline{v}(\alpha)=1\mbox{ iff }\exists k\in\mathbb{N}\colon \mathcal{M},x_k\Vdash \alpha.\label{alphaEQ}\end{equation} For induction on $\alpha$, the base case of $p$ follows from the definition of $v$. For the $\neg$ case, simply observe that 
\begin{eqnarray}
\overline{v}(\neg\alpha)=1 & \Leftrightarrow & \overline{v}(\alpha)=0\quad\mbox{by definition of }\overline{v}\nonumber\\
 &\Leftrightarrow & \forall j\in\mathbb{N}\colon\mathcal{M},x_j\nVdash \alpha\quad\mbox{by the inductive hypothesis}\nonumber \\
&\Rightarrow &  \mathcal{M},x_{m}\Vdash\neg \alpha  \quad\mbox{where }\alpha=\psi_m\mbox{ in the enumeration}\nonumber\\
&\Rightarrow & \exists k\in\mathbb{N}\colon \mathcal{M},x_{k}\Vdash \neg \alpha\nonumber
\end{eqnarray}
and the last line implies the second by Persistence and the semantic clause for $\neg$ in Definition \ref{pmtruth1}. 

For the $\wedge$ case, simply observe that
\begin{eqnarray}
\overline{v}(\alpha\wedge\beta)=1 &\Leftrightarrow & \overline{v}(\alpha)=1\mbox{ and }\overline{v}(\beta)=1\quad\mbox{by definition of }\overline{v}\nonumber\\
&\Leftrightarrow & \exists j,j'\in\mathbb{N}\colon\mathcal{M},x_j\Vdash \alpha\mbox{ and }\mathcal{M},x_{j'}\Vdash \beta\quad\mbox{by the inductive hypothesis}\nonumber\\
&\Rightarrow &\mathcal{M},x_{\mathrm{max}(j,j')}\Vdash \alpha\mbox{ and }\mathcal{M},x_{\mathrm{max}(j,j')}\Vdash \beta\quad\mbox{by Persistence}\nonumber\\
&\Rightarrow &\mathcal{M},x_{\mathrm{max}(j,j')}\Vdash \alpha\wedge \beta\quad\mbox{by Definition \ref{pmtruth1}}\nonumber \\
&\Rightarrow &\exists k\in\mathbb{N}\colon\mathcal{M},x_k\Vdash \alpha\wedge \beta\nonumber 
\end{eqnarray}
and the last line implies the second. Thus, we have established (\ref{alphaEQ}). 

Since $\mathcal{M},x_0\Vdash \neg\varphi$, we have $\overline{v}(\neg\varphi)=1$ by (\ref{alphaEQ}), so $\varphi$ is not a classical tautology.\end{proof} 
 
Our second perspective on classical partial-state frames is a topological one, which is well-known in the literature on forcing in set theory (see, e.g., \citealt[\S~5.1.3]{Jech2008} and \citealt{Takeuti1973}).\footnote{The history of forcing in \citealt{Moore1988} (p. 163) attributes this topological perspective on forcing to Dana Scott.}

\begin{remark}[Perspective 2 -- Regular Open Truth Sets]\label{Persp2} Let $\mathcal{O}(S,\sqsubseteq)$ be the set of all downsets in the poset $\langle S,\sqsubseteq\rangle$. Then $\mathcal{O}(S,\sqsubseteq)$ is an open set topology on $S$, the \textit{downset topology} or \textit{Alexandrov topology} (though the latter term is also often used for the topology of upsets rather than downsets); and by Remark \ref{Persp1}, $\adm\subseteq \mathcal{O}(S,\sqsubseteq)$ for any classical partial-state frame $\mathcal{F}=\langle S, \sqsubseteq , \{R_i\}_{i\in \ind},\adm\rangle$. Since for any downset in $\langle S,\sqsubseteq\rangle$, its complement is an upset, and vice versa, the closed sets in $\mathcal{O}(S,\sqsubseteq)$ are just the upsets in $\langle S,\sqsubseteq\rangle$. Thus, the \textit{closure} $\mathrm{cl}(X)$ of a set $X$, the smallest upset that includes $X$, is the set $\mathord{\Uparrow}X=\{y\in S\mid \exists x\sqsubseteq y\colon x\in X\}$ of all states with refinements in $X$; and the \textit{interior} $\mathrm{int}(X)$ of $X$, the largest downset  that is included in $X$, is the set $\{y\in S\mid \forall x\sqsubseteq y \colon x\in X\}$ ($=S\setminus \mathord{\Uparrow} (S\setminus X)$) of all states all of whose refinements are in $X$. (Note that in what follows, we sometimes apply the operations $\mathrm{int}$ and $\mathrm{cl}$ to sets $X\subseteq S$ that are not downsets.) From this perspective, we can rewrite the $\neg$ and $\rightarrow$ clauses of Definition \ref{pmtruth1} equivalently as:
\begin{itemize}
\item $\llbracket \neg\varphi\rrbracket^\mathcal{M}=\mathrm{int}(S\setminus \llbracket \varphi\rrbracket^\mathcal{M})$;
\item $\llbracket \varphi\rightarrow\psi\rrbracket^\mathcal{M}= \mathrm{int}((S\setminus\llbracket \varphi\rrbracket^\mathcal{M})\cup \llbracket \psi\rrbracket^\mathcal{M})$.
\end{itemize}
Since the interior of the complement is the complement of the closure, $\mathrm{int}(S\setminus \llbracket \varphi\rrbracket^\mathcal{M})=S\setminus \mathrm{cl}(\llbracket\varphi\rrbracket^\mathcal{M})$. So by the $\neg$ clause,  $\llbracket \neg\neg\varphi\rrbracket^\mathcal{M} = \mathrm{int}(S\setminus (S\setminus \mathrm{cl}(\llbracket \varphi\rrbracket^\mathcal{M})))=\mathrm{int}(\mathrm{cl}(\llbracket \varphi\rrbracket^\mathcal{M}))$. Thus, the classical requirement that $\llbracket \varphi\rrbracket^\mathcal{M}=\llbracket \neg\neg\varphi\rrbracket^\mathcal{M}$ is equivalent to the requirement that for all admissible propositions $X\in\adm$, $X=\mathrm{int}(\mathrm{cl}(X))$.  

Also note that from this perspective, the clause for $\vee$ from Fact \ref{ForcingDis} is equivalent to:
\begin{itemize}
\item $\llbracket \varphi\vee\psi\rrbracket^\mathcal{M}=\mathrm{int}(\mathrm{cl}(\llbracket \varphi \rrbracket^\mathcal{M}\cup \llbracket\psi\rrbracket^\mathcal{M}))$.
\end{itemize}

Open sets $X$ with the property that $X=\mathrm{int}(\mathrm{cl}(X))$ are called \textit{regular open}. What the Refinability principle of Remark \ref{Persp1} adds to Persistence topologically is the idea that the admissible propositions in $\adm$ should be not only open sets in $\mathcal{O}(S,\sqsubseteq)$, as per Persistence, but \textit{regular open} sets in $\mathcal{O}(S,\sqsubseteq)$. 
 
For any partial-state frame $\mathcal{F}=\langle S,\sqsubseteq,\{R_i\}_{i\in\ind},\adm\rangle$, the condition that $\adm$ be a set of regular open sets in $\mathcal{O}(S,\sqsubseteq)$ (still satisfying the required closure under $\cap$, $\supset$, and $\blacksquare_i$ from Definition \ref{PosetMod}) is not only necessary but also sufficient for $\mathcal{F}$ to be a classical frame. As observed by Tarski \citeyearpar{Tarski1937,Tarski1938,Tarski1956}, for any topological space $\mathcal{S}=\langle S,\mathcal{O}\rangle$, the structure $\langle \mathrm{RO}(\mathcal{S}), \meet, -, \top\rangle$ where $\mathrm{RO}(\mathcal{S})$ is the set of \textit{all} regular open sets in $\mathcal{S}$,  $X\meet Y= X\cap Y$, $-X=\mathrm{int}(S\setminus X)$, and $\top=S$ is a complete Boolean algebra, with the meet of an arbitrary $\mathcal{X}\subseteq\mathrm{RO}(\mathcal{S})$ given by $\bigmeet \mathcal{X}=\mathrm{int}(\bigcap \mathcal{X})$ and the join by $\bigvee\mathcal{X}=\mathrm{int}(\mathrm{cl}(\bigcup\mathcal{X}))$.\footnote{Note that where $\mathcal{O}$ is the Alexandrov topology $\mathcal{O}(S,\sqsubseteq)$, we have  $\bigmeet \mathcal{X}=\bigcap \mathcal{X}$.} This is called the \textit{regular open algebra} of $\mathcal{S}$. Another way to the same point is that the \textit{open} sets ordered by inclusion form a complete \textit{Heyting} algebra, and the \textit{regular} elements of a Heyting algebra, those $X$ for which $\mathord{-}\mathord{-}X=X$, form a complete Boolean algebra (cf. \citealt{Glivenko1929}). Now even if $\adm$ in $\mathcal{F}=\langle S,\sqsubseteq,\{R_i\}_{i\in\ind},\adm\rangle$ does not include \textit{all} regular open sets in $\mathcal{O}(S,\sqsubseteq)$, the closure conditions on $\adm$ from Definition \ref{PosetMod} ensure that $\adm$ is closed under the operations $\meet$ and $-$ just defined, and $S\in \adm$. Thus, if $\mathcal{F}=\langle S,\sqsubseteq,\{R_i\}_{i\in\ind},\adm\rangle$ is a partial-state frame such that all $X\in\adm$ are regular open, then $\langle \adm,\meet, -,\top\rangle$ is a \textit{subalgebra} of the regular open algebra arising from $\mathcal{O}(S,\sqsubseteq)$, so $\langle \adm,\meet, -,\top\rangle$ is a Boolean algebra (though not necessarily a complete Boolean algebra).

From these observations it is a short step, which we leave to the reader, to Lemma \ref{RegClass}.\hfill $\triangleleft$
\end{remark}

\begin{lemma}[Regular Opens and Classicality]\label{RegClass} If $\mathcal{F}=\langle S,\sqsubseteq,\{R_i\}_{i\in\ind},\adm\rangle$ is a partial-state frame in which every $X\in\adm$ is a regular open set in $\mathcal{O}(S,\sqsubseteq)$, and if $\varphi\in\mathcal{L}(\sig,\emptyset)$ is a classical propositional tautology, then $\varphi$ is valid over $\mathcal{F}$; so by Fact \ref{TruthSub}.\ref{TruthSub3}, every substitution instance of $\varphi$ in $\mathcal{L}(\sig,\ind)$ is valid over $\mathcal{F}$.
\end{lemma}

The perspectives of Remarks \ref{Persp1} and \ref{Persp2} come together easily in part \ref{RefReg3} of the following fact. For part \ref{RefReg2.5}, let $\mathord{\Downarrow}X=\{y\in S\mid \exists x\in X\colon y\sqsubseteq x\}$, so $\mathord{\Downarrow}X=X$ iff $X$ is a downset.

\begin{fact}[Regular Opens, Persistence, and Refinability]\label{RefReg} For any poset $\langle S,\sqsubseteq\rangle$ and $X\subseteq S$:
\begin{enumerate} 
\item\label{RefReg1} $\mathrm{int}(\mathrm{cl}(X))=\{x\in S\mid \forall x'\sqsubseteq x\,\exists x''\sqsubseteq x'\colon x''\in X\}$;
\item\label{RefReg2.5} $\mathrm{int}(\mathrm{cl}(\mathord{\Downarrow}X))$ is the smallest regular open set that includes $X$;
\item\label{RefReg3} $X$ satisfies \textit{persistence} and \textit{refinability} as in Remark \ref{Persp1} iff $X$ is regular open in $\mathcal{O}(S,\sqsubseteq)$.
\end{enumerate}
\end{fact}
\begin{proof} Part \ref{RefReg1} is immediate from the first-order definitions of $\mathrm{int}(X)$ and $\mathrm{cl}(X)$ in Remark \ref{Persp2}; part \ref{RefReg2.5} follows straightforwardly from part \ref{RefReg1}; and part \ref{RefReg3} follows from part \ref{RefReg1} and the combined formulation of \textit{persistence} and \textit{refinability} in Remark \ref{Persp1}.
\end{proof}

To get a sense of the constraint that admissible propositions must be regular open sets, it helps to consider some simple examples, such as the following.
 
\begin{example}[Beth Comb]\label{BethEx} Consider the infinite poset $\langle S,\sqsubseteq\rangle$ depicted three times in Figure \ref{Comb}, dubbed the \textit{Beth comb} due to its use in Beth semantics for intuitionistic logic \citep{Beth1956} (see \citealt[p. 898]{Humberstone2011}). As usual, arrows implied by reflexivity and transitivity are omitted. Let us call the $w_n$ states \textit{worlds} and the $s_n$ states \textit{non-worlds}. In addition to $S$ and $\emptyset$, there are two kinds of regular open sets in the downset topology on the Beth comb: 
\begin{itemize}
\item[(i)] if a downset $X$ contains only worlds, then $X$ is regular open iff  $\forall w_k\in X$ $\exists l>k$: $w_l\not\in X$.
\item[(ii)] if a downset $X$ contains a non-world, then where $n=\min\{k\mid s_k\in X\}$, $X$ is regular open iff $\mathord{\downarrow}s_n\subseteq X$ and $w_{n-1}\not\in X$.
\end{itemize}
Thus, the infinite downsets indicated by the ellipses in Figure \ref{Comb} are \textit{not} regular open. The smallest regular open set including the set in the middle (resp. the right) of Figure \ref{Comb} is  $S$ (resp. $\mathord{\downarrow}s_1$).\hfill$\triangleleft$

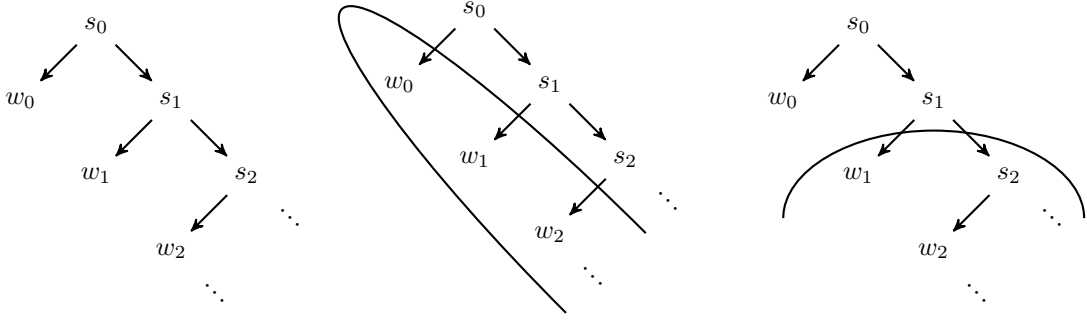
\begin{figure}[h]
 \begin{center}
\begin{tikzpicture}[-,>=stealth',shorten >=1pt,shorten <=1pt, auto,node
distance=2cm,thick,every loop/.style={<-,shorten <=1pt}]
\tikzstyle{every state}=[fill=gray!20,draw=none,text=black]

\node (x_0) at (0,0) {{$s_0$}};
\node (w_0) at (-1,-1) {{$w_0$}};
\node (x_1) at (1,-1) {{$s_1$}};
\node (w_1) at (0,-2) {{$w_1$}};
\node (x_2) at (2,-2) {{$s_2$}};
\node (w_2) at (1,-3) {{$w_2$}};

\node (dots) at (2.5,-2.5) {{\rotatebox[origin=c]{45}{$\boldsymbol{\vdots}$}}};

\node (dots) at (1.5,-3.5) {{\rotatebox[origin=c]{45}{$\boldsymbol{\vdots}$}}};

\path (x_0) edge[->] node {{}} (w_0);
\path (x_0) edge[->] node {{}} (x_1);
\path (x_1) edge[->] node {{}} (w_1);
\path (x_1) edge[->] node {{}} (x_2);
\path (x_2) edge[->] node {{}} (w_2);
\end{tikzpicture}
\begin{tikzpicture}[-,>=stealth',shorten >=1pt,shorten <=1pt, auto,node
distance=2cm,thick,every loop/.style={<-,shorten <=1pt}]
\tikzstyle{every state}=[fill=gray!20,draw=none,text=black]

\node (x_0) at (0,0) {{$s_0$}};
\node (w_0) at (-1,-1) {{$w_0$}};
\node (x_1) at (1,-1) {{$s_1$}};
\node (w_1) at (0,-2) {{$w_1$}};
\node (x_2) at (2,-2) {{$s_2$}};
\node (w_2) at (1,-3) {{$w_2$}};

\node (dots) at (2.5,-2.5) {{\rotatebox[origin=c]{45}{$\boldsymbol{\vdots}$}}};

\node (dots) at (1.5,-3.5) {{\rotatebox[origin=c]{45}{$\boldsymbol{\vdots}$}}};

\path (x_0) edge[->] node {{}} (w_0);
\path (x_0) edge[->] node {{}} (x_1);
\path (x_1) edge[->] node {{}} (w_1);
\path (x_1) edge[->] node {{}} (x_2);
\path (x_2) edge[->] node {{}} (w_2);
\draw[rotate=45] (-2,-3.75) arc (180:0: .75cm and 5cm);

\end{tikzpicture}
\begin{tikzpicture}[-,>=stealth',shorten >=1pt,shorten <=1pt, auto,node
distance=2cm,thick,every loop/.style={<-,shorten <=1pt}]
\tikzstyle{every state}=[fill=gray!20,draw=none,text=black]

\node (space) at (-2,0) {{$\,$}};
\node (x_0) at (0,0) {{$s_0$}};
\node (w_0) at (-1,-1) {{$w_0$}};
\node (x_1) at (1,-1) {{$s_1$}};
\node (w_1) at (0,-2) {{$w_1$}};
\node (x_2) at (2,-2) {{$s_2$}};
\node (w_2) at (1,-3) {{$w_2$}};

\node (dots) at (2.5,-2.5) {{\rotatebox[origin=c]{45}{$\boldsymbol{\vdots}$}}};

\node (dots) at (1.5,-3.5) {{\rotatebox[origin=c]{45}{$\boldsymbol{\vdots}$}}};

\path (x_0) edge[->] node {{}} (w_0);
\path (x_0) edge[->] node {{}} (x_1);
\path (x_1) edge[->] node {{}} (w_1);
\path (x_1) edge[->] node {{}} (x_2);
\path (x_2) edge[->] node {{}} (w_2);

\draw (-1,-2.6) arc (180:00: 2cm and 1.2cm);
\end{tikzpicture}

\end{center}
\caption{the Beth comb (left) and two examples of non-regular-open sets (middle and right).}\label{Comb}
\end{figure}
\end{example}

Since we will refer to sets as in Fact \ref{RefReg} so frequently, we use the following notation.

\begin{notation}[$\mathrm{RO}$]\label{ROnotation} For a poset $\langle S,\sqsubseteq\rangle$, let $\mathrm{RO}(S,\sqsubseteq)$ be the set of all $X\subseteq S$ satisfying \textit{persistence} and \textit{refinability}, or equivalently, the set of all regular open sets in $\mathcal{O}(S,\sqsubseteq)$. 

For a partial-state frame $\mathcal{F}=\langle S,\sqsubseteq,\{R_i\}_{i\in\ind},\adm\rangle$, let $\mathrm{RO}(\mathcal{F})=\mathrm{RO}(S,\sqsubseteq)$.\hfill $\triangleleft$
\end{notation}

Let us summarize the perspectives on classicality from this section with the following proposition. 

\begin{proposition}[Characterizations of Classicality]\label{characterizations} For any partial-state frame $\mathcal{F}=\langle S, \sqsubseteq , \{R_i\}_{i\in \ind},\adm\rangle$, the following are equivalent:
\begin{enumerate}
\item\label{Char1} the set of $\varphi\in\mathcal{L}(\sig,\ind)$ valid over $\mathcal{F}$ is a classical normal modal logic as in Definition \ref{NML};
\item\label{Char2} for every $\varphi\in\mathcal{L}(\sig,\ind)$, $\neg\neg\varphi$ is equivalent to $\varphi$ over $\mathcal{F}$;
\item\label{Char4} $\adm\subseteq \mathrm{RO}(\mathcal{F})$.
\end{enumerate}
\end{proposition}
\begin{proof} Part \ref{Char1} obviously implies part \ref{Char2}, and we observed in Remarks \ref{Persp1} and \ref{Persp2} that part \ref{Char2} is equivalent to part \ref{Char4}. Thus, to complete the proof, it suffices to show that part \ref{Char4} implies part \ref{Char1}. By Fact \ref{TruthSub} and \ref{Nec&K}, the set of formulas valid over $\mathcal{F}$ is closed under Uniform Substitution, Modus Ponens, and Necessitation, and contains the K axiom; and by Lemma \ref{RegClass}, for a partial-state frame $\mathcal{F}$ satisfying \ref{Char4}, the set of formulas valid over $\mathcal{F}$ contains all classical propositional tautologies. Thus, the set of formulas valid over such a frame is a classical normal modal logic as in Definition \ref{NML}.\end{proof}

Proposition \ref{characterizations} motivates the following central definition of the paper.

\begin{definition}[Possibility Frames]\label{PosFrames} A \textit{possibility frame} is a partial-state frame $\mathcal{F}=\langle S, \sqsubseteq , \{R_i\}_{i\in\ind},\adm\rangle$ as in Definition \ref{PosetMod} in which $\adm\subseteq \mathrm{RO}(\mathcal{F})$. A \textit{full} possibility frame is a possibility frame in which $\adm=\mathrm{RO}(\mathcal{F})$.\hfill $\triangleleft$
\end{definition}

Full possibility frames are the possibility semantic analogue of Kripke frames, which are \textit{full} world frames (recall Example \ref{KripkeExample} or see Appendix \S~\ref{GFS}). As previewed in \S~\ref{intro}, in \S~\ref{NoKripke}, \S~\ref{VtoPossSection}, and \S~\ref{CompFull} we will show that full possibility frames provide a more general semantics for normal modal logics than Kripke frames; and in \S~\ref{GFPF}, we will show that general possibility frames provide a fully general semantics. 

Two of our three examples of partial-state frames from \S~\ref{PSFramesSem} are examples of possibility frames.

\begin{example}[World Frames Cont.]\label{KripkeAgain} One can check that every world frame, viewed as a partial-state frame as in Example \ref{KripkeExample}, is a possibility frame, and every full world frame is a full possibility frame. \hfill $\triangleleft$
\end{example}

\begin{example}[Powerset Possibilization Cont.]\label{transform} One can check that for any world frame $\mathfrak{F}$, its powerset possibilization $\mathfrak{F}^\pow$ as in Example \ref{PowerPoss} is a possibility frame, and if $\mathfrak{F}$ is a full world frame, then $\mathfrak{F}^\pow$ is a full possibility frame (cf. Facts \ref{transform2} and \ref{FullTFAE}). \hfill $\triangleleft$
\end{example} 

Recall that the operation of powerset possibilization captures the view that the space of possibilities is the space of nonempty sets of worlds. This view builds in strong assumptions about the nature of possibilities, e.g., that the space of possibilities ordered by refinement has the structure of a \textit{complete} and \textit{atomic} Boolean lattice (minus the minimum element). Using these properties and one other property concerning the interplay of accessibility and refinement, we can characterize the class of possibility frames that are isomorphic to the powerset possibilization of some Kripke frame, as in \S~\ref{RichFrames}. In \S~\ref{RichFrames} we will also explain that any \textit{full} possibility frame can be transformed into a semantically equivalent possibility frame that shares several of the properties of powerset possibilizations, with the exception of \textit{atomicity}. In \S\S~\ref{VtoPossSection}-\ref{DualEquiv}, we will show that it is precisely this difference with respect to atomicity that makes full possibility frames more general than Kripke frames and their equivalent powerset possibilizations. In \S~\ref{PrincFrames} and \S~\ref{VtoPossSection}, we will generalize even further away from powerset possibilizations with our \textit{principal} possibility frames, which drop \textit{lattice-completeness} as well.

Using Example \ref{KripkeAgain} or Example \ref{transform}, we get an easy completeness proof for \textbf{K}.

\begin{corollary}[Soundness and Completeness of \textbf{K}]\label{Ksound2} \textbf{K} is sound with respect to the class of all possibility frames and complete with respect to the class of full possibility frames.
\end{corollary}

\begin{proof}
Soundness is given by Proposition \ref{characterizations}. Completeness follows from the completeness of \textbf{K} with respect to the class of full world frames together with Example \ref{KripkeAgain}, or Example \ref{transform} and Fact \ref{WtoP1}.
\end{proof}

An important property of a \textit{full} possibility frame $\mathcal{F}$ is that $\mathrm{RO}(\mathcal{F})$ is closed under $\blacksquare_i$ from Definition \ref{PosetMod}. This follows from the requirement of a partial-state frame that $\adm$ be closed under $\blacksquare_i$ plus the requirement of a full possibility frame that $\adm=\mathrm{RO}(\mathcal{F})$; and this is not trivial, for there are possibility frames $\mathcal{F}$ that lack the property. By contrast, it is easy to check that for any $\mathcal{F}$, $\mathrm{RO}(\mathcal{F})$ is closed under $\cap$ and $\supset$. 

The fact that not every possibility frame is such that $\mathrm{RO}(\mathcal{F})$ is closed under $\blacksquare_i$ means that not every possibility frame can be turned into a full possibility frame simply by replacing its set of admissible propositions $\adm$ by $\mathrm{RO}(\mathcal{F})$. It is helpful to put this point in terms of the following terminology and notation.

\begin{definition}[Foundations of Frames and Associated Full Frames]\label{Foundation} A \textit{foundation} is a tuple $F=\langle S,\sqsubseteq, \{R_i\}_{i\in\ind}\rangle$ where $\langle S,\sqsubseteq\rangle$ is a nonempty poset and each $R_i$ is a binary relation on $S$. 

Given a possibility frame $\mathcal{F}=\langle S,\sqsubseteq,\{R_i\}_{i\in\ind},\adm\rangle$, the \textit{foundation of} $\mathcal{F}$ is the tuple $\mathcal{F}_\found=\langle S,\sqsubseteq, \{R_i\}_{i\in\ind}\rangle$. We say that $\mathcal{F}$ is \textit{based on} $\mathcal{F}_\found$. 

Given a foundation $F=\langle S,\sqsubseteq, \{R_i\}_{i\in\ind}\rangle$, let $F^\full=\langle S,\sqsubseteq, \{R_i\}_{i\in\ind},\mathrm{RO}(S,\sqsubseteq)\rangle$. Given a possibility frame $\mathcal{F}=\langle S,\sqsubseteq,\{R_i\}_{i\in\ind},\adm\rangle$, if $\mathrm{RO}(S,\sqsubseteq)$ is closed under $\blacksquare_i$ for each $i\in\ind$, then we call $(\mathcal{F}_\found)^\full= {\langle S,\sqsubseteq,\{R_i\}_{i\in\ind},\mathrm{RO}(S,\sqsubseteq)\rangle}$ the \textit{associated full frame} of $\mathcal{F}$. \hfill $\triangleleft$
\end{definition}

Many possibility frames may be based on the same foundation. Yet at most one \textit{full} possibility frame may be based on a given foundation $F$. If $F=\langle S,\sqsubseteq, \{R_i\}_{i\in\ind}\rangle$ is such that $\mathrm{RO}(S,\sqsubseteq)$ is closed under $\blacksquare_i$, then $F^\full$ is the unique full possibility frame based on $F$. If $\mathrm{RO}(S,\sqsubseteq)$ is not closed under $\blacksquare_i$, then $F^\full$ is not even a partial-state frame as in Definition \ref{PosetMod}. In \S~\ref{FullFrames}, we will give first-order conditions on the interplay of $\sqsubseteq$ and $R_i$ that are equivalent to $\mathrm{RO}(S,\sqsubseteq)$ being closed under $\blacksquare_i$; and in \S~\ref{FullFrames} and \S~\ref{GFPF}, we will see that every possibility frame can be turned into a modally equivalent one satisfying even stronger interplay conditions.

The fact that closure of  $\mathrm{RO}(\mathcal{F})$ under $\cap$ and $\supset$ is guaranteed whereas closure under $\blacksquare_i$ is an extra requirement is reflected in a syntactic fact: to embed classical modal logic into intuitionistic modal logic by an extension of the G\"{o}del-Gentzen negative translation, the translation can simply ``pass through'' conjunctions and implications (and negations), but we must double-negate the box formulas, as follows.

\begin{definition}[Modal Negative Translation]\label{NegTrans} Let $G\colon \mathcal{L}(\sig, \ind)\to\mathcal{L}(\sig, \ind)$ be recursively defined as follows:
\begin{enumerate}
\item\label{NegTrans1} $p^{G}=\neg\neg p$;
\item $(\neg\varphi)^{G}=\neg\varphi^{G}$;
\item $(\varphi\wedge\psi)^{G} = \varphi^{G}\wedge\psi^{G}$;
\item $(\varphi\rightarrow\psi)^{G} = (\varphi^{G}\rightarrow\psi^{G})$;
\item\label{NegTrans6} $(\Box_i\varphi)^G=\neg\neg \Box_i \varphi^G$. \hfill $\triangleleft$
\end{enumerate}
\end{definition}
\noindent Since we defined $(\varphi\vee\psi):=\neg (\neg\varphi\wedge\neg\psi)$, we can think of $(\varphi\vee\psi)^G=\neg (\neg \varphi^G\wedge\neg\psi^G)$.

If we simply translated $\Box_i\varphi$ as $\Box_i \varphi^G$, we would not get the modal extension of the famous result that $\varphi$ is a theorem of classical logic iff $\varphi^G$ is a theorem of intuitionistic logic \citep{Godel1933,Gentzen1933,Gentzen1936,Gentzen1974}. For example, although $\vdash_\mathbf{K} \neg\neg\Box_i p\rightarrow\Box_i p$, one can check that $\neg\neg\Box_i\neg\neg p \rightarrow \Box_i\neg\neg p$ is not valid over intuitionistic frames and hence $\nvdash_\mathbf{HK} \neg\neg\Box_i\neg\neg p \rightarrow \Box_i\neg\neg p$ by Theorem \ref{Bozic}. But if we double-negate box formulas, we obtain the following  (cf. \citealt[pp. 231-2]{Bozic1984}).  

\begin{proposition}[Full and Faithful Translation]\label{IntTrans}
For all $\varphi\in\mathcal{L}(\sig, \ind)$: $\vdash_\mathbf{K} \varphi$ iff $\vdash_\mathbf{HK} \varphi^G$.
\end{proposition}

\begin{proof} From left to right, first, the famous result for propositional logic gives us that the translations of classical propositional tautologies are theorems of \textbf{HK}. Second, it is easy to check that the translation of the \textbf{K} axiom is a theorem of \textbf{HK} and that we can match in \textbf{HK} applications of Modus Ponens and Necessitation in \textbf{K}. Finally, that we can match in \textbf{HK} applications of Uniform Substitution in \textbf{K} follows from the fact that for any $\varphi\in\mathcal{L}(\sig,\ind)$ and substitution $\sigma\colon\sig\to \mathcal{L}(\sig,\ind)$, we have $\vdash_\mathbf{HK}(\widehat{\sigma}(\varphi))^G\leftrightarrow\widehat{\sigma^G}(\varphi^G)$, where $\widehat{\sigma}\colon \mathcal{L}(\sig,\ind)\to \mathcal{L}(\sig,\ind)$ is the usual extension of $\sigma$ to all formulas, and $\sigma^G$ is the substitution defined by $\sigma^G(p)=(\sigma(p))^G$. (Note: this part of the proof would fail if we had translated $\Box_i\varphi$ as $\Box_i \varphi^G$.) 

From right to left, if $\nvdash_\mathbf{K}\varphi$, then by the completeness of \textbf{K} with respect to the class of world models according to the standard satisfaction relation $\vDash$ of Kripke semantics (Definition \ref{pwmtruth1}), there is a world model $\mathfrak{M}$ that falsifies $\varphi$ according to $\vDash$. Thus, $\mathfrak{M}$ falsifies $\varphi^G$ according to $\vDash$, since $\varphi$ and $\varphi^G$ are equivalent according to $\vDash$. Viewing $\mathfrak{M}$ as a partial-state model as in Example \ref{KripkeExample}, it is easy to see that it counts as an intuitionistic modal model as in Example \ref{IntMod}, and over world models the forcing relation $\Vdash$ reduces to the relation $\vDash$. Thus, we have an intuitionistic modal model that falsifies $\varphi^G$ according to $\Vdash$, so $\nvdash_\mathbf{HK}\varphi^G$ by the soundness of \textbf{HK} with respect to intuitionistic modal models according to $\Vdash$ (Theorem \ref{Bozic}).
\end{proof}

Returning to our semantics, the reason that closure of $\mathrm{RO}(\mathcal{F})$ under $\blacksquare_i$ is nontrivial is because such closure depends on the interplay of the $R_i$ and $\sqsubseteq$ relations in $\mathcal{F}$, to which we turn in \S~\ref{FullFrames}.

Before proceeding further, let us consider a simple example of a full possibility frame.

\begin{figure}[h]
 \begin{center}
\begin{tikzpicture}[-,>=stealth',shorten >=1pt,shorten <=1pt, auto,node
distance=2cm,thick,every loop/.style={<-,shorten <=1pt}]
\tikzstyle{every state}=[fill=gray!20,draw=none,text=black]

\node (x_0) at (0,0) {{$s_0$}};
\node (w_0) at (-1,-1) {{$w_0$}};
\node (x_1) at (1,-1) {{$s_1$}};
\node (w_1) at (0,-2) {{$w_1$}};
\node (x_2) at (2,-2) {{$s_2$}};
\node (w_2) at (1,-3) {{$w_2$}};
\node (x_3) at (3,-3) {{}};

\node (dots) at (2.5,-2.5) {{\rotatebox[origin=c]{45}{$\boldsymbol{\vdots}$}}};

\node (dots) at (1.5,-3.5) {{\rotatebox[origin=c]{45}{$\boldsymbol{\vdots}$}}};

\path (x_0) edge[->] node {{}} (w_0);
\path (x_0) edge[->] node {{}} (x_1);
\path (x_1) edge[->] node {{}} (w_1);
\path (x_1) edge[->] node {{}} (x_2);
\path (x_2) edge[->] node {{}} (w_2);

\path (w_0) edge[->,dashed] node {{}} (x_1);
\path (w_1) edge[->,dashed] node {{}} (x_2);

\path (x_0) edge[->,dashed, bend left] node {{}} (x_1);
\path (x_1) edge[->,dashed, bend left] node {{}} (x_2);
\end{tikzpicture}

\end{center}
\caption{a temporal flow on the Beth comb.}\label{BethFlowFig}
\end{figure}
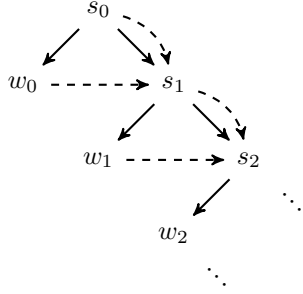

\begin{example}[A Temporal Flow on the Beth Comb]\label{BethFlowEx} Returning to the Beth comb of Example \ref{BethEx}, let us think of the $w_n$ as \textit{instants} of time, the $s_n$ as unending \textit{stretches} of time, and $s_{n+1}$ as the whole of \textit{the future} relative to $w_n$. Define an accessibility relation on the Beth comb by: $xR_< y$ iff $x$ is the $n$-th world or stretch of time and $y=s_{n+1}$. The result is shown in Figure \ref{BethFlowFig}, in which successive instants of time peel off of the future. First observe that if $X\in\mathrm{RO}(S,\sqsubseteq)$, then $\blacksquare_< X\in \mathrm{RO}(S,\sqsubseteq)$. For if $X$ is a regular open set of type (i) in Example \ref{BethEx}, then $\blacksquare_< X=\emptyset$, which is regular open; and if $X$ is a regular open set of type (ii) in Example \ref{BethEx}, then where $n=\min\{k\mid s_k\in X\}$,  $\blacksquare_< X=\mathord{\downarrow}s_{n-1}$, which is a regular open set of type (ii). Thus, the structure $\mathcal{F}=\langle S,\sqsubseteq, R_<,\mathrm{RO}(S,\sqsubseteq)\rangle$ is a full possibility frame. It is easy to see that this frame validates exactly the same formulas as the Kripke frame $\langle \mathbb{N},<\rangle$. Indeed, any possibility frame, like this one, in which every state is refined by an endpoint (world, instant) is semantically equivalent to a Kripke frame based on those endpoints (see \S~\ref{AtomicSection}). Yet this frame is a good example of how assumptions about relations and modal axioms from Kripke semantics must be reconsidered in possibility semantics. For example, our $R_<$ is \textit{functional}, and yet $\Diamond_< p\rightarrow \Box_< p$ is not valid. (Later we will see that for any full possibility frame, there is a semantically equivalent one in which the relations are partially functional.) Our $R_<$ is not transitive, and yet it validates $\Box_< p\rightarrow \Box_< \Box_< p$. This raises obvious questions about how modal axioms correspond to possibility frame properties, which we treat in \S~\ref{LemmScottCorr}. For now, note in connection with $\Box_< p\rightarrow \Box_< \Box_< p$ that the frame in Figure \ref{BethFlowFig} is such that $f_<(f_<(x))\sqsubseteq f_<(x)$, where $f_<(x)$ is the unique $y$ such that $xR_< y$.

Below we will see a possibility frame without endpoints, based on the infinite complete binary tree (Example \ref{BinaryEx}). We will also see a possibility frame for which there is no equivalent Kripke frame in \S~\ref{NoKripke}. \hfill $\triangleleft$
\end{example}

\subsection{The Interplay of Accessibility and Refinement}\label{FullFrames}

In this section, we will show that two first-order conditions on the interplay of $R_i$ and $\sqsubseteq$  are \textit{necessary} and \textit{sufficient} for $\mathrm{RO}(\mathcal{F})$ to be closed under $\blacksquare_i$ (Proposition \ref{ROtoRO}), so every full possibility frame satisfies these conditions. We will then show that every full possibility frame can be turned into a semantically equivalent one that satisfies a stronger and simpler condition on the interplay of $R_i$ and $\sqsubseteq$ (Proposition \ref{Representation}).

For $\mathrm{RO}(\mathcal{F})$ to be closed under $\blacksquare_i$, it must be that if $Y$ satisfies \textit{persistence} and \textit{refinability}, then so does $\blacksquare_iY$. Let us first consider \textit{persistence} for $\blacksquare_i Y$: if $x'\sqsubseteq x$ and $x\in\blacksquare_i Y$, then $x'\in\blacksquare_i Y$. The necessary and sufficient condition to ensure this, as shown by Proposition \ref{ROtoRO}, is that if a state $x$ has ``ruled out'' a state $z$, then $z$ remains ruled out by every refinement $x'$ of $x$ (recall Notation \ref{notation}):
\begin{itemize}
\item \Rrule{} -- if ${x}'\sqsubseteq{x}$, and for all $y\in R_i(x)$, $y\incomp z$, then for all $y'\in R_i(x')$, $y'\incomp z$.
\end{itemize}
Or equivalently: 
\begin{itemize}
\item \Rrule{} -- if ${x}'\sqsubseteq{x}$ and ${x}' R_i{y}'\comp{z}$, then $\exists {y}$: ${x} R_i{y}\comp{z}$ (see Figure \ref{RruleFig}).
\end{itemize}
Note that \Rrule{} is implied by the \Rcomm{} condition and hence the \upR{} condition from intuitionistic frames in Example \ref{IntMod}.\footnote{Also note that since $\sqsubseteq$ is transitive, \Rrule{} is equivalent to: if $x'\sqsubseteq x$, $x'R_iy'$, and $z'\sqsubseteq y'$, then $\exists y$: $xR_iy\comp z'$.} In increasing strength, the order is: \Rrule{}, \Rcomm{}, and \upR{} (see Figure \ref{InterplayTable}).

 \begin{figure}[h]
\begin{center}
\begin{tikzpicture}[->,>=stealth',shorten >=1pt,shorten <=1pt, auto,node
distance=2cm,thick,every loop/.style={<-,shorten <=1pt}]
\tikzstyle{every state}=[fill=gray!20,draw=none,text=black]
\node (x-up) at (0,0) {{$x$}};
\node (x) at (0,-1.5) {{$x'$}};
\node (y) at (2,-1.5) {{$y'$}};
\node (z) at (4,-1.5) {{$z$}};
\node (comp1a) at (3.15, -2.25) {{}};
\node (comp1b) at (2.85, -2.25) {{}};

\node at (5,-1.5) {{\textit{$\Rightarrow$}}};

\path (x) edge[dashed,->] node {{}} (y);
\path (x-up) edge[->] node {{}} (x);

\path (y) edge[->] node {{}} (comp1a);
\path (z) edge[->] node {{}} (comp1b);

\node at (7,.3) {{$\exists$}};
\node (x-up') at (6,0) {{$x$}};
\node (x') at (6,-1.5) {{$x'$}};
\node (y') at (8,-1.5) {{$y'$}};
\node (y-up') at (8,0) {{$y$}};

\node (z') at (10,0) {{$z$}};
\node (comp1a') at (9.15, -.75) {{}};
\node (comp1b') at (8.85, -.75) {{}};

\path (x') edge[dashed,->] node {{}} (y');
\path (x-up') edge[->] node {{}} (x');

\path (x-up') edge[dashed,->] node {{}} (y-up');

\path (y-up') edge[->] node {{}} (comp1a');
\path (z') edge[->] node {{}} (comp1b');

\end{tikzpicture}
\end{center}
\caption{the \Rrule{} condition.}\label{RruleFig}
\end{figure}

Next, consider \textit{refinability} for $\blacksquare_i Y$: if $x\not\in \blacksquare_i Y$, then $\exists x'\sqsubseteq x$ $\forall x''\sqsubseteq x'$, $x''\not\in \blacksquare_i Y$. The condition to ensure this can be understood by adopting the game perspective of Remark \ref{Game}.
  
\begin{remark}[Accessibility Game]\label{Game} Given a partial-state frame $\mathcal{F}=\langle S, \sqsubseteq , \{R_i\}_{i\in \ind},\adm\rangle$, $x,y\in S$, and $i\in\ind$,  the \textit{accessibility game} $\mathrm{G}(\mathcal{F},x,y,i)$ for players \textbf{A} and \textbf{E} has the following rounds, depicted in Figure \ref{AccGame}:
\begin{enumerate} 
\item \textbf{A} chooses a $y'\sqsubseteq y$; 
\item \textbf{E} chooses an $x'\sqsubseteq x$;
\item \textbf{A} chooses an $x''\sqsubseteq x'$;
\item[] if $R_i(x'')=\emptyset$, then \textbf{A} wins, otherwise play continues;
\item  \textbf{E} chooses a $y''\in R_i(x'')$;
\item[] if $y''\comp y'$, then \textbf{E} wins, otherwise \textbf{A} wins.
\end{enumerate}
One can think of \textbf{A} and \textbf{E} as arguing about whether $y$ is  \textit{accessible} to $x$: if it is, then for any way $y'$ of further specifying $y$, there should be some way $x'$ of further specifying $x$ that ``locks in'' access to states compatible with $y'$, i.e., such that all refinements $x''$ of $x'$ have access to some state $y''$ compatible with $y'$. If refinements of $x$ cannot keep up with refinements of $y$ in this way, then $y$ is not accessible to $x$. Thus, player \textbf{A} is trying to show that $y$ is not accessible to $x$, while player \textbf{E} is trying to block \textbf{A}'s argument.

Now consider the following condition on a partial-state frame $\mathcal{F}$:
\begin{itemize}
\item \Rwinweak{} -- if ${x}R_i{y}$, then $\forall {y'}\sqsubseteq {y}$ $\exists{x'}\sqsubseteq {x}$ $\forall {x''}\sqsubseteq {x'}$ $\exists {y''}\comp y'$: ${x''} R_i{y''}$.
\end{itemize}
This condition says that if $xR_iy$, then \textbf{E} has a \textit{winning strategy} in the accessibility game $\mathrm{G}(\mathcal{F},x,y,i)$, in line with our way of thinking about the accessibility game above. $\hfill$ $\triangleleft$
\end{remark}

\begin{figure}[h]
 \begin{center}
\begin{tikzpicture}[->,>=stealth',shorten >=1pt,shorten <=1pt, auto,node
distance=2cm,thick,every loop/.style={<-,shorten <=1pt}]
\tikzstyle{every state}=[fill=gray!20,draw=none,text=black]

\node (x) at (0,0) {{$x$}};
\node (y) at (3,0) {{$y$}};

\node (x') at (0,-1.5) {{$x'$}};
\node (x'') at (0,-3) {{$x''$}};
\node (y') at (3,-1.5) {{$y'$}};
\node (y'') at (3,-3) {{$y''$}};

\node (5) at (3.1,-2.3) {{\scalebox{1.4}{$\comp$} ?}};

\node (4) at (1.4,-3.4) {{4.\,\textbf{E} chooses}};

\path (x') edge[<-] node {{ 2.\,\textbf{E} chooses\;}} (x);

\path (y) edge[->] node {{ 1.\,\textbf{A} chooses}} (y');
\path (y'') edge[dashed,<-] node {{}} (x'');
\path (x'') edge[<-] node {{ 3.\,\textbf{A} chooses\;}} (x');

\end{tikzpicture}
\end{center}
\caption{the accessibility game $\mathrm{G}$ -- if $y''\comp y'$, \textbf{E} wins, otherwise \textbf{A} wins.}\label{AccGame}
\end{figure}

Now we will show that  \Rrule{} and \Rwinweak{} characterize closure of $\mathrm{RO}(S,\sqsubseteq)$ under $\blacksquare_i$.

\begin{proposition}[First-order Characterization of Closure of $\mathrm{RO}(S,\sqsubseteq)$ under $\blacksquare_i$]\label{ROtoRO} For any poset $\langle S,\sqsubseteq\rangle$ and binary relation $R_i$ on $S$, the following are equivalent: 
\begin{enumerate}
\item\label{ROtoRO1} $\mathrm{RO}(S,\sqsubseteq)$ is closed under $\blacksquare_i$;
\item\label{ROtoRO2} $R_i$ and $\sqsubseteq$ satisfy \Rrule{} and \Rwinweak{}.
\end{enumerate}
\end{proposition}

\begin{proof} We begin with the implication from \ref{ROtoRO2} to \ref{ROtoRO1}. Assume that $X\subseteq S$ satisfies \textit{persistence} and \textit{refinability}.

To show that $\blacksquare_iX$ satisfies \textit{persistence}, suppose $x'\sqsubseteq x$ and $x'\not\in\blacksquare_i X$, so there is a $y'$ such that $x' R_iy'$ and $y'\not\in X$. Then since $X$ satisfies \textit{refinability}, there is a $z\sqsubseteq y'$ such that (i) for all $z'\sqsubseteq z$, $z'\not\in X$. Since $z\sqsubseteq y'$ gives us $y'\comp z$, from $x' R_iy'$ and \Rrule{} we have a $ y$ such that $ x R_i y\comp z$. Then from $ y\comp z$, (i), and  \textit{persistence} for $X$, it follows that $y\not\in X$, which with $ x R_i y$ implies $x\not\in\blacksquare_i X$. Thus, $\blacksquare_iX$ satisfies \textit{persistence}.

For \textit{refinability}, suppose that $x\not\in\blacksquare_i X$, so there is a $y$ such that $xR_i y$ and $y\not\in X$. Then since $X$ satisfies \textit{refinability}, there is a $y'\sqsubseteq y$ such that (ii) for all $ z\sqsubseteq y'$, $z\not\in X$. Assuming \Rwinweak{}, from $xR_iy$ and $y'\sqsubseteq y$ we have that $\exists x'\sqsubseteq x$ $\forall x''\sqsubseteq x'$ $\exists y''$: $y''\comp y'$ and $x'' R_i y''$. From $y''\comp y'$, (ii), and \textit{persistence} for $X$, it follows that $y''\not\in X$, which with $x'' R_i y''$ implies $x''\not\in\blacksquare_i X$. So we have shown that if $x\not\in\blacksquare_i X$, then $\exists x'\sqsubseteq x$ $\forall x''\sqsubseteq x'$,  $x''\not\in \blacksquare_i X$. Thus, $\blacksquare_i X$ satisfies \textit{refinability}.

Now let us prove the implication from \ref{ROtoRO1} to \ref{ROtoRO2}. 

First, suppose that \Rrule{} does not hold, so we have ${x}'\sqsubseteq{x} $ and ${x}' R_i{y}'\comp{z}$, but for all ${y} $ with ${x}  R_i{y} $, we have ${y} \incomp{z}$. Let $V=\{v\in S\mid v\incomp z\}$, so ${x} \in \blacksquare_i V$. One can check that $V$ satisfies \textit{persistence} and \textit{refinability} (see Fact \ref{ConVal}), so $V\in\mathrm{RO}(S,\sqsubseteq)$. Since ${y}' \comp {z}$, ${y}'\not\in V$, which with ${x}' R_i{y}'$ implies  ${x}'\not\in\blacksquare_i V$. But then since ${x}'\sqsubseteq{x} $, ${x} \in \blacksquare_i V$, and ${x}'\not\in\blacksquare_i V$, $\blacksquare_iV$ does not satisfy \textit{persistence}, so $\blacksquare_iV\not\in \mathrm{RO}(S,\sqsubseteq)$. But then since $V\in\mathrm{RO}(S,\sqsubseteq)$, $\mathrm{RO}(S,\sqsubseteq)$ is not closed under $\blacksquare_i$.

Second, suppose that \Rwinweak{} does not hold, so we have ${x}R_i{y}$ and $\exists {y'}\sqsubseteq {y}$ $\forall{x'}\sqsubseteq {x}$ $\exists {x''}\sqsubseteq {x'}$ $\forall {y''}$: ${x''} R_i{y''}$ implies ${y''}\incomp {y'}$. Let $V=\{v\in S\mid  v\incomp y'\}$, so $V\in\mathrm{RO}(S,\sqsubseteq)$ by the same reasoning as above. Since $y'\sqsubseteq y$, $y\not\in V$, which with $xR_iy$ implies $x\not\in\blacksquare_i V$. But it also follows from our supposition that $\forall{x'}\sqsubseteq {x}$ $\exists {x''}\sqsubseteq {x'}$: $x''\in\blacksquare_i V$. Thus, $\blacksquare_i V$ does not satisfy \textit{refinability}, so $\blacksquare_i V\not\in \mathrm{RO}(S,\sqsubseteq)$. But then since $V\in\mathrm{RO}(S,\sqsubseteq)$, $\mathrm{RO}(S,\sqsubseteq)$ is not closed under $\blacksquare_i$.\end{proof}

As an immediate corollary of Proposition \ref{ROtoRO}, we have the following (recall Definitions \ref{PosFrames} and \ref{Foundation}).
 
\begin{corollary}[Interplay of $R_i$ and $\sqsubseteq$ for Full Possibility Frames]\label{Fullinterplay} $\,$
\begin{enumerate}
\item\label{Fullinterplay1} If $\mathcal{F}$ is a full possibility frame, then $\mathcal{F}$ satisfies \Rrule{} and \Rwinweak{};
\item\label{Fullinterplay2} For any foundation $F=\langle S,\sqsubseteq, \{R_i\}_{i\in \ind}\rangle$, $F^\full = \langle S,\sqsubseteq, \{R_i\}_{i\in \ind},\mathrm{RO}(S,\sqsubseteq)\rangle $  is a full possibility frame iff  $F$ satisfies \Rrule{} and \Rwinweak{}.
\end{enumerate}
\end{corollary}

In addition to identifying conditions so that $\mathrm{RO}(\mathcal{F})$ is closed under $\blacksquare_i$, let us identify conditions so that for every state $x$, its set $R_i(x)$ of accessible states is in $\mathrm{RO}(\mathcal{F})$. 

\begin{fact}[$R_i(x)$ and $\mathrm{RO}(\mathcal{F})$]\label{R(x)RO} For any partial-state frame $\mathcal{F}=\langle S, \sqsubseteq , \{R_i\}_{i\in \ind},\adm\rangle$, the following are equivalent:
\begin{enumerate}
\item for all $x\in S$, $R_i(x)\in \mathrm{RO}(\mathcal{F})$;
\item $\mathcal{F}$ satisfies both
\begin{enumerate}
\item \Rdown{} -- if $xR_iy$ and $y'\sqsubseteq y$, then $xR_iy'$ (recall Figure \ref{RdownFig});
\item \Rdense{} -- $xR_iy$ if $\forall y'\sqsubseteq y$ $\exists y''\sqsubseteq y'$: $xR_iy''$ (see Figure \ref{RdenseFig}).
\end{enumerate}
\end{enumerate}
\end{fact}
\begin{proof} \Rdown{} corresponds to \textit{persistence} of $R_i(x)$ and \Rdense{} corresponds to \textit{refinability} of $R_i(x)$.
\end{proof}

 \begin{figure}[h]
\begin{center}
\begin{tikzpicture}[->,>=stealth',shorten >=1pt,shorten <=1pt, auto,node
distance=2cm,thick,every loop/.style={<-,shorten <=1pt}]
\tikzstyle{every state}=[fill=gray!20,draw=none,text=black]

\node (x) at (0,0) {{$x$}};
\node (y) at (2,0) {{$y$}};
\node (y') at (2,-1.5) {{$y'$}};
\node (y'') at (2,-3) {{$y''$}};

\node at (3,0) {{\textit{$\Rightarrow$}}};

\path (x) edge[dashed,->] node {{}} (y'');
\path (y) edge[->] node {{$\forall$}} (y');
\path (y') edge[->] node {{$\exists$}} (y'');

\node (xR) at (4,0) {{$x$}};
\node (yR) at (6,0) {{$y$}};

\path (xR) edge[dashed,->] node {{}} (yR);

\end{tikzpicture}
\end{center}
\caption{\Rdense{}}\label{RdenseFig}
\end{figure}

Note how the conditions above relate to the \upR{} condition from intuitionistic frames (Example \ref{IntMod}). 

\begin{fact}[Deriving \upR{}]\label{DerivingUp} For any partial-state frame $\mathcal{F}$, if $\mathcal{F}$ satisfies \Rrule{}, \Rdown{}, and \Rdense{}, then $\mathcal{F}$ satisfies \upR{}.
\end{fact}

\begin{proof} To show that $x'\sqsubseteq x$ and $x'R_iy$ together imply $xR_iy$, suppose $x'\sqsubseteq x$ but $y\not\in R_i(x)$. Then by \Rdense{}, there is a $y'\sqsubseteq y$ such that (i) for all $y''\sqsubseteq y'$, $y''\not\in R_i(x)$. We claim that for all $z\in R_i(x)$, $z\incomp y'$. For if $z\comp y'$, then there is a $y''\sqsubseteq z$ with $y''\sqsubseteq y'$; and then by \Rdown{}, $z\in R_i(x)$ and $y''\sqsubseteq z$ together imply $y''\in R_i(x)$; but by (i), $y''\sqsubseteq y'$ implies $y''\not\in R_i(x)$. Let $V=\{v\in S\mid v\incomp y'\}$, so $x\in\blacksquare_i V$. One can check that $V\in\mathrm{RO}(\mathcal{F})$ (see Fact \ref{ConVal}), so by the proof of Proposition \ref{ROtoRO}, \Rrule{} implies that $\blacksquare_iV$ satisfies \textit{persistence}. Thus, from $x\in\blacksquare_i V$ and $x'\sqsubseteq x$ we have $x'\in\blacksquare_i V$, which implies $y'\not\in R_i(x')$ by the definition of $V$, which with $y'\sqsubseteq y$ implies $y\not\in R_i(x')$ by \Rdown{}. This shows that $\mathcal{F}$ satisfies \upR{}.\end{proof}

Next, note how the \Rdown{} condition from intuitionistic frames simplifies other conditions.

\begin{fact}[Simplifying with \Rdown{}]\label{Simp} For any partial-state frame satisfying \Rdown{}, the following are equivalent:
\begin{enumerate}
\item \Rwinweak{} -- if ${x}R_i{y}$, then $\forall {y'}\sqsubseteq {y}$ $\exists{x'}\sqsubseteq {x}$ $\forall {x''}\sqsubseteq {x'}$ $\exists {y''}\comp y'$: ${x''} R_i{y''}$;
\item \Rwin{} -- if $xR_iy$, then $\forall y'\sqsubseteq y$ $\exists x'\sqsubseteq x$ $\forall x''\sqsubseteq x'$ $\exists y''\sqsubseteq y'$: $x'' R_i y''$;
\item \Rref{} -- if $xR_iy$, then $\exists x'\sqsubseteq x$ $\forall x''\sqsubseteq x'$ $\exists y'\sqsubseteq y$: $x'' R_i y'$.
\end{enumerate}
\end{fact}

\begin{figure}[h]
 \begin{center}
\begin{tikzpicture}[->,>=stealth',shorten >=1pt,shorten <=1pt, auto,node
distance=2cm,thick,every loop/.style={<-,shorten <=1pt}]
\tikzstyle{every state}=[fill=gray!20,draw=none,text=black]

\node (x) at (0,0) {{$x$}};
\node (y) at (3,0) {{$y$}};

\node (x') at (0,-1.5) {{$x'$}};
\node (x'') at (0,-3) {{$x''$}};
\node (y') at (3,-1.5) {{$y'$}};
\node (y'') at (3,-3) {{$y''$}};

\node (4) at (1.4,-3.4) {{4.\,\textbf{E} chooses}};

\path (x') edge[<-] node {{ 2.\,\textbf{E} chooses\;}} (x);

\path (y) edge[->] node {{ 1.\,\textbf{A} chooses}} (y');
\path (y') edge[->] node {{?}} (y'');
\path (y'') edge[dashed,<-] node {{}} (x'');
\path (x'') edge[<-] node {{ 3.\,\textbf{A} chooses\;}} (x');

\end{tikzpicture}
\end{center}
\caption{the accessibility game $\underline{\mathrm{G}}$ -- if $y''\sqsubseteq y'$, \textbf{E} wins, otherwise \textbf{A} wins.}\label{AccGame2}
\end{figure}
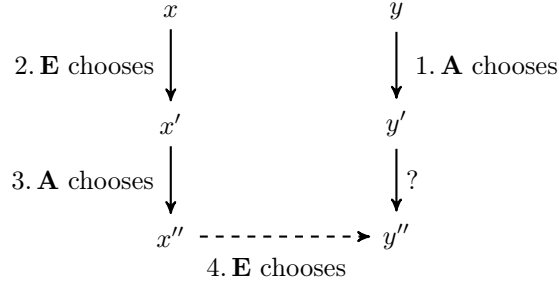

Corresponding to \Rwin{} is a modified accessibility game $\underline{\mathrm{G}}(\mathcal{F},x,y,i)$ that differs from the accessibility game $\mathrm{G}(\mathcal{F},x,y,i)$ of Remark \ref{Game} by changing the winning condition (see Figure \ref{AccGame2}):
\begin{itemize}
\item if $y''\sqsubseteq y'$, then \textbf{E} wins, otherwise \textbf{A} wins.
\end{itemize}
As shown by Fact \ref{Simp}, in frames satisfying \Rdown{} the accessibility games $\underline{\mathrm{G}}$ and $\mathrm{G}$ are equivalent, i.e., \textbf{E} has a winning strategy in the one iff \textbf{E} has a winning strategy in the other.

Having considered the condition that $xR_iy$ \textit{implies} that \textbf{E} has a winning strategy in the accessibility game $\underline{\mathrm{G}}(\mathcal{F},x,y,i)$, it is natural to consider the stronger condition that $xR_iy$ is \textit{equivalent} to \textbf{E} having a winning strategy in the accessibility game $\underline{\mathrm{G}}(\mathcal{F},x,y,i)$:
\begin{itemize}
\item \RWin{} -- $xR_iy$ iff $\forall y'\sqsubseteq y$ $\exists x'\sqsubseteq x$ $\forall x''\sqsubseteq x'$ $\exists y''\sqsubseteq y'$: $x'' R_i y''$ (see Figure \ref{Rwinfig}).
\end{itemize}
Remarkably, not only does this one condition imply \textit{all} of the others, and not only is it equivalent to the conjunction of the conditions from Proposition \ref{ROtoRO} and Fact \ref{R(x)RO}, but also, we may assume \RWin{}  without loss of generality when working with full possibility frames. We prove these claims in turn.

\begin{proposition}[Ultimate Condition]\label{MasterCon} The condition \RWin{} implies  all of the conditions above and is implied by the conjunction of the conditions from Proposition \ref{ROtoRO} and Fact \ref{R(x)RO}: \Rrule{}, \Rwinweak{}, \Rdown{}, and \Rdense{}. Thus, by Proposition \ref{ROtoRO} and Fact \ref{R(x)RO}, for any partial-state frame $\mathcal{F}=\langle S,\sqsubseteq,\{R_i\}_{i\in\ind},\adm\rangle$, the following are equivalent:
\begin{enumerate}
\item $\mathcal{F}$ satisfies \RWin{};
\item for all $X\in \mathrm{RO}(\mathcal{F})$ and $x\in S$, $\blacksquare_iX\in\mathrm{RO}(\mathcal{F})$ and $R_i(x)\in\mathrm{RO}(\mathcal{F})$.
\end{enumerate}\end{proposition}

\begin{proof} First, we show that \RWin{} implies all the other conditions. It is easy to see that \RWin{} implies \upR{} and therefore \Rrule{} and that \RWin{} implies \Rdown{} and therefore the conditions in Fact \ref{Simp}. To see that \RWin{} implies \Rdense{}, suppose it is not the case that $xR_iy$. Then by the right-to-left direction of \RWin{}, there is a $y'\sqsubseteq y$ such that (i) $\forall x'\sqsubseteq x$ $\exists x''\sqsubseteq x'$ $\forall y''\sqsubseteq y'$, \textit{not} $x'' R_i y''$. Now we claim that  for all $y''\sqsubseteq y'$, it is not the case that $xR_iy''$. For if $xR_iy''$, then by the left-to-right direction of \RWin{}, $\exists x'\sqsubseteq x$ $\forall x''\sqsubseteq x'$ $\exists y'''\sqsubseteq y''$: $x''R_i y'''$, which with $y'''\sqsubseteq y'$ contradicts (i). Thus, we have shown that if \textit{not} $xR_iy$, then $\exists y'\sqsubseteq y$ $\forall y''\sqsubseteq y'$, \textit{not} $xR_iy''$, which is the contrapositive of \Rdense.

Next, we show that the stated conditions imply \RWin{}. The left-to-right direction of \RWin{}, \Rwin{}, follows from \Rwinweak{} and \Rdown{} by Fact \ref{Simp}. For the right-to-left direction, suppose $\forall y'\sqsubseteq y$ $\exists x'\sqsubseteq x$ $\forall x''\sqsubseteq x'$ $\exists y''\sqsubseteq y'$: $x'' R_i y''$, which implies $xR_iy''$ by \upR{}, which follows from the other stated conditions by Fact \ref{DerivingUp}. Thus, $\forall y'\sqsubseteq y$ $\exists y''\sqsubseteq y'$: $xR_iy''$, which implies $xR_iy$ by \Rdense.
\end{proof}

Given the strength and importance of the \RWin{} condition, we introduce the following terminology.

\begin{definition}[Strong Possibility Frames]\label{StrongPoss} A \textit{strong} possibility frame is a possibility frame as in Definition \ref{PosFrames} that satisfies \RWin{}. \hfill $\triangleleft$
\end{definition}

\begin{figure}[h]
 \begin{center}
\begin{tikzpicture}[->,>=stealth',shorten >=1pt,shorten <=1pt, auto,node
distance=2cm,thick,every loop/.style={<-,shorten <=1pt}]
\tikzstyle{every state}=[fill=gray!20,draw=none,text=black]

\node (x0) at (-5,0) {{$x$}};
\node (y0) at (-2,0) {{$y$}};

\path (x0) edge[dashed,->] node {{}} (y0);

\node at (-1,0) {{\textit{$\Leftrightarrow$}}};

\node (x) at (0,0) {{$x$}};
\node (y) at (3,0) {{$y$}};

\node (x') at (0,-1.5) {{$x'$}};
\node (x'') at (0,-3) {{$x''$}};
\node (y') at (3,-1.5) {{$y'$}};
\node (y'') at (3,-3) {{$y''$}};

\path (x') edge[<-] node {{ b.$\,\exists\;$}} (x);
\path (x) edge[dashed,->] node {{}} (y);
\path (y) edge[->] node {{ a.$\,\forall$}} (y');
\path (y') edge[->] node {{ d.$\,\exists$}} (y'');
\path (y'') edge[dashed,<-] node {{}} (x'');
\path (x'') edge[<-] node {{ c.$\,\forall\;$ }} (x');

\end{tikzpicture}
\end{center}
\caption{\RWin{}}\label{Rwinfig}
\end{figure}

We will now prove that any full possibility frame can be turned into a semantically equivalent \textit{strong} and full possibility frame, simply by modifying the accessibility relations in the frame. 

\begin{proposition}[From Full Possibility Frames to Strong Possibility Frames]\label{Representation} For any partial-state frame $\mathcal{F}=\langle S,\sqsubseteq, \{R_i\}_{i\in\ind}, \adm\rangle$,  define $\mathcal{F}^\tight=\langle S,\sqsubseteq, \{R_i^\tight\}_{i\in\ind}, \adm\rangle$ by $xR_i^\tight y$ iff for all $Z\in\adm$, $x\in\blacksquare_i^\mathcal{F} Z$ implies $y\in Z$, where $\blacksquare_i^\mathcal{F}$ is the $\blacksquare_i$ operator for $\mathcal{F}$. Then:
\begin{enumerate}
\item\label{Representation1} if every $X\in\adm$ satisfies \textit{persistence}, then $\mathcal{F}^\tight$ satisfies \upR{} and \Rdown{};
\item\label{Representation2} if $\mathcal{F}$ is a \textit{possibility frame}, then $\mathcal{F}^\tight$ is a possibility frame satisfying \Rdense{};
\item\label{Representation3} if $\mathcal{F}$ is a \textit{full} possibility frame, then $\mathcal{F}^\tight$ is a \textit{strong} and full possibility frame; 
\item\label{Representation4} for all $\pi\colon \sig\to\adm$, $x\in S$, and $\varphi\in\mathcal{L}(\sig,\ind)$: $\langle \mathcal{F},\pi\rangle, x\Vdash \varphi$ iff $\langle \mathcal{F}^\tight,\pi\rangle, x\Vdash \varphi$.
\end{enumerate}
\end{proposition}

\begin{proof} For \ref{Representation1}, we prove \upR{} and \Rdown{} at the same time. Suppose $x'\sqsubseteq x$, $y'\sqsubseteq y$ and $x'R_i^\tight y$. To show that $xR_i^\tight y'$, consider any $Z\in\adm$ with $x\in\blacksquare_i^\mathcal{F}Z$. Then since $\mathcal{F}$ is a partial-state frame, $\blacksquare_i^\mathcal{F} Z\in \adm$, so $\blacksquare_i^\mathcal{F} Z$ satisfies \textit{persistence} by the assumption of part \ref{Representation1}, so from $x'\sqsubseteq x$  and $x\in\blacksquare_i^\mathcal{F}Z$ we have $x'\in\blacksquare_i^\mathcal{F}Z$. Then since $x'R_i^\tight y$, we have $y\in Z$, which with $y'\sqsubseteq y$ and \textit{persistence} for $Z$ implies $y'\in  Z$. Thus, $xR_i^\tight y'$. 

For part \ref{Representation2}, we first show that $\mathcal{F}^\tight$ is a possibility frame, assuming that $\mathcal{F}$ is. Since $\mathrm{RO}(\mathcal{F})=\mathrm{RO}(\mathcal{F}^\tight)$, and the set $\adm$ of admissible propositions is the same in $\mathcal{F}$ and $\mathcal{F}^\tight$, we need only show that $\adm$ is closed under $\blacksquare_i^{\mathcal{F}^\tight}$, the $\blacksquare_i$ operator for $\mathcal{F}^\tight$. Since $\adm$ is closed under $\blacksquare_i^\mathcal{F}$ by assumption, it suffices to show that for all $Z\in \adm$, $\blacksquare_i ^\mathcal{F}Z = \blacksquare_i^{\mathcal{F}^\tight}Z$. To see that $\blacksquare_i ^\mathcal{F}Z \supseteq \blacksquare_i^{\mathcal{F}^\tight}Z$, suppose $x\not\in \blacksquare_i ^\mathcal{F}Z$, so there is a $y$ such that $xR_iy$ but $y\not\in Z$.  Then since $xR_iy$ clearly implies $xR_i^\tight y$, we have $x\not\in \blacksquare_i^{\mathcal{F}^\tight} Z$. To see that $\blacksquare_i ^\mathcal{F}Z \subseteq \blacksquare_i^{\mathcal{F}^\tight}Z$, suppose $x\not\in  \blacksquare_i^{\mathcal{F}^\tight}Z$, so there is a $y$ such that $xR_i^\tight y$ but $y\not\in Z$. Then by the definition of $R_i^\tight$, $x\not\in\blacksquare_i^\mathcal{F} Z$. 

For \Rdense{} in part \ref{Representation2}, assume that $\forall y'\sqsubseteq y$ $\exists y''\sqsubseteq y'$: $x R_i^\tight y''$. For reductio, suppose \textit{not} $xR_i^\tight y$, so there is some $Z\in\adm$ such that $x\in\blacksquare_i^\mathcal{F} Z$ but $y\not\in Z$. Then since $Z$ satisfies \textit{refinability}, $\exists y'\sqsubseteq y$ $\forall y''\sqsubseteq y'$, $y''\not\in Z$, which with $x\in\blacksquare_i^\mathcal{F} Z$ implies that \textit{not} $xR_i^\tight y''$.  But this contradicts our initial assumption. 

For part \ref{Representation3}, we have already shown that $\mathcal{F}^\tight$ is a possibility frame if $\mathcal{F}$ is. Then since $\mathrm{RO}(\mathcal{F})=\mathrm{RO}(\mathcal{F}^\tight)$, and the set $\adm$ of admissible propositions is the same in $\mathcal{F}$ and $\mathcal{F}^\tight$, $\mathcal{F}^\tight$ is a full possibility frame if $\mathcal{F}$ is.

To show that $\mathcal{F}^\tight$ satisfies \RWin{}, by Proposition \ref{MasterCon} it suffices to show \upR{}, \Rdown{}, \Rdense{}, and \Rref{}, the first three of which we have already shown. For \Rref{}, suppose for reductio that $xR_i^\tight y$ but $\forall x'\sqsubseteq x$ $\exists x''\sqsubseteq x'$ such that (i) $\forall y'$, if $x''R_i^\tight y'$ then $y'\not\sqsubseteq y$. It follows that (ii) $\forall y'$, if $x''R^\tight_iy'$, then $y'\incomp y$. For otherwise there is a $y''\sqsubseteq y'$ with $y''\sqsubseteq y$, and by \Rdown{},  $x''R^\tight_iy'$ and $y''\sqsubseteq y'$ together imply $x''R^\tight_i y''$, which with $y''\sqsubseteq y$ contradicts (i). Then since $x''R_i y'$ implies $x''R_i^\tight y'$, (ii) implies (iii) $\forall y'$, if $x''R_iy'$, then $y'\incomp y$. Now define $V=\{v\in S\mid v\incomp y\}$. One can check that $V$ satisfies \textit{persistence} and \textit{refinability} (see Fact \ref{ConVal}), so $V\in\mathrm{RO}(\mathcal{F})$. Then because $\mathcal{F}$ is a \textit{full} possibility frame, so $\adm=\mathrm{RO}(\mathcal{F})$, we have $V\in\adm$.  Since $xR_i^\tight y$ and $y\not\in V$, $x\not\in\blacksquare_i^\mathcal{F} V$,  but $\forall x'\sqsubseteq x$ $\exists x''\sqsubseteq x'$: $x''\in \blacksquare_i^\mathcal{F} V$ by (iii). Thus, $\blacksquare_i^\mathcal{F} V$ does not satisfy \textit{refinability}, so $\blacksquare_i^\mathcal{F} V\not\in\mathrm{RO}(\mathcal{F})$ and hence $\blacksquare_i^\mathcal{F} V\not\in\adm$. But if $V\in\adm$ and $\blacksquare_i^\mathcal{F} V\not\in\adm$, then by Definition \ref{PosetMod}, $\mathcal{F}$ is not a partial-state frame, which contradicts the assumption that $\mathcal{F}$ is a possibility frame.

The proof of part \ref{Representation4} is by induction on $\varphi$. Since $\mathcal{F}^\tight$ differs from $\mathcal{F}$ only with respect to the accessibility relations, the only inductive case to check is the $\Box_i$ case, which follows from the fact established above that for all $Z\in \adm$, $\blacksquare_i ^\mathcal{F}Z = \blacksquare_i^{\mathcal{F}^\tight}Z$, together with Fact \ref{TruthSub}. \end{proof} 

Not only can any \textit{full} possibility frame be transformed into a semantically equivalent \textit{strong} possibility frame, but also in \S~\ref{GFPF} we will show that \textit{any} possibility frame can be transformed into a semantically equivalent strong possibility frame. Thus, one may assume without loss of generality that \RWin{} gets at the essence of the interplay between accessibility and refinement in possibility frames.

We have already seen two examples of strong possibility frames. Every world frame, viewed as a possibility frame as in Example \ref{KripkeAgain}, is trivially a strong possibility frame. Less trivially, the \textit{powerset possibilization} of any world frame as in Example \ref{PowerPoss} is also a strong possibility frame.

\begin{fact}[Powerset Possibilization Cont.]\label{transform2} For any world frame $\mathfrak{F}=\langle \wo{W},\{\wo{R}_i\}_{i\in\ind},\wo{A}\rangle$, its powerset possibilization $\mathfrak{F}^\pow$ is a strong possibility frame.
\end{fact}

\begin{proof} By Proposition \ref{MasterCon}, it suffices to show that $\mathfrak{F}^\pow$ satisfies \upR{} and \Rdown{}, which is easy to check, and: 
\begin{itemize}
\item \Rdense{} -- $XR_i^\pow Y$ if $\forall Y'\sqsubseteq Y$ $\exists Y''\sqsubseteq Y'$: $XR_i^\pow Y''$;
\item \Rref{} --  if $X R_i^\pow Y$, then $\exists X'\sqsubseteq X$ $\forall X''\sqsubseteq X'$ $\exists Y'\sqsubseteq Y$: $X''R_i^\pow Y'$.
\end{itemize} 
For \Rdense{}, given nonempty $X,Y\subseteq\wo{W}$, assume that for all nonempty $Y'\subseteq Y$ there is a nonempty $Y''\subseteq Y'$ with $XR_i^\pow Y''$, i.e., $Y''\subseteq \wo{R}_i[X]$. Then for any $y\in Y$, taking $Y'=\{y\}$ implies $Y''=\{y\}\subseteq\wo{R}_i[X]$. Since this holds for every $y\in Y$, we have $Y\subseteq \wo{R}_i[X]$, i.e., $XR_i^\pow Y$. 

For \Rref{}, we must check that if $\emptyset\not=Y\subseteq \mathrm{R}_i[X]$, then there is a nonempty $X'\subseteq X$ such that for all nonempty $X''\subseteq X'$, there is a nonempty $Y_{X''}\subseteq Y$ with $Y_{X''}\subseteq \mathrm{R}_i[X'']$. Since $\emptyset\not=Y\subseteq \mathrm{R}_i[X]$, $X\cap\mathrm{R}_i^{-1}[Y]\not = \emptyset$, so pick $x\in X\cap\mathrm{R}_i^{-1}[Y]$ and $y\in \mathrm{R}_i(x)\cap Y$. Setting $X'=\{x\}$, we have $X'\subseteq X$, and there is only one nonempty $X''\subseteq X'$, namely $X'$ itself. Then setting $Y_{X'}=\{y\}$, we have $Y_{X'}\subseteq Y$ and $Y_{X'}\subseteq \mathrm{R}_i[X']$, so we are done.
\end{proof}

This proof actually shows that $\mathfrak{F}^\pow$ satisfies the following stronger condition (see Figure \ref{InterplayTable}):
\begin{itemize}
\item \RrefPlus{} -- if $xR_iy$, then $\underline{\exists y'\sqsubseteq y}$ $\exists x'\sqsubseteq x$ $\forall x''\sqsubseteq x'$: $x''R_iy'$.
\end{itemize}
We will see this condition satisfied in another frame in \S~\ref{NoKripke}. (Whether every frame is modally equivalent to one satisfying \RrefPlus{} depends on the assumption of the ultrafilter axiom. See Appendix \S~\ref{Strengths}.)

We are now in a good position to consider Humberstone's \citeyearpar{Humberstone1981} original frames for possibility semantics.  

\begin{remark}[Humberstone Frames]\label{HumbFrame} While we defined full possibility frames as partial-state frames $\mathcal{F}$ in which $\adm=\mathrm{RO}(\mathcal{F})$, Humberstone \citeyearpar{Humberstone1981} built strong conditions on the interplay of $R_i$ and $\sqsubseteq$ into his definition of frames. A \textit{Humberstone frame} is a tuple $\mathcal{F}=\langle S,\sqsubseteq,\{R_i\}_{i\in\ind},\adm\rangle$, where $\langle S,\sqsubseteq\rangle$ is a nonempty poset, $R_i$ is a binary relation on $S$, and $\adm=\mathrm{RO}(S,\sqsubseteq)$, such that $\mathcal{F}$ satisfies \upR{}, \Rdown{}, and:
\begin{itemize}
\item \RrefPlusPlus{} -- if $xR_iy$, then $\exists x'\sqsubseteq x$ $\forall x''\sqsubseteq x'$, $x''R_iy$. 
\end{itemize}
A \textit{Humberstone model} is a tuple $\mathcal{M}=\langle \mathcal{F},\pi\rangle$ where $\mathcal{F}$ is a Humberstone frame and $\pi\colon \sig\to \adm$. (In fact, Humberstone \citeyearpar{Humberstone1981} took $\pi$ to be a \textit{partial} function, but that approach is equivalent to the approach of this paper, as explained in Remark \ref{ThreeVals}. Later Humberstone \citeyearpar[p.~900]{Humberstone2011} took the total function approach.)

Since \upR{} and \RrefPlusPlus{} together imply the conditions of Proposition \ref{ROtoRO}, the set $\adm$ of admissible propositions in a Humberstone frame is closed under $\blacksquare_i$, and as previously noted, $\mathrm{RO}(S,\sqsubseteq)$ is automatically closed under $\cap$ and $\supset$. Thus, Humberstone frames are partial-state frames as in Definition \ref{PosetMod}, and since they satisfy $\adm=\mathrm{RO}(S,\sqsubseteq)$, they are full possibility frames as in Definition \ref{PosFrames}.

Humberstone's \RrefPlusPlus{} condition is not implied by any combination of the other conditions discussed in this section. In contrast to Fact \ref{transform2}, the powerset possibilization of a Kripke frame is not necessarily a Humberstone frame (see Appendix \S~\ref{Strengths}). It is an open question whether Humberstone frames are as general as the full possibility frames of Definition \ref{PosFrames}, or even whether they are more general than Kripke frames, for the purposes of characterizing normal modal logics (see Problem \ref{HumProb} in \S~\ref{OpenProb}). 

In \S~\ref{AtomicSection}, we will see that all \textit{finite} full possibility frames, and hence all finite Humberstone frames, can be turned into modally equivalent Kripke frames. By contrast, in \S~\ref{NoKripke} we will construct \textit{infinite} full possibility frames for which there are no modally equivalent Kripke frames. The question of whether every infinite Humberstone frame is modally equivalent to a Kripke frame is an open question.\hfill $\triangleleft$\end{remark}

Let us finally consider a concrete example of the interplay of accessibility and refinement.

\begin{figure}[h]
\begin{center}
\begin{tikzpicture}[->,>=stealth',shorten >=1pt,shorten <=1pt, auto,node
distance=2cm,thick,every loop/.style={<-,shorten <=1pt}]
\tikzstyle{every state}=[fill=gray!20,draw=none,text=black]

\node (empty) at (0,0) {{$\emptyset$}};
\node (0) at (-4,-1) {{$0$}};
\node (1) at (4,-1) {{$1$}};

\node (00) at (-6,-2) {{$00$}};
\node (01) at (-2,-2) {{$01$}};
\node (10) at (2,-2) {{$10$}};
\node (11) at (6,-2) {{$11$}};

\node (000) at (-7,-3) {{$000$}};
\node (001) at (-5,-3) {{$001$}};

\node (010) at (-3,-3) {{$010$}};
\node (011) at (-1,-3) {{$011$}};

\node (100) at (1,-3) {{$100$}};
\node (101) at (3,-3) {{$101$}};

\node (110) at (5,-3) {{$110$}};
\node (111) at (7,-3) {{$111$}};

\path (empty) edge[->] node {{}} (0);
\path (empty) edge[->] node {{}} (1);

\path (0) edge[->] node {{}} (00);
\path (0) edge[->] node {{}} (01);

\path (1) edge[->] node {{}} (10);
\path (1) edge[->] node {{}} (11);

\path (00) edge[->] node {{}} (000);
\path (00) edge[->] node {{}} (001);

\path (01) edge[->] node {{}} (010);
\path (01) edge[->] node {{}} (011);

\path (10) edge[->] node {{}} (100);
\path (10) edge[->] node {{}} (101);

\path (11) edge[->] node {{}} (110);
\path (11) edge[->] node {{}} (111);

\path (empty) edge[->] node {{}} (1);

\path (0) edge[dashed,->, bend right=7.5] node {{}} (1);
\path (00) edge[dashed,->, bend right=7.5] node {{}} (01);
\path (00) edge[dashed,->, bend right=10] node {{}} (10);
\path (000) edge[dashed,->, bend right=7.5] node {{}} (001);
\path (000) edge[dashed,->, bend right=12.5] node {{}} (010);
\path (000) edge[dashed,->, bend right=15] node {{}} (100);
 
\end{tikzpicture}
\end{center}
\caption{the infinite complete binary tree with outgoing $R_1$ arrows as in Example \ref{BinaryEx} shown only for states along the leftmost branch. Reflexive accessibility loops are omitted.}\label{BinaryTreeFig}
\end{figure}
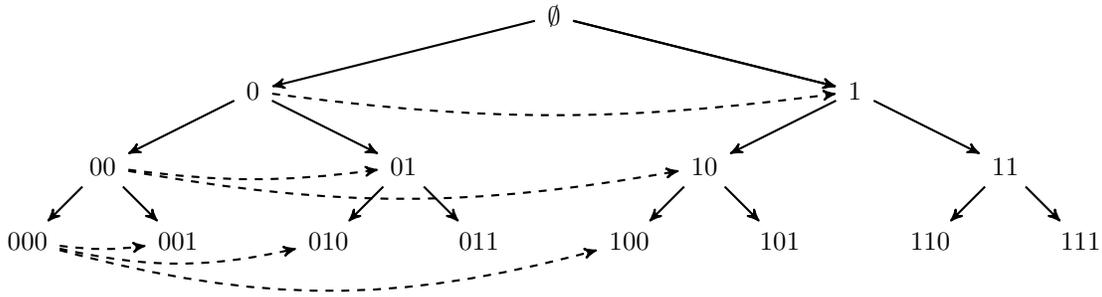

\begin{example}[Accessibility Relations on the Infinite Complete Binary Tree]\label{BinaryEx}
Consider the set $2^{<\omega}$ of all finite binary strings. For $x,x'\in 2^{<\omega}$, let $x'\sqsubseteq x$ iff  $x'$ extends $x$, i.e., $x$ is an initial segment of $x'$. We can view $\langle 2^{<\omega},\sqsubseteq\rangle$ as the infinite complete binary tree with $\sqsubseteq$ as the reflexive transitive closure of the child relation---a simple example of a possibility space in which every possibility can be further refined. Observe that for each $x\in 2^{<\omega}$, $\mathord{\downarrow}x$ satisfies \textit{refinability} and hence $\mathord{\downarrow}x\in \mathrm{RO}(2^{<\omega},\sqsubseteq)$, so $\mathrm{RO}(2^{<\omega},\sqsubseteq)$ is an atomless Boolean algebra.  As an exercise, consider various definitions of accessibility relations on $\langle 2^{<\omega},\sqsubseteq\rangle$ and then check which, if any, of the interplay conditions are satisfied. For example, for $n\in\mathbb{N}$, define an accessibility relation $R_n$ by: $xR_ny$ iff $x$ and $y$ have the same length and differ in no more than $n$ places, i.e., where $x=\langle x_1,\dots,x_k\rangle$ and $y=\langle y_1,\dots,y_k\rangle$, we have $|\{i\mid 1\leq i\leq k,\, x_i\not=y_i\}|\leq n$. See Figure \ref{BinaryTreeFig} for $R_1$. Since $x$ and $y$ must have the same length, $R_n$ does not satisfy  \upR{}, \Rdown{}, or \RrefPlus{}. However, $R_n$ does satisfy \Rcomm{} and \Rwin{}. For \Rwin{}, if $xR_ny$ and $y'=\langle y_1,\dots, y_l\rangle$ is an extension $y=\langle y_1,\dots, y_k\rangle$, let $x'$ be the result of concatenating $\langle y_{k+1},\dots, y_l\rangle$ on to the end of $x$. Then clearly for every extension $x''$ of $x'$, there is an extension $y''$ of $y'$ such that $x''R_n y''$, since we can use the same move of copying and concatenating. Thus, the structure $\mathcal{F}=\langle 2^{<\omega},\sqsubseteq, \{R_n\}_{n\in\mathbb{N}}, \mathrm{RO}(2^{<\omega},\sqsubseteq) \rangle$ is a full possibility frame. Alternatively, suppose we define $R_n^\leq$ by: $xR_n^\leq y$ iff $y$ is at least as long as $x$ and $xR_nz$ where $z$ is the initial segment of $y$ of the same length as $x$. Then $R_n^\leq$ satisfies \upR{}, \Rdown{}, and \Rref{}, but still not \RrefPlus{}.  \hfill $\triangleleft$
\end{example}

Although we have seen that we can assume without loss of generality stronger conditions than \Rrule{} and \Rwinweak{}, we do not build these conditions into the definition of possibility frames. The greatest payoff in terms of simplifying our theory seems to come from assuming the \Rdown{} property as in intuitionistic frames and Humberstone frames. We give frames with this property the following honorific title.

\begin{definition}[Standard Possibility Frames]\label{Standard} A \textit{standard possibility frame} is a possibility frame satisfying \Rdown{}. \hfill $\triangleleft$
\end{definition}

The topic of this section has been the interplay of two \textit{relations}: accessibility and refinement. It is worth mentioning how the situation changes if we consider the possibility semantic analogue not of Kripke frames but of \textit{neighborhood}  frames \citep{Montague1970,Scott1970}, which can characterize non-normal modal logics. 

\begin{remark}[Neighborhood Possibility Frames] A \textit{full neighborhood possibility frame} may be defined as a tuple $\mathcal{F}=\langle S,\sqsubseteq, \{N_i\}_{i\in \ind},\adm\rangle$ where $\langle S,\sqsubseteq\rangle$ is a nonempty poset, $N_i\colon S\to \wp(\adm)$, $\adm=\mathrm{RO}(S,\sqsubseteq)$, and $\sqsubseteq$ and $N_i$ satisfy the following interplay conditions:
\begin{itemize}
\item if $x'\sqsubseteq x$, then $N_i(x')\supseteq N_i(x)$;
\item if $X\not\in N_i(x)$, then $\exists x'\sqsubseteq x$ $\forall x''\sqsubseteq x'$, $X\not\in N_i(x'')$.
\end{itemize}
The standard neighborhood semantics clause for $\Box_i$ now applies:
\begin{itemize}
\item $\mathcal{M},x\Vdash \Box_i\varphi$ iff $\llbracket \varphi\rrbracket^\mathcal{M}\in N_i(x)$.
\end{itemize}
The two interplay conditions above ensure that $\llbracket \Box_i\varphi\rrbracket^\mathcal{M}$ satisfies \textit{persistence} and \textit{refinability} if $\llbracket \varphi\rrbracket^\mathcal{M}$ does.

The logic of any class of neighborhood possibility frames is a \textit{congruential} modal logic, i.e., such that if $\varphi\leftrightarrow\psi\in\mathbf{L}$, then $\Box_i\varphi\leftrightarrow\Box_i\psi\in\mathbf{L}$.  To characterize normal modal logics, we can use \textit{normal} neighborhood possibility frames in which for each $x\in S$, $N_i(x)$ is a \textit{filter} in $\adm$: $N_i(x)\neq\emptyset$, and $X, Y\in N_i(x)$ iff $X\cap Y\in N_i(x)$. We mention this only to state the following fact: full normal neighborhood possibility frames are to \textit{complete} Boolean algebras with operators ($\mathcal{C}$-BAOs) as our full (relational) possibility frames are to \textit{complete} and \textit{completely additive} BAOs ($\mathcal{CV}$-BAOs). (One can also define morphisms between neighborhood possibility frames to play a role parallel to that of our possibility morphisms in \S~\ref{morphisms}, but we will not go into the details here.) The previous fact should make sense after the duality theory of \S\S~\ref{PossToBAO}-\ref{DualEquiv}.\hfill $\triangleleft$
\end{remark}

For easy reference, all of the interplay conditions discussed in this section and elsewhere in the paper are collected in Figure \ref{InterplayTable}.   

\begin{figure}
\begin{center}
\begin{tabular}{lll}
\Rrule{} & if $x'\sqsubseteq x$ and $x'R_iy'\comp z$, then $\exists y$: $xR_i y\comp z$ & \S~\ref{FullFrames} \\ 
\Rcomm{} & if $x'\sqsubseteq x$ and $x'R_iy'$, then $\exists y$: $xR_i y$ and $y'\sqsubseteq y$ & \S~\ref{PSFramesSem}, \S~\ref{FullFrames} \\ 
\upR{} & if $x'\sqsubseteq x$ and $x'R_i y'$, then $xR_i y'$ &  \S~\ref{PSFramesSem}, \S~\ref{FullFrames} \\
\\
\Rdown{} & if $y'\sqsubseteq y$ and $xR_iy$, then $xR_i y'$ &  \S~\ref{PSFramesSem}, \S~\ref{FullFrames} \\
\\
\Rwinweak{} & if $xR_iy$, then $\forall y'\sqsubseteq y$ $\exists x'\sqsubseteq x$ $\forall x''\sqsubseteq x'$ $\exists y''\comp y'$: $x''R_iy''$ & \S~\ref{FullFrames} \\
\Rwin{} & if $xR_iy$, then $\forall y'\sqsubseteq y$ $\exists x'\sqsubseteq x$ $\forall x''\sqsubseteq x'$ $\exists y''\sqsubseteq y'$: $x''R_iy''$ & \S~\ref{FullFrames} \\
\RWin{} & $xR_i y$ iff $\forall y'\sqsubseteq y$ $\exists x'\sqsubseteq x$ $\forall x''\sqsubseteq x'$ $\exists y''\sqsubseteq y'$: $x''R_iy''$ & \S~\ref{FullFrames} \\
\\
\Rref{} & if $xR_iy$, then $\exists x'\sqsubseteq x$ $\forall x''\sqsubseteq x'$ $\underline{\exists y'\sqsubseteq y}$: $x''R_iy'$ & \S~\ref{FullFrames}, \S~\ref{Strengths} \\
\RrefPlus{} & if $xR_iy$, then $\underline{\exists y'\sqsubseteq y}$ $\exists x'\sqsubseteq x$ $\forall x''\sqsubseteq x'$\; $x''R_iy'$ & \S~\ref{FullFrames}, \S~\ref{Strengths}  \\
\RrefPlusPlus{} & if $xR_iy$, then $\exists x'\sqsubseteq x$ $\forall x''\sqsubseteq x'$\; $x''R_iy$ & \S~\ref{FullFrames}, \S~\ref{Strengths} \\
\\
\Rdense{} & $xR_iy$ if $\forall y'\sqsubseteq y$ $\exists y''\sqsubseteq y'$: $xR_i y''$ & \S~\ref{FullFrames}  \\ 
\\
\Rmax{} & if $R_i(x)\not=\emptyset$, then $R_i(x)$ has a maximum in $\langle S, \sqsubseteq\rangle$ & \S~\ref{FuncFrames} \\
\Rmaxe{} & $R_i(x)$ has a maximum in $\langle S, \sqsubseteq\rangle$ & \S~\ref{FuncFrames} \\
\Rprinc{} & if $R_i(x)\not=\emptyset$, then $R_i(x)$ is a principal downset in $\langle S, \sqsubseteq\rangle$ & \S~\ref{FuncFrames} \\
\\
\end{tabular}
\caption{all of the interplay conditions relating accessibility and refinement mentioned in the paper. Each block of conditions is ordered from weaker to stronger conditions, except that \Rmaxe{} and \Rprinc{} are incomparable.}\label{InterplayTable}
\end{center}
\end{figure} 

\subsection{Accessibility and Possibility} \label{Acc&Poss}

So far we have said nothing about the semantics of $\Diamond_i$.  Since we use the classical definition $\Diamond_i\varphi:=\neg\Box_i\neg\varphi$, the semantics of $\Diamond_i$ is derived directly from that of $\neg$ and $\Box_i$ as follows.

\begin{figure}[h]
\begin{center}
\begin{tikzpicture}[->,>=stealth',shorten >=1pt,shorten <=1pt, auto,node
distance=2cm,thick,every loop/.style={<-,shorten <=1pt}]
\tikzstyle{every state}=[fill=gray!20,draw=none,text=black]

\node (x-up) at (5,0) {{$x$}};
\node (x) at (5,-1.75) {{$x'$}};
\node (forall) at (4.75,-.75) {{$\forall$}};
\node (exists) at (5.85,-1.4) {{$\exists$}};
\node (y) at (7,-1.75) {{$y'$}};

\node (Diamond) at (5.7,0) {{$\Vdash\Diamond_i \varphi$}};
\node (phi) at (7.55,-1.75) {{$\Vdash\varphi$}};

\path (x) edge[dashed,->] node {{}} (y);
\path (x-up) edge[->] node {{}} (x);

\node (x-up') at (0,0) {{$x$}};
\node (x') at (0,-1.75) {{$x'$}};
\node (y') at (2,-1.75) {{$y'$}};
\node (y'') at (2,-3.5) {{$y''$}};

\node (forall') at (-.25,-.75) {{$\forall$}};
\node (exists') at (.85,-1.4) {{$\exists$}};
\node (exists'') at (1.75,-2.5) {{$\exists$}}; 

\node (Diamond') at (.7,0) {{$\Vdash\Diamond_i \varphi$}};
\node (phi') at (2.6,-3.5) {{$\Vdash\varphi$}};

\path (x-up') edge[->] node {{}} (x');
\path (x') edge[dashed,->] node {{}} (y');
\path (y') edge[->] node {{}} (y'');

\end{tikzpicture}
\end{center}
\caption{the semantic clause for $\Diamond_i$, without (left) and with (right) the \Rdown{} condition.}\label{DiamondFig}
\end{figure}
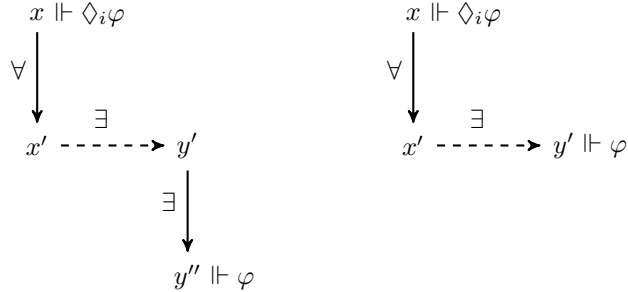

\begin{fact}[Forcing $\Diamond_i$]\label{ForcingDiamond} Given a partial-state frame $\mathcal{F}=\langle S, \sqsubseteq , \{R_i\}_{i\in\ind},\adm\rangle$, $x\in S$, and $Y\subseteq S$, define:
\begin{enumerate}[label=\arabic*.,ref=\arabic*]
\item $x\in\blacklozenge_i Y$ iff $\forall x'\sqsubseteq x$ $\exists y'$: $x'R_iy'$ and $\exists y''\sqsubseteq y'$: $y''\in Y$.
\end{enumerate}
Then for any possibility model $\mathcal{M}$ and $\varphi\in\mathcal{L}(\sig,\ind)$, $\llbracket \Diamond_i\varphi\rrbracket^\mathcal{M}=\blacklozenge_i \llbracket\varphi\rrbracket^\mathcal{M}$, i.e.:
\begin{enumerate}[label=\arabic*.,ref=\arabic*,resume]
\item\label{ForcingDiamond2} $\mathcal{M},x\Vdash\Diamond_i\varphi$ iff $\forall x'\sqsubseteq x$ $\exists y'$: $x'R_iy'$ and $\exists y''\sqsubseteq y'$: $\mathcal{M},y''\Vdash\varphi$ (see Figure \ref{DiamondFig}, left).
\end{enumerate}
For standard frames satisfying \Rdown{} as in \S~\ref{FullFrames}, these conditions simplify to:
\begin{enumerate}[label=\arabic*.,ref=\arabic*,resume]
 \item $x\in\blacklozenge_i Y$ iff $\forall x'\sqsubseteq x$ $\exists y'$: $x'R_iy'$ and $y'\in Y$;
\item $\mathcal{M},x\Vdash\Diamond_i\varphi$ iff $\forall x'\sqsubseteq x$ $\exists y'$: $x'R_iy'$ and $\mathcal{M},y'\Vdash\varphi$ (see Figure \ref{DiamondFig}, right).
\end{enumerate}
\end{fact}

Although unfamiliar, the clause for $\Diamond_i$ is quite intuitive: for whatever sense of `possible' is at issue, the clause says that $x$ \textit{forces} that $\varphi$ is \textit{possible} iff for every refinement $x'$ of $x$ (think of this as the forcing part), $x'$ has access to a state that forces $\varphi$, or if we are not assuming \Rdown{}, then $x'$ has access to state that can be refined to force $\varphi$ (think of this as the possibility part).

It is useful to know some shortcuts for thinking about sequences of modal operators involving diamonds. Since $\Diamond_{i_1}\dots\Diamond_{i_n}\varphi$ is equivalent to $\neg\Box_{i_1}\dots \Box_{i_n}\neg\varphi$, we have $\mathcal{M},x\Vdash \Diamond_{i_1}\dots\Diamond_{i_n}\varphi$ iff $\forall x'\sqsubseteq x$ $\exists y'_1,\dots,y'_n$: $x'R_{i_1}y'_1\dots y'_{n-1}R_{i_n}y'_n$ and $\exists y''_n\sqsubseteq y'_n$: $\mathcal{M},y''_n\Vdash \varphi$. Over standard frames satisfying \Rdown{}, we can simplify this further as follows.

\begin{fact}[Iterated Modalities]\label{ItMod} For any possibility model $\mathcal{M}$ based on a standard possibility frame and $n\geq 1$:
\begin{enumerate}
\item\label{ItMod1} $\mathcal{M},x\Vdash \Diamond_{i_1}\dots\Diamond_{i_n}\varphi$ iff $\forall x'\sqsubseteq x$ $\exists y'_1,\dots,y'_n$: $x'R_{i_1}y'_1\dots y'_{n-1}R_{i_n}y'_n$ and $\mathcal{M},y'_n\Vdash \varphi$;
\item\label{ItMod2} $\mathcal{M},x\Vdash \Box_j \Diamond_{i_1}\dots\Diamond_{i_n}\varphi$ iff $\forall y$, if $xR_j y$, then $\exists y_1,\dots,y_n$: $yR_{i_1}y_1\dots y_{n-1}R_{i_n}y_n$ and $\mathcal{M},y_n\Vdash \varphi$.
\end{enumerate}
\end{fact}
\begin{proof} For part \ref{ItMod1}, in the truth condition for $\neg\Box_{i_1}\dots \Box_{i_n}\neg\varphi$ given above, together $y'_{n-1}R_{i_n}y'_n$ and $\exists y''_n\sqsubseteq y'_n$ imply $y'_{n-1}R_{i_n} y''_n$ by the \Rdown{} property of standard frames.

For part \ref{ItMod2}, by part \ref{ItMod1} we have that $\mathcal{M},x\Vdash \Box_j \Diamond_{i_1}\dots\Diamond_{i_n}\varphi$ iff $\forall y$, if $xR_jy$, then $\forall y'\sqsubseteq y$ $\exists y'_1,\dots,y'_n$: $y'R_{i_1}y'_1\dots y'_{n-1}R_{i_n}y'_n$ and $\mathcal{M},y'_n\Vdash \varphi$. This clearly implies the clause given in part \ref{ItMod2}, and the converse implication also holds due to \Rdown{}: for $xR_jy$ and $y'\sqsubseteq y$ together imply $xR_jy'$.
\end{proof}

Fact \ref{ItMod} shows that for any sequence $M_1\dots M_n$ of modal operators beginning with a box, we can think of the semantic clause for $M_1\dots M_n \, p$  in standard models exactly as in Kripke semantics. If $M_1$ is a diamond, we must add an initial $\forall x'\sqsubseteq x$, but then the rest of the clause is as in Kripke semantics.

Note that we cannot conclude from $xR_iy$ and $\mathcal{M},y\Vdash\varphi$ that $\mathcal{M},x\Vdash \Diamond_i\varphi$. This shows that in possibility semantics, $y$ being \textit{accessible to} $x$, which guarantees that $\mathcal{M},x\Vdash \Box_i\varphi\Rightarrow \mathcal{M},y\Vdash\varphi$, is not the same as $y$ being \textit{possible relative to} $x$, in the sense that would guarantee $\mathcal{M},y\Vdash\varphi\Rightarrow \mathcal{M},x\Vdash \Diamond_i\varphi$ (cf.~\citealt[p.~327f]{Humberstone1981}). We can, however, define a relation of relative possibility that will guarantee the latter implication, as in Remark \ref{RelPoss}.

\begin{remark}[Relative Possibility]\label{RelPoss}  Given a partial-state frame $\mathcal{F}=\langle S, \sqsubseteq , \{R_i\}_{i\in\ind},\adm\rangle$ and $x,y\in S$, define 
\begin{itemize} 
\item $xR_{i\Diamond}y$ iff $x\in \blacklozenge_i \mathord{\downarrow}y$,
\end{itemize}
or equivalently, 
\begin{itemize}
\item $xR_{i\Diamond}y$ iff  $\forall x'\sqsubseteq x$ $\exists y'$: $x'R_iy'\comp y$,
\end{itemize}
which for standard frames simplifies to
\begin{itemize}
\item $xR_{i\Diamond}y$ iff  $\forall x'\sqsubseteq x$ $\exists y'$: $x'R_iy'\sqsubseteq y$.
\end{itemize}
Let us make three observations concerning the relation $R_{i\Diamond}$.

First, if $xR_{i\Diamond}y$, then for any downset $Z\in\adm$, $y\in Z\Rightarrow x\in\blacklozenge_i Z$. Thus, assuming each $Z\in \adm$ is a downset,  for any model $\mathcal{M}$ based on $\mathcal{F}$, if $xR_{i\Diamond}y$, then for any $\varphi\in\mathcal{L}(\sig,\ind)$, $\mathcal{M},y\Vdash\varphi\Rightarrow\mathcal{M},x\Vdash \Diamond_i\varphi$.

Second, using the $R_{i\Diamond}$ relation, we can rewrite \Rwinweak{} from \S~\ref{FullFrames} equivalently as follows:
\begin{itemize}
\item \Rwinweak{} -- if ${x}R_i{y}$, then $\forall {y'}\sqsubseteq {y}$ $\exists{x'}\sqsubseteq {x}$: $x'R_{i\Diamond}y'$.
\end{itemize}
Thus, \Rwinweak{} relates the fundamental notions of \textit{accessibility} and \textit{relative possibility}: if $y$ is \textit{accessible to} $x$, then for every refinement $y'$ of $y$ there is a refinement $x'$ of $x$ such that $y'$ is \textit{possible relative to} $x'$.

Third, assuming only $\mathcal{M},x\Vdash\Diamond_i\varphi$, we cannot conclude that there is a $y$ such that $xR_{i\Diamond }y$ and $\mathcal{M},y\Vdash\varphi$. The reason is that a partial state $x$ might determine that $\varphi$ is possible without yet determining any \textit{particular} witness $y$ for the possibility of $\varphi$. More formally: it may be that for every state $y$ that forces $\varphi$, $x$ can still be refined to an $x_y$ that ``rules out'' that particular $y$, in the sense that every state in $R_i(x_y)$ is incompatible with $y$, i.e., $\forall y\in\llbracket\varphi\rrbracket^\mathcal{M}$ $\exists x_y\sqsubseteq x$ $\forall z\in R_i(x_y)$, $z\incomp y$. This is consistent with $\mathcal{M},x\Vdash\Diamond_i\varphi$. \hfill $\triangleleft$
\end{remark}

In \S~\ref{NoKripke}, we will see the semantics for $\Diamond_i$ in action in a concrete example.

\subsection{Full Possibility Frames with No Kripke Equivalents}\label{NoKripke}

We will conclude \S~\ref{FromPartToPoss} by constructing a full possibility frame $\mathcal{F}$ that validates a modal formula that is not valid over any Kripke frame. Thus, the logic of $\mathcal{F}$ will be a normal modal logic that is \textit{Kripke-frame inconsistent}---it is not sound with respect to any Kripke frame---and hence \textit{Kripke-frame incomplete}---it is not sound and complete with respect to any class of Kripke frames. From $\mathcal{F}$ we will generate continuum many full possibility frames with distinct Kripke-inconsistent logics. This requires a polymodal language, since every syntactically consistent normal unimodal logic is Kripke-frame consistent \citep{Makinson1971}. In \S~\ref{CompFull}, we will generate continuum many full possibility frames with distinct Kripke-\textit{incomplete} \textit{uni}modal logics.

Suppose $\varphi$ and $\psi$ are modal formulas such that the propositional variable $p$ does not occur in $\psi$ (or at least it only occurs in the form $\top:=(p\vee\neg p)$). Then consider the following formula:
\begin{equation}\Diamond_i (p\wedge\psi)\rightarrow \big(\Diamond_i(p\wedge \varphi)\wedge \Diamond_i(p\wedge\neg\varphi)\big)\tag{\textsc{Split}}.\label{Splitting}\end{equation}
We claim that any Kripke frame $\mathfrak{F}$ that validates (\ref{Splitting}) must also validate $\neg\Diamond_i\psi$. Suppose $\neg\Diamond_i\psi$ is not valid over $\mathfrak{F}$, so there is a model $\mathfrak{M}$ based on $\mathfrak{F}$ and a world $w$ such that $\mathfrak{M},w\vDash \Diamond_i\psi$. Hence there is some world $v$ such that $w\mathrm{R}_iv$ and $\mathfrak{M},v\vDash \psi$. Let $\mathfrak{M}'$ be the model based on $\mathfrak{F}$ that differs from $\mathfrak{M}$ only in that $\llbracket p\rrbracket^{\mathfrak{M}'}=\{v\}$. Then since $\psi$ does not contain $p$, we still have $\mathfrak{M}', v\vDash \psi$. Now observe that the antecedent of (\ref{Splitting}) is true at $w$ in $\mathfrak{M}'$, but clearly the consequent of (\ref{Splitting}) is false since $\llbracket p\rrbracket^{\mathfrak{M}'}$ is a singleton set. 

Worlds cannot split, but possibilities can. We will construct a full possibility frame that validates an instance of (\ref{Splitting}) \textit{and} $\Diamond_i\psi$.\footnote{Since we saw that this conjunction is Kripke-inconsistent, the example requires multiple modal operators by the point above about Makinson's theorem. Another way to see this (thanks here to Lloyd Humberstone) is that no syntactically consistent normal unimodal logic contains a pair of formulas of the form $\Diamond \alpha$ and $\Diamond\neg\alpha$ (see \citealt{French2015}), but any normal modal logic containing (\ref{Splitting}) and $\Diamond_i\psi$ contains both $\Diamond\varphi(\top/p)$ and $\Diamond\neg\varphi(\top/p)$ (substitute $\top$ for $p$ in (\ref{Splitting})).} The construction is a possibility frame version of the construction in \citealt{Litak2005} of a $\mathcal{CV}$-BAO that generates a variety of BAOs with no atomic members. (We will precisely relate full possibility frames and $\mathcal{CV}$-BAOs in \S\S~\ref{PossToBAO}-\ref{DualEquiv}.) The main idea of the following is exactly as in \citealt{Litak2005} (also cf.~\citealt{Venema2002}), but where Litak defines operators on an algebra, we define relations on a frame. One may compare the two constructions to see the relative benefits of thinking in terms of relations vs. operators.

Consider the standard topology on $(0,1)$ generated by the basis of open intervals $(a,b)$. Let $\mathrm{RO}(0,1)$ be the set of regular open sets. Recall the fact that a subset of $(0,1)$ is open iff it is the union of a countable set of pairwise disjoint open intervals. A subset of $(0,1)$ is regular open only if it is the union of a countable set of pairwise disjoint and \textit{non-adjacent} open intervals---for if adjacent intervals $(a,b)$ and $(b,c)$ are in the set, then the interior of the closure of the union will contain $b$, which is not in the union. For any regular open $O$, since
\[O=\bigcup \{(a,b)\mid (a,b)\subseteq O \mbox{ and }\neg\exists (a',b')\colon (a,b)\subsetneq (a',b')\subseteq O\},\]
we can canonically ``encode'' $O$ as the following set of pairs:
\begin{equation}\sigma_O=\{\langle a,b\rangle\mid (a,b)\subseteq O\mbox{ and }\neg \exists (a',b')\colon (a,b)\subsetneq (a',b')\subseteq O\}.\label{EncodeEq}\end{equation}
This gives us a convenient way of shrinking a given regular open set $O$, as follows:
\begin{equation}O_{-}=\bigcup\{(a,b)\mid  \langle a,b+\frac{|a-b|}{2}\rangle\in \sigma_O\}.\label{Shrink}\end{equation}
If we consider $\langle \mathrm{RO}(0,1)\setminus\{\emptyset\},\subseteq\rangle$ as a possibility space, then since $O_{-}\subsetneq O$ and $O_{-}$ is regular open, the possibility $O_{-}$ is a strict refinement of $O$. As in Remark \ref{Persp2}, the regular open sets of any topology, ordered by inclusion, form a complete Boolean lattice. In addition, $\langle \mathrm{RO}(0,1),\subseteq\rangle$ is \textit{atomless}. Thus, in $\langle \mathrm{RO}(0,1)\setminus\{\emptyset\},\subseteq\rangle$, every possibility can be further refined. This is the key to validating an instance of~(\ref{Splitting}).

Building a possibility frame on a complete Boolean lattice---minus the bottom element---makes it easy to deal with the set $\adm$ of admissible propositions. For if $\langle S,\sqsubseteq\rangle$ is a poset obtained from a complete Boolean lattice by deleting the bottom element, then the regular open sets in the downset topology on $\langle S,\sqsubseteq\rangle$ are just $\emptyset$ and each \textit{principal downset} $\mathord{\downarrow}x=\{x'\in S\mid x'\sqsubseteq x\}$ for $x\in S$ (see Fact \ref{FullTFAE}). Thus, the regular open sets in the downset topology on $\langle \mathrm{RO}(0,1)\setminus\{\emptyset\},\subseteq\rangle$ are $\emptyset$ and each $\mathord{\downarrow}O=\{O'\in\mathrm{RO}(0,1)\setminus\{\emptyset\}\mid O'\subseteq O\}$.

We noted above how a nonempty regular open set can be canonically encoded as a nonempty set of pairs of real numbers. We will take as the possibilities of our frame $\mathcal{F}$ not only every nonempty regular open set of reals, but also every nonempty set of pairs of reals $\langle a,b\rangle$ with $0\leq a<b\leq 1$. Thus, where $[0,1]^2_<={\{\langle a,b\rangle \in [0,1]^2\mid a<b\}}$, let
\[S= (\mathrm{RO}(0,1) \cup \wp ([0,1]^2_<))\setminus \{\emptyset\}.\]
We will write `$O$', `$O'$', etc., for elements of $\mathrm{RO}(0,1)\setminus\{\emptyset\}$, and `$\sigma$', `$\sigma'$', `$\tau$', etc., for elements of $\wp([0,1]^2_<)\setminus\{\emptyset\}$.

The refinement relation $\sqsubseteq$ in our frame $\mathcal{F}$ will be $\subseteq$. So regular open sets can refine regular open sets, and sets of pairs can refine sets of pairs, but regular open sets cannot refine sets of pairs or vice versa. Since $\langle \wp([0,1]^2_<),\subseteq\rangle$ is a complete (and atomic) Boolean lattice, we have taken a disjoint union of two complete Boolean lattices---minus the bottom elements. Again, this makes it easy to deal with the set $\adm$ of admissible propositions in our full possibility frame $\mathcal{F}$. The regular open sets in the downset topology on our $\langle S,\sqsubseteq\rangle$ are just $\emptyset$ and for each $O\in\mathrm{RO}(0,1)\setminus\{\emptyset\}$ and $\sigma\in \wp([0,1]^2_<)\setminus\{\emptyset\}$, the sets $\mathord{\downarrow} O$, $\mathord{\downarrow}\sigma$, and $\mathord{\downarrow}O\cup\mathord{\downarrow}\sigma$.

All that remains is to define the accessibility relations in $\mathcal{F}$. Let $R_i$ be the universal relation on $S$. Before defining the other relations, let us explain the strategy to ensure that $\mathcal{F}$ validates an instance of  (\ref{Splitting}).

\textit{Strategy}. We will define a formula $\psi$ such that $\llbracket \psi\rrbracket^\mathcal{M}=\mathrm{RO}(0,1)\setminus\{\emptyset\}$ for any model $\mathcal{M}$ based on $\mathcal{F}$. Thus, $\Diamond_i(p\wedge\psi)$ will say that $p$ is true at some regular open set, which with our observation above about the admissible sets in $\adm$ implies that there is some regular open $O$ such that $\llbracket p\rrbracket^\mathcal{M}\cap\mathrm{RO}(0,1)=\mathord{\downarrow}O$. Then we will define a formula $\varphi $ such that for any model $\mathcal{M}$ based on $\mathcal{F}$, if $\llbracket p\rrbracket^\mathcal{M}\cap\mathrm{RO}(0,1)=\mathord{\downarrow}O$, then $\llbracket \varphi \rrbracket^\mathcal{M}=\mathord{\downarrow} O_{-}$ for $O_{-}$ as in (\ref{Shrink}).  Then since $\emptyset\not=\mathord{\downarrow} O_{-}\subsetneq \mathord{\downarrow}O\subseteq \llbracket p\rrbracket^\mathcal{M}$, we have $\emptyset\not=\llbracket \varphi \rrbracket^\mathcal{M}\subsetneq \llbracket p\rrbracket^\mathcal{M}$. Hence there is an $x$ with $\mathcal{M},x\Vdash p\wedge\varphi $ and a $y$ with $\mathcal{M},y\Vdash p$ and $\mathcal{M},y\nVdash \varphi $, which by Refinability and Persistence implies that there is a $y'\sqsubseteq y$ with $\mathcal{M},y'\Vdash p\wedge \neg \varphi $. Thus, the formula $\Diamond_i (p\wedge\psi)\rightarrow \big(\Diamond_i (p\wedge\varphi )\wedge\Diamond_i (p\wedge\neg\varphi ) \big)$ will be valid over~$\mathcal{F}$. 

Our strategy to define $\varphi$ will be to define a polymodal formula $\alpha$ and appropriate accessibility relations of $\mathcal{F}$ for the operators in $\alpha$ such that if $\llbracket p\rrbracket^\mathcal{M}\cap \mathrm{RO}(0,1)=\mathord{\downarrow} O$, then $\llbracket \alpha\rrbracket^\mathcal{M}=\mathord{\downarrow}\sigma_{O_{-}}$. We will also define an accessibility relation $R_\triangleright$ such that $\llbracket \alpha\rrbracket^\mathcal{M}=\mathord{\downarrow}\sigma_{O_{-}}$ implies $\llbracket \Diamond_\triangleright\alpha\rrbracket^\mathcal{M}=\mathord{\downarrow}O_{-}$. We can then take $\varphi:=\Diamond_\triangleright\alpha$.

In addition to the universal relation $R_i$, we will define four other accessibility relations for $\mathcal{F}$. The relation  $R_\triangleright$ just mentioned only relates regular open sets to sets of pairs:
\begin{itemize} 
\item $OR_\triangleright\sigma$ iff $\forall \langle a',b'\rangle\in \sigma$ $\exists (a,b)\subseteq O$:  $(a,b)\subseteq (a',b')$.
\end{itemize}
For any $O\in \mathrm{RO}(0,1)\setminus\{\emptyset\}$,  $OR_\triangleright\sigma_O$ for $\sigma_O$ as in (\ref{EncodeEq}); and $xR_\triangleright y$ only if $x\in\mathrm{RO}(0,1)\setminus\{\emptyset\}$. So for any model $\mathcal{M}$ based on our frame, we will have $\llbracket \Diamond_\triangleright\top\rrbracket^\mathcal{M}=\mathrm{RO}(0,1)\setminus\{\emptyset\}$. Thus, we can take $\psi$ in the \textit{Strategy} to be~$\Diamond_\triangleright\top$.

Since we want our frame $\mathcal{F}$ to be a \textit{full} possibility frame, we will check that each accessibility relation we define for $\mathcal{F}$ satisfies the \Rrule{} and \Rwinweak{} conditions from \S~\ref{FullFrames} (recall Proposition \ref{ROtoRO}). In fact, we will show that each relation $R$ satisfies the following stronger set of conditions from \S~\ref{FullFrames}:
\begin{itemize}
\item \upR{} -- if $x'\subseteq x$ and $x'Ry'$, then $xRy'$; 
\item \Rdown{} -- if $y'\subseteq y$ and $xRy$, then $xRy'$;
\item \RrefPlus{} -- if $xRy$, then $\exists y'\subseteq y$ $\exists x'\subseteq x$ $\forall x''\subseteq x'$, $x''Ry'$.
\end{itemize}
Clearly the relation $R_\triangleright$ defined above satisfies \upR{} and \Rdown{}. For \RrefPlus{}, given $OR_\triangleright\sigma$, let $\sigma'=\{\langle a',b'\rangle\}$ for one of the $\langle a',b'\rangle\in\sigma$. Then $OR_\triangleright\sigma$ implies that there is an $O'=(a,b)\subseteq O$ such that $O'\subseteq (a',b')$. Then for all $O''\subseteq O'$, $O'' R_\triangleright \sigma'$. This establishes \RrefPlus{}.

The following lemma is the motivation for defining $R_\triangleright$ as above. Informally, it says that $\Diamond_\triangleright$ can take us from the canonical encoding $\sigma_O$ of a regular open set $O$ back to $O$ itself. Although we have not yet defined all of $\mathcal{F}$, the lemma holds no matter what further accessibility relations we define.

\begin{lemma}\label{ExLem1} For any formula $\chi$ and model $\mathcal{M}$ based on the frame $\mathcal{F}$, if $\llbracket \chi\rrbracket^\mathcal{M}=\mathord{\downarrow}\sigma_O$, then $\llbracket \Diamond_\triangleright \chi\rrbracket^\mathcal{M}=\mathord{\downarrow} O$.
\end{lemma}
\begin{proof}Suppose $\llbracket \chi\rrbracket^\mathcal{M}=\mathord{\downarrow}\sigma_O$. First, we show $\mathcal{M},O\Vdash \Diamond_\triangleright\chi$, which implies $\llbracket \Diamond_\triangleright \chi\rrbracket^\mathcal{M}\supseteq\mathord{\downarrow} O$ by Persistence. Recall from \S~\ref{Acc&Poss} that $\mathcal{M},O\Vdash \Diamond_\triangleright\chi$ if $\forall O'\subseteq O$ $\exists \sigma'$: $O'R_\triangleright\sigma'$ and $\mathcal{M},\sigma'\Vdash\chi$. Given $O'\subseteq O$, let 
\[\sigma'=\{\langle a',b'\rangle\mid \langle a',b'\rangle\in\sigma_O\mbox{ and }\exists (a,b)\subseteq O'\colon (a,b)\subseteq (a',b')\}.\] 
Observe that $\sigma'\not=\emptyset$ and $O'R_\triangleright\sigma'$. Since $\sigma'\subseteq\sigma_O$ and $\llbracket \chi\rrbracket^\mathcal{M}=\mathord{\downarrow}\sigma_O$, we have $\mathcal{M},\sigma'\Vdash \chi$. Thus, $\mathcal{M},O\Vdash \Diamond_\triangleright\chi$.

Next, sets of pairs do not force $\Diamond_\triangleright\chi$, since they have no $R_\triangleright$-successors. We also show that if $O'\not\subseteq O$, then $\mathcal{M},O'\nVdash \Diamond_\triangleright\chi$, so $\llbracket \Diamond_\triangleright \chi\rrbracket^\mathcal{M}\subseteq\mathord{\downarrow} O$.  If $O'\not\subseteq O$, then since $O'$ is open and $O$ is regular open, it follows that there is an $(a,b)\subseteq O'$ such that $(a,b)\cap O=\emptyset$.  We claim that for any $\sigma$ with $(a,b)R_\triangleright\sigma$, $\mathcal{M},\sigma\nVdash \chi$. Then since $(a,b)\subseteq O'$, we have $\mathcal{M},O'\nVdash \Diamond_\triangleright \chi$ by the truth clause for $\Diamond_\triangleright$ given \Rdown{} (recall \S~\ref{Acc&Poss}). To prove the claim, suppose $(a,b)R_\triangleright\sigma$. Since $\llbracket \chi\rrbracket^\mathcal{M}=\mathord{\downarrow}\sigma_O$, $\mathcal{M},\sigma\Vdash \chi$ only if $\sigma\subseteq\sigma_O$. Suppose for reductio that $\sigma\subseteq \sigma_O$, so there is an $\langle a',b'\rangle\in\sigma$ such that $(a',b')\subseteq O$. Together $\langle a',b'\rangle\in\sigma$ and  $(a,b)R_\triangleright\sigma$ imply  $(a,b)\cap (a',b')\not=\emptyset$ by the definition of $R_\triangleright$. But since $(a',b')\subseteq O$, $(a,b)\cap (a',b')\not=\emptyset$ contradicts $(a,b)\cap O=\emptyset$ from above. Thus, $\sigma\not\subseteq\sigma_O$, which completes the proof.\end{proof}

Next, we define a relation that only relates sets of pairs to regular open sets: 
\begin{itemize}
\item $\sigma R_\triangleleft O$ iff $\forall (a,b)\subseteq O$ $\exists \langle a',b'\rangle\in\sigma$: $(a,b)\subseteq (a',b')$.
\end{itemize}
Note that since each $O$ is a union of open intervals, $\sigma R_\triangleleft O$ implies $O\subseteq \underset{\langle a,b\rangle\in\sigma}\bigcup (a,b)$. 

Clearly $R_\triangleleft$ satisfies \upR{} and \Rdown{}. For \RrefPlus{}, given $\sigma R_\triangleleft O$, let $O'=(a,b)$ for one of the $(a,b)\subseteq O$. Then $\sigma R_\triangleleft O$ implies that there is an $\langle a',b'\rangle\in\sigma$ such that $O'\subseteq (a',b')$. Let $\sigma'=\{\langle a',b'\rangle\}$. Then for all nonempty $\sigma''\subseteq\sigma'$, i.e, $\sigma''=\sigma'$, we have $\sigma'' R_\triangleleft O'$. This establishes \RrefPlus{}.

The following lemma is the motivation for defining $R_\triangleleft$ as above.

\begin{lemma}\label{SubLem} For any model $\mathcal{M}$ based on the frame $\mathcal{F}$, if $\llbracket p\rrbracket^\mathcal{M}\cap\mathrm{RO}(0,1)=\mathord{\downarrow}O$, then $\mathcal{M},\{\langle a,b\rangle\}\Vdash \Box_\triangleleft p$ iff $(a,b)\subseteq O$. 
\end{lemma}
\begin{proof} Suppose $\mathcal{M},\{\langle a,b\rangle\}\nVdash \Box_\triangleleft p$, so there is an $O'$ with $\{\langle a,b\rangle\} R_\triangleleft O'$ and $\mathcal{M},O'\nVdash p$. Since $\mathcal{M},O'\nVdash p$ and $\llbracket p\rrbracket^\mathcal{M}\cap\mathrm{RO}(0,1)=\mathord{\downarrow}O$, we have $O'\not\subseteq O$. Since $\{\langle a,b\rangle\} R_\triangleleft O'$, we have $O'\subseteq (a,b)$, which with $O'\not\subseteq O$ implies $(a,b)\not\subseteq O$. Conversely, if $(a,b)\not\subseteq O$, so $\mathcal{M},(a,b)\nVdash p$, then since $\{\langle a,b\rangle\} R_\triangleleft (a,b)$, $\mathcal{M},\{\langle a,b\rangle\}\nVdash \Box_\triangleleft p$.
\end{proof} 

Finally, we define two relations that only relate sets of pairs to sets of pairs:
\begin{itemize}
\item $\sigma R_\subsetneq \tau$ iff $\tau$ is a singleton $\{\langle c,d\rangle\}$ and $\exists \langle a,b\rangle \in \sigma$: $(a,b)\subsetneq (c,d)$;
\item $\sigma R_+\tau$ iff  $\tau$ is a singleton $\{\langle c,d\rangle\}$ and $\exists \langle a,b\rangle \in \sigma$: $\langle a,b+\frac{|a-b|}{2}\rangle=\langle c,d\rangle$.
\end{itemize}
Both relations clearly satisfy \upR{} and \Rdown{}. For \RrefPlus{}, if $\sigma R_\subsetneq \tau$, so $\tau$ is a singleton $\{\langle c,d\rangle\}$ and there is an $\langle a,b\rangle\in\sigma$ with $(a,b)\subsetneq (c,d)$, let $\sigma'=\{\langle a,b\rangle\}$. Then  for all nonempty $\sigma''\subseteq \sigma'$, i.e., $\sigma''=\sigma'$, we have $\sigma'' R_\subsetneq \tau$. So $R_\subsetneq$ satisfies \RrefPlus{}, and the same form of argument applies to $R_+$.\footnote{In fact, this shows that $R_\subsetneq $ and $R_+$ satisfy the stronger \RrefPlusPlus{} condition from Remark \ref{HumbFrame}.}

We have now shown that $\mathcal{F}=\langle S,\sqsubseteq, \{R_i,R_\triangleright,R_\triangleleft,R_\subsetneq, R_+\}, \mathrm{RO}(S,\sqsubseteq)\rangle$ is a full possibility frame.

The following lemma is the motivation for defining $R_\subsetneq$ and $R_+$ as above. Informally, it says that we can write a formula that takes us from a regular open  $O$ to the canonical encoding $\sigma_{O_{-}}$ of the shrunken $O_{-}\subsetneq O$.

\begin{lemma}\label{ExLem2} Let $\alpha:=\Diamond_+\top\wedge \Box_+(\Box_\triangleleft p\wedge \neg\Diamond_\subsetneq \Box_\triangleleft p)$. For any model $\mathcal{M}$ based on $\mathcal{F}$, if $\llbracket p\rrbracket^\mathcal{M}\cap \mathrm{RO}(0,1)=\mathord{\downarrow} O$, then $\llbracket \alpha\rrbracket^\mathcal{M}=\mathord{\downarrow}\sigma_{O_{-}}$.
\end{lemma}
\begin{proof} By Fact \ref{ForcingDiamond}.\ref{ForcingDiamond2} and the definition of $R_+$, $\mathcal{M},x\Vdash\Diamond_+\top$ iff $x=\sigma$ for some $\sigma\subseteq [0,1]^2_<$ such that for every $\langle a,b\rangle\in \sigma$, $\langle a,b+\frac{|a-b|}{2}\rangle\in [0,1]^2_{<}$. Now consider $\mathcal{M},\sigma\Vdash \Box_+(\Box_\triangleleft p\wedge \neg\Diamond_\subsetneq \Box_\triangleleft p)$. By definition of $R_+$, together with $\mathcal{M},\sigma\Vdash\Diamond_+\top$, this is equivalent to: for all $\langle a,b\rangle\in\sigma$, $\mathcal{M},\{\langle a,b+\frac{|a-b|}{2}\rangle \}\Vdash\Box_\triangleleft p\wedge \neg\Diamond_\subsetneq \Box_\triangleleft p $. By Lemma~\ref{SubLem}, ${\mathcal{M},\{\langle a,b+\frac{|a-b|}{2}\rangle \}\Vdash \Box_\triangleleft p}$ iff $(a,b+\frac{|a-b|}{2})\subseteq O$. Also observe that $\mathcal{M},\{\langle a,b+\frac{|a-b|}{2}\rangle \}\Vdash \neg\Diamond_\subsetneq \Box_\triangleleft p$ iff there is no $(c,d)$ such that $(a,b+\frac{|a-b|}{2})\subsetneq (c,d)\subseteq O$. Putting all of this together with the definition of $\sigma_O$ from (\ref{EncodeEq}), we have ${\mathcal{M},\sigma\Vdash \Box_+(\Box_\triangleleft p\wedge \neg\Diamond_\subsetneq \Box_\triangleleft p)}$ iff for each $\langle a,b\rangle\in\sigma$, $\langle a,b+\frac{|a-b|}{2}\rangle\in \sigma_O$. It follows by the definition of $O_{-}$ in (\ref{Shrink}) that $\mathcal{M},\sigma\Vdash \alpha$ iff $\sigma\subseteq \sigma_{O_{-}}$.\end{proof}

Together Lemmas \ref{ExLem1} and \ref{ExLem2} immediately imply the following final piece of the argument.

\begin{proposition}\label{ValProp} For the formula $\alpha$ from Lemma \ref{ExLem2} and any model $\mathcal{M}$ based on $\mathcal{F}$, if $\llbracket p\rrbracket^\mathcal{M}\cap \mathrm{RO}(0,1)=\mathord{\downarrow} O$, then $\llbracket \Diamond_\triangleright \alpha \rrbracket^\mathcal{M}=\mathord{\downarrow}O_{-}$.
\end{proposition}

Where $\varphi :=\Diamond_\triangleright\alpha$ and $\psi:=\Diamond_\triangleright\top$, based on the \textit{Strategy} outlined above we have completed the proof that our full possibility frame $\mathcal{F}$ validates $\Diamond_i (p\wedge\psi)\rightarrow \big(\Diamond_i (p\wedge\varphi )\wedge\Diamond_i (p\wedge\neg\varphi  )\big)$.

From this one example of a full possibility frame whose logic is Kripke-frame inconsistent, we can easily generate continuum many others, as in Theorem~\ref{NoEquiv}. The proof of Theorem~\ref{NoEquiv} uses the fact, proved in \S~\ref{OpFrameAlg}, that given any two full possibility frames $\mathcal{H}$ and $\mathcal{G}$, there is a full possibility frame $\mathcal{H}\biguplus\mathcal{G}$, the \textit{disjoint union} of $\mathcal{H}$ and $\mathcal{G}$, such that a formula is valid over $\mathcal{H}\biguplus\mathcal{G}$ iff it is valid over both $\mathcal{H}$ and $\mathcal{G}$.

\begin{theorem}[Polymodal Possibility Frames with No Kripke Equivalents]\label{NoEquiv} There are continuum many full possibility frames for the polymodal language above whose logics are pairwise distinct and Kripke-frame inconsistent.
\end{theorem}
 
\begin{proof} We first prove the weaker claim that there are continuum many full possibility frames whose logics are pairwise distinct and Kripke-frame \textit{incomplete}. We then show how to modify the argument to obtain the stated theorem for Kripke-frame \textit{inconsistent} logics.

Our frame $\mathcal{F}$ above is a frame for the language $\mathcal{L}(\sig,\{i,\triangleright,\triangleleft,\subsetneq,+ \})$.  Take any continuum-sized set $\{\mathcal{G}_j\}_{j\in J}$ of Kripke frames for the unimodal language $\mathcal{L}(\sig,\{+\})$, viewed as full possibility frames, such that the logics of the $\mathcal{G}_j$ are pairwise distinct. That such a set of frames exists is a standard fact. Extend each $\mathcal{G}_j$ to a frame $\mathcal{G}'_j$ for  $\mathcal{L}(\sig,\{i,\triangleright,\triangleleft,\subsetneq,+\})$ such that the accessibility relations for $i$, $\triangleright$, $\triangleleft$, and $\subsetneq$ are empty in $\mathcal{G}'_j$. Then it is easy to see that each $\mathcal{G}'_j$ is still a full possibility frame, and the polymodal logics of the $\mathcal{G}'_j$ are still pairwise distinct. Finally, consider the continuum-sized set $\{\mathcal{F}\biguplus\mathcal{G}_j'\}_{j\in J}$ where $\mathcal{F}\biguplus\mathcal{G}_j'$ is the disjoint union of $\mathcal{F}$ and $\mathcal{G}_j'$ as in Definition \ref{DisUn}. We claim that (i) the logic of each $\mathcal{F}\biguplus\mathcal{G}_j'$ is Kripke-frame incomplete and (ii) the logics of the $\mathcal{F}\biguplus\mathcal{G}_j'$ are pairwise distinct. For (i), for our chosen $\varphi$ and $\psi$ above, we showed that (\ref{Splitting}) is valid over $\mathcal{F}$, and (\ref{Splitting}) is valid over $\mathcal{G}'_j$ since the accessibility relation for $i$ is empty in $\mathcal{G}'_j$, so (\ref{Splitting}) is valid over $\mathcal{F}\biguplus\mathcal{G}_j'$ by Proposition \ref{DisPres}. However, $\neg\Diamond_i \psi$ is not valid over $\mathcal{F}$, so by Proposition \ref{DisPres}, it is not valid over $\mathcal{F}\biguplus\mathcal{G}_j'$. Then since we observed at the beginning of this section that any Kripke frame whose logic includes (\ref{Splitting}) also includes $\neg\Diamond_i \psi$, we have established (i). For (ii), by the initial description of $\{\mathcal{G}_j\}_{j\in J}$, for any distinct $\mathcal{G}_j$ and $\mathcal{G}_k$, there is a $\chi\in \mathcal{L}(\sig,\{+\})$ that is valid over one but not the other. Suppose that $\mathcal{G}_j\Vdash\chi$ but $\mathcal{G}_k\nVdash\chi$. Then since $\mathcal{G}_j'$ and $\mathcal{G}_k'$ are obtained from $\mathcal{G}_j$ and $\mathcal{G}_k$ by adding empty accessibility relations for $\triangleright$ and $\triangleleft$ (and $i$ and $\subsetneq$), it follows that $\mathcal{G}_j'\Vdash (\Box_\triangleright\bot\wedge\Box_\triangleleft\bot)\rightarrow\chi$ but $\mathcal{G}_k'\nVdash (\Box_\triangleright\bot\wedge\Box_\triangleleft\bot)\rightarrow \chi$. By our construction of $\mathcal{F}$, every regular open set has an $R_\triangleright$-successor and every set of pairs has an $R_\triangleleft$-successor, so $\mathcal{F}\Vdash \neg (\Box_\triangleright\bot\wedge\Box_\triangleleft\bot)$ and hence $\mathcal{F}\Vdash(\Box_\triangleright\bot\wedge\Box_\triangleleft\bot)\rightarrow\chi$. Combining $\mathcal{G}_j'\Vdash (\Box_\triangleright\bot\wedge\Box_\triangleleft\bot)\rightarrow\chi$,  $\mathcal{G}_k'\nVdash (\Box_\triangleright\bot\wedge\Box_\triangleleft\bot)\rightarrow \chi$, and $\mathcal{F}\Vdash(\Box_\triangleright\bot\wedge\Box_\triangleleft\bot)\rightarrow\chi$, it follows by Proposition \ref{DisPres} that $\mathcal{F}\biguplus\mathcal{G}_j'\Vdash (\Box_\triangleright\bot\wedge\Box_\triangleleft\bot)\rightarrow \chi$ but $\mathcal{F}\biguplus\mathcal{G}_k'\nVdash (\Box_\triangleright\bot\wedge\Box_\triangleleft\bot)\rightarrow\chi$. Thus, we have established (ii). 

Now we prove the claim about Kripke-frame inconsistent logics. For each $j\in J$, let $(\mathcal{F}\biguplus\mathcal{G}_j')'$ be obtained from $\mathcal{F}\biguplus\mathcal{G}_j'$ by making the accessibility relation for $i$ the universal relation in $(\mathcal{F}\biguplus\mathcal{G}_j')'$ and changing nothing else. First, it is easy to check that  $(\mathcal{F}\biguplus\mathcal{G}_j')'$ is still a full possibility frame. Second, from the fact that  (\ref{Splitting}) is valid over $\mathcal{F}\biguplus\mathcal{G}_j'$, it is easy to check that (\ref{Splitting}) is valid over $(\mathcal{F}\biguplus\mathcal{G}_j')'$. This uses the fact that no state from $\mathcal{G}_j'$ has an $R_\triangleright$-successor, so our formula $\psi$, i.e., $\Diamond_\triangleright\top$, that appears in the antecedent of (\ref{Splitting}) cannot be true at a state from $\mathcal{G}_j'$. Not only does $(\mathcal{F}\biguplus\mathcal{G}_j')'$ validate (\ref{Splitting}), but unlike $\mathcal{F}\biguplus\mathcal{G}_j'$, it also clearly validates $\Diamond_i\psi$, i.e., $\Diamond_i \Diamond_\triangleright\top$. Thus, by the reasoning at the beginning of this section, the logic of $(\mathcal{F}\biguplus\mathcal{G}_j')'$ is Kripke-frame \textit{inconsistent}. It only remains to show that for $j\not=k$, the logics of $(\mathcal{F}\biguplus\mathcal{G}_j')'$ and $(\mathcal{F}\biguplus\mathcal{G}_k')'$ are distinct. This follows from the fact above that $\mathcal{F}\biguplus\mathcal{G}_j'\Vdash (\Box_\triangleright\bot\wedge\Box_\triangleleft\bot)\rightarrow \chi$ and $\mathcal{F}\biguplus\mathcal{G}_k'\nVdash (\Box_\triangleright\bot\wedge\Box_\triangleleft\bot)\rightarrow\chi$, which implies $(\mathcal{F}\biguplus\mathcal{G}_j')'\Vdash (\Box_\triangleright\bot\wedge\Box_\triangleleft\bot)\rightarrow \chi$ and $(\mathcal{F}\biguplus\mathcal{G}_k')'\nVdash (\Box_\triangleright\bot\wedge\Box_\triangleleft\bot)\rightarrow\chi$ because $\chi\in \mathcal{L}(\sig,\{+\})$ does not contain the $\Box_i$ modality, and $(\mathcal{F}\biguplus\mathcal{G}_j')'$ and $(\mathcal{F}\biguplus\mathcal{G}_k')'$ differ from $\mathcal{F}\biguplus\mathcal{G}_j'$ and $\mathcal{F}\biguplus\mathcal{G}_k'$, respectively, only in the accessibility relation for $\Box_i$. Thus, $\{(\mathcal{F}\biguplus\mathcal{G}_j')'\}_{j\in J}$ is our desired continuum-sized set of frames. 
\end{proof}

In \S~\ref{CompletenessTheory}, we will use Theorem~\ref{NoEquiv}, our duality theory for full possibility frames in \S~\ref{DualityTheory}, and known results about polymodal-to-unimodal reduction to prove the following result for the unimodal case.

\begin{restatable}[Unimodal Full Possibility Frames with No Kripke Equivalents]{theorem}{UNCOUNT}\label{uncount} There are continuum many full possibility frames for the unimodal language whose logics are pairwise distinct and Kripke-frame incomplete.
\end{restatable}
 
The syntactic form of (\ref{Splitting}) is a direct way of getting at the distinction between worlds and possibilities. It remains to be seen in what indirect ways the differences between Kripke frames and full possibility frames may show up syntactically in the basic polymodal language. It also remains to be seen what other kinds of mathematical structures can produce full possibility frames with no equivalent Kripke frame.

In what follows, we will be interested not only in full possibility frames, but possibility frames generally. \S~\ref{SpecialClasses} contains a catalogue of some of the most important other classes of possibility frames. 

With our introduction to possibility semantics now complete, we proceed in \S~\ref{morphisms} to our first major topic in the model theory of modal logic based on possibilities: morphisms between possibility frames.

\section{Possibility Morphisms}\label{morphisms}

Morphisms between frames and models are of fundamental importance in possible world semantics, as well as the more general setting of possibility semantics. Recall that given \textit{world} frames $\mathfrak{F}=\langle \wo{W},\{\wo{R}_i\}_{i\in\ind},\wo{A}\rangle$ and $\mathfrak{F}'=\langle \wo{W}',\{\wo{R}_i'\}_{i\in\ind},\wo{A}'\rangle$ (see Appendix \S~\ref{GFS}), a \textit{p-morphism} or \textit{bounded morphism} \citep[p.~309]{Blackburn2001} from $\mathfrak{F}$ to $\mathfrak{F}'$ is a function $f\colon \wo{W}\to\wo{W}'$ such that for all $w,v\in \wo{W}$ and $v'\in\wo{W}'$:
\begin{itemize}
\item[(a)] if $w\wo{R}_i v$, then $f(w)\wo{R}_i'f(v)$; 
\item[(b)] if $f(w)\wo{R}_i'v'$, then $\exists v\in \wo{W}$: $w\wo{R}_i v$ and $f(v)=v'$;
\item[(c)] $\forall X\in\wo{A}'$: $f^{-1}[X]\in \wo{A}$.
\end{itemize}
Note that (a) is equivalent to $f[\mathrm{R}_i(w)]\subseteq\mathrm{R}_i'(f(w))$ and (b) is equivalent to $f[\mathrm{R}_i(w)]\supseteq\mathrm{R}_i'(f(w))$. Also note that if $\mathfrak{F}$ and $\mathfrak{F}'$ are \textit{full} world frames, so $\wo{A}=\wp(\wo{W})$ and $\wo{A}'=\wp(\wo{W}')$, then (c) is trivially satisfied. 

Given world \textit{models} $\mathfrak{M}=\langle \wo{W},\{\wo{R}_i\}_{i\in\ind}, \wo{V}\rangle$ and $\mathfrak{M}'=\langle \wo{W}',\{\wo{R}_i'\}_{i\in\ind}, \wo{V}'\rangle$ for $\mathcal{L}(\sig,\ind)$, a p-morphism from $\mathfrak{M}$ to $\mathfrak{M}'$ is a function $f\colon \wo{W}\to\wo{W}'$ satisfying (a)-(b) above such that for all $p\in\sig$:
\begin{itemize}
\item[(d)] $\mathrm{V}(p)= f^{-1}[\mathrm{V}'(p)]$.
\end{itemize}

The key semantic preservation facts about p-morphisms are the following. First, if $f$ is a p-morphism from $\mathfrak{M}$ to $\mathfrak{M}'$, then for all $w\in \wo{W}$ and $\varphi\in\mathcal{L}(\sig,\ind)$, $\mathfrak{M},w\vDash \varphi$ iff $\mathfrak{M}',f(w)\vDash \varphi$. Second, if $f$ is a \textit{surjective} p-morphism from $\mathfrak{F}$ to $\mathfrak{F}'$, then $\mathfrak{F}\vDash \varphi$ implies $\mathfrak{F}'\vDash \varphi$. Third, say that $f$ is an \textit{embedding} iff $f$ is an \textit{injective} p-morphism such that for all $X\in \wo{A}$, there is an $X'\in\wo{A}'$ such that $f[X]=f[\wo{W}]\cap X'$ \citep[p.~309]{Blackburn2001}; then if there is an embedding of $\mathfrak{F}$ into $\mathfrak{F}'$, we have that $\mathfrak{F}'\vDash \varphi$ implies $\mathfrak{F}\vDash \varphi$.

Below we present analogous definitions and preservation facts suitably generalized for possibility semantics. We will define three different grades of \textit{possibility morphisms}, which requires some explanation.

\begin{remark}[Three Grades of Possibility Morphisms] While in possible world semantics there is one central notion of p-morphism, for possibility semantics we will define three different grades of morphisms: possibility morphisms, strict possibility morphisms, and p-morphisms. They are related as follows:
\begin{enumerate}
\item the concept of a possibility morphism\footnote{Cf. Goldblatt's \citeyearpar{Goldblatt2006b} concept of a modal map between general frames.} is the most general concept, singled out by the fact that if $h$ is a function from a possibility frame $\mathcal{F}$ to a possibility frame $\mathcal{G}$, then $h^{-1}[\cdot]$ is a BAO-homomorphism from $\mathcal{G}^\under$ to $\mathcal{F}^\under$ iff $h$ is a possibility morphism (Theorem \ref{PossMorphBAOMorph}.\ref{MtoM1});
\item over the class of all possibility frames, every p-morphism is a strict possibility morphism, \textbf{but not vice versa}, and every strict possibility morphism is a possibility morphism, \textbf{but not vice versa};
\item over the class of \textit{full} possibility frames satisfying additional conditions of strength (Definition \ref{StrongPoss}) and separativity (Definition \ref{SepFrames}) that can be assumed without loss of generality, \textbf{every possibility morphism is a strict possibility morphism} (Proposition \ref{SepStrongStrict}), but not every strict possibility morphism is a p-morphism;
\item over the classes of \textit{rich} possibility frames and \textit{filter-descriptive} possibility frames that will be important in our duality theory, \textbf{every possibility morphism is a p-morphism} (Propositions \ref{PossToP1} and \ref{PossToP2});
\item over the class of world frames regarded as possibility frames, i.e., in which $\sqsubseteq$ is discrete, every strict possibility morphism is a p-morphism. \hfill $\triangleleft$
\end{enumerate}
\end{remark}
 
\begin{definition}[Possibility Morphisms]\label{PossMorph} Given possibility frames $\mathcal{F}=\langle S,\sqsubseteq,\{R_i\}_{i\in \ind},\adm\rangle$ and $\mathcal{F}'={\langle S',\sqsubseteq',\{R_i'\}_{i\in \ind},\adm'\rangle}$,  a \textit{possibility morphism} from $\mathcal{F}$ to $\mathcal{F'}$ is a function $h\colon S\to S'$  such that for all $ x \in S$:
\begin{enumerate}[label=\arabic*.,ref=\arabic*]
\item\label{Sqmatch} \SqMatch{} -- $\forall X'\in \adm'$: $\mathord{\downarrow}'h(x)\cap X'=\emptyset$ iff $\mathord{\downarrow}x\cap h^{-1}[X']=\emptyset$;
\item\label{Rmatch} \RMatch{} -- $\forall X'\in \adm'$ $\forall i\in\ind$: $R_i'(h(x))\subseteq X'$ iff $R_i(x)\subseteq h^{-1}[X']$;
\item\label{PossMorph5} \PullBack{} -- $\forall X'\in \adm'$: $h^{-1}[X']\in \adm$.\footnote{If $\mathcal{F}$ and $\mathcal{F}'$ are full possibility frames, then \PullBack{} says that $h\colon S\to S'$ is such that the inverse image of each regular open subset of $\mathcal{F}'$ is a regular open subset of $\mathcal{F}$. A map between topological spaces such that the inverse image of each regular open set is regular open is called an \textit{R-map} \citep{Carnahan1973}.}
\end{enumerate}
A \textit{strict} possibility morphism from $\mathcal{F}$ to $\mathcal{F}'$ is a function $h\colon S\to S'$ satisfying \PullBack{} such that the following conditions hold for all $x,y\in S$ and $y',z'\in S'$, and the latter two hold for all $i\in \ind$: 
\begin{enumerate}[label=\arabic*.,ref=\arabic*,resume]
\item\label{SqForthCon} \SqForth{} -- if $y\sqsubseteq x$, then $h(y)\sqsubseteq' h(x)$;
\item \SqBack{} -- if $y'\sqsubseteq' h( x )$, then $\exists y$: $y\sqsubseteq x $ and $h(y)\sqsubseteq' y'$ (see Figure \ref{SqBackFig});
\item \RForth{} -- if $xR_i y$, then $h(x)R_i' h(y)$;
\item\label{SRBackCon} \SRBack{} -- if $h(x)R_i' y'$ and $z'\sqsubseteq ' y'$, then $\exists y$: $xR_i y$ and $h(y)\comp' z'$ (see Figure \ref{RBackFig}).
\end{enumerate}
As stated in Facts \ref{TautCond}-\ref{PullFact} below, these \textit{forth} and \textit{back} conditions jointly imply \SqMatch{} and \RMatch{} above, and they imply \PullBack{} whenever $\mathcal{F}$ is a full possibility frame.   

A \textit{p-morphism} is defined in the same way as a strict possibility morphism, but with strengthened versions of the two back conditions:
\begin{enumerate}[label=\arabic*.,ref=\arabic*,resume]
\item \pSqBack{} -- if $y'\sqsubseteq' h( x )$, then $\exists y$: $y\sqsubseteq x $ and $h(y)= y'$;
\item \pRBack{} -- if $h(x)R_i' y'$, then $\exists y$: $xR_i y$ and $h(y)=y'$.
\end{enumerate}
We also highlight the following special classes of possibility morphisms: 
\begin{enumerate}[label=\arabic*.,ref=\arabic*,resume]
\item\label{DenseMorph} a possibility morphism $h$ is \textit{dense}  iff $\forall x'\in S'$ $\exists x\in S$: $h(x)\sqsubseteq' x'$;
\item\label{pre-embed} a possibility morphism $h$ is \textit{robust}  iff $\forall X\in \adm$: $X=h^{-1}[h[X]]$ and $\exists X'\in \adm'$: $h[X]=h[S]\cap X'$;\footnote{See Definition 104 of \citealt{Goranko2007} for the  notion of a bounded \textit{strong} morphism between general frames, which requires that $\forall X\in P$: $X=h^{-1}[h[X]]$ and $h[X]\in P'$. Note that if $h$ is robust and surjective, then it is strong in this sense.} 
\item\label{embed} a possibility morphism $h$ is a \textit{strong embedding} iff for all $x,y\in S$:  (i) $y\sqsubseteq x$ iff $h(y)\sqsubseteq' h(x)$;\footnote{Condition (i) implies that $h$ is injective given the antisymmetry of $\sqsubseteq$.} (ii) for all $i\in \ind$, $xR_i y$ iff $h(x)R_i' h(y)$; and (iii) $\forall X\in \adm$ $\exists X'\in \adm'$: $h[X]=h[S]\cap X'$; 
\item a possibility morphism $h$ is a \textit{$\sqsubseteq$-strong embedding} iff it satisfies (i) and (iii) of part \ref{embed}; 
\item $h\colon S\to S'$ is an \textit{isomorphism} iff it is a bijection satisfying (i) and (ii) of part \ref{embed}, \PullBack{}, and $\forall X\in \adm$: $h[X]\in \adm'$.\hfill $\triangleleft$
\end{enumerate}
Given possibility models $\mathcal{M}=\langle S,\sqsubseteq,\{R_i\}_{i\in\ind},\pi\rangle$ and $\mathcal{M}'=\langle S',\sqsubseteq',\{R_i'\}_{i\in\ind},\pi'\rangle$ for $\mathcal{L}(\sig,\ind)$, in the weakest sense a possibility morphism from $\mathcal{M}$ to $\mathcal{M}'$ is an $h\colon S\to S'$ satisfying clauses \ref{Sqmatch}-\ref{Rmatch} above with $\adm'$ replaced by $\{\llbracket \varphi\rrbracket^{\mathcal{M}'}\mid \varphi\in\mathcal{L}(\sig,\ind)\}$ such that for all $p\in\sig$:
\begin{enumerate}[label=\arabic*.,ref=\arabic*,resume]
\item\label{PossMorphAtomic} $\pi(p)= h^{-1}[\pi'(p)]$.
\end{enumerate}
Strict possibility morphisms/p-morphisms between models are defined using the same \textit{forth} and \textit{back} conditions as above.\hfill $\triangleleft$
\end{definition}

 \begin{figure}[h]
\begin{center}
\begin{tikzpicture}[->,>=stealth',shorten >=1pt,shorten <=1pt, auto,node
distance=2cm,thick,every loop/.style={<-,shorten <=1pt}]
\tikzstyle{every state}=[fill=gray!20,draw=none,text=black]
\node (x) at (0,0) {{$x$}};
\node (hx) at (3,0) {{$h(x)$}};
\node (y') at (3,-1.5) {{$y'$}};
\draw (0,-1.5) ellipse (1.25cm and 2.25cm);
\draw (3,-1.5) ellipse (1.25cm and 2.25cm);

\node (F) at (0,1.25) {{$\mathcal{F}$}};
\node (F) at (3,1.25) {{$\mathcal{F}'$}};

\node at (5,-1.5) {{\textit{$\Rightarrow$}}};

\path (x) edge[dotted,->] node {{}} (hx);
\path (hx) edge[->] node {{}} (y');

\node (xR) at (7,0) {{$x$}};
\node (hxR) at (10,0) {{$h(x)$}};
\node (y'R) at (10,-1.5) {{$y'$}};
\node (yR) at (7,-1.5) {{$y$}};
\node (hyR) at (10,-3) {{$h(y)$}};

\node (exists) at (6.75,-.65) {{$\exists$}};

\draw (7,-1.5) ellipse (1.25cm and 2.25cm);
\draw (10,-1.5) ellipse (1.25cm and 2.25cm);

\node (F) at (7,1.25) {{$\mathcal{F}$}};
\node (F) at (10,1.25) {{$\mathcal{F}'$}};

\path (xR) edge[dotted,->] node {{}} (hxR);
\path (yR) edge[dotted,->] node {{}} (hyR);
\path (xR) edge[->] node {{}} (yR);
\path (y'R) edge[->] node {{}} (hyR);
\path (hxR) edge[->] node {{}} (y'R);

\end{tikzpicture}
\end{center}
\caption{the \SqBack{} condition of strict possibility morphisms. Dotted lines indicate the possibility morphism $h$, while a solid line from $s$ down to $t$ indicates that $t$ is a refinement of $s$.}\label{SqBackFig}
\end{figure}
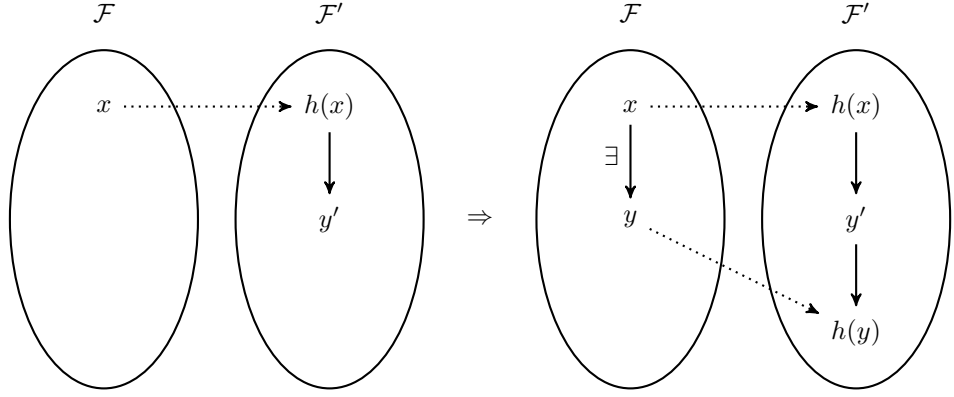 

 \begin{figure}[h]
\begin{center}
\begin{tikzpicture}[->,>=stealth',shorten >=1pt,shorten <=1pt, auto,node
distance=2cm,thick,every loop/.style={<-,shorten <=1pt}]
\tikzstyle{every state}=[fill=gray!20,draw=none,text=black]

\node (x) at (.3,.75) {{$x$}};

\node (hx) at (4,.75) {{$h(x)$}};
\node (Bhx) at (3.75,.75) {{}};
\node (Shx) at (4.25,.75) {{}};

\node (y') at (5.3,.75) {{$y'$}};
\node (Sy') at (5.2,.75) {{}};

\node (z') at (5.3,-.3) {{$z'$}};

\draw (1,0) ellipse (1.5cm and 1.5cm);
\draw (4.5,0) ellipse (1.5cm and 1.5cm);

\node (F) at (1,1.9) {{$\mathcal{F}$}};
\node (F') at (4.5,1.9) {{$\mathcal{F}'$}};

\path (x) edge[dotted,->] node {{}} (Bhx);
\path (y') edge[->] node {{}} (z');
\path (Shx) edge[dashed,->] node {{}} (Sy');

\node at (6.75,0) {{\textit{$\Rightarrow$}}};

\node (xR) at (8.3,.75) {{$x$}};

\node (BhxR) at (11.75,.75) {{}};
\node (hxR) at (12,.75) {{$h(x)$}};
\node (ShxR) at (12.25,.75) {{}};

\node (y'R) at (13.3,.75) {{$y'$}};
\node (Sy'R) at (13.2,.75) {{}};

\draw (9,0) ellipse (1.5cm and 1.5cm);
\draw (12.5,0) ellipse (1.5cm and 1.5cm);

\node (FR) at (9,1.9) {{$\mathcal{F}$}};
\node (F'R) at (12.5,1.9) {{$\mathcal{F}'$}};

\path (xR) edge[dotted,->] node {{}} (BhxR);
\path (ShxR) edge[dashed,->] node {{}} (Sy'R);

\node (yR) at (9.7,-.3) {{$y$}};
\path (xR) edge[dashed,->] node {{}} (yR);

\node (hyR) at (12,-.3) {{$h(y)$}};
\node (BhyR) at (11.75,-.3) {{}};
\path (yR) edge[dotted,->] node {{}} (BhyR);

\node (z'R) at (13.3,-.3) {{$z'$}};

\node (cL) at (12.9,-1.1) {{}};
\node (cR) at (12.65,-1.125) {{}};
\path (hyR) edge[->] node {{}} (cL);
\path (z'R) edge[->] node {{}} (cR);

\path (y'R) edge[->] node {{}} (z'R);

\node (exists) at (8.8,.05) {{$\exists$}};

\end{tikzpicture}
\end{center}
\caption{the \SRBack{} condition of strict possibility morphisms. Dashed lines indicate the accessibility relations $R_i$ and $R_i'$. Note that if $R_i'$ satisfies \Rdown{} from \S~\ref{FullFrames}, then \SRBack{} is equivalent to the following simpler condition: if $h(x)R_i'y'$, then $\exists y$: $xR_iy$ and $h(y)\comp' y'$.}\label{RBackFig}
\end{figure}
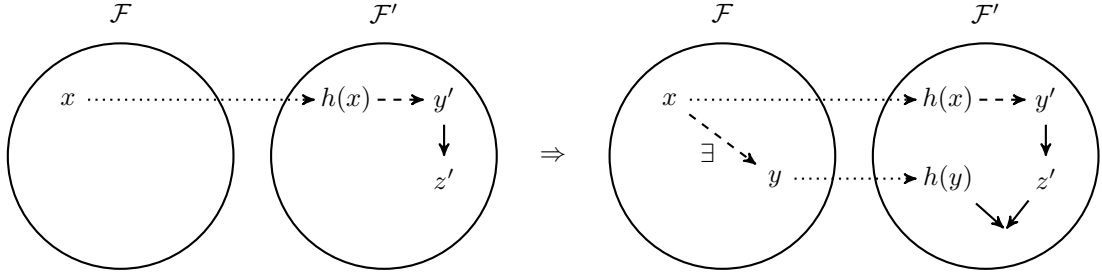

Let us make a number of clarificatory comments on Definition \ref{PossMorph}.

First, it can sometimes be useful to think of \SqMatch{} and \RMatch{} with the right-hand sides of the `iff' written as $h[\mathord{\downarrow}x]\cap X'=\emptyset$  and $h[R_i(x)]\subseteq X'$, respectively.

Second, note that if $h$ is a possibility morphism from a frame $\mathcal{F}=\langle S,\sqsubseteq, \{R_i\}_{i\in\ind},\adm\rangle$ to a frame $\mathcal{F}'=\langle S',\sqsubseteq', \{R_i'\}_{i\in\ind},\adm'\rangle$, and we consider admissible models $\mathcal{M}=\langle S,\sqsubseteq, \{R_i\}_{i\in\ind},\pi\rangle$ and $\mathcal{M}'={\langle S',\sqsubseteq', \{R_i'\}_{i\in\ind},\pi'\rangle}$ such that for all $p\in\sig$, $\pi(p)=h^{-1}[\pi'(p)]$, then by Definition \ref{PossMorph}, $h$ is a possibility morphism from $\mathcal{M}$ to $\mathcal{M}'$, since $\{\llbracket \varphi\rrbracket^{\mathcal{M}'}\mid \varphi\in\mathcal{L}(\sig,\ind)\}\subseteq\adm'$ by Fact \ref{TruthSub}.\ref{TruthSub2}.

Third, the difference between \SRBack{} and \pRBack{}  can be understood as follows. Where $\langle S',\sqsubseteq'\rangle$ is the poset of the target frame $\mathcal{F}'$ and $X\subseteq S'$, recall the definitions of $\mathord{\Downarrow}X$, $\mathrm{cl}(X)$, and $\mathrm{int}(X)$ from Remark \ref{Persp2} and following. Also recall from Fact \ref{RefReg}.\ref{RefReg2.5} that for any $X\subseteq S'$, $\mathrm{int}(\mathrm{cl}(\mathord{\Downarrow}X))$ is the smallest regular open set in the downset topology on $\langle S',\sqsubseteq'\rangle$ that includes $X$. Now consider the following back conditions:
\begin{itemize}
\item[B1] $R_i'(h(x))\subseteq h[R_i(x)]$;
\item[B2] $\mathrm{cl}(\mathord{\Downarrow} R_i'(h(x)))\subseteq \mathrm{cl}(\mathord{\Downarrow} h[R_i(x)])$;
\item[B3] $\mathrm{int}(\mathrm{cl}(\mathord{\Downarrow} R_i'(h(x))))\subseteq \mathrm{int}(\mathrm{cl}(\mathord{\Downarrow} h[R_i(x)]))$.
\end{itemize}
B1-B3 clarify the earlier back conditions for $R_i$ as follows.

\begin{fact}[Back Conditions for $R_i$]\label{BCC} For any possibility frames $\mathcal{F}$ and $\mathcal{F}'$ and function $h\colon \mathcal{F}\to\mathcal{F}'$:
\begin{enumerate}
\item\label{BCC1} B1 implies B2, which implies B3;
\item\label{BCC2} \pRBack{} is equivalent to B1;
\item\label{BCC4} \SRBack{} is equivalent to $\mathord{\Downarrow} R_i'(h(x))\subseteq \mathrm{cl}(\mathord{\Downarrow} h[R_i(x)])$, which is equivalent to B2;
\item\label{BCC5} B3 implies the right-to-left direction of \RMatch{};
\item\label{BCC6} if $\mathcal{F}'$ is a \textit{full} possibility frame, then B3 is equivalent to the right-to-left direction of \RMatch{}.
\end{enumerate}
\end{fact}
\begin{proof} Part \ref{BCC1} is immediate from the monotonicity property of the operators $\mathord{\Downarrow}$, $\mathrm{cl}$, and $\mathrm{int}$.  Part \ref{BCC2} and the first equivalence of part \ref{BCC4} are easy to verify by inspection of the conditions. For the second equivalence of part \ref{BCC4}, from left to right, if $\mathord{\Downarrow} R_i'(h(x))\subseteq \mathrm{cl}(\mathord{\Downarrow} h[R_i(x)])$, then $\mathrm{cl}(\mathord{\Downarrow} R_i'(h(x)))\subseteq \mathrm{cl}(\mathrm{cl}(\mathord{\Downarrow} h[R_i(x)]))=\mathrm{cl}(\mathord{\Downarrow} h[R_i(x)])$, and the right-to-left direction follows from $\mathord{\Downarrow} R_i'(h(x))\subseteq \mathrm{cl}(\mathord{\Downarrow} R_i'(h(x)))$.

For part \ref{BCC5}, to establish \RMatch{}, suppose that for an $X'\in \adm'$, $R_i(x)\subseteq h^{-1}[X']$, i.e., $h[R_i(x)]\subseteq X'$. We must show that $R_i'(h(x))\subseteq X'$. From $h[R_i(x)]\subseteq X'$, we have $\mathrm{int}(\mathrm{cl}(\mathord{\Downarrow}h[R_i(x)]))\subseteq \mathrm{int}(\mathrm{cl}(\mathord{\Downarrow}X'))$, and from $X'\in\adm'\subseteq \mathrm{RO}(S',\sqsubseteq')$, we have $X'=\mathrm{int}(\mathrm{cl}(\mathord{\Downarrow}X'))$, so $\mathrm{int}(\mathrm{cl}(\mathord{\Downarrow}h[R_i(x)]))\subseteq X'$. Then since $R_i'(h(x))\subseteq \mathrm{int}(\mathrm{cl}(\mathord{\Downarrow} R_i'(h(x))))$, we have $R_i'(h(x))\subseteq X'$ by B3.

For part \ref{BCC6}, we assume that for any $X'\in \adm'$, $h[R_i(x)]\subseteq X'$ implies $R_i'(h(x))\subseteq X'$. By Fact \ref{RefReg}.\ref{RefReg2.5}, $h[R_i(x)]\subseteq\mathrm{int}(\mathrm{cl}(\mathord{\Downarrow} h[R_i(x)]))\in\mathrm{RO}(S',\sqsubseteq')$, so since $\mathcal{F}'$ is \textit{full}, we have $\mathrm{int}(\mathrm{cl}(\mathord{\Downarrow} h[R_i(x)]))\in\adm'$, so we can use our previous assumption with $X'=\mathrm{int}(\mathrm{cl}(\mathord{\Downarrow} h[R_i(x)]))$ to conclude that $R_i'(h(x))\subseteq \mathrm{int}(\mathrm{cl}(\mathord{\Downarrow} h[R_i(x)]))$, which implies that $\mathrm{int}(\mathrm{cl}(\mathord{\Downarrow} R_i'(h(x))))\subseteq\mathrm{int}(\mathrm{cl}(\mathord{\Downarrow} h[R_i(x)]))$, which is B3.
\end{proof}

Thus, \textit{p-morphisms} use B1 and \textit{strict} possibility morphisms use B2. We could also give a special name to morphisms satisfying B3, but we will not need to in this paper. 

An analysis similar to that of the back conditions for $R_i$ above applies to back conditions for $\sqsubseteq$. Note that our \SqBack{} is equivalent to $\mathord{\downarrow}'h(x)\subseteq\mathrm{cl}(h[\mathord{\downarrow}x])$, which is equivalent to $\mathrm{cl}(\mathord{\downarrow}'h(x))\subseteq\mathrm{cl}(h[\mathord{\downarrow}x])$.

Observe how the other conditions on strict possibility morphisms imply the other \textit{matching} conditions on possibility morphisms, with Fact \ref{TautCond}.\ref{Taut4} following from Fact \ref{BCC}.

\begin{fact}[Strict Conditions \& Matching Conditions]\label{TautCond} For any possibility frames $\mathcal{F}$ and $\mathcal{F}'$ and $h\colon \mathcal{F}\to\mathcal{F}'$:
\begin{enumerate}
\item\label{Taut1} \SqForth{} implies the left-to-right direction of \SqMatch{};
\item\label{Taut2} \SqBack{} implies the right-to-left direction of \SqMatch{};
\item\label{Taut3} \RForth{} implies the left-to-right direction of \RMatch{};
\item\label{Taut4} \SRBack{} implies the right-to-left direction of \RMatch{}.
\end{enumerate}
\end{fact}

As in the case of p-morphisms mentioned at the beginning of this section, so too in the case of strict possibility morphisms, the \PullBack{} condition comes for free for \textit{full} frames, as in Fact \ref{PullFact}.

\begin{fact}[\PullBack{} to Full Frames]\label{PullFact} If $\mathcal{F}$ is a full possibility frame, $\mathcal{F}'$ is any possibility frame, and $h\colon \mathcal{F}\to\mathcal{F}'$ satisfies \SqForth{} and \SqBack{}, then $h$ satisfies \PullBack{}.
\end{fact}

\begin{proof} Where $\mathcal{F}=\langle S,\sqsubseteq,\{R_i\}_{i\in \ind},\adm\rangle$ and $\mathcal{F}'=\langle S',\sqsubseteq',\{R_i'\}_{i\in\ind},\adm'\rangle$, suppose $X'\in \adm'$, so $X'$ satisfies \textit{persistence} and \textit{refinability} with respect to $\langle S',\sqsubseteq'\rangle$. Then we will show that $h^{-1}[X']$ satisfies \textit{persistence} and \textit{refinability} with respect to $\langle S,\sqsubseteq\rangle$, so $h^{-1}[X']\in \mathrm{RO}(S,\sqsubseteq)=\adm$ by the assumption that $\mathcal{F}$ is full. 

For \textit{persistence}, suppose $x\in h^{-1}[X']$, so $h(x)\in X'$, and $y\sqsubseteq x$. Then by \SqForth{}, $h(y)\sqsubseteq' h(x)$, so $h(x)\in X'$ implies $h(y)\in X'$ by \textit{persistence} for $X'$, so $y\in h^{-1}[X']$ as desired. For \textit{refinability}, suppose $x\not\in h^{-1}[X']$, so $h(x)\not\in X'$. Then by \textit{refinability} for $X'$, there is a $y'\sqsubseteq' h(x)$ such that (i) for all $z'\sqsubseteq' y'$, $z'\not\in X'$. By \SqBack{}, $y'\sqsubseteq' h(x)$ implies there is a $y\sqsubseteq x$ such that $h(y)\sqsubseteq' y'$. Now we claim that for all $z\sqsubseteq y$, $z\not\in h^{-1}[X']$. For if $z\sqsubseteq y$, then $h(z)\sqsubseteq' h(y)$ by \SqForth{}, which with $h(y)\sqsubseteq' y'$ and (i) implies $h(z)\not\in X'$, so $z\not\in h^{-1}[X']$. Thus, $x\not\in h^{-1}[X']$ implies $\exists y\sqsubseteq x$ $\forall z\sqsubseteq y$: $z\not\in h^{-1}[X']$ as desired.\end{proof}

As a final clarification on Definition \ref{PossMorph}, let us note several aspects of our definition of  \textit{strong embeddings}.

\begin{remark}[Strong Embeddings] $\,$ 
\begin{enumerate} 
\item All strong embeddings are \textit{robust} possibility morphisms, as injectivity gives $h^{-1}[h[X]]=X$;
\item \textit{surjective} strong embeddings are equivalent to isomorphisms;
\item without surjectivity, strong embeddings are not guaranteed to be strict possibility morphisms, so we may speak of \textit{strict} strong embeddings. In \S~\ref{OpFrameAlg}, we consider the kind of \textit{subframes} that are images of strict strong embeddings. \hfill $\triangleleft$    
\end{enumerate}
\end{remark}

Now the following result demonstrates the importance of Definition \ref{PossMorph}.
 
\begin{proposition}[Preservation by Possibility Morphisms]\label{Preservation} For any possibility models $\mathcal{M}$ and $\mathcal{M}'$ and possibility frames $\mathcal{F}$ and $\mathcal{F}'$:
\begin{enumerate}
\item\label{Preservation1} if there is a possibility morphism $h$ from $\mathcal{M}$ to $\mathcal{M}'$, then for all $ x \in\mathcal{M}$ and $\varphi\in\mathcal{L}(\sig,\ind)$, $\mathcal{M}, x \Vdash\varphi$ iff $\mathcal{M}',h( x )\Vdash \varphi$;
\item\label{Preservation3} if there is a \textit{dense} possibility morphism from $\mathcal{F}$ to $\mathcal{F}'$, then for all $\varphi\in\mathcal{L}(\sig,\ind)$, $\mathcal{F}\Vdash\varphi$ implies $\mathcal{F}'\Vdash\varphi$;
\item\label{Preservation2} if there is a \textit{robust} possibility morphism from $\mathcal{F}$ to $\mathcal{F}'$, then for all $\varphi\in\mathcal{L}(\sig,\ind)$, $\mathcal{F}'\Vdash\varphi$ implies $\mathcal{F}\Vdash\varphi$.
\end{enumerate}
\end{proposition}

\begin{proof} Part \ref{Preservation1} is by induction on $\varphi$. The atomic case is by Definition \ref{PossMorph}.\ref{PossMorphAtomic}; the $\wedge$ case is routine; and since $\varphi\rightarrow\psi$ is equivalent to $\neg (\varphi\wedge\neg\psi)$ over possibility frames (see proof of Lemma \ref{classicality}), we do not need a $\rightarrow$~case.
 
For the $\neg$ and $\Box_i$ cases, the inductive hypothesis gives us $h^{-1}[\llbracket \varphi\rrbracket^{\mathcal{M}'}]=\llbracket \varphi\rrbracket^\mathcal{M}$. Thus, by \SqMatch{} with $\adm'$ replaced by $\{\llbracket \varphi\rrbracket^{\mathcal{M}'}\mid \varphi\in\mathcal{L}(\sig,\ind)\}$, $\mathord{\downarrow}'h(x)\cap \llbracket \varphi\rrbracket^{\mathcal{M}'}=\emptyset$ iff $\mathord{\downarrow}x\cap \llbracket \varphi\rrbracket^\mathcal{M}=\emptyset$, so $\mathcal{M}',h(x)\Vdash \neg\varphi$ iff $\mathcal{M},x\Vdash\neg\varphi$. Similarly, by \RMatch{}, $R_i'(h(x))\subseteq \llbracket \varphi\rrbracket^{\mathcal{M}'}$ iff $R_i(x)\subseteq \llbracket \varphi\rrbracket^\mathcal{M}$, so $\mathcal{M}',h(x)\Vdash\Box_i\varphi$ iff $\mathcal{M},x\Vdash \Box_i\varphi$.

For part \ref{Preservation3}, if $\mathcal{F}'\nVdash\varphi$, then there is a possibility model $\mathcal{M}'=\langle \mathcal{F}',\pi'\rangle$ and $y'\in\mathcal{M}'$ such that $\mathcal{M}',y'\nVdash\varphi$, in which case Refinability implies that there is an $x'\sqsubseteq' y'$ such that $\mathcal{M}',x'\Vdash\neg\varphi$. Given our morphism $h$ from $\mathcal{F}$ to $\mathcal{F}'$, define a valuation $\pi$ on $\mathcal{F}$ by $\pi(p)=h^{-1}[\pi'(p)]$. Then $h^{-1}[\pi'(p)]\in \adm$ by \PullBack{}, so $\mathcal{M}=\langle \mathcal{F},\pi\rangle$ is an admissible model based on $\mathcal{F}$, and $h$ is a possibility morphism from $\mathcal{M}$ to $\mathcal{M}'$ according to Definition \ref{PossMorph}.\ref{PossMorphAtomic}. Finally, since $h$ is a \textit{dense} possibility morphism from $\mathcal{F}$ to $\mathcal{F}'$, there is an $x\in \mathcal{M}$ such that $h(x)\sqsubseteq' x'$, which with $\mathcal{M}',x'\Vdash\neg\varphi$ implies $\mathcal{M}',h(x)\nVdash\varphi$, which with part \ref{Preservation1} implies $\mathcal{M},x\nVdash\varphi$, so $\mathcal{F}\nVdash\varphi$. 

For part \ref{Preservation2}, suppose $\mathcal{F}\nVdash\varphi$, so there is a possibility model $\mathcal{M}=\langle \mathcal{F},\pi\rangle$ and $x\in\mathcal{M}$ such that $\mathcal{M},x\nVdash \varphi$. Since our morphism $h$ from $\mathcal{F}$ to $\mathcal{F}'$ is \textit{robust}, we can choose for each $p\in \sig$ a $\pi'(p)\in P'$ such that $h[\pi(p)]=h[S]\cap \pi'(p)$. Let this define  a valuation $\pi'$ on $\mathcal{F}'$, so $\mathcal{M}'=\langle\mathcal{F}',\pi'\rangle$ is an admissible model based on $\mathcal{F}'$. Then for all $y\in S$, from the equation $h[\pi(p)]=h[S]\cap \pi'(p)$ we have that $y\in\pi(p)$ implies $h(y)\in\pi'(p)$; and from the same equation we have that $h(y)\in \pi'(p)$ implies $h(y)\in h[\pi(p)]$, which with $\pi(p)=h^{-1}[h[\pi(p)]]$ from the robustness of $h$ implies $y\in \pi(p)$. Thus, $h$ is a possibility morphism from $\mathcal{M}$ to $\mathcal{M}'$ according to Definition \ref{PossMorph}.\ref{PossMorphAtomic}, in which case from $\mathcal{M},x\nVdash\varphi$ and part \ref{Preservation1} we have $\mathcal{M}',h(x)\nVdash \varphi$, so $\mathcal{F}'\nVdash\varphi$.\end{proof} 

We will use Proposition \ref{Preservation} in \S~\ref{SpecialClasses} and \S~\ref{DualityTheory} to show that various frame constructions preserve validity and non-validity. We have already seen one frame construction that preserves validity and non-validity, namely the construction of $\mathcal{F}^\tight$ from $\mathcal{F}$ in Proposition \ref{Representation}. The proof of Proposition \ref{Representation}.\ref{Representation4} gives us the following.

\begin{fact}[$\mathcal{F}^\tight$ Construction]\label{TightCon} For any possibility frame $\mathcal{F}=\langle S,\sqsubseteq,\{R_i\}_{i\in\ind},\adm\rangle$, the identity map on $S$ is a surjective robust possibility morphism from $\mathcal{F}$ to the frame $\mathcal{F}^\tight=\langle S,\sqsubseteq, \{R_i^\tight\}_{i\in\ind},\adm\rangle$ in Proposition \ref{Representation}.
\end{fact}
\noindent A surjective robust morphism is also a dense and robust morphism, so by Proposition \ref{Preservation}, $\mathcal{F}\Vdash \varphi$ iff $\mathcal{F}^\tight\Vdash \varphi$.

An important fact about our morphisms is that they compose to form morphisms of the same type.

\begin{fact}[Composition of Morphisms]\label{Composition} For any possibility frames $\mathcal{F}$, $\mathcal{G}$, and $\mathcal{H}$ and functions $f\colon \mathcal{F}\to\mathcal{G}$ and $g\colon \mathcal{G}\to\mathcal{H}$:
\begin{enumerate}
\item\label{Composition1} if $f$ and $g$ are possibility morphisms, then $g\circ f$ is a possibility morphism;
\item\label{Composition2} if $f$ and $g$ are strict possibility morphisms, then $g\circ f$ is a strict possibility morphism;
\item\label{Composition3} if $f$ and $g$ are p-morphisms, then $g\circ f$ is a p-morphism.
\end{enumerate}
\end{fact}

\begin{proof}  First, the \PullBack{} condition of possibility morphisms is preserved by composition: if $X^\mathcal{H}\in\adm^\mathcal{H}$, then
  $g^{-1}[X^\mathcal{H}]\in\adm^\mathcal{G}$ by \PullBack{} for $g$, so $(g\circ f)^{-1}[X^\mathcal{H}]=f^{-1}
  [g^{-1}[X^\mathcal{H}]]\in\adm^\mathcal{F}$ by \PullBack{} for $f$.

For part \ref{Composition1}, the proofs that $g\circ f$ satisfies \SqMatch{} and \RMatch{} follow exactly the same pattern, so we include only the latter. We must show that (a) $\forall x^\mathcal{F}\in\mathcal{F}$ $\forall X^\mathcal{H}\in \adm^\mathcal{H}$: $R_i^\mathcal{H}(g(f(x^\mathcal{F})))\subseteq X^\mathcal{H}$ iff $R_i^\mathcal{F}(x^\mathcal{F})\subseteq (g\circ f)^{-1}[X^\mathcal{H}]$. By \textit{$R$-matching} for $g$, we have that $\forall x^\mathcal{G}\in\mathcal{G}$ $\forall X^\mathcal{H}\in \adm^\mathcal{H}$: $R_i^\mathcal{H}(g(x^\mathcal{G}))\subseteq X^\mathcal{H}$ iff $R_i^\mathcal{G}(x^\mathcal{G})\subseteq g^{-1}[X^\mathcal{H}]$. Thus, we have (b) $\forall x^\mathcal{F}\in\mathcal{F}$ $\forall X^\mathcal{H}\in \adm^\mathcal{H}$: $R_i^\mathcal{H}(g(f(x^\mathcal{F})))\subseteq X^\mathcal{H}$ iff $R_i^\mathcal{G}(f(x^\mathcal{F}))\subseteq g^{-1}[X^\mathcal{H}]$. By \textit{pull back} for $g$, $X^\mathcal{H}\in\adm^\mathcal{H}$ implies $g^{-1}[X^\mathcal{H}]\in\adm^\mathcal{G}$, so by \textit{$R$-matching} for $f$, we have (c) $R_i^\mathcal{G}(f(x^\mathcal{F}))\subseteq g^{-1}[X^\mathcal{H}]$ iff $R_i^\mathcal{F}(x^\mathcal{F})\subseteq f^{-1}[g^{-1}[X^\mathcal{H}]]=(g\circ f)^{-1}[X^\mathcal{H}]$. Together (b) and (c) imply (a). 

For part \ref{Composition2}, that \SqForth{} (resp. \RForth{}) for $f$ and $g$ implies \SqForth{} (resp. \RForth{}) for $g\circ f$ is obvious. To show that $g\circ f$ satisfies \SqBack{}, we must show that if $y^\mathcal{H}\sqsubseteq^\mathcal{H} g(f(x))$, then there is a $y^\mathcal{F}\in\mathcal{F}$ such that $y^\mathcal{F}\sqsubseteq^\mathcal{F} x$ and $g(f(y^\mathcal{F}))\sqsubseteq^\mathcal{H} y^\mathcal{H}$. So suppose $y^\mathcal{H}\sqsubseteq^\mathcal{H} g(f(x))$. Then by \SqBack{} for $g$, there is a $y^\mathcal{G}\in\mathcal{G}$ such that $y^\mathcal{G}\sqsubseteq^\mathcal{G}f(x)$ and  $g(y^\mathcal{G})\sqsubseteq^\mathcal{H} y^\mathcal{H}$. Given $y^\mathcal{G}\sqsubseteq^\mathcal{G}f(x)$ and  \SqBack{} for $f$, there is a $y^\mathcal{F}\in\mathcal{F}$ such that $y^\mathcal{F}\sqsubseteq^\mathcal{F} x$ and $f(y^\mathcal{F})\sqsubseteq^\mathcal{G} y^\mathcal{G}$. By \SqForth{} for $g$,  $f(y^\mathcal{F})\sqsubseteq^\mathcal{G} y^\mathcal{G}$ implies $g(f(y^\mathcal{F}))\sqsubseteq^\mathcal{H} g(y^\mathcal{G})$, which with $g(y^\mathcal{G})\sqsubseteq^\mathcal{H} y^\mathcal{H}$ implies $g(f(y^\mathcal{F}))\sqsubseteq^\mathcal{H} y^\mathcal{H}$, which with $y^\mathcal{F}\sqsubseteq^\mathcal{F} x$ means that $y^\mathcal{F}$ is our desired witness. 

Next, to show that $g\circ f$ satisfies \SRBack{}, we must show that if $g(f(x))R_i^\mathcal{H} y^\mathcal{H}$ and $z^\mathcal{H}\sqsubseteq^\mathcal{H} y^\mathcal{H}$, then there is a $y^\mathcal{F}\in\mathcal{F}$ such that $xR_i^\mathcal{F} y^\mathcal{F}$ and $g(f(y^\mathcal{F}))\comp^\mathcal{H} z^\mathcal{H}$. So suppose $g(f(x))R_i^\mathcal{H} y^\mathcal{H}$ and $z^\mathcal{H}\sqsubseteq^\mathcal{H} y^\mathcal{H}$. Then by \SRBack{} for $g$, there is $y^\mathcal{G}\in\mathcal{G}$ such that $f(x)R_i^\mathcal{G} y^\mathcal{G}$ and  $g(y^\mathcal{G})\comp^\mathcal{H}z^\mathcal{H}$. Since $g(y^\mathcal{G})\comp^\mathcal{H}z^\mathcal{H}$, there is a $u^\mathcal{H}\in\mathcal{H}$ such that $u^\mathcal{H}\sqsubseteq^\mathcal{H}g(y^\mathcal{G})$ and $u^\mathcal{H}\sqsubseteq^\mathcal{H}z^\mathcal{H}$. By \SqBack{} for $g$, $u^\mathcal{H}\sqsubseteq^\mathcal{H}g(y^\mathcal{G})$ implies that there is a $u^\mathcal{G}\in\mathcal{G}$ such that $u^\mathcal{G}\sqsubseteq^\mathcal{G} y^\mathcal{G}$ and $g(u^\mathcal{G})\sqsubseteq^\mathcal{H} u^\mathcal{H}$. Given $f(x)R_i^\mathcal{G} y^\mathcal{G}$ and $u^\mathcal{G}\sqsubseteq^\mathcal{G} y^\mathcal{G}$, \SRBack{} for $f$ implies that there is a $y^\mathcal{F}\in\mathcal{F}$ such that $xR_i^\mathcal{F}y^\mathcal{F}$ and $f(y^\mathcal{F})\comp^\mathcal{G} u^\mathcal{G}$. Thus, there is a $v^\mathcal{G}\in\mathcal{G}$ such that $v^\mathcal{G}\sqsubseteq^\mathcal{G} f(y^\mathcal{F})$ and $v^\mathcal{G}\sqsubseteq^\mathcal{G} u^\mathcal{G}$, which with \SqForth{} for $g$ implies $g(v^\mathcal{G})\sqsubseteq^\mathcal{H} g(f(y^\mathcal{F}))$ and $g(v^\mathcal{G})\sqsubseteq^\mathcal{H} g(u^\mathcal{G})$, which with $g(u^\mathcal{G})\sqsubseteq^\mathcal{H} u^\mathcal{H}\sqsubseteq^\mathcal{H} z^\mathcal{H}$ from above gives us $g(f(y^\mathcal{F}))\comp^\mathcal{H} z^\mathcal{H}$. This completes the proof of part \ref{Composition2}.

Part \ref{Composition3} is well known.\end{proof} 

The importance of Fact \ref{Composition} is that it allows us to think in categorical terms as follows.

\begin{remark}[Categories]\label{Categories} By Fact \ref{Composition}, any class $\mathsf{F}$ of possibility frames together with all possibility morphisms (resp.~strict possibility morphisms, p-morphisms) between frames in $\mathsf{F}$ constitutes a \textit{category}, where the objects are the frames, the morphisms are the possibility morphisms, the identity morphism for each frame is the identity function, and composition of morphisms is functional composition. This is the categorical perspective we will adopt for the duality theory of \S~\ref{DualityTheory}. All of the concepts from category theory that we will use can be found in, e.g., \S\S~3-4 of \citealt{Adamek2009}. \end{remark}

\section{Special Classes of Frames}\label{SpecialClasses} 

In this section, we survey classes of possibility frames that are important for understanding the relations between possibility frames and world frames and between possibility frames and Boolean algebras with operators. Here is a brief statement of the importance of each class of frames to be considered:
 
\begin{itemize}
\item \textbf{separative} frames (\S~\ref{SepSec}) -- these are important because \textit{separativity} simplifies reasoning about frames, without loss of generality, and is related to the \textit{tight} frames of \S~\ref{TightSection} and the \textit{principal} frames of \S~\ref{PrincFrames}.
\item \textbf{atomic} frames (\S~\ref{AtomicSection}) -- these are important because any atomic (full) possibility frame can be easily transformed into a semantically equivalent (full) \textit{world} frame.
\item \textbf{extended} frames (\S~\ref{ExtFrames}) -- these frames have a distinguished minimum element $\bot$, which can be useful when working with functional frames as in \S~\ref{FuncFrames} or principal frames as in \S~\ref{PrincFrames}.
\item \textbf{functional} frames (\S~\ref{FuncFrames}) -- these frames support the functional semantics for $\Box_i$ mentioned in \S~\ref{intro} (cf.~\citealt{Holliday2014}), where $\mathcal{M},x\Vdash \Box_i \varphi$ iff $\mathcal{M},f_i(x)\Vdash \varphi$ (or $f_i$ is undefined at $x$), and they are important in the duality theory relating possibility frames and so-called $\mathcal{T}$-BAOs in \S~\ref{DualityTheory}.
\item \textbf{tight} frames (\S~\ref{TightSection}) -- the notion of tightness will be used in characterizing the \textit{rich} frames of \S~\ref{RichFrames} and the \textit{filter-descriptive} frames of \S~\ref{Fdes}, both of which are central to the duality theory of \S~\ref{DualityTheory}.
\item \textbf{principal} frames (\S~\ref{PrincFrames}) -- these will be important in the duality theory relating possibility frames and \textit{completely additive} Boolean algebras with operators ($\mathcal{V}$-BAOs) in \S~\ref{VtoPossSection}.
\item \textbf{rich} frames (\S~\ref{RichFrames}) -- this subclass of principal frames will be important in providing a categorical duality with \textit{complete} Boolean algebras with  \textit{completely additive} operators ($\mathcal{CV}$-BAOs) in \S~\ref{DualEquiv}.
\end{itemize} 
For most of these frame classes, the \textit{powerset possibilization} of a Kripke frame (Example \ref{PowerPoss}) will provide an example of a frame in the class, but we wish to abstract away from some properties of such powerset possibilizations---especially their atomicity in light of \S~\ref{AtomicSection}. A diagram illustrating the relations between some of the above frame classes will appear in Figure \ref{FrameClasses} at the end of our tour in \S~\ref{RichFrames}.

Recall that we have already encountered some special classes of possibility frames in \S~\ref{FullFrames}, namely \textit{strong} possibility frames (Definition \ref{StrongPoss}) and \textit{standard} possibility frames (Definition \ref{Standard}). These will reappear in \S\S~\ref{FuncFrames}-\ref{RichFrames} and when we study correspondence theory for possibility semantics in \S~\ref{LemmScottCorr}.
 
\subsection{Separative Frames}\label{SepSec}

The following relation on states behaves much like the refinement relation $\sqsubseteq$ in possibility frames.

\begin{definition}[$\cof$]\label{CoRef} Given a partial-state frame $\mathcal{F}=\langle S,\sqsubseteq ,\{R_i\}_{i\in\ind},\adm\rangle$ and $x,y\in S$, define 
\[x\cof y\mbox{ iff }\forall x'\sqsubseteq x,\, x'\comp y,\]
i.e., $x\cof y$ iff $\forall x'\sqsubseteq x\, \exists x''\sqsubseteq x'$: $x''\sqsubseteq y$. Let $x\simeq_\cofsub y$ iff $x\cof y$ and $y\cof x$.\hfill $\triangleleft$ 
\end{definition}  

If $\langle S,\sqsubseteq\rangle$ is such that every state is refined by a \textit{minimal} point (see Definition \ref{AtDef}), then $x\cof y$ iff every minimal point that refines $x$ also refines $y$.  Note that $\cof$ is a preorder,\footnote{For transitivity, suppose $x_1\cof x_2$ and $x_2\cof x_3$. Toward showing $x_1\cof x_3$, suppose $ z_1 \sqsubseteq x_1$. Then given $x_1\cof x_2$, there is some $ z'_1 \sqsubseteq  z_1 $ such that $ z'_1 \sqsubseteq x_2$. Then given $x_2\cof x_3$, there is some $ z'_2\sqsubseteq  z'_1 $ such that $ z'_2\sqsubseteq x_3$. Given $ z'_2\sqsubseteq  z'_1 \sqsubseteq  z_1 $, by the transitivity of $\sqsubseteq$ we have $ z'_2\sqsubseteq   z_1 $. Thus, for any $ z_1 \sqsubseteq x_1$ there is a $ z'_2\sqsubseteq   z_1 $ such that $ z'_2\sqsubseteq x_3$, which implies $x_1\cof x_3$.} but not necessarily antisymmetric. 

We first observe that sets of possibilities satisfying \textit{persistence} and \textit{refinability} are closed under $\cof$.

\begin{fact}[\textit{$\cof$-persistence}]\label{CofClose} Given a poset $\langle S,\sqsubseteq\rangle$, if $X\in\mathrm{RO}(S,\sqsubseteq)$ (recall Notation \ref{ROnotation}), then $X$ satisfies \textit{$\cof$-persistence}: if $x\in X$ and $x'\cof x$, then $x'\in X$.
\end{fact}

\begin{proof} Suppose that $x'\cof x$. If $x'\not\in X$, then by \textit{refinability} there is a $y'\sqsubseteq  x'$ such that (i) for all $y''\sqsubseteq y'$, $y''\not\in X$. By Definition \ref{CoRef}, together $x'\cof x$ and $y'\sqsubseteq  x'$ imply that there is a $y''\sqsubseteq y'$ such that $y''\sqsubseteq x$. By (i), $y''\sqsubseteq y'$ implies $y''\not\in X$, which with $y''\sqsubseteq x$ and \textit{persistence} implies $x\not\in X$.\end{proof}

Next, we observe that taking the $\cof$-successors of a given possibility is a way of generating a set of possibilities that satisfies \textit{persistence} and \textit{refinability}.

\begin{fact}[$\cof$-generated Propositions]\label{CofGenerated} Given a poset $\langle S,\sqsubseteq\rangle$ and $x\in S$, $\{x'\in S\mid x'\cof x\}\in\mathrm{RO}(S,\sqsubseteq)$.
\end{fact}

\begin{proof}
Since $x''\sqsubseteq x'\cof x$ implies $x''\cof x$, the set $\{x'\in S\mid x'\cof x\}$ satisfies \textit{persistence}. 

For \textit{refinability}, suppose $y\not\in \{x'\in S\mid x'\cof x\}$, so $y\not\cof x$, so there is a $y'\sqsubseteq y$ such that for all $y''\sqsubseteq y'$, $y''\not\sqsubseteq x$. It follows that for all $y''\sqsubseteq y'$ and $y'''\sqsubseteq y''$, $y'''\not\sqsubseteq x$, so $y''\not\cof x$ and hence $y''\not\in  \{x'\in S\mid x'\cof x\}$. Thus, $ \{x'\in S\mid x'\cof x\}$ satisfies \textit{refinability}.
\end{proof}

Finally, we observe that $\cof$ behaves like $\sqsubseteq$ with respect to the forcing relation in possibility models.

\begin{fact}[Forcing and $\cof$]\label{PrecPerRef}\textnormal{For any possibility model $\mathcal{M}=\langle S,\sqsubseteq ,\{R_i\}_{i\in\ind},\pi\rangle$,  $x, y\in S$, and $\varphi\in\mathcal{L}(\sig,\ind)$: 
\begin{enumerate}
\item\label{PrecPer} $\cof$-Persistence: if $\mathcal{M},x\Vdash\varphi$ and $ y\cof x$, then $\mathcal{M}, y\Vdash\varphi$;
\item\label{PrecRef} $\cof$-Refinability: if $\mathcal{M},x\nVdash \varphi$, then $\exists x'\cof x$: $\mathcal{M},x'\Vdash\neg\varphi$;
\item\label{PrecNeg} $\cof$-Negation: $\mathcal{M},x\Vdash\neg\varphi$ iff $\forall x'\cof x$, $\mathcal{M},x'\nVdash\varphi$;
\item\label{PrecDup} $\cof$-Duplication: if $x\simeq_\cofsub  y$, then $\mathcal{M},x\Vdash \varphi$ iff $\mathcal{M}, y\Vdash \varphi$.
\end{enumerate}} 
\end{fact}
\begin{proof} By Fact \ref{TruthSub} and Definition \ref{PosFrames}, the truth set of any formula in a possibility model satisfies \textit{persistence} and \textit{refinability}. Thus, part \ref{PrecPer} follows from  Fact \ref{CofClose}.

Part \ref{PrecRef} is immediate from \textit{refinability} and the fact that $x'\sqsubseteq x$ implies $x'\cof x$.

For part \ref{PrecNeg}, the right-to-left direction holds because $x'\sqsubseteq x$ implies $x'\cof x$. For the left-to-right direction, suppose there is a $ x'\cof  x $ with $\mathcal{M}, x'\Vdash\varphi$. Since $ x'\cof  x $, there is a $ x''\sqsubseteq x'$ with $ x''\sqsubseteq x $. By \textit{persistence}, $\mathcal{M}, x'\Vdash\varphi$ and $ x''\sqsubseteq x'$ together imply $\mathcal{M}, x''\Vdash\varphi$, which with $ x''\sqsubseteq x $ implies $\mathcal{M}, x \nVdash\neg\varphi$.

Part \ref{PrecDup} follows from part \ref{PrecPer}.
\end{proof}

In light of the similarities between $\cof$ and $\sqsubseteq$ observed above, it is natural to consider frames in which $\langle S,\sqsubseteq\rangle$ is, in the terminology of set-theoretic forcing (e.g., \citealt[p. 4]{Jech1986}), a \textit{separative poset}.

\begin{definition}[Separative Frames]\label{SepFrames} A poset $\langle S,\sqsubseteq\rangle$ is \textit{separative} iff $\sqsubseteq \,= \, \cof$, i.e., $x\sqsubseteq y$ iff $\forall x'\sqsubseteq x$ $\exists x''\sqsubseteq x'$: $x''\sqsubseteq y$. A  partial-state frame $\mathcal{F}=\langle S,\sqsubseteq,\{R_i\}_{i\in \ind},\adm\rangle$ is separative iff $\langle S,\sqsubseteq\rangle$ is separative. \hfill $\triangleleft$
\end{definition} 

\noindent Note that every world frame and powerset possibilization thereof (Example \ref{PowerPoss}) is separative. 

Separativity can be characterized in other useful ways.

\begin{fact}[Separativity]\label{Sep&Princ} For any poset $\langle S,\sqsubseteq\rangle$:
\begin{enumerate}
\item\label{Sep&Princ1} $\langle S,\sqsubseteq\rangle$ is separative iff for all $y\in S$, $\mathord{\downarrow}y$ satisfies \textit{refinability}, so $\mathord{\downarrow}y\in\mathrm{RO}(S,\sqsubseteq)$;
\item\label{Sep&Princ2} if $\langle S,\sqsubseteq\rangle$ is separative, then for any $x,y\in S$, if $x\not\sqsubseteq y$, then $x$ and $y$ are distinguishable by a set in $\mathrm{RO}(S,\sqsubseteq)$: $\mathord{\downarrow}y\in\mathrm{RO}(S,\sqsubseteq)$ and $y\in \mathord{\downarrow}y$ but $x\not\in \mathord{\downarrow}y$.\end{enumerate}
\end{fact}
\begin{proof} For part \ref{Sep&Princ1},  \textit{refinability} for $\mathord{\downarrow}y$ says that if $\forall x'\sqsubseteq x$ $\exists x''\sqsubseteq x'$: $x''\in\mathord{\downarrow}y$, then $x\in\mathord{\downarrow}y$. But this is just to say that if $\forall x'\sqsubseteq x$ $\exists x''\sqsubseteq x'$: $x''\sqsubseteq y$, then $x\sqsubseteq y$, which is the nontrivial direction of separativity.

Part \ref{Sep&Princ2} is immediate from part \ref{Sep&Princ1}.\end{proof}

Recall the notion of \textit{differentiation} of general frames \citep[Def. 5.6]{Blackburn2001}.

\begin{definition}[Differentiated Frames]\label{DiffFrames} A partial-state frame $\mathcal{F}=\langle S,\sqsubseteq,\{R_i\}_{i\in \ind},\adm\rangle$ is \textit{differentiated} iff for all $x,y\in S$: $x=y$ iff for all $Z\in \adm$, $x\in Z$ iff $y\in Z$.
\end{definition}

By Fact \ref{Sep&Princ}.\ref{Sep&Princ2}, we have the following.

\begin{fact}[Separativity and Differentiation] Every separative full possibility frame is differentiated.
\end{fact} 

Another useful fact concerns possibility morphisms to separative frames.

\begin{fact}[Separativity and Morphisms]\label{SepOrd} For any possibility frames $\mathcal{F}=\langle S,\sqsubseteq,\{R_i\}_{i\in\ind},\adm\rangle$ and \\ $\mathcal{F}'=\langle S',\sqsubseteq',\{R_i'\}_{i\in\ind},\adm'\rangle$ and possibility morphism $h\colon S\to S'$:
\begin{enumerate}
\item\label{SepOrd1} if for every $x\in S$, $\mathord{\downarrow}'h(x)\in\adm'$, then $h$ is such that for all $x,y\in S$, $y\cof x$ implies $h(y)\sqsubseteq' h(x)$;
\item\label{SepOrd2} if $\mathcal{F}'$ is a full separative frame, then $h$ is such that for all $x,y\in S$, $y\cof x$ implies $h(y)\sqsubseteq' h(x)$.
\end{enumerate}
\end{fact}
\begin{proof} For part \ref{SepOrd1}, we use the \PullBack{} property of possibility morphisms: for any $X'\in\adm'$, $h^{-1}[X']\in \adm$. Suppose that for every $x\in S$, $\mathord{\downarrow}'h(x)\in\adm'$, so $h^{-1}[\mathord{\downarrow}'h(x)]\in \adm$ by \PullBack{}, which means that $h^{-1}[\mathord{\downarrow}'h(x)]\in\mathrm{RO}(S,\sqsubseteq)$ since $\mathcal{F}$ is a possibility frame. Then by Fact \ref{CofClose}, together $x\in h^{-1}[\mathord{\downarrow}'h(x)]$ and $y\cof x$ imply $y\in h^{-1}[\mathord{\downarrow}'h(x)]$, which means $h(y)\in \mathord{\downarrow}'h(x)$, so $h(y)\sqsubseteq' h(x)$. Thus, $y\cof x$ implies $h(y)\sqsubseteq' h(x)$.

Part \ref{SepOrd2} follows from part \ref{SepOrd1} and Fact \ref{Sep&Princ}.\ref{Sep&Princ1}.
\end{proof}

A key fact motivating the definition of \textit{strict} possibility morphisms in Definition \ref{PossMorph} is that every possibility morphism to a full, separative, and strong (Definition \ref{StrongPoss}) possibility frame is strict.

\begin{proposition}[Separative, Strong, Strict]\label{SepStrongStrict} For any possibility frame $\mathcal{F}=\langle S,\sqsubseteq,\{R_i\}_{i\in \ind},\adm\rangle$ and full, separative, and strong possibility frame $\mathcal{F}'=\langle S',\sqsubseteq',\{R_i'\}_{i\in\ind},\adm'\rangle$, if $h\colon S\to S'$ is a possibility morphism, then $h$ is a strict possibility morphism.
\end{proposition}
\begin{proof} The \SqForth{}  property follows from the stronger property in Fact \ref{SepOrd}.\ref{SepOrd2}. 

For \SqBack{}, suppose $y'\sqsubseteq' h(x)$, so $y'\in \mathord{\downarrow}'h(x)$. Since $\mathcal{F}'$ is separative and full, $\mathord{\downarrow}'y'\in P'$ by  Fact \ref{Sep&Princ}.\ref{Sep&Princ1}. Then by the \SqMatch{} condition of possibility morphisms,  $\mathord{\downarrow}'h(x)\cap \mathord{\downarrow}' y'\neq\emptyset$ implies $\mathord{\downarrow}x\cap h^{-1}[\mathord{\downarrow}'y']\neq\emptyset$, so there is a $y\sqsubseteq x$ such that $h(y)\sqsubseteq' y'$, which establishes \SqBack{}.

For \RForth{}, since $\mathcal{F}'$ is strong and full, we have $R_i'(h(x))\in P'$ by Proposition \ref{MasterCon}. Then by the \PullBack{} property of possibility morphisms, $h^{-1}[R_i'(h(x))]\in P$, so by the left-to-right direction of the \RMatch{} property of possibility morphisms, $R_i'(h(x))\subseteq R_i'(h(x))$ implies $R_i(x)\subseteq h^{-1}[R_i'(h(x))]$, which is equivalent to the \RForth{} condition that $xR_iy$ implies $h(x)R_i'h(y)$.

Finally, for \SRBack{}, suppose $h(x)R_i'y'$ and $z'\sqsubseteq' y'$. Since $\mathcal{F}'$ is separative and full, $\mathord{\downarrow}'z'\in P'$ by Fact \ref{Sep&Princ}.\ref{Sep&Princ1} and hence $\neg \mathord{\downarrow}'z'\in P'$. From $h(x)R_i'y'$ and $z'\sqsubseteq' y'$, we have $R_i'(h(x))\not\subseteq \neg \mathord{\downarrow}'z'$, which by the right-to-left direction of \RMatch{} implies  $R_i(x)\not\subseteq h^{-1}[\neg \mathord{\downarrow}'z']$, so there is a $y\in S$ such that $xR_iy$ and $h(y)\not\in \neg \mathord{\downarrow}'z'$, which implies $h(y)\comp ' z'$, which establishes \SRBack{}.
\end{proof}

It will follow from Proposition \ref{UpDownRich} and Theorem \ref{AlmostBack} that for any full possibility frame $\mathcal{F}$, there is a separative full possibility frame $\mathcal{F}'$ and a dense and robust possibility morphism from $\mathcal{F}$ to $\mathcal{F}'$, so by Proposition \ref{Preservation}, $\mathcal{F}$ and $\mathcal{F}'$ validate the same formulas. But we can also prove the following stronger result more directly. 

\begin{restatable}[Separative Quotient]{proposition}{SepQuo}\label{sep-quo} For every possibility frame $\mathcal{F}$, there is a separative possibility frame $\mathcal{F}^\simeq$ (such that if $\mathcal{F}$ is full, so is $\mathcal{F}^\simeq$) and a surjective robust possibility morphism from $\mathcal{F}$ to $\mathcal{F}^\simeq$. Thus, by Proposition \ref{Preservation}, for all $\varphi\in\mathcal{L}(\sig,\ind)$,  $\mathcal{F}\Vdash\varphi$ iff $\mathcal{F}^\simeq\Vdash\varphi$.
\end{restatable}

\begin{proof} See Appendix \S~\ref{SepProof}. The claim also follows from Proposition \ref{Tightening} together with Fact \ref{TightSepDiff}.\ref{TightSepDiff0}.
\end{proof}

\begin{figure}[h]
\begin{center}
\begin{tikzpicture}[yscale=1, ->,>=stealth',shorten >=1pt,shorten <=1pt, auto,node
distance=2cm,thick,every loop/.style={<-,shorten <=1pt}]
\tikzstyle{every state}=[fill=gray!20,draw=none,text=black]

\node[circle,draw=black!100,fill=black!100,inner sep=0pt,minimum size=.175cm]  (x) at (0,0) {{}};
\node[circle,draw=black!100,fill=black!100,inner sep=0pt,minimum size=.175cm]  (y) at (-1,-1) {{}};
\node[circle,draw=black!100,fill=black!100,inner sep=0pt,minimum size=.175cm] (z) at (1,-1) {{}};
\node[circle,draw=black!100,fill=black!100,inner sep=0pt,minimum size=.175cm]  (y') at (-1,-2) {{}};
\node[circle,draw=black!100,fill=black!100,inner sep=0pt,minimum size=.175cm] (z') at (1,-2) {{}};

\path (x) edge[->] node {{}} (y);
\path (x) edge[->] node {{}} (z);
\path (y) edge[->] node {{}} (y');
\path (y) edge[->] node {{}} (z');
\path (z) edge[->] node {{}} (y');
\path (z) edge[->] node {{}} (z');

\node[circle,draw=black!100,fill=black!100,inner sep=0pt,minimum size=.175cm]  (a) at (4,-1) {{}};
\node[circle,draw=black!100,fill=black!100,inner sep=0pt,minimum size=.175cm]  (b) at (3,-2) {{}};
\node[circle,draw=black!100,fill=black!100,inner sep=0pt,minimum size=.175cm]  (c) at (5,-2) {{}};
\path (a) edge[->] node {{}} (b);
\path (a) edge[->] node {{}} (c);

\path (x) edge[dotted,->,bend left] node {{}} (a);
\path (y) edge[dotted,->,bend left] node {{}} (a);
\path (z) edge[dotted,->,bend left] node {{}} (a);

\path (y') edge[dotted,->,bend right] node {{}} (b);
\path (z') edge[dotted,->,bend right] node {{}} (c);

\end{tikzpicture}
\end{center}
\caption{A possibility frame (left) and its separative quotient (right). For both frames, assume that the accessibility relation $R_i$, not shown, is the universal relation.}\label{SepQuotFig}
\end{figure}
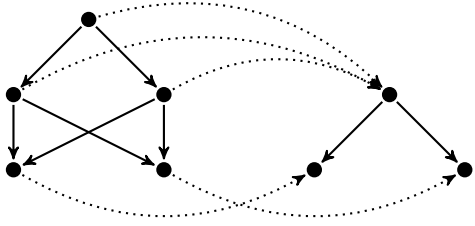

The separative quotient $\mathcal{F}^\simeq$ simply collapses possibilities that are equivalent according to the  $\cof$ relation, as in Figure \ref{SepQuotFig}.

Since the construction of \textit{strong} possibility frames from full possibility frames for Proposition \ref{Representation} preserves separativity, and the composition of two surjective robust possibility morphisms is also a surjective robust possibility morphism, from Propositions \ref{sep-quo} and \ref{Representation} and Fact \ref{TightCon} we have the following.

\begin{corollary}[Separative Strong Frames] For any full possibility frame $\mathcal{F}$, there is a separative, strong, and full possibility frame $(\mathcal{F}^\simeq)^\tight$ and a surjective robust possibility morphism from $\mathcal{F}$ to $(\mathcal{F}^\simeq)^\tight$.
\end{corollary} 

\subsection{Atomic Frames}\label{AtomicSection}

World frames and their powerset possibilizations (Examples \ref{KripkeExample} and \ref{PowerPoss}) are examples of what we will call \textit{atomic} possibility frames in Definition \ref{AtDef}, deviating slightly from the standard definition for posets.

\begin{definition}[Atomic Poset]\label{AtomicPoset} Given a poset $\langle S,\sqsubseteq\rangle$, an \textit{atom} in $\langle S,\sqsubseteq\rangle$ is an $a\in S$ that is not the minimum of $\langle S,\sqsubseteq\rangle$ (if there is one) such that for all $ x\in S$, if $x\sqsubseteq a$, then either $x=a$ or $x$ is the minimum of $\langle S,\sqsubseteq\rangle$. A poset is \textit{atomic} iff for every non-minimum element $x\in S$, there is an atom $a$ such that $a\sqsubseteq x$. \hfill $\triangleleft$
\end{definition}
\noindent This is not quite the notion we want for possibility frames (though it would be fine for the extended frames of \S~\ref{ExtFrames}), so we define atomic possibility frames a bit differently.

\begin{definition}[Minimal Points and Atomic Frames]\label{AtDef} A \textit{minimal point} in a poset $\langle S,\sqsubseteq\rangle$ is an $a\in S$ such that for all $ x\in S$, if $x\sqsubseteq a$, then $x=a$. Let $\mathsf{min}\langle S,\sqsubseteq\rangle$ be the set of minimal points in $\langle S,\sqsubseteq\rangle$. A partial-state frame $\mathcal{F}=\langle S, \sqsubseteq, \{R_i\}_{i\in\ind},\adm\rangle$ is \textit{atomic} iff for every $x\in S$, there is an $a\in \mathsf{min}\langle S,\sqsubseteq\rangle$ such that $a\sqsubseteq x$. \hfill $\triangleleft$
\end{definition} 

\noindent This notion of atomic is the flipped version of what is called the \textit{McKinsey condition} \citep[p.~82]{Chagrov1997}, the condition that for all $x\in S$, there is a $y\in \mathsf{max}\langle S,\sqsubseteq\rangle$ such that $x\sqsubseteq y$. 

Note that the condition of atomicity concerns only the poset $\langle S,\sqsubseteq\rangle$, not the set $\adm$ of admissible propositions. Following standard terminology \citep[Def.~5.65]{Blackburn2001}, we could say that an atomic possibility frame is \textit{discrete} iff for each of its minimal points $a$, $\{a\}\in \adm$.

Using the following construction (cf.~\citealt{Venema1998}), we will show in Proposition \ref{AtToWorld} that atomic possibility frames are semantically equivalent to world frames, and atomic full possibility frames are semantically equivalent to full world frames (Kripke frames).

\begin{definition}[Atom Structure]\label{AtSt} Given an atomic possibility frame $\mathcal{F}=\langle S, \sqsubseteq, \{R_i\}_{i\in\ind},\adm\rangle$ and possibility model $\mathcal{M}=\langle\mathcal{F},\pi\rangle$, define $\mathfrak{At}\mathcal{F}=\langle S',\sqsubseteq', \{R_i'\}_{i\in \ind}, \adm'\rangle$ and $\mathfrak{At}\mathcal{M}=\langle \mathfrak{At}\mathcal{F},\pi'\rangle$ by:
\begin{enumerate}
\item\label{AtSt1} $S'= \mathsf{min}\langle S,\sqsubseteq\rangle$ and $\sqsubseteq'$ is the identity relation on $S'$;
\item\label{AtSt2} for all $a,b\in S'$, $aR_i'b$ iff $\exists x\in S$: $aR_i x$ and $b\sqsubseteq x$;
\item\label{AtSt3} $\adm'=\{\mathsf{min}\langle S,\sqsubseteq\rangle \cap X\mid X\in \adm\}$;
\item\label{AtSt4} for all $a\in S'$, $a\in \pi'(p)$ iff $a\in \pi(p)$. \hfill $\triangleleft$
\end{enumerate}
\end{definition}

Note that if $\mathcal{F}$ satisfies \Rdown{} (\S~\ref{FullFrames}), then the definition of $R_i'$ is equivalent to: $aR_i'b$ iff $aR_ib$.

\begin{proposition}[From Atomic Possibility Frames to World Frames]\label{AtToWorld} For any atomic possibility frame $\mathcal{F}$:
\begin{enumerate}
\item\label{AtToWorld1} $\mathfrak{At}\mathcal{F}$ is a world frame as in Example \ref{KripkeExample};
\item\label{AtToWorld2} if $\mathcal{F}$ is a \textit{full} possibility frame, then $\mathfrak{At}\mathcal{F}$ is a \textit{full} world frame (Kripke frame) as in Example \ref{KripkeExample};
\item\label{AtToWorld3} the identity map on $\mathfrak{At}\mathcal{F}$ is a dense and robust possibility morphism (and a strong embedding if $\mathcal{F}$ satisfies \Rdown{}) from $\mathfrak{At}\mathcal{F}$ into $\mathcal{F}$. Thus, by Proposition \ref{Preservation}, for all $\varphi\in\mathcal{L}(\sig,\ind)$, $\mathfrak{At}\mathcal{F}\Vdash\varphi$ iff $\mathcal{F}\Vdash\varphi$.
\end{enumerate}
\end{proposition} 
 
\begin{proof} For part \ref{AtToWorld1}, the only fact to check is that $\adm'$ satisfies the closure conditions of a (general) world frame as in Definition \ref{GenSem}, i.e., that it is closed under intersection, complement, and $\blacksquare_i'$. Since $\adm$ is closed under intersection, $\adm'$ is clearly closed under intersection as well. For the other two closure conditions, suppose $X'\in \adm'$, so there is an $X\in \adm$ such that $X'=\mathsf{min}\langle S,\sqsubseteq\rangle\cap X$. Since $X\in \adm$, it follows by Definition \ref{PosetMod} and Remark \ref{Persp2} that $\mathrm{int}(S\setminus X)=\{x\in S\mid \mathord{\downarrow}x\cap X=\emptyset\}\in\adm$. Now we claim that $S'\setminus X'= \mathsf{min}\langle S,\sqsubseteq\rangle \cap \mathrm{int}(S\setminus X)$. For the right-to-left inclusion, if $a\in \mathsf{min}\langle S,\sqsubseteq\rangle \cap \mathrm{int}(S\setminus X)$, then from $a\in \mathsf{min}\langle S,\sqsubseteq\rangle$ we have $a\in S'$ and from $a\in \mathrm{int}(S\setminus X)$ we have $a\not\in X$, so $a\not\in X'$. From left to right, if $a\in S'\setminus X'$, then $a\in\mathsf{min}\langle S,\sqsubseteq\rangle$ but $a\not\in X$, which together imply $a\in \mathrm{int}(S\setminus X)$, so $a\in \mathsf{min}\langle S,\sqsubseteq\rangle\cap \mathrm{int}(S\setminus X)$. Thus, $S'\setminus X'= \mathsf{min}\langle S,\sqsubseteq\rangle \cap \mathrm{int}(S\setminus X)$, so $\adm'$ is closed under complement because $\adm$ is closed under taking the interior of the complement.

Finally, we claim that $\blacksquare_i' X'=\mathsf{min}\langle S,\sqsubseteq\rangle \cap \blacksquare_i X$. For the right-to-left inclusion, if $a\not\in\blacksquare_i' X'$, then there is a $b\in S'=\mathsf{min}\langle S,\sqsubseteq\rangle$ such that $aR_i'b$ and $b\not\in X'$. From $b\in\mathsf{min}\langle S,\sqsubseteq\rangle$ and $b\not\in X'=\mathsf{min}\langle S,\sqsubseteq\rangle\cap X$, it follows that $b\not\in X$. Since $aR_i'b$, there is an $x\in S$ such that $aR_ix$ and $b\sqsubseteq x$. By \textit{persistence} for $X$, together $b\not\in X$ and $b\sqsubseteq x$ imply $x\not\in X$, which with $aR_ix$ implies $a\not\in\blacksquare_i X$. For the left-to-right inclusion, suppose $x\not\in\mathsf{min}\langle S,\sqsubseteq\rangle\cap\blacksquare_i X$. If $x\not\in \mathsf{min}\langle S,\sqsubseteq\rangle$, then $x\not\in S'$, so $x\not\in\blacksquare_i'X'$. So suppose $x\in \mathsf{min}\langle S,\sqsubseteq\rangle$ but $x\not\in \blacksquare_i X$, so there is a $y$ such that $xR_i y$ and $y\not\in X$. Then by \textit{refinability} for $X$, there is a minimal point $b\sqsubseteq y$ such that $b\not\in X$ and hence $b\not\in X'$. From $xR_iy$ and $b\sqsubseteq y$, we have $xR_i' b$, which with $b\not\in X'$ implies $x\not\in\blacksquare_i'X'$. Thus, we have shown that $\blacksquare_i' X'=\mathsf{min}\langle S,\sqsubseteq\rangle \cap \blacksquare_i X$, so $\adm'$ is closed under $\blacksquare_i'$ because $\adm$ is closed under $\blacksquare_i$.

For part \ref{AtToWorld2}, if $\mathcal{F}$ is \textit{full}, then for any set $U$ of minimal points in $\mathcal{F}$, we have $\mathrm{int}(\mathrm{cl}(U))\in \adm$ and ${\mathsf{min}\langle S,\sqsubseteq\rangle } \cap \mathrm{int}(\mathrm{cl}(U))=U$, so $U\in \adm'$ by Definition \ref{AtSt}.\ref{AtSt3} and hence $\mathfrak{At}\mathcal{F}$ is full.

For part \ref{AtToWorld3}, since $\sqsubseteq'$ is the restriction of $\sqsubseteq$ to $S'$, the identity map $h$ on $S'$ is such that $a\sqsubseteq' b$ iff $h(a)\sqsubseteq h(b)$. Since $\sqsubseteq'$ is the identity relation and all elements in $h[S']$ are minimal, $h$ satisfies \SqBack{} and \SqMatch{}. Since $h$ is injective, the requirement $X = h^{-1}[h[X]]$ of \textit{robust} possibility morphisms holds. In addition, if $\mathcal{F}$ satisfies \Rdown{}, then the requirement of \textit{strong embeddings} that $aR_i'b$ iff $h(a)R_ih(b)$ also holds. It only remains to show: 
\begin{itemize}
\item \RMatch{} -- $\forall a\in S'\,\forall X\in \adm$: $R_i(h(a))\subseteq X$ iff $R_i'(a)\subseteq h^{-1}[X]$;
\item \PullBack{} -- $\forall X\in \adm$: $h^{-1}[X]\in\adm'$;
\item \textit{embedding} -- $\forall X'\in \adm'$ $\exists X\in \adm$: $h[X']=h[S']\cap X$;
\item \textit{dense} -- $\forall x\in S$ $\exists a\in S'$: $h(a)\sqsubseteq x$.
\end{itemize}
For \PullBack{}, for any $X\in \adm$, $h^{-1}[X]=\mathsf{min}\langle S,\sqsubseteq\rangle \cap X\in \adm '$. For \RMatch{}, since $h(a)=a$ and $h^{-1}[X]=\mathsf{min}\langle S,\sqsubseteq\rangle\cap X$, we must show that $R_i(a)\subseteq X$ iff $R_i'(a)\subseteq \mathsf{min}\langle S,\sqsubseteq\rangle\cap X$. This follows from the fact, established in the proof of part \ref{AtToWorld1}, that $\mathsf{min}\langle S,\sqsubseteq\rangle \cap \blacksquare_i X = \blacksquare_i' (\mathsf{min}\langle S,\sqsubseteq\rangle\cap X)$.

Finally, from the facts that $h$ is the identity map and $S'=\mathsf{min}\langle S,\sqsubseteq\rangle$, the \textit{embedding} condition is just that for all $X'\in \adm'$ there is an $X\in \adm$ such that $X'=\mathsf{min}\langle S,\sqsubseteq\rangle\cap X$, which is immediate from the definition of $\adm'$ in Definition \ref{AtSt}.\ref{AtSt3}. From the same facts, the \textit{dense} condition is just the condition that $\mathcal{F}$ is atomic.\end{proof}
 
Although atomic possibility frames and their atom structures are semantically equivalent by Proposition \ref{AtToWorld}.\ref{AtToWorld3}, many non-isomorphic atomic possibility frames can have the same atom structure, so we do lose information when going from an atomic possibility frame to its atom structure, in a way that we do not lose information when going from a world frame to its powerset possibilization. The following proposition records the relationship between atom structures and powerset possibilizations.

\begin{proposition}[Atom Structures and Powerset Possibilizations]\label{AtomProp} For any world frame $\mathfrak{F}$ (regarded as a possibility frame as in Examples \ref{KripkeExample} and \ref{KripkeAgain}) and atomic possibility frame $\mathcal{F}$:
\begin{enumerate}
\item\label{AtomProp1} $\mathfrak{At}(\mathfrak{F}^\pow)$ is isomorphic to $\mathfrak{F}$;
\item\label{AtomProp2} the function $h\colon \mathcal{F}\to (\mathfrak{At}\mathcal{F})^\pow$ defined by $h(x)=\{a\in\mathfrak{At}\mathcal{F}\mid a\sqsubseteq^\mathcal{F} x\}$ is a dense and robust possibility morphism from $\mathcal{F}$ to $(\mathfrak{At}\mathcal{F})^\pow$, so by Proposition \ref{Preservation}, for all $\varphi\in\mathcal{L}(\sig,\ind)$, $\mathcal{F}\Vdash\varphi$ iff $(\mathfrak{At}\mathcal{F})^\pow\Vdash\varphi$;
\item\label{AtomProp3} if $\mathcal{F}$ is separative, then the $h$ from part \ref{AtomProp2} is a $\sqsubseteq$-strong embedding of $\mathcal{F}$ into $(\mathfrak{At}\mathcal{F})^\pow$. 
\end{enumerate}
\end{proposition}
\begin{proof} Part \ref{AtomProp1} is easy to check using the definitions.

For part \ref{AtomProp2}, let $\mathcal{F}=\langle S, \sqsubseteq, \{R_i\}_{i\in\ind},\adm\rangle$, $\mathfrak{At}\mathcal{F}=\langle S', \sqsubseteq', \{R_i'\}_{i\in\ind},\adm'\rangle$, and $(\mathfrak{At}\mathcal{F})^\pow={\langle S^{\prime\pow}, \sqsubseteq^{\prime\pow}, \{R_i^{\prime\pow}\}_{i\in\ind},\adm^{\prime\pow}\rangle}$. First, we show that $h$ satisfies:
\begin{itemize}
\item \PullBack{} -- $\forall \mathcal{X}\in \adm^{\prime\pow}$: $h^{-1}[\mathcal{X}]\in \adm$;
\item \textit{robust} -- $\forall X\in P$: $X=h^{-1}[h[X]]$ and $\exists \mathcal{X}\in \adm^{\prime\pow}$: $h[X]=h[S]\cap \mathcal{X}$.
\end{itemize} 
For \PullBack{}, by definition of $(\mathfrak{At}\mathcal{F})^{\pow}$, $\adm^{\prime\pow}= \{\mathord{\downarrow} X\mid X\in \adm'\}$, where $\mathord{\downarrow} X = \{Y\in S^{\prime\pow}\mid Y\sqsubseteq^{\prime\pow} X\}=\{Y\subseteq \mathsf{min}\langle S,\sqsubseteq\rangle\mid \emptyset\not=Y\subseteq X\}$, and by definition of $\mathfrak{At}\mathcal{F}$, $\adm'=\{\mathsf{min}\langle S,\sqsubseteq\rangle\cap X\mid X\in \adm\}$. So for $\mathcal{X}\in \adm^{\prime\pow}$, there is an $X\in\adm$ such that (i) $\mathcal{X}=\mathord{\downarrow}(\mathsf{min}\langle S,\sqsubseteq\rangle\cap X)=\{Y\subseteq \mathsf{min}\langle S,\sqsubseteq\rangle\mid \emptyset\not=Y\subseteq X\}$. We claim that $h^{-1}[\mathcal{X}]=X$. For the left-to-right inclusion, if $x\in h^{-1}[\mathcal{X}]$, so $h(x)\in\mathcal{X}$, then by (i), $h(x)\subseteq X$. By \textit{refinability} for $X$, if $x\not\in X$, then there is a minimal point $a\sqsubseteq x$ such that $a\not\in X$, so $h(x)\not\subseteq X$, contradicting what we just showed. Thus, $x\in X$. For the right-to-left inclusion, by \textit{persistence} for $X$, if $x\in X$, then for every minimal point $a\sqsubseteq x$, $a\in X$, so $h(x)\subseteq X$, which implies $h(x)\in\mathcal{X}$ and hence $x\in h^{-1}[\mathcal{X}]$. Thus, $h^{-1}[\mathcal{X}]=X\in\adm$.

For \textit{robust}, to show $X=h^{-1}[h[X]]$, we show that if $h(x)\in h[X]$, then $x\in X$. If $h(x)\in h[X]$, then there is an $x'\in X$ such that $h(x)=h(x')$. From $h(x)=h(x')$ we have that $x\cof x'$, which with $x'\in X$ implies $x\in X$ by Fact \ref{CofClose}. Next, we must show that there is an $\mathcal{X}\in \adm^{\prime\pow}$ such that $h[X]=h[S]\cap\mathcal{X}$. Taking $\mathcal{X}=\{Y\subseteq \mathsf{min}\langle S,\sqsubseteq\rangle\mid \emptyset\not=Y\subseteq X\}$, we have $\mathcal{X}\in \adm^{\prime\pow}$ by the unpacking of definitions above. For any $x\in X$, $h(x)\subseteq X$ by \textit{persistence} for $X$, so $h(x)\in\mathcal{X}$. Thus, $h[X]\subseteq\mathcal{X}$. Now suppose $Y\in h[S]\cap\mathcal{X}$. Since $Y\in\mathcal{X}$, $Y\subseteq X$, and since $Y\in h[S]$, there is a $y\in S$ such that $h(y)=Y$. If $y\not\in X$, then by \textit{refinability} for $X$ there is a minimal point $a\sqsubseteq y$ such that $a\not\in X$, so $h(y)\not\subseteq X$, contradicting the fact from the previous sentence that $h(y)=Y\subseteq X$. Thus, $y\in X$, so $Y=h(y)\in h[X]$. This shows that $h[S]\cap\mathcal{X}\subseteq h[X]$.

Next, we show that $h$ satisfies:
\begin{itemize}
\item \SqForth{} -- if $y\sqsubseteq x$, then $h(y)\sqsubseteq^{\prime\pow} h(x)$;
\item \SqBack{} -- if $Y\sqsubseteq^{\prime\pow} h( x )$, then $\exists y$: $y\sqsubseteq x $ and $h(y)\sqsubseteq^{\prime\pow} Y$;
\item \RMatch{} -- $\forall \mathcal{X}\in \adm^{\prime\pow}$: $R^{\prime\pow}_i(h(x))\subseteq \mathcal{X}$ iff $R_i(x)\subseteq h^{-1}[\mathcal{X}]$;
\item \textit{dense} -- $\forall X\in S^{\prime\pow}$ $\exists x\in S$: $h(x)\sqsubseteq^{\prime\pow} X$.
\end{itemize}

For \SqForth{}, if $y\sqsubseteq x$, then $h(y)\subseteq h(x)$, so $h(y)\sqsubseteq^{\prime\pow} h(x)$ by the definition of $(\mathfrak{At}\mathcal{F})^{\pow}$. 

For \SqBack{}, suppose $Y\sqsubseteq^{\prime\pow} h( x )$, so $Y\subseteq h( x )$. Since every state in $(\mathfrak{At}\mathcal{F})^{\pow}$ is a nonempty set of minimal points from $\mathcal{F}$, there is a minimal point $y\in Y$.  Then $Y\subseteq h( x )$ implies $y\sqsubseteq x$, and $h(y)=\{y\}$, so $y\in Y$ implies $h(y)\subseteq Y$ and hence $h(y)\sqsubseteq^{\prime\pow} Y$. 

For \RMatch{}, suppose $R_i(x)\not\subseteq h^{-1}[\mathcal{X}]$, so $x\not\in\blacksquare_i h^{-1}[\mathcal{X}]$. By \PullBack{}, $h^{-1}[\mathcal{X}]\in \adm$, so $\blacksquare_i h^{-1}[\mathcal{X}]\in\adm$. Then by \textit{refinability} for $\blacksquare_i h^{-1}[\mathcal{X}]$ and the fact that $\mathcal{F}$ is atomic, $x\not\in\blacksquare_i h^{-1}[\mathcal{X}]$ implies that there is a minimal point $a\sqsubseteq x$ such that $a\not\in\blacksquare_i h^{-1}[\mathcal{X}]$, so there is a $y$ such that $aR_iy$ and $y\not\in h^{-1}[\mathcal{X}]$, so $h(y)\not\in\mathcal{X}$. Since $a\sqsubseteq x$, $a\in h(x)$. Now we claim that $h(x)R_i^{\prime\pow} h(y)$, which by the definition of $(\mathfrak{At}\mathcal{F})^{\pow}$ is equivalent to $h(y)\subseteq R_i'[h(x)]$, where $R_i'$ is the accessibility relation in $\mathfrak{At}\mathcal{F}$. To prove the claim, take a $b\in h(y)$, so $b\sqsubseteq y$. Then since $aR_iy$, it follows by the definition of $R_i'$ that $aR_i'b$, which with $a\in h(x)$ implies $b\in R_i'[h(x)]$. Hence $h(y)\subseteq R_i'[h(x)]$, so $h(x)R_i^{\prime\pow} h(y)$, which with $h(y)\not\in\mathcal{X}$ implies $R^{\prime\pow}_i(h(x))\not\subseteq \mathcal{X}$.

Conversely, suppose $R^{\prime\pow}_i(h(x))\not\subseteq \mathcal{X}$, so there is a $Y$ such that $h(x)R_i^{\prime\pow} Y$, i.e., $Y\subseteq R_i'[h(x)]$, and $Y\not\in\mathcal{X}$. Then by \textit{refinability} for $\mathcal{X}$ and the fact that $(\mathfrak{At}\mathcal{F})^{\pow}$ is atomic, $Y\not\in\mathcal{X}$ implies that there is a minimal point $B$ in $(\mathfrak{At}\mathcal{F})^{\pow}$ such that $B\sqsubseteq^{\prime\pow} Y$, i.e., $B\subseteq Y$, and $B\not\in\mathcal{X}$. That $B$ is a minimal point in $(\mathfrak{At}\mathcal{F})^{\pow}$ means that $B=\{b\}$ for a $b$ in $\mathfrak{At}\mathcal{F}$. Now given $Y\subseteq R_i'[h(x)]$, we have $b\in R_i'[h(x)]$, so there is a minimal point $a\in h(x)$, i.e., $a\sqsubseteq x$, such that $aR_i'b$. Since $h(b)=B$ and $B\not\in\mathcal{X}$, $b\not\in h^{-1}[\mathcal{X}]$, which with $aR_i'b$ implies $a\not\in\blacksquare_i h^{-1}[\mathcal{X}]$, which with $a\sqsubseteq x$ and \textit{persistence} for  $\blacksquare_i h^{-1}[\mathcal{X}]$ implies $x\not\in\blacksquare_i h^{-1}[\mathcal{X}]$, so $R_i(x)\not\subseteq h^{-1}[\mathcal{X}]$. 

Finally, to show that $h$ is \textit{dense}, since every $X\in S^{\prime\pow}$ is a nonempty set of minimal points from $\mathcal{F}$, simply take a minimal point $x\in X$, so $h(x)=\{x\}\subseteq X$, which means $h(x)\sqsubseteq^{\prime\pow}X$.

For part \ref{AtomProp3}, if $\mathcal{F}$ is atomic and separative, then as noted after Definition \ref{CoRef}, $x\sqsubseteq y$ iff every minimal point that refines $x$ also refines $y$, i.e., $h(x)\subseteq h(y)$, which is equivalent to $h(x)\sqsubseteq^{\prime\pow} h(y)$. This implies, together with the fact from above that $h$ is a robust possibility morphism, that $h$ is a $\sqsubseteq$-strong embedding.\end{proof} 

This proof provides an example of how our \SqBack{} clause may apply when the standard back clause for a p-morphism does not, i.e., $Y\sqsubseteq^{\prime\pow} h( x )$ does not imply there is a $y\sqsubseteq x$ such that $h(y)=Y$, as required by a p-morphism. For if all we assume is that $\mathcal{F}$ is atomic, there is no guarantee that for each set $Y$ of minimal points in $\mathcal{F}$, there is a $y$ in $\mathcal{F}$ such that the set of minimal points refining $y$, our $h(y)$ above, is exactly $Y$.

In \S~\ref{RichFrames}, we will identify the possibility frames $\mathcal{F}$ for which $(\mathfrak{At}\mathcal{F})^\pow$ is \textit{isomorphic} to $\mathcal{F}$.

\subsection{Extended Frames}\label{ExtFrames}

When we defined the powerset possibilization of a world frame (Example \ref{PowerPoss}), we chopped off the bottom element $\emptyset$ of the poset $\langle \wp(\wo{W}),\subseteq\rangle$. However, it would sometimes be convenient---especially when dealing with the \textit{functional} frames of \S~\ref{FuncFrames} or the \textit{principal} frames of \S~\ref{PrincFrames}---to allow in our frames a distinguished minimum element $\inc$, which may be thought of as the ``impossible state.'' 

To be clear: a possibility frame as in Definition \ref{PosFrames} is already allowed to have a minimum in its poset $\langle S,\sqsubseteq\rangle$, but such frames are somewhat uninteresting for the following reason.

\begin{fact}[Collapse]\label{collapse} If a possibility model $\mathcal{M}=\langle S, \sqsubseteq, \{R_i\}_{i\in\ind},\pi\rangle$ is such that $\langle S,\sqsubseteq\rangle$ has a minimum element, then for every formula $\varphi\in\mathcal{L}(\sig,\ind)$, either $\llbracket \varphi\rrbracket^\mathcal{M}=S$ or $\llbracket \varphi\rrbracket^\mathcal{M}=\emptyset$.
\end{fact} 

\begin{proof} By Persistence, if $\varphi$ is true anywhere in $\mathcal{M}$, then it is true at the minimum. But then $\varphi$ must be true everywhere, for if $\mathcal{M},x\nVdash\varphi$, then by Refinability, there is an $x'\sqsubseteq x$ with $\mathcal{M},x'\Vdash\neg\varphi$, which contradicts the fact that $\varphi$ is true at the minimum. (From the topological perspective of Remark \ref{Persp2}, the point is that in a poset $\langle S,\sqsubseteq\rangle$ with a minimum element, the only regular open sets in $\mathcal{O}(S,\sqsubseteq)$ are $S$ and $\emptyset$.)
\end{proof}

What we want to allow is a \textit{distinguished} minimum $\inc$ that does not lead to Fact \ref{collapse}.
 
\begin{definition}[Extended Possibility Frames and Models]\label{ExtendedFrames} An \textit{extended possibility frame} is a tuple $\mathcal{E}=\langle S, \sqsubseteq,\inc, \{R_i\}_{i\in\ind},\adm\rangle$ where: $\langle S,\sqsubseteq\rangle$ is a poset with minimum $\inc$; $R_i$ is a binary relation on $S$ such that $R_i(\inc)=\{\inc\}$ and $xR_i\inc$ for all $x\in S$; $\adm$ is a subset of $\wp(S)$ such that $\inc\in\bigcap \adm$; and the structure $\mathcal{E}_{-}=\langle S_{-}, \sqsubseteq_{-}, \{R_{i_{-}}\}_{i\in\ind},\adm_{-}\rangle$ defined as follows is a \textit{possibility frame} as in Definition \ref{PosFrames}:
\begin{enumerate}
\item $S_{-}=S\setminus\{\inc\}$; 
\item $\sqsubseteq_{-}$ and $R_{i_{-}}$ are the restrictions of $\sqsubseteq$ and $R_i$ to $S_{-}$; 
\item $\adm_{-}=\{X\setminus \{\inc\}\mid X\in\adm\}$. 
\end{enumerate}

An \textit{extended possibility model $\mathcal{M}$ based on $\mathcal{E}$} is a tuple $\mathcal{M}=\langle \mathcal{E},\pi\rangle$ where $\pi\colon \sig\to \adm$. For Fact~\ref{ExtLem} below, we define $\pi_{-}\colon\sig\to \adm_{-}$ by $\pi_{-}(p)=\pi(p)\setminus\{\inc\}$.\hfill $\triangleleft$
\end{definition}

The semantics for extended models essentially ignores the impossible state $\inc$ as follows.

\begin{definition}[Forcing for Extended Models]\label{ExtendedSemantics}
Given an extended possibility model $\mathcal{M}$, $x\in \mathcal{M}$, and $\varphi\in\mathcal{L}(\sig,\ind)$, we define $\mathcal{M},x\Vdash\varphi$ as in Definition \ref{pmtruth1} except with a modified clause for $\neg$:
\begin{enumerate}[label=\arabic*.,ref=\arabic*,resume]
\item $\mathcal{M},x\Vdash\neg\varphi$ iff $\forall x'\sqsubseteq x$: if $x'\not=\inc$, then $\mathcal{M},x'\nVdash\varphi$. \hfill $\triangleleft$
\end{enumerate}
\end{definition}

Now an easy induction shows that all formulas are true at the impossible state $\inc$, using the fact that in an extended model, every $p\in\sig$ is true at $\inc$ given the requirement that $\inc\in \bigcap \adm$.

\begin{fact}[Incoherence] For any extended possibility model $\mathcal{M}$ and $\varphi\in\mathcal{L}(\sig,\ind)$, $\mathcal{M},\inc\Vdash\varphi$.
\end{fact}

\noindent It is also easy to see that Fact \ref{collapse} does not hold for extended possibility models.

As a natural example of an extended possibility frame, we have the following.

\begin{example}[Extended Powerset Possibilization]\label{ExtPowerPoss} Given a world frame $\mathfrak{F}=\langle \wo{W},\{\wo{R}_i\}_{i\in\ind},\wo{A}\rangle$ and a world model $\mathfrak{M}=\langle \mathfrak{F},\wo{V}\rangle$, the \textit{extended powerset possibilizations} of $\mathfrak{F}$ and $\mathfrak{M}$, $\mathfrak{F}^\pow_\inc=\langle S,\sqsubseteq,\inc,\{R_i\}_{i\in\ind},\adm\rangle$ and $\mathfrak{M}^\pow_\inc=\langle \mathfrak{F}^\pow_\inc,\pi\rangle$, are defined by: $S=\wp(\wo{W})$; $X\sqsubseteq Y$ iff $X\subseteq Y$; $\inc=\emptyset$; $XR_iY$ iff $Y\subseteq \mathrm{R}_i[X]$; $\adm = \{\mathord{\downarrow}X\mid X\in \wo{A}\}$; and $\pi(p)=\{X\in S\mid X\subseteq\mathrm{V}(p)\}$. Note that if $\mathfrak{F}$ is a Kripke frame, then $P=\{\mathord{\downarrow}X\mid X\in S\}$. \hfill $\triangleleft$
\end{example}

We can switch back and forth between extended and non-extended frames whenever convenient, by restriction as in Definition \ref{ExtendedFrames} and extension as in Definition \ref{Extending}.

\begin{definition}[Extending Frames]\label{Extending} Given a partial-state frame $\mathcal{F}=\langle S, \sqsubseteq, \{R_i\}_{i\in\ind},\adm\rangle$ and $\inc\not\in S$, define $\mathcal{F}_\inc=\langle S_\inc, \sqsubseteq_\inc,\inc, \{R_{i_\inc}\}_{i\in\ind},\adm_\inc\rangle$ as follows: $S_\inc = S\cup \{\inc\}$; $x\sqsubseteq_\inc y$ iff $x\sqsubseteq y$ or $x=\inc$; $xR_{i_\inc}y$ iff $xR_iy$ or $y=\inc$; and $\adm_\inc = \{X\cup\{\inc\}\mid X\in \adm\}$. 

Given a valuation $\pi\colon \sig\to \adm$, define $\pi_\inc\colon \sig\to \adm_\inc$ by $\pi_\inc (p)=\pi(p)\cup\{\inc\}$. \hfill $\triangleleft$
\end{definition}

The following fact records that restriction and extension work as desired.

\begin{fact}[Equivalence of Extended and Restricted Frames]\label{ExtLem} For any possibility frame $\mathcal{F}$ and extended possibility frame $\mathcal{E}$:
\begin{enumerate}
\item $\mathcal{F}_\inc$ is an extended possibility frame such that $(\mathcal{F}_\inc)_{-}=\mathcal{F}$, and for any possibility model $\langle \mathcal{F},\pi\rangle$, $x\in \mathcal{F}$, and $\varphi\in\mathcal{L}(\sig,\ind)$, $\langle\mathcal{F},\pi\rangle,x\Vdash\varphi$ iff $\langle\mathcal{F}_\inc,\pi_\inc\rangle,x\Vdash\varphi$;
\item $\mathcal{E}_{-}$ is a possibility frame such that $(\mathcal{E}_{-})_\inc$ is isomorphic to $\mathcal{E}$, and for any extended possibility model $\langle \mathcal{E},\pi\rangle$, $x\in\mathcal{E}_{-}$, and $\varphi\in\mathcal{L}(\sig,\ind)$, $\langle\mathcal{E},\pi\rangle,x\Vdash\varphi$ iff $\langle\mathcal{E}_{-},\pi_{-}\rangle,x\Vdash\varphi$.
\end{enumerate}
\end{fact}

\subsection{Functional Frames}\label{FuncFrames} 

For the powerset possibilization $\mathfrak{F}^\pow$ of a world frame $\mathfrak{F}=\langle \wo{W},\{\wo{R}_i\}_{i\in\ind},\wo{A}\rangle$ (Example \ref{PowerPoss}), we defined its accessibility relations by $XR_i^\pow Y$ iff $Y\subseteq \mathrm{R}_i[X]$. As a result, for any possibility $X\in \mathfrak{F}^\pow$, if the set $R_i^\pow(X)=\{Y\in S^\pow\mid XR_i^\pow Y\}$ is nonempty, it has a \textit{maximum} in $\langle S^\pow,\sqsubseteq^\pow\rangle$, i.e., a single possibility $f_i(X)$ of which all possibilities accessible from $X$ are refinements, namely $f_i(X)=\mathrm{R}_i[X]$. This makes possible a \textit{functional} semantics for the modality $\Box_i$, which we may retain even as we generalize away from powerset possibilizations.

\begin{definition}[Quasi-Functional and Functional Possibility Frames]\label{FuncFramesDef} A \textit{quasi-functional} possibility frame is a possibility frame  $\mathcal{F}=\langle S, \sqsubseteq, \{R_i\}_{i\in \ind},\adm\rangle$ satisfying:
\begin{enumerate}[label=\arabic*.,ref=\arabic*]
\item \Rmax{} -- if $R_i(x)\not=\emptyset$, then $R_i(x)$ has a maximum in $\langle S,\sqsubseteq\rangle$.
\end{enumerate}
An \textit{extended quasi-functional} possibility frame is an extended possibility frame $\mathcal{F}=\langle S, \sqsubseteq,\inc, \{R_i\}_{i\in\ind},\adm\rangle$ (Definition \ref{ExtendedFrames}) satisfying:
\begin{enumerate}[label=\arabic*.,ref=\arabic*,resume]
\item \Rmaxe{} -- $R_i(x)$ has a maximum in $\langle S,\sqsubseteq\rangle$.
\end{enumerate}
A \textit{functional} possibility frame is a possibility frame in which each $R_i$ is partially functional, i.e., $xR_iy$ and $xR_iy'$ together imply $y=y'$. An \textit{extended functional} possibility frame is an extended possibility frame $\mathcal{E}$ whose restriction $\mathcal{E}_{-}$ as in Definition~\ref{ExtendedFrames} is functional.

For each $i\in\ind$, let $f_i\colon S\pfun S$ be the partial function such that for all $x\in S$, if $R_i(x)\not=\emptyset$, then $f_i(x)$ is the maximum of $R_i(x)$. We write `$f_i(x)\mathord{\downarrow}$' to indicate that $f_i$ is defined at $x$. In an extended quasi-functional possibility frame, $f_i$ is a total function.

A (quasi-)functional possibility model is a possibility model based on a (quasi-)functional possibility frame, and similarly for extended models.\hfill $\triangleleft$\end{definition} 

Note that every functional frame is a quasi-functional frame. Also note the following.

\begin{fact}[\Rprinc{}]\label{RprincFact} If a frame is quasi-functional and satisfies \Rdown{}, then it satisfies the following condition, and vice versa: \Rprinc{} -- if $R_i(x)\not=\emptyset$, then $R_i(x)$ is a \textit{principal downset} in $\langle S,\sqsubseteq\rangle$.\footnote{As usual, a principal downset in $\langle S,\sqsubseteq\rangle$ is an $X\subseteq S$ such that $X=\mathord{\downarrow}x$ for some $x\in S$.}\end{fact}

Over quasi-functional models, we obtain the following simple semantic clause for $\Box_i$.

\begin{fact}[Functional Semantics]\label{FuncSem} For any extended quasi-functional possibility model $\mathcal{M}$, $x\in \mathcal{M}$, and $\varphi\in\mathcal{L}(\sig,\ind)$:
\[\mbox{$\mathcal{M},x\Vdash \Box_i\varphi$ iff $\mathcal{M},f_i(x)\Vdash \varphi$}.\]
For quasi-functional models that are not extended, $\mathcal{M},x\Vdash \Box_i\varphi$ iff $\mathcal{M},f_i(x)\Vdash \varphi$ or $f_i$ is undefined at $x$.
\end{fact}

 As suggested by Fact \ref{FuncSem}, we can always go from a quasi-functional to an equivalent functional frame.

\begin{proposition}[From Quasi-Functional to Functional Frames]\label{QtoF} For any quasi-functional possibility frame $\mathcal{F}=\langle S, \sqsubseteq, \{R_i\}_{i\in \ind},\adm\rangle$, define its \textit{functionalization} $\mathcal{F}_f=\langle S, \sqsubseteq, \{f_i\}_{i\in \ind},\adm\rangle$ where $f_i$ is the partial function (partially functional relation) given in Definition \ref{FuncFramesDef}. Then:
\begin{enumerate}
\item\label{QtoF1} $\mathcal{F}_f$ is a functional possibility frame;
\item\label{QtoF2} if $\mathcal{F}$ is a full possibility frame, then $\mathcal{F}_f$ is a full possibility frame satisfying the following conditions (cf.~\Rrule{} and \Rwinweak{} from \S~\ref{FullFrames}):
\begin{enumerate}
\item \Frule{} -- if $x'\sqsubseteq x$ and $f_i(x')\comp z$, then $f_i(x)\comp z$;\footnote{We intend this to mean: if $x'\sqsubseteq x$, $f_i(x')\mathord{\downarrow}$, and $f_i(x')\comp z$, then $f_i(x)\mathord{\downarrow}$ and $f_i(x)\comp z$, but we leave the definedness implicit.}
\item \Fwinweak{} -- $\forall y\sqsubseteq f_i(x)$ $\exists x'\sqsubseteq x$ $\forall x''\sqsubseteq x'$: $y\comp f_i(x'')$;\footnote{In \citealt{Holliday2014}, this condition was called \textit{$f$-refinability}. Like in the previous footnote, we take it to mean that if $f_i(x)\mathord{\downarrow}$, then the condition holds with $f_i(x'')\mathord{\downarrow}$.}
\end{enumerate}
\item\label{QtoF3} if $\mathcal{F}$ is separative and satisfies \Rrule{}, then $\mathcal{F}_f$ satisfies:
\begin{enumerate}
\item \Fmon{} -- if $x'\sqsubseteq x$ and $f_i(x')\mathord{\downarrow}$, then $f_i(x)\mathord{\downarrow}$ and $f_i(x')\sqsubseteq f_i(x)$;\footnote{In \citealt{Holliday2014}, this condition was called \textit{$f$-persistence}.}
\end{enumerate} 
\item\label{QtoF4} for all $\pi\colon \sig\to\adm$, $x\in S$, and $\varphi\in\mathcal{L}(\sig,\ind)$, we have $\langle \mathcal{F},\pi\rangle,x\Vdash \varphi$ iff $\langle \mathcal{F}_f,\pi\rangle,x\Vdash \varphi$;  
\item\label{QtoF5} if a logic \textbf{L} is sound and complete with respect to a class $\mathsf{F}$ of quasi-functional possibility frames, then \textbf{L} is sound and complete with respect to the class of functionalizations of frames from $\mathsf{F}$.
\end{enumerate}
\end{proposition}
\begin{proof} For part \ref{QtoF1}, we must first verify that $\mathcal{F}_f$ is indeed a possibility frame. Since we have only modified the accessibility relations, it suffices to check that $\adm$ is still closed under $\blacksquare_i^{\mathcal{F}_f}$ for each $i\in \ind$, as required of a partial-state frame (Definition \ref{PosetMod}). Since $\mathcal{F}$ is quasi-functional, for any $x\in S$ and $Y\in\adm$, we have $x\in \blacksquare_i^\mathcal{F} Y$ iff either $f_i(x)\in Y$ or $f_i(x)$ is undefined (the same point as in Fact \ref{FuncSem}), which is then equivalent to $x\in\blacksquare_i^{\mathcal{F}_f} Y$. Thus, $\blacksquare_i^\mathcal{F} Y=\blacksquare_i^{\mathcal{F}_f}Y$, so $\adm$ is closed under $\blacksquare_i^{\mathcal{F}_f}$ by virtue of being closed under $\blacksquare_i^\mathcal{F}$.

For part \ref{QtoF2}, if $\mathcal{F}$ is full, then by Proposition \ref{ROtoRO}, $\mathcal{F}$ satisfies \Rrule{} and \Rwinweak{}. It is then easy to check that $\mathcal{F}_f$ satisfies \Frule{} and \Fwinweak{}.

For part \ref{QtoF3}, as just noted, \Rrule{} for $\mathcal{F}$ implies \Frule{} for $\mathcal{F}_f$, which gives us that if $x'\sqsubseteq x$ and $f_i(x')\mathord{\downarrow}$, then $f_i(x)\mathord{\downarrow}$. Now if $f_i(x')\not\sqsubseteq f_i(x)$, then by separativity there is a $z\sqsubseteq f_i(x')$, so $f_i(x')\comp z$, such that \textit{not} $f_i(x)\comp z$, which implies $x'\not\sqsubseteq x$ by \Frule{}. This establishes \Fmon{}.

Part \ref{QtoF4} has an obvious proof by induction using Fact \ref{FuncSem} in the $\Box_i$ case.

Part \ref{QtoF5} is immediate from part \ref{QtoF4}.
\end{proof}

From Example \ref{transform} we know that the powerset possibilization $\mathfrak{F}^\pow$ of a world frame $\mathfrak{F}$ is a possibility frame, and $\mathfrak{F}^\pow$ is full if $\mathfrak{F}$ is full, so the observation at the beginning of this section gives us the following.

\begin{example}[Powerset Possibilization Cont.]\label{PowFun} The powerset possibilization $\mathfrak{F}^\pow$ of any world frame (resp.~full world frame) $\mathfrak{F}$ is a \textit{quasi-functional} possibility frame (resp.~full possibility frame).
\end{example}

Putting together Example \ref{PowFun} and Proposition \ref{QtoF}, we can define the \textit{functional powerset possibilization} of a world frame $\mathfrak{F}=\langle \wo{W},\{\wo{R}_i\}_{i\in\ind},\wo{A}\rangle$ as the functionalization $(\mathfrak{F}^\pow)_f$ of the powerset possibilization $\mathfrak{F}^\pow$ of $\mathfrak{F}$. More directly, each (partial) accessibility function $f_i$ in  $(\mathfrak{F}^\pow)_f$ is defined by $f_i(X)=\mathrm{R}_i[X]$ (when $\mathrm{R}_i[X]\neq\varnothing$). For an application of this functional powerset possibilization construction, see \citealt{Benthem2015}.

Together Fact \ref{WtoP1}, Example \ref{PowFun}, and Proposition \ref{QtoF} show that the functional powerset possibilization of a world frame validates exactly the same formulas as the original world frame. Thus, (full) functional possibility frames are as general as (full) world frames in the following sense.

\begin{corollary}[Completeness for Functional Frames] If a logic \textbf{L} is sound and complete with respect to a class $\mathsf{F}$ of world frames, then \textbf{L} is sound and complete with respect to a class of functional possibility frames, viz., the class of functional powerset possibilizations of frames from $\mathsf{F}$. Moreover, this statement holds for \textit{full} world/possibility frames.
\end{corollary} 

Not only can every full world frame be transformed into a semantically equivalent \textit{functional} full possibility frame, but more remarkably, so can every \textit{full possibility frame}. Thus, we could assume without modal-logical loss of generality that our full possibility frames are always functional.

\begin{proposition}[From Relations to Functions]\label{FromReltoFun} For any full possibility frame $\mathcal{F}$, there is a \textit{functional} full possibility frame $\mathcal{F}'$ and a dense and robust possibility morphism $h\colon\mathcal{F}\to\mathcal{F}'$. Thus, by Proposition \ref{Preservation}, for all $\varphi\in\mathcal{L}(\sig,\ind)$, $\mathcal{F}\Vdash\varphi$ iff $\mathcal{F}'\Vdash\varphi$.
\end{proposition}

Proposition \ref{FromReltoFun} will follow from Proposition \ref{UpDownRich}.\ref{UpDownRich3} and Theorem \ref{AlmostBack}.\ref{Inverses0} in \S~\ref{DualEquiv}.

\subsection{Tight Frames}\label{TightSection} 

Kripke frames and their powerset possibilizations are examples of \textit{tight} possibility frames in the following standard sense  \citep[p. 251]{Chagrov1997}, which will be important in \S~\ref{DualityTheory}.

\begin{definition}[Tight Frames]\label{TightFrames} Let $\mathcal{F}$ be a partial-state frame $\mathcal{F}=\langle S,\sqsubseteq,\{R_i\}_{i\in\ind},\adm\rangle$.
\begin{enumerate}
\item\label{Tight1} $\mathcal{F}$ is $R$-\textit{tight} iff for all $i\in\ind$ and $x,y\in S$, if $\forall Z\in\adm$, $x\in\blacksquare_i Z\Rightarrow y\in Z$, then $xR_iy$;
\item\label{Tight2} $\mathcal{F}$ is \textit{$\sqsubseteq$-tight} iff for all $x,y\in S$, if $\forall Z\in\adm$, $x\in Z\Rightarrow y\in Z$, then $y\sqsubseteq x$;
\item\label{Tight3} $\mathcal{F}$ is \textit{tight} iff it is $R$-\textit{tight} and \textit{$\sqsubseteq$-tight}.
 \hfill $\triangleleft$
 \end{enumerate}
\end{definition}
\noindent Note that for any partial-state frame, $xR_iy$ implies $\forall Z\in\adm$, $x\in\blacksquare_i Z\Rightarrow y\in Z$; and for any possibility frame, $y\sqsubseteq x$ implies $\forall Z\in\adm$, $x\in Z\Rightarrow y\in Z$. We may also call a poset $\langle S,\sqsubseteq\rangle$ \textit{tight with respect to $\adm$} when for all $x,y\in S$, if $\forall Z\in\adm$, $x\in Z\Rightarrow y\in Z$, then $y\sqsubseteq x$.

 The notion of \textit{$\sqsubseteq$-tightness} is related to the notions of separativity and differentiation from \S~\ref{SepSec} as follows.

\begin{fact}[$\sqsubseteq$-tightness, Separativity, and Differentiation]\label{TightSepDiff} $\,$
\begin{enumerate}
\item\label{TightSepDiff0} Every $\sqsubseteq$-tight possibility frame is separative;
\item\label{TightSepDiff1} Every separative full possibility frame is $\sqsubseteq$-tight;
\item\label{TightSepDiff2} Every $\sqsubseteq$-tight possibility frame is differentiated. 
\end{enumerate}
\end{fact}
\begin{proof} For part \ref{TightSepDiff0}, if $y\cof x$, then by Fact \ref{CofClose}, $\forall Z\in\adm$,  $x\in Z\Rightarrow y\in Z$, so $y\sqsubseteq x$ by $\sqsubseteq$-\textit{tightness}.

Parts \ref{TightSepDiff1} and \ref{TightSepDiff2} follow from Fact \ref{Sep&Princ} and Definition \ref{DiffFrames}, respectively.
\end{proof} 

The notion of \textit{$R$-tightness} is related to the interplay conditions discussed in \S~\ref{FullFrames} as follows.

\begin{lemma}[$R$-tightness and Interplay Conditions]\label{TightStrong} For any possibility frame $\mathcal{F}$:
\begin{enumerate}
\item\label{TightStrong1} if $\mathcal{F}$ is $R$-tight, then $\mathcal{F}$ satisfies \upR{}, \Rdown{}, and \Rdense{};
\item\label{TightStrong1.5} if $\mathcal{F}$ is $R$-tight and satisfies \Rref{}, then $\mathcal{F}$ is strong, i.e., satisfies \RWin{};
\item\label{TightStrong1.75} if $\mathcal{F}$ is \textit{full}, then $\mathcal{F}$ is $R$-tight iff $\mathcal{F}$ is strong.
\end{enumerate}
\end{lemma}

\begin{proof} The proof of part \ref{TightStrong1} is the same as the proof for Proposition \ref{Representation} that $\mathcal{F}^\tight$ satisfies \upR{}, \Rdown{}, and \Rdense{}. Part \ref{TightStrong1.5} follows from part \ref{TightStrong1}, Proposition \ref{MasterCon}, and Fact \ref{Simp}. 

For part \ref{TightStrong1.75} from left to right, if $\mathcal{F}$ is full, then it satisfies \Rwinweak{} by Corollary \ref{Fullinterplay}, which combines with \upR{}, \Rdown{}, and \Rdense{} from part \ref{TightStrong1} to give us that $\mathcal{F}$ is strong by Proposition \ref{MasterCon}. From right to left, if $\mathcal{F}=\langle S, \sqsubseteq , \{R_i\}_{i\in \ind},\adm\rangle$ satisfies \Rdown{} and \Rdense{}, then by Fact \ref{R(x)RO}, for all $x\in S$, $R_i(x)$ satisfies \textit{persistence} and \textit{refinability}, so if $\mathcal{F}$ is also full, then $R_i(x)\in \adm$. Now suppose that $\forall Z\in \adm$, $R_i(x)\subseteq Z\Rightarrow y\in Z$. Then since $R_i(x)\in\adm$, we have $y\in R_i(x)$, so $xR_iy$. Thus, $\mathcal{F}$ is $R$-tight.
\end{proof}

Recall that for Proposition \ref{Representation}, we took a possibility frame $\mathcal{F}=\langle S,\sqsubseteq,\{R_i\}_{i\in\ind},\adm\rangle$ and constructed a new possibility frame $\mathcal{F}^\tight=\langle S,\sqsubseteq,\{R_i^\tight\}_{i\in\ind},\adm\rangle$ by setting $xR_i^\tight y$ iff $\forall Z\in \adm$,  $x\in\blacksquare_i^\mathcal{F} Z\Rightarrow y\in Z$. Since we showed in the proof of Proposition \ref{Representation}.\ref{Representation2} that for all $Z\in\adm$, $\blacksquare_i^\mathcal{F} Z=\blacksquare_i^{\mathcal{F}^\tight} Z$, it follows that $xR_i^\tight y$ iff $\forall Z\in \adm$,  $x\in\blacksquare_i^{\mathcal{F}^\tight} Z\Rightarrow y\in Z$. Thus, $\mathcal{F}^\tight$ is \textit{$R$-tight}. By applying the same idea to $\sqsubseteq$ in addition to $R_i$, i.e., defining $x'\sqsubseteq^t x$ iff $\forall Z\in\adm$,  $x\in Z \Rightarrow x'\in Z$, and then taking the quotient of the frame with respect to the equivalence relation defined by $x\equiv^t y$ if $x\sqsubseteq^t y$ and $y\sqsubseteq^t x$, it is straightforward to prove the following.

\begin{figure}[h]
\begin{center}
\begin{tikzpicture}[yscale=1, ->,>=stealth',shorten >=1pt,shorten <=1pt, auto,node
distance=2cm,thick,every loop/.style={<-,shorten <=1pt}]
\tikzstyle{every state}=[fill=gray!20,draw=none,text=black]

\node[circle,draw=black!100,fill=black!100,inner sep=0pt,minimum size=.175cm]  (x) at (0,1) {{}};
\node[circle,draw=black!100,fill=black!100,inner sep=0pt,minimum size=.175cm]  (y) at (-1.5,-1.5) {{}};
\node[circle,draw=black!100,fill=black!100,inner sep=0pt,minimum size=.175cm] (z) at (0,-1.5) {{}};
\node[circle,draw=black!100,fill=black!100,inner sep=0pt,minimum size=.175cm] (w) at (1.5,-1.5) {{}};
\node[circle,draw=black!100,fill=black!100,inner sep=0pt,minimum size=.175cm] (v) at (3,-1.5) {{}};

\node[circle,draw=black!100,fill=black!100,inner sep=0pt,minimum size=.175cm]  (x') at (1.5,0) {{}};

\path (x) edge[->] node {{}} (y);
\path (x) edge[->] node {{}} (z);
\path (x) edge[->] node {{}} (w);

\path (x') edge[->] node {{}} (z);
\path (x') edge[->] node {{}} (w);
\path (x') edge[->] node {{}} (v);

\node[rectangle,draw, minimum width = 1cm, minimum height = 1cm] at (-1.5,-1.5) {};
\node[rectangle,draw, minimum width = 1cm, minimum height = 1cm] at (0,-1.5) {};
\node[rectangle,draw, minimum width = 2.5cm, minimum height = 1cm] at (2.25,-1.5) {};

\node[circle,draw=black!100,fill=black!100,inner sep=0pt,minimum size=.175cm]  (a) at (7,1) {{}};
\node[circle,draw=black!100,fill=black!100,inner sep=0pt,minimum size=.175cm]  (b) at (5.5,-1.5) {{}};
\node[circle,draw=black!100,fill=black!100,inner sep=0pt,minimum size=.175cm] (c) at (7,-1.5) {{}};
\node[circle,draw=black!100,fill=black!100,inner sep=0pt,minimum size=.175cm] (d) at (8.5,-1.5) {{}};

\node[circle,draw=black!100,fill=black!100,inner sep=0pt,minimum size=.175cm]  (a') at (8.5,0) {{}};

\path (a) edge[->] node {{}} (b);
\path (a) edge[->] node {{}} (c);
\path (a) edge[->] node {{}} (d);

\path (a') edge[->] node {{}} (c);
\path (a') edge[->] node {{}} (d);

\path (a) edge[->] node {{}} (a');

\node[rectangle,draw, minimum width = 1cm, minimum height = 1cm] at (5.5,-1.5) {};
\node[rectangle,draw, minimum width = 1cm, minimum height = 1cm] at (7,-1.5) {};
\node[rectangle,draw, minimum width = 1cm, minimum height = 1cm] at (8.5,-1.5) {};

\path (w) edge[dotted,->,bend right] node {{}} (d);
\path (v) edge[dotted,->,bend right] node {{}} (d);

\end{tikzpicture}
\end{center}
\caption{A possibility frame (left) and its tightening (right). For both frames, assume that the accessibility relation $R_i$, not shown, is the universal relation, and let $\adm$ be the subalgebra of $\mathrm{RO}(\mathcal{F})$ generated by the sets represented by rectangles. Then the frame on the left is separative but not tight.}\label{TighteningFig}
\end{figure}
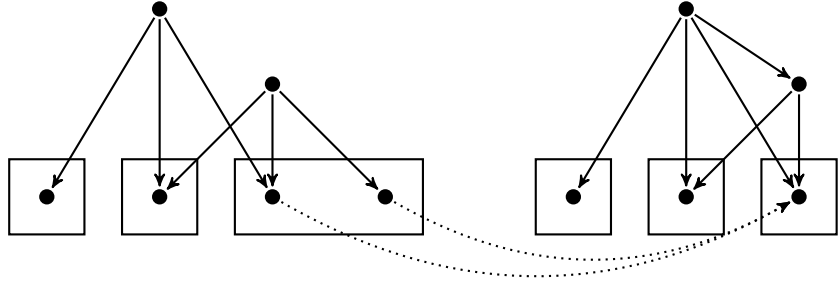

\begin{proposition}[Tightening]\label{Tightening} For any possibility frame $\mathcal{F}$, there is a tight possibility frame $\mathcal{F}^t$ (which is full if $\mathcal{F}$ is) and a surjective robust possibility morphism from $\mathcal{F}$ to $\mathcal{F}^t$. Thus, by Proposition \ref{Preservation}, for all $\varphi\in\mathcal{L}(\sig,\ind)$, $\mathcal{F}\Vdash\varphi$ iff $\mathcal{F}^t\Vdash\varphi$.
\end{proposition}
\noindent For example, Figure \ref{TighteningFig} shows a possibility frame and its tightening with respect to $\sqsubseteq$.

That every possibility frame $\mathcal{F}$ is semantically equivalent to a tight possibility frame will also follow from Corollary \ref{Arbitrary-to-FD} in \S~\ref{Fdes}. 

\subsection{Principal Frames}\label{PrincFrames}  
 
A special feature of the powerset possibilization $\mathfrak{F}^\pow$ of a Kripke frame $\mathfrak{F}$ from Example \ref{PowerPoss} is that the set $\adm$ of admissible propositions in $\mathfrak{F}^\pow$ is the set of all \textit{principal downsets} in the poset $\langle S,\sqsubseteq\rangle$ underlying  $\mathfrak{F}^\pow$, plus the empty set. If we consider the extended powerset possibilization $\mathfrak{F}^\pow_\bot$ as in Example \ref{ExtPowerPoss}, then we can simply say that $\adm$ is the set of all principal downsets. Thus, every possibility $x\in S$ gives rise to a proposition $\mathord{\downarrow}x\in \adm$, expressing \textit{that the possibility $x$ obtains}; and for every nonempty proposition $X\in\adm$, there is a least specific possibility $x\in S$ where that proposition is true, which could be thought of as \textit{the possibility that $X$}. 

We can detach this feature of $\mathfrak{F}^\pow$ from its other features, e.g., that its poset is a complete and atomic Boolean lattice (minus the minimum), or that it is a quasi-functional possibility frame. Let us  consider possibility frames in which $\adm$ is the set of all principal downsets in $\langle S,\sqsubseteq\rangle$ plus $\emptyset$. 

\begin{figure}[h]
\begin{center}
\begin{tikzpicture}[yscale=1, ->,>=stealth',shorten >=1pt,shorten <=1pt, auto,node
distance=2cm,thick,every loop/.style={<-,shorten <=1pt}]
\tikzstyle{every state}=[fill=gray!20,draw=none,text=black]

\node[circle,draw=black!100,fill=black!100,inner sep=0pt,minimum size=.175cm]  (a) at (4,-1) {{}};
\node[circle,draw=black!100,fill=black!100,inner sep=0pt,minimum size=.175cm]  (b) at (3,-2) {{}};
\node[circle,draw=black!100,fill=black!100,inner sep=0pt,minimum size=.175cm]  (x) at (4,-2) {{}};
\node[circle,draw=black!100,fill=black!100,inner sep=0pt,minimum size=.175cm]  (c) at (5,-2) {{}};
\path (a) edge[->] node {{}} (b);
\path (a) edge[->] node {{}} (c);
\path (a) edge[->] node {{}} (x);

\end{tikzpicture}
\end{center}
\caption{A non-principal possibility frame. Assume that $R_i$ is the universal relation and that $\adm=\mathrm{RO}(\mathcal{F})$. The set containing the bottom left and bottom middle possibilities belongs to $\adm$ but is not a principal downset.}\label{NonPrincFrame}
\end{figure}
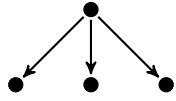

\begin{definition}[Principal Possibility Frame]\label{PrincPoss} A \textit{principal possibility frame} is a possibility frame \\ $\mathcal{F}=\langle S,\sqsubseteq,\{R_i\}_{i\in\ind},\adm\rangle$ in which $\adm$ is the set of all principal downsets in $\langle S,\sqsubseteq\rangle$ plus $\emptyset$.

An \textit{extended principal possibility frame} is an extended possibility frame $\mathcal{F}=\langle S,\sqsubseteq,\bot,\{R_i\}_{i\in\ind},\adm\rangle$ as in \S~\ref{ExtFrames} in which $\adm$ is the set of all principal downsets in $\langle S,\sqsubseteq\rangle$. \hfill $\triangleleft$
\end{definition}
\noindent In contrast to the powerset possibilization of a Kripke frame, a simple example of a \textit{non-principal} frame is the full possibility frame depicted in Figure \ref{NonPrincFrame}.

In order for $\adm$ to be the set of all principal downsets (plus $\emptyset$) \textit{and} satisfy the closure conditions on $\adm$ required of a partial-state frame in Definition \ref{PosetMod}, $\langle S,\sqsubseteq\rangle$ must have a particular form, which we will describe.

Recall that a poset $\langle A,\leq\rangle$ is a lower semilattice iff for all $x,y\in A$, $\{x,y\}$ has a greatest lower bound $x\meet y$ in $\langle A,\leq\rangle$. A lower semilattice $\langle A,\leq\rangle$ with a minimum $\bot$ is a \textit{$p$-semilattice} (for \textit{p}seudocomplemented) iff for all $x\in A$, there is an $x^*\in A$ such that for all $z\in A$, $z\meet x=\bot$ iff $z\leq x^*$. A lower semilattice $\langle A,\leq\rangle$ is an \textit{implicative semilattice} \citep{Nemitz1965} (or a \textit{Brouwerian semilattice} in \citealt{Birkhoff1940} and \citealt{Kohler1981}, or a \textit{relatively pseudocomplemented} semilattice) iff for all $x,y\in A$, there is an $x*y\in A$ such that for all $z\in A$, $z\meet x\leq y$ iff $z\leq x*y$. A \textit{bounded} implicative semilattice is an implicative semilattice with a minimum $\bot$. Any bounded implicative semilattice is a p-semilattice, with $x^*=x*\bot$.\footnote{Of course, we are treating all of these as classes of posets, rather than algebraic structures of particular similarity types.} Finally, note that a Heyting algebra may be defined as a bounded implicative \textit{lattice}. Every finite bounded implicative semilattice is a Heyting algebra, but not every infinite bounded implicative semilattice is a Heyting algebra.  

When we add an accessibility relation $R_i$ to $\langle A,\leq\rangle$, we need one more definition.

\begin{definition}[$\Rlatbox_i$ Operation]\label{DotOp} Given a poset $\langle A,\leq\rangle$ with minimum $\bot$ and a binary relation $R_i$ on $A$, define a partial operation $\Rlatbox_i$ on $A$ as follows: for $y\in A$, if the set $\blacksquare_i \mathord{\downarrow}y$, i.e., $\{x\in A\mid R_i(x)\subseteq \mathord{\downarrow}y\}$ (Definition~\ref{PosetMod}), is a nonempty downset and contains a \textit{maximum} element, then $\Rlatbox_iy=\mathrm{max}(\blacksquare_i \mathord{\downarrow}y)$. If $\blacksquare_i \mathord{\downarrow}y$ is empty, then $\Rlatbox_i y=\bot$. If $\blacksquare_i \mathord{\downarrow}y$ is nonempty but is not a downset or does not contain a maximum element, then $\Rlatbox_i y$ is undefined. \hfill $\triangleleft$
\end{definition}

The following fact is a straightforward consequence of the definitions just given.
 
\begin{fact}[Principal Partial-State Frames]\label{ArePsem} Let $\mathcal{F}=\langle S,\sqsubseteq,\{R_i\}_{i\in\ind},\adm\rangle$ be such that $\langle S,\sqsubseteq\rangle$ is a nonempty poset, $R_i$ is a binary relation on $S$, and $\adm$ is the set of all principal downsets in $\langle S,\sqsubseteq\rangle$ plus $\emptyset$. 

Then the following are equivalent:
\begin{enumerate}
\item $\mathcal{F}$ is a partial-state frame;
\item $\adm$ is closed under the operations $\cap$, $\supset$, and $\blacksquare_i$ from Definition \ref{PosetMod};
\item the extension $\langle S_\bot,\sqsubseteq_\bot\rangle$ of $\langle S,\sqsubseteq\rangle$ (Definition \ref{Extending}) is a \textit{bounded implicative semilattice} on which the operation $\Rlatbox_i$ defined from $R_{i_\bot}$ (Definition \ref{Extending}) is a \textit{total} operation.
\end{enumerate} 
\end{fact}

In addition to the closure requirements on $\adm$ in partial-state frames, there is the requirement in possibility frames that $\adm\subseteq\mathrm{RO}(S,\sqsubseteq)$ (Definition \ref{PosFrames}). The following is immediate from Definition \ref{TightFrames}.\ref{Tight2} and Fact~\ref{Sep&Princ}.

\begin{fact}[Principal Possibility Frames, Separativity, and $\sqsubseteq$-tightness]\label{PrincPossSep} If $\adm$ is the set of all principal downsets in $\langle S,\sqsubseteq\rangle$ plus $\emptyset$, then $\langle S,\sqsubseteq\rangle$ is $\sqsubseteq$-tight with respect to $\adm$, and the following are equivalent:
\begin{enumerate}
\item $\adm\subseteq \mathrm{RO}(S,\sqsubseteq)$;
\item $\langle S,\sqsubseteq\rangle$ is separative.
\end{enumerate}
\end{fact}

The requirements on principal possibility frames from Facts \ref{ArePsem} and \ref{PrincPossSep} are related as follows.

\begin{fact}[Boolean Lattices]\label{Boolean} The following are equivalent:
\begin{enumerate}
\item $\langle A,\leq\rangle$ is a p-semilattice such that the restriction of $\leq$ to $A\setminus\{\bot\}$ is a separative partial order;
\item $\langle A,\leq\rangle$ is a Boolean lattice.
\end{enumerate}
\end{fact}

\begin{proof} Frink \citeyearpar{Frink1962} showed that $\langle A,\leq\rangle$ is a Boolean lattice iff $\langle A,\leq\rangle$ is a p-semilattice such that for all $x\in A$, $x^{**}= x$. Thus, we need only show that a p-semilattice $\langle A,\leq\rangle$ obeys $x^{**}= x$ iff the restriction of $\leq$ to $A\setminus\{\bot\}$ is separative. Let $\sqsubseteq$ be the restriction.

For any p-semilattice $\langle A,\leq\rangle$ and $x\in A$, $x\leq x^{**}$, so we need only show that $x^{**}\leq x$. If $x^{**}=\bot$, then $x^{**}\leq x$, so suppose $x^{**}\not=\bot$.  Then we claim that $x^{**} \sqsubseteq_\cofsub x$ (Definition \ref{CoRef}). For suppose $z\sqsubseteq x^{**}$, so $z\leq x^{**}$ and $z\not=\bot$. Then $z\meet x\not=\bot$, for otherwise $z\leq x^*$, contradicting $z\leq x^{**}$ and $z\not=\bot$. Thus, $z\meet x\sqsubseteq x$. So we have shown that for every $z\sqsubseteq x^{**}$ there is a $z'\sqsubseteq z$, namely $z'=z\meet x$, such that $z'\sqsubseteq x$, which means $x^{**} \sqsubseteq_\cofsub x$. Then by Definition \ref{SepFrames}, if $\sqsubseteq$ is separative, $x^{**} \cof x$ implies $x^{**}\sqsubseteq x$ and hence $x^{**}\leq x$. 

In the other direction, to show that $\sqsubseteq$ is separative, for $x,y\in A\setminus \{\bot\}$, assume that $x\not\sqsubseteq y$. It follows that $x\meet y^*\not=\bot$, for otherwise $x\leq y^{**}$, which with $y^{**}\leq y$ and the transitivity of $\leq$ implies $x\leq y$ and hence $x\sqsubseteq y$, contradicting our assumption. Thus, $x\meet y^*\sqsubseteq x$. Moreover, for any  $x''\sqsubseteq x\meet y^*$, so $x''\not=\bot$, we have $x''\not\sqsubseteq y$. So from the assumption that $x\not\sqsubseteq y$, we have shown that there is an $x'\sqsubseteq x$, namely $x'=x\meet y^*$, such that for all $x''\sqsubseteq x'$, $x''\not\sqsubseteq y$. Hence $\sqsubseteq$ is separative.\end{proof} 

Combining the previous three facts, we can characterize principal possibility frames as follows. 

\begin{fact}[Characterization of Principal Possibility Frames]\label{PrincEquiv} Let $\mathcal{F}=\langle S,\sqsubseteq,\{R_i\}_{i\in\ind},\adm\rangle$ be such that $\langle S,\sqsubseteq\rangle$ is a nonempty poset, $R_i$ is a binary relation on $S$, and $\adm\subseteq \wp (S)$. Then the following are equivalent:
\begin{enumerate}
\item\label{PrincEquiv1} $\mathcal{F}$ is a principal possibility frame;
\item\label{PrincEquiv2} $\adm$ is the set of all principal downsets in $\langle S,\sqsubseteq\rangle$ plus $\emptyset$, and $\langle S_\bot,\sqsubseteq_\bot\rangle$ is a Boolean lattice on which the $\Rlatbox_i$ defined from $R_{i_\bot}$ is a \textit{total} operation.
\end{enumerate}
\end{fact}

Now let us consider the interplay conditions on $R_i$ and $\sqsubseteq$ that hold for any principal possibility frame. Recall that in Proposition \ref{ROtoRO} we showed the following conditions to be necessary and sufficient for a poset $\langle S,\sqsubseteq\rangle$ with accessibility relation $R_i$ to be such that $\mathrm{RO}(S,\sqsubseteq)$ is closed under $\blacksquare_i$:
\begin{itemize}
\item \Rrule{} -- if ${x}'\sqsubseteq{x}$ and ${x}' R_i{y}'\comp{z}$, then $\exists {y}$: ${x} R_i{y}\comp{z}$;
\item \Rwinweak{} -- if ${x}R_i{y}$, then $\forall {y'}\sqsubseteq {y}$ $\exists{x'}\sqsubseteq {x}$ $\forall {x''}\sqsubseteq {x'}$ $\exists {y''}\comp y'$: ${x''} R_i{y''}$.
\end{itemize}
Thus, every full possibility frame satisfies \Rrule{} and \Rwinweak{}. Also recall that in Lemma \ref{TightStrong}.\ref{TightStrong1.75}, we observed that if $\mathcal{F}$ is full, then $\mathcal{F}$ is $R$-tight iff $\mathcal{F}$ is strong, i.e., iff it satisfies \RWin{} (Definition \ref{StrongPoss}).

\begin{proposition}[Interplay of Accessibility and Refinement in Principal Frames]\label{InterPrinc} For any principal possibility frame $\mathcal{F}$:
\begin{enumerate}
\item\label{InterPrinc1} $\mathcal{F}$ satisfies \Rrule{} and \Rwinweak{}; 
\item\label{InterPrinc2} $\mathcal{F}$ is $R$-tight iff $\mathcal{F}$ is strong.
\end{enumerate}
\end{proposition}
\begin{proof}  Let $\mathcal{F}=\langle S,\sqsubseteq, \{R_i\}_{i\in\ind},\adm\rangle$ be a principal possibility frame.

For part \ref{InterPrinc1}, suppose $\mathcal{F}$ does not satisfy \Rrule{}, so we have ${x}'\sqsubseteq{x}$ and ${x}' R_i{y}'\comp{z}$, but for all $y$, if $xR_iy$, then $y\incomp z$, so $y\sqsubseteq z^*$. Hence $x\in \blacksquare_i \mathord{\downarrow} z^*$, but $x'\not\in  \blacksquare_i \mathord{\downarrow} z^*$, so  $\blacksquare_i \mathord{\downarrow} z^*$ does not satisfy \textit{persistence}. Then since $\mathcal{F}$ is a possibility frame, $\blacksquare_i \mathord{\downarrow} z^*\not\in\adm$.  Yet $\mathord{\downarrow} z^*$ is a principal downset and hence in $\adm$, since $\mathcal{F}$ is a principal possibility frame. Thus, $\adm$ is not closed under $\blacksquare_i$, so $\mathcal{F}$ is not a partial-state frame, a contradiction.

Next, suppose $\mathcal{F}$ does not satisfy \Rwinweak{}, so we have $xR_iy$ and $\exists y'\sqsubseteq y$ $\forall x'\sqsubseteq x$ $\exists x''\sqsubseteq x'$ $\forall y''$, if $x''R_iy''$, then $y''\incomp y'$, so $y''\sqsubseteq y'^*$. Hence $\forall x'\sqsubseteq x$ $\exists x''\sqsubseteq x'$: $x''\in\blacksquare_i \mathord{\downarrow} y'^*$. But since $xR_iy$ and $y'\sqsubseteq y$, $x\not\in \blacksquare_i \mathord{\downarrow} y'^*$. It follows that $\blacksquare_i\mathord{\downarrow} y'^*$ does not satisfy \textit{refinability}, so $\blacksquare_i \mathord{\downarrow} y'^*\not\in \adm$. Yet $ \mathord{\downarrow} y'^*$ is a principal downset and hence in $\adm$. Thus, $\adm$ is not closed under $\blacksquare_i$, so $\mathcal{F}$ is not a partial-state frame, a contradiction.

For part \ref{InterPrinc2}, the proof of the left-to-right direction is the same as the proof of the left-to-right direction of Lemma \ref{TightStrong}.\ref{TightStrong1.75}, but using the fact just shown that any principal frame satisfies \Rwinweak{}. From right to left, to show that $\mathcal{F}$ is $R$-tight, suppose that \textit{not} $xR_iy$. Then we must find a $Z\in \adm$ such that $R_i(x)\subseteq Z$ but $y\not\in Z$. Since \textit{not} $xR_iy$, by the \Rdense{} property that follows from \RWin{} (Proposition \ref{MasterCon}) there is a $y'\sqsubseteq y$ such that (i) for all $y''\sqsubseteq y'$, \textit{not} $xR_iy''$. Now we claim that for all $z\in R_i(x)$, $z\incomp y'$. For if $xR_iz$ and $z\comp y'$, so there is a $z'\sqsubseteq z$ such that $z'\sqsubseteq y'$, then by the \Rdown{} property that follows from \RWin{} (Proposition \ref{MasterCon}), $xR_iz$ implies $xR_iz'$, which with $z'\sqsubseteq y'$ contradicts (i). Thus, for all $z\in R_i(x)$, $z\incomp y'$, which implies $R_i(x)\subseteq \mathord{\downarrow} y'^*$. But since $y'\sqsubseteq y$, $y\not\in \mathord{\downarrow} y'^*$. Finally, since $\mathcal{F}$ is principal, $\mathord{\downarrow} y'^*\in\adm$, so our desired $Z$ is $\mathord{\downarrow} y'^*$.
\end{proof}

For \textit{quasi-functional} frames as in \S~\ref{FuncFrames}, we can say more about the interplay of accessibility and refinement. By the same reasoning as in the proof of Proposition \ref{QtoF}.\ref{QtoF3}, we have the following.

\begin{fact}[Quasi-Functional Principal Frames]\label{FPrincMon} If $\mathcal{F}$ is a quasi-functional principal possibility frame, then $\mathcal{F}$ satisfies \Fmon{}: if $x'\sqsubseteq x$ and $f_i(x')\mathord{\downarrow}$, then $f_i(x)\mathord{\downarrow}$ and $f_i(x')\sqsubseteq f_i(x)$.\footnote{This also holds for \textit{extended} quasi-functional principal frames (see Definitions \ref{ExtendedFrames}, \ref{FuncFramesDef}, and \ref{PrincPoss}), for which $f_i(\bot)=\bot$.}
\end{fact}

From the assumption that $\mathcal{F}$ is a possibility frame in which $P$ is the set of all principal downsets in $\langle S,\sqsubseteq\rangle$ plus $\emptyset$, we can deduce facts not only about the structure of $\langle S,\sqsubseteq\rangle$ (Fact \ref{PrincEquiv}) and the interplay of $R_i$ and $\sqsubseteq$ (Proposition \ref{InterPrinc} and Fact \ref{FPrincMon}), but also about the nature of possibility morphisms to $\mathcal{F}$. 

\begin{fact}[Possibility Morphisms to Principal Frames]\label{PrincOrd} For any possibility frame $\mathcal{F}$ and principal possibility frame $\mathcal{F}'$, if $h\colon \mathcal{F}\to\mathcal{F}'$ is a possibility morphism, then $h$ satisfies:
\begin{enumerate}
\item\label{PrincOrd1} \SqForth{} -- if $y\sqsubseteq x$, then $h(y)\sqsubseteq' h(x)$;
\item\label{PrincOrd2} \SqBack{} -- if $y'\sqsubseteq' h(x)$, then  $\exists y\sqsubseteq x$: $h(y)\sqsubseteq' y'$.
\end{enumerate}
\end{fact} 
\begin{proof} For part \ref{PrincOrd1}, since all principal downsets in $\mathcal{F}'$ are admissible propositions in $\mathcal{F}'$, Fact \ref{SepOrd}.\ref{SepOrd1} implies that $h$ satisfies \SqForth{}. For part \ref{PrincOrd2}, by the right-to-left direction of \SqMatch{}, we have that $\forall x\in\mathcal{F}$ $\forall X'\in\adm'$, if $\mathord{\downarrow}'h(x)\cap X'\neq\emptyset$, then $\mathord{\downarrow}x\cap h^{-1}[X']\neq\emptyset$.  Since $\mathcal{F}'$ is a principal frame, for any $y'\in\mathcal{F}'$, $\mathord{\downarrow}'y'\in\adm'$, so by \SqMatch{}, if $\mathord{\downarrow}'h(x)\cap \mathord{\downarrow}'y'\neq\emptyset$, then $\mathord{\downarrow}x\cap h^{-1}[\mathord{\downarrow}'y']\neq\emptyset$. Then since $y'\sqsubseteq ' h(x)$ implies  $\mathord{\downarrow}'h(x)\cap \mathord{\downarrow}'y'\neq\emptyset$, it implies $\mathord{\downarrow}x\cap h^{-1}[\mathord{\downarrow}'y']\neq\emptyset$, which means there is a $y\sqsubseteq x$ such that $h(y)\sqsubseteq' y'$.
\end{proof}

Now that we have an idea of the structure of principal possibility frames, let us return to the point made at the beginning of this section: in an extended principal possibility frame, for every proposition $X\in\adm$, there is a least specific possibility $x\in X$. Then since the truth set $\llbracket\varphi\rrbracket^\mathcal{M}$ of a formula $\varphi$ in a partial-state model $\mathcal{M}$ always belongs to $\adm$ (Fact \ref{TruthSub}), we have the following.

\begin{fact}[Principal Truth Sets]\label{propositions1} If $\mathcal{F}$ is an extended principal possibility frame, then for any extended possibility model $\mathcal{M}$ based on $\mathcal{F}$ and $\varphi\in\mathcal{L}(\sig,\ind)$, there is an element $\Vert \varphi\rVert^\mathcal{M}\in S$ such that $\llbracket\varphi\rrbracket^\mathcal{M}=\mathord{\downarrow}\Vert \varphi\rVert^\mathcal{M}$. For any $\varphi,\psi\in\mathcal{L}(\sig,\ind)$: $\lVert\neg\varphi\rVert=\lVert \varphi\rVert^*$, $\lVert \varphi\wedge\psi\rVert=\lVert\varphi\rVert\meet\lVert \psi\rVert$, and $\lVert \Box_i\varphi\rVert=\Rlatbox_i\lVert\varphi\rVert$.
\end{fact}

Fact \ref{propositions1} shows that principal possibility frames will easily transform into modal \textit{algebras} (see \S~\ref{DualityTheory}).

In \S~\ref{VtoPossSection}, we will identify a vast source of principal possibility frames. In particular, we will show that any $\mathcal{V}$-BAO can be turned into a semantically equivalent principal possibility frame, and any $\mathcal{T}$-BAO can be turned into a semantically equivalent \textit{quasi-functional} principal possibility frame (and hence into a semantically equivalent \textit{functional} principal possibility frame by Proposition \ref{QtoF}).

\subsubsection{Lattice-Complete Principal Frames}\label{LatCompSec}

Since principal frames are always based on Boolean lattices, an important special case is the following.

\begin{definition}[Lattice-Complete Principal Frames]\label{LatCompPrinc} A principal possibility frame $\mathcal{F}=\langle S,\sqsubseteq,\{R_i\}_{i\in\ind},\adm\rangle$ is \textit{lattice-complete} iff $\langle S_\bot,\sqsubseteq_\bot\rangle$ is a complete Boolean lattice. \hfill $\triangleleft$
\end{definition}

When $\langle S_\bot,\sqsubseteq_\bot\rangle$ is a \textit{complete} Boolean lattice, the property characterizing principal frames in Fact \ref{PrincEquiv} that $\boxdot_i$ is a total operation on $\langle S_\bot,\sqsubseteq_\bot\rangle$ not only implies the interplay conditions \Rrule{} and \Rwinweak{} as in Proposition \ref{InterPrinc}, but it is also \textit{implied by} these interplay conditions, as shown by the following.

\begin{fact}[Characterization of Lattice-Complete Principal Frames]\label{CharLatComp} Let $\mathcal{F}=\langle S,\sqsubseteq,\{R_i\}_{i\in\ind},\adm\rangle$ be such that $\langle S,\sqsubseteq\rangle$ is a nonempty poset, $R_i$ is a binary relation on $S$, and $\adm\subseteq \wp (S)$. Then the following are equivalent:
\begin{enumerate}
\item\label{CharLatComp1} $\mathcal{F}$ is a lattice-complete principal possibility frame;
\item\label{CharLatComp2} $\adm$ is the set of all principal downsets in $\langle S,\sqsubseteq\rangle$ plus $\emptyset$, $\langle S_\bot,\sqsubseteq_\bot\rangle$ is a complete Boolean lattice, and $R_i$ satisfies \Rrule{} and \Rwinweak{}.
\end{enumerate}
\end{fact}

\begin{proof} The implication from part \ref{CharLatComp1} to part \ref{CharLatComp2} follows from Definitions \ref{PrincPoss} and \ref{LatCompPrinc} and Proposition \ref{InterPrinc}.

From part \ref{CharLatComp2} to part \ref{CharLatComp1}, given the characterization of principal frames in terms of $\boxdot_i$ being a total operation in Fact \ref{PrincEquiv}, it suffices to show that if $R_i$ satisfies \Rrule{} and \Rwinweak{}, then $\boxdot_i$ is such a total operation, which is to say that for any $y\in S_\bot$, if $\{x\in S_\bot\mid R_{i_\bot}(x)\subseteq \mathord{\downarrow} y\}$ is nonempty (where $\mathord{\downarrow} y=\{y'\in S_\bot\mid y'\sqsubseteq_\bot y\}$), then it is a downset and has a \textit{maximum} in $\langle S_\bot,\sqsubseteq_\bot\rangle$. For the downset condition, \Rrule{} for $R_i$ implies that $\{x\in S_\bot\mid R_{i_\bot}(x)\subseteq \mathord{\downarrow} y\}$ is a downset by reasoning similar to that in the proof of Proposition~\ref{ROtoRO} (using Facts~\ref{Boolean} and \ref{Sep&Princ}.\ref{Sep&Princ1} for the \textit{refinability} of $\mathord{\downarrow}y$ restricted to $S$). As for its having a maximum, by the completeness assumption, this is equivalent to 
\begin{equation}
\bigvee\{x\in S_\bot\mid R_{i_\bot}(x)\subseteq \mathord{\downarrow} y\}\in \{x\in S_\bot\mid R_{i_\bot}(x)\subseteq \mathord{\downarrow} y\}.\label{LatCompEquiv}
\end{equation}
Recall from the definition of $R_{i_\bot}$ in Definition \ref{Extending} that $R_{i_\bot}(\bot)=\{\bot\}$ and for $x\in S$, $R_{i_\bot}(x)=R_i(x)\cup\{\bot\}$. If $\bigvee\{x\in S_\bot\mid R_{i_\bot}(x)\subseteq \mathord{\downarrow} y\}=\bot$, then (\ref{LatCompEquiv}) holds, so suppose it is not $\bot$. Then to establish (\ref{LatCompEquiv}), suppose for \textit{reductio} that $\bigvee\{x\in S_\bot\mid R_{i_\bot}(x)\subseteq \mathord{\downarrow}y\}R_{i_\bot}z$ but $z\not\sqsubseteq_\bot y$. Then $z\not=\bot$, so $\bigvee\{x\in S_\bot\mid R_{i_\bot}(x)\subseteq \mathord{\downarrow} y\}R_i z$, and $z\meet y^*\not=\bot$, so  $z\meet y^*\sqsubseteq z$. Let $z'=z\meet y^*$. Then by \Rwinweak{} for $R_i$, there is an $x'\sqsubseteq  \bigvee\{x\in S\mid R_i(x)\subseteq \mathord{\downarrow}y\}$ such that (i) $\forall x''\sqsubseteq x'$ $\exists z''\comp z'$: $x''R_iz''$. Given $x'\sqsubseteq \bigvee\{x\in S_\bot\mid R_{i_\bot}(x)\subseteq \mathord{\downarrow}y\}$, there is an $x\in \{x\in S_\bot\mid R_{i_\bot}(x)\subseteq \mathord{\downarrow}y\}$ such that $x'\meet x\not=\bot$ (cf.~Fact \ref{Useful}), so $x'\meet x\sqsubseteq x'$. Let $x''=x'\meet x$. Then by (i), there is a $z''\comp z'$ such that $x''R_i z''$. Given $x''R_iz''\comp z'$ and $x''\sqsubseteq x$, it follows by \Rrule{} that there is a $u$ such that $xR_iu\comp z'$. Then since $z'=z\meet y^*$, we have $xR_iu\comp y^*$, which contradicts the fact that $R_{i_\bot}(x)\subseteq\mathord{\downarrow}y$.\end{proof}

Recall from Proposition \ref{InterPrinc}.\ref{InterPrinc2} that a principal frame is \textit{$R$-tight} iff it is \textit{strong}. For lattice-complete principal frames, we can add the following observation concerning $R$-tightness. 

\begin{fact}[\textit{R-tight} and \Rprinc{}]\label{TightPrinc} If $\mathcal{F}$ is a lattice-complete principal possibility frame, then $\mathcal{F}$ satisfies the \textit{R-tight} condition iff $\mathcal{F}$ satisfies the \Rprinc{} condition from \S~\ref{FuncFrames}: $R_i(x)$ is a principal downset if $R_i(x)\not=\emptyset$.
\end{fact}

\begin{proof} From left to right, assuming $\mathcal{F}$ satisfies \textit{R-tight} and $R_i(x)\not=\emptyset$, we show that $R_i(x)=\mathord{\downarrow}\bigvee R_i(x)$, where $\bigvee$ is the join operation in the complete Boolean lattice minus $\bot$ that underlies $\mathcal{F}$. By Lemma \ref{TightStrong}.\ref{TightStrong1}, \textit{R-tight} implies \Rdown{}, so $R_i(x)$ is a downset. Thus, to show that $R_i(x)=\mathord{\downarrow}\bigvee R_i(x)$, it suffices to show that $\bigvee R_i(x)\in R_i(x)$. Consider any $Z\in\adm$ such that $x\in\blacksquare_i Z$, so $R_i(x)\subseteq Z$. If $\bigvee R_i(x)\not\in Z$, then by \textit{refinability} for $Z$, there is a $z\sqsubseteq \bigvee R_i(x)$ such that for all $z'\sqsubseteq z$, $z'\not\in Z$. Since $z\sqsubseteq \bigvee R_i(x)$, there is some $y\in R_i(x)$ such that $z\meet y\not=\bot$ (cf. Fact \ref{Useful}), so $z\meet y\sqsubseteq y$, but $z\meet y\not\in Z$ by the previous sentence. Then by \textit{persistence} for $Z$, $y\not\in Z$, which contradicts the fact that $y\in R_i(x)\subseteq Z$. Thus, $\bigvee R_i(x)\in Z$. Since this holds for any $Z\in \adm$ such that $x\in\blacksquare_i Z$, \textit{R-tight} implies $\bigvee R_i(x)\in R_i(x)$.

From right to left, assume $\mathcal{F}$ satisfies \Rprinc{} and that for all $Z\in\adm$, $x\in\blacksquare_i Z$ implies $y\in Z$. It follows that $R_i(x)\not=\emptyset$, for $R_i(x)=\emptyset$ implies $x\in\blacksquare_i\emptyset$, which with our assumption implies the contradiction $y\in\emptyset$. Thus, by \Rprinc{}, $R_i(x)=\mathord{\downarrow}f_i(x)$ for some $f_i(x)$. Now suppose \textit{not} $xR_iy$, so $y\not\sqsubseteq f_i(x)$. Then since $\mathcal{F}$ is separative by Fact \ref{Boolean}, there is a $z\sqsubseteq y$ such that \textit{not} $z\comp f_i(x)$. Then since $R_i(x)=\mathord{\downarrow}f_i(x)$, it follows that $x\in\blacksquare_i \mathord{\downarrow}z^*$, and $\mathord{\downarrow}z^*\in\adm$ since $\mathcal{F}$ is a principal frame, so our initial assumption implies that $y\in \mathord{\downarrow}z^*$, which contradicts $z\sqsubseteq y$. Thus, $xR_iy$, which shows that $\mathcal{F}$ satisfies \textit{$R$-tight}.
\end{proof}

Finally, it is important to note that lattice-complete principal possibility frames are a special case of \textit{full} possibility frames, as shown by the following. 
 
\begin{fact}[Completeness and Fullness]\label{FullTFAE} The following are equivalent:
\begin{enumerate}
\item\label{FullTFAE1} $\mathcal{F}$ is a lattice-complete principal possibility frame;
\item\label{FullTFAE2} $\mathcal{F}$ is a full possibility frame in which $\langle S_\bot,\sqsubseteq_\bot\rangle$ is a complete Boolean lattice.
\end{enumerate}
\end{fact}

\begin{proof} It suffices to show that if $\langle S_\bot,\sqsubseteq_\bot\rangle$ is a complete Boolean lattice, then the set of principal downsets in $\langle S,\sqsubseteq\rangle$ plus $\emptyset$, which is $\adm$ in a principal frame, is exactly the set of \textit{all} $X\subseteq S$ satisfying \textit{persistence} and \textit{refinability} in $\langle S,\sqsubseteq\rangle$, which is $\adm$ in a full frame. If $\langle S_\bot,\sqsubseteq_\bot\rangle$ is a Boolean lattice, then $\langle S,\sqsubseteq\rangle$ is separative by Fact \ref{Boolean}, so every principal downset satisfies not only \textit{persistence} but also \textit{refinability} by Fact \ref{Sep&Princ}.\ref{Sep&Princ1}; and of course $\emptyset$ satisfies \textit{persistence} and \textit{refinability}. In the other direction, suppose $X\subseteq S$ satisfies \textit{persistence} and \textit{refinability}. If $X=\emptyset$, we are done, so suppose $X\not=\emptyset$. Since $\langle S_\bot,\sqsubseteq_\bot\rangle$ is \textit{complete}, $\bigvee X$ exists in $\langle S_\bot,\sqsubseteq_\bot\rangle$, and since $X\subseteq S$, $\bot\not\in X$, so $\bigvee X$ exists in $\langle S,\sqsubseteq\rangle$. Now we claim that $X=\mathord{\downarrow}\bigvee X$, so $X$ is a principal downset. Since  $X\subseteq \mathord{\downarrow}\bigjoin X$ is immediate, we need only show that $\mathord{\downarrow}\bigjoin X\subseteq  X$. By \textit{persistence}, it suffices to show that  $\bigjoin X\in X$.  If $\bigjoin X\not\in X$, then \textit{refinability} in $\langle S,\sqsubseteq\rangle$ implies that there is an $x\in S$ such that $x\sqsubseteq \bigjoin X$ and for all  $x'\sqsubseteq x$, $x'\not\in  X$. Now for any $y\in X$, if $x\wedge y\neq\inc$, then $x\wedge y\sqsubseteq y$, so $x\wedge y\in X$ by \textit{persistence}, but also $x\wedge y\sqsubseteq x$, so $x\wedge y\not\in X$ by the previous sentence; hence for all $y\in X$, $x\wedge y=\inc$, so $y\sqsubseteq x^*$. Hence $\bigjoin X\sqsubseteq x^*$, which contradicts the fact that $x\sqsubseteq \bigjoin X$. Thus, $\bigjoin X\in X$, as desired.\end{proof}

\subsection{Rich Frames}\label{RichFrames} 

If we combine the idea of lattice-complete principal possibility frames from \S~\ref{LatCompSec} with the idea of \textit{strong} possibility frames from \S~\ref{FullFrames}, we arrive at a special class of full possibility frames that will be of great importance in the duality theory of \S~\ref{DualityTheory}. We will begin with a very direct definition of this class of frames, which does not refer back to further definitions from previous sections. As in \S~\ref{ExtFrames}, if $\langle S,\sqsubseteq\rangle$ is a poset, then $\langle S_\bot,\sqsubseteq_\bot\rangle$ is the result of extending $\langle S,\sqsubseteq\rangle$ with a new minimum element $\bot$.

\begin{definition}[Rich Frames]\label{RichPossFrames} A \textit{rich} frame is a tuple $\mathcal{F}=\langle S,\sqsubseteq,\{R_i\}_{i\in\ind},\adm\rangle$ where $\langle S,\sqsubseteq\rangle$ is a nonempty poset such that $\langle S_\bot,\sqsubseteq_\bot\rangle$ is a complete Boolean lattice, $\adm$ is the set of all principal downsets in $\langle S,\sqsubseteq\rangle$ plus $\emptyset$, and each $R_i$ is a binary relation on $S$ satisfying:
\begin{itemize}
\item \RWin{} -- $xR_iy$ iff $\forall y'\sqsubseteq y$ $\exists x'\sqsubseteq x$ $\forall x''\sqsubseteq x'$ $\exists y''\sqsubseteq y'$: $x'' R_i y''$ (see Figure \ref{Rwinfig}).  \hfill $\triangleleft$
\end{itemize}
\end{definition} 

\noindent Recall from \S~\ref{FullFrames} that possibility frames satisfying \RWin{} are the \textit{strong} possibility frames. 

The first thing to prove about rich frames is that they are in fact possibility frames in the sense of Definition \ref{PosFrames} and therefore strong possibility frames. The following fact characterizes rich frames as a special kind of principal possibility frames or equivalently as a special kind of full possibility frames.

\begin{fact}[Characterizations of Rich Frames]\label{RichFullStrong} For any $\mathcal{F}=\langle S,\sqsubseteq,\{R_i\}_{i\in\ind},\adm\rangle$, the following are equivalent:
\begin{enumerate}
\item\label{RichFullStrongA} $\mathcal{F}$ is a rich frame;
\item\label{RichFullStrongB} $\mathcal{F}$ is a principal possibility frame that is lattice-complete and satisfies $R$-tight;
\item\label{RichFullStrongB2} $\mathcal{F}$ is a principal possibility frame that is lattice-complete and satisfies \Rprinc;
\item\label{RichFullStrongC} $\mathcal{F}$ is a full possibility frame in which $\langle S_\bot,\sqsubseteq_\bot\rangle$ is a complete Boolean lattice, satisfying $R$-tight;
\item\label{RichFullStrongC2} $\mathcal{F}$ is a full possibility frame in which $\langle S_\bot,\sqsubseteq_\bot\rangle$ is a complete Boolean lattice, satisfying   \Rprinc.
\end{enumerate} 
\end{fact}

\begin{proof} First observe that if $\mathcal{F}$ is rich and therefore satisfies \RWin{}, then $\mathcal{F}$ satisfies the weaker \Rrule{} and \Rwinweak{} from \S~\ref{FullFrames}, so Definition \ref{RichPossFrames} and Fact \ref{CharLatComp} together imply that $\mathcal{F}$ is a \textit{principal} possibility frame, and then $\mathcal{F}$ is lattice-complete by Definition \ref{RichPossFrames}. Thus, $\mathcal{F}$ being rich implies that $\mathcal{F}$ is a lattice-complete principal possibility frame satisfying \RWin{}. The converse implication is clear, so $\mathcal{F}$ being rich is equivalent to $\mathcal{F}$ being a lattice-complete principal possibility frame satisfying \RWin{}. 

Then to see that parts \ref{RichFullStrongA}, \ref{RichFullStrongB}, and \ref{RichFullStrongB2} are equivalent, it suffices to recall from Proposition \ref{InterPrinc}.\ref{InterPrinc2} that if $\mathcal{F}$ is a principal possibility frame, then $\mathcal{F}$ satisfies \RWin{} iff it satisfies $R$-tight, and from Fact \ref{TightPrinc} that if $\mathcal{F}$ is a lattice-complete principal possibility frame, then $\mathcal{F}$ satisfies $R$-tight iff it satisfies \Rprinc{}.

To complete the proof, it suffices to observe that parts \ref{RichFullStrongB} and \ref{RichFullStrongC} are equivalent, and that parts \ref{RichFullStrongB2} and \ref{RichFullStrongC2} are equivalent. This follows from Fact \ref{FullTFAE}, which showed that $\mathcal{F}$ is a lattice-complete principal possibility frame iff $\mathcal{F}$ is a full possibility frame in which $\langle S_\bot,\sqsubseteq_\bot\rangle$ is a complete Boolean lattice. \end{proof}

\begin{figure}[h]\begin{center}
\begin{tikzpicture}[semithick]

      \draw (-2,-7.6) rectangle (2,0) ;
      \draw (0,0) circle (2) ;
      \draw (0,0) circle (3) ;
       \draw[dashed] (0,0) circle (4) ;
       \draw (0,0) circle (5) ;

      \node at (0,2.4) {{\textit{principal}}};
       \node at (0,-5.5) {{\textit{$R$-tight}}};
       \node at (0,-5.9) {{(equivalent to \RWin{}}};
        \node at (0,-6.3) {{for full or principal $\mathcal{F}$)}};
      \node at (0,-6.75) {{(equivalent to \Rprinc{}}};
      \node at (0,-7.15) {{for lattice-complete $\mathcal{F}$)}};
       \node at (0,1.4) {{\textit{lattice-}}};
       \node at (0,1) {{\textit{complete}}};
      \node at (0,.5) {{(all full)}};
       \node at (0,3.4) {{\textit{$\sqsubseteq$-tight}}};
        \node at (0,4.4) {{\textit{separative}}};

\fill[pattern=dots] (2,0) arc (0:-180:2);

\end{tikzpicture}
\end{center}
\caption{\textit{rich frames} (dotted region) shown inside other frame classes from \S~\ref{SpecialClasses}. Each label applies to everything inside the smallest circle (or rectangle) that contains the label. The dashed circle reflects the fact that the distinction between \textit{separative} and \textit{$\sqsubseteq$-tight} disappears for full possibility frames (Fact \ref{TightSepDiff}).}\label{FrameClasses}
\end{figure}
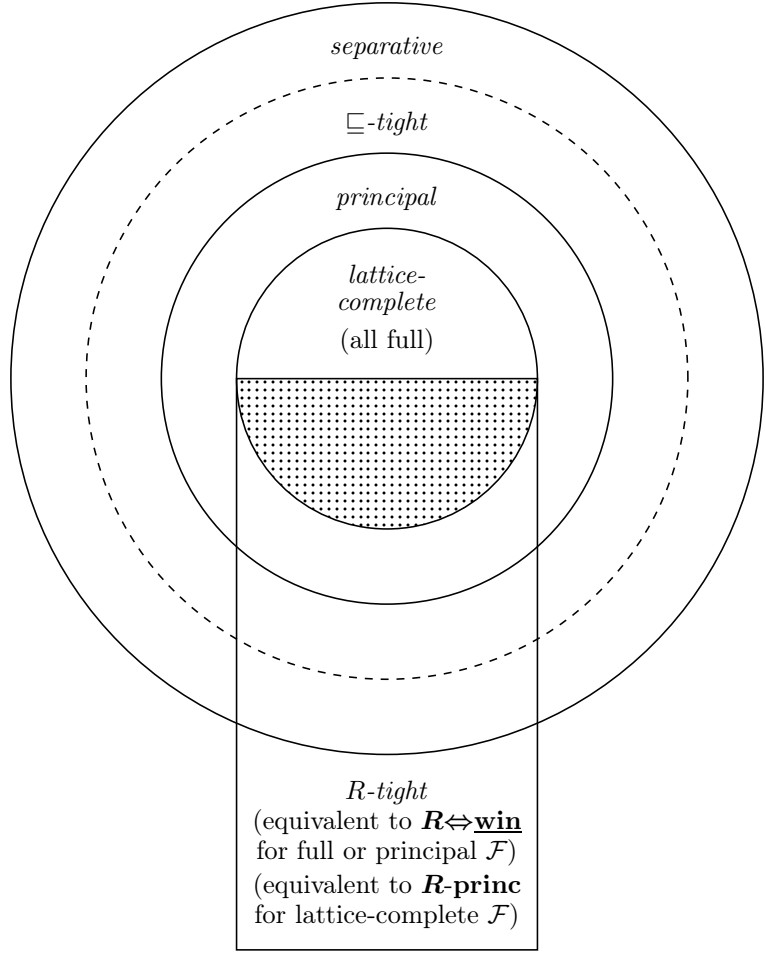

Figure \ref{FrameClasses} shows how the class of rich frames sits inside other frame classes we have discussed in \S~\ref{SpecialClasses}.

The next significant fact about rich frames is that \textit{atomic} (\S~\ref{AtomicSection}) rich frames are exactly the possibility frames that are isomorphic to the \textit{powerset possibilization} of some Kripke frame (Example \ref{PowerPoss}), as shown by the following proposition. For part \ref{PowChar2}, recall the notion of the atom structure $\mathfrak{At}\mathcal{F}$ from Definition \ref{AtSt}. 

\begin{proposition}[Characterization of Powerset Possibilizations]\label{PowChar}$\,$
\begin{enumerate}
\item\label{PowChar1} If $\mathfrak{F}$ is a Kripke frame, then $\mathfrak{F}^\pow$ is an atomic rich frame.
\item\label{PowChar2} If $\mathcal{F}$ is an atomic rich frame, then $\mathcal{F}$ is isomorphic to $(\mathfrak{At}\mathcal{F})^\pow$.
\end{enumerate}
\end{proposition}
 
\begin{proof} Part \ref{PowChar1} is straightforward (see Fact \ref{transform2}). For part \ref{PowChar2}, we showed in Proposition \ref{AtomProp}.\ref{AtomProp3} that for any atomic and separative possibility frame $\mathcal{F}$, which includes any atomic rich frame, the function $h\colon \mathcal{F}\to (\mathfrak{At}\mathcal{F})^\pow$ defined by $h(x)=\{a\in\mathfrak{At}\mathcal{F}\mid a\sqsubseteq^\mathcal{F} x\}$ is a $\sqsubseteq$-strong embedding of $\mathcal{F}$ into $(\mathfrak{At}\mathcal{F})^\pow$. We will now prove that it is a possibility isomorphism as in Definition \ref{PossMorph} when $\mathcal{F}=\langle S,\sqsubseteq, \{R_i\}_{i\in\ind},\adm\rangle$ is an atomic rich frame. First, since we are assuming that $\langle S_\bot,\sqsubseteq_\bot\rangle$ is a complete Boolean lattice, every nonempty $A\subseteq\mathsf{min}\langle S,\sqsubseteq\rangle$ has a least upper bound $x_A$ in $\langle S,\sqsubseteq\rangle$, and $h(x_A)=A$. Thus, $h$ is surjective. Then since $h$ is a surjective $\sqsubseteq$-strong embedding, we have that for all $X\in \adm^\mathcal{F}$, $h[X]\in\adm^{(\mathfrak{At}\mathcal{F})^\pow}$, as required for an isomorphism.

It only remains to show that $xR_i^\mathcal{F} y$ iff $h(x)R_i^{(\mathfrak{At}\mathcal{F})^\pow}h(y)$. By the definition of powerset possibilization, $h(x)R_i^{(\mathfrak{At}\mathcal{F})^\pow}h(y)$ iff $h(y)\subseteq R_i^{\mathfrak{At}\mathcal{F}}[h(x)]$. Since $\mathcal{F}$ satisfies \RWin{}, it satisfies \Rdown{}. Thus, as noted after Definition \ref{AtSt}, the  relation $R_i^{\mathfrak{At}\mathcal{F}}$ in $\mathfrak{At}\mathcal{F}$ is such that for all $a,b\in \mathfrak{At}\mathcal{F}$, $aR_i^{\mathfrak{At}\mathcal{F}}b$ iff $aR_i^\mathcal{F} b$. Thus, $h(y)\subseteq R_i^{\mathfrak{At}\mathcal{F}}[h(x)]$ iff $h(y)\subseteq R_i^\mathcal{F}[h(x)]$. Putting all of this together with the definition of $h$, we have
\begin{equation}h(x)R_i^{(\mathfrak{At}\mathcal{F})^\pow}h(y)\mbox{ iff } \{b\in\mathfrak{At}\mathcal{F}\mid b\sqsubseteq^\mathcal{F} y\}\subseteq R_i^\mathcal{F}[\{a\in\mathfrak{At}\mathcal{F}\mid a\sqsubseteq^\mathcal{F} x\}],\label{NeedRight}\end{equation} 
so we must show that the right side is equivalent to $xR_i^\mathcal{F}y$. If $xR_i^\mathcal{F}y$, then by the left-to-right direction of \RWin{}, for any minimal point $b=y'\sqsubseteq y$, there is a minimal $a=x''\sqsubseteq x$ and a $y''\sqsubseteq b$ such that $aR_i^\mathcal{F}y''$. Then since $b$ is minimal, $y''=b$, so $aR_i^\mathcal{F}b$. Hence the right side of (\ref{NeedRight}) holds. In the other direction, suppose \textit{not} $xR_i^\mathcal{F}y$. If $R_i^\mathcal{F}(x)=\emptyset$, then for any minimal point $a\sqsubseteq x$, we have $R_i^\mathcal{F}(a)=\emptyset$ by \upR{}, which follows from \RWin{}, but there is at least one minimal point $b\sqsubseteq y$, so the right side of (\ref{NeedRight}) does not hold. On the other hand if $R_i^\mathcal{F}(x)\not=\emptyset$, then since $\mathcal{F}$ satisfies \Rprinc{} by Fact \ref{RichFullStrong}, where $R_i^\mathcal{F}(x)=\mathord{\downarrow}f_i(x)$, \textit{not} $xR_i^\mathcal{F}y$ implies \textit{not} $y\sqsubseteq f_i(x)$, which with separativity implies there is a minimal point $b\sqsubseteq y$ such that $b\not\sqsubseteq f_i(x)$, which implies \textit{not} $xR_i^\mathcal{F}b$. Then there can be no minimal point $a\sqsubseteq x$ such that $aR_i^\mathcal{F}b$, for that would imply $xR_i^\mathcal{F}b$ by \upR{}. Thus, $b$ shows that the right side of (\ref{NeedRight}) does not hold.\end{proof}

That rich frames are not required to be atomic is the key to their importance in possibility semantics. The following theorem demonstrates the special relation of rich frames to full possibility frames in general. We prove part \ref{full-to-princ1} in \S~\ref{DualEquiv} as Proposition \ref{UpDownRich}.\ref{UpDownRich1} and Theorem \ref{AlmostBack}.\ref{Inverses0}, and part \ref{full-to-princ2} in \S~\ref{DualEquiv} as Theorem \ref{RefSub}.

\begin{restatable}[From Full Frames to Rich Frames]{theorem}{FulltoPrinc}\label{full-to-princ}$\,$
\begin{enumerate}
\item\label{full-to-princ1} For any full possibility frame $\mathcal{F}$, there is a rich possibility frame $\mathcal{F}'$ and a strict, dense, and robust possibility morphism  $h\colon \mathcal{F}\to\mathcal{F}'$. Thus, by Proposition \ref{Preservation}, for all $\varphi\in\mathcal{L}(\sig,\ind)$, $\mathcal{F}\Vdash \varphi$ iff $\mathcal{F}'\Vdash\varphi$.
\item\label{full-to-princ2}  The category of rich possibility frames with (strict) possibility morphisms is a \textit{reflective subcategory} of the category of full possibility frames with (strict) possibility morphisms.
\end{enumerate}
\end{restatable}
\noindent For a reminder of the definition of reflective subcategory, see the discussion preceding Theorem \ref{RefSub} below.

In \S~\ref{DualEquiv}, we will also show that the category of rich possibility frames with (strict) possibility morphisms is \textit{dually equivalent} to the category of $\mathcal{CV}$-BAOs with complete BAO-homomorphisms (Theorem \ref{CVDuals}). 

\section{Beginnings of Duality Theory}\label{DualityTheory}

The duality theory relating possible world frames and Boolean algebras with operators (BAOs) has been one of the most fruitful areas of mathematical modal logic (see, e.g.,  \citealt{Goldblatt2006} and \citealt[Ch. 5]{Blackburn2001}). As noted in \S~\ref{intro}, two of the fundamental results of this theory are:\footnote{For a detailed study of the following two dualities, as well as a third ``hybrid'' duality, see \citealt{Givant2014}.}
\begin{itemize}
\item \citealt{Thomason1975}: the category of $\mathcal{CAV}$-BAOs (see Definition \ref{BAOclasses}) with complete BAO-homomorphisms is dually equivalent to the category of \textit{full} world frames (Kripke frames) with p-morphisms.
\item \citealt{Goldblatt1974}: the category \textbf{BAO} of BAOs with BAO-homomorphisms is dually equivalent to the category of \textit{descriptive} world frames with p-morphisms. (For a topological version of this duality, see \citealt{Esakia1974}.)
\end{itemize}
We will prove analogues for possibility semantics of both of these results in \S~\ref{DualEquiv} and \S~\ref{Fdes}, respectively:
\begin{itemize}
\item Theorems \ref{CVDuals} and \ref{RefSub}: the category \textbf{$\mathcal{CV}$-BAO} of $\mathcal{CV}$-BAOs (see Definition \ref{BAOclasses}) with complete BAO-homomorphisms is dually equivalent to a reflective subcategory of the category \textbf{FullPoss} of \textit{full} possibility frames with strict possibility morphisms.
\item Theorem \ref{Dual2}: \textbf{BAO} is dually equivalent to the category \textbf{FiltPoss}  of \textit{filter-descriptive} possibility frames (Definition \ref{F-desc}) with p-morphisms.
\end{itemize}
In addition, we will prove results relating classes of BAOs (in particular, $\mathcal{V}$-BAOs and $\mathcal{T}$-BAOs as in Definition \ref{BAOclasses}) to possibility frames that have no obvious analogues on the side of world semantics. 

Figure \ref{CatFig} presents our main categories, functors, and categorical relationships. The functors $(\cdot)^\under$, $(\cdot)_\rela$, and $(\cdot)_\gff$ will be defined in \S~\ref{PossToBAO}, \S~\ref{VtoPossSection}, and \S~\ref{GFPF}, respectively.

 \begin{figure}[h]
\begin{center}
\begin{tabular}{l|c|c|c|c|c}
 && functor && functor & \\ 
 && $\to$ && $\to$ & \\ 
\hline
category: &\textbf{FullPoss} & $(\cdot)^\under$ &\textbf{$\mathcal{CV}$-BAO} & $(\cdot)_\rela$ & \textbf{RichPoss}\\
\hline
objects: & full & & $\mathcal{CV}$-BAOs && rich \\
& possibility frames &&&& possibility frames  \\
\hline
morphisms: & strict  & & complete  & & p-morphisms \\
& possibility morphisms && BAO-homomorphisms && \\
&&& &\\
category: &\textbf{Poss} & $(\cdot)^\under$ & \textbf{BAO} & $(\cdot)_\gff$ &  \textbf{FiltPoss}\\
\hline
objects: & possibility frames & & BAOs && filter-descriptive \\
&&& & & possibility frames\\
\hline
morphisms:&possibility morphisms & & BAO-homomorphisms && p-morphisms \\
\end{tabular}\end{center}

\begin{center}
\begin{tabular}{ccccc}
&&&&\\
\textbf{$\mathcal{CV}$-BAO} & is dually equivalent to & \textbf{RichPoss} & which is a reflective subcategory of & \textbf{FullPoss} \\ &&&&\\
\textbf{BAO} & is dually equivalent to & \textbf{FiltPoss} & which is a reflective subcategory of & \textbf{Poss} 
\end{tabular}
\end{center}

\caption{main categories, functors, and categorical relationships.}\label{CatFig}
\end{figure}

The duality we develop in \S\S~\ref{PossToBAO}-\ref{DualEquiv} between full possibility frames and complete and completely additive BAOs parallels the connection between forcing posets and Boolean valued models in set theory. As Takeuti and Zaring \citeyearpar[p.~1]{Takeuti1973} nicely explain the connection: 
\begin{quote}
One feature [of the theory developed in this book] is that it establishes a relationship between Cohen's method of forcing and Scott-Solovay's method of Boolean valued models.
The key to this theory is found in a rather simple correspondence between partial order structures and complete Boolean algebras\dots. With each partial order structure \textbf{P}, we associate the complete Boolean algebra of regular open sets determined by the order topology on \textbf{P}. With each Boolean algebra \textbf{B}, we associate the partial order structure whose universe is that of \textbf{B} minus the zero element and whose order is the natural order on \textbf{B}.
\end{quote}
This is exactly how we proceed in \S\S~\ref{PossToBAO}-\ref{DualEquiv}, except the key additional step now is to connect the accessibility relations on our partial order structures with the modal operators on our Boolean algebras.

The duality we develop in \S\S~\ref{GFPF}-\ref{Fdes} between arbitrary BAOs and filter-descriptive possibility frames is very similar to the theory of descriptive frames developed by Goldblatt \citeyearpar[\S1.9]{Goldblatt1974} (published in \citealt{Goldblatt1976a,Goldblatt1976b,Goldblatt1993}), except where the theory of descriptive frames relies on assumptions about ultrafilters, we will make do with reasoning about proper filters. As a result, unlike the theory of descriptive frames, the theory of filter-descriptive frames does not require going beyond ZF set theory.  

For a review of BAOs and the semantics for normal modal logics using BAOs, see Appendix \S~\ref{AlgSem}. We will proceed in this section assuming familiarity with that material. For the purposes of possibility semantics, it is more convenient to take the multiplicative operators $\blacksquare_i$ as the primitive operators of the BAO, rather than the additive operators $\blacklozenge_i$, which we define by $\blacklozenge_i x:=\mathord{-}\blacksquare_i \mathord{-}x$. 

For a BAO $\mathbb{A}=\langle A, \meet, -, \top, \{\blacksquare_i\}_{i\in \ind}\rangle$, let $\langle A,\leq\rangle$ be the associated Boolean lattice, with $x\leq y$  iff $x\wedge y= x$. We will focus on the following classes of BAOs, with labels from \citealt{Litak2005}.

\begin{definition}[Classes of BAOs]\label{BAOclasses} Let $\mathbb{A}=\langle A, \meet, -, \top, \{\blacksquare_i\}_{i\in \ind}\rangle$ be a nontrivial BAO. 
\begin{enumerate}
\item\label{C-BAO} $\mathbb{A}$ is a $\mathcal{C}$-BAO iff $\langle A,\leq\rangle$ is a \textit{complete} lattice;
\item $\mathbb{A}$ is an $\mathcal{A}$-BAO iff $\langle A,\leq\rangle$ is an \textit{atomic} lattice as in Definition \ref{AtomicPoset};
\item\label{V-BAO} $\mathbb{A}$ is a $\mathcal{V}$-BAO iff each $\blacksquare_i$ operator is \textit{completely multiplicative}, i.e., for every $i\in \ind$ and $X\subseteq A$, if $\underset{x\in X}{\bigmeet} x$ exists in $\langle A,\leq\rangle$, then $\underset{x\in X}{\bigmeet}\blacksquare_i x$ exists in $\langle A,\leq\rangle$ and $\underset{x\in X}{\bigmeet}\blacksquare_i x=\blacksquare_i\underset{x\in X}{\bigmeet} x$;
\item\label{T-BAO} $\mathbb{A}$ is a $\mathcal{T}$-BAO iff each $\blacksquare_i$ operator admits a \textit{left adjoint}, i.e., there is a function $f_i\colon A\to A$ such that for all $x,y\in A$, $x\leq \blacksquare_i y$ iff $f_i(x)\leq y$. \hfill $\triangleleft$
\end{enumerate}
\end{definition}
Since $\mathbb{A}$ is a BAO, parts \ref{V-BAO}-\ref{T-BAO} could be equivalently stated in terms of $\bigjoin$ and $\blacklozenge_i$ instead of $\bigmeet$ and $\blacksquare_i$. For part \ref{V-BAO}, the equivalent condition is that each $\blacklozenge_i$ is \textit{completely additive}: if $\underset{x\in X}{\bigjoin} x$ exists in $\langle A,\leq\rangle$, then $\underset{x\in X}{\bigjoin}\blacklozenge_i x$ exists in $\langle A,\leq\rangle$ and $\underset{x\in X}{\bigjoin}\blacklozenge_i x=\blacklozenge_i\underset{x\in X}{\bigjoin} x$.  For part \ref{T-BAO}, the equivalent condition is that there is a function $g_i\colon A\to A$ such that for all $x,y\in A$, $\blacklozenge_i x\leq y$ iff $x\leq g_i(y)$. The `$\mathcal{T}$' in `$\mathcal{T}$-BAO' suggests \textit{tense algebras}, which are usually taken to be BAOs with two operators, one of which is the \textit{conjugate} of the other (see, e.g., \citealt[Def. 1.11]{Jonsson1952a}, \citealt[\S4.7]{Jonsson1993}, \citealt[Prop. 8.5]{Venema2007}), i.e., such that $\blacklozenge_1 x\wedge y=\bot$ iff $x\wedge \blacklozenge_2 y=\bot$, in which case $\blacksquare_1$ and $\blacklozenge_2$ form an adjoint pair as above: $x\leq \blacksquare_1 y$ iff $\blacklozenge_2 x\leq y$. We use `$f_i$' instead of some kind of diamond notation as a reminder that this $f_i$ need not be definable by a term from the signature of the BAO, as well as to suggest a connection between $\mathcal{T}$-BAOs and the (quasi-)functional possibility frames of \S~\ref{FuncFrames}, which will be made explicit in Theorems \ref{PtoB}.\ref{PtoB4} and \ref{VtoPossFrames}.\ref{VtoPoss2} below.

The following is a well-known and useful fact.

\begin{fact}[$\mathcal{T}$ and $\mathcal{V}$]\label{TandV} Any $\mathcal{T}$-BAO is a $\mathcal{V}$-BAO, and any $\mathcal{CV}$-BAO is a $\mathcal{CT}$-BAO. Thus, $\mathcal{CV}=\mathcal{CT}$.
\end{fact}
\begin{proof} For the first part, if $\bigmeet X$ exists, then $\blacksquare_i\bigmeet X$ is a lower bound of $\{\blacksquare_i x\mid x\in X\}$. In a $\mathcal{T}$-BAO it is also the greatest: if $z\leq \blacksquare_i x$ for all $x\in X$, then $f_i(z)\leq x$ for all $x\in X$, so $f_i(z)\leq \bigmeet X$ and hence $z\leq \blacksquare_i \bigmeet X$. For the second part, in a $\mathcal{CV}$-BAO $\mathbb{A}$ the left adjoint $f_i$ of $\blacksquare_i$ is given by $f_i(x) =\bigmeet \{ y\in\mathbb{A}\mid x\leq\blacksquare_i y\}$.
\end{proof}

Recall from \S~\ref{intro} that where $\mathcal{X}$ is a class of BAOs, $\mathrm{ML}(\mathcal{X})$ is the set of modal logics \textbf{L} such that \textbf{L} is the logic of some subclass of $\mathcal{X}$, and $\mathcal{ALG}$ is the class of all BAOs, we have:
\[\mathrm{ML}(\mathcal{CAV})\subsetneq \mathrm{ML}(\mathcal{CV})\subsetneq \mathrm{ML}(\mathcal{T})\subsetneq \mathrm{ML}(\mathcal{V})\subsetneq\mathrm{ML}(\mathcal{ALG}).\]

Finally, we fix our terminology for discussing homomorphisms between BAOs.

\begin{definition}[Homomorphisms]\label{Homomorph} Given BAOs $\mathbb{A}=\langle A, \meet, -, \top, \{\blacksquare_i\}_{i\in \ind}\rangle$ and $\mathbb{A}'=\langle A', \meet', -', \top', \{\blacksquare_i'\}_{i\in \ind}\rangle$, a function $h\colon A\to A'$ is a \textit{Boolean algebra homomorphism} iff  for every $x,y\in A$: $h(x\meet y)=h(x)\meet' h(y)$; $h(-x)=-'h(x)$; and $h(\top)=\top'$. For a \textit{BAO-homomorphism}, we additionally require that $h(\blacksquare_ix)=\blacksquare_i'h(x)$. For a \textit{complete} BAO-homomorphism, we  require that if $\underset{j\in J}{\bigmeet} x_j$ exists in $\mathbb{A}$, then $\underset{j\in J}{\bigmeet'} h(x_j)$ exists in $\mathbb{A}'$ and $h(\underset{j\in J}{\bigmeet} x_j)=\underset{j\in J}{\bigmeet'} h(x_j)$. A (complete) \textit{BAO-embedding} is an injective (complete) BAO-homomorphism. A \textit{BAO-isomorphism} is a bijective BAO-homomorphism (and hence complete).\hfill$\triangleleft$\end{definition}

\subsection{From Possibility Frames to BAOs}\label{PossToBAO}

We begin our duality theory with the simpler direction, going from frames to BAOs. As motivation for the following definition, recall the discussion of regular open algebras and subalgebras thereof in Remark \ref{Persp2}.

\begin{definition}[Underlying BAO]\label{RegOpAlg} Given a possibility frame $\mathcal{F}=\langle S, \sqsubseteq , \{R_i\}_{i\in \ind},\adm\rangle$, its \textit{underlying BAO} is the algebra $\mathcal{F}^\under=\langle A,\wedge,-,\top,\{\blacksquare_i\}_{i\in\ind}\rangle$ where:
\begin{enumerate}
\item $A=\adm$;
\item $X\wedge Y=X\cap Y$;
\item $-X=\mathrm{int}(S\setminus X)$;\footnote{Recall from Remark \ref{Persp2} that $\mathrm{int}(X)=\{y\in S\mid \forall x\sqsubseteq y,\, x\in X\}$.}
\item $\top = S$;
\item $\blacksquare_i Y = \{x\in S\mid R_i(x)\subseteq Y\}$. \hfill $\triangleleft$
\end{enumerate}
\end{definition}

In what follows, several notions and results for possibility frames will specialize to well-known notions and results when the possibility frames in question are Kripke frames. Definition \ref{RegOpAlg} is already an example.

\begin{remark}[Full Complex Algebra] If $\mathcal{F}$ is a Kripke frame, regarded as a possibility frame as in Examples \ref{KripkeExample} and \ref{KripkeAgain}, then its underlying BAO as in Definition \ref{RegOpAlg} is its \textit{full complex algebra} (see Appendix \S~\ref{WRM}).
\end{remark}

Let us now verify that the ``underlying BAO'' of a possibility frame is indeed a BAO and moreover that a formula $\varphi$ is satisfiable over a possibility frame according to possibility semantics iff $\varphi$ is satisfiable over its underlying BAO according to algebraic semantics, so the two are semantically equivalent.  As explained in \S~\ref{AlgSem}, we say that $\varphi$ is \textit{satisfiable} over a BAO $\mathbb{A}$ iff there is an algebraic model $\mathbb{M}=\langle \mathbb{A},\theta\rangle$ based on $\mathbb{A}$ such that $\tilde{\theta}(\varphi)\not=\bot$, where $\tilde{\theta}$ is the meaning function extending $\theta$, and $\bot$ is the bottom of $\mathbb{A}$; and $\varphi$ is \textit{valid} over $\mathbb{A}$ iff for all algebraic models $\mathbb{M}=\langle \mathbb{A},\theta\rangle$ based on $\mathbb{A}$, $\tilde{\theta}(\varphi)=\top$, where $\top$ is the top of $\mathbb{A}$. Thus, as usual, $\varphi$ is valid iff $\neg\varphi$ is not satisfiable, so if the same formulas are satisfiable over a possibility frame and over a BAO, then the same formulas are valid over each, so they have the same logic.

\begin{theorem}[From Possibility Frames to BAOs]\label{PtoB} For any possibility frame $\mathcal{F}$ and model $\mathcal{M}=\langle\mathcal{F},\pi\rangle$:
\begin{enumerate}[label=\arabic*.,ref=\arabic*]
\item\label{PtoB1} $\mathcal{F}^\under$ is a BAO;
\item\label{PtoB2} if $\mathcal{F}$ is a full possibility frame, then $\mathcal{F}^\under$ is a $\mathcal{CV}$-BAO;
\item\label{PtoB3} if $\mathcal{F}$ is a principal possibility frame, then $\mathcal{F}^\under$ is a $\mathcal{V}$-BAO;
\item\label{PtoB4} if $\mathcal{F}$ is a quasi-functional principal possibility frame, then $\mathcal{F}^\under$ is a $\mathcal{T}$-BAO;
\item\label{PtoB5} where $\tilde{\pi}$ is the meaning function on $\mathcal{F}^\under$ extending $\pi$, for all $\varphi\in\mathcal{L}(\sig,\ind)$, $\llbracket \varphi\rrbracket^\mathcal{M}=\tilde{\pi}(\varphi)$;
\item\label{PtoB6} for all $\varphi\in\mathcal{L}(\sig,\ind)$, $\varphi$ is satisfiable over $\mathcal{F}$ iff $\varphi$ is satisfiable over $\mathcal{F}^\under$.
\end{enumerate}
\end{theorem}

\begin{proof} For part \ref{PtoB1}, we already explained in Remark \ref{Persp2} that if $\langle A,\wedge,-,\top\rangle$ is defined as in Definition \ref{RegOpAlg}, then it is a Boolean algebra, and it is easy to check that $\blacksquare_i \top=\top$ and $\blacksquare_i (X\meet Y)=\blacksquare_i X\meet \blacksquare_i Y$, so $\mathcal{F}^\under$ is a BAO. For part \ref{PtoB2}, we already explained in Remark \ref{Persp2} that if all regular open sets in $\mathcal{O}(S,\sqsubseteq)$ are admissible propositions, as in a full possibility frame, then $\langle A,\wedge,-,\top\rangle$ is a complete Boolean algebra, and since arbitrary meets are given by intersection in this regular open algebra of the Alexandrov topology $\mathcal{O}(S,\sqsubseteq)$, clearly for any $\mathcal{X}\subseteq A$ we have $\underset{X\in \mathcal{X}}{\bigmeet}\blacksquare_i X=\blacksquare_i\underset{X\in \mathcal{X}}{\bigmeet} X$. Thus, $\mathcal{F}^\under$ is a $\mathcal{CV}$-BAO. 
 
For parts \ref{PtoB3} and \ref{PtoB4}, suppose $\mathcal{F}=\langle S, \sqsubseteq, \{R_i\}_{i\in \ind},\adm\rangle$ is a \textit{principal} possibility frame. Let $\langle S',\sqsubseteq'\rangle$ be the Boolean lattice associated with $\mathcal{F}$, so $S'=S\cup\{\bot'\}$, where $\bot'$ is the minimum of $\langle S',\sqsubseteq'\rangle$. Let $R_i'$ be the extension of $R_i$ to $\langle S',\sqsubseteq'\rangle$ defined by $xR_i'y$ iff $xR_iy$ or $y=\bot'$ (recall Definition \ref{Extending}). Where $\langle S',\meet',-',\top'\rangle$ is the Boolean algebra obtained from $\langle S',\sqsubseteq'\rangle$, consider the structure $\mathcal{F}'=\langle S',\meet',-',\top',\{\latbox_i'\}_{i\in \ind}\rangle$ where $\latbox_i'$ is the operation on $S'$ from Definition \ref{DotOp}, so $\latbox_i' y=\mbox{max}\{x\in S'\mid  R_i'(x)\subseteq\mathord{\downarrow}' y\}$. We claim that $\mathcal{F}'=\langle S',\meet',-',\top',\{\latbox_i'\}_{i\in \ind}\rangle$ is a BAO that is BAO-isomorphic to $\mathcal{F}^\under=\langle \adm,\meet,-,\top,\{\blacksquare_i\}_{i\in\ind}\rangle$ from Definition \ref{RegOpAlg}. It is easy to check that $\latbox_i' \top' =  \top'$ and $\latbox_i' (x\meet' y)=\latbox_i' x\meet' \latbox_i' y$. Since $\mathcal{F}$ is a principal frame, $\adm$ is the set of all principal downsets in $\langle S,\sqsubseteq\rangle$ plus $\emptyset$. Define a function $h\colon S'\to \adm$ by $h(x)=\{x'\in S\mid x'\sqsubseteq' x\}$, noting that $h(\bot')=\emptyset$. It is straightforward to check that $h$ is a BAO-isomorphism, using the fact that $h(\latbox_i' x)=\{y\in S\mid y\sqsubseteq'\latbox_i' x\}=\{y\in S\mid R_i(y)\subseteq \{x'\in S\mid x'\sqsubseteq' x\} \}= \blacksquare_i \{x'\in S\mid x'\sqsubseteq' x\} = \blacksquare_i h(x)$.
 
For part \ref{PtoB3}, it suffices to show that  $\mathcal{F}'$ is a $\mathcal{V}$-BAO. To reduce clutter, we will drop some of the primes. For reductio, suppose that for an ${X}\subseteq S'$, $\bigmeet X$ exists in $\mathcal{F}'$, but $\latbox_i \bigmeet X$ is not the greatest lower bound of $\{\latbox_i x\mid x\in{X}\}$. Since $\mathcal{F}'$ is a BAO, for any $y\in{X}$,  $\bigmeet X\sqsubseteq y$ implies $\latbox_i\bigmeet X\sqsubseteq \latbox_i y$, so  $\latbox_i \bigmeet X$ is a lower bound of $\{\latbox_i x\mid x\in{X}\}$. Hence if $\latbox_i \bigmeet X$ is not the \textit{greatest}, then there is a $u\in S'$ such that $u\sqsubseteq\latbox_i x$ for all $x\in{X}$, but $u\not\sqsubseteq \latbox_i \bigmeet X=\mbox{max}\{y\in S'\mid R_i'(y)\subseteq\mathord{\downarrow} \bigmeet X \}$. Thus, $u\not\in \{y\in S'\mid R_i'(y)\subseteq\mathord{\downarrow} \bigmeet X \}$, so there is a $z\in R_i'(u)$ with $z\not\sqsubseteq \bigmeet X$, so for some $x\in X$, $z\not\sqsubseteq x$ and hence $z\meet - x\not=\bot$. Given $uR_i'z$, $z\wedge - x\not=\bot$, and $u\sqsubseteq\latbox_i x$, by \Rrule{} from Proposition \ref{InterPrinc} there is a $v\in S'$ with $\latbox_i xR_i' v$ and $v\meet -x\not=\bot$. But $\Rlatbox_i x=\mathrm{max}\{y\in S'\mid R_i'(y)\subseteq\mathord{\downarrow}x\}$, so $\latbox_i xR_i' v$ implies $v\sqsubseteq x$, contradicting $v\meet -x\not=\bot$.

For part \ref{PtoB4}, it suffices to show that $\mathcal{F}'$ is a $\mathcal{T}$-BAO. Since $\mathcal{F}$ is quasi-functional, $\mathcal{F}_\bot$ is an extended quasi-functional frame (recall Definition \ref{FuncFramesDef}), so for every $z\in S'$, there is an $f_i'(z)\in S'$ such that $f_i'(z)=\mathrm{max}(R_i'(z))$. We must show that $y\sqsubseteq \latbox_i' x$ iff $f_i'(y)\sqsubseteq x$. The left hand side says that $y\sqsubseteq\mathrm{max}\{z\in S'\mid R_i'(z)\subseteq\mathord{\downarrow} x\}$. If $f_i'(y)\sqsubseteq x$, then since $f_i'(y)=\mathrm{max}(R_i'(y))$, $R_i'(y)\subseteq\mathord{\downarrow} x$, so $y\sqsubseteq \mathrm{max}\{z\in S'\mid R_i'(z)\subseteq\mathord{\downarrow} x\}$. In the other direction, if $y\sqsubseteq \mathrm{max}\{z\in S'\mid R_i'(z)\subseteq\mathord{\downarrow} x\}$, then by Fact \ref{FPrincMon}, $f_i'(y)\sqsubseteq f_i'(\mathrm{max}\{z\in S'\mid R_i'(z)\subseteq\mathord{\downarrow} x\})$, and given $f_i'(z)=\mathrm{max}(R_i'(z))$,  we have $f_i'(\mathrm{max}\{z\in S'\mid R_i'(z)\subseteq\mathord{\downarrow} x\})\sqsubseteq x$, so $f_i'(y)\sqsubseteq x$. 

Part \ref{PtoB5} is by an easy induction on $\varphi$ using Fact \ref{TruthSub}.\ref{TruthSub1}, which is left to the reader.

For part \ref{PtoB6}, we already have from part \ref{PtoB5} that if $\varphi$ is satisfiable over $\mathcal{F}$, then $\varphi$ is satisfiable over $\mathcal{F}^\under$. For the other direction, suppose $\varphi$ is satisfied in some algebraic model $\langle \mathcal{F}^\under,\theta\rangle$. Since $\theta\colon\sig\to\adm$, $\theta$ is an admissible valuation for $\mathcal{F}$, and then by part \ref{PtoB5}, since $\varphi$ is satisfied in $\langle \mathcal{F}^\under,\theta\rangle$,  $\varphi$ is satisfied in $\langle \mathcal{F},\theta\rangle$.\end{proof}

As an aside concerning Theorem \ref{PtoB}.\ref{PtoB2}, recall from Fact \ref{TandV} that every $\mathcal{CV}$-BAO is a $\mathcal{CT}$-BAO, which leads to the following observation.

\begin{observation}[Adjoints in Full Possibility Frames]\label{ResidFull} When $\mathcal{F}=\langle S,\sqsubseteq, \{R_i\}_{i\in\ind},\adm\rangle$ is a full possibility frame, the left adjoint of $\blacksquare_i$ in $\mathcal{F}^\under$ is the $f_i\colon \adm\to\adm$ defined by $f_i(X)=\mathrm{int}(\mathrm{cl}(\mathord{\Downarrow} R_i[X]))$ as in Fact \ref{RefReg}.\ref{RefReg2.5}.
\end{observation}

By Theorem \ref{PtoB}, we can transfer completeness results from possibility semantics to algebraic semantics as follows.

\begin{corollary}[From Frame Completeness to BAO Completeness]\label{TransferPossAlg} For any normal modal logic \textbf{L}, if $\mathbf{L}$ is sound and complete with respect to a class of full possibility frames (resp. principal possibility frames, quasi-functional principal possibility frames), then $\mathbf{L}$ is sound and complete with respect to a class of $\mathcal{CV}$-BAOs (resp. $\mathcal{V}$-BAOs, $\mathcal{T}$-BAOs).\end{corollary}
 
Not only are possibility frames semantically equivalent to their underlying BAOs, but also \textit{morphisms} between possibility frames transform into morphisms between their underlying BAOs in the other direction. Compare this with the analogous standard results for world frames and their underlying BAOs in, e.g., \citealt[Thm.~1.5.9]{Goldblatt1974} or \citealt[Props.~5.51, 5.79]{Blackburn2001}.

\begin{theorem}[From Possibility Morphisms to BAO-Homomorphisms]\label{PossMorphBAOMorph}  
For any possibility frames $\mathcal{F}$ and $\mathcal{F}'$ and function $h\colon\mathcal{F}\to\mathcal{F}'$:
\begin{enumerate}[label=\arabic*.,ref=\arabic*]
\item\label{MtoM1} $h$ is a possibility morphism iff the function $h^\under$ defined for $X'\in \mathcal{F}'^\under$ by $h^\under(X')=h^{-1}[X']$ is a BAO-homomorphism $h^\under\colon \mathcal{F}'^\under\to\mathcal{F}^\under$;
\item\label{MtoM1.5} if $h$ is a possibility morphism and $\mathcal{F}$ and $\mathcal{F}'$ are full possibility frames, then $h^\under$ is a complete BAO-homomorphism;
\item\label{MtoM2} if $h$ is dense as in Definition \ref{PossMorph}.\ref{DenseMorph}, then $h^\under$ is injective;
\item\label{MtoM3} if $h$ is robust as in Definition \ref{PossMorph}.\ref{pre-embed}, then $h^\under$ is surjective;
\item\label{MtoM4} if $f\colon \mathcal{F}\to\mathcal{F}$ is the identity map on $\mathcal{F}$, then $f^\under$ is the identity map on $\mathcal{F}^\under$;
\item\label{MtoM5} for any possibility morphisms $f\colon \mathcal{F}\to\mathcal{G}$ and $g\colon \mathcal{G}\to\mathcal{H}$, $(g\circ f)^\under = f^\under \circ g^\under$.
\end{enumerate}
\end{theorem}
  
\begin{proof} For part \ref{MtoM1}, first observe that the \textit{pull back} condition on possibility morphisms (Definition \ref{PossMorph}.\ref{PossMorph5}) is equivalent to the condition that $h^\under$ is a well-defined function to the domain of $\mathcal{F}^\under$, i.e., $\adm$ in $\mathcal{F}$. Second, clearly we have $h^\under(X'\wedge Y')=h^\under(X')\wedge h^\under(Y')$. Next, by our definitions:
\begin{eqnarray} 
&&h^\under(-{X'})=- h^\under({X}')\nonumber\\
&\Leftrightarrow&  h^\under(\mathrm{int}(S'\setminus {X}'))=\mathrm{int}(S\setminus h^\under({X}'))\nonumber\\
&\Leftrightarrow& h^{-1}[\mathrm{int}(S'\setminus {X}')]=\mathrm{int}(S\setminus h^{-1}[{X}'])\nonumber\\
&\Leftrightarrow& \forall {x}\in{S},\, \mathord{\downarrow}' h(x)\cap X'=\emptyset\mbox{ iff }\mathord{\downarrow}x\cap h^{-1}[X']=\emptyset,\nonumber
\end{eqnarray}
and the `iff' is exactly the \SqMatch{} condition on possibility morphisms. Finally, by our definitions:
\begin{eqnarray}
&&h^\under(\blacksquare_{R_i'}{X}')=\blacksquare_{R_i} h^\under({X}')\nonumber\\
&\Leftrightarrow&  h^\under(\{{y}'\in {S}'\mid R_i'({y}')\subseteq{X}'\})=\{{y}\in {S}\mid R_i({y})\subseteq h^\under({X}')\}\nonumber\\
&\Leftrightarrow& h^{-1}[\{{y}'\in {S}'\mid R_i'({y}')\subseteq{X}'\}]=\{{y}\in {S}\mid R_i({y})\subseteq h^{-1}[{X}']\}\nonumber\\
&\Leftrightarrow& \forall {y}\in{S},\, R_i' (h({y}))\subseteq{X}'\mbox{ iff }R_i({y})\subseteq h^{-1}[{X}'],\nonumber
\end{eqnarray}
and the `iff' is exactly the \RMatch{} condition on possibility morphisms.

For part \ref{MtoM1.5}, if $\mathcal{F}$ and $\mathcal{F}'$ are full possibility frames, then arbitrary meets coincide with intersections, so $h^\under(\underset{j\in J}{\bigmeet}{X}'_j)=\underset{j\in J}{\bigmeet }h^\under({X}'_j)$.

For part \ref{MtoM2}, assume $h$ is dense and that $h^\under (X')=h^\under(Y')$, so $h^{-1}[X']=h^{-1}[Y']$. To show that $X'=Y'$, consider an $x'\in X'$. We claim that for all $y'\sqsubseteq' x'$ there is a $z'\sqsubseteq' y'$ such that $z'\in Y'$, so by \textit{refinability} for $Y'$, we have $x'\in Y'$. Suppose $y'\sqsubseteq' x'$. Since $h$ is dense, there is a $z\in\mathcal{F}$ such that $h(z)\sqsubseteq' y'$.  Since $x'\in X'$ and $h(z)\sqsubseteq' y'\sqsubseteq' x'$, we have $h(z)\in X'$ by \textit{persistence} for $X'$. Then since $h^{-1}[X']=h^{-1}[Y']$, $h(z)\in X'$ implies $h(z)\in Y'$. So where $z'=h(z)$, we have proven the claim needed to apply \textit{refinability} and conclude $x'\in Y'$. Hence $X'\subseteq Y'$, and the other direction is the same, so $h^\under$ is injective.

For part \ref{MtoM3}, suppose $X\in\mathcal{F}^\under$, so $X\in \adm^\mathcal{F}$. Then assuming $h$ is robust, there is an $X'\in \adm^{\mathcal{F}'}$ such that $h[X]=h[S^\mathcal{F}]\cap X'$. We claim that $h^{-1}[X']=X$. Since $h[X]=h[S^\mathcal{F}]\cap X'$, we have $h^{-1}[h[X]]= h^{-1}[h[S^\mathcal{F}]\cap X']=h^{-1}[h[S^\mathcal{F}]]\cap h^{-1}[X']=S^\mathcal{F}\cap h^{-1}[X']=h^{-1}[X']$. Then since $h$ is robust, $h^{-1}[h[X]]=X$, which with the previous sentence gives us $h^{-1}[X']=X$, so $h^\under(X')=X$. This shows that $h^\under$ is surjective.

Parts \ref{MtoM4}-\ref{MtoM5} are easy to check.\end{proof}

Thus, we arrive at the first part of the categorical picture that we prepared for in Remark \ref{Categories}.

\begin{corollary}[The $(\cdot)^\under$ Functor]\label{FuncCor1} The $(\cdot)^\under$ operation given in Definition \ref{RegOpAlg} and Theorem \ref{PossMorphBAOMorph} is a \textit{contravariant functor} from \textbf{Poss} to \textbf{BAO} (resp. \textbf{FullPoss} to \textbf{$\mathcal{CV}$-BAO}).
\end{corollary}

\subsection{From \texorpdfstring{$\mathcal{V}$}{V}-BAOs to Possibility Frames}\label{VtoPossSection}
 
To go from BAOs to possibility frames, we will consider two different routes, one in this section and \S~\ref{DualEquiv}, and the other in \S\S~\ref{GFPF}-\ref{Fdes}. In this section, we show how $\mathcal{V}$-BAOs can be turned into semantically equivalent \textit{principal} possibility frames and $\mathcal{CV}$-BAOs can be turned into semantically equivalent \textit{full} possibility frames. 

In addition, in the case of $\mathcal{CV}$-BAOs, we will establish a categorical connection with full possibility frames. We will show that there is a contravariant functor $(\cdot)_\rela$ from the category \textbf{$\mathcal{CV}$-BAO} of $\mathcal{CV}$-BAOs with complete BAO-homomorphisms to the category \textbf{RichPoss} of rich possibility frames with p-morphisms. Thus, $(\cdot)_\rela$  and $(\cdot)^\under$ from Corollary \ref{FuncCor1} will form a pair of contravariant functors between these categories.

We begin by introducing the following definition of a class of BAOs, inspired by possibility frames, which can be easily transformed into semantically equivalent possibility frames.

\begin{definition}[$\mathcal{R}$-BAOs]\label{plentiful} Given a BAO $\mathbb{A}=\langle A, \meet, -, \top, \{\blacksquare_i\}_{i\in \ind}\rangle$, for each $i\in\ind$, define a relation $R_i\subseteq A\times A$ by:
\begin{itemize}
\item[1.] $x R_i y$ iff for all $y'\in A$, if $\bot\not=y'\leq y$, then $x\meet\blacklozenge_i y'\not=\bot$.
\end{itemize}
Then $\mathbb{A}$ is an $\mathcal{R}$-BAO iff for every $i\in\ind$ and $x,y\in A$:
\begin{itemize}
\item[2.] if $x\meet\blacklozenge_i y\not=\bot$, then there is a $y'\in A$ such that $\bot\not=y'\leq y$ and $x R_i y'$. \hfill $\triangleleft$
\end{itemize}
\end{definition}

It turns out that the class of $\mathcal{R}$-BAOs is exactly the class of $\mathcal{V}$-BAOs, a result which Holliday and Litak \citeyearpar{Litak2015} use to show that not all normal modal logics are sound and complete with respect to a class of $\mathcal{V}$-BAOs.\footnote{Around the same time that Litak and I proved that complete additivity of an operator in a BAO is equivalent to the $\forall\exists\forall$ first-order condition $\mathcal{R}$, in January 2015, Hajnal Andr\'{e}ka, Zal\'{a}n Gyenis, and Istv\'{a}n N\'{e}meti independently proved that complete additivity of an operator on a poset is preserved under ultraproducts. After they learned of our result on the first-orderness of complete additivity in BAOs from Steven Givant, Andr\'{e}ka et al. \citeyearpar{Andreka2016} extended it from BAOs to arbitrary posets.} Here we use it to show that any $\mathcal{V}$-BAO can be transformed into a semantically equivalent principal possibility frame. To establish the result in Lemma \ref{VandH} below, the following straightforward fact is useful.

\begin{fact}[Least Upper Bounds]\label{Useful} For any BAO $\mathbb{A}=\langle A, \meet, -, \top, \{\blacksquare_i\}_{i\in\ind}\rangle$, $x\in A$, and $Y\subseteq A$, the following are equivalent:
\begin{enumerate}
\item $x=\bigvee Y$;
\item $x$ is an upper bound of $Y$, and for any $z\in A$, if $z\meet x\not=\bot$, then for some $y\in Y$, $z\meet y\not=\bot$.
\end{enumerate}
\end{fact}

\begin{lemma}[\citealt{Litak2015}]\label{VandH} $\,$
\begin{enumerate}
\item\label{VtoH} Every $\mathcal{V}$-BAO is an $\mathcal{R}$-BAO;
\item\label{HtoV} Every $\mathcal{R}$-BAO is a $\mathcal{V}$-BAO.
\end{enumerate}
\end{lemma} 
\begin{proof} For part \ref{VtoH}, suppose $\mathbb{A}$ is not an $\mathcal{R}$-BAO, so for some $x,y\in A$, $x\meet \blacklozenge_i y\not=\bot$ but 
\begin{equation}
\forall y'\in A, \mbox{ if }\bot\not=y'\leq y\mbox{ then }\exists y''\in A\colon\bot\not=y''\leq y'\mbox{ and }x\meet\blacklozenge_i y''=\bot.\label{VtoHeq1}
\end{equation} 
Let $Y'$ be the set of all $y'\in A$ such that $\bot\not=y'\leq y$, and let $Y''$ be the set of all $y''\in A$ such that $\bot\not= y''\leq y$ and $x\meet\blacklozenge_i y''=\bot$. Then from (\ref{VtoHeq1}), we have
\begin{equation}
\forall y'\in Y'\, \exists y''\in Y''\colon y''\leq y'.\label{VtoHeq2}
\end{equation}
Now we claim that
\begin{equation}y=\bigvee Y''.\label{CompAd1}\end{equation} 
First, by definition of $Y''$, $y$ is an upper bound of $Y''$. Second, observe that for all $z\in A$, if $z\meet y\not=\bot$, then $z\meet y\in Y'$, which with (\ref{VtoHeq2}) implies that for some $y''\in Y''$, $\bot\not =y''\leq z\meet y$ and hence $z\meet y''\not=\bot$. Thus, we have (\ref{CompAd1}) by Fact \ref{Useful}. Now if 
\begin{equation}\blacklozenge_i y= \bigvee \{\blacklozenge_i y''\mid y''\in Y'' \},\label{CompAd2}\end{equation}
then since $x\meet\blacklozenge_i y\not=\bot$, by Fact \ref{Useful} there is a $y''\in Y''$ such that $x\meet\blacklozenge_i y''\not=\bot$. But by the definition of $Y''$, there is no such $y''$. Thus, although (\ref{CompAd1}) holds, (\ref{CompAd2}) does not, so $\mathbb{A}$ is not a $\mathcal{V}$-BAO.

For part \ref{HtoV}, suppose $\mathbb{A}$ is not a $\mathcal{V}$-BAO, so there is some $X\subseteq A$ such that $\bigvee X$ exists in $\langle  A,\leq\rangle$, but $\blacklozenge_i \bigvee X$ is not the least upper bound of $\{ \blacklozenge _i x\mid x\in X\}$. For every $x\in X$, since $x\leq \bigvee X$, we have $\blacklozenge_i x\leq \blacklozenge_i \bigvee X$, so $\blacklozenge_i \bigvee X$ is an upper bound of $\{ \blacklozenge _i x\mid x\in X\}$. Then by the assumption that it is not a \textit{least} upper bound, Fact \ref{Useful} implies that there is a $z\in A$ such that (i) $z\meet \blacklozenge_i \bigvee X\not=\bot$, but (ii) for all $x\in X$, $z\meet \blacklozenge_i x=\bot$. Now if $\mathbb{A}$ is an $\mathcal{R}$-BAO, then (i) implies that there is a $u\in A$ such that (iii) $\bot\not= u\leq \bigvee X$ and (iv) $z R_i u$, i.e., for all $u'\in A$ such that $\bot\not= u'\leq u$, $z\meet \blacklozenge_i u'\not=\bot$. Given (iii), by Fact \ref{Useful} there is an $x\in X$ such that $u\meet x\not=\bot$. Then where $u'=u\meet x$, we have $\bot\not= u'\leq u$ and hence $z\meet \blacklozenge_i u'\not=\bot$ by (iv). But then since $\blacklozenge_i u'\leq \blacklozenge_i x$, we have $z\meet\blacklozenge_i x\not=\bot$, which contradicts (ii). Thus, $\mathbb{A}$ is not an $\mathcal{R}$-BAO.\end{proof} 

Thomason \citeyearpar{Thomason1975} observed that any $\mathcal{CAV}$-BAO can be turned into an equivalent Kripke frame whose domain is the set of atoms in the BAO and whose accessibility relations $R_i$ are defined by: $aR_ib$ iff $a\meet \blacklozenge_i b\not=\bot$. Note how this definition relates to part \ref{Hposs3} of the following definition for turning $\mathcal{V}$-BAOs into possibility frames.
 
\begin{definition}[Full Frame and Principal Frame of a $\mathcal{V}$-BAO]\label{Hposs} Given a $\mathcal{V}$-BAO $\mathbb{A}=\langle A,\meet,-,\top, \{\blacksquare_i\}_{i\in\ind}\rangle$ and algebraic model $\mathbb{M}=\langle \mathbb{A},\theta\rangle$, we define the \textit{full frame} $\mathbb{A}_\fullV=\langle S,\sqsubseteq, \{R_i\}_{i\in\ind},\adm_\fullV\rangle$ \textit{of} $\mathbb{A}$, the \textit{principal frame} $\mathbb{A}_\rela=\langle S,\sqsubseteq, \{R_i\}_{i\in\ind},\adm_\rela\rangle$ \textit{of} $\mathbb{A}$, and $\mathbb{M}_\rela =\langle S,\sqsubseteq, \{R_i\}_{i\in\ind} ,\pi \rangle$ as follows:
\begin{enumerate}
\item $S=A\setminus \{\bot\}$;
\item $\sqsubseteq$ is the restriction of $\leq$ to $S$;
\item\label{Hposs3} $xR_iy$ iff for all $y'\in A$: if $\bot\not=y'\sqsubseteq y$, then $x\meet \blacklozenge_i y'\not=\bot$;
\item $\adm_\fullV = \mathrm{RO}(S,\sqsubseteq)$ (recall Notation \ref{ROnotation});
\item\label{Hposs4} $\adm_\rela=\{\mathord{\downarrow}x\mid x\in S\}\cup\{\emptyset\}$, where $\mathord{\downarrow}x=\{x'\in S\mid x'\sqsubseteq x\}$;
\item\label{Hposs5} $\pi(p)=\mathord{\downarrow}\theta(p)=\{x\in S\mid x\sqsubseteq \theta(p)\}$. \hfill $\triangleleft$ 
\end{enumerate} 
\end{definition}

The following property of $\mathbb{A}_\fullV$/$\mathbb{A}_\rela$ when $\mathbb{A}$ is a $\mathcal{T}$-BAO will be useful later.

\begin{lemma}[Full/Principal Frame of a $\mathcal{T}$-BAO]\label{RelT} For any $\mathcal{T}$-BAO $\mathbb{A}$, $\mathbb{A}_\fullV$/$\mathbb{A}_\rela$ is such that for all $x,y\in S$, $xR_iy$ iff $y\sqsubseteq f_i(x)$, where $f_i$ is the left adjoint of $\blacksquare_i$ (recall Definition \ref{BAOclasses}.\ref{T-BAO}). 
\end{lemma}
\begin{proof} If we do not have $xR_iy$, then there is a $y'$ such that $\bot\not= y'\sqsubseteq y$ and $x\meet \blacklozenge_i y'=\bot$, so $x\sqsubseteq \blacksquare_i \mathord{-}y'$. Hence $f_i(x)\sqsubseteq \mathord{-}y'$, which with $\bot\not=y'$ implies $y'\not\sqsubseteq f_i(x)$, which with $y'\sqsubseteq y$ implies $y\not\sqsubseteq f_i(x)$.

If $y\not\sqsubseteq f_i(x)$, then $y'=y\meet \mathord{-}f_i(x)\not=\bot$.  Since $f_i$ is the left adjoint of $\blacksquare_i$, we have $x\sqsubseteq \blacksquare_i f_i(x)$ from $f_i(x)\sqsubseteq f_i(x)$, so $x\meet\blacklozenge_i \mathord{-}f_i(x)=\bot$, which implies $x\meet\blacklozenge_i y'=\bot$. Thus, we do not have $xR_iy$.
\end{proof}

Another important fact is that for a $\mathcal{CV}$-BAO, its full frame and principal frame are the same.

\begin{lemma}[Full = Principal Frame of a $\mathcal{CV}$-BAO] For any $\mathcal{CV}$-BAO $\mathbb{A}$, $\mathbb{A}_\fullV=\mathbb{A}_\rela$.
\end{lemma}

\begin{proof} As shown in the proof of Fact \ref{FullTFAE}, if $\langle A,\leq\rangle$ is a \textit{complete} Boolean lattice and $\langle S,\sqsubseteq\rangle$ is the result of deleting the bottom element as in Definition \ref{Hposs}, then the set $\mathrm{RO}(S,\sqsubseteq)$ of regular open sets in the downset topology on $\langle S,\sqsubseteq\rangle$ is exactly the set of principal downsets in $\langle S,\sqsubseteq\rangle$ plus $\emptyset$, so $\adm_\fullV=\adm_\rela$.
\end{proof}

Theorem \ref{VtoPossFrames} records the crucial properties of the $(\cdot)_\rela$ and $(\cdot)_\fullV$ transformations. For part \ref{VtoPoss5} of the theorem, we generalize the notion of satisfiability over a BAO from a single formula to a set of formulas: $\Sigma\subseteq\mathcal{L}(\sig,\ind)$ is satisfiable over a BAO $\mathbb{A}$ iff there is an algebraic model $\mathbb{M}=\langle \mathbb{A},\theta\rangle$ and an $x\in\mathbb{A}$, $x\not=\bot$, such that for all $\sigma\in\Sigma$, $x\leq \tilde{\theta}(\sigma)$. This is not the only generalization that makes sense, but it fits here.

\begin{theorem}[From $\mathcal{V}$-BAOs to Possibility Frames]\label{VtoPossFrames} For any $\mathcal{V}$-BAO $\mathbb{A}$ and algebraic model $\mathbb{M}=\langle \mathbb{A},\theta\rangle$:
\begin{enumerate}[label=\arabic*.,ref=\arabic*]
\item\label{FromGenToPos1} $\mathbb{A}_\rela$ is a tight principal possibility frame, and $\mathbb{A}_\fullV$ is a tight full possibility frame;
\item\label{VtoPoss2} if $\mathbb{A}$ is a $\mathcal{T}$-BAO, then  $\mathbb{A}_\rela$ and $\mathbb{A}_\fullV$ are quasi-functional;
\item\label{VtoPoss3} if $\mathbb{A}$ is a $\mathcal{CV}$-BAO, then $\mathbb{A}_\rela=\mathbb{A}_\fullV$ is a rich possibility frame;
\item\label{VtoPoss1.5} $\mathbb{M}_\rela$ is a possibility model based on $\mathbb{A}_\rela$ and $\mathbb{A}_\fullV$;
\item\label{VtoPoss4} for all $\varphi\in\mathcal{L}(\sig,\ind)$, $\llbracket \varphi\rrbracket^{\mathbb{M}_\rela}=\mathord{\downarrow}\tilde{\theta}(\varphi)$;
\item\label{VtoPoss4.5} for all $\Sigma\subseteq \mathcal{L}(\sig,\ind)$, if $\Sigma$ is satisfiable over $\mathbb{A}$, then $\Sigma$ is satisfiable over $\mathbb{A}_\fullV$;
\item\label{VtoPoss5} for all $\Sigma\subseteq \mathcal{L}(\sig,\ind)$, $\Sigma$ is satisfiable over $\mathbb{A}$ iff $\Sigma$ is satisfiable over $\mathbb{A}_\rela$.
\end{enumerate}
\end{theorem}
 
\begin{proof} For part \ref{FromGenToPos1}, we first show that $\mathbb{A}_\rela=\langle S,\sqsubseteq, \{R_i\}_{i\in\ind},\adm \rangle$ is a principal possibility frame. By Fact \ref{PrincEquiv}, it suffices to show: that $\adm$ is the set of all principal downsets in $\langle S,\sqsubseteq\rangle$ plus $\emptyset$, which is immediate from Definition \ref{Hposs}; that $\langle S_\bot,\sqsubseteq_\bot\rangle$ is a Boolean lattice, which is also immediate, since $\langle S_\bot,\sqsubseteq_\bot\rangle$ is isomorphic to the underlying lattice $\langle A,\leq\rangle$ of $\mathbb{A}$, which we will now identify with $\langle S_\bot,\sqsubseteq_\bot\rangle$; and that the operation $\Rlatbox_i$ defined on $\langle A,\leq\rangle$ using $R_i$ as in Definition \ref{DotOp} is a total operation. For the last property, we claim that for all $y\in A$, $\blacksquare_i^\mathbb{A} y=\Rlatbox_iy=\mbox{max}\{x\in A\mid  R_i(x)\subseteq\mathord{\downarrow} y\}$, so $\Rlatbox_i$ is total since $\blacksquare_i^\mathbb{A}$ is total. (To reduce clutter, we will drop the superscript for $\mathbb{A}$, but remember that $\blacksquare_i$ is the operator in $\mathbb{A}$, not the defined operation $\blacksquare_i^{\mathbb{A}_\rela}$ on $\adm$ in the frame $\mathbb{A}_\rela$.) 
First note that $\{x\in A\mid R_i(x)\subseteq\mathord{\downarrow}y\}$ is a downset in $\langle A,\leq\rangle$, as required by Definition \ref{DotOp}, since $x'\leq x$ implies $R_i(x')\subseteq R_i(x)$ by Definition \ref{Hposs}.\ref{Hposs3}. Now suppose $x\not\leq \blacksquare_i y$, so $x\meet \blacklozenge_i \mathord{-} y\not=\bot$. Then since $\mathbb{A}$ is a $\mathcal{V}$-BAO and hence an $\mathcal{R}$-BAO by Lemma \ref{VandH}, it follows that there is a $z\leq -y$ such that $z\not=\bot$ and $xR_i z$, which implies $R_i(x)\not\subseteq\mathord{\downarrow} y$. This shows that $\blacksquare_i y$ is an upper bound of the set $\{x\in A\mid  R_i(x)\subseteq\mathord{\downarrow} y\}$. To show that it is the maximum, we show that $\blacksquare_i y$ belongs to the set. Suppose $\blacksquare_i y R_i z$, so by the definition of $R_i$ in $\mathbb{A}_\rela$, for all $z'\leq z$ with $z'\not=\bot$, we have $\blacksquare_i y\wedge \blacklozenge_i z'\not=\bot$. Now if $z\not\leq y$, then $\bot\not=z\meet -y \leq z$, so by the previous step, $\blacksquare_i y\meet \blacklozenge_i(z\meet -y)\not=\bot$, which is a contradiction given the properties of $\blacksquare_i$ (Definition~\ref{BAOs}). Thus, $R_i(\blacksquare_i y)\subseteq \mathord{\downarrow} y$, as desired.

To show that $\mathbb{A}_\rela$ is $R$-tight, let $\boxtimes_i$ be the modal operator on $\adm$ in $\mathbb{A}_\rela$, in order to distinguish it from $\blacksquare_i$/$\boxdot_i$. We need to show that if $\forall Z\in\adm$,  $x\in\boxtimes_i Z\Rightarrow y\in Z$, then $xR_iy$. Suppose that we do not have $xR_iy$, so there is a $y'\in A$ with $\bot\not=y'\leq y$ and $x\meet \blacklozenge_i y'=\bot$. Then we claim that $x\in \boxtimes_i \mathord{\downarrow} \mathord{-}y'$, i.e., $R_i(x)\subseteq \mathord{\downarrow} \mathord{-}y'$. If $xR_iv$, then for all $v'\in A$ such that $\bot\not=v'\leq v$, we have $x\meet \blacklozenge_i v'\not=\bot$, which with $x\meet \blacklozenge_i y'=\bot$ implies $v'\not\leq y'$. Since this holds for all $v'\not=\bot$ with $v'\leq v$, we have $v\leq \mathord{-}y'$. Thus, $x\in \boxtimes_i \mathord{\downarrow} \mathord{-}y'$.  But since $\bot\not=y'\leq y$, we have $y\not\leq -y'$, so $y\not\in\mathord{\downarrow} \mathord{-}y'$. Thus, we have a $Z$ for which $x\in\boxtimes_i Z$ but $y\not\in Z$. Hence $\mathbb{A}_\rela$ is $R$-tight. Finally, since $\mathbb{A}_\rela$ is a principal possibility frame, it is $\sqsubseteq$-tight by Fact \ref{PrincPossSep}.

Next we show that $\mathbb{A}_\fullV$ is a full possibility frame that is tight. Since $\mathbb{A}_\rela$ is $R$-tight, it is \textit{strong} by Proposition \ref{InterPrinc}.\ref{InterPrinc2}, i.e., $\sqsubseteq$ and $R_i$ satisfy \RWin{}. Then since $\mathbb{A}_\fullV=((\mathbb{A}_\rela)_\found)^\full$ (Definition \ref{Foundation}), it follows by Corollary \ref{Fullinterplay}.\ref{Fullinterplay2} that $\mathbb{A}_\fullV$ is indeed a full possibility frame. Then since $\mathbb{A}_\fullV$ is also strong, it is $R$-tight by Lemma \ref{TightStrong}.\ref{TightStrong1.75}. Finally, since $\mathbb{A}_\fullV$ is a full frame and is separative by Fact \ref{Boolean}, it is $\sqsubseteq$-tight by Fact \ref{TightSepDiff}.\ref{TightSepDiff1}.

For part \ref{VtoPoss2}, by Lemma \ref{RelT}, if $\mathbb{A}$ is a $\mathcal{T}$-BAO, then $\mathbb{A}_\fullV$ and $\mathbb{A}_\rela$ satisfy \Rprinc{} as in Fact \ref{RprincFact} with $R_i(x)=\mathord{\downarrow}f_i(x)$, so they are quasi-functional by Definition \ref{FuncFramesDef}.

For part \ref{VtoPoss3}, any $\mathcal{CV}$-BAO $\mathbb{A}$ is a $\mathcal{T}$-BAO by Fact \ref{TandV}, so $\mathbb{A}_\rela$ satisfies \Rprinc{} as above. By Fact \ref{RichFullStrong}, being a rich frame is equivalent to being a full frame that satisfies \Rprinc{} and for which $\langle S_\bot,\sqsubseteq_\bot\rangle$ is a complete Boolean lattice. It follows that if $\mathbb{A}$ is a $\mathcal{CV}$-BAO, then $\mathbb{A}_\fullV=\mathbb{A}_\rela$ is a rich frame.

For part \ref{VtoPoss1.5}, that $\pi(p)\in \adm_\rela$ is immediate from Definition \ref{Hposs}.\ref{Hposs4}-\ref{Hposs5}, so $\mathbb{M}_\rela$ is based on $\mathbb{A}_\rela$, and $\adm_\rela\subseteq \adm_\fullV$, so $\mathbb{M}_\rela$ is based on $\mathbb{A}_\fullV$ as well.

For part \ref{VtoPoss4}, $\llbracket \varphi\rrbracket^{\mathbb{M}_\rela}=\mathord{\downarrow}\tilde{\theta}(\varphi)$ can be rewritten as: for all $x\in S$, $\mathbb{M}_\rela,x\Vdash \varphi$ iff $x\sqsubseteq \tilde{\theta}(\varphi)$. The proof is by induction on $\varphi$. The atomic case is by definition of $\mathbb{M}_\rela$, and the $\wedge$ case is routine. For the $\neg$ case, $\tilde{\theta}(\neg\varphi)=-\tilde{\theta}(\varphi)$, so we show that $\mathbb{M}_\rela,x\Vdash\neg \varphi$ iff $x\sqsubseteq -\tilde{\theta}(\varphi)$. Now $x\sqsubseteq -\tilde{\theta}(\varphi)$ iff for all $x'\sqsubseteq x$ (so $x'\not=\inc$), $x'\not\sqsubseteq \tilde{\theta}(\varphi)$. Then since $x'\not\sqsubseteq \tilde{\theta}(\varphi)$ is equivalent to $\mathbb{M}_\rela,x'\nVdash \varphi$ by the inductive hypothesis, the right side of the previous `iff' is equivalent to $\mathbb{M}_\rela,x\Vdash\neg\varphi$. For the $\Box_i$ case, $\tilde{\theta}(\Box_i\varphi)=\blacksquare_i\tilde{\theta}(\varphi)$, so we show that $\mathbb{M}_\rela,x\Vdash \Box_i\varphi$ iff $x\sqsubseteq \blacksquare_i\tilde{\theta}(\varphi)$. From left to right, suppose $x\not\sqsubseteq \blacksquare_i\tilde{\theta}(\varphi)$, so $x\meet \blacklozenge_i \mathord{-}\tilde{\theta}(\varphi)\not=\bot$. Then since $\mathbb{A}$ is a $\mathcal{V}$-BAO and hence an $\mathcal{R}$-BAO by Lemma \ref{VandH}, there is a $y\in A$ such that $\bot\not=y\sqsubseteq -\tilde{\theta}(\varphi)$ and $xR_i y$. Hence $y\not\sqsubseteq\tilde{\theta}(\varphi)$, which with the inductive hypothesis implies $\mathbb{M}_\rela,y\nVdash \varphi$, which with $xR_iy$ implies $\mathbb{M}_\rela,x\nVdash \Box_i\varphi$. In the other direction, suppose $\mathbb{M}_\rela,x\nVdash \Box_i\varphi$, so there is a $y\in A$ such that $xR_iy$ and $\mathbb{M}_\rela,y\nVdash \varphi$. Then $y\not\sqsubseteq\tilde{\theta}(\varphi)$ by the inductive hypothesis, so $y\meet -\tilde{\theta}(\varphi)\not=\bot$. It follows by the definition of $R_i$ in $\mathbb{M}_\rela$ that $x\meet \blacklozenge_i (y\meet -\tilde{\theta}(\varphi))\not=\inc$, which by the properties of $\blacksquare_i$ (Definition \ref{BAOs}) implies $x\meet \blacklozenge_i \mathord{-}\tilde{\theta}(\varphi)\not=\bot$ and hence $x\not\sqsubseteq \blacksquare_i\tilde{\theta}(\varphi)$.

Part \ref{VtoPoss4.5} is immediate from parts \ref{VtoPoss1.5} and \ref{VtoPoss4}.

For part \ref{VtoPoss5}, the left-to-right direction is immediate from part \ref{VtoPoss4}. For the right-to-left direction, given a possibility model $\langle \mathbb{A}_\rela,\pi\rangle$, define $\pi_{-\rela}\colon \sig\to A$ (where $A$ is the domain of $\mathbb{A}$) such that $\pi_{-\rela}(p)=\bot$ if $\pi(p)=\emptyset$, and otherwise $\pi_{-\rela}(p)=x$ for the $x$ such that $\pi(p)=\mathord{\downarrow}x$, which is guaranteed to exist since $\mathbb{A}_\rela$ is a \textit{principal} possibility frame for which $\pi$ is an admissible valuation. Then $\langle\mathbb{A},\pi_{-\rela}\rangle$ is an algebraic model such that $\langle\mathbb{A}_\rela,(\pi_{-\rela})_\rela\rangle = \langle\mathbb{A}_\rela,\pi\rangle$, so if $\Sigma$ is satisfied at a point $x$ in  $ \langle\mathbb{A}_\rela,\pi\rangle$ and hence in $\langle\mathbb{A}_\rela,(\pi_{-\rela})_\rela\rangle$, then by part \ref{VtoPoss4}, $\Sigma$ is satisfied at $x$ in $\langle\mathbb{A},\pi_{-\rela}\rangle$ as well.
\end{proof}

Note the contrast between parts \ref{VtoPoss4.5} and \ref{VtoPoss5} of Theorem \ref{VtoPossFrames}: we cannot always turn a possibility model based on $\mathbb{A}_\fullV$ into a semantically equivalent algebraic model based on $\mathbb{A}$. We will return to this point in \S~\ref{MacNeille}.

While Theorem \ref{VtoPossFrames} applies to all $\mathcal{V}$-BAOs, if we focus in on $\mathcal{CV}$-BAOs we can go further: not only are $\mathcal{CV}$-BAOs semantically equivalent to full possibility frames, but also morphisms between $\mathcal{CV}$-BAOs transform into morphisms between their associated possibility frames in the other direction. Thomason \citeyearpar{Thomason1975} observed that given $\mathcal{CAV}$-BAOs $\mathbb{A}$ and $\mathbb{A}'$ and a complete BAO-homomorphism $h\colon \mathbb{A}'\to\mathbb{A}$, for every atom $a$ in $\mathbb{A}$, there is a unique atom $a'$ in $\mathbb{A}'$ such that $a\leq h(a')$, and the function that sends each atom $a$ in $\mathbb{A}$ to its $a'$  with $a\leq h(a')$ is a p-morphism from the atom structure of $\mathbb{A}$ to that of $\mathbb{A}'$. Since Thomason's construction depends on atoms, we cannot use it for $\mathcal{CV}$-BAOs in general. The construction of $h_\rela$ we use for Theorem \ref{BAOMorphPossMorph} has been used by Bezhanishvili \citeyearpar[Thm.~30]{Bezhanishvili1999} to establish a duality between, on the one hand, certain frames for intuitionistic modal logic together with p-morphism-like maps, and on the other hand, \textit{complete} and \textit{completely join-prime generated} Heyting algebras with operators (HAOs) together with complete HAO-homomorphisms. Although this construction does not depend on atoms, it does depend on lattice-completeness for the arbitrary meets in the definition of $h_\rela$, so we cannot use it for $\mathcal{V}$-BAOs in general.

\begin{theorem}[From BAO-Homomorphisms to Possibility Morphisms]\label{BAOMorphPossMorph}
For any $\mathcal{CV}$-BAOs $\mathbb{A}$ and $\mathbb{A}'$ and complete BAO-homomorphism $h\colon\mathbb{A}'\to\mathbb{A}$, define $h_\rela\colon \mathbb{A}_\rela\to\mathbb{A}'_\rela$ by $h_\rela(x)=\bigmeet' \{x'\in \mathbb{A}'\mid x\leq h(x') \}$. Then: 
\begin{enumerate}
\item\label{BAOMorphPossMorph1} $h_\rela$ is a p-morphism;
\item\label{BAOMorphPossMorph2} if $h$ is surjective, then $h_\rela$ is a strong embedding;
\item\label{BAOMorphPossMorph3} if $h$ is injective, then $h_\rela$ is surjective;
\item\label{BAOMorphPossMorph4} if $f\colon \mathbb{A}\to\mathbb{A}$ is the identity map on $\mathbb{A}$, then $f_\rela$ is the identity map on $\mathbb{A}_\rela$;
\item\label{BAOMorphPossMorph5} for any $\mathcal{CV}$-BAOs $\mathbb{A}$, $\mathbb{B}$, and $\mathbb{C}$ and complete BAO-homomorphisms $f\colon \mathbb{A}\to\mathbb{B}$ and $g\colon \mathbb{B}\to\mathbb{C}$, $(g\circ f)_\rela = f_\rela \circ g_\rela$.
\end{enumerate}
\end{theorem}
 
\begin{proof} We begin by verifying that $h_\rela$ is indeed a function to the domain of $\mathbb{A}_\rela'$, which is that of $\mathbb{A}'$ minus $\bot'$. This requires that for $x\in\mathbb{A}_\rela$, so $x\not=\bot$, we have $h_\rela(x)\not=\bot'$. Indeed, if $h_\rela(x)=\bigmeet' \{x' \mid x\leq h(x') \}=\bot'$, then since $h$ is complete, $\bigmeet \{h(x')\mid x\leq h(x')\}=\bot$, which implies $x\leq\bot$, contradicting $x\not=\bot$. Since $\bot$ is not in the domain of $h_\rela$ and $\bot'$ is not in its range, we may write $h_\rela(x)=\bigmeet' \{x'\in \mathbb{A}'_\rela\mid x\sqsubseteq h(x') \}$. Here $\sqsubseteq$ and $\sqsubseteq'$ are the refinement relations in $\mathbb{A}_\rela$ and $\mathbb{A}_\rela'$, respectively, as in Definition \ref{Hposs}.

Now we make several observations about $h_\rela$, assuming that any element to which $h_\rela$ is applied is not $\bot$:
\begin{itemize}
\item[(a)] $y\sqsubseteq x$ implies $h_\rela(y)\sqsubseteq' h_\rela(x)$. For if $y\sqsubseteq x$, then $x\sqsubseteq h(z')$ implies $y\sqsubseteq h(z')$, so $\{z'\mid x\sqsubseteq h(z')\}\subseteq\{z'\mid y\sqsubseteq h(z')\}$, so $\bigmeet' \{z'\mid y\sqsubseteq h(z')\}\sqsubseteq' \bigmeet' \{z'\mid x\sqsubseteq h(z')\}$, which means $h_\rela(y)\sqsubseteq' h_\rela(x)$. 
\item[(b)] $h_\rela(h(y'))\sqsubseteq' y'$. By definition of $h_\rela$, $h_\rela(h(y'))=\bigmeet'\{z'\mid h(y')\sqsubseteq h(z')\}$, and $y'\in \{z'\mid h(y')\sqsubseteq h(z')\}$.
\item[(c)] $x\sqsubseteq h(h_\rela(x))$. By definition of $h_\rela$, $h(h_\rela(x))=h(\bigmeet'\{z'\mid x\sqsubseteq h(z')\})$; since $h$ is a complete homomorphism,  $h(\bigmeet'\{z'\mid x\sqsubseteq h(z')\})=\bigmeet\{h(z')\mid x\sqsubseteq h(z')\}$; and $x\sqsubseteq \bigmeet\{h(z')\mid x\sqsubseteq h(z')\}$.
\item[(d)] $h_\rela(x)\sqsubseteq' y'$ iff $x\sqsubseteq h(y')$. Since $h_\rela(x)=\bigmeet'\{z'\mid x\sqsubseteq h(z')\}$, clearly $x\sqsubseteq h(y')$ implies $h_\rela(x)\sqsubseteq' y'$.  In the other direction, $h_\rela(x)\sqsubseteq' y'$ implies $h(h_\rela(x))\sqsubseteq h(y')$, which with (c) implies $x\sqsubseteq h(y')$.
\item[(e)] $h_\rela(f_i(x))=f'_i(h_\rela(x))$, where $f_i$ and $f_i'$ are the left adjoints of $\blacksquare_i$ and $\blacksquare_i'$, respectively, which exist since $\mathbb{A}$ and $\mathbb{A}'$ are $\mathcal{CV}$-BAOs and hence $\mathcal{CT}$-BAOs (Fact \ref{TandV}), and we assume $f_i(x)\neq \bot$. Note that $f_i(x)=\bot$ iff $f'_i(h_\rela(x))=\bot'$ by adjointness and the argument for (d), using $h(\blacksquare_i'\bot')=\blacksquare_i\bot$. Given $f_i(x)\neq \bot$, by (d), $h_\rela(f_i(x))\sqsubseteq' y'$ iff $f_i(x)\sqsubseteq h(y')$. Since $f_i$ is the left adjoint of $\blacksquare_i$, the right side of the `iff' is equivalent to  $x\sqsubseteq \blacksquare_i h(y')=h(\blacksquare_i' y')$, which by (d) is equivalent to $h_\rela(x)\sqsubseteq ' \blacksquare_i' y'$, which by the fact that $f_i'$ is the left adjoint of $\blacksquare_i'$ is equivalent to $f_i'(h_\rela(x))\sqsubseteq' y'$. Thus, for any $y'\in\mathbb{A}_\rela'$, $h_\rela(f_i(x))\sqsubseteq' y'$ iff $f_i'(h_\rela(x))\sqsubseteq' y'$, so $h_\rela(f_i(x))=f'_i(h_\rela(x))$. 
\end{itemize}
Observation (a) shows that $h_\rela$ satisfies \SqForth{}. Next, we show:
\begin{itemize}
\item if $y'\sqsubseteq' h_\rela( x )$, then $\exists y$: $y\sqsubseteq x $ and $h_\rela(y)= y'$ (\pSqBack{});
\item if $ x R_i y$, then $h_\rela(x)R_i' h_\rela(y)$ (\RForth{}); 
\item if $h_\rela( x )R_i' y'$, then $\exists y$: $ x R_i y$ and $h_\rela(y)=y'$ (\pRBack{}).
\end{itemize}
Since $\mathbb{A}_\rela$ is a full possibility frame by Theorem \ref{VtoPossFrames}.\ref{VtoPoss3}, the \PullBack{} property of $h_\rela$ follows from \SqForth{} and \SqBack{} by Fact \ref{PullFact}.

For \pSqBack{}, suppose $y'\sqsubseteq' h_\rela( x )$. First, we claim that $h(y')\meet x\not=\bot$. If not, then $x\sqsubseteq -h(y')=h(-y')$, which by (a) and (b) implies $h_\rela(x)\sqsubseteq' h_\rela(h(-y'))\sqsubseteq'-y'$. But then $y'\sqsubseteq' h_\rela(x)\sqsubseteq' -y'$, so $y'=\bot'$, contradicting $y'\in\mathbb{A}_\rela'$. Second, we claim that $h_\rela (h(y')\meet x)\sqsubseteq ' y'$. By (a), $h_\rela (h(y')\meet x)\sqsubseteq' h_\rela(h(y'))$, and by (b), $h_\rela(h(y'))\sqsubseteq' y'$, so $h_\rela (h(y')\meet x)\sqsubseteq ' y'$.  Third, we claim that $y'\sqsubseteq' h_\rela (h(y')\meet x)$. By definition of $h_\rela$, it suffices to show that $y'\sqsubseteq' x'$ for all $x'\in \mathbb{A}'_\rela$ such that $h(y')\wedge x\sqsubseteq h(x')$. From $h(y')\wedge x\sqsubseteq h(x')$, we have $x\sqsubseteq -h(y')\vee h(x')$ and hence $x\sqsubseteq h(-y'\vee x')$, which by (a) and (b) implies $h_\rela (x)\sqsubseteq' h_\rela (h(-y'\vee x'))\sqsubseteq' -y'\vee x'$. Then since $y'\sqsubseteq' h_\rela( x )$, it follows that $y'\sqsubseteq' x'$, as desired. Thus, we have shown that $h_\rela (h(y')\meet x)= y'$, so setting $y=h(y')\meet x$ completes the proof of \pSqBack{}.

For \RForth{}, if $xR_i y$, then $y\sqsubseteq f_i(x)$ by Lemma \ref{RelT}, so $h_\rela(y)\sqsubseteq' h_\rela(f_i(x))$ by (a), so $h_\rela(y)\sqsubseteq' f_i'(h_\rela(x))$ by (e), so $h_\rela(x)R_i' h_\rela(y)$ by Lemma \ref{RelT}.

For \pRBack{}, suppose $h_\rela(x)R_i'y'$. First, we claim that $h(y')\wedge f_i(x)\neq\bot$. If not, then $f_i(x)\sqsubseteq - h(y')$, which by adjointness implies $x\sqsubseteq \blacksquare_i \mathord{-}h(y')$ and hence $x\sqsubseteq h(\blacksquare_i ' \mathord{-}y')$, which by (a) and (b) implies $h_\rela (x)\sqsubseteq' h_\rela (h(\blacksquare_i ' \mathord{-}y')) \sqsubseteq' \blacksquare_i'\mathord{-}y'$, which contradicts $h_\rela(x)R_i'y'$.  Second, we have $h_\rela (h(y')\wedge f_i(x))\sqsubseteq' y'$ by (d). Third, we claim that $y'\sqsubseteq' h_\rela(h(y')\wedge f_i(x))$. Since $h_\rela(x)R_i'y'$, we have $y'\sqsubseteq' f_i'(h_\rela(x))$ by Lemma \ref{RelT}, so $y'\sqsubseteq' h_\rela(f_i(x))$ by (e). The rest of the proof is analogous to that used to show $y'\sqsubseteq' h_\rela(h(y')\meet x)$ above for \pSqBack{}, substituting $f_i(x)$ for $x$.  Thus, we have shown that $h_\rela (h(y')\wedge f_i(x))=y'$, so setting $y=h(y')\wedge f_i(x)$ completes the proof of \pRBack{}.

For part \ref{BAOMorphPossMorph2}, assuming $h$ is surjective, we must show that (i) $h_\rela(y)\sqsubseteq' h_\rela(x)$ implies $y\sqsubseteq x$ (so $h_\rela$ is injective), (ii) $h_\rela(x)R_i' h_\rela(y)$ implies $xR_iy$, and (iii)  for all $X\in \adm$ there is an $X'\in \adm'$ such that $h_\rela[X]=h_\rela[S]\cap X'$, where $S$ is the domain of $\mathbb{A}_\rela$. For (i), if $h_\rela(y)\sqsubseteq' h_\rela(x)$, then
\[\begin{array}{llll}
&&\bigmeet' \{y'\in \mathbb{A}'\mid y\leq h(y') \}\leq'\bigmeet' \{x'\in \mathbb{A}'\mid x\leq h(x') \}&\mbox{by definition of }h_\rela\mbox{ and }\sqsubseteq' \\
&\Rightarrow& h(\bigmeet' \{y'\in \mathbb{A}'\mid y\leq h(y') \})\leq h(\bigmeet' \{x'\in \mathbb{A}'\mid x\leq h(x') \})&\mbox{since }h\mbox{ is a homomorphism} \\
&\Rightarrow& \bigmeet \{h(y')\mid  y\leq h(y') \}\leq \bigmeet \{h(x')\mid x\leq h(x') \}&\mbox{since }h\mbox{ is complete}\\
&\Rightarrow& y\leq x&\mbox{since }h\mbox{ is surjective,}
\end{array}\] 
so $y\sqsubseteq x$ since $y\in\mathbb{A}_\rela$. For (ii), if $h_\rela(x)R_i' h_\rela(y)$, then by Lemma \ref{RelT}, $h_\rela(y)\sqsubseteq' f_i'(h_\rela(x))$, which with (e) implies $h_\rela(y)\sqsubseteq' h_\rela (f_i(x))$, which with (i) implies $y\sqsubseteq f_i(x)$, which with Lemma \ref{RelT} implies $xR_iy$. Finally, for (iii), suppose $X\in\adm$, so by Definition \ref{Hposs}.\ref{Hposs4}, either $X=\emptyset$ or $X=\mathord{\downarrow}x$ for some non-minimum element $x$ of $\mathbb{A}$. If $X=\emptyset$, we can take $X'=\emptyset$. So suppose $X=\mathord{\downarrow}x$. Let $X'=\mathord{\downarrow}' h_\rela(x)$. Then $X'\in \adm'$ by Definition \ref{Hposs}.\ref{Hposs4}. It only remains to show that $h_\rela[\mathord{\downarrow} x]=h_\rela[S]\cap \mathord{\downarrow}'h_\rela(x)$. Suppose $y'\in h_\rela[\mathord{\downarrow} x]$, so there is a $y\sqsubseteq x$ such that $h_\rela(y)=y'$. Then by \SqForth{}, $y'=h_\rela(y)\sqsubseteq' h_\rela(x)$, so $y'\in  h_\rela[S]\cap \mathord{\downarrow}'h_\rela(x)$. In the other direction, suppose $y'\in h_\rela[S]\cap \mathord{\downarrow}'h_\rela(x)$, so there is a $y\in S$ such that $y'=h_\rela(y)\sqsubseteq' h_\rela(x)$. Then by (i) above we have $y\sqsubseteq x$, which with $y'=h_\rela(y)$ implies $y'\in h_\rela[\mathord{\downarrow} x]$. This completes the proof of part \ref{BAOMorphPossMorph2}. 

For part \ref{BAOMorphPossMorph3}, assuming $h$ is injective and $y'\in\mathbb{A}_\rela'$, we claim that $h_\rela(h(y'))=y'$. As in (b) above, by definition of $h_\rela$, $h_\rela(h(y'))=\bigmeet'\{z'\mid h(y')\sqsubseteq h(z')\}$, and $y'\in \{z'\mid h(y')\sqsubseteq h(z')\}$. Now we claim that $y'$ is the minimum of $\{z'\mid h(y')\sqsubseteq h(z')\}$, which implies $h_\rela(h(y'))=y'$. For suppose $z'$ is a member of the set, so $h(y')\sqsubseteq h(z')$. Then since an injective Boolean
homomorphism reflects order, $y'\sqsubseteq' z'$, as desired. Hence $h_\rela(h(y'))=y'$, which shows that $h_\rela$ is surjective.

Parts \ref{BAOMorphPossMorph4}-\ref{BAOMorphPossMorph5} are easy to check.\end{proof}

From Theorems \ref{VtoPossFrames}.\ref{VtoPoss3} and \ref{BAOMorphPossMorph}, we obtain the second piece of our categorical picture.

\begin{corollary}[The $(\cdot)_\rela$ Functor]\label{FuncCor2} The $(\cdot)_\rela$ operation given by Definition \ref{Hposs} and Theorem \ref{BAOMorphPossMorph} is a contravariant functor from \textbf{$\mathcal{CV}$-BAO} to \textbf{RichPoss}.
\end{corollary} 

\subsection{\texorpdfstring{$(\cdot)_\rela$ and $(\cdot)^\under$}{(.)\_p and (.)\^{}b}, and Dual Equivalence with Rich Frames}\label{DualEquiv}
  
Let us now consider the relation between $(\cdot)_\rela$ from \S~\ref{VtoPossSection} and $(\cdot)^\under$ from \S~\ref{PossToBAO}. We begin by relating an arbitrary $\mathcal{V}$-BAO $\mathbb{A}$ to the underlying BAO $ (\mathbb{A}_\rela)^\under$ of its principal frame $\mathbb{A}_\rela$. In \S~\ref{MacNeille}, we will relate an arbitrary $\mathcal{V}$-BAO to the underlying BAO $ (\mathbb{A}_\fullV)^\under$ of its full frame $\mathbb{A}_\fullV$. Recall that for a $\mathcal{CV}$-BAO $\mathbb{A}$, $\mathbb{A}_\rela=\mathbb{A}_\fullV$.

\begin{theorem}[From BAOs to Frames and Back]\label{BAOstoFrames} Given a $\mathcal{V}$-BAO $\mathbb{A}$ whose associated Boolean lattice is $\langle A,\leq\rangle$ with minimum $\bot$, define $\zeta_\mathbb{A}\colon \mathbb{A}\to (\mathbb{A}_\rela)^\under$ by $\zeta_\mathbb{A}(x)=\{x'\in A\setminus\{\bot\}\mid x'\leq x\}$. Then:
\begin{enumerate}
\item\label{BAOstoFrames1} If $\mathbb{A}$ is a $\mathcal{V}$-BAO, then $\zeta_\mathbb{A}$ is a complete BAO-isomorphism.
\item\label{BAOstoFrames2} If $\mathbb{A}$ and $\mathbb{B}$ are $\mathcal{CV}$-BAOs, and $h\colon \mathbb{A}\to\mathbb{B}$ is a complete BAO-homomorphism, then $(h_\rela)^\under \circ \zeta_\mathbb{A}= \zeta_\mathbb{B}\circ h$, so the following diagram commutes:
\end{enumerate}
\begin{center}
\begin{tikzpicture}[->,>=stealth',shorten >=1pt,shorten <=1pt, auto,node
distance=2cm,thick,every loop/.style={<-,shorten <=1pt}]
\tikzstyle{every state}=[fill=gray!20,draw=none,text=black]

\node (A) at (0,2) {{$\mathbb{A}$}};
\node (A') at (0,0) {{$(\mathbb{A}_\rela)^\under$}};
\node (B) at (3,2) {{$\mathbb{B}$}};
\node (B') at (3,0) {{$(\mathbb{B}_\rela)^\under$}};

\path (A) edge[->] node {{$h$}} (B);
\path (A') edge[<-] node {{$\zeta_\mathbb{A}$}} (A);
\path (B) edge[->] node {{$\zeta_\mathbb{B}$}} (B');
\path (B') edge[<-] node {{$(h_\rela)^\under$}} (A');

\end{tikzpicture}
\end{center}
\end{theorem}

\begin{proof} For part \ref{BAOstoFrames1}, it is straightforward to check that $\zeta_\mathbb{A}$ is a complete Boolean algebra isomorphism, so we only show that it respects the operators. Let $\mathbb{A}=\langle A,\meet,-,\top, \{\blacksquare_i\}_{i\in\ind}\rangle$. For $x\in \mathbb{A}$, $\zeta_\mathbb{A}(\blacksquare_i^\mathbb{A} x)  =  \{y\in A\setminus \{\bot\}\mid y\leq \blacksquare_i^\mathbb{A} x\} $ by the definition of $\zeta_\mathbb{A}$, and $\blacksquare_i^{ (\mathbb{A}_\rela)^\under} \zeta_\mathbb{A} (x)=\{y\in A\setminus\{\bot\}\mid R_i^{\mathbb{A}_\rela}(y)\subseteq \zeta_\mathbb{A}(x)\}$ by Definition \ref{RegOpAlg} and the fact from Definition \ref{Hposs} that $\mathbb{A}_\rela=A\setminus\{\bot\}$. So we must show that for all $y\in A\setminus\{\bot\}$, $y\leq \blacksquare_i^\mathbb{A}x$ iff $R_i^{\mathbb{A}_\rela}(y)\subseteq \zeta_\mathbb{A}(x)$. Suppose $y\leq \blacksquare_i^\mathbb{A}x$ and consider a $z$ such that $yR_i^{\mathbb{A}_\rela}z$, which means that for all $z'$ such that $\bot\not=z'\leq z$, $y\meet\blacklozenge_i^\mathbb{A} z'\not=\bot$. It follows that $z\leq x$, for if $z\not\leq x$, then $z\meet -x\not=\bot$, so by the previous step, $y\meet\blacklozenge_i^\mathbb{A} (z\meet -x)\not=\bot$, which implies $y\meet\blacklozenge_i^\mathbb{A} \mathord{-}x\not=\bot$, which contradicts $y\leq \blacksquare_i^\mathbb{A}x$. In the other direction, if $y\not\leq \blacksquare_i^\mathbb{A}x$, then $y\meet -\blacksquare_i^\mathbb{A}x\not=\bot$, so $y\meet \blacklozenge_i^\mathbb{A}\mathord{-}x\not=\bot$. Then since $\mathbb{A}$ is a $\mathcal{V}$-BAO and hence an $\mathcal{R}$-BAO by Lemma \ref{VandH}, by Definition \ref{plentiful} there is a $z$ such that $\bot\not=z\leq -x$ and $yR_i^{\mathbb{A}_\rela} z$. Then $z\not\leq x$, so $R_i^{\mathbb{A}_\rela}(y)\not\subseteq \zeta_\mathbb{A}(x)$.
 
For part \ref{BAOstoFrames2}, let $\bot_\mathbb{A}$ and $\bot_\mathbb{B}$ be the minimums of $\mathbb{A}$ and $\mathbb{B}$, respectively. For $a\in\mathbb{A}$, we have
\[\begin{array}{llll}
(h_\rela)^\under(\zeta_\mathbb{A}(a))&=&(h_\rela)^\under(\{a'\in A\setminus\{\bot_\mathbb{A}\}\mid a'\leq a\}) & \mbox{by definition of $\zeta_\mathbb{A}$}\\
&=&h_\rela^{-1}[\{a'\in A\setminus\{\bot_\mathbb{A}\}\mid a'\leq a\}] & \mbox{by definition of $(\cdot)^\under$}\\
&=& \{b\in B\setminus\{\bot_\mathbb{B}\}\mid  h_\rela (b)\leq a \} & \mbox{by definition of $(\cdot)^{-1}$} \\
&=& \{b\in B\setminus\{\bot_\mathbb{B}\}\mid \bigmeet \{x\in\mathbb{A} \mid b\leq h(x)\}\leq a\} & \mbox{by definition of $(\cdot)_\rela$} \\
&=& \{b\in B\setminus\{\bot_\mathbb{B}\}\mid b\leq h(a)\}& (\dagger)\\
&=& \zeta_\mathbb{B}(h(a)) & \mbox{by definition of $\zeta_\mathbb{B}$}.
\end{array}\] 
For $(\dagger)$, if $b\leq h(a)$, then $a\in \{x\in\mathbb{A}\mid b\leq h(x)\}$, so $\bigmeet \{x\in\mathbb{A}\mid b\leq h(x)\}\leq a$. In the other direction, since $h$ is a complete BAO-homomorphism, we have
\begin{eqnarray*}
&& \bigmeet \{x\in\mathbb{A}\mid b\leq h(x)\}\leq a\\
\Rightarrow &&h(\bigmeet \{x\in\mathbb{A}\mid b\leq h(x)\})\leq h(a) \\
\Rightarrow &&  \bigmeet \{h(x)\mid x\in\mathbb{A}\mbox{ and }b\leq h(x)\}\leq h(a)\\
\Rightarrow && b\leq h(a),
\end{eqnarray*}
which completes the proof.\end{proof}

In the other direction, going from a possibility frame $\mathcal{F}$ to $(\mathcal{F}^\under)_\rela$, this transformation is only well-defined if $\mathcal{F}^\under$ is a $\mathcal{V}$-BAO (recall Definition \ref{Hposs}), which is guaranteed if $\mathcal{F}$ is a full or principal possibility frame (Theorem \ref{PtoB}.\ref{PtoB2}-\ref{PtoB3}). If $\mathcal{F}=\langle S,\sqsubseteq,\{R_i\}_{i\in\ind},\adm\rangle$ is full or principal, then by Definitions \ref{RegOpAlg} and \ref{Hposs}, $(\mathcal{F}^\under)_\rela$ is such that  $S^{(\mathcal{F}^\under)_\rela}=\adm\setminus\{\emptyset\}$, $X\sqsubseteq^{(\mathcal{F}^\under)_\rela}Y$ iff $X\subseteq Y$, and $\adm^{(\mathcal{F}^\under)_\rela}$ is the set of all principal downsets in $\langle S^{(\mathcal{F}^\under)_\rela},\sqsubseteq^{(\mathcal{F}^\under)_\rela}\rangle$ plus $\emptyset$. In contrast to Theorem \ref{BAOstoFrames}, clearly $(\mathcal{F}^\under)_\rela$ is not guaranteed to be isomorphic to $\mathcal{F}$, since $(\mathcal{F}^\under)_\rela$ is always a principal possibility frame, with its poset a Boolean lattice minus $\bot$. Even if $\mathcal{F}$ is a principal frame, it might not be isomorphic to $(\mathcal{F}^\under)_\rela$, because the relation $R_i^{(\mathcal{F}^\under)_\rela}$ may relate more possibilities than $R_i$ does. To get an isomorphism we need to start with principal frames that are \textit{tight} as in \S~\ref{TightSection}. (Since all principal frames are \textit{$\sqsubseteq$-tight} by Fact \ref{PrincPossSep}, the requirement we have in mind here is \textit{$R$-tight}.)

\begin{proposition}[Tight Principal Frames and BAOs]\label{TightBAO} A principal possibility frame $\mathcal{F}$ is isomorphic to $(\mathcal{F}^\under)_\rela$ iff $\mathcal{F}$ is a tight principal possibility frame.
\end{proposition}
 
\begin{proof} For the left-to-right direction, by Theorem \ref{PtoB}.\ref{PtoB3}, for any principal possibility frame $\mathcal{F}$, $\mathcal{F}^\under$ is a $\mathcal{V}$-BAO, and by Theorem \ref{VtoPossFrames}.\ref{FromGenToPos1}, for any $\mathcal{V}$-BAO $\mathbb{A}$, $\mathbb{A}_\rela$ is a tight principal possibility frame. Thus, $(\mathcal{F}^\under)_\rela$ is a tight principal possibility frame, so any isomorphic frame is too. 

The right-to-left direction follows from Theorem \ref{Inverses}.\ref{Inverses1} below.
\end{proof}

Recall from \S~\ref{RichFrames} that \textit{rich} possibility frames are equivalent to tight principal possibility frames that are \textit{lattice-complete} (Fact \ref{RichFullStrong}). Rich frames arise now thanks to Proposition \ref{UpDownRich}.\ref{UpDownRich1}.

\begin{proposition}[From Full Frames to Rich Frames and Functional Frames]\label{UpDownRich} If $\mathcal{F}=\langle S,\sqsubseteq, \{R_i\}_{i\in\ind},\adm\rangle$ is a \textit{full} possibility frame, then:
\begin{enumerate}
\item\label{UpDownRich1} $(\mathcal{F}^\under)_\rela$ is a rich possibility frame;
\item\label{UpDownRich2} the relation $R_i^{(\mathcal{F}^\under)_\rela}$ in $(\mathcal{F}^\under)_\rela$ is given by $XR_i^{(\mathcal{F}^\under)_\rela}Y$ iff $Y\subseteq \mathrm{int}(\mathrm{cl}(\mathord{\Downarrow}R_i[X]))$; 
\item\label{UpDownRich3} $(\mathcal{F}^\under)_\rela$ can be functionalized (Proposition \ref{QtoF}) by $f_i(X)=\mathrm{int}(\mathrm{cl}(\mathord{\Downarrow}R_i[X]))$ (when this set is nonempty).
\end{enumerate} 
\end{proposition}

\begin{proof} Part \ref{UpDownRich1} is immediate from Theorems \ref{PtoB}.\ref{PtoB2} and \ref{VtoPossFrames}.\ref{VtoPoss3}. Part \ref{UpDownRich2} is immediate from Observation \ref{ResidFull} and Lemma \ref{RelT}. Part \ref{UpDownRich3} is immediate from part \ref{UpDownRich2} and Proposition \ref{QtoF}.
\end{proof}

Before proving the right-to-left direction of Proposition \ref{TightBAO}, we will show that if $\mathcal{F}$ is an arbitrary full or principal frame, then although we are not guaranteed that $\mathcal{F}$ is \textit{isomorphic} to $(\mathcal{F}^\under)_\rela$, we are guaranteed a rather close connection. For the following, recall the terminology for morphisms from Definition \ref{PossMorph}. 

\begin{theorem}[From Frames to BAOs and Almost Back]\label{AlmostBack} Given a full or principal possibility frame $\mathcal{F}$ whose associated poset is $\langle S,\sqsubseteq\rangle$, define $\zeta_\mathcal{F}\colon \mathcal{F}\to(\mathcal{F}^\under)_\rela$ by $\zeta_\mathcal{F}(x)=\{x'\in S\mid x'\cof x\}$ (recall Definition \ref{CoRef}). Then:
\begin{enumerate}
\item\label{Inverses0} $\zeta_\mathcal{F}$ is a strict, dense, and robust possibility morphism from $\mathcal{F}$ to $(\mathcal{F}^\under)_\rela$;
\item\label{Inverses.5} if $\mathcal{F}$ is \textit{separative}, then $\zeta_\mathcal{F}$ is a $\sqsubseteq$-strong embedding; 
\item\label{Inverses.75} if $\mathcal{F}$ is separative and \textit{strong}, then $\zeta_\mathcal{F}$ is a strong embedding.
\end{enumerate}
\end{theorem}
Note that if $\mathcal{F}$ is full, then Fact \ref{CofGenerated} implies that $\zeta_\mathcal{F}(x)\in \adm^\mathcal{F}$; and if $\mathcal{F}$ is principal, then Fact \ref{Boolean} implies that $\cof\,=\,\sqsubseteq$, so $\zeta_\mathcal{F}(x)=\{x'\in S\mid x'\sqsubseteq x\}\in \adm^\mathcal{F}$. Thus, $\zeta_\mathcal{F}(x)$ is indeed in the domain of $(\mathcal{F}^\under)_\rela$. 
\begin{proof} Let $\mathcal{F}=\langle S,\sqsubseteq, \{R_i\}_{i\in\ind},\adm\rangle$ and $(\mathcal{F}^\under)_\rela=\langle S^{(\mathcal{F}^\under)_\rela},\sqsubseteq^{(\mathcal{F}^\under)_\rela}, \{R_i^{(\mathcal{F}^\under)_\rela}\}_{i\in\ind},\adm^{(\mathcal{F}^\under)_\rela}\rangle$. We begin with the observation ($\star$) that if $X\in S^{(\mathcal{F}^\under)_\rela}$ and $x\in X$, then by Fact \ref{CofClose}, $\zeta_\mathcal{F}(x)\subseteq X$, so $\zeta_\mathcal{F}(x)\sqsubseteq^{(\mathcal{F}^\under)_\rela} X$. Since each $X\in S^{(\mathcal{F}^\under)_\rela}$ is nonempty, it follows that there is an $x\in S$ such that $\zeta_\mathcal{F}(x)\sqsubseteq^{(\mathcal{F}^\under)_\rela} X$, so $\zeta_\mathcal{F}$ is \textit{dense} in the sense of Definition \ref{PossMorph}.\ref{DenseMorph}. For the special features of \textit{robust} morphisms, we must show that for all $X\in\adm$, $X=\zeta^{-1}_\mathcal{F}[\zeta_\mathcal{F}[X]]$ and there is an $\mathcal{X}'\in \adm^{(\mathcal{F}^\under)_\rela}$ such that $\zeta_\mathcal{F}[X]=\zeta_\mathcal{F}[S]\cap \mathcal{X}'$. Since $X\subseteq \zeta^{-1}_\mathcal{F}[\zeta_\mathcal{F}[X]]$ is immediate, we begin with $X\supseteq \zeta^{-1}_\mathcal{F}[\zeta_\mathcal{F}[X]]$. If $\zeta_\mathcal{F}(y)\in \zeta_\mathcal{F}[X]$, then for some $x\in X$, $\zeta_\mathcal{F}(y)=\{y'\in S\mid y'\cof y\}=\zeta_\mathcal{F}(x)=\{x'\in S\mid x'\cof x\}$. Hence $y\cof x$, which with $x\in X$ and Fact \ref{CofClose} implies $y\in X$. Thus, $X\supseteq \zeta^{-1}_\mathcal{F}[\zeta_\mathcal{F}[X]]$. For the required $\mathcal{X}'$, take $\mathcal{X}'=\{Y\in  (\mathcal{F}^\under)_\rela\mid Y\subseteq X\}$, so by ($\star$), $\zeta_\mathcal{F}[X]\subseteq \mathcal{X}'$ and hence  $\zeta_\mathcal{F}[X]\subseteq \zeta_\mathcal{F}[S]\cap \mathcal{X}'$. In the other direction, if $Y\in \zeta_\mathcal{F}[S] \cap \mathcal{X}'$, then there is a $y\in S$ such that $Y=\zeta_\mathcal{F}(y)\in \mathcal{X}'$, so $\zeta_\mathcal{F}(y)\subseteq X$. Then $y\in X$, so $\zeta_\mathcal{F}(y)\in \zeta_\mathcal{F}[X]$ and hence $Y\in \zeta_\mathcal{F}[X]$. Thus, $\zeta_\mathcal{F}[X]\supseteq \zeta_\mathcal{F}[S]\cap \mathcal{X}'$. Next, note that $\zeta_\mathcal{F}$ satisfies \PullBack{}, since $\zeta_\mathcal{F}^{-1}[\emptyset]=\emptyset$ and for every nonempty $\mathcal{X}\in \adm^{(\mathcal{F}^\under)_\rela}$, there is an $X\in \adm$ such that $\mathcal{X}=\mathord{\downarrow}X$ by Definition \ref{Hposs}.\ref{Hposs4}, and $\zeta_\mathcal{F}^{-1}[\mathcal{X}]=X$.

For the rest of part \ref{Inverses0}, we show that $\zeta_\mathcal{F}$ satisfies: 
\begin{itemize}
\item if $y\sqsubseteq x $, then $\zeta_\mathcal{F}(y)\sqsubseteq^{(\mathcal{F}^\under)_\rela} \zeta_\mathcal{F}(x)$ (\SqForth{}); 
\item if $Y\sqsubseteq^{(\mathcal{F}^\under)_\rela} \zeta_\mathcal{F}( x )$, then $\exists y$: $y\sqsubseteq  x $ and $\zeta_\mathcal{F}(y)\sqsubseteq^{(\mathcal{F}^\under)_\rela} Y$ (\SqBack{});
\item if $x R_i y$, then $\zeta_\mathcal{F}(x) R_i^{(\mathcal{F}^\under)_\rela} \zeta_\mathcal{F}(y)$ (\RForth{}); 
\item if $\zeta_\mathcal{F}(x)R^{(\mathcal{F}^\under)_\rela}_i Y$ and $Z\sqsubseteq^{(\mathcal{F}^\under)_\rela} Y$, then $\exists y$: $xR_i y$ and $\zeta_\mathcal{F}(y)\comp^{(\mathcal{F}^\under)_\rela} Z$ (\SRBack{}).
\end{itemize}
For \SqForth{}, if $y\sqsubseteq x$, then $\zeta_\mathcal{F}(y)\subseteq \zeta_\mathcal{F}(x)$, so $\zeta_\mathcal{F}(y)\sqsubseteq^{(\mathcal{F}^\under)_\rela} \zeta_\mathcal{F}(x)$. For \SqBack{}, if $Y\sqsubseteq^{(\mathcal{F}^\under)_\rela} \zeta_\mathcal{F}( x )$, so $Y\subseteq \{x'\in S\mid x'\cof x\}$, then there is a $y'\in Y$ such that $y'\cof x$, which implies that there is a $y\sqsubseteq y'$ such that $y\sqsubseteq x$. Since $Y$ satisfies \textit{persistence} and \textit{refinability}, from $y'\in Y$ and $y\sqsubseteq y'$ we have that $\{x'\in S\mid x'\cof y\}\subseteq Y$ by Fact \ref{CofClose}, so $\zeta_\mathcal{F}(y)\sqsubseteq^{(\mathcal{F}^\under)_\rela}  Y$.

To show that $\zeta_\mathcal{F}$ satisfies \RForth{}, we use the fact from Corollary \ref{Fullinterplay} and Proposition \ref{InterPrinc} that any full or principal possibility frame satisfies  \Rwinweak{} from \S~\ref{FullFrames}, plus the definition of $\blacklozenge_i Y$ from Fact \ref{ForcingDiamond}:
\[\begin{array}{llll}
& & xR_i y & \\
&\Rightarrow & \forall y'\sqsubseteq y\; \exists x'\sqsubseteq x\; \forall x''\sqsubseteq x' \;\exists y''\colon x''R_i y''\mbox{ and }y'\comp y''&\mbox{by \Rwinweak{}}\\
&\Rightarrow & \forall Y\in \adm, \mbox{ if }\emptyset\not = Y\subseteq \zeta_\mathcal{F}(y)\mbox{, then } & \mbox{since every $Y\in \adm$}\\
&& \zeta_\mathcal{F}(x)\cap \{x'\in S\mid \forall x''\sqsubseteq x'\,\exists y''\colon x''R_i y''\mbox{ and }\exists y'''\sqsubseteq y''\colon y'''\in Y\}\not=\emptyset & \mbox{is a downset} \\
&\Leftrightarrow & \forall Y\in \mathcal{F}^\under\colon\mbox{if } \bot^{\mathcal{F}^\under}\not = Y\leq^{\mathcal{F}^\under} \zeta_\mathcal{F}(y)\mbox{, then }\zeta_\mathcal{F}(x)\meet^{\mathcal{F}^\under}\blacklozenge_i Y\not= \bot^{\mathcal{F}^\under} &\mbox{by definition of $(\cdot)^\under$ and $\blacklozenge_i$}\\
&\Leftrightarrow & \zeta_\mathcal{F}(x)R_i^{(\mathcal{F}^\under)_\rela} \zeta_\mathcal{F}(y) & \mbox{by definition of $(\cdot)_\rela$}.
\end{array}\]

For \SRBack{}, suppose that $\zeta_\mathcal{F}(x)R_i^{(\mathcal{F}^\under)_\rela}Y$ and $Z\sqsubseteq^{(\mathcal{F}^\under)_\rela}Y$. By definition of $(\cdot)_\rela$ (Definition \ref{Hposs}), it follows that  $\zeta_\mathcal{F}(x)\meet \blacklozenge_i Z\not=\bot$, i.e., 
\[\{x'\in S\mid x'\cof x\}\cap \{x'\in S\mid \forall x''\sqsubseteq x'\,\exists z\colon x''R_i z\mbox{ and }\exists z'\sqsubseteq z\colon z'\in Z\}\not=\emptyset.\] 
So there is an $x'\cof x$ such that for all $x''\sqsubseteq x'$ there is a $z$ such that $x''R_i z$ and a $z'\sqsubseteq z$ such that $z'\in Z$; since $x'\cof x$, there is an $x''\sqsubseteq x'$ such that $x''\sqsubseteq x$; and since $x''\sqsubseteq x'$, there is a $z$ such that $x''R_i z$ and a $z'\sqsubseteq z$ such that $z'\in Z$. Then since $z'\in Z$, it follows by ($\star$) above that $\zeta_\mathcal{F}(z')\sqsubseteq^{(\mathcal{F}^\under)_\rela} Z$. By \Rrule{} (Corollary \ref{Fullinterplay}, Proposition \ref{InterPrinc}), together $x''\sqsubseteq x$, $x''R_iz$ and $z'\sqsubseteq z$ imply there is a $y$ such that $xR_iy$ and $y\comp z'$, so there is a $z''$ with $z''\sqsubseteq y$ and $z''\sqsubseteq z'$. Thus, by \SqForth{}, $\zeta_\mathcal{F}(z'')\sqsubseteq^{(\mathcal{F}^\under)_\rela} \zeta_\mathcal{F}(y)$ and $\zeta_\mathcal{F}(z'')\sqsubseteq^{(\mathcal{F}^\under)_\rela} \zeta_\mathcal{F}(z')$, which with $\zeta_\mathcal{F}(z')\sqsubseteq^{(\mathcal{F}^\under)_\rela} Z$ implies $\zeta_\mathcal{F}(y)\comp^{(\mathcal{F}^\under)_\rela} Z$. This establishes \SRBack{}. 

For part \ref{Inverses.5}, if $\mathcal{F}$ is separative, then $\zeta_\mathcal{F}(x)=\mathord{\downarrow}x=\{x'\in S\mid x'\sqsubseteq x\}$. Then since $\sqsubseteq$ is a partial order, we have $y\sqsubseteq x$ iff $\{y'\in S\mid y'\sqsubseteq y\}\subseteq \{x'\in S\mid x'\sqsubseteq x\}$ iff $\zeta_\mathcal{F}(y)\sqsubseteq^{(\mathcal{F}^\under)_\rela}\zeta_\mathcal{F}(x)$, which also shows that $\zeta_\mathcal{F}$ is injective. Since we already showed that $\zeta_\mathcal{F}$ is robust, it follows that $\zeta_\mathcal{F}$ is a $\sqsubseteq$-strong embedding.

For part \ref{Inverses.75}, all we need to add to part \ref{Inverses.5} is that $xR_iy$ iff $\zeta_\mathcal{F}(x)R_i^{(\mathcal{F}^\under)_\rela} \zeta_\mathcal{F}(y)$. If $\mathcal{F}$ is not only separative but also strong, i.e., satisfies \RWin{} from \S~\ref{FullFrames}, then using the definition of $\blacklozenge_i Y$ assuming \Rdown{} from Fact \ref{ForcingDiamond}, we have:
\[\begin{array}{llll}
& & xR_i y & \\
&\Leftrightarrow & \forall y'\sqsubseteq y\; \exists x'\sqsubseteq x\; \forall x''\sqsubseteq x' \;\exists y''\sqsubseteq y'\colon x''R_i y''&\mbox{by \RWin{}}\\
&\Leftrightarrow & \forall Y\in \adm, \mbox{ if }\emptyset\not = Y\subseteq \mathord{\downarrow}y\mbox{, then } &\mbox{since every $Y\in \adm$ is a downset}\\
&& \mathord{\downarrow}x\cap \{x'\in S\mid \forall x''\sqsubseteq x'\,\exists y''\in Y\colon x''R_i y''\}\not=\emptyset &\mbox{and $\forall y'\sqsubseteq y$: $\mathord{\downarrow}y'\in\adm$ by Fact \ref{Sep&Princ}.\ref{Sep&Princ1}}\\
&\Leftrightarrow & \forall Y\in \mathcal{F}^\under, \mbox{ if } \bot^{\mathcal{F}^\under}\not = Y\leq^{\mathcal{F}^\under}\mathord{\downarrow}y\mbox{, then }\mathord{\downarrow}x\meet^{\mathcal{F}^\under}\blacklozenge_i Y\not= \bot^{\mathcal{F}^\under}  &\mbox{by definition of $(\cdot)^\under$ and $\blacklozenge_i$}\\
&\Leftrightarrow &\mathord{\downarrow}xR_i^{(\mathcal{F}^\under)_\rela} \mathord{\downarrow}y& \mbox{by definition of $(\cdot)_\rela$}\\
&\Leftrightarrow & \zeta_\mathcal{F}(x)R_i^{(\mathcal{F}^\under)_\rela} \zeta_\mathcal{F}(y)&\mbox{by definition of $\zeta_\mathcal{F}$,}
\end{array}\]
which completes the proof.\end{proof}
 
We are now ready to prove the analogue of Theorem \ref{BAOstoFrames} going from $\mathcal{F}$ to $(\mathcal{F}^\under)_\rela$. 

\begin{theorem}[From Frames to BAOs and Back]\label{Inverses} $\,$
\begin{enumerate}
\item\label{Inverses1} If $\mathcal{F}$ is a \textit{tight principal} possibility frame, then $\zeta_\mathcal{F}$ is an isomorphism.
\item\label{Inverses2} If $\mathcal{F}$ and $\mathcal{G}$ are \textit{rich} possibility frames, and $h\colon \mathcal{F}\to\mathcal{G}$ is a possibility morphism, then $(h^\under)_\rela \circ \zeta_\mathcal{F}=\zeta_\mathcal{G}\circ h$, so the following diagram commutes:
\begin{center}
\begin{tikzpicture}[->,>=stealth',shorten >=1pt,shorten <=1pt, auto,node
distance=2cm,thick,every loop/.style={<-,shorten <=1pt}]
\tikzstyle{every state}=[fill=gray!20,draw=none,text=black]

\node (F) at (0,2) {{$\mathcal{F}$}};
\node (F') at (0,0) {{$(\mathcal{F}^\under)_\rela$}};
\node (G) at (3,2) {{$\mathcal{G}$}};
\node (G') at (3,0) {{$(\mathcal{G}^\under)_\rela$}};

\path (F) edge[->] node {{$h$}} (G);
\path (F') edge[<-] node {{$\zeta_\mathcal{F}$}} (F);
\path (G) edge[->] node {{$\zeta_\mathcal{G}$}} (G');
\path (G') edge[<-] node {{$(h^\under)_\rela$}} (F');

\end{tikzpicture}
\end{center}
\end{enumerate}
\end{theorem}

\begin{proof} Since we are dealing with principal frames, we have $\zeta_\mathcal{F}(x)=\mathord{\downarrow}x=\{x'\in S\mid x'\sqsubseteq x\}$.

For part \ref{Inverses1}, by Fact \ref{PrincPossSep} and Proposition \ref{InterPrinc}.\ref{InterPrinc2}, every tight principal frame is separative and strong, so by Theorem \ref{AlmostBack}.\ref{Inverses.75} we have that $\zeta_\mathcal{F}$ is a strong embedding. Moreover, since the domain of $(\mathcal{F}^\under)_\rela$ is the set of principal downsets in $\mathcal{F}$, and $\zeta_\mathcal{F}(x)=\mathord{\downarrow}x=\{x'\in S\mid x'\sqsubseteq x\}$ since $\mathcal{F}$ is principal, we have that $\zeta_\mathcal{F}$ is surjective. Then since a surjective strong embedding is equivalent to an isomorphism, $\zeta_\mathcal{F}$ is an isomorphism. 

For part \ref{Inverses2}, assuming $\mathcal{F}$ and $\mathcal{G}$ are \textit{rich} possibility frames, $\mathcal{F}^\under$ and $\mathcal{G}^\under$ are $\mathcal{CV}$-BAOs by Theorem \ref{PtoB}.\ref{PtoB2}. Thus, for any possibility morphism $h\colon \mathcal{F}\to\mathcal{G}$, which gives us the complete BAO-homomorphism $h^\under\colon {\mathcal{G}^\under\to\mathcal{F}^\under}$ as in Theorem \ref{PossMorphBAOMorph}, we can form the possibility morphism $(h^\under)_\rela\colon (\mathcal{F}^\under)_\rela\to(\mathcal{G}^\under)_\rela$ as in Theorem \ref{BAOMorphPossMorph}. Now for $x\in S^{\mathcal{F}}$, we have:
\[\begin{array}{llll}
(h^\under)_\rela \circ \zeta_\mathcal{F}(x)&=& (h^\under)_\rela(\{x'\in S^{\mathcal{F}}\mid x'\sqsubseteq^\mathcal{F} x\})&\mbox{by definition of $\zeta_\mathcal{F}$}\\
&=& \bigmeet \{X\in (\mathcal{G}^\under)_\rela\mid \{x'\in S^{\mathcal{F}}\mid x'\sqsubseteq^\mathcal{F} x\} \sqsubseteq^{(\mathcal{F}^\under)_\rela} h^\under(X)\}&\mbox{by definition of $(\cdot)_\rela$} \\
&=& \bigmeet \{X\in (\mathcal{G}^\under)_\rela\mid \{x'\in S^{\mathcal{F}}\mid x'\sqsubseteq^\mathcal{F} x\}\sqsubseteq^{(\mathcal{F}^\under)_\rela} h^{-1}[X]\} &\mbox{by definition of $(\cdot)^\under$}\\
&=& \bigmeet \{X\in (\mathcal{G}^\under)_\rela\mid \{x'\in S^{\mathcal{F}}\mid x'\sqsubseteq^\mathcal{F} x\} \subseteq h^{-1}[X]\} &\mbox{by definition of $(\mathcal{F}^\under)_\rela$} \\
&=& \bigmeet \{X\in (\mathcal{G}^\under)_\rela\mid \forall x'\sqsubseteq^\mathcal{F} x,\, h(x')\in X\} &\mbox{by definition of }h^{-1}\\
&=& \{y\in S^\mathcal{G}\mid y\sqsubseteq^\mathcal{G} h(x)\} &  (\ddagger) \\
&=& \zeta_\mathcal{G}(h(x)) &\mbox{by definition of $\zeta_\mathcal{G}$.}
\end{array}\]
For $(\ddagger)$, since $\mathcal{G}$ is a rich and hence principal frame, Definitions \ref{Hposs} and \ref{RegOpAlg} imply that the domain of $(\mathcal{G}^\under)_\rela$ is the set of all principal downsets in $\mathcal{G}$, so $\mathord{\downarrow}^\mathcal{G}h(x)=\{y\in S^\mathcal{G}\mid y\sqsubseteq^\mathcal{G} h(x)\}$ is in $(\mathcal{G}^\under)_\rela$. Moreover, Fact~\ref{PrincOrd} implies that for all $x'\sqsubseteq^\mathcal{F} x$, $h(x')\sqsubseteq^\mathcal{G} h(x)$, so $h(x')\in \mathord{\downarrow}^\mathcal{G}h(x)$, which implies that $\mathord{\downarrow}^\mathcal{G}h(x)\in \{X\in (\mathcal{G}^\under)_\rela\mid \forall x'\sqsubseteq^\mathcal{F} x,\, h(x')\in X\}$. Now consider any $X\in (\mathcal{G}^\under)_\rela$ such that $\forall x'\sqsubseteq^\mathcal{F} x$, $h(x')\in X$, so in particular, $h(x)\in X$. Then since $X$ is a downset in $\mathcal{G}$, $\mathord{\downarrow}^\mathcal{G}h(x)\subseteq X$, so $\mathord{\downarrow}^\mathcal{G}h(x)\sqsubseteq^{(\mathcal{G}^\under)_\rela} X$. Putting it all together, we have shown that $\mathord{\downarrow}^\mathcal{G}h(x)$ is the greatest lower bound of $\{X\in (\mathcal{G}^\under)_\rela\mid \forall x'\sqsubseteq^\mathcal{F} x,\, h(x')\in X\}$ in $(\mathcal{G}^\under)_\rela$, as $(\ddagger)$~claims.\end{proof}

Putting together Theorems \ref{BAOstoFrames} and \ref{Inverses} and Corollaries \ref{FuncCor1} and \ref{FuncCor2}, we have the following analogue of Thomason's \citeyearpar{Thomason1975} dual equivalence result for the categories of $\mathcal{CAV}$-BAOs with complete BAO-homomorphism and of Kripke frames with p-morphisms.

\begin{theorem}[Dual Equivalence]\label{CVDuals} \textbf{$\mathcal{CV}$-BAO} is dually equivalent to \textbf{RichPoss}.\end{theorem}  

\subsection{Reflection with Rich Frames}\label{RichReflectionSection}

Next we observe that \textbf{RichPoss}, whose morphisms are p-morphisms, is a full subcategory of \textbf{FullPoss}, whose morphisms are strict possibility morphisms. 

\begin{proposition}[Full Subcategory]\label{PossToP1} Every possibility morphism between rich possibility frames is a p-morphism.
\end{proposition}

\begin{proof} Immediate from Theorems \ref{BAOMorphPossMorph}.\ref{BAOMorphPossMorph1} and \ref{Inverses}.
\end{proof}

In fact, we will show that \textbf{RichPoss} is a \textit{reflective subcategory} of \textbf{FullPoss}.  Recall that this means it is a full subcategory such that for any full frame $\mathcal{F}$, there is a rich frame $\mathbf{R}(\mathcal{F})$ and a \textbf{FullPoss}-morphism $r\colon \mathcal{F}\to \mathbf{R}(\mathcal{F})$  such that for any rich frame $\mathcal{G}$ and \textbf{FullPoss}-morphism $g\colon \mathcal{F}\to\mathcal{G}$, there is a unique \textbf{RichPoss}-morphism $\overline{g}\colon \mathbf{R}(\mathcal{F})\to\mathcal{G}$ such that $g=\overline{g}\circ r$.\footnote{Sometimes this is not taken to be the definition of being a reflective subcategory, but rather a condition that is proved equivalent to being a reflective subcategory (e.g., \citealt[\S I.18]{Balbes1974}).} The pair $\langle \mathbf{R}(\mathcal{F}),r\rangle$ is called the \textit{reflection} of $\mathcal{F}$ in \textbf{RichPoss}. We will take as the reflection of $\mathcal{F}$ the pair $\langle (\mathcal{F}^\under)_\rela,\zeta_\mathcal{F}\rangle$ with $\zeta_\mathcal{F}$ from Theorem \ref{AlmostBack}. We already showed in Theorem \ref{AlmostBack}.\ref{Inverses0} that $\zeta_\mathcal{F}$ is a strict possibility morphism and hence a morphism in \textbf{FullPoss}.
    
\begin{theorem}[Rich Reflections]\label{RefSub} For any full possibility frame $\mathcal{F}$, rich possibility frame $\mathcal{G}$, and possibility morphism $g\colon \mathcal{F}\to\mathcal{G}$, define $\overline{g}\colon (\mathcal{F}^\under)_\rela\to\mathcal{G}$ such that for $X\in  (\mathcal{F}^\under)_\rela$, $\overline{g}(X)=\bigvee g[X]$, where $\bigvee$ is the join operation in the complete Boolean lattice underlying the rich frame $\mathcal{G}$. Then:
\begin{enumerate}
\item\label{RefSub1} $\overline{g}$ is a p-morphism;
\item\label{RefSub2} $\overline{g}$ is the unique possibility morphism from $(\mathcal{F}^\under)_\rela$ to $\mathcal{G}$ such that $g=\overline{g}\circ \zeta_\mathcal{F}$, so the following diagram commutes:
\begin{center}
\begin{tikzpicture}[->,>=stealth',shorten >=1pt,shorten <=1pt, auto,node
distance=2cm,thick,every loop/.style={<-,shorten <=1pt}]
\tikzstyle{every state}=[fill=gray!20,draw=none,text=black]

\node (F) at (0,2) {{$\mathcal{F}$}};
\node (F') at (0,0) {{$(\mathcal{F}^\under)_\rela$}};
\node (G) at (3,2) {{$\mathcal{G}$}};

\path (F) edge[->] node {{$g$}} (G);
\path (F') edge[<-] node {{$\zeta_\mathcal{F}$}} (F);
\path (G) edge[<-] node {{$\overline{g}$}} (F');

\end{tikzpicture}
\end{center}
\end{enumerate}
\end{theorem} 

\begin{proof} For part \ref{RefSub1}, we use the fact that since $\mathcal{G}$ is rich, the map $\zeta_\mathcal{G}\colon \mathcal{G}\to (\mathcal{G}^\under)_\rela$ is an isomorphism by Theorem \ref{Inverses}.\ref{Inverses1}, and the map $(g^\under)_\rela: (\mathcal{F}^\under)_\rela\to (\mathcal{G}^\under)_\rela$ is a p-morphism by Theorem \ref{BAOMorphPossMorph}.\ref{BAOMorphPossMorph1}, which applies since $g^\under$ is a complete BAO-homomorphism by Theorem \ref{PossMorphBAOMorph}.\ref{MtoM1.5} and Fact \ref{RichFullStrong}. Thus, to show that $\overline{g}$ is a p-morphism, it suffices to prove that $\overline{g}=\zeta_\mathcal{G}^{-1}\circ (g^\under)_\rela$, or equivalently, $\zeta_\mathcal{G}\circ \overline{g}=(g^\under)_\rela$.  On one hand, for each $X\in (\mathcal{F}^\under)_\rela$:
\[\begin{array}{llll}
\zeta_\mathcal{G}(\overline{g}(X))&=&\{x'\in \mathcal{G}\mid x'\cof^\mathcal{G}\overline{g}(X)\} &\mbox{by definition of }\zeta_\mathcal{G}\\
&=&\{x'\in \mathcal{G}\mid x'\sqsubseteq^\mathcal{G}\overline{g}(X)\} &\mbox{since }\mathcal{G}\mbox{ is rich and hence separative}\\
&=&\{x'\in \mathcal{G}\mid x'\sqsubseteq^\mathcal{G}\bigvee g[X]\}  &\mbox{by definition of }\overline{g}.
\end{array}
\]
 On the other hand, for each $X\in (\mathcal{F}^\under)_\rela$:
\[\begin{array}{llll}
(g^\under)_\rela(X)&=&\bigwedge\{Y\in (\mathcal{G}^\under)_\rela\mid X\leq^{\mathcal{F}^\under} g^\under (Y)\} &\mbox{by definition of }(\cdot)_\rela \\
&=&\bigwedge\{Y\in (\mathcal{G}^\under)_\rela\mid X\subseteq g^{-1}[Y]\} &\mbox{by definition of }(\cdot)^\under\\
&=&\bigcap\{Y\in (\mathcal{G}^\under)_\rela\mid X\subseteq g^{-1}[Y]\} &\mbox{because }\mathcal{G}\mbox{ is full}\\
&=&\bigcap\{Y\in (\mathcal{G}^\under)_\rela\mid g[X]\subseteq Y\}.
\end{array}
\]
If $Y\in (\mathcal{G}^\under)_\rela$, then since $\mathcal{G}$ is rich and hence principal, $Y=\{y'\in\mathcal{G}\mid y'\sqsubseteq^\mathcal{G}y\}$ for some $y\in\mathcal{G}$, so $g[X]\subseteq Y$ implies $\bigvee g[X]\sqsubseteq^\mathcal{G} y$, which in turn implies $\zeta_\mathcal{G}(\overline{g}(X))\subseteq Y$. Thus, $\zeta_\mathcal{G}(\overline{g}(X))\subseteq (g^\under)_\rela(X)$. Conversely, $g[X]\subseteq \zeta_\mathcal{G}(\overline{g}(X))\in (\mathcal{G}^\under)_\rela$, so $\zeta_\mathcal{G}(\overline{g}(X))\in \{Y\in (\mathcal{G}^\under)_\rela\mid g[X]\subseteq Y\}$ and hence $\zeta_\mathcal{G}(\overline{g}(X))\supseteq (g^\under)_\rela(X)$.

For part \ref{RefSub2}, for each $x\in\mathcal{F}$, the reflexivity of $\cof^\mathcal{F}$ implies $x\in\zeta_\mathcal{F}(x)$. In addition, since $\mathcal{G}$ is rich and hence full and separative, Fact \ref{SepOrd}.\ref{SepOrd2} implies $g(y)\sqsubseteq^\mathcal{G}g(x)$ for every $y\in\zeta_\mathcal{F}(x)$. Hence $\overline{g}(\zeta_\mathcal{F}(x))=\bigvee g[\zeta_\mathcal{F}(x)]=g(x)$, so $\overline{g}\circ\zeta_\mathcal{F}=g$, as desired.

Finally, we prove that $\overline{g}$ is the \textit{unique} possibility morphism from $(\mathcal{F}^\under)_\rela$ to $\mathcal{G}$ such that $g=\overline{g}\circ \zeta_\mathcal{F}$. Suppose $\overline{h}$ is a possibility morphism from $(\mathcal{F}^\under)_\rela$ to $\mathcal{G}$ such that $g=\overline{h}\circ \zeta_\mathcal{F}$. First, we claim that for all $X\in (\mathcal{F}^\under)_\rela$, $\bigvee g[X]\sqsubseteq^\mathcal{G} \overline{h}(X)$. If not, then there is an $x\in X$ such that $g(x)\not\sqsubseteq^\mathcal{G} \overline{h}(X)$. Then since $g=\overline{h}\circ \zeta_\mathcal{F}$, we have $\overline{h}(\zeta_\mathcal{F}(x))\not\sqsubseteq^\mathcal{G} \overline{h}(X)$. Since $X\in (\mathcal{F}^\under)_\rela$, $X$ is closed under $\cof^\mathcal{F}$ by Fact \ref{CofClose}, so $x\in X$ implies $\zeta_\mathcal{F}(x)=\{x'\in S\mid x'\cof^\mathcal{F} x\}\subseteq X$, so $\zeta_\mathcal{F}(x)\sqsubseteq^{(\mathcal{F}^\under)_\rela} X$. By Fact \ref{PrincOrd}, $\overline{h}$ satisfies \SqForth{}, so $\zeta_\mathcal{F}(x)\sqsubseteq^{(\mathcal{F}^\under)_\rela} X$ implies $\overline{h}(\zeta_\mathcal{F}(x))\sqsubseteq^\mathcal{G} \overline{h}(X)$, contradicting what we obtained above. Thus, $\bigvee g[X]\sqsubseteq^\mathcal{G} \overline{h}(X)$. Second, we claim that $\overline{h}(X)\sqsubseteq^\mathcal{G} \bigvee g[X]$. If not, then there is a $Y'\sqsubseteq^\mathcal{G} \overline{h}(X)$ such that in the Boolean algebra associated with $\mathcal{G}$'s poset, $Y'\meet  \bigvee g[X]=\bot$.  By Fact \ref{PrincOrd}, $\overline{h}$ satisfies \SqBack{}, so $Y'\sqsubseteq^\mathcal{G} \overline{h}(X)$ implies that there is a $Y\in (\mathcal{F}^\under)_\rela$ such that $Y\sqsubseteq^{(\mathcal{F}^\under)_\rela} X$ and $\overline{h}(Y)\sqsubseteq^\mathcal{G} Y'$. Thus, $Y'\meet  \bigvee g[X]=\bot$ implies $\overline{h}(Y)\meet  \bigvee g[X]=\bot$, which implies that for all $x\in X$, $\overline{h}(Y)\meet g(x)=\bot$. Since $Y\sqsubseteq^{(\mathcal{F}^\under)_\rela} X$, $Y\cap X\not=\emptyset$, so take an $x\in Y\cap X$.  As above, since $Y\in (\mathcal{F}^\under)_\rela$, $x\in Y$ implies $\zeta_\mathcal{F}(x)=\{x'\in S\mid x'\cof^\mathcal{F} x\}\subseteq Y$, so $\zeta_\mathcal{F}(x)\sqsubseteq^{(\mathcal{F}^\under)_\rela} Y$, which with \SqForth{} for $\overline{h}$ gives us $\overline{h}(\zeta_\mathcal{F}(x))\sqsubseteq^{\mathcal{G}} \overline{h}(Y)$ and hence $g(x)\sqsubseteq^{\mathcal{G}} \overline{h}(Y)$, which contradicts $\overline{h}(Y)\meet g(x)=\bot$ above. Thus, $\overline{h}(X)\sqsubseteq^\mathcal{G} \bigvee g[X]$. We have shown that $\overline{h}(X)=\bigvee g[X]$, so $\overline{h}=\overline{g}$ by the definition of $\overline{g}$.\end{proof}

From Theorem \ref{AlmostBack}.\ref{Inverses0}, Proposition \ref{PossToP1}, and Theorem \ref{RefSub} and the definition of reflective subcategories, we obtain our desired result.

\begin{theorem}[Reflective Subcategories]  \textbf{RichPoss} is a reflective subcategory of \textbf{FullPoss}.
\end{theorem}

\subsection{From Arbitrary BAOs to Possibility Frames}\label{GFPF}

Going beyond $\mathcal{V}$-BAOs, \textit{any} BAO can be transformed into a semantically equivalent world frame and hence possibility frame, namely its \textit{general ultrafilter frame} (see Appendix \S~\ref{AlgSem}). Below we give another way of transforming any BAO into a semantically equivalent possibility frame, namely its general \textit{filter} frame. While worlds must come from ultrafilters, possibilities can come from any proper filters. 

Recall that a \textit{proper filter} in a BAO $\mathbb{A}=\langle A, \meet, -, \top, \{\blacksquare_i\}_{i\in \ind}\rangle$ is a nonempty $F\subseteq A$ such that for all $x,y\in A$: $x,y\in F$ implies $x\meet y\in F$ ($F$ is downward directed); if $x\leq y$ and $x\in F$, then $y\in F$ ($F$ is an upset); and $\bot\not\in F$, which with the previous condition is equivalent to $F\subsetneq A$ ($F$ is proper).

\begin{definition}[Filter Frames and General Filter Frames]\label{FiltEx} Given a BAO $\mathbb{A}=\langle A, \meet, -, \top, \{\blacksquare_i\}_{i\in \ind}\rangle$ and algebraic model $\mathbb{M}=\langle\mathbb{A},\theta\rangle$, we define the \textit{filter frame} $\mathbb{A}_{\ff}=\langle S,\sqsubseteq,\{R_i\}_{i\in\ind},\adm_{\ff}\rangle$, the \textit{general filter frame} $\mathbb{A}_{\gff}=\langle S,\sqsubseteq,\{R_i\}_{i\in\ind},\adm_{\gff}\rangle$, and $\mathbb{M}_{\gff}=\langle S,\sqsubseteq,\{R_i\}_{i\in\ind},\pi\rangle$ as follows: 
\begin{enumerate}[label=\arabic*.,ref=\arabic*]
\item $S$ is the set of proper filters in $\mathbb{A}$; 
\item\label{FiltEx2} $X\sqsubseteq Y$ iff $X\supseteq Y$;
\item\label{FiltEx3} $XR_iY$ iff for all $x\in A$, if $\blacksquare_i x\in X$ then $x\in Y$;
\item\label{FiltEx4} $\adm_{\ff}=\mathrm{RO}(S,\sqsubseteq)$ (recall Notation \ref{ROnotation});
\item\label{FiltEx4.5} where $\widehat{x}=\{X\in S\mid x\in X\}$, $\adm_{\gff}=\{\widehat{x}\mid x\in A\}$;
\item\label{FiltEx5} $\pi(p)= \{X\in S\mid \tilde{\theta}(p)\in X\}$. \hfill $\triangleleft$
\end{enumerate} 
\end{definition}
As observed by Stone \citeyearpar{Stone1936} and Tarski \citeyearpar{Tarski1937b}, the set of \textit{ideals} of a Boolean algebra ordered by \textit{inclusion} ($I\leq I'$ iff $I\subseteq I'$) forms a \textit{complete Heyting algebra}.\footnote{Since the lattice of open sets of any topological space is a complete Heyting algebra, also known as a \textit{locale}, the observation also follows from the fact that the lattice of ideals of a Boolean algebra is isomorphic to the lattice of open sets in the dual Stone space \citep{Stone1937}. The lattice of ideals of a Boolean algebra is therefore a \textit{compact} and \textit{zero-dimensional} locale.} The same holds for the set of filters of a Boolean algebra ordered by inclusion, which is isomorphic as a lattice to the set of ideals ordered by inclusion. But note that in Definition \ref{FiltEx}.\ref{FiltEx2} we order the (proper) filters by \textit{reverse} inclusion. The set of filters of a Boolean algebra ordered by reverse inclusion is therefore what is known as a complete \textit{co-Heyting} algebra. Thus, the (general) filter frame of a BAO is always based on a complete co-Heyting algebra minus its bottom element.  

The following notation and fact about filters will be useful in what follows.

\begin{fact}[Filter Generated by a Subset]\label{FiltGen} Given a BAO $\mathbb{A}=\langle A, \meet, -, \top, \{\blacksquare_i\}_{i\in \ind}\rangle$ and $X\subseteq A$, let $[X)$ be the smallest filter in $\mathbb{A}$ that contains $X$, which must exist because the intersection of filters is a filter. Then  $[X)=\{x\in A\mid \exists x_1,\dots,x_n\in X\mid x_1\meet\dots\meet x_n\leq x\}$. For $y\in A$, $[\{y\})=\mathord{\uparrow} y=\{x\in A\mid y\leq x\}$.
\end{fact}
 
While the claim made above about the relationship between a BAO and its general \textit{ultrafilter} frame requires going beyond ZF set theory, since it requires the use of the ultrafilter axiom, our claim about the relationship between a BAO and its general \textit{filter} frame does not go beyond ZF. For the following result, recall the notions of \textit{strong} frames from Definition \ref{StrongPoss} and of \textit{tight} frames from Definition \ref{TightFrames}. Also recall that we say $\varphi$ is \textit{satisfiable} over a BAO $\mathbb{A}$ iff there is some algebraic model $\langle \mathbb{A},\theta\rangle$ with $\tilde{\theta}(\varphi)\not=\bot$.

\begin{theorem}[From BAOs to Possibility Frames]\label{BAOtoGenPos} For any BAO $\mathbb{A}$ and algebraic model $\mathbb{M}=\langle \mathbb{A},\theta\rangle$:
\begin{enumerate} 
\item\label{IsRelGen} $\mathbb{A}_{\gff}$ is a \textit{strong} and \textit{tight} possibility frame, and $\mathbb{A}_{\ff}$ is a strong, tight, and \textit{full} possibility frame;
\item\label{IsBasedOn} $\mathbb{M}_{\gff}$ is a possibility model based on $\mathbb{A}_{\gff}$ and $\mathbb{A}_{\ff}$;
\item\label{BAOtoGenPos3}  for all $\varphi\in\mathcal{L}(\sig,\ind)$ and $X\in \mathbb{A}_{\gff}$: $\mathbb{M}_{\gff},X\Vdash \varphi$ iff $\tilde{\theta}(\varphi)\in X$; 
\item\label{BAOtoGenPos4} for all $\varphi\in\mathcal{L}(\sig,\ind)$ and nonzero $x\in \mathbb{A}$,  $x\leq\tilde{\theta}(\varphi)$ iff $\mathbb{M}_{\gff},\mathord{\uparrow}x\Vdash \varphi$; 
\item\label{BAOtoGenPos5} for all $\varphi\in\mathcal{L}(\sig,\ind)$, if $\varphi$ is satisfiable over $\mathbb{A}$, then $\varphi$ is satisfiable over $\mathbb{A}_\gff$ and $\mathbb{A}_\ff$.
\end{enumerate}
\end{theorem}
\begin{proof} For part \ref{IsRelGen}, we first show that $\adm_{\gff}$ is closed under $\cap$, $\supset$, and $\blacksquare_i^{\mathbb{A}_\gff}$ as in Definition \ref{PosetMod}, so $\mathbb{A}_\gff$ is a partial-state frame. (Note the distinction between $\blacksquare_i^{\mathbb{A}_\gff}$ and $\blacksquare_i$, the latter being the operator in $\mathbb{A}$.)  Since $\widehat{\bot}=\emptyset$, $\emptyset\in \adm_{\gff}$. Now consider $\mathcal{X},\mathcal{Y}\in \adm_{\gff}$, so there are $x,y\in A$ such that $\mathcal{X}=\widehat{x}$ and $\mathcal{Y}=\widehat{y}$. We claim that:
\begin{itemize}
\item[(i)] $\widehat{x}\cap \widehat{y}=\widehat{x\meet y}$;
\item[(ii)] $\widehat{x}\supset\widehat{y}=\{Z\in S\mid \forall Z'\sqsubseteq Z,\, Z'\in \widehat{x}\Rightarrow Z'\in \widehat{y}\}=\widehat{-x\join y\;\;}$\hspace{-.05in};
\item[(iii)] $\blacksquare_i^{\mathbb{A}_\gff}\widehat{x}=\{Z\in S\mid R_i(Z)\subseteq \widehat{x}\}= \widehat{\blacksquare_i x}$.
\end{itemize}
For part (i), for any $Z\in S$, we have: $Z\in \widehat{x}\cap \widehat{y}$ iff $x,y\in Z$ iff $x\meet y\in Z$ (since $Z$ is a filter) iff $Z\in \widehat{x\meet y}$. Part (ii) is also easy to check. For part (iii), we have $Z\in \widehat{\blacksquare_i x}$ iff $\blacksquare_i x\in Z$, and $\blacksquare_i x\in Z$ only if (by the definition of $R_i$) for every $Z'\in R_i(Z)$, $x\in Z'$, i.e., $Z'\in \widehat{x}$, which is equivalent to $R_i(Z)\subseteq \widehat{x}$. Conversely, if $\blacksquare_i x\not\in Z$, then $Y=\{y \mid\blacksquare_i y\in Z\}$ is a proper filter such that $ZR_iY$ and $x\not\in Y$, so $R_i(Z)\not\subseteq \widehat{x}$.

Second, we show that $\adm_\gff\subseteq\mathrm{RO}(S,\sqsubseteq)$, i.e., every $\mathcal{X}\in \adm_{\gff}$ satisfies \textit{persistence} and \textit{refinability}, so $\mathbb{A}_\gff$ is a possibility frame. By definition of $\adm_{\gff}$, $\mathcal{X}=\widehat{x}$ for some $x\in A$. For \textit{persistence}, if $X'\sqsubseteq X$, so $X'\supseteq X$, and $X\in \widehat{x}$, so $x\in X$,  then $x\in X'$ and hence $X'\in \widehat{x}$.  For \textit{refinability}, if $X\not\in \widehat{x}$, so $x\not\in X$, then $X'=[X\cup\{-x\})$ is a proper filter, and $X'\supseteq X$, so $X'\sqsubseteq X$. Moreover, for every proper filter $X''\sqsubseteq X'$, i.e., $X''\supseteq X'$, we have $x\not\in X''$, so $X''\not\in \widehat{x}$. Thus, $\widehat{x}$ satisfies \textit{refinability}.

Third, to see that $\mathbb{A}_\gff$ is \textit{$R$-tight}, suppose that \textit{not} $XR_iY$, so by Definition \ref{FiltEx}.\ref{FiltEx3} there is some $x\in A$ such that $\blacksquare_i x\in X$ but $x\not\in Y$. Then $X\in \widehat{\blacksquare_i x}$ but $Y\not\in\widehat{x}$. It follows by (iii) above that $X\in \blacksquare_i^{\mathbb{A}_\gff}\widehat{x}$ but $Y\not\in \widehat{x}$, and by definition of $\adm_\gff$, $\widehat{x}\in\adm_\gff$. Thus, if for all $\mathcal{Z}\in \adm_\gff$, $X\in\blacksquare^{\mathbb{A}_\gff}_i\mathcal{Z}$ implies $Y\in\mathcal{Z}$, then $XR_iY$. So $\mathbb{A}_\gff$ is \textit{$R$-tight}. For \textit{$\sqsubseteq$-tight}, if $Y\not\sqsubseteq X$, so $Y\not\supseteq X$, then taking an $x\in X$ such that $x\not\in Y$, we have $X\in\widehat{x}$ but $Y\not\in\widehat{x}$, and $\widehat{x}\in \adm_\gff$. Thus, if for all $\mathcal{Z}\in\adm_\gff$, $X\in \mathcal{Z}$ implies $Y\in\mathcal{Z}$, then $Y\sqsubseteq X$. Hence $\mathbb{A}_\gff$ is \textit{$\sqsubseteq$-tight}.

Finally, since $\mathbb{A}_\gff$ is tight, to show that $\mathbb{A}_\gff$ is strong it suffices by Lemma \ref{TightStrong}.\ref{TightStrong1.5} to show that it satisfies \Rref{}. To that end, consider proper filters ${X}$ and ${Y}$ such that ${X}R_i{Y}$. Where 
\begin{equation}\mathrm{X}'= {X}\cup \{\blacklozenge_i y \mid y\in {Y}\},\label{Udef}\end{equation}
suppose for reductio that $[\mathrm{X}')$ is not a proper filter, i.e., $\inc\in[\mathrm{X}')$. Then by (\ref{Udef}) and Fact \ref{FiltGen}, there are $x_1,\dots, x_m\in{X}$ and $y_1,\dots,y_k\in{Y}$ such that
\[x_1\meet \dots \meet x_m \meet \blacklozenge_i y_1\meet\dots\meet \blacklozenge_i y_k \leq\inc,\]
which implies 
\begin{equation}
x_1\meet\dots\wedge x_m\leq \blacksquare_i \mathord{-}(y_1\meet\dots\meet y_k) \label{GivesBox}
\end{equation}
by the properties of $\blacklozenge_i$ and $\blacksquare_i$ (see Definition \ref{BAOs}). Since ${X}$ is a filter, $x_1,\dots,x_m\in{X}$ implies $x_1\meet\dots\meet x_m\in {X}$, which with (\ref{GivesBox}) implies $\blacksquare_i \mathord{-}(y_1\meet\dots\meet y_k) \in{X}$, which with ${X}R_i{Y}$ implies $\mathord{-}(y_1\meet\dots\meet y_k) \in{Y}$, which contradicts the facts that $y_1,\dots,y_k\in{Y}$ and ${Y}$ is a \textit{proper} filter. Thus, ${X'}=[\mathrm{X}')$ is a proper filter.

Now consider any proper filter ${X''}\sqsubseteq{X'}$, i.e., ${X''}\supseteq{X'}$. Where
\begin{equation}\mathrm{Y}'={Y}\cup \{x \mid \blacksquare_i x\in {X''}\},\label{Vdef}\end{equation}
suppose for reductio that $[\mathrm{Y}')$ is not a proper filter, i.e., $\inc \in [\mathrm{Y}')$.  Then since $Y$ and $ \{x \mid \blacksquare_i x\in {X''}\}$ are filters, it follows by (\ref{Vdef}) and Fact \ref{FiltGen} that there are $y\in{Y}$ and $\blacksquare_i x \in{X''}$ such that $y \meet  x \leq \inc$, which implies 
\begin{equation}\blacksquare_i x \leq \blacksquare_i \mathord{-} y\label{GivesBox2}\end{equation}
by the properties of $\blacksquare_i$. Then ${X'}=[\mathrm{X}')$ and (\ref{Udef}) together imply $\blacklozenge_i y\in{X'}$. Then since ${X''}\supseteq{X'}$, we have $\blacklozenge_i y\in {X''}$. On the other hand, since ${X''}$ is a filter, $\blacksquare_i x \in{X''}$ and (\ref{GivesBox2}) together imply $\blacksquare_i \mathord{-} y \in {X''}$. The previous two points contradict the fact that ${X''}$ is a proper filter. Thus, ${Y'}=[\mathrm{Y}')$ is a proper filter. Moreover, by (\ref{Vdef}), ${X''}R_i{Y'}$. So we have shown that $\exists{X'}\sqsubseteq{X}$ $\forall{X''}\sqsubseteq{X'}$ $\exists{Y'}\sqsubseteq{Y}$: ${X''}R_i{Y'}$, which establishes \Rref{}. 

For the claim about $\mathbb{A}_\ff$ in part \ref{IsRelGen}: since $\mathbb{A}_\gff$ is strong, $\mathbb{A}_\ff$ is also strong, since they have the same $\sqsubseteq$ and $R_i$ relations; then by Proposition \ref{ROtoRO}, $\mathrm{RO}(S,\sqsubseteq)$ is closed in the ways required for a partial-state frame, which with $\adm_{\ff}=\mathrm{RO}(S,\sqsubseteq)$ means that $\mathbb{A}_\ff$ is a full possibility frame; and then since $\mathbb{A}_\ff$ is full and strong, it is \textit{$R$-tight} by Lemma \ref{TightStrong}.\ref{TightStrong1.75}; and since $\mathbb{A}_\ff$ is full and clearly \textit{separative}, it is \textit{$\sqsubseteq$-tight} by Fact \ref{TightSepDiff}.\ref{TightSepDiff1}.

For part \ref{IsBasedOn}, that $\pi(p)\in \adm_{\gff}$ is immediate from Definition \ref{FiltEx}.\ref{FiltEx4.5}-\ref{FiltEx5}, so $\mathbb{M}_{\gff}$ is based on $\mathbb{A}_{\gff}$, and $P_{\gff}\subseteq P_{\ff}$, so $\mathbb{M}_{\gff}$ is based on $\mathbb{A}_{\ff}$ as well.

For part \ref{BAOtoGenPos3}, the proof is by induction on $\varphi$. The base case is immediate from the definition of $\pi$ in Definition \ref{FiltEx}.\ref{FiltEx5}. The $\neg$ and $\wedge$ cases are also straightforward. For the $\Box_i$ case, if $\mathbb{M}_{\gff},X\nVdash \Box_i\varphi$, then there is a $Y\in S$ such that $XR_iY$ and $\mathbb{M}_{\gff},Y\nVdash \varphi$, which with the inductive hypothesis implies $\tilde{\theta}(\varphi)\not\in Y$, which with $XR_iY$ implies $\blacksquare_i \tilde{\theta}(\varphi)\not\in X$ and hence $\tilde{\theta}(\Box_i\varphi)\not\in X$. In the other direction, suppose $\tilde{\theta}(\Box_i\varphi)\not\in X$. Where
\begin{equation}
\mathrm{Y}=\{y\mid \blacksquare_i y\in X\}\cup \{-\tilde{\theta}(\varphi)\},
\label{Sdef}\end{equation}
suppose for reductio that $[\mathrm{Y})$ is not a proper filter, i.e., $\inc\in [\mathrm{Y})$. Then since $\{y\mid \blacksquare_i y\in X\}$ is a filter, by (\ref{Sdef}) and Fact \ref{FiltGen} there is a $\blacksquare_i y\in X$ such that $y\meet -\tilde{\theta}(\varphi)\leq\inc$, which implies
\begin{equation} \blacksquare_i y\leq \blacksquare_i \tilde{\theta}(\varphi)
\label{GivesBox3}\end{equation}
by the properties of $\blacksquare_i$. Then since $\blacksquare_i y\in X$ and $X$ is a filter,  (\ref{GivesBox3}) implies $\blacksquare_i \tilde{\theta}(\varphi)\in X$ and hence $\tilde{\theta}(\Box_i\varphi)\in X$, contradicting our initial supposition. Hence ${Y}=[\mathrm{Y})$ is a proper filter, which with (\ref{Sdef}) implies $\tilde{\theta}(\varphi)\not\in {Y}$ and hence $\mathbb{M}_{\gff},{Y}\nVdash\varphi$ by the inductive hypothesis. Also by (\ref{Sdef}), we have ${X}R_i{Y}$, which with $\mathbb{M}_{\gff},{Y}\nVdash\varphi$ implies  $\mathbb{M}_{\gff},{X}\nVdash \Box_i\varphi$.

Part \ref{BAOtoGenPos4} is immediate from part \ref{BAOtoGenPos3}: $\mathbb{M}_{\gff},\mathord{\uparrow}x\Vdash \varphi$ iff $\tilde{\theta}(\varphi)\in \mathord{\uparrow}x$ iff $x\leq \tilde{\theta}(\varphi)$.

Part \ref{BAOtoGenPos5} is immediate from parts \ref{IsBasedOn} and \ref{BAOtoGenPos3}.\end{proof}

While Theorem \ref{BAOtoGenPos} shows that satisfiability of formulas is preserved in moving from a BAO to its filter frame or general filter frame, adding to Theorem \ref{BAOtoGenPos} the following obvious lemma shows that unsatisfiability of formulas is also preserved in moving from a BAO to its \textit{general} filter frame. Note, by contrast, that we cannot always turn a possibility model based on $\mathbb{A}_{\ff}$ into an equivalent algebraic model based on $\mathbb{A}$.

\begin{lemma}\label{ObvLem} Given a BAO $\mathbb{A}=\langle A, \meet, -, \top, \{\blacksquare_i\}_{i\in \ind}\rangle$ and a possibility model $\mathcal{M}=\langle \mathbb{A}_{\gff},\pi\rangle$ based on $\mathbb{A}_{\gff}$, define $\mathcal{M}_{-\gff}=\langle \mathbb{A},\pi_{-\gff}\rangle$ with $\pi_{-\gff}\colon \sig\to A$ given by $\pi_{-\gff}(p)=x$ for the $x\in A$ such that $\pi(p)=\{X\in \mathbb{A}_{\gff}\mid x\in X\}$, which must exist by Definition \ref{FiltEx}.\ref{FiltEx4.5} and the fact that $\mathcal{M}$ is based on $\mathbb{A}_{\gff}$. Then $(\mathcal{M}_{-\gff})_{\gff}=\mathcal{M}$.
\end{lemma}

Thus, we arrive at our desired result on the modal equivalence of a BAO and its general filter frame.

\begin{theorem}[Semantic Equivalence of BAOs and General Filter Frames]\label{SatGenFil} For any BAO $\mathbb{A}$ and $\varphi\in\mathcal{L}(\sig,\ind)$, $\varphi$ is satisfiable over $\mathbb{A}$ iff $\varphi$ is satisfiable over $\mathbb{A}_{\gff}$.
\end{theorem}

\begin{proof} The left-to-right direction is given by Theorem \ref{BAOtoGenPos}.\ref{BAOtoGenPos5}. From right to left, if $\varphi$ is satisfied in a possibility model $\mathcal{M}$ based on $\mathbb{A}_\gff$, then by Lemma \ref{ObvLem}, it is satisfied in the possibility model $(\mathcal{M}_{-\gff})_{\gff}$, in which case by Theorem \ref{BAOtoGenPos}.\ref{BAOtoGenPos3} it is satisfied in the algebraic model $\mathcal{M}_{-\gff}$ based on $\mathbb{A}$.
\end{proof}

Not only is every BAO semantically equivalent to its general filter frame, but also homomorphisms between BAOs transform into possibility morphisms between their general filter frames as follows.

\begin{theorem}[From BAO-homomorphisms to Possibility Morphisms II]\label{BAOMorphPossMorphII}
For any BAOs $\mathbb{A}$ and $\mathbb{A}'$ and BAO-homomorphism $h\colon\mathbb{A}'\to\mathbb{A}$, define $h_\gff\colon \mathbb{A}_\gff\to\mathbb{A}'_\gff$ by $h_\gff(X)=h^{-1}[X]$. Then:
\begin{enumerate}
\item\label{BAOMorphPossMorphIIa} $h_\gff$ is a p-morphism;
\item\label{BAOMorphPossMorphIIb} if $h$ is surjective, then $h_\gff$ is a strong embedding;
\item\label{BAOMorphPossMorphIIc} if $h$ is injective, then $h_\gff$ is surjective;
\item\label{BAOMorphPossMorphIId} as a function from $\mathbb{A}_\ff$ to $\mathbb{A}_\ff'$, $h_\gff$ still satisfies parts \ref{BAOMorphPossMorphIIa}-\ref{BAOMorphPossMorphIIc};
\item\label{BAOMorphPossMorphIIe} if $f\colon \mathbb{A}\to\mathbb{A}$ is the identity map on $\mathbb{A}$, then $f_\gff$ is the identity map on $\mathbb{A}_\gff$;
\item\label{BAOMorphPossMorphIIf} for any BAO-homomorphisms $f\colon \mathbb{A}\to\mathbb{B}$ and $g\colon \mathbb{B}\to\mathbb{C}$, $(g\circ f)_\gff = f_\gff \circ g_\gff$.
\end{enumerate}
\end{theorem}
 
\begin{proof} For part \ref{BAOMorphPossMorphIIa}, that $h_\gff$ is a function from $\mathbb{A}_\gff$ to $\mathbb{A}_\gff'$, i.e., sending proper filters from $\mathbb{A}$ to proper filters from $\mathbb{A}'$, follows from the definition of $h_\gff$ and the fact that $h$ is a homomorphism. Where $\mathbb{A}_{\gff}=\langle S,\sqsubseteq,\{R_i\}_{i\in\ind},\adm_{\gff}\rangle$ and $\mathbb{A}_{\gff}'=\langle S',\sqsubseteq',\{R_i'\}_{i\in\ind},\adm_{\gff}'\rangle$, we will show that for all proper filters $X,Y$ from $\mathbb{A}$ and $Y'$ from $\mathbb{A}'$:
\begin{itemize}
\item  if $Y\sqsubseteq  X $, then $h_\gff(Y)\sqsubseteq' h_\gff( X )$ (\SqForth{}); 
\item if $Y'\sqsubseteq' h_\gff( X )$, then $\exists Y$: $Y\sqsubseteq X $ and $h_\gff(Y)= Y'$ (\pSqBack{});
\item if $ X R_i  Y$, then $h_\gff(X)R_i 'h_\gff(Y)$ (\RForth{}); 
\item if $h_\gff( X )R_i' Y'$, then $\exists Y$: $ X R_i  Y$ and $ h_\gff(Y)=Y'$ (\pRBack{});
\item $\forall \mathcal{X}'\in \adm'_\gff$, $h_\gff^{-1}[\mathcal{X}']\in \adm_\gff$ (\textit{pull back}).
\end{itemize}

For \SqForth{}, if $Y\sqsubseteq X$, so $Y\supseteq X$, then $h_\gff(Y)=h^{-1}[Y]\supseteq h^{-1}[X]=h_\gff(X)$, so $h_\gff (Y)\sqsubseteq' h_\gff (X)$. 

For \SqBack{}, suppose $Y'\sqsubseteq' h_\gff( X )$, so $Y'\supseteq h_\gff( X )$. We claim that $[h[Y']\cup X)$ is a proper filter. If not, so $\bot\in [h[Y']\cup X)$, then by Fact \ref{FiltGen} and the facts that $h[Y']$ is closed under finite meets and $X$ is a filter, there are $y\in h[Y']$ and $x\in X$ such that $y\meet x\leq \bot$, so $x\leq - y$ and hence $- y\in X$. Since $y\in h[Y']$, there is a $y'\in Y'$ with $h(y')=y$, so $-h(y')\in X$. Then since $h$ is a homomorphism, $h(-y')\in X$, so $-y'\in h^{-1}[X]=h_\gff (X)$, which with $Y'\supseteq h_\gff (X)$ gives us $-y' \in Y'$, which contradicts the fact that $y'\in Y'$ and $Y'$ is a proper filter. Thus, $[h[Y']\cup X)$ is indeed a proper filter.  Let $Y=[h[Y']\cup X)$. Then since $Y\supseteq X$, we have $Y\sqsubseteq X$. In addition, clearly $h_\gff(Y)=h^{-1}[[h[Y']\cup X)]\supseteq Y'$. Finally, we claim that $h_\gff(Y)\subseteq Y'$.\footnote{Thanks to David Gabelaia and Mamuka Jibladze for pointing out this strengthening of my original proof.} For if $z'\in h_\gff(Y)$, so $h(z')\in [h[Y']\cup X)$, then there is a $y'\in Y'$ and $x\in X$ such that $h(y')\wedge x\leq h(z')$, which implies $x\leq -h(y')\vee h(z')$ and hence $x\leq h(-y'\vee z')$. Then since $x\in X$, we have $h(-y'\vee z')\in X$, so $-y'\vee z'\in h^{-1}[X]=h_\gff(X)$.  Since $Y'\supseteq h_\gff( X )$, it follows that $-y'\vee z'\in Y'$, which with $y'\in Y'$ implies $z'\in Y'$, which completes the proof that $h_\gff(Y)\subseteq Y'$. Thus, $h_\gff(Y)= Y'$, so $h_\gff$ satisfies \pSqBack{}.

For \RForth{}, suppose $XR_iY$. For any $z'\in\mathbb{A}'$, if $\blacksquare_i'z'\in h_\gff(X)=h^{-1}[X]$, then $h(\blacksquare_i'z')\in X$, so $\blacksquare_ih(z')\in X$, which implies $h(z')\in Y$ by $XR_iY$, so $z'\in h^{-1}[Y]=h_\gff(Y)$. Thus, $h_\gff(X)R_i'h_\gff(Y)$.

For \pRBack{}, suppose $h_\gff(X)R_i ' Y'$. Since $X$ is a filter, so is $Z=\{x\mid \blacksquare_i x\in X\}$. We claim that $[h[Y']\cup Z)$ is also a proper filter. If not, then there is a $y'\in Y'$ and $\blacksquare_i x\in X$ such that $h(y')\wedge x\leq \bot$, so $x\leq -h(y')$ and hence $\blacksquare_i x\leq \blacksquare_i \mathord{-} h(y')=h(\blacksquare_i' \mathord{-} y')$, using the fact that $h$ is a BAO-homomorphism. Then since $\blacksquare_i x\in X$, we have $h(\blacksquare_i' \mathord{-}y')\in X$ and hence $\blacksquare_i' \mathord{-}y'\in h_\gff(X)$. Then since $h_\gff(X)R_i ' Y'$, we have $-y'\in Y'$, contradicting the fact that $y'\in Y'$ and $Y'$ is a proper filter. Thus, $[h[Y']\cup Z)$ is a proper filter. Let $Y=[h[Y']\cup Z)$. Then clearly $XR_iY$ and $h_\gff (Y)\supseteq Y'$. Finally, we claim that $h_\gff (Y)\subseteq Y'$. For if $z'\in h_\gff (Y)$, so $h(z')\in [h[Y']\cup Z)$, then there is a $y'\in Y'$ and $\blacksquare_i x\in X$ such that $h(y')\wedge x\leq h(z')$, which implies $x\leq -h(y')\vee h(z')$ and hence $\blacksquare_i x\leq \blacksquare_i (-h(y')\vee h(z'))=h(\blacksquare_i' (-y'\vee z'))$. Then since $\blacksquare_i x\in X$, we have $h(\blacksquare_i' (-y'\vee z'))\in X$, so $\blacksquare_i' (-y'\vee z')\in h_\gff (X)$. Since $h_\gff(X)R_i ' Y'$, it follows that $-y'\vee z'\in Y'$, which with $y'\in Y'$ implies $z'\in Y'$, which completes the proof that $h_\gff (Y)\subseteq Y'$. Thus, $h_\gff (Y)= Y'$, so $h_\gff$ satisfies \pRBack{}.

For \textit{pull back}, recall that $\adm_\gff=\{\widehat{x}\mid x\in A\}$ and $\adm'_\gff=\{\widehat{x'}\mid x'\in A'\}$. Now suppose $\widehat{x'}\in \adm'_\gff$, so $x'\in A'$. Then $h(x')\in A$, so $\widehat{h(x')}\in P_\gff$. We claim that $h^{-1}_\gff[\widehat{x'}]=\widehat{h(x')}$, so $h^{-1}_\gff[\widehat{x'}]\in\adm_\gff$, as desired. For $h^{-1}_\gff[\widehat{x'}]$ is the set of proper filters $X$ in $\mathbb{A}$ such that $h_\gff(X)=h^{-1}[X]\in\widehat{x'}$, which means $x'\in h^{-1}[X]$, which means $h(x')\in X$, and $\widehat{h(x')}$ is the set of proper filters $X$ in $\mathbb{A}$ such that $h(x')\in X$.

For part \ref{BAOMorphPossMorphIIb}, we must show that if $h$ is surjective, then (i) $h_\gff(Y)\sqsubseteq' h_\gff(X)$ implies $Y\sqsubseteq X$ (so $h_\gff$ is injective), (ii) $h_\gff(X)R_i'h_\gff(Y)$ implies $XR_iY$, and (iii) for every $\widehat{x}\in\adm_\gff$ there is a $\widehat{x'}\in\adm_\gff'$ such that $h_\gff[\widehat{x}]=h_\gff[S]\cap\widehat{x'}$. For (i), assume $h_\gff(Y)\sqsubseteq' h_\gff(X)$, so $h^{-1}[Y]\supseteq h^{-1}[X]$. Since $h$ is surjective, for any $x\in X$ there is an $x'\in A'$ such that $h(x')=x$, so $x'\in h^{-1}[X]$ and hence $x'\in h^{-1}[Y]$ by our assumption, so $h(x')=x\in Y$. Thus, $Y\supseteq X$, so $Y\sqsubseteq X$. For (ii), to show $XR_iY$, we must show that for all $\blacksquare_i x\in X$, $x\in Y$. Given $\blacksquare_i x\in X$, since $h$ is surjective, there is an $x'\in A'$ such that $h(x')=x$, and then since $h$ is a BAO-homomorphism, $\blacksquare_i x=\blacksquare_i h(x')=h(\blacksquare_i' x')$. Hence $h(\blacksquare_i'x')\in X$, so $\blacksquare_i'x'\in h^{-1}[X]=h_\gff(X)$, which with $h_\gff(X)R_i'h_\gff(Y)$ implies $x'\in h_\gff(Y)=h^{-1}[Y]$, so $h(x')=x\in Y$. Thus, $XR_iY$. For (iii), since $h$ is surjective, given $\widehat{x}\in \adm_\gff$ there is an $x'\in A'$, so $\widehat{x'}\in \adm_\gff'$, such that $h(x')=x$, so $h_\gff[\widehat{x}]=h_\gff[\widehat{h(x')}]$. Now we claim that $h_\gff[\widehat{h(x')}]=h_\gff[S]\cap\widehat{x'}$. For the left-to-right inclusion, suppose $X'\in h_\gff[\widehat{h(x')}]$, so there is an $X\in S$ such that $X\in\widehat{h(x')}$, so $h(x')\in X$, and $h_\gff(X)=h^{-1}[X]=X'$. It follows that $x'\in X'$, so $X'\in\widehat{x'}$. For the right-to-left inclusion, suppose that $X'\in\widehat{x'}$, so $X'\in S'$ and $x'\in X'$, and that $X'\in h_\gff[S]$, so there is an $X\in S$ such that $h_\gff(X)=h^{-1}[X]=X'$. Hence $h(x')\in X$, so $X\in \widehat{h(x')}$, which with $h_\gff(X)=X'$ implies $X'\in h_\gff[\widehat{h(x')}]$. 

For part \ref{BAOMorphPossMorphIIc}, assuming $h$ is injective, for any $X'\in\mathbb{A}'_\gff$, $h_\gff\big([h[X'])\big)=h^{-1}\big[[h[X'])\big]=X'$, so $h_\gff$ is surjective.
 
For part \ref{BAOMorphPossMorphIId}, to show that $h_\gff$ is a p-morphism from $\mathbb{A}_\ff$ to $\mathbb{A}_\ff'$, we need only add to the proof of part \ref{BAOMorphPossMorphIIa} above that $h_\gff$ satisfies \PullBack{} with respect to $\mathbb{A}_\ff$ and $\mathbb{A}'_\ff$: for all $\mathcal{X}'\in \adm_\ff'=\mathrm{RO}(S',\sqsubseteq')$, we have $h_\gff^{-1}[\mathcal{X}']\in \adm_\ff=\mathrm{RO}(S,\sqsubseteq)$. To show that $h_\gff^{-1}[\mathcal{X}']$ satisfies \textit{persistence} with respect to $\sqsubseteq$, suppose $X\in h_\gff^{-1}[\mathcal{X}']$, so $h_\gff(X)\in\mathcal{X}'$, and $Y\sqsubseteq X$. Then $h_\gff(Y)\sqsubseteq' h_\gff(X)$ by \SqForth{}, which with $h_\gff(X)\in\mathcal{X}'$ and the \textit{persistence} of $\mathcal{X}'$ with respect to $\sqsubseteq'$ implies $h_\gff(Y)\in\mathcal{X}'$, so $Y\in h_\gff^{-1}[\mathcal{X}']$. To show that $h_\gff^{-1}[\mathcal{X}']$ satisfies \textit{refinability} with respect to $\sqsubseteq$, observe that if $X\not\in h_\gff^{-1}[\mathcal{X}']$, so $h_\gff(X)\not\in\mathcal{X}'$, then by \textit{refinability} for $\mathcal{X}'$ with respect to $\sqsubseteq'$, there is a $Y'\sqsubseteq' h_\gff(X)$ such that (a) for all $Y''\sqsubseteq' Y'$, $Y''\not\in \mathcal{X}'$. Given $Y'\sqsubseteq' h_\gff(X)$ and \SqBack{}, there is a $Y\sqsubseteq X$ such that $h_\gff(Y)\sqsubseteq ' Y'$. Then for any $Z\sqsubseteq Y$, we have $h_\gff(Z)\sqsubseteq' h_\gff(Y)\sqsubseteq' Y'$ by \SqForth{}, so $h_\gff(Z)\not\in\mathcal{X}'$ by (a), so $Z\not\in h_\gff^{-1}[\mathcal{X}']$. Thus, we have shown that if $X\not\in h_\gff^{-1}[\mathcal{X}']$, then there is a $Y\sqsubseteq X$ such that for all $Z\sqsubseteq Y$, $Z\not\in h_\gff^{-1}[\mathcal{X}']$, so $h_\gff^{-1}[\mathcal{X}']$ satisfies \textit{refinability} with respect to $\sqsubseteq$.

Next, to show that if $h$ is surjective, then $h_\gff$ is a strong embedding from $\mathbb{A}_\ff$ to $\mathbb{A}_\ff'$, we need only add to the proof of part \ref{BAOMorphPossMorphIIb} above that for all $\mathcal{X}\in \adm_\ff=\mathrm{RO}(S,\sqsubseteq)$, there is an $\mathcal{X}'\in\adm_\ff'=\mathrm{RO}(S',\sqsubseteq')$ such that $h_\gff[\mathcal{X}]=h_\gff[S]\cap \mathcal{X}'$. Where $\mathord{\Downarrow}h_\gff[\mathcal{X}]=\{Y'\in S'\mid \exists X'\in h_\gff[\mathcal{X}]\colon Y'\sqsubseteq' X'\}$, let $\mathcal{X}'=\mathrm{int}(\mathrm{cl}(\mathord{\Downarrow}h_\gff[\mathcal{X}]))$, recalling that $\mathrm{int}(\mathcal{Y}')=\{Y'\in S'\mid \forall Z'\sqsubseteq' Y',\, Z'\in \mathcal{Y}'\}$ and $\mathrm{cl}(\mathcal{Y}')=\{Y'\in S'\mid \exists Z'\sqsubseteq' Y'\colon Z'\in \mathcal{Y}'\}$. Then by Fact \ref{RefReg}.\ref{RefReg2.5}, $\mathcal{X}'\in \adm'_\ff$ and $h_\gff[\mathcal{X}]\subseteq h_\gff[S]\cap \mathcal{X}'$. To show $h_\gff[S]\cap \mathcal{X}'\subseteq h_\gff[\mathcal{X}]$, suppose $X'\in h_\gff[S]$ but $X'\not\in h_\gff[\mathcal{X}]$.  Since $X'\in h_\gff[S]$, there is an $X\in S$ such that $h_\gff(X)=X'$, which with $X'\not\in h_\gff[\mathcal{X}]$ implies $X\not\in\mathcal{X}$. Then by \textit{refinability} for $\mathcal{X}$ with respect to $\sqsubseteq$, there is a $Y\sqsubseteq X$ such that (b) for all $Z\sqsubseteq Y$, $Z\not\in\mathcal{X}$. By \SqForth{}, $Y\sqsubseteq X$ implies $h_\gff(Y)\sqsubseteq' h_\gff(X)=X'$.
Now for any $Z'\sqsubseteq' h_\gff(Y)$, by \SqBack{} there is a $Z\sqsubseteq Y$ with $h_\gff(Z)\sqsubseteq' Z'$.  Then by (b), $Z\not\in\mathcal{X}$. If $Z'\in \mathord{\Downarrow}h_\gff[\mathcal{X}]$, so there is a $V\in \mathcal{X}$ such that $Z'\sqsubseteq' h_\gff(V)$, then from $h_\gff(Z)\sqsubseteq' Z'$ above we have $h_\gff(Z)\sqsubseteq' h_\gff(V)$, which with (i) above implies $Z\sqsubseteq V$, which with $V\in\mathcal{X}$ and \textit{persistence} for $\mathcal{X}$ implies $Z\in\mathcal{X}$, contradicting what we just deduced from (b). Thus, $Z'\not\in \mathord{\Downarrow}h_\gff[\mathcal{X}]$. Then since $Z'$ was an arbitrary refinement of $h_\gff(Y)$, we have $h_\gff(Y)\not\in\mathrm{cl}(\mathord{\Downarrow}h_\gff[\mathcal{X}])$, which with $h_\gff(Y)\sqsubseteq' X'$ from above implies $X'\not\in \mathrm{int}(\mathrm{cl}(\mathord{\Downarrow}h_\gff[\mathcal{X}]))=\mathcal{X}'$. This shows that $h_\gff[S]\cap \mathcal{X}'\subseteq h_\gff[\mathcal{X}]$.
       
 Finally, since the domains of $\mathbb{A}'_\gff$ and $\mathbb{A}'_\ff$ are the same, part \ref{BAOMorphPossMorphIIc} implies that $h_\gff$ is onto $\mathbb{A}'_\ff$ if $h$ is injective.

Parts \ref{BAOMorphPossMorphIIe}-\ref{BAOMorphPossMorphIIf} are easy to check.\end{proof} 

From Theorems \ref{BAOtoGenPos}.\ref{IsRelGen} and \ref{BAOMorphPossMorphII}, we obtain the next piece of our categorical picture.

\begin{corollary}[The $(\cdot)_\gff$ Functor]\label{GffFunc} The $(\cdot)_\gff$ operation given in Definition \ref{FiltEx} and Theorem \ref{BAOMorphPossMorphII} is a contravariant functor from \textbf{BAO} to the category of possibility frames with p-morphisms. Thus, together $(\cdot)_\gff$ and $(\cdot)^\under$ from Corollary \ref{FuncCor1} form a pair of contravariant functors between these categories.
\end{corollary}

We will be more specific about the type of possibility frames in the image of $(\cdot)_\gff$ in the next section.

\subsection{\texorpdfstring{$(\cdot)_\gff$ and $(\cdot)^\under$}{(.)\_g and (.)\^{}b}, and Dual Equivalence with Filter-Descriptive Frames}\label{Fdes} 

Let us now consider the relation between the functor $(\cdot)_\gff$ from \S~\ref{GFPF} and the functor $(\cdot)^\under$ from \S~\ref{PossToBAO}.

\begin{proposition}[From BAOs to Frames and Back II]\label{BAOsFramesII} Given a BAO $\mathbb{A}$, define $\eta_\mathbb{A}\colon \mathbb{A}\to(\mathbb{A}_\gff)^\under$ by $\eta_\mathbb{A}(x)=\widehat{x}$, where $\widehat{x}$ is the set of proper filters in $\mathbb{A}$ that contain $x$, as in Definition \ref{FiltEx}. Then:
\begin{enumerate}
\item\label{BAOsFramesIIa} $\eta_\mathbb{A}$ is a BAO-isomorphism; 
\item\label{BAOsFramesIIb} if $g\colon \mathbb{A}\to\mathbb{B}$ is a BAO-homomorphism, then $(g_\gff)^\under\circ \eta_\mathbb{A}= \eta_\mathbb{B}\circ g$, so the following diagram commutes:
\begin{center}
\begin{tikzpicture}[->,>=stealth',shorten >=1pt,shorten <=1pt, auto,node
distance=2cm,thick,every loop/.style={<-,shorten <=1pt}]
\tikzstyle{every state}=[fill=gray!20,draw=none,text=black]

\node (A) at (0,2) {{$\mathbb{A}$}};
\node (A') at (0,0) {{$(\mathbb{A}_\gff)^\under$}};
\node (B) at (3,2) {{$\mathbb{B}$}};
\node (B') at (3,0) {{$(\mathbb{B}_\gff)^\under$}};

\path (A) edge[->] node {{$g$}} (B);
\path (A') edge[<-] node {{$\eta_\mathbb{A}$}} (A);
\path (B) edge[->] node {{$\eta_\mathbb{B}$}} (B');
\path (B') edge[<-] node {{$(g_\gff)^\under$}} (A');

\end{tikzpicture}
\end{center}
\end{enumerate}
\end{proposition}
\begin{proof} For part \ref{BAOsFramesIIa}, to see that $\eta_\mathbb{A}$ is order-reflecting and hence injective, observe that  if $x\not\leq^\mathbb{A}y$, then $x\in\mathord{\uparrow}x$ but $y\not\in\mathord{\uparrow}x$, so $\widehat{x}\not\subseteq\widehat{y}$. To see that $\eta_\mathbb{A}$ is surjective, recall that by Definition \ref{FiltEx}, the set $P_\gff$ of admissible propositions in $\mathbb{A}_\gff$ is $\{\widehat{x}\mid x\in\mathbb{A}\}$, and by Definition \ref{RegOpAlg}, $P_\gff$ is the domain of $(\mathbb{A}_\gff)^\under$. Then since $x\leq^\mathbb{A}y$ implies $\widehat{x}\subseteq\widehat{y}$, it follows that $\eta_\mathbb{A}$ is an order-isomorphism and hence a Boolean isomorphism. Finally, for the modal operations, we have: 
\begin{eqnarray}
\eta_\mathbb{A}(\blacksquare_i^\mathbb{A}x)&=&\widehat{\blacksquare_i^\mathbb{A}x}\nonumber\\
&=& \{Y\in\mathbb{A}_\gff\mid \blacksquare_i^\mathbb{A}x\in Y\}\nonumber\\
&=& \{Y\in\mathbb{A}_\gff\mid R_i^{\mathbb{A}_\gff}(Y)\subseteq\widehat{x}\}\label{AlreadyShowed}\\
&=& \blacksquare_i^{(\mathbb{A}_\gff)^\under} \widehat{x} = \blacksquare_i^{(\mathbb{A}_\gff)^\under} \eta_\mathbb{A}(x).\nonumber
\end{eqnarray}
The proof of (\ref{AlreadyShowed}) is essentially the same as the $\Box_i$ case of the proof of Theorem \ref{BAOtoGenPos}.\ref{BAOtoGenPos3}. 
 
For part \ref{BAOsFramesIIb}, given $g\colon \mathbb{A}\to\mathbb{B}$, $g_\gff\colon \mathbb{B}_\gff\to\mathbb{A}_\gff$, $(g_\gff)^\under\colon (\mathbb{A}_\gff)^\under\to (\mathbb{B}_\gff)^\under$, $\eta_\mathbb{A}\colon \mathbb{A}\to(\mathbb{A}_\gff)^\under$, and $\eta_\mathbb{B}\colon \mathbb{B}\to(\mathbb{B}_\gff)^\under$, we have: 
\[\begin{array}{llll}
(g_\gff)^\under (\eta_\mathbb{A}(x)) & = & (g_\gff)^\under (\widehat{x}^\mathbb{A})&\mbox{by definition of $\eta_\mathbb{A}$}\\
&=& (g_\gff)^{-1}[\widehat{x}^\mathbb{A}]=\{X\in \mathbb{B}_\gff \mid g_\gff(X)\in \widehat{x}^\mathbb{A}\} & \mbox{by definition of $(\cdot)^\under$} \\
 &=& \{X\in\mathbb{B}_\gff\mid g^{-1}[X]\in\widehat{x}^\mathbb{A}\} & \mbox{by definition of $(\cdot)_\gff$}\\
&=& \widehat{g(x)}^\mathbb{B} &(\triangle)\\
&=&\eta_\mathbb{B}(g(x)) &  \mbox{by definition of $\eta_\mathbb{B}$.} 
\end{array}\]
For $(\triangle)$, if $X\in \widehat{g(x)}^\mathbb{B}$, then $X$ is a proper filter in $\mathbb{B}$, so $X\in\mathbb{B}_\gff$, and $g(x)\in X$, so $x\in g^{-1}[X]$. Then since $g$ is a homomorphism from $\mathbb{A}$ to $\mathbb{B}$, that $X$ is a proper filter in $\mathbb{B}$ implies that $g^{-1}[X]$ is a proper filter in $\mathbb{A}$, so $x\in g^{-1}[X]$ implies $g^{-1}[X]\in\widehat{x}^\mathbb{A}$. In the other direction, if $g^{-1}[X]\in\widehat{x}^\mathbb{A}$, then $x\in g^{-1}[X]$, so $g(x)\in X$, and since $X\in \mathbb{B}_\gff$, $X$ is a proper filter in $\mathbb{B}$, so $g(x)\in X$ implies $X\in \widehat{g(x)}^\mathbb{B}$.
\end{proof}

Let us now go in the other direction, from a frame $\mathcal{F}$ to $(\mathcal{F}^\under)_\gff$. Later we will also consider $(\mathcal{F}^\under)_\ff$.

\begin{definition}[Filter Extension and General Filter Extension]\label{FiltExt} For a possibility frame $\mathcal{F}$, its \textit{filter extension} is the possibility frame $(\mathcal{F}^\under)_\ff$, and its \textit{general filter extension} is the possibility frame $(\mathcal{F}^\under)_\gff$. \hfill $\triangleleft$
\end{definition}

In contrast to Proposition \ref{BAOsFramesII}, it is clear that many possibility frames $\mathcal{F}$ will not be isomorphic to $(\mathcal{F}^\under)_\gff$. The same point applies in the case of taking the general \textit{ultrafilter} frame of the underlying BAO of a \textit{world} frame (see \S~\ref{AlgSem}), which is isomorphic to the original world frame iff the original frame is \textit{descriptive} (see \citealt[Thm.~5.76]{Blackburn2001}). Goldblatt \citeyearpar[p.~33f]{Goldblatt1974} originally defined a \textit{descriptive} frame to be a world frame $\mathcal{F}$ that is differentiated (Axiom I), tight (Axiom II), and such that for every ultrafilter $u$ in the underlying BAO of $\mathcal{F}$, there is a world $w$ in $\mathcal{F}$ such that $u$ is the set of admissible propositions from $\mathcal{F}$ that contain $w$ (Axiom III). Descriptive world frames are a special case of possibility frames (Example \ref{KripkeAgain}). But we would like a notion analogous to descriptive so that the general \textit{filter} frame of a BAO will qualify. 

\begin{definition}[Filter-Descriptive]\label{F-desc} A possibility frame $\mathcal{F}=\langle S,\sqsubseteq,\{R_i\}_{i\in\ind},\adm\rangle$ is \textit{filter-descriptive} iff it is tight and for every proper filter $F$ in $\mathcal{F}^\under$, there is an $x\in S$ such that $F=P(x)=\{X\in \adm\mid x\in X\}$. \hfill $\triangleleft$
\end{definition}
\noindent Recall that \textit{tight} implies \textit{$\sqsubseteq$-tight} (Definition \ref{TightFrames}), which implies \textit{differentiated} (Fact \ref{TightSepDiff}.\ref{TightSepDiff2}).

The following two propositions show that \textit{filter-descriptive} is indeed the notion we want.

\begin{proposition}[Filter-Descriptive Frames and BAOs]\label{GenAb}  For any BAO $\mathbb{A}$ and possibility frame $\mathcal{F}$:
\begin{enumerate}
\item\label{GenAb1} $\mathbb{A}_\gff$ is filter-descriptive;
\item\label{GenAb2} $\mathcal{F}$ is isomorphic to $(\mathcal{F}^\under)_\gff$ iff $\mathcal{F}$ is filter-descriptive.
\end{enumerate}
\end{proposition}

\begin{proof} For part \ref{GenAb1}, we have already shown for Theorem \ref{BAOtoGenPos}.\ref{IsRelGen} that $\mathbb{A}_\gff$ is tight. Let us show that $\mathbb{A}_\gff$ satisfies the condition about filters. Consider a proper filter $F$ in $(\mathbb{A}_\gff)^\under$. Since the domain of $(\mathbb{A}_\gff)^\under$ is the set $P_\gff$ of admissible propositions in $\mathbb{A}_\gff$, we have $F\subseteq P_\gff=\{\widehat{x}\mid x\in \mathbb{A}\}$. Recall that $\widehat{x}$ is the set of all proper filters in $\mathbb{A}$ that contain $x$. Now let $Z=\{x\in\mathbb{A}\mid \widehat{x}\in F\}$. We claim that $Z$ is a proper filter in $\mathbb{A}$, so $Z\in\mathbb{A}_\gff$. Clearly $Z\neq\emptyset$. Suppose $x\in Z$, so $\widehat{x}\in F$, and $x\leq_\mathbb{A} y$. Since $\widehat{x}$ and $\widehat{y}$ are the sets of proper filters in $\mathbb{A}$ containing $x$ and $y$, respectively, $x\leq_\mathbb{A} y$ implies that any proper filter that contains $x$ also contains $y$, so $\widehat{x}\subseteq\widehat{y}$ and hence $\widehat{x}\leq_{(\mathbb{A}_\gff)^\under} \widehat{y}$. Then since $\widehat{x}\in F$ and $F$ is a filter in $(\mathbb{A}_\gff)^\under$, we have $\widehat{y}\in F$ and hence $y\in Z$. Thus, $Z$ is an upset in $\mathbb{A}$. Next, if $x,y\in Z$, then $\widehat{x},\widehat{y}\in F$, so $\widehat{x}\cap\widehat{y}\in F$ because $F$ is a filter in $(\mathbb{A}_\gff)^\under$, so $\widehat{x\meet y}\in F$ by the fact that $\widehat{x}\cap\widehat{y}=\widehat{x\meet y}$ (recall the proof of Theorem \ref{BAOtoGenPos}), so $x\meet y\in Z$. Thus, $Z$ is also downward directed. So $Z$ is a filter. Moreover, if $\bot\in Z$, then $\widehat{\bot}=\emptyset\in F$, in which case $F$ would not be a proper filter in $(\mathbb{A}_\gff)^\under$. Thus, $Z$ is a proper filter in $\mathbb{A}$. Finally, by our definitions, $P_\gff(Z)=\{\widehat{x}\mid Z\in\widehat{x}\}=\{\widehat{x}\mid x\in Z\}=F$.

For part \ref{GenAb2}, the left-to-right direction follows from part \ref{GenAb1}. The right-to-left direction follows from the following Proposition \ref{FramesBAOsII}. \end{proof}

Note that since every normal modal logic is sound and complete with respect to a BAO (Theorem \ref{AdAlg}), it follows from Theorem \ref{SatGenFil} together with Proposition \ref{GenAb}.\ref{GenAb1} that every consistent normal modal logic is sound and complete with respect to a filter-descriptive possibility frame. We will return to this point in \S~\ref{Canonical}.

\begin{proposition}[From Frames to BAOs and Back II]\label{FramesBAOsII} For any possibility frame $\mathcal{F}=\langle S,\sqsubseteq, \{R_i\}_{i\in\ind},\adm\rangle$, define $\eta_\mathcal{F}\colon \mathcal{F}\to (\mathcal{F}^\under)_\gff$ by $\eta_\mathcal{F}(x)=P(x)$, where $P(x)=\{X\in \adm\mid x\in X\}$. Then:
\begin{enumerate}
\item\label{FramesBAOsII0} $\eta_\mathcal{F}$ is a possibility morphism;
\item\label{FramesBAOsII1} if $\mathcal{F}$ is filter-descriptive, then $\eta_\mathcal{F}$ is an isomorphism;
\item\label{FramesBAOsII2} if $\mathcal{H}$ is a possibility frame and $g\colon \mathcal{F}\to\mathcal{H}$ is a possibility morphism, then $(g^\under)_\gff\circ \eta_\mathcal{F}=\eta_\mathcal{H}\circ g$, so the following diagram commutes:
\begin{center}
\begin{tikzpicture}[->,>=stealth',shorten >=1pt,shorten <=1pt, auto,node
distance=2cm,thick,every loop/.style={<-,shorten <=1pt}]
\tikzstyle{every state}=[fill=gray!20,draw=none,text=black]

\node (F) at (0,2) {{$\mathcal{F}$}};
\node (F') at (0,0) {{$(\mathcal{F}^\under)_\gff$}};
\node (G) at (3,2) {{$\mathcal{H}$}};
\node (G') at (3,0) {{$(\mathcal{H}^\under)_\gff$}};

\path (F) edge[->] node {{$g$}} (G);
\path (F') edge[<-] node {{$\eta_\mathcal{F}$}} (F);
\path (G) edge[->] node {{$\eta_\mathcal{H}$}} (G');
\path (G') edge[<-] node {{$(g^\under)_\gff$}} (F');

\end{tikzpicture}
\end{center}
\end{enumerate}
\end{proposition}

\begin{proof} Recall that the set of proper filters in $\mathcal{F}^\under$ is the domain of $(\mathcal{F}^\under)_\gff$, and $\adm_\gff=\{\widehat{X}\mid X\in \mathcal{F}^\under\}$ is the set of admissible sets in $(\mathcal{F}^\under)_\gff$, where $\widehat{X}$ is the set of proper filters in $\mathcal{F}^\under$ that contain $X$. Let $\sqsubseteq'$ and $R_i'$ be the refinement and accessibility relations in $(\mathcal{F}^\under)_\gff$.

For part \ref{FramesBAOsII0}, for every $x\in S$, $P(x)$ is a proper filter in $\mathcal{F}^\under$, so $\eta_\mathcal{F}$ is indeed a map from $\mathcal{F}$ to $(\mathcal{F}^\under)_\gff$. 

For \SqMatch{}, we must show that for all $x\in S$ and all $\widehat{X}\in P_\gff$, $\mathord{\downarrow}'\eta_\mathcal{F}(x)\cap \widehat{X}=\emptyset$ iff $\mathord{\downarrow}x\cap \eta_\mathcal{F}^{-1}[\widehat{X}]=\emptyset$. From left-to-right, suppose there is a $y\in \mathord{\downarrow}x\cap \eta_\mathcal{F}^{-1}[\widehat{X}]$, which means $y\sqsubseteq x$ and $P(y)\in \widehat{X}$. Since $y\sqsubseteq x$, we have $P(y)\supseteq P(x)$ by \textit{persistence}. Thus, $\eta_\mathcal{F}(y)\in \mathord{\downarrow}'\eta_\mathcal{F}(x)\cap \widehat{X}$. Conversely, suppose there is a proper filter $F\in \mathord{\downarrow}'\eta_\mathcal{F}(x)\cap \widehat{X}$, which means $F\supseteq P(x)$ and $X\in F$. Then $\neg X\not\in P(x)$, so $x\not\in\neg X$, which implies that there is a $y\sqsubseteq x$ such that $y\in X$, so $X\in P(y)$ and hence $P(y)\in \widehat{X}$. Therefore, $y\in \mathord{\downarrow}x\cap \eta_\mathcal{F}^{-1}[\widehat{X}]$. 

For \RMatch{}, we must show that for all $x\in S$ and $\widehat{X}\in\adm_\gff$, $R_i'(\eta_\mathcal{F}(x))\subseteq \widehat{X}$ iff $R_i(x)\subseteq \eta_\mathcal{F}^{-1}[\widehat{X}]$. From left-to-right,  suppose there is a $y\in S$ such that $xR_iy$ and $y\not\in \eta_\mathcal{F}^{-1}[\widehat{X}]$. From $xR_iy$, it follows by the definition of $R_i'$ in $(\mathcal{F}^\under)_\gff$ that $P(x)R_i'P(y)$. From $y\not\in \eta_\mathcal{F}^{-1}[\widehat{X}]$, we have $P(y)\not\in\widehat{X}$. Thus, $R_i'(\eta_\mathcal{F}(x))\not\subseteq \widehat{X}$. Conversely, suppose there is a proper filter $F$ in $\mathcal{F}^\under$ such that $P(x)R_i'F$ but $F\not\in\widehat{X}$, so $X\not\in F$. Then by the definition of $R_i'$, we have $\blacksquare_i X\not\in P(x)$, so $x\not\in \blacksquare_i X$, which implies that there is a $y\in S$ such that $xR_iy$ and $y\not\in X$, so $X\not\in P(y)$ and hence $P(y)\not\in\widehat{X}$. Thus, $R_i(x)\not\subseteq \eta_\mathcal{F}^{-1}[\widehat{X}]$.

Finally, for \PullBack{}, we must show that for any $\widehat{X}\in \adm_\gff$, we have $\eta_\mathcal{F}^{-1}[\widehat{X}]\in \adm$.  Since $X\in \adm$, it suffices to show that $\eta_\mathcal{F}^{-1}[\widehat{X}]=X$. If $x\in \eta_\mathcal{F}^{-1}[\widehat{X}]$, then $P(x)\in \widehat{X}$, so $X\in P(x)$ and hence $x\in X$. Conversely, if $x\in X$, so $X\in P(x)$, then since $P(x)$ is a filter in $\mathcal{F}^\under$, we have $P(x)\in \widehat{X}$ and hence $x\in \eta_\mathcal{F}^{-1}[\widehat{X}]$.

For part \ref{FramesBAOsII1}, since $\mathcal{F}$ is \textit{filter-descriptive}, every proper filter in $\mathcal{F}^\under$ is $P(x)=\eta_\mathcal{F}(x)$ for some $x\in S$, so $\eta_\mathcal{F}$ is a surjective map from $\mathcal{F}$ onto $(\mathcal{F}^\under)_\gff$. That $\eta_\mathcal{F}$ is injective follows from the assumption that $\mathcal{F}$ is \textit{$\sqsubseteq$-tight} and hence \textit{differentiated} (Fact \ref{TightSepDiff}.\ref{TightSepDiff2}). The equivalence of $xR_i y$ and $P(x)R_i'P(y)$ follows from the assumption that $\mathcal{F}$ is \textit{$R$-tight}. Next, by definition of $(\mathcal{F}^\under)_\gff$, $P(x)\sqsubseteq ' P(y)$ iff $P(x)\supseteq P(y)$, and the equivalence of $x\sqsubseteq y$ and $P(x)\supseteq P(y)$ follows from the assumption that $\mathcal{F}$ is \textit{$\sqsubseteq$-tight}. The final part of establishing an isomorphism is to show that for all $X\in \adm$, $\eta_\mathcal{F}[X]\in \adm_\gff$.  It suffices to show that for all $X\in\adm$, $\eta_\mathcal{F}[X]=\widehat{X}$. Since $\eta_\mathcal{F}[X]=\{P(x)\mid x\in X\}$, if $F\in \eta_\mathcal{F}[X]$, then $F=P(x)$ for some $x\in X$, so $X\in P(x)=F$ and hence $F\in\widehat{X}$; in the other direction, if $F\in\widehat{X}$, so $F$ is a proper filter in $\mathcal{F}^\under$ such that $X\in F$, then since $\mathcal{F}$ is filter-descriptive, there is some $x\in S$ such that $F=P(x)$, and then $X\in F$ implies $x\in X$, so $F\in \eta_\mathcal{F}[X]$.

For part \ref{FramesBAOsII2}, given $g\colon \mathcal{F}\to\mathcal{H}$, $g^\under\colon \mathcal{H}^\under\to \mathcal{F}^\under$, and $(g^\under)_\gff\colon (\mathcal{F}^\under)_\gff\to(\mathcal{H}^\under)_\gff$, $\eta_\mathcal{F}\colon \mathcal{F}\to(\mathcal{F}^\under)_\gff$, and $\eta_\mathcal{H}\colon \mathcal{H}\to(\mathcal{H}^\under)_\gff$, we have:
\[\begin{array}{llll}
(g^\under)_\gff(\eta_\mathcal{F}(x))&=& (g^\under)_\gff (\adm^\mathcal{F}(x))&\mbox{by definition of $\eta_\mathcal{F}$}\\
&=& (g^\under)^{-1}[\adm^\mathcal{F}(x)]= \{X\in \mathcal{H}^\under \mid  g^\under(X)\in \adm^\mathcal{F}(x)\}&\mbox{by definition of $(\cdot)_\gff$}\\
&=& \{X\in \mathcal{H}^\under \mid  g^{-1}[X]\in \adm^\mathcal{F}(x)\}&\mbox{by definition of $(\cdot)^\under$}\\
&=& \adm^\mathcal{H}(g(x)) & (*)\\
&=& \eta_\mathcal{H}(g(x)) &\mbox{by definition of $\eta_\mathcal{H}$.}
\end{array}\]
For $(*)$, if $X\in \adm^\mathcal{H}(g(x))$, then $X\in \adm^\mathcal{H}$, so $X\in\mathcal{H}^\under$, and by \textit{pull back} for $g$, $g^{-1}[X]\in \adm^\mathcal{F}$. From $X\in \adm^\mathcal{H}(g(x))$ we also have $g(x)\in X$,  so $x\in g^{-1}[X]$, which with $g^{-1}[X]\in \adm^\mathcal{F}$ implies $g^{-1}[X]\in \adm^\mathcal{F}(x)$. In the other direction, if $g^{-1}[X]\in \adm^\mathcal{F}(x)$, then $x\in g^{-1}[X]$, so $g(x)\in X$. Then since $X\in \mathcal{H}^\under$, $X\in \adm^\mathcal{H}$, so $g(x)\in X$ implies $X\in \adm^\mathcal{H}(g(x))$. This completes the proof of $(*)$.
\end{proof}

Putting together Propositions \ref{BAOsFramesII} and \ref{FramesBAOsII} and Corollaries \ref{FuncCor1} and \ref{GffFunc}, we have the following analogue of Goldblatt's \citeyearpar{Goldblatt1974} dual equivalence result for the categories of BAOs with BAO-homomorphism and of descriptive frames with p-morphisms. \textbf{FiltPoss} is the category of filter-descriptive possibility frames with p-morphisms.

\begin{theorem}[Dual Equivalence II]\label{Dual2} \textbf{BAO} is dually equivalent to \textbf{FiltPoss}.
\end{theorem}

\subsection{Reflection with Filter-Descriptive Frames}\label{FilterReflectionSection}

Next we observe that \textbf{FiltPoss} is a full subcategory of the category \textbf{Poss} of possibility frames with possibility morphisms.

\begin{proposition}[Full Subcategory II]\label{PossToP2} Every possibility morphism between filter-descriptive possibility frames is a p-morphism.
\end{proposition}
\begin{proof} Immediate from Theorem \ref{BAOMorphPossMorphII}.\ref{BAOMorphPossMorphIIa} and Proposition \ref{FramesBAOsII}.
\end{proof}

In fact, we will show that \textbf{FiltPoss} is a reflective subcategory of \textbf{Poss}. This is analogous to Goldblatt's \citeyearpar{Goldblatt2006b} result that the category of descriptive frames with p-morphisms is a reflective subcategory of the category of general world frames with what he calls \textit{modal maps}.

Before proving the reflective subcategory claim, let us prove a preliminary lemma.

\begin{lemma}[Filter Transfer]\label{FiltTransfer} For any possibility frame $\mathcal{F}$, filter-descriptive possibility frame $\mathcal{G}$, possibility morphism $g\colon \mathcal{F}\to\mathcal{G}$, and proper filter $F$ in $\mathcal{F}^\under$, the set $\{X\in \mathcal{G}^\under\mid g^{-1}[X]\in F\}$ is a proper filter in $\mathcal{G}^\under$.
\end{lemma}
\begin{proof} Clearly $\{X\in \mathcal{G}^\under\mid g^{-1}[X]\in F\}\neq\emptyset$. Next we show it is closed under taking supersets in $\mathcal{G}^\under$. Suppose $X\in \mathcal{G}^\under$ and  $g^{-1}[X]\in F$. For any $Y\in \mathcal{G}^\under$ such that $X\subseteq Y$, we have $g^{-1}[X]\subseteq g^{-1}[Y]$. Then since $F$ is a filter in $\mathcal{F}^\under$ and $g^{-1}[Y]\in\mathcal{F}^\under$ by the \PullBack{} property of possibility morphisms, we conclude that $g^{-1}[Y]\in F$. For closure under intersection, if $X,X'\in \mathcal{G}^\under$, $g^{-1}[X]\in F$, and $g^{-1}[X']\in F$, then since $F$ is a filter in $\mathcal{F}^\under$, $g^{-1}[X\cap X']=g^{-1}[X]\cap g^{-1}[X']\in F$. Thus, $\{X\in \mathcal{G}^\under\mid g^{-1}[X]\in F\}$ is a filter in $\mathcal{G}^\under$. To see that it is a proper filter, since $F$ is proper and $g^{-1}[\emptyset]=\emptyset$, we have $g^{-1}[\emptyset]\not\in F$, so $\emptyset\not\in \{X\in \mathcal{G}^\under\mid g^{-1}[X]\in F\}$.
\end{proof}

We now prove an analogue of Theorem \ref{RefSub} for filter-descriptive frames.

\begin{theorem}[Filter-Descriptive Reflections]\label{FiltRefSub} For any possibility frame $\mathcal{F}$, filter-descriptive possibility frame $\mathcal{G}$, and possibility morphism $g\colon \mathcal{F}\to\mathcal{G}$, define $\overline{g}\colon (\mathcal{F}^\under)_\gff\to\mathcal{G}$ such that for $F\in  (\mathcal{F}^\under)_\gff$, $\overline{g}(F)$ is the unique $x\in \mathcal{G}$ such that $P(x)=\{X\in \mathcal{G}^\under\mid g^{-1}[X]\in F\}$, which exists by Lemma \ref{FiltTransfer} and the fact that $\mathcal{G}$ is filter-descriptive. Then:
\begin{enumerate}
\item\label{RefSub1'} $\overline{g}$ is a p-morphism;
\item\label{RefSub2'} $\overline{g}$ is the unique possibility morphism from $(\mathcal{F}^\under)_\gff$ to $\mathcal{G}$ such that $g=\overline{g}\circ \eta_\mathcal{F}$, so the following diagram commutes:
\begin{center}
\begin{tikzpicture}[->,>=stealth',shorten >=1pt,shorten <=1pt, auto,node
distance=2cm,thick,every loop/.style={<-,shorten <=1pt}]
\tikzstyle{every state}=[fill=gray!20,draw=none,text=black]

\node (F) at (0,2) {{$\mathcal{F}$}};
\node (F') at (0,0) {{$(\mathcal{F}^\under)_\gff$}};
\node (G) at (3,2) {{$\mathcal{G}$}};

\path (F) edge[->] node {{$g$}} (G);
\path (F') edge[<-] node {{$\eta_\mathcal{F}$}} (F);
\path (G) edge[<-] node {{$\overline{g}$}} (F');

\end{tikzpicture}
\end{center}
\end{enumerate}
\end{theorem} 

\begin{proof} For part \ref{RefSub1'}, we use the fact that since $\mathcal{G}$ is filter-descriptive, the map $\eta_\mathcal{G}\colon \mathcal{G}\to (\mathcal{G}^\under)_\gff$ is an isomorphism by Proposition \ref{FramesBAOsII}.\ref{FramesBAOsII1}, and the map $(g^\under)_\gff: (\mathcal{F}^\under)_\gff\to (\mathcal{G}^\under)_\gff$ is a p-morphism by Theorem \ref{BAOMorphPossMorphII}.\ref{BAOMorphPossMorphIIa}. Thus, to show that $\overline{g}$ is a p-morphism, it suffices to prove that $\overline{g}=\eta_\mathcal{G}^{-1}\circ (g^\under)_\gff$, or equivalently, $\eta_\mathcal{G}\circ \overline{g}=(g^\under)_\gff$. On one hand, for each $F\in (\mathcal{F}^\under)_\gff$ we have:
\[\begin{array}{llll}
\eta_\mathcal{G}(\overline{g}(F))&=&\{X\in \mathcal{G}^\under\mid \overline{g}(F)\in X\} &\mbox{by definition of }\eta_\mathcal{G}\\
&=&\{X\in \mathcal{G}^\under\mid X\in P(\overline{g}(F))\}&\mbox{by definition of }P(\cdot)\\
&=&\{X\in \mathcal{G}^\under\mid g^{-1}[X]\in F\}&\mbox{by definition of }\overline{g}.\\
\end{array}
\]
 On the other hand,  for each $F\in (\mathcal{F}^\under)_\gff$ we have:
\[\begin{array}{llll}
(g^\under)_\gff(F)&=& (g^\under)^{-1}[F] = \{X\in\mathcal{G}^\under \mid g^\under(X)\in F\}& \mbox{by definition of }(\cdot)_\gff\\
&=&   \{X\in\mathcal{G}^\under \mid g^{-1}[X]\in F\} &\mbox{by definition of }(\cdot)^\under. \\

\end{array}
\]

For part \ref{RefSub2'}, by Proposition \ref{FramesBAOsII}.\ref{FramesBAOsII2} we have $(g^\under)_\gff \circ \eta_\mathcal{F} = \eta_\mathcal{G}\circ g$. Since $\mathcal{G}$ is filter-descriptive, by Proposition \ref{FramesBAOsII}.\ref{FramesBAOsII1}, $\eta_\mathcal{G}$ is an isomorphism. Hence $(g^\under)_\gff \circ \eta_\mathcal{F} = \eta_\mathcal{G}\circ g$ implies $\eta_\mathcal{G}^{-1}\circ (g^\under)_\gff \circ \eta_\mathcal{F} = g$. Then since we showed above that  $\overline{g}=\eta_\mathcal{G}^{-1}\circ (g^\under)_\gff$, we conclude that $\overline{g}\circ \eta_\mathcal{F}=g$, as desired.

Finally, we prove that $\overline{g}$ is the \textit{unique} possibility morphism from $(\mathcal{F}^\under)_\gff$ to $\mathcal{G}$ such that $g=\overline{g}\circ \eta_\mathcal{F}$. Suppose $\overline{h}$ is a possibility morphism from $(\mathcal{F}^\under)_\gff$ to $\mathcal{G}$ such that $g=\overline{h}\circ \eta_\mathcal{F}$. By the definition of $\overline{g}$, to show that $\overline{h}=\overline{g}$, it suffices to show that for each $F\in (\mathcal{F}^\under)_\gff$, we have $P(\overline{h}(F))=\{X\in\mathcal{G}^\under \mid g^{-1}[X]\in F\}$. Recall that $P(\overline{h}(F))=\{X\in P^\mathcal{G}\mid \overline{h}(F)\in X\}=\{X\in \mathcal{G}^\under\mid \overline{h}(F)\in X\}$. For any $X\in\mathcal{G}^\under$, i.e., $X\in \adm^\mathcal{G}$, it follows by \PullBack{} for $\overline{h}$ that $\overline{h}^{-1}[X]\in P^{(\mathcal{F}^\under)_\gff}$, which by definition of $(\cdot)_\gff$ implies that $\overline{h}^{-1}[X]= \widehat{Y}$ for some $Y\in \mathcal{F}^\under$. Since $g=\overline{h}\circ \eta_\mathcal{F}$, we have $g^{-1}[X]= \eta_\mathcal{F}^{-1}[\overline{h}^{-1}[X]]$ and hence $g^{-1}[X]=\eta_\mathcal{F}^{-1}[\widehat{Y}]$. Then since $\eta_\mathcal{F}^{-1}[\widehat{Y}]=Y$, we have $g^{-1}[X]=Y$. Thus, $\overline{h}(F)\in X$ iff $F\in \overline{h}^{-1}[X]=\widehat{Y}$ iff $Y\in F$ iff $g^{-1}[X]\in F$, which completes the proof.\end{proof}

From Propositions \ref{FramesBAOsII}.\ref{FramesBAOsII0} and \ref{PossToP2}, Theorem \ref{FiltRefSub}, and the definition of reflective subcategories, we obtain our desired result.

\begin{theorem}[Reflective Subcategories II]\label{Reflective2}  \textbf{FiltPoss} is a reflective subcategory of \textbf{Poss}.
\end{theorem}

However, we do not obtain the analogous result with \textbf{Poss} replaced by the category of possibility frames with \textit{strict} possibility morphisms. For unlike the case of the $\zeta_\mathcal{F}\colon \mathcal{F}\to (\mathcal{F}^\under)_\rela$ in Theorem \ref{AlmostBack}, the $\eta_\mathcal{F}\colon \mathcal{F}\to (\mathcal{F}^\under)_\gff$ in Proposition \ref{FramesBAOsII} is not guaranteed to be a strict possibility morphism. Consider the \SqBack{} condition: if $Y\sqsubseteq^{(\mathcal{F}^\under)_\gff} \eta_\mathcal{F}(x)$, then $\exists y\colon y\sqsubseteq^\mathcal{F} x$ and $\eta_\mathcal{F}(y)\sqsubseteq^{(\mathcal{F}^\under)_\gff} Y$. To see why this is not guaranteed, suppose $Y$ is an \textit{ultrafilter} in $\mathcal{F}^\under$, so it is a minimal point in $(\mathcal{F}^\under)_\gff$. Then $\eta_\mathcal{F}(y)\sqsubseteq^{(\mathcal{F}^\under)_\gff} Y$ iff $\eta_\mathcal{F}(y)= Y$. But it is not guaranteed that there is a $y$ in the possibility frame $\mathcal{F}$ such that the set $\eta_\mathcal{F}(y)=P(y)$ of admissible propositions containing $y$ is an ultrafilter, let alone the particular ultrafilter $Y$. 

Compare this to the relation between a (general) \textit{world} frame $\mathfrak{F}$ and the general ultrafilter frame of its underlying BAO, $(\mathfrak{F}^*)_*$ (see Appendix \S~\ref{AlgSem}): if $\mathfrak{F}$ is an arbitrary world frame, not necessarily descriptive in Goldblatt's sense, then it is not guaranteed that the function $f$ that sends each world $x$ to the ultrafilter of admissible propositions containing $x$ (see, e.g., \citealt[Thm.~5.76(iii)]{Blackburn2001}) is a p-morphism from $\mathfrak{F}$ to $(\mathfrak{F}^*)_*$. The back clause of a p-morphism requires that if $f(x)R^{(\mathfrak{F}^*)_*}Y$, so $Y$ is an ultrafilter in the underlying BAO $\mathfrak{F}^*$ of $\mathfrak{F}$, then there is a $y\in\mathfrak{F}$ such that $xR^\mathfrak{F}y$ and $f(y)=Y$, so the set of admissible propositions containing $y$ is exactly the ultrafilter $Y$. But if $\mathfrak{F}$ is an arbitrary world frame, then it is not guaranteed that there is such a $y\in\mathfrak{F}$ (cf.~\citealt[p.~94-95]{Blackburn2001}).

Although $\eta_\mathcal{F}$ is not guaranteed to be a strict possibility morphism from $\mathcal{F}$ to $(\mathcal{F}^\under)_\gff$, we have the following corollary of Theorem \ref{PtoB} and Proposition \ref{BAOsFramesII}.\ref{BAOsFramesIIa}.

\begin{corollary}[From Arbitrary Frames to Filter-Descriptive Frames]\label{Arbitrary-to-FD} For any possibility frame $\mathcal{F}$:
\begin{enumerate}
\item $\mathcal{F}^\under$ is isomorphic to $((\mathcal{F}^\under)_\gff)^\under$;
\item for all $\varphi\in\mathcal{L}(\sig,\ind)$, $\mathcal{F}\Vdash\varphi$ iff $(\mathcal{F}^\under)_\gff\Vdash\varphi$.
\end{enumerate}
\end{corollary}

If we instead consider the filter extension $(\mathcal{F}^\under)_\ff$, then we cannot expect an isomorphism of BAOs and validity preservation in both directions as in Corollary \ref{Arbitrary-to-FD}. However, by the proof of Proposition \ref{BAOsFramesII}.\ref{BAOsFramesIIa}, the map $\eta_\mathbb{A}$ sending each $x\in \mathbb{A}$ to $\eta_\mathbb{A}(x)=\widehat{x}\in (\mathbb{A}_\ff)^\under$ is a BAO-embedding of $\mathbb{A}$ into $(\mathbb{A}_\ff)^\under$, so we may take $\mathbb{A}=\mathcal{F}^\under$ to obtain part 1 of the following result. Part 2 also follows directly from Theorems \ref{PtoB}.\ref{PtoB6} and \ref{BAOtoGenPos}.\ref{BAOtoGenPos5}.

\begin{corollary}[From Arbitrary Frames to Filter Extensions]\label{FilterPres} For any possibility frame $\mathcal{F}$:
\begin{enumerate}
\item there is a BAO-embedding of $\mathcal{F}^\under$ into $((\mathcal{F}^\under)_\ff)^\under$;
\item for all $\varphi\in\mathcal{L}(\sig,\ind)$, if $(\mathcal{F}^\under)_\ff\Vdash\varphi$, then $\mathcal{F}\Vdash\varphi$.
\end{enumerate}
\end{corollary}

\subsection{MacNeille Completions and Canonical Extensions of BAOs}\label{MacNeille}

We saw in \S~\ref{DualEquiv} that any $\mathcal{V}$-BAO $\mathbb{A}$ is isomorphic to $(\mathbb{A}_\rela)^\under$, the underlying BAO of the \textit{principal frame} of $\mathbb{A}$, and in \S~\ref{Fdes} that any BAO $\mathbb{A}$ is isomorphic to $(\mathbb{A}_\gff)^\under$, the underlying BAO of the \textit{general filter frame} of $\mathbb{A}$. Next we consider the relation of a $\mathcal{V}$-BAO $\mathbb{A}$ to $(\mathbb{A}_\fullV)^\under$, the underlying BAO of the \textit{full frame} of $\mathbb{A}$ (Definition \ref{Hposs}), and the relation of a BAO $\mathbb{A}$ to $(\mathbb{A}_\ff)^\under$, the underlying BAO of the \textit{filter frame} of $\mathbb{A}$ (Definition \ref{FiltEx}). 
  
To characterize the relation between $\mathbb{A}$ and $(\mathbb{A}_\fullV)^\under$, recall the result (due to \citealt{MacNeille1937,Tarski1937}) that for any Boolean algebra $\mathfrak{A}$, there is a \textit{unique}---up to isomorphism---\textit{complete} Boolean algebra $\mathfrak{A}'$ such that $\mathfrak{A}$ embeds \textit{densely} into $\mathfrak{A}'$, i.e., there is an injective homomorphism $\xi\colon \mathfrak{A}\to\mathfrak{A}'$ such that for every non-minimum $a'\in \mathfrak{A}'$, there is a non-minimum $a\in \mathfrak{A}$ such that $\xi(a)\leq' a'$.\footnote{This $\mathfrak{A}'$ is sometimes called the \textit{completion} of $\mathfrak{A}$ (\citealt[p.~214]{Givant2000}; \citealt[p.~82]{Jech2002}); but in lattice theory, a ``completion'' of a poset is given by any order-embedding of it into a complete lattice (\citealt[p.~165]{Davey2002}).} This $\mathfrak{A}'$ can be characterized and constructed in several ways. For one, it can be constructed by an application of the general MacNeille completion of a poset. (For the moment let us blur the distinction between Boolean algebras and Boolean lattices.) Given a poset $\langle P,\leq\rangle$ and $X\subseteq P$, consider the sets of lower bounds $X^l=\{y\in P\mid \forall x\in X,\, y\leq x\}$ and upper bounds $X^u=\{y\in P\mid \forall x\in X,\, x\leq y\}$ of $X$. The \textit{MacNeille completion} of $\langle P,\leq\rangle$ is the poset $\mathbf{M}( P,\leq)=\langle P',\leq'\rangle$ where $P'=\{X\subseteq P\mid X^{ul}=X\}$ (the set of ``normal'' ideals) and $X\leq' Y$ iff $X\subseteq Y$. For any poset $\langle P,\leq\rangle$, $\mathbf{M}(P,\leq)$ is a complete lattice with $\bigmeet \mathcal{X}=\bigcap\mathcal{X}$ and $\bigjoin \mathcal{X}=(\bigcup \mathcal{X})^{ul}$ for $\mathcal{X}\subseteq P'$, and $\langle P,\leq\rangle$ embeds densely into $\mathbf{M}(P,\leq)$ by $x\mapsto \mathord{\downarrow}x$. For an arbitrary poset $\langle P,\leq\rangle$, there is no guarantee that $\mathbf{M}(P,\leq)$ is a Boolean lattice. But if $\langle P,\leq\rangle$ is a Boolean lattice to begin with, then so is $\mathbf{M}(P,\leq)$.

Compare the above with the construction of the full \textit{regular open algebra} based on $\langle P,\leq\rangle$ as in Remark~\ref{Persp2}, which for any poset $\langle P,\leq\rangle$ produces a complete Boolean lattice, which we will call $\mathbf{R}(P,\leq)$. Recall from \S~\ref{ExtFrames} that if $\langle P,\leq\rangle$ has a minimum element $\bot$, then $\mathbf{R}(P,\leq)$ contains only $P$ and $\emptyset$; so the better comparison is between $\mathbf{M}( P,\leq)$ and $\mathbf{R}(P_{-},\leq_{-})$ where $P_{-}=P\setminus \{\bot\}$ and $\leq_{-}$ is the restriction of $\leq$ to $P_{-}$. If $\langle P_{-},\leq_{-}\rangle$ is \textit{separative} as in \S~\ref{SepSec}, then $\langle P,\leq\rangle$ embeds densely into $\mathbf{R}(P_{-},\leq_{-})$ by $x\mapsto \mathord{\downarrow}_{-}x=\{y\in P_{-}\mid y\leq x\}$. ($\mathbf{R}(P_{-},\leq_{-})$ is often called the \textit{completion by regular cuts} of $\langle P,\leq\rangle$.) Then since $\langle P,\leq\rangle$ being Boolean implies that $\langle P_{-},\leq_{-}\rangle$ is separative (Fact \ref{Boolean}), $\langle P,\leq\rangle$ being Boolean implies that $\langle P,\leq\rangle$ densely embeds into $\mathbf{R}(P_{-},\leq_{-})$. These points show that although for an arbitrary poset $\langle P,\leq\rangle$, there is no guarantee that $\mathbf{M}(P,\leq)$ and $\mathbf{R}(P_{-},\leq_{-})$ are isomorphic, if $\langle P,\leq\rangle$ is Boolean, then $\mathbf{M}(P,\leq)$ and $\mathbf{R}(P_{-},\leq_{-})$ are both complete Boolean lattices into which our Boolean $\langle P,\leq\rangle$ densely embeds, so they are isomorphic.\footnote{Yet another way to get an isomorphic copy of this completion is by taking the regular open algebra of the \textit{Stone space} \citep{Stone1937} of the original Boolean algebra.}

Let us apply the foregoing points to the relation between $\mathbb{A}$ and $(\mathbb{A}_\fullV)^\under$. Recall that the poset in $\mathbb{A}_\fullV$ is the Boolean lattice of $\mathbb{A}$ minus its bottom element, and the Boolean algebra reduct of  $(\mathbb{A}_\fullV)^\under$ is the full regular open algebra based on the poset in $\mathbb{A}_\fullV$. Thus, by the previous paragraph, the Boolean reduct of the BAO  $(\mathbb{A}_\fullV)^\under$ is isomorphic to the \textit{MacNeille completion} of the Boolean reduct of the BAO $\mathbb{A}$. 

To relate $\mathbb{A}$ and $(\mathbb{A}_\fullV)^\under$ as BAOs, recall that the \textit{Monk completion} or \textit{lower} MacNeille completion of a  BAO $\mathbb{A}$ \citep{Monk1970} is obtained by extending the operators $\blacklozenge_i$ of $\mathbb{A}$ to operators $\blacklozenge_i^\circ$ on the MacNeille completion of the Boolean reduct of $\mathbb{A}$ as follows, where $\xi$ is the embedding into the completion with order $\leq$:
\begin{equation}\blacklozenge^\circ_i Y = \bigvee \{\xi(\blacklozenge_i x)\mid x\in\mathbb{A}\mbox{ and }\xi(x)\leq Y\}.\label{Monk1}\end{equation}
If $\blacklozenge_i$ is completely additive, then so is  $\blacklozenge_i^\circ$ \citep[Thm.~1.2]{Monk1970}; otherwise the Monk completion of a BAO is not even guaranteed to be a BAO. For the dual box operator, we have:
\begin{equation}\blacksquare^\circ_i Y = \bigmeet \{\xi(\blacksquare_i x)\mid x\in\mathbb{A}\mbox{ and }Y\leq \xi(x)\}.\label{Monk2}\end{equation}
\noindent Properties of the Monk completion have been investigated in \citealt{Monk1970,Givant1999,Gehrke2005,Harding2007} and \citealt{Theunissen2007}. With Theorem \ref{MonkComp} below, we obtain a new characterization of the Monk completion of a $\mathcal{V}$-BAO.  The key point is that instead of thinking of extending the operators from $\mathbb{A}$ to $\mathbb{A}^\circ$ in terms of (\ref{Monk1})/(\ref{Monk2}), we can equivalently think in terms of first defining an accessibility relation $R_i$ on $\mathbb{A}$ (minus its bottom $\bot$) as in Definition \ref{Hposs}.\ref{Hposs3}, and then defining the extended $\blacksquare_i^\circ$ using the relation $R_i$ as usual, so $\blacksquare_i^\circ Y$ is the set of $x$ for which $R_i(x)\subseteq Y$.

\begin{theorem}[Monk Completion]\label{MonkComp} The Monk completion of a $\mathcal{V}$-BAO $\mathbb{A}$ is isomorphic to $(\mathbb{A}_\fullV)^\under$.
\end{theorem}

\begin{proof} Since we already saw above that the MacNeille completion of the Boolean reduct of $\mathbb{A}$ is isomorphic to the Boolean reduct of $(\mathbb{A}_\fullV)^\under$, we need only show that extending the operator $\blacksquare_i^\mathbb{A}$ to $\blacksquare^\circ_i$ as in (\ref{Monk2}) is equivalent to extending $\blacksquare_i^\mathbb{A}$ to $\blacksquare_i^{(\mathbb{A}_\fullV)^\under}$, where $\blacksquare_i^{(\mathbb{A}_\fullV)^\under} Y=\{x\in \mathbb{A}_\fullV\mid R_i^{\mathbb{A}_\fullV}(x)\subseteq Y\}$. Recall that the poset $\langle S,\sqsubseteq\rangle$ in $\mathbb{A}_\fullV$ is the poset of $\mathbb{A}$ minus its bottom $\bot$, and the poset of $(\mathbb{A}_\fullV)^\under$ is the set of all regular open sets in the downset topology on $\langle S,\sqsubseteq\rangle$,  ordered by inclusion. For $z\in\mathbb{A}_\fullV$,  let $\mathord{\downarrow}^{\mathbb{A}_\fullV}z=\{z'\in\mathbb{A}_\fullV\mid z'\sqsubseteq z\}$. As above, the Boolean reduct of $\mathbb{A}$ embeds into that of $(\mathbb{A}_\fullV)^\under$ by $\xi(x)=\mathord{\downarrow}^{\mathbb{A}_\fullV}x$. Then we first observe that for any $Y\in (\mathbb{A}_\fullV)^\under$:
\[\begin{array}{lcll}
\blacksquare^\circ_i Y &=& \bigmeet \{\xi(\blacksquare_i^\mathbb{A} y)\mid y\in\mathbb{A}\mbox{ and }Y\leq \xi(y)\} &\mbox{by (\ref{Monk2})}\\
&=& \bigcap \{\mathord{\downarrow}^{\mathbb{A}_\fullV}\blacksquare_i^\mathbb{A} y\mid y\in\mathbb{A} \mbox{ and }Y\subseteq \mathord{\downarrow}^{\mathbb{A}_\fullV}y\}&\mbox{by definition of }(\mathbb{A}_\fullV)^\under\mbox{ and }\xi \\
&=& \{x\in\mathbb{A}_\fullV\mid \forall y\in \mathbb{A}\colon \mbox{if }Y\subseteq \mathord{\downarrow}^{\mathbb{A}_\fullV}y\mbox{, then }x\sqsubseteq \blacksquare_i^\mathbb{A} y\}&\mbox{by definition of }\mathord{\downarrow}^{\mathbb{A}_\fullV}\\
&=& \{x\in\mathbb{A}_\fullV\mid \forall y\in \mathbb{A}\colon \mbox{if }Y\subseteq \mathord{\downarrow}^{\mathbb{A}_\fullV}y\mbox{, then }x\in \blacksquare_i^{(\mathbb{A}_\fullV)^\under} \mathord{\downarrow}^{\mathbb{A}_\fullV} y\}\\
&\supseteq &  \blacksquare_i^{(\mathbb{A}_\fullV)^\under}Y,
\end{array}\]
where the last equation holds by the equivalence of $x\sqsubseteq \blacksquare_i^\mathbb{A} y$ and $x\in \blacksquare_i^{(\mathbb{A}_\fullV)^\under} \mathord{\downarrow}^{\mathbb{A}_\fullV} y$ already proved for Theorem \ref{VtoPossFrames}.\ref{VtoPoss4}, and the $\supseteq$ inclusion holds because $Y\subseteq \mathord{\downarrow}^{\mathbb{A}_\fullV}y$ implies $\blacksquare_i^{(\mathbb{A}_\fullV)^\under}Y\subseteq \blacksquare_i^{(\mathbb{A}_\fullV)^\under}\mathord{\downarrow}^{\mathbb{A}_\fullV}y$. 

In the other direction, suppose $x\not\in  \blacksquare_i^{(\mathbb{A}_\fullV)^\under}Y$, so $\exists z$: $xR_i^{\mathbb{A}_\fullV}z$ and $z\not\in Y$.  Since $z\not\in Y$, it follows by \textit{refinability} for $Y$ that $\exists z'\sqsubseteq z$ $\forall z''\sqsubseteq z'$, $z''\not\in Y$.  Then by \textit{persistence} for $Y$, it follows that for all $y\in Y$, $y\meet z'\not\in Y$, so $y\meet z'=\bot$, so $y\sqsubseteq -z'$. Thus, $Y\subseteq \mathord{\downarrow}^{\mathbb{A}_\fullV} \mathord{-}z'$. Finally, by definition of the relation $R_i^{\mathbb{A}_\fullV}$ in $\mathbb{A}_\fullV$ (Definition \ref{Hposs}.\ref{Hposs3}), $xR_i^{\mathbb{A}_\fullV}z$ and $z'\sqsubseteq z$ together imply that $x\meet \blacklozenge_i^\mathbb{A}z'\not=\bot$, so $x\not\sqsubseteq \blacksquare_i^\mathbb{A} \mathord{-}z'$. Thus, we have found a $y\in\mathbb{A}$, namely $y=\mathord{-}z'$, such that $Y\subseteq \mathord{\downarrow}^{\mathbb{A}_\fullV}y$ but $x\not\sqsubseteq \blacksquare_i^\mathbb{A} y$, so $x\not\in\blacksquare^\circ_iY$ by the equations above.\end{proof}

Using Theorem~\ref{MonkComp}, we can transfer results about Monk completions to results about $(\mathbb{A}_\fullV)^\under$ and vice versa. We will see an example of the utility of this connection in \S~\ref{AtomlessFull}. 

Let us now turn to the relation between $\mathbb{A}$ and $(\mathbb{A}_\ff)^\under$. Following \citealt{Jonsson1952a}, a BAO $\mathbb{B}$ is a \textit{perfect extension} of a BAO $\mathbb{A}$ iff: (i) $\mathbb{B}$ is a $\mathcal{CAV}$-BAO, and there is a BAO-embedding $e$ of $\mathbb{A}$ into $\mathbb{B}$; (ii) if $a$ and $b$ are distinct atoms in $\mathbb{B}$, then there is an $x\in\mathbb{A}$ such that $a\leq^\mathbb{B} e(x)$ and $e(x)\leq^\mathbb{B}-b$; and (iii) if $X$ is a set of elements from $\mathbb{A}$ such that $\bigvee^\mathbb{B} e[X]=\top$, then there is a finite $X'\subseteq X$ such that $\bigvee X'=\top$. J\'{o}nsson and Tarski showed that every BAO has a perfect extension, assuming the ultrafilter axiom or an equivalent axiom, and any two perfect extensions are isomorphic, so we may speak of \textit{the} perfect extension of a BAO. The perfect extension of $\mathbb{A}$ can be constructed as the full complex algebra of the ultrafilter frame of $\mathbb{A}$ (see \S~\ref{AlgSem}). If we think of the ultrafilter frame as a full world frame or full possibility frame (recall Example \ref{KripkeExample}) instead of a Kripke frame, then we can say that the perfect extension of $\mathbb{A}$ arises as the underlying BAO of the ultrafilter frame of $\mathbb{A}$. By contrast, our $(\mathbb{A}_\ff)^\under$ is the underlying BAO of the \textit{filter} frame of $\mathbb{A}$.  

Assuming the ultrafilter axiom, the filter frame $\mathbb{A}_\ff$ of $\mathbb{A}$ is an \textit{atomic} full possibility frame as in \S~\ref{AtomicSection}, and the \textit{atom structure} $\mathfrak{At}(\mathbb{A}_\ff)$ of $\mathbb{A}_\ff$ (Definition \ref{AtSt}) is the \textit{ultra}filter frame of $\mathbb{A}$ considered as a full possibility frame. Thus, the underlying BAO $(\mathfrak{At}(\mathbb{A}_\ff))^\under$ of $\mathfrak{At}(\mathbb{A}_\ff)$ is isomorphic to the perfect extension of $\mathbb{A}$. By Proposition \ref{AtToWorld}.\ref{AtToWorld3}, there is a dense and robust possibility morphism from $\mathfrak{At}(\mathbb{A}_\ff)$ to $\mathbb{A}_\ff$, so by Theorem \ref{PossMorphBAOMorph}, $(\mathfrak{At}(\mathbb{A}_\ff))^\under$ is isomorphic to $(\mathbb{A}_\ff)^\under$. Therefore, assuming the ultrafilter axiom, $(\mathbb{A}_\ff)^\under$ is a perfect extension of~$\mathbb{A}$. 

The perfect extension of a BAO has come to be called the \textit{canonical extension} of the BAO. However, there is a different definition of the canonical extension, due to Gehrke and Harding \citeyearpar[Def.~2.5]{Gehrke2001} (in the more general setting of lattices with additional operations), that does not require the canonical extension to be atomic, as required for a perfect extension; with this definition, one can prove in ZF set theory that the canonical extension exists and is unique up to isomorphism (op.~cit., Props.~2.6-2.7), and then one can prove in ZF plus the ultrafilter axiom that the canonical extension is a perfect extension (cf.~op.~cit., Lem.~3.4). In addition to the construction of this canonical extension given in \citealt{Gehrke2001}, there is another way of constructing it as the MacNeille completion of a certain intermediate structure (see \citealt{Ghilardi1997,Dunn2005,Gehrke2008}). According to the Gehrke-Harding definition applied to BAOs, a BAO $\mathbb{B}$ is a canonical extension of a BAO $\mathbb{A}$ iff: (i$'$) $\mathbb{B}$ is a $\mathcal{CV}$-BAO, and there is a BAO-embedding $e$ of $\mathbb{A}$ into $\mathbb{B}$; (ii$'$) every element of $\mathbb{B}$ is a join of meets of $e$-images of elements of $\mathbb{A}$---or equivalently in this Boolean case, every element of $\mathbb{B}$ is a meet of joins of $e$-images of elements of $\mathbb{A}$; and (iii$'$) for any sets $X,Y$ of elements of $\mathbb{A}$, if $\bigwedge^\mathbb{B} e[X]\leq^\mathbb{B} \bigvee^\mathbb{B} e[Y]$, then there are finite $X'\subseteq X$ and $Y'\subseteq Y$ such that $\bigwedge X'\leq \bigvee Y'$. Following the terminology of \citealt{Conradie2016}, let us call such a $\mathbb{B}$ a \textit{constructive canonical extension} of $\mathbb{A}$ to emphasize that its existence can be proved without nonconstructive choice principles. Indeed, the following is provable without use of the ultrafilter axiom.

\begin{theorem}[Canonical Extension]\label{CanExtProp} For any BAO $\mathbb{A}$, $(\mathbb{A}_\ff)^\under$ is a constructive canonical extension of $\mathbb{A}$.
\end{theorem} 

\begin{proof} For condition (i$'$), by Theorem \ref{BAOtoGenPos}.\ref{IsRelGen}, $\mathbb{A}_\ff$ is a full possibility frame, so by Theorem \ref{PtoB}.\ref{PtoB2}, $(\mathbb{A}_\ff)^\under$ is a $\mathcal{CV}$-BAO. By the proof of Proposition \ref{BAOsFramesII}, the map $\eta_\mathbb{A}$ sending each $x\in \mathbb{A}$ to $\eta_\mathbb{A}(x)=\widehat{x}\in (\mathbb{A}_\ff)^\under$ is a BAO-embedding. Recall that $\widehat{x}$ is the set of all proper filters in $\mathbb{A}$ that contain $x$.

For condition (ii$'$), each $\mathcal{X}\in (\mathbb{A}_\ff)^\under$ is a regular open downset of filters in $\mathbb{A}_\ff$, so we have 
\[\mathcal{X}=\mathrm{int}(\mathrm{cl}(\mathcal{X}))=\mathrm{int}(\mathrm{cl} (\underset{F\in \mathcal{X}}\bigcup \mathord{\downarrow}F))=\mathrm{int}(\mathrm{cl} (\underset{F\in \mathcal{X}}\bigcup \,\underset{x\in F}{\bigcap}\hat{x}))=\underset{F\in \mathcal{X}}{\bigvee}\,\underset{x\in F}{\bigwedge}\hat{x},\]
where $\bigvee$ and $\bigmeet$ are the join and meet in the regular open algebra $(\mathbb{A}_\ff)^\under$. 

For condition (iii$'$), suppose that for sets $X,Y$ of elements of $\mathbb{A}$, $\bigwedge \{\widehat{x}\mid x\in X\}\leq^{(\mathbb{A}_\ff)^\under}\bigvee \{\widehat{y}\mid y\in Y\}$, so 
\begin{equation}\bigcap \{\widehat{x}\mid x\in X\}\subseteq \mathrm{int}(\mathrm{cl}(\bigcup\{\widehat{y}\mid y\in Y\})),\label{WhichMeans}\end{equation}
which means that for any proper filter $F$ in $\mathbb{A}$ such that $X\subseteq F$, we have that for all proper filters $F'\supseteq F$, there is a proper filter $F''\supseteq F'$ such that $y\in F''$ for some $y\in Y$. Now suppose for reductio that for every finite $X'\subseteq X$ and $Y'\subseteq Y$, $\bigwedge X'\not\leq \bigvee Y'$. It follows by Fact \ref{FiltGen} that $F=[X \cup \{-y\mid y\in Y\})$ is a proper filter. Taking $F'=F$ in the unpacking of (\ref{WhichMeans}) above, there is a proper filter $F''\supseteq F$ such that $y\in F''$ for some $y\in Y$, which is impossible given our choice of $F$.
\end{proof}

We will return to this connection between $(\mathbb{A}_\ff)^\under$ and the constructive canonical extension of $\mathbb{A}$ in \S~\ref{Canonical}.

\subsection{Frame Constructions and Algebraic Constructions}\label{OpFrameAlg}

The next step in developing the duality theory for possibility frames and BAOs is to relate \textit{frame constructions} that preserve the validity of modal formulas with \textit{algebraic constructions} that preserve the validity of modal formulas, or in more algebraic terms, that preserve universally quantified algebraic equations. The story here parallels the story for world frames and BAOs \citep[\S\S~1.4-1.6]{Goldblatt1974} but with some twists.

Let us recall the standard algebraic notions of homomorphic images, subalgebras, and direct products. A BAO $\mathbb{A}$ is a \textit{homomorphic image} of a BAO $\mathbb{B}$ iff there is a surjective BAO-homomorphism as in Definition \ref{Homomorph} from $\mathbb{B}$ to $\mathbb{A}$. Where $A$ and $B$ are the domains of $\mathbb{A}$ and $\mathbb{B}$, respectively, $\mathbb{A}$ is a  \textit{subalgebra} of $\mathbb{B}$ iff $A\subseteq B$, $A$ is closed under the operations of $\mathbb{B}$, and the operations of $\mathbb{A}$ are the restrictions to $A$ of the operations of $\mathbb{B}$. Finally, given a family $\{\mathbb{A}_j\}_{j\in J}$ of BAOs $\mathbb{A}_j=\langle A_j, \meet_j, -_j, \top_j, \{\blacksquare_{i,j}\}_{i\in \ind}\rangle$, their \textit{direct product} is the BAO $\underset{j\in J}{\prod}\mathbb{A}_j=\langle A, \meet, -, \top, \{\blacksquare_i\}_{i\in \ind}\rangle$ where $A$ is the Cartesian product $\underset{j\in J}\prod A_j$ and the operations of $\underset{j\in J}{\prod}\mathbb{A}_j$ are defined coordinatewise from those of the $\mathbb{A}_j$, i.e., for functions $x,y\in\underset{j\in J}\prod A_j$, $x\meet y$ is the function in  $\underset{j\in J}\prod A_j$ such that for all $j\in J$, $(x\meet y)_j=x_j\meet_jy_j$, where $f_j$ is the value of function $f$ at $j$, and similarly for the other operations. One can easily check that $\underset{j\in J}{\prod}\mathbb{A}_j$ is indeed a BAO.

Taking homomorphic images, subalgebras, and direct products of algebras are ways of preserving universally quantified algebraic equations. Let us now compare these with ways of preserving the validity of modal formulas over possibility frames. We recall the following from Proposition \ref{Preservation}.\ref{Preservation3}.

\begin{proposition}[Preservation Under Dense Possibility Morphisms]\label{DenseRecall} For any possibility frames $\mathcal{F}$ and $\mathcal{F}'$, if there is a \textit{dense} possibility morphism from $\mathcal{F}$ to $\mathcal{F}'$, then for all $\varphi\in\mathcal{L}(\sig,\ind)$, $\mathcal{F}\Vdash\varphi$ implies $\mathcal{F}'\Vdash\varphi$.
\end{proposition}
\noindent A special case of dense possibility morphisms are \textit{surjective} possibility morphisms, but surjectivity is not required to preserve validity.

We have already done the work with Theorems \ref{BAOMorphPossMorphII} and \ref{PossMorphBAOMorph} to relate possibility morphisms and subalgebras as follows.

\begin{proposition}[Possibility Morphisms and Subalgebras]\label{Poss&Sub} For any BAOs $\mathbb{A}$ and $\mathbb{B}$:
\begin{enumerate}[label=\arabic*.,ref=\arabic*]
\item\label{Poss&Sub1} if $\mathbb{A}$ is isomorphic to a subalgebra of $\mathbb{B}$, then there is a surjective p-morphism from $\mathbb{B}_\gff$ to $\mathbb{A}_\gff$;
\item\label{Poss&Sub2} if $\mathbb{A}$ is isomorphic to a subalgebra of $\mathbb{B}$, then there is a surjective p-morphism from $\mathbb{B}_\ff$ to $\mathbb{A}_\ff$.
\end{enumerate}
Conversely, for any possibility frames $\mathcal{F}$ and $\mathcal{G}$:
\begin{enumerate}[label=\arabic*.,ref=\arabic*,resume]
\item\label{Poss&Sub3} if there is a dense possibility morphism from $\mathcal{F}$ to $\mathcal{G}$, then $\mathcal{G}^\under$ is isomorphic to a subalgebra of $\mathcal{F}^\under$.
\end{enumerate} 
\end{proposition}

\begin{proof} For part \ref{Poss&Sub1}, if $\mathbb{A}$ is isomorphic to a subalgebra of $\mathbb{B}$, then the isomorphism is an injective homomorphism from $\mathbb{A}$ to $\mathbb{B}$, so by parts \ref{BAOMorphPossMorphIIa} and \ref{BAOMorphPossMorphIIc} of Theorem \ref{BAOMorphPossMorphII}, there is a surjective p-morphism from $\mathbb{B}_\gff$ to $\mathbb{A}_\gff$.  The proof of part \ref{Poss&Sub2} is the same but using part \ref{BAOMorphPossMorphIId} of Theorem  \ref{BAOMorphPossMorphII} instead of parts \ref{BAOMorphPossMorphIIa} and \ref{BAOMorphPossMorphIIc}.

For part \ref{Poss&Sub3}, if there is a dense possibility morphism from $\mathcal{F}$ to $\mathcal{G}$, then by Theorem \ref{PossMorphBAOMorph}, there is an injective homomorphism from $\mathcal{G}^\under$ to $\mathcal{F}^\under$, so $\mathcal{G}^\under$ is isomorphic to a subalgebra of $\mathcal{F}^\under$.
\end{proof}

Another validity preserving construction for possibility frames is given by the notion of \textit{generated subframes}, which parallels the standard definition for world frames \citep[\S~1.4]{Goldblatt1974}. In fact, for possibility frames there is also a more liberal notion of \textit{selective subframe}. (Compare this to the definition of \textit{cofinal subframes} in \citealt[p.~295]{Chagrov1997} and \citealt[p.~61]{Bezhanishvili2006}.)

\begin{definition}[Subframes]\label{Subframes} Given a possibility frame $\mathcal{F}=\langle S,\sqsubseteq, \{R_i\}_{i\in\ind},\adm\rangle$, a \textit{subframe} of $\mathcal{F}$ is a tuple $\mathcal{F}'=\langle S',\sqsubseteq ', \{R_i'\}_{i\in\ind},\adm'\rangle$ such that:
\begin{enumerate}[label=\arabic*.,ref=\arabic*]
\item\label{SubRestrict} $\emptyset\neq S'\subseteq S$; $\sqsubseteq'\,=\,\sqsubseteq \cap \,( S'\times S')$, and $R_i'=R_i\cap (S'\times S')$;
\item\label{SubAd} $\adm'=\{X\cap S'\mid X\in\adm\}$.
\end{enumerate} 
A \textit{generated subframe} of $\mathcal{F}$ is a subframe $\mathcal{F}'$ of $\mathcal{F}$ such that $S'$ is closed under $\sqsubseteq$ and $R_i$ from $\mathcal{F}$:
\begin{enumerate}[label=\arabic*.,ref=\arabic*,resume]
\item if $x\in S'$ and $y\sqsubseteq x$, then $y\in S'$; if $x\in S'$ and $xR_iy$, then $y\in S'$.
\end{enumerate}
A \textit{selective subframe} of $\mathcal{F}$ is a subframe $\mathcal{F}'$ of $\mathcal{F}$ such that:
\begin{enumerate}[label=\arabic*.,ref=\arabic*,resume]
\item\label{SelSub1} if $x\in S'$ and $y\sqsubseteq x$, then there is a $z\in S'$ such that $z\sqsubseteq y$.
\item\label{SelSub2} if $x\in S'$, $xR_iy$, and $u\sqsubseteq y$, then there is a $z\in S'$ such that $z\comp u$ and $xR_i z$. \hfill $\triangleleft$
\end{enumerate} 
\end{definition}

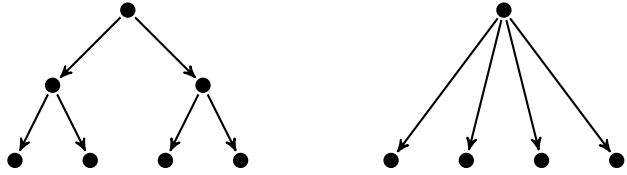
\begin{figure}[h]
\begin{center}
\begin{tikzpicture}[yscale=1, ->,>=stealth',shorten >=1pt,shorten <=1pt, auto,node
distance=2cm,thick,every loop/.style={<-,shorten <=1pt}]
\tikzstyle{every state}=[fill=gray!20,draw=none,text=black]

\node[circle,draw=black!100,fill=black!100,inner sep=0pt,minimum size=.175cm]  (x) at (0,0) {{}};
\node[circle,draw=black!100,fill=black!100,inner sep=0pt,minimum size=.175cm]  (y) at (-1,-1) {{}};
\node[circle,draw=black!100,fill=black!100,inner sep=0pt,minimum size=.175cm] (z) at (1,-1) {{}};

\node[circle,draw=black!100,fill=black!100,inner sep=0pt,minimum size=.175cm]  (y') at (-1.5,-2) {{}};
\node[circle,draw=black!100,fill=black!100,inner sep=0pt,minimum size=.175cm] (y'') at (-.5,-2) {{}};

\node[circle,draw=black!100,fill=black!100,inner sep=0pt,minimum size=.175cm]  (z') at (.5,-2) {{}};
\node[circle,draw=black!100,fill=black!100,inner sep=0pt,minimum size=.175cm] (z'') at (1.5,-2) {{}};

\path (x) edge[->] node {{}} (y);
\path (x) edge[->] node {{}} (z);
\path (y) edge[->] node {{}} (y');
\path (y) edge[->] node {{}} (y'');
\path (z) edge[->] node {{}} (z');
\path (z) edge[->] node {{}} (z'');

\node[circle,draw=black!100,fill=black!100,inner sep=0pt,minimum size=.175cm]  (a) at (5,0) {{}};

\node[circle,draw=black!100,fill=black!100,inner sep=0pt,minimum size=.175cm]  (b') at (3.5,-2) {{}};
\node[circle,draw=black!100,fill=black!100,inner sep=0pt,minimum size=.175cm] (b'') at (4.5,-2) {{}};

\node[circle,draw=black!100,fill=black!100,inner sep=0pt,minimum size=.175cm]  (c') at (5.5,-2) {{}};
\node[circle,draw=black!100,fill=black!100,inner sep=0pt,minimum size=.175cm] (c'') at (6.5,-2) {{}};

\path (a) edge[->] node {{}} (b');
\path (a) edge[->] node {{}} (b'');
\path (a) edge[->] node {{}} (c');
\path (a) edge[->] node {{}} (c'');

\end{tikzpicture}
\end{center}

\caption{A possibility frame (left) and one of its selective subframes (right) that is not a generated subframe. For both frames, assume that the accessibility relation $R_i$ is the universal relation.}\label{SelectiveFig}
\end{figure}

Note the following facts about Definition \ref{Subframes}. First, a subframe of $\mathcal{F}$ is not guaranteed to be a possibility frame. Second, every generated subframe of $\mathcal{F}$ is a selective subframe of $\mathcal{F}$, but not vice versa, as shown in Figure \ref{SelectiveFig}. Third, if $\mathcal{F}$ is a world frame, so $\sqsubseteq$ is identity, then every selective subframe of $\mathcal{F}$ is a generated subframe of $\mathcal{F}$. 

\begin{proposition}[Preservation under Selective Subframes]\label{GenPres} If $\mathcal{F}'$ is a selective subframe of a possibility frame $\mathcal{F}$, then:
\begin{enumerate}
\item\label{GenPres1} $\mathcal{F}'$ is a possibility frame;
\item\label{GenPres2} if $\mathcal{F}$ is full, then $\mathcal{F}'$ is full;
\item\label{GenPres3} for all $\varphi\in\mathcal{L}(\sig,\ind)$, $\mathcal{F}\Vdash\varphi$ implies $\mathcal{F}'\Vdash\varphi$.
\end{enumerate}
\end{proposition}

\begin{proof} For part \ref{GenPres1}, we must first show that $\adm'$ is closed under the operations $\cap$, $\supset'$, and $\blacksquare_i'$ as in Definition \ref{PosetMod}, and $\emptyset\in\adm'$.  Since $\emptyset\in\adm$, we have $\emptyset=\emptyset\cap S'\in \adm'$ by Definition \ref{Subframes}.\ref{SubAd}. For the closure conditions, suppose $X',Y'\in \adm'$, so by Definition \ref{Subframes}.\ref{SubAd} there are $X,Y\in\adm$ such that $X'=X\cap S'$ and $Y'=Y\cap S'$. Where $\supset$ and $\blacksquare_i$ are the operations in $\mathcal{F}$, we claim that:
\begin{itemize}
\item[(a)] $X'\cap Y' = (X\cap Y)\cap S'$;
\item[(b)] $X'\supset' Y'=(X\supset Y)\cap S'$; 
\item[(c)] $\blacksquare_i'X'=(\blacksquare_i X)\cap S'$.
\end{itemize}
From (a)-(c), $X,Y\in\adm$, the fact that $\adm$ is closed under $\cap$, $\supset$, and $\blacksquare_i$, and Definition \ref{Subframes}.\ref{SubAd}, it follows that $X'\cap Y'\in\adm'$, $X'\supset'Y'\in\adm'$, and $\blacksquare_i'X'\in\adm'$, as desired. Hence $\mathcal{F}'$ is a partial-state frame.

Equation (a) is immediate from $X'=X\cap S'$ and $Y'=Y\cap S'$. For (b), the right-to-left inclusion is straightforward. For the left-to-right inclusion, suppose $x\not\in (X\supset Y)\cap S'$. If $x\not\in S'$, then $x\not\in X'\supset' Y'$, so suppose $x\in S'$. From $x\not\in X\supset Y$, there is a $z\sqsubseteq x$ such that $z\in X$ but $z\not\in Y$. Then by \textit{refinability} for $Y$, there is a $z'\sqsubseteq z$ such that for all $z''\sqsubseteq z'$, $z''\not\in Y$. By Definition \ref{Subframes}.\ref{SelSub1}, since $x\in S'$ and $z'\sqsubseteq x$, there is a $z''\in S'$ such that $z''\sqsubseteq z'$. Then by the \textit{refinability} step, $z''\not\in Y$, so $z''\not\in Y'=Y\cap S'$, but by \textit{persistence} for $X$, $z''\sqsubseteq z\in X$ implies $z''\in X$, which with $z''\in S'$ implies $z''\in X'=X\cap S'$. Finally, from $z''\in S'$ and $z''\sqsubseteq z'\sqsubseteq z\sqsubseteq x$, we have $z''\sqsubseteq' x$. Hence we have a $z''\sqsubseteq' x$ such that $z''\in X'$ and $z''\not\in Y'$, so $x\not\in X'\supset' Y'$.

 For (c), the right-to-left inclusion is again straightforward. For the left-to-right inclusion, suppose $x\in S'$ but $x\not\in\blacksquare_i X$, so there is a $y$ such that $xR_i y$ and $y\not\in X$. Then by \textit{refinability} for $X$, there is a $u\sqsubseteq y$ such that for all $u'\sqsubseteq u$, $u'\not\in X$. By Definition \ref{Subframes}.\ref{SelSub2}, since $x\in S'$, $xR_i y$, and $u\sqsubseteq y$, there is a $z\in S'$ such that $z\comp u$ and $xR_i z$, so $xR_i'z$. Since for all $u'\sqsubseteq u$, $u'\not\in X$, and $z\comp u$, \textit{persistence} for $X$ implies $z\not\in X$, so $z\not\in X'=X\cap S'$, which with $xR_i'z$ implies $x\not\in\blacksquare_i' X'$.

To show that $\mathcal{F}'$ is a possibility frame, it only remains to show that $P'\subseteq\mathrm{RO}(S',\sqsubseteq')$. Suppose $X'\in \adm'$, so there is an $X\in \adm$ such that $X'=X\cap S'$. To show that $X'$ satisfies \textit{persistence} with respect to $\langle S',\sqsubseteq'\rangle$, suppose $x\in X'=X\cap S'$ and $y\sqsubseteq' x$. Then $y\sqsubseteq x$, which with $x\in X$ and \textit{persistence} for $X$ implies $y\in X$, which with $y\in S'$ implies $y\in X'=X\cap S'$. To show that $X'$ satisfies \textit{refinability} with respect to $\langle S',\sqsubseteq'\rangle$, suppose $x\in S'$ but $x\not\in X'=X\cap S'$. Then $x\not\in X$, so by \textit{refinability} for $X$, there is a $y\sqsubseteq x$ such that for all $z\sqsubseteq y$, $z\not\in X$.  By Definition \ref{Subframes}.\ref{SelSub1}, since $x\in  S'$ and $y\sqsubseteq x$, there is a $z\in S'$ such that $z\sqsubseteq y$ and hence $z\sqsubseteq' x$. Now for any $u\sqsubseteq' z$, we have $u\sqsubseteq y$ and hence $u\not\in X$ by the \textit{refinability} step, so $u\not\in X'$. Thus, we have shown that if $x\not\in X'$, then there is a $z\sqsubseteq' x$ such that for all $u\sqsubseteq' z$, $u\not\in X'$, as desired.
  
For part \ref{GenPres2}, assuming $\mathcal{F}$ is full, so $P=\mathrm{RO}(S,\sqsubseteq)$, we will show that $\mathrm{RO}(S',\sqsubseteq')\subseteq P'$, which with the previous paragraph shows that $\mathcal{F}'$ is full. Suppose $X'\in \mathrm{RO}(S',\sqsubseteq')$.  By Definition \ref{Subframes}.\ref{SubAd}, to show $X'\in P'$, it suffices to show that there is an $X\in\adm=\mathrm{RO}(S,\sqsubseteq)$ such that $X'=X\cap S'$. Where $\mathord{\Downarrow}X'=\{x\in S\mid \exists y\in X'\colon x\sqsubseteq y\}$, let $X=\mathrm{int}(\mathrm{cl}(\mathord{\Downarrow}X'))$, so $X'\subseteq X\in\mathrm{RO}(S,\sqsubseteq)$ by Fact \ref{RefReg}.\ref{RefReg2.5}. It remains to show that $X\cap S'\subseteq X'$. Suppose $x\in S'$. Then if $x\not\in X'$,  \textit{refinability} for $X'$ with respect to $\sqsubseteq'$ gives us a $y\sqsubseteq' x$ such that (i) for all $z\sqsubseteq' y$, $z\not\in X'$. Now suppose for reductio that $y\in  \mathrm{cl}(\mathord{\Downarrow}X')$, so there is a $u\sqsubseteq y$ and a $y'\in X'$ such that $u\sqsubseteq y'$. Then by Definition \ref{Subframes}.\ref{SelSub1}, there is a $z\in S'$ such that $z\sqsubseteq u$. Since $z\sqsubseteq u\sqsubseteq y$, $z\sqsubseteq u\sqsubseteq y'$, $z\in S'$, and $y,y'\in S'$, we have $z\sqsubseteq' y$ and $z\sqsubseteq' y'$. By (i), $z\sqsubseteq' y$ implies $z\not\in X'$, which with $z\sqsubseteq' y'$ and \textit{persistence} for $X'$ with respect to $\sqsubseteq'$ implies $y'\not\in X'$, contradicting $y'\in X'$ from above. Hence $y\not\in \mathrm{cl}(\mathord{\Downarrow}X')$, which with $y\sqsubseteq' x$ implies $x\not\in\mathrm{int}(\mathrm{cl}(\mathord{\Downarrow}X'))$, so $x\not\in X$, which completes the proof of part \ref{GenPres2}.

For part \ref{GenPres3}, the observations made in the proof of part \ref{GenPres1} imply that the identity map on $S'$ is a strong embedding from $\mathcal{F}'$ to $\mathcal{F}$ (Definition \ref{PossMorph}.\ref{embed}), so by Proposition \ref{Preservation}.\ref{Preservation2}, $\mathcal{F}\Vdash\varphi$ implies $\mathcal{F}'\Vdash\varphi$.\end{proof}
 
The following proposition shows that being a selective subframe is equivalent to being the image of a strict strong embedding as in Definition \ref{PossMorph}. 

\begin{proposition}[Embeddings and Subframes]\label{TautEm} For any possibility frames $\mathcal{F}$ and $\mathcal{G}$, the following are equivalent:
\begin{enumerate}
\item\label{TautEm1} there is a strict strong embedding from $\mathcal{F}$ into $\mathcal{G}$; 
\item\label{TautEm2} $\mathcal{F}$ is isomorphic to a selective subframe of $\mathcal{G}$.
\end{enumerate}
In addition, for any possibility frames $\mathcal{F}$ and $\mathcal{G}$, the following are equivalent:
\begin{enumerate}
\item[3.] there is a p-morphism that is a strong embedding from $\mathcal{F}$ to $\mathcal{G}$;
\item[4.] $\mathcal{F}$ is isomorphic to a generated subframe of $\mathcal{G}$.
\end{enumerate}
\end{proposition}
\begin{proof} From \ref{TautEm1} to \ref{TautEm2}, where $h$ is the strict strong embedding of $\mathcal{F}$ into $\mathcal{G}$, we first claim that the subframe $\mathcal{G}'$ of $\mathcal{G}$ with domain $h[S^\mathcal{F}]$ is a selective subframe of $\mathcal{G}$. Suppose $x^\mathcal{G}\in h[S^\mathcal{F}]$ and $y^\mathcal{G}\sqsubseteq^\mathcal{G} x^\mathcal{G}$. Thus, $x^\mathcal{G}=h(x^\mathcal{F})$ for some $x^\mathcal{F}\in \mathcal{F}$. Then since $y^\mathcal{G}\sqsubseteq^\mathcal{G} h(x^\mathcal{F})$ and $h$ is a strict possibility morphism, by \SqBack{} there is a $y^\mathcal{F}\in\mathcal{F}$ such that $h(y^\mathcal{F})\sqsubseteq^\mathcal{G} y^\mathcal{G}$. Then taking $z^\mathcal{G}=h(y^\mathcal{F})$, we have $z^\mathcal{G}\in h[S^\mathcal{F}]$ and $z^\mathcal{G}\sqsubseteq^\mathcal{G} y^\mathcal{G}$, so $\mathcal{G}'$ satisfies part \ref{SelSub1} of Definition \ref{Subframes}. Next, suppose $x^\mathcal{G}=h(x^\mathcal{F})R_i^\mathcal{G} y^\mathcal{G}$ and $u^\mathcal{G}\sqsubseteq^\mathcal{G} y^\mathcal{G}$. Then by \SRBack{}, there is a $y^\mathcal{F}\in\mathcal{F}$ such that $x^\mathcal{F}R_i^\mathcal{F}y^\mathcal{F}$ and $h(y^\mathcal{F})\between^\mathcal{G} u^\mathcal{G}$. By \RForth{}, $x^\mathcal{F}R_i^\mathcal{F}y^\mathcal{F}$ implies $h(x^\mathcal{F})R_i^\mathcal{G}h(y^\mathcal{F})$. Then taking $z^\mathcal{G}=h(y^\mathcal{F})$, we have $z^\mathcal{G}\in h[S^\mathcal{F}]$, $z^\mathcal{G}\between^\mathcal{G}u^\mathcal{G}$, and $x^\mathcal{G}R_i^\mathcal{G} z^\mathcal{G}$, so $\mathcal{G}'$ satisfies part \ref{SelSub2} of Definition \ref{Subframes}.
 
Since $h$ is an injection from $\mathcal{F}$ to $\mathcal{G}$, it is a bijection between $\mathcal{F}$ and our subframe $\mathcal{G}'=\langle h[S^\mathcal{F}],\sqsubseteq',\{R_i'\}_{i\in \ind},\adm'\rangle$. We claim that $h$ is an isomorphism between $\mathcal{F}$ and $\mathcal{G}'$. Since $h$ is a strong embedding, we already have that $y\sqsubseteq x$ iff $h(y)\sqsubseteq' h(x)$, and $xR_iy$ iff $h(x)R_i' h(y)$. Next, we must show that $h\colon \mathcal{F}\to\mathcal{G}'$ satisfies \PullBack{}, so for all $X'\in \adm'$, $h^{-1}[X']\in\adm^\mathcal{F}$. If $X'\in\adm'$, then since $\mathcal{G}'$ is a subframe of $\mathcal{G}$, and the domain of $\mathcal{G}'$ is $h[S^\mathcal{F}]$, by Definition \ref{Subframes}.\ref{SubAd} there is an $X\in \adm^\mathcal{G}$ such that $X'=X\cap h[S^\mathcal{F}]$. Since $h$ is a possibility morphism from $\mathcal{F}$ to $\mathcal{G}$, $X\in \adm^\mathcal{G}$ implies $h^{-1}[X]\in\adm^\mathcal{F}$. Since $\mathcal{F}$ is a possibility frame, we also have $S^\mathcal{F}\in\adm^\mathcal{F}$. Then since $\adm^\mathcal{F}$ is closed under $\cap$, we have $h^{-1}[X]\cap S^\mathcal{F}\in \adm^\mathcal{F}$, which with $h^{-1}[X']=h^{-1}[X\cap h [S^\mathcal{F}]]=h^{-1}[X]\cap h^{-1}[h[S^\mathcal{F}]]=h^{-1}[X]\cap S^\mathcal{F}$ gives us $h^{-1}[X']\in \adm^\mathcal{F}$.

Finally, we must show that for all $X\in \adm^\mathcal{F}$, $h[X]\in \adm'$. Since $h$ is a strong embedding from $\mathcal{F}$ to $\mathcal{G}$, for all $X\in \adm^\mathcal{F}$, there is an $X^\mathcal{G}\in \adm^\mathcal{G}$ such that $h[X]=X^\mathcal{G}\cap h[S^\mathcal{F}]$. Then since $\mathcal{G}'$ is a subframe of $\mathcal{G}$ with domain $h[S^\mathcal{F}]$, we have $h[X]\in\adm'$. This completes the proof that $h$ is an isomorphism.

We leave the proof from \ref{TautEm2} to \ref{TautEm1} to the reader. Note how Definition \ref{Subframes}.\ref{SelSub1}-\ref{SelSub2} is used to show that the isomorphism from $\mathcal{F}$ to the selective subframe of $\mathcal{G}$ satisfies \SqBack{} and \SRBack{} with respect to $\mathcal{G}$.

The proof that 3 and 4 are equivalent is an easy adaptation of the arguments for 1 and 2.\end{proof} 

Now we can relate selective subframes and homomorphic images as follows. 

\begin{proposition}[Selective Subframes and Homomorphic Images]\label{Gen&Hom} For any BAOs $\mathbb{A}$ and $\mathbb{B}$:
\begin{enumerate}[label=\arabic*.,ref=\arabic*]
\item\label{Gen&Hom1} if $\mathbb{A}$ is a homomorphic image of $\mathbb{B}$, then $\mathbb{A}_\gff$ is isomorphic to a generated subframe of $\mathbb{B}_\gff$;
\item\label{Gen&Hom2} if $\mathbb{A}$ is a homomorphic image of $\mathbb{B}$, then $\mathbb{A}_\ff$ is isomorphic to a generated subframe of $\mathbb{B}_\ff$.
\end{enumerate}
Conversely, for any possibility frames $\mathcal{F}$ and $\mathcal{G}$:
\begin{enumerate}[label=\arabic*.,ref=\arabic*,resume]
\item\label{Gen&Hom3} if $\mathcal{F}$ is isomorphic to a selective subframe of $\mathcal{G}$, then $\mathcal{F}^\under$ is a homomorphic image of $\mathcal{G}^\under$.
\end{enumerate}
\end{proposition}
 
\begin{proof} For part \ref{Gen&Hom1}, if there is a surjective homomorphism from $\mathbb{B}$ to $\mathbb{A}$, then by Theorem \ref{BAOMorphPossMorphII}, there is a p-morphism that is a strong embedding from $\mathbb{A}_\gff$ into $\mathbb{B}_\gff$, so by Proposition \ref{TautEm}, $\mathbb{A}_\gff$ is isomorphic to a generated subframe of $\mathbb{B}_\gff$. The proof of part \ref{Gen&Hom2} is the same but using part \ref{BAOMorphPossMorphIId} of Theorem  \ref{BAOMorphPossMorphII} instead of parts \ref{BAOMorphPossMorphIIa} and \ref{BAOMorphPossMorphIIb}.

For part \ref{Gen&Hom3}, if $\mathcal{F}$ is isomorphic to a selective subframe of $\mathcal{G}$, then by Proposition \ref{TautEm}, there is a strong embedding from $\mathcal{F}$ into $\mathcal{G}$, so by Theorem \ref{PossMorphBAOMorph}, there is a surjective homomorphism from $\mathcal{G}^\under$ to $\mathcal{F}^\under$.
\end{proof}

A third validity preserving construction for possibility frames is given by the notion of \textit{disjoint unions} of possibility frames, which parallels the standard definition for world frames \citep[\S~1.6]{Goldblatt1974}.
 
\begin{definition}[Disjoint Union]\label{DisUn} Given a nonempty indexed family $\{\mathcal{F}_j\}_{j\in J}$ of possibility frames $\mathcal{F}_j=\langle S_j,\sqsubseteq_j,\{R_{i,j}\}_{i\in\ind},P_j\rangle$, let $\mathcal{F}_j'=\langle S_j',\sqsubseteq_j',\{R_{i,j}'\}_{i\in\ind},P_j'\rangle$ be the isomorphic copy of $\mathcal{F}_j$ with domain $S'_j=S_j\times \{j\}$, so that $\{\mathcal{F}'_j\}_{j\in J}$ is a family of pairwise disjoint possibility frames. Then the \textit{disjoint union} of $\{\mathcal{F}_j\}_{j\in J}$ is the tuple $\underset{j\in J}{\biguplus}\mathcal{F}_j=\langle S,\sqsubseteq, \{R_i\}_{i\in\ind},P\rangle$ defined by: 
\begin{enumerate}
\item\label{DisUn1} $S= \underset{j\in J}{\bigcup} S_j'$; $\sqsubseteq \,= \underset{j\in J}{\bigcup}\sqsubseteq_j'$; $R_i=\underset{j\in J}{\bigcup}R_{i,j}'$; 
\item\label{DisUn2} $P= \{X\subseteq S\mid \forall j\in J\colon X\cap S_j'\in P_j'\}$.
\end{enumerate}
\end{definition}

\begin{proposition}[Preservation Under Disjoint Unions]\label{DisPres} For any nonempty indexed family $\{\mathcal{F}_j\}_{j\in J}$ of possibility frames:
\begin{enumerate}
\item\label{DissPoss1} $\underset{j\in J}{\biguplus}\mathcal{F}_j$ is a possibility frame; 
\item\label{DissPoss2} if each $\mathcal{F}_j$ is full, then $\underset{j\in J}{\biguplus}\mathcal{F}_j$ is full;
\item\label{DissPoss3} for all $\varphi\in\mathcal{L}(\sig,\ind)$, $\underset{j\in J}{\biguplus}\mathcal{F}_j\Vdash\varphi$ iff for all $j\in J$, $\mathcal{F}_j\Vdash \varphi$.
\end{enumerate}
\end{proposition}

\begin{proof} For part \ref{DissPoss1}, we must first show that the set $\adm$ of admissible propositions in $\underset{j\in J}{\biguplus}\mathcal{F}_j$ contains $\emptyset$, which is immediate from Definition \ref{DisUn}.\ref{DisUn2} and $\emptyset\in \adm'_j$, and is closed under the operations $\cap$, $\supset$, and $\blacksquare_i$ from Definition \ref{PosetMod}. By the disjointness of the $S_j'$, for all $X,Y\in\adm$ we have:
\begin{itemize}
\item[(i)] $X\cap Y = \underset{j\in J}{\bigcup} (X\cap Y \cap S'_j)$;
\item[(ii)] $X\supset Y =  \underset{j\in J}{\bigcup} ((X\cap S'_j)\supset_j' (Y\cap S'_j))$;
\item[(iii)] $\blacksquare_i X = \underset{j\in J}{\bigcup} \blacksquare_{i,j}' (X\cap S_j')$.
\end{itemize}
If $X,Y\in \adm$, then by definition of $\adm$, for all $j\in J$, $X\cap S'_j\in \adm'_j$ and $Y\cap S'_j\in \adm'_j$. Since $\adm_j'$ is closed under $\cap$, $\supset'_j$, and $\blacksquare_{i,j}'$, it follows that $X\cap Y\cap S'_j\in \adm'_j$, $(X\cap S'_j)\supset_j' (Y\cap S'_j)\in \adm'_j$, and $ \blacksquare_{i,j}' (X\cap S_j')\in\adm'_j$. Since $S_j'$ is the whole domain of $\mathcal{F}'_j$, $(X\cap S'_j)\supset_j' (Y\cap S'_j)=(X\supset'_jY)\cap S'_j$ and $ \blacksquare_{i,j}' (X\cap S_j')=(\blacksquare_{i,j}' X)\cap S_j'$.  The previous two sentences, equations (i)-(iii), and the definition of $\adm$ jointly imply that $X\cap Y,X\supset Y,\blacksquare_iX\in\adm$.

Next, we must show that each $X\in \adm$ satisfies \textit{persistence} and \textit{refinability} with respect to $\sqsubseteq$. For \textit{persistence}, suppose $x\in X$ and $y\sqsubseteq x$. Then since $X=\underset{j\in J}{\bigcup} (X\cap S_j')$, we have $x\in X\cap S_j'$ for some $j\in J$. By definition of $\adm$, $X\cap S_j'\in P'_j$, so $X\cap S_j'$ satisfies \textit{persistence} with respect to $\sqsubseteq_j'$. Since $\sqsubseteq \,= \underset{j\in J}{\bigcup}\sqsubseteq_j'$, $y\sqsubseteq x$, and $x\in S_j'$, it follows by the disjointness of the $\sqsubseteq'_j$ relations that $y\sqsubseteq'_j x$, which with the previous sentence implies that $y\in X\cap S_j'$, so $y\in X$. Thus, we have shown that $X$ satisfies \textit{persistence} with respect to $\sqsubseteq$.

For \textit{refinability}, assume that for all $y\sqsubseteq x$ there is a $z\sqsubseteq y$ such that $z\in X$. Consider such a $z$.  As above, $z\in X\cap S_j'$ for some $j\in J$, which with $z\sqsubseteq x$ implies that for all $y\sqsubseteq x$, $y\sqsubseteq'_j x$ and $y\in S'_j$. Thus, from our initial assumption, it follows that for all $y\sqsubseteq'_j x$ there is a $z\sqsubseteq'_j y$ such that $z\in X\cap S_j'$. Since $X\cap S_j'\in \adm'_j$, $X\cap S_j'$ satisfies \textit{refinability} with respect to $\sqsubseteq'_j$, so the previous sentence implies $x\in X\cap S_j'$, so $x\in X$. Thus, we have shown that $X$ satisfies \textit{refinability} with respect to $\sqsubseteq$.

For part \ref{DissPoss2}, assume that each $\mathcal{F}_j$ is full, which implies that each $\mathcal{F}_j'$ is full. Then the fact that $\underset{j\in J}{\biguplus}\mathcal{F}_j$ is full is a consequence of the following, which is easy to check: if $X\in \mathrm{RO}(\underset{j\in J}{\biguplus}\mathcal{F}_j)$, then for all $j\in J$, $X\cap S'_j\in \mathrm{RO}(\mathcal{F}_j')$, so $X\cap S'_j\in P'_j$ since $\mathcal{F}'_j$ is full. Then by the definition of $\adm$ in $\underset{j\in J}{\biguplus}\mathcal{F}_j$, $X\in \adm$.

For part \ref{DissPoss3}, the obvious embedding of $\mathcal{F}_j$ into $\underset{j\in J}{\biguplus}\mathcal{F}_j$ gives us that $\underset{j\in J}{\biguplus}\mathcal{F}_j\Vdash\varphi$ implies $\mathcal{F}_j\Vdash\varphi$ by Proposition \ref{Preservation}.\ref{Preservation2}. In the other direction, suppose $\varphi$ is falsified at an $x$ in $\underset{j\in J}{\biguplus}\mathcal{F}_j$ with an admissible valuation $\pi$. Then $x\in\mathcal{F}_j'$ for some $j\in J$; the restriction $\pi_j$ of $\pi$ to $\mathcal{F}_j'$ is an admissible valuation for $\mathcal{F}_j'$ by Definition \ref{DisUn}.\ref{DisUn2}; and the identity map on $\mathcal{F}_j'$ is a possibility morphism from $\langle\mathcal{F}_j',\pi_j\rangle$ to $\langle\underset{j\in J}{\biguplus}\mathcal{F}_j,\pi\rangle$ as in Definition \ref{PossMorph}.\ref{PossMorphAtomic}, so $\langle\underset{j\in J}{\biguplus}\mathcal{F}_j ,\pi\rangle, x\nVdash\varphi$ implies $\langle\mathcal{F}_j',\pi_j\rangle, x\nVdash \varphi$ by Proposition \ref{Preservation}.\ref{Preservation1}, so $\mathcal{F}_j\nVdash\varphi$ since $\mathcal{F}_j$ and $\mathcal{F}_j'$ are isomorphic.\end{proof} 

As usual, we can relate disjoint unions and direct products as follows. 

\begin{proposition}[Disjoint Unions and Direct Products]\label{Dis&Prod} For any nonempty indexed family $\{\mathcal{F}_j\}_{j\in J}$ of possibility frames, $(\underset{j\in J}{\biguplus}\mathcal{F}_j)^\under$ is BAO-isomorphic to $\underset{j\in J}\prod \mathcal{F}_j^\under$.
\end{proposition}
\begin{proof} As in Definition \ref{DisUn}, $\{\mathcal{F}_j'\}_{j\in J}$ is the family of pairwise disjoint frames such that for each $j\in J$, $\mathcal{F}_j'$ is isomorphic to $\mathcal{F}_j$, which is used to construct $\underset{j\in J}{\biguplus}\mathcal{F}_j$. Since $\underset{j\in J}\prod \mathcal{F}_j'^\under$ is obviously BAO-isomorphic to $\underset{j\in J}\prod \mathcal{F}_j^\under$, it suffices to show that $(\underset{j\in J}{\biguplus}\mathcal{F}_j)^\under$  is BAO-isomorphic to $\underset{j\in J}\prod \mathcal{F}_j'^\under$. Where $\adm$ is the set of admissible propositions in $\underset{j\in J}{\biguplus}\mathcal{F}_j$, and $\adm_j'$ is the set of admissible propositions in $\mathcal{F}_j'$, we define a map $h\colon \adm \to \underset{j\in J}{\prod}\adm_j'$ so that for $X\in \adm$, $h(X)$ is the element of $\underset{j\in J}{\prod}\adm_j'$ whose value at $j$ is $X\cap S'_j$, so $h$ is clearly a bijection by Definition \ref{DisUn}.\ref{DisUn2}. That $h$ is a BAO-homomorphism follows from equations (i)-(iii) in the proof of Proposition \ref{DisPres}.\ref{DissPoss1}. For example, for $\blacksquare_i$, where $f_k$ is the value of function $f$ at $k$, $\mathbb{A}=(\underset{j\in J}{\biguplus}\mathcal{F}_j)^\under$, and $\mathbb{B}= \underset{j\in J}\prod \mathcal{F}_j'^\under$: 
\begin{eqnarray*}
h(\blacksquare_i^\mathbb{A} X)_k & = & \blacksquare_i^\mathbb{A} X\cap S'_k\quad\mbox{by definition of }h \\
 & = & \underset{j\in J}{\bigcup} \blacksquare_{i,j}' (X\cap S_j') \cap S'_k\quad\mbox{by (iii)}\\
 & = & \blacksquare_{i,k}' (X\cap S_k') \quad\mbox{by the pairwise disjointness of the $S_j'$}\\
&=& \blacksquare_{i,k}' h(X)_k\quad\mbox{by definition of }h\\
&=& (\blacksquare_i^\mathbb{B}h(X))_k\quad\mbox{by definition of the direct product},
\end{eqnarray*}
so $h(\blacksquare_i^\mathbb{A} X)=\blacksquare_i^\mathbb{B}h(X)$. The other cases are similar.\end{proof} 

We will see in \S\S~\ref{DefViaDual}-\ref{DefViaDual2} how the above connections between frame constructions and algebraic constructions can be applied to the study of modally definable classes of frames.

\section{Beginnings of Definability \& Correspondence Theory}\label{special}

In this section, we turn to modal \textit{definability theory}, including \textit{correspondence theory}, in the context of possibility semantics (cf.~\citealt{Benthem1983,Benthem1980} on definability theory for possible world semantics). Three classic questions concerning modal definability are the following from \citealt{Benthem1983} (p.~13):
\begin{enumerate}
\item When does a given modal formula define a first-order property of the relations in frames? 
\item When can a given first-order property of frames be defined by means of modal formulas?
\item Which classes of frames are definable at all by means of modal formulas?
\end{enumerate}
We will take these questions in reverse order: we discuss question 3 for possibility frames in \S~\ref{DefViaDual}, question 2 for full possibility frames in \S~\ref{DefViaDual2}, and question 1 for full possibility frames in \S~\ref{LemmScottCorr}. 

To make these questions precise, we need to fix the relevant notions of definability and of first-orderness. We begin with definability, for which the general notion is that of relative definability.

\begin{definition}[Relative Modal Frame Definability]\label{RelDef} Let $\mathsf{F}$ and $\mathsf{G}$ be classes of possibility frames. 

A set $\Sigma\subseteq\mathcal{L}(\sig,\ind)$ of modal formulas \textit{defines $\mathsf{F}$ relative to $\mathsf{G}$} iff for all $\mathcal{F}\in\mathsf{G}$, $\mathcal{F}\in\mathsf{F}$ iff every $\varphi\in\Sigma$ is valid over $\mathcal{F}$. $\mathsf{F}$ is \textit{modally definable} (or \textit{modal axiomatic}) relative to $\mathsf{G}$ iff there is some $\Sigma\subseteq\mathcal{L}(\sig,\ind)$ that defines $\mathsf{F}$ relative to $\mathsf{G}$. \hfill $\triangleleft$
\end{definition} 

For question 3 above, in \S~\ref{DefViaDual} we give a characterization of the classes of possibility frames that are modally definable relative to the class of all possibility frames. Theorem \ref{ModDef1} in \S~\ref{DefViaDual} is the analogue for possibility semantics of Goldblatt's \citeyearpar[Thm.~1.12.11]{Goldblatt1974} characterization of modally definable classes of world frames. 

Turning to the notion of first-orderness, let $\mathcal{L}^1(\ind)$ be the first-order language with equality that contains binary relation symbols $\dot{\sqsubseteq}$ and $\dot{R_i}$ for each $i\in\ind$. We interpret $\mathcal{L}^1(\ind)$ in possibility frames $\mathcal{F}=\langle S,\sqsubseteq,\{R_i\}_{i\in\ind},\adm\rangle$ in the obvious way, interpreting $\dot{\sqsubseteq}$ as the refinement relation $\sqsubseteq$ and $\dot{R}_i$ as the accessibility relation $R_i$. The set $\adm$ of admissible propositions plays no role in interpreting $\mathcal{L}^1(\ind)$, so instead of talking of \textit{frames} we can talk of \textit{foundations} $F=\langle S,\sqsubseteq, \{R_i\}_{i\in\ind}\rangle$ as in Definition \ref{Foundation}.

Foundations serve as standard first-order structures for $\mathcal{L}^1(\ind)$, so we can apply to foundations standard notions of $\mathcal{L}^1(\ind)$-definability, $\mathcal{L}^1(\ind)$-elementary equivalence, etc. From Corollary \ref{Fullinterplay}.\ref{Fullinterplay2}, we have the following definability result. 

\begin{corollary}[First-Order Definability of Foundations of Full Frames]\label{DefFullFrames} The class of foundations $F=\langle S,\sqsubseteq, \{R_i\}_{i\in\ind}\rangle $ for which $F^\full =\langle S,\sqsubseteq, \{R_i\}_{i\in\ind},\mathrm{RO}(S,\sqsubseteq)\rangle$ is a full possibility frame is definable by a set of $\mathcal{L}^1(\ind)$-sentences relative to the class of all $\mathcal{L}^1(\ind)$-structures, viz., by the $\mathcal{L}^1(\ind)$-sentences expressing that $\sqsubseteq$ is a partial order and that $\sqsubseteq$ and $R_i$ for each $i\in\ind$ satisfy the conditions \Rrule{} and \Rwinweak{} from \S~\ref{FullFrames}. \hfill $\triangleleft$
\end{corollary}

Using the notion of foundations, we can also talk about $\mathcal{L}^1(\ind)$-definability and $\mathcal{L}^1(\ind)$-elementary equivalence for \textit{frames} as follows.

\begin{definition}[$\mathcal{L}^1(\ind)$-Definability and Equivalence of Frames]\label{FOdef} Let $\mathsf{F}$ and $\mathsf{G}$ be classes of possibility frames.  

A set $\Sigma$ of $\mathcal{L}^1(\ind)$-sentences \textit{defines $\mathsf{F}$ relative to $\mathsf{G}$} iff for all $\mathcal{F}\in\mathsf{G}$, $\mathcal{F}\in\mathsf{F}$ iff every $\varphi\in\Sigma$ is true in $\mathcal{F}_\found$.

Possibility frames $\mathcal{F}$ and $\mathcal{G}$ are \textit{$\mathcal{L}^1(\ind)$-elementarily equivalent} iff $\mathcal{F}_\found$ and $\mathcal{G}_\found$ are $\mathcal{L}^1(\ind)$-elementarily equivalent. \textsf{F} is \textit{closed under $\mathcal{L}^1(\ind)$-elementary equivalence relative to $\mathsf{G}$} iff whenever $\mathcal{F}\in\mathsf{F}$, $\mathcal{G}\in\mathsf{G}$, and $\mathcal{F}$ is $\mathcal{L}^1(\ind)$-elementarily equivalent to $\mathcal{G}$, then $\mathcal{G}\in\mathsf{F}$.\hfill $\triangleleft$\end{definition}

For succinctness, we will adopt the following convention: when we say that a class of \textit{full} possibility frames is simply `definable', in the modal or $\mathcal{L}^1(\ind)$ sense, or that it is simply `closed under $\mathcal{L}^1(\ind)$-elementary equivalence', we mean \textit{relative to the class of full possibility frames}. 

In \S~\ref{DefViaDual2}, we answer question 2 above for full possibility frames. We give a characterization of when a class of full possibility frames that is closed under $\mathcal{L}^1(\ind)$-elementary equivalence is modally definable. Theorem \ref{ModDef2} in \S~\ref{DefViaDual2} is the analogue for possibility semantics of the well-known Goldblatt-Thomason Theorem (\citealt[Thm.~8]{Goldblatt1975}, \citealt[Thm.~3.19]{Blackburn2001}) about when a class of full world frames closed under elementary equivalence is modally definable (relative to the class of full world frames). Stating the result for classes closed under elementary equivalence is of course stronger than what question 2 asks for, because while first-order definability by a set of sentences implies closure under elementary equivalence, the converse implication does not hold. (Recall that closure under elementary equivalence is equivalent to the weaker condition of being the union of classes each of which is first-order definable by a set of sentences.) 

We will postpone discussion of question 1 and our relevant results until \S~\ref{LemmScottCorr}.
 
\subsection{Modally Definable Classes of Possibility Frames}\label{DefViaDual}

In this section, we use the results of \S~\ref{OpFrameAlg} to prove an analogue for possibility frames of Goldblatt's \citeyearpar[Thm.~1.12.11]{Goldblatt1974} characterization of modally definable classes of world frames. 

Goldblatt's \citeyearpar[Thm.~1.12.11]{Goldblatt1974} theorem states that a class $\mathsf{K}$ of world frames is modally definable \textit{if and only} if it is closed under surjective p-morphisms, generated subframes, and disjoint unions, while both $\mathsf{K}$ and its complement are closed under general ultrafilter extensions (i.e., taking the general ultrafilter frame of the underlying BAO of a frame, as in \S~\ref{AlgSem}). The left-to-right direction simply uses the fact that the validity of modal formulas is preserved by surjective p-morphisms, generated subframes, and disjoint unions, while a frame validates exactly the same formulas as its general ultrafilter extension. For the right-to-left direction, Goldblatt's strategy for showing that the modal logic of \textsf{K} \textit{defines} \textsf{K} involves switching from the universe of frames to the universe of BAOs, where we can then apply Birkhoff's \citeyearpar{Birkhoff1935} HSP theorem: the smallest equationally definable class of algebras containing a given class $\mathsf{C}$ of algebras is $\mathbf{HSP}(\mathsf{C})$, the class of all \textbf{h}omomorphic images of \textbf{s}ubalgebras of direct \textbf{p}roducts of algebras from $\mathsf{C}$. Given suitable connections between algebraic constructions and frame constructions (recall \S~\ref{OpFrameAlg}), one can then use the HSP theorem to show that any frame validating the modal logic of $\mathsf{K}$ in fact belongs to $\mathsf{K}$. This now standard strategy is exactly the strategy we will follow below. Compare the proof of Theorem \ref{ModDef1} below to van Benthem's \citeyearpar[Thm.~16.1]{Benthem1983} proof of Goldblatt's \citeyearpar[Thm.~1.12.11]{Goldblatt1974} characterization of modally definable classes of world frames. Also note that the HSP theorem is a theorem of ZF set theory \citep{Andreka1981}, so Theorem \ref{ModDef1} does not require the axiom of choice or even the ultrafilter axiom.

\begin{theorem}[Modal Definability of Possibility Frames]\label{ModDef1} $\,$
\begin{enumerate}
\item\label{ModDef1a} If a class $\mathsf{F}$ of possibility frames is definable by modal formulas, then $\mathsf{F}$ is closed under dense possibility morphisms, selective subframes, and disjoint unions, while both $\mathsf{F}$ and its complement are closed under general filter extensions (see \S~\ref{Fdes}).
\item\label{ModDef1b} If a class $\mathsf{F}$ of possibility frames is closed under surjective p-morphisms, generated subframes, and disjoint unions, while both $\mathsf{F}$ and its complement are closed under general filter extensions, then $\mathsf{F}$ is modally definable.
\end{enumerate}
\end{theorem}

\begin{proof} For part \ref{ModDef1a}, if $\mathsf{F}$ is defined by a set $\Sigma$ of modal formulas, then since the validity of modal formulas is preserved by dense possibility morphisms (Proposition \ref{DenseRecall}), selective subframes (Proposition \ref{GenPres}), disjoint unions (Proposition \ref{DisPres}), and general filter extensions (Corollary \ref{Arbitrary-to-FD}), $\mathsf{F}$ is closed under these constructions. Since validity is also reflected by general filter extensions, i.e., if $\varphi$ is valid over the general filter extension of a frame, then $\varphi$ is valid over the frame (Corollary \ref{Arbitrary-to-FD}), it follows that the complement of $\mathsf{F}$ is closed under general filter extensions. For if $\mathcal{F}\not\in\mathsf{F}$, so it does not validate some formula from $\Sigma$, then $(\mathcal{F}^\under)_\gff$ does not validate that formula from $\Sigma$, so  $(\mathcal{F}^\under)_\gff\not\in\mathsf{F}$.

For part \ref{ModDef1b}, assume $\mathsf{F}\neq\emptyset$; otherwise $\mathsf{F}$ is defined by $\{\bot\}$. Let $\Sigma$ be the modal logic of $\mathsf{F}$, i.e., the set of all $\varphi\in\mathcal{L}(\sig,\ind)$ valid over every frame in $\mathsf{F}$. We claim that $\Sigma$ defines $\mathsf{F}$: a frame $\mathcal{F}$ is in $\mathsf{F}$ iff $\mathcal{F}$ validates $\Sigma$. The left-to-right direction is immediate. For the right-to-left direction, as usual, we can turn the set $\Sigma$ of modal formulas into a set $\Sigma^\under$ of algebraic equations, and then since $\Sigma$ is the modal logic of $\mathsf{F}$, $\Sigma^\under$ is the algebraic equational theory of $\mathsf{F}^\under =\{\mathcal{F}^\under\mid \mathcal{F}\in\mathsf{F}\}$ by Theorem \ref{PtoB}.\ref{PtoB6}. Then the class of all BAOs in which the equations of $\Sigma^\under$ hold is the smallest equationally definable class of BAOs containing $\mathsf{F}^\under$, which by  Birkhoff's theorem is $\mathbf{HSP}(\mathsf{F}^\under)$, the class of all homomorphic images of subalgebras of direct products of algebras from~$\mathsf{F}^\under$. 

Now suppose  $\mathcal{F}$ is a possibility frame validating $\Sigma$, so by the previous paragraph, the equations of $\Sigma^\under$ hold in $\mathcal{F}^\under$ and hence $\mathcal{F}^\under\in \mathbf{HSP}(\mathsf{F}^\under)$. Thus, there is a nonempty family $\{\mathcal{F}_j\}_{j\in J}$ of possibility frames from $\mathsf{F}$ such that $\mathcal{F}^\under$ is a homomorphic image of a subalgebra $\mathbb{S}$ of $\underset{j\in J}\prod \mathcal{F}_j^\under$. Since $\mathcal{F}^\under$ is a homomorphic image of $\mathbb{S}$, by Proposition \ref{Gen&Hom}.\ref{Gen&Hom1}, $(\mathcal{F}^\under)_\gff$ is isomorphic to a generated subframe of $\mathbb{S}_\gff$. Next, from Proposition \ref{Dis&Prod}, $\underset{j\in J}\prod \mathcal{F}_j^\under$ is isomorphic to $(\underset{j\in J}{\biguplus}\mathcal{F}_j)^\under$. Thus, the subalgebra $\mathbb{S}$ of $\underset{j\in J}\prod \mathcal{F}_j^\under$ is isomorphic to a subalgebra of $(\underset{j\in J}{\biguplus}\mathcal{F}_j)^\under$, whence by Proposition \ref{Poss&Sub}.\ref{Poss&Sub1}, there is a surjective p-morphism from  $((\underset{j\in J}{\biguplus}\mathcal{F}_j)^\under)_\gff$ to $\mathbb{S}_\gff$.

Now we argue that $\mathcal{F}\in\mathsf{F}$. Since each $\mathcal{F}_j$ is in $\mathsf{F}$, $\underset{j\in J}{\biguplus}\mathcal{F}_j\in\mathsf{F}$ by closure under disjoint unions; then $((\underset{j\in J}{\biguplus}\mathcal{F}_j)^\under)_\gff\in\mathsf{F}$ by closure under general filter extensions; then 
since there is a surjective p-morphism from $((\underset{j\in J}{\biguplus}\mathcal{F}_j)^\under)_\gff$ to $\mathbb{S}_\gff$, $\mathbb{S}_\gff\in\mathsf{F}$ by closure under surjective p-morphisms; then since $(\mathcal{F}^\under)_\gff$ is isomorphic to a generated subframe of $\mathbb{S}_\gff$, $(\mathcal{F}^\under)_\gff\in\mathsf{F}$ by closure under generated subframes and isomorphisms; and then finally $\mathcal{F}\in\mathsf{F}$ by the closure of the complement of $\mathsf{F}$ under general filter extensions. \end{proof}

\subsection{Modally Definable Classes of Full Possibility Frames}\label{DefViaDual2}

Next we will prove an analogue for \textit{full} possibility frames of the Goldblatt-Thomason Theorem for full world frames. Goldblatt and Thomason's \citeyearpar{Goldblatt1975} result concerns classes $\mathsf{K}$ of full world frames that are \textit{closed under elementary equivalence} in the sense that if $\langle \wo{W},\{\wo{R}_i\}_{i\in\ind}\rangle\in \mathsf{K}$, and $\langle \wo{W},\{\wo{R}_i\}_{i\in\ind}\rangle$ and $\langle \wo{W}',\{\wo{R}_i'\}_{i\in\ind}\rangle$ satisfy the same sentences of the first-order language with equality and a binary predicate $\dot{\mathrm{R}}_i$ for each $i\in \ind$, then $\langle \wo{W}',\{\wo{R}_i'\}_{i\in\ind}\rangle\in\mathsf{K}$. The theorem states that any class $\mathsf{K}$ of full world frames that is closed under elementary equivalence is modally definable relative to the class of full world frames iff $\mathsf{K}$ is closed under surjective p-morphisms, generated subframes, and disjoint unions, while the complement of $\mathsf{K}$ is closed under ultrafilter extensions (i.e., taking the ultrafilter frame of the underlying BAO of the frame, as in \S~\ref{AlgSem}). Note the differences between this result and Goldblatt's result on modally definable classes of (general) world frames. There is no closure under general ultrafilter extensions, because such extensions are not full frames (except in the trivial case of the general ultrafilter extension of a \textit{finite} full frame, which is always isomorphic to the initial frame). Thus, we consider ultrafilter extensions instead of general ultrafilter extensions. Recall that if $\varphi$ is valid over the ultrafilter extension of a frame, then $\varphi$ is valid over the frame, but the converse is not guaranteed, as it was in the case of general ultrafilter extensions. However, if $\mathsf{K}$ is closed under elementary equivalence as well as surjective p-morphisms, then $\mathsf{K}$ is closed under ultrafilter extensions \citep[Lem.~3.8]{Benthem1979b}. This fact allows the same form of argument as used in the proof of Goldblatt's theorem for world frames to be used in a proof of the Goldblatt-Thomason Theorem for full world frames. 

Our analogue of Goldblatt-Thomason concerns classes of full possibility frames that are closed under $\mathcal{L}^1(\ind)$-elementary equivalence as in Definition \ref{FOdef} (relative to the class of full possibility frames). The restriction to classes closed under elementary equivalence helps in a way analogous to the way it helps in the case of full world frames: any class of full possibility frames closed under $\mathcal{L}^1(\ind)$-elementary equivalence and dense (and strict) possibility morphisms is also closed under \textit{filter extensions}, as in Lemma \ref{JvBLem}. This fact will allow us to easily adapt the proof of Theorem \ref{ModDef1} above to give a proof of Theorem \ref{ModDef2} below.

To prove the key Lemma \ref{JvBLem}, we need some help from first-order model theory. Given a first-order structure $\mathfrak{A}$ for a language $\mathcal{L}$ and an element $a$ in $\mathfrak{A}$, let $(\mathcal{L},\dot{a})$ be the expansion of $\mathcal{L}$ with a new constant symbol $\dot{a}$, and let $(\mathfrak{A},a)$ be the $(\mathcal{L},\dot{a})$-expansion of $\mathfrak{A}$ that interprets $\dot{a}$ as $a$. A structure $\mathfrak{A}$ for $\mathcal{L}$ is \textit{2-saturated} iff for any element $a$ in $\mathfrak{A}$ and any set $\Sigma(\mathrm{x})$ of formulas of $(\mathcal{L},\dot{a})$ containing at most the variable $\mathrm{x}$ free, if every finite subset of $\Sigma(\mathrm{x})$ is satisfied by some object or other in $(\mathfrak{A},a)$, then $\Sigma(\mathrm{x})$ is satisfied by an object in $(\mathfrak{A},a)$. It is a standard result in first-order model theory that every model has a 2-saturated (indeed, $\omega$-saturated) elementary extension \citep[\S5.1]{Chang1990}. Here we assume the axiom of choice.\footnote{In fact, the result that every model has an $\omega$-saturated elementary extension does not require full choice. It suffices for this result to assume the ultrafilter axiom and the axiom of dependent choice (thanks to Dan Appel and Nick Ramsey for discussion of this point). Pincus \citeyearpar{Pincus1977} showed that ZF plus the ultrafilter axiom and dependent choice is strictly weaker than ZFC.}

\begin{lemma}[2-Saturated Elementary Extensions]\label{SatEx} For every full possibility frame $\mathcal{F}=\langle S,\sqsubseteq, \{R_i\}_{i\in\ind},\adm\rangle$, viewed as a structure for the first-order language $\mathcal{L}^1(\ind,\adm)$ with identity, binary predicate symbols $\dot{\sqsubseteq}$ and $\dot{R}_i$ for each $i\in\ind$, and a unary predicate symbol $\dot{X}$ for each $X\in \adm$, there is a $\mathcal{L}^1(\ind,\adm)$-structure $\mathcal{F}^\mathrm{s}=\langle S^\mathrm{s},\sqsubseteq^\mathrm{s}, \{R_i^\mathrm{s}\}_{i\in\ind},\adm^\mathrm{s}\rangle$  such that:
\begin{enumerate}
\item\label{SatEx1} $\mathcal{F}^\mathrm{s}$ is a 2-saturated $\mathcal{L}^1(\ind,\adm)$-elementary extension of $\mathcal{F}$;
\item\label{SatEx2} $\langle S^\mathrm{s},\sqsubseteq^\mathrm{s} , \{R_i^\mathrm{s}\}_{i\in\ind}\rangle$ is an $\mathcal{L}^1(\ind)$-elementary extension of $\langle S,\sqsubseteq, \{R_i\}_{i\in \ind}\rangle$;
\item\label{SatEx3} $\mathcal{F}'=\langle S^\mathrm{s},\sqsubseteq^\mathrm{s}, \{R_i^\mathrm{s}\}_{i\in\ind},\mathrm{RO}(S^\mathrm{s},\sqsubseteq^\mathrm{s})\rangle$ is a full possibility frame.
\end{enumerate}
\end{lemma}

\begin{proof} We have already noted that part \ref{SatEx1} follows from standard results in first-order model theory. Part \ref{SatEx2} is immediate from part \ref{SatEx1}. For part \ref{SatEx3}, by Proposition \ref{ROtoRO}, $\mathcal{F}'$ is a full possibility frame iff it satisfies the $\mathcal{L}^1(\ind)$-sentences expressing \Rrule{} and \Rwinweak{} for each $i\in\ind$ and that $\sqsubseteq$ is a partial order. Since $\mathcal{F}$ is a full possibility frame, $\mathcal{F}$ satisfies those sentences by Proposition \ref{ROtoRO}, so by part \ref{SatEx2}, $\mathcal{F}'$ also satisfies them.
\end{proof}

We can now prove the key lemma that will take us from Theorem \ref{ModDef1} to Theorem \ref{ModDef2}. It is based on \textsc{Lemma 3.8} of \citealt{Benthem1979b}, which is in turn based on \textsc{Lemma 9} of \citealt{Fine1975c}.

\begin{lemma}[From Frames to Filter Extensions]\label{JvBLem} For every full possibility frame $\mathcal{F}=\langle S,\sqsubseteq, \{R_i\}_{i\in\ind},\adm\rangle$, there is a full possibility frame $\mathcal{F}'=\langle S',\sqsubseteq', \{R_i'\}_{i\in\ind},P'\rangle$ such that $\langle S',\sqsubseteq',\{R_i'\}_{i\in\ind}\rangle$ is an $\mathcal{L}^1(\ind)$-elementary extension of $\langle S,\sqsubseteq, \{R_i\}_{i\in\ind}\rangle$ and there is a dense and strict possibility morphism from $\mathcal{F}'$ to the filter extension $(\mathcal{F}^\under)_\ff$ of $\mathcal{F}$.
\end{lemma} 

\begin{proof} Given $\mathcal{F}$, let $\mathcal{F}^\mathrm{s}$ and $\mathcal{F}'$ be as in Lemma \ref{SatEx}. Define a function $g$ with domain $S'$ by $g(x)=\{{X\in \adm} \mid x\in \dot{X}^{\mathcal{F}^\mathrm{s}}\}$, where $\dot{X}^{\mathcal{F}^\mathrm{s}}$ is the interpretation of the predicate symbol $\dot{X}$ in the $\mathcal{L}^1(\ind,\adm)$-structure $\mathcal{F}^\mathrm{s}$. We claim that $g$ is a dense and strict possibility morphism from $\mathcal{F}'$ to $(\mathcal{F}^\under)_\ff$. 

First, we check that $g(x)$ is in the domain of $(\mathcal{F}^\under)_\ff$, i.e., that $g(x)$ is a proper filter in $\mathcal{F}^\under$. For every $X,Y\in \adm$, where $\dot{(X\cap Y)}$ and $\dot{(-X)}$ are the predicate symbols of $\mathcal{L}^1(\ind,\adm)$ corresponding to the sets $X\cap Y$ and $-X=X\supset \emptyset$ in $\adm$, the $\mathcal{L}^1(\ind,\adm)$ sentences $\forall \mathrm{z}((\dot{X}(\mathrm{z})\wedge \dot{Y}(\mathrm{z}))\leftrightarrow \dot{(X\cap Y)}(\mathrm{z}))$ and $\forall\mathrm{z}( \dot{(-X)}(\mathrm{z})\rightarrow \neg \dot{X}(\mathrm{z}))$, as well as $\forall \mathrm{z} \dot{S}(\mathrm{z})$, are true in the possibility frame $\mathcal{F}$, regarded as an $\mathcal{L}^1(\ind,\adm)$-structure, and hence true in its $\mathcal{L}^1(\ind,\adm)$-elementary extension $\mathcal{F}^\mathrm{s}$ that is used to define $g$. It follows that $X,Y\in g(x)$ iff $X\cap Y\in g(x)$, and $g(x)\neq\emptyset$, so $g(x)$ is a \textit{filter} in $\mathcal{F}^\under$, and that $-X\in g(x)$ implies $X\not\in g(x)$, so $g(x)$ is a \textit{proper} filter.

Next, we show that $g$ satisfies the \textit{dense} condition, i.e., that for every $y'\in (\mathcal{F}^\under)_\ff$, there is a $y \in \mathcal{F}'$ such that $g(y)\sqsubseteq^{(\mathcal{F}^\under)_\ff} y'$. Consider any element of $(\mathcal{F}^\under)_\ff$, i.e., any filter $F$ in $\mathcal{F}^\under$. Let $\Sigma=\{\dot{X}(\mathrm{x})\mid X\in F\}$. For any finite $\Sigma_0\subseteq\Sigma$, let $F_0=\{X\mid \dot{X}(\mathrm{x})\in\Sigma_0\}$, so $F_0$ is a finite subset of $F$. Hence $\bigcap F_0\in F$ by the fact that $F$ is a filter, so then $\bigcap F_0\not=\emptyset$ by the fact that $F$ is a proper filter. Thus, there is an $x\in \bigcap F_0$, which means that $x$ satisfies $\Sigma_0$ in $\mathcal{F}$. Hence $\exists \mathrm{x} (\bigwedge \Sigma_0)$ is true in $\mathcal{F}$ and therefore in its $\mathcal{L}^1(\ind,\adm)$-elementary extension $\mathcal{F}^\mathrm{s}$. Thus, we have shown that every finite subset of $\Sigma$ is satisfied by some object or other in $\mathcal{F}^\mathrm{s}$. Then since $\mathcal{F}^\mathrm{s}$ is 2-saturated, it follows that $\Sigma$ is satisfied in $\mathcal{F}^\mathrm{s}$ by an object $y$, so $y\in\mathcal{F}'$ as well. Then by the definition of $g$,  $g(y)\supseteq F$, so $g(y)\sqsubseteq^{(\mathcal{F}^\under)_\ff} F$ by definition of $(\mathcal{F}^\under)_\ff$. Hence $g$ satisfies the \textit{dense} condition.

Finally, we show that $g$ is a strict possibility morphism:
\begin{itemize}
\item  if $y\sqsubseteq'  x $, then $g(y)\sqsubseteq^{(\mathcal{F}^\under)_\ff} g( x )$ (\SqForth{}); 
\item if $F\sqsubseteq^{(\mathcal{F}^\under)_\ff} g( x )$, then $\exists y$: $y\sqsubseteq' x $ and $g(y)\sqsubseteq^{(\mathcal{F}^\under)_\ff} F$ (\SqBack{});
\item if $ x R_i'y$, then $g(x)R_i^{(\mathcal{F}^\under)_\ff} g(y)$ (\RForth{}); 
\item if $g(x)R^{(\mathcal{F}^\under)_\ff}_i F$ and $G\sqsubseteq ^{(\mathcal{F}^\under)_\ff} F$, then $\exists y$: $xR_i' y$ and $g(y)\comp^{(\mathcal{F}^\under)_\ff} G$ (\SRBack{}).
\end{itemize}
Since $\mathcal{F}'$ is a full possibility frame, the \PullBack{} property of $g$ follows from \SqForth{} and \SqBack{} by Fact \ref{PullFact}.

For \SqForth{}, since $\mathcal{F}$ is a possibility frame, we have that for all $X\in\adm$, $\mathcal{F}$ satisfies the $\mathcal{L}^1(\ind,\adm)$ sentence $\forall \mathrm{x}\forall\mathrm{y}((\mathrm{y}\,\dot{\sqsubseteq}\,\mathrm{x}\wedge \dot{X}(\mathrm{x}))\rightarrow \dot{X}(\mathrm{y}))$, so its $\mathcal{L}^1(\ind,\adm)$-elementary extension $\mathcal{F}^\mathrm{s}$ does as well. Now if $y\sqsubseteq ' x$, so $y\sqsubseteq^\mathrm{s} x$, then since $\mathcal{F}^\mathrm{s}$ satisfies the given sentence, it follows by the definition of $g$ that $g(y)\supseteq g(x)$, so $g(y)\sqsubseteq^{(\mathcal{F}^\under)_\ff} g( x )$. 

The proof of \RForth{} is analogous, using for every $X\in \adm$ the sentence $\forall \mathrm{x}\forall\mathrm{y}((\mathrm{x}\dot{R}_i\mathrm{y}\wedge \dot{(\blacksquare_i X)}(\mathrm{x}))\rightarrow \dot{X}(\mathrm{y}))$, where $\dot{(\blacksquare_i X)}$ is the predicate symbol corresponding to $\blacksquare_i X$ in $\adm$, and then the definition of $R_i^{(\mathcal{F}^\under)_\ff}$.

For \SqBack{}, suppose that $F\sqsubseteq^{(\mathcal{F}^\under)_\ff} g( x )$, so $F$ is a filter in $\mathcal{F}^\under$ such that $F\supseteq g( x )$. Expand the language $\mathcal{L}^1(\ind,\adm)$ with a new constant symbol $\dot{x}$. Let $\Sigma =\{\dot{X}(\mathrm{y})\mid X\in F\}\cup \{\mathrm{y}\,\dot{\sqsubseteq}\, \dot{x}\}$. For any finite  $\Sigma_0\subseteq
\Sigma$, let $F_0=\{X\mid \dot{X}(\mathrm{y})\in \Sigma_0\}$ as above, so $\bigcap F_0\in F$ as above. Let $Q=\bigcap F_0$, so $Q\in F$. Then since $F$ is a proper filter and $F\supseteq g(x)$, it follows that $-Q\not\in g(x)$. Hence $x$ is not in the interpretation of the predicate symbol $\dot{(-Q)}$ in $\mathcal{F}^\mathrm{s}$. Now since the $\mathcal{L}^1(\ind,\adm)$ sentence $\forall \mathrm{x} (\neg (\dot{-Q})(\mathrm{x})\rightarrow \exists \mathrm{y} (\mathrm{y}\,\dot{\sqsubseteq}\,\mathrm{x}\wedge \dot{Q}(y)))$ is true in the possibility frame $\mathcal{F}$, it is true in its $\mathcal{L}^1(\ind,\adm)$-elementary extension $\mathcal{F}^\mathrm{s}$. Combining the previous two steps, we have that there is a $y\in \mathcal{F}^\mathrm{s}$ such that $y\sqsubseteq^\mathrm{s}x$, so $y\sqsubseteq ' x$, and $y\in Q=\bigcap F_0$. Hence $y$ satisfies the set $\Sigma_0\cup \{\mathrm{y} \,\dot{\sqsubseteq}\, \dot{x}\}$ in the expansion of $\mathcal{F}^\mathrm{s}$ that interprets the new constant $\dot{x}$ as $x$. Since $\Sigma_0$ was an arbitrary finite subset of $\Sigma$, we can apply the fact that $\mathcal{F}^\mathrm{s}$ is 2-saturated to conclude that there is an object $y$ that satisfies the whole set $\Sigma$ in the expansion of $\mathcal{F}^\mathrm{s}$. Then by the definition of $\Sigma$ and $g$, we have $g(y)\supseteq F$, so $g(y)\sqsubseteq^{(\mathcal{F}^\under)_\ff} F$. Finally, since $\mathrm{y}\,\dot{\sqsubseteq}\,\dot{x}\in\Sigma$, we have that $y\sqsubseteq^\mathrm{s} x$, so $y\sqsubseteq' x$. This completes the proof of \SqBack{}. 

For \SRBack{}, suppose $g( x )R^{(\mathcal{F}^\under)_\ff}_i F$ and $G\sqsubseteq^{(\mathcal{F}^\under)_\ff} F$, so for every $X\in\adm$, $\blacksquare_i X\in g(x)$ implies $X\in G$ by the definitions of $R^{(\mathcal{F}^\under)_\ff}_i$ and $\sqsubseteq^{(\mathcal{F}^\under)_\ff}$. As above, expand the language with a new constant $\dot{x}$, but this time let $\Sigma = \{\dot{X}(\mathrm{z})\mid X\in G\}\cup \{\exists\mathrm{y}(\dot{x}\dot{R}_i\mathrm{y}\wedge \mathrm{z}\,\dot{\sqsubseteq}\,\mathrm{y})\}$. For any finite $\Sigma_0\subseteq\Sigma$, let $G_0=\{X\mid \dot{X}(\mathrm{z})\in \Sigma_0\}$, so as above, where $Q=\bigcap G_0$, $Q\in G$. Then since $G$ is a proper filter, the first sentence of this paragraph implies that $\blacksquare_i \mathord{-}Q\not\in g(x)$, so $x$ is not in the interpretation of the predicate symbol $(\dot{\blacksquare_i \mathord{-}Q})$ in $\mathcal{F}^\mathrm{s}$. Now since the $\mathcal{L}^{1}(\ind,\adm)$ sentence $\forall \mathrm{x} (\neg (\dot{\blacksquare_i\mathord{-}Q})(\mathrm{x})\rightarrow \exists \mathrm{y} \exists\mathrm{z}(\mathrm{x}R_i\mathrm{y}\wedge \mathrm{z}\,\dot{\sqsubseteq}\,\mathrm{y}\wedge \dot{Q}(z)))$ is true in $\mathcal{F}$, it is true in $\mathcal{F}^\mathrm{s}$. As above, we then deduce that there is an object $z$ in $\mathcal{F}^\mathrm{s}$ that satisfies the whole of $\Sigma$ in the expansion of $\mathcal{F}^\mathrm{s}$ that interprets $\dot{x}$ as $x$. Thus, there are $y,z\in \mathcal{F}^\mathrm{s}$ such that $xR^\mathrm{s}_iy$ and $z\sqsubseteq^\mathrm{s} y$, so $xR_i'y$ and $z\sqsubseteq' y$, and $g(z)\supseteq G$, so $g(z)\sqsubseteq^{(\mathcal{F}^\under)_\ff}G$. Since $z\sqsubseteq ' y$, we have $g(z)\sqsubseteq^{(\mathcal{F}^\under)_\ff} g(y)$ by \SqForth{}, which with $g(z)\sqsubseteq^{(\mathcal{F}^\under)_\ff}G$ implies $g(y)\comp^{(\mathcal{F}^\under)_\ff} G$. This completes the proof of \SRBack{}.\end{proof}

We now obtain our possibility-semantic analogue of the Goldblatt-Thomason Theorem.

\begin{theorem}[Modal Definability of Full Possibility Frames]\label{ModDef2} $\,$
\begin{enumerate}
\item\label{ModDef2a} If a class $\mathsf{F}$ of full possibility frames is definable by modal formulas, then $\mathsf{F}$ is closed under dense possibility morphisms, selective subframes, and disjoint unions, while its complement is closed under filter extensions. 
\item\label{ModDef2b} If a class $\mathsf{F}$ of full possibility frames is closed under $\mathcal{L}^1(\ind)$-elementary equivalence, dense and strict possibility morphisms, generated subframes, and disjoint unions, while its complement is closed under filter extensions, then $\mathsf{F}$ is modally definable. 
\end{enumerate}
\end{theorem}

\begin{proof} For part \ref{ModDef2a}, all we need to add to the proof of Theorem \ref{ModDef1}.\ref{ModDef1a} is the fact, given in Corollary \ref{FilterPres}, that if $\varphi$ is valid over the filter extension of a frame, then $\varphi$ is valid over the frame. By the same reasoning as before, that fact and the assumption that $\mathsf{F}$ is modally definable together imply that the complement of $\mathsf{F}$ is closed under filter extensions. 
 
For part \ref{ModDef2b}, since $\mathsf{F}$ is closed under elementary equivalence and dense and strict possibility morphisms, it is closed under filter extensions by Lemma \ref{JvBLem}. Thus, we have a setup analogous to that of Theorem \ref{ModDef1}.\ref{ModDef1b}: \textsf{F} is closed under surjective p-morphisms, generated subframes, and disjoint unions, while both $\mathsf{F}$ and its complement are closed under filter extensions. Now reproduce the proof of Theorem \ref{ModDef1}.\ref{ModDef1b} but with $(\mathcal{F}^\under)_\ff$ in place of $(\mathcal{F}^\under)_\gff$, $((\underset{j\in J}{\biguplus}\mathcal{F}_j)^\under)_\ff$ in place of $((\underset{j\in J}{\biguplus}\mathcal{F}_j)^\under)_\gff$, Proposition \ref{Gen&Hom}.\ref{Gen&Hom2} in place of Proposition \ref{Gen&Hom}.\ref{Gen&Hom1}, and Proposition \ref{Poss&Sub}.\ref{Poss&Sub2} in place of Proposition \ref{Poss&Sub}.\ref{Poss&Sub1}. The resulting argument is a proof of part \ref{ModDef2b} of the current theorem.\end{proof}

\subsection{Modal Formulas with First-Order Correspondents}\label{LemmScottCorr}

In this section, we address question 1 from the beginning of \S~\ref{special}, going in the opposite direction relative to \ref{DefViaDual2}: given a modal formula $\varphi$, is the class of full possibility frames that $\varphi$ defines also definable by a first-order sentence of $\mathcal{L}^1(\ind)$? Let us phrase this in the standard terms of first-order \textit{correspondence}.

\begin{definition}[Relative Frame Correspondence]\label{RelCor} A formula $\varphi\in\mathcal{L}(\sig,\ind)$ \textit{globally corresponds} to a sentence $\psi\in\mathcal{L}^1(\ind)$ \textit{relative to} a class $\mathsf{F}$ of possibility frames iff for any  $\mathcal{F}\in\mathsf{F}$, $\varphi$ is valid over $\mathcal{F}$ as a possibility frame iff $\psi$ is true in $\mathcal{F}_\found$ (recall Definition \ref{Foundation}) as a structure for $\mathcal{L}^1(\ind)$. 

A formula $\varphi\in\mathcal{L}(\sig,\ind)$ \textit{locally corresponds} to a formula $\psi(\mathrm{x})\in\mathcal{L}^1(\ind)$ with exactly one free variable $\mathrm{x}$, relative to $\mathsf{F}$, iff for any $\mathcal{F}\in\mathsf{F}$ and $s\in\mathcal{F}$, $\mathcal{F},s\Vdash \varphi$ (recall Definition \ref{pmtruth1}) iff $\psi(\mathrm{x})$ is satisfied by $s$ in $\mathcal{F}_\found$. \hfill $\triangleleft$  
\end{definition} 

Note that $\varphi\in\mathcal{L}(\sig,\ind)$ globally corresponds to $\psi\in\mathcal{L}^1(\ind)$ relative to $\mathsf{F}$ iff $\varphi$ and $\psi$ define the same class of possibility frames relative to $\mathsf{F}$, as in Definitions \ref{RelDef} and \ref{FOdef}.

In \S~\ref{intro}, we quoted Goldblatt's \citeyearpar[p.~51]{Goldblatt2006} remark that a ``substantial reason for the great success'' of Kripke semantics is the way in which many natural modal axioms correspond to first-order properties of the accessibility relations in Kripke frames (full world frames), e.g., seriality ($\Box\varphi\rightarrow\Diamond\varphi$), reflexivity ($\Box\varphi\rightarrow\varphi$), transitivity ($\Box\varphi\rightarrow\Box\Box\varphi$), symmetry ($\varphi\rightarrow\Box\Diamond\varphi$), etc.  When moving to a more general semantics, we may lose such nice correspondences. For example, although for any full world frame $\mathfrak{F}$, $\mathfrak{F}$ validates the 4 axiom, $\Box_i\varphi\rightarrow \Box_i\Box_i\varphi$, iff $\wo{R}_i$ is transitive, it is \textit{not} the case that for any (general) world frame $\mathfrak{g}$ (see \S~\ref{GFS}), $\mathfrak{g}$ validates the 4 axiom iff $\wo{R}_i$ is transitive. The reason is that even if $\wo{R}_i$ is not transitive, $\mathfrak{g}$ may validate the 4 axiom because of the limitations on admissible valuations over $\mathfrak{g}$.\footnote{As Kracht \citeyearpar{Kracht1993} discusses, familiar correspondences may be restored relative to certain classes of general world frames defined by conditions on the set of admissible propositions.} Analogous points apply to other axioms. 

Given a modal formula $\varphi$, we will ask whether it has a first-order correspondent over full possibility frames. One can seek an answer in terms of semantic or syntactic features of $\varphi$. Over full \textit{world} frames, a semantic feature that is necessary and sufficient for $\varphi$ to have a global first-order correspondent is that the validity of $\varphi$ be preserved under taking \textit{ultrapowers} of full world frames, and a related feature characterizes local first-order correspondence (see \citealt[\S~VIII]{Benthem1976b,Benthem1983}). A syntactic feature that is sufficient for local and global correspondence over full world frames is that $\varphi$ be a \textit{Sahlqvist formula} (see \citealt[\S~3.6]{Blackburn2001}, \citealt{Sahlqvist1975}, \citealt[\S~I.4]{Benthem1976}, \citealt[\S~IX]{Benthem1983}). 
 
For full possibility frames, on the semantic side an analogous result with ultrapowers holds, as shown in \citealt{Yamamoto2016}. Given a full possibility frame $\mathcal{F}$, consider its foundation $\mathcal{F}_\found$ (Definition \ref{Foundation}) as a first-order structure for $\mathcal{L}^{1}(\ind)$, and then take an ultrapower $\prod_U \mathcal{F}_\found$ of $\mathcal{F}_\found$ in the standard sense. By Corollary \ref{DefFullFrames} and the preservation of first-order sentences under taking ultrapowers, $\prod_U \mathcal{F}_\found$ is the foundation of a full possibility frame. Thus, there is a full possibility frame $(\prod_U \mathcal{F}_\found)^\full$ based on  $\prod_U \mathcal{F}_\found$. We say that $(\prod_U \mathcal{F}_\found)^\full$ is an ultrapower of $\mathcal{F}$. Then the analogue of van Benthem's result is that $\varphi$ has a global first-order correspondent over full possibility frames iff the validity of $\varphi$ is preserved under taking ultrapowers of full possibility frames.

On the syntactic side, Yamamoto \citeyearpar{Yamamoto2016} shows that the possibility-semantic version of the Sahlqvist correspondence theorem also holds: every Sahlqvist formula has a local (resp.~global) first-order correspondent over full possibility frames. Yamamoto's proof uses the connection between full possibility frames and $\mathcal{CV}$-BAOs (\S\S~\ref{PossToBAO}-\ref{DualEquiv}) and the methods of algebraic modal correspondence \citep{Conradie2014}, further supporting the theme of the present paper of the importance of duality theory.

Below we will give some simpler methods for calculating first-order correspondents for restricted classes of Sahlqvist formulas. First, we discuss the extent to which the ``minimal valuation'' heuristic for first-order correspondence in Kripke semantics can be applied to possibility semantics. Second, we give an analogue for full possibility frames of one of the most elegant first-order correspondence results for full world frames, namely Lemmon and Scott's \citeyearpar[\S~4]{Lemmon1977} result (also covered in \citealt[\S~3.3, \S~5.5]{Chellas1980}, \citealt[\S~6]{Popkorn1994}, and \citealt[\S~9]{Garson2014}) for formulas of the form $\Diamond_\alpha\Box_\beta p\rightarrow \Box_\delta\Diamond_\gamma p$, where $\Diamond_\alpha$ is a sequence of diamond operators, $\Box_\beta$ is a sequence of box operators, etc., which covers many familiar modal axioms.

As in the case of possible world semantics, so too in the case of possibility semantics, the problem of establishing first-order correspondence can be viewed as the problem of showing that certain second-order formulas have first-order equivalents. Every modal formula has a \textit{second-order correspondent} over full frames. To get to this second-order perspective, we begin with a first-order language $\mathcal{L}^1(\sig,\ind)$ that extends $\mathcal{L}^1(\ind)$ with a predicate $\dot{Q}$ for each $q\in\sig$. While we interpreted $\mathcal{L}^1(\ind)$ in possibility \textit{frames}, we interpret  $\mathcal{L}^1(\sig,\ind)$ in possibility \textit{models} $\mathcal{M}=\langle S,\sqsubseteq,\{R_i\}_{i\in\ind},\pi\rangle$, interpreting $\dot{Q}$ as $\pi(q)$. We relate the modal  $\mathcal{L}(\sig,\ind)$ to the first-order  $\mathcal{L}^1(\sig,\ind)$ over possibility models using the following definition and proposition.

\begin{definition}[Standard Translation]\label{StanTran}
For each first-order variable $\var{x}$, define the \textit{standard translation} $ST_\var{x}\colon \mathcal{L}(\sig,\ind)\to \mathcal{L}^1(\sig,\ind)$ as follows:
\begin{enumerate}
\item $ST_\var{x}(q)=\dot{Q}\var{x}$ for $q\in\sig$;
\item $ST_\var{x}(\neg\varphi)=\forall \var{y}(\var{y}\,\dot{\sqsubseteq}\,\var{x}\rightarrow \neg ST_\var{y}(\varphi))$, where $\var{y}$ is a fresh variable;
\item $ST_\var{x}((\varphi\wedge\psi))=(ST_\var{x}(\varphi)\wedge ST_\var{x}(\psi))$;
\item $ST_\var{x}((\varphi\rightarrow\psi))=\forall \var{y}((\var{y}\,\dot{\sqsubseteq}\,\var{x}\wedge ST_\var{y}(\varphi))\rightarrow ST_\var{y}(\psi))$, where $\var{y}$ is a fresh variable;
\item $ST_\var{x}(\Box_i\varphi)=\forall \var{y}(\var{x}\dot{R_i}\var{y}\rightarrow ST_\var{y}(\varphi))$, where $\var{y}$ is a fresh variable.
\end{enumerate}
Let $ST(\varphi)=\forall \var{x} ST_\var{x}(\varphi)$.  \hfill $\triangleleft$
\end{definition}

\begin{proposition}[Correspondence on Models] For any possibility model $\mathcal{M}$, $s\in\mathcal{M}$, and $\varphi\in\mathcal{L}(\sig,\ind)$:
\begin{enumerate}
\item $\mathcal{M},s\Vdash \varphi$ iff the formula $ST_\mathrm{x}(\varphi)$ is satisfied by $s$ in $\mathcal{M}$ regarded as a first-order structure for $\mathcal{L}^1(\sig,\ind)$;
\item $\mathcal{M}\Vdash \varphi$ iff the sentence $ST(\varphi)$ is true in $\mathcal{M}$ regarded as a first-order structure for $\mathcal{L}^1(\sig,\ind)$.
\end{enumerate}
\end{proposition}

 To have this kind of result at the level of \textit{frame validity}, we move to the monadic \textit{second-order} language $\mathcal{L}^2(\sig,\ind)$ obtained from $\mathcal{L}^1(\ind)$ by adding a predicate \textit{variable} $\dot{Q}$ for each $q\in\sig$. We interpret $\mathcal{L}^2(\sig,\ind)$ in possibility frames---or rather, foundations $\mathcal{F}_\found=\langle S,\sqsubseteq,\{R_i\}_{i\in\ind}\rangle$ of possibility frames---with $\wp(S)$ as the domain for the monadic second-order quantifiers. We then relate the modal language $\mathcal{L}(\sig,\ind)$ to the second-order language $\mathcal{L}^2(\sig,\ind)$ over full possibility frames using the following definition and proposition.

\begin{definition}[Second-Order Translation] For each first-order variable $\mathrm{x}$, define the \textit{second-order translation} $SOT_\mathrm{x}\colon \mathcal{L}(\sig,\ind)\to \mathcal{L}^2(\sig,\ind)$ as follows, where $\dot{Q}$ is a second-order variable:
\begin{enumerate}
\item $RO(\dot{Q})= \forall \var{v} (\dot{Q}\var{v}\leftrightarrow\forall \var{v}' (\var{v}'\,\dot{\sqsubseteq}\,\var{v}\rightarrow \exists \var{v}'' (\var{v}''\,\dot{\sqsubseteq}\,\var{v}'\wedge \dot{Q}\var{v}'')))$;
\item $SOT_{\mathrm{x}} (\varphi)= \forall \dot{Q_1}\dots \forall \dot{Q_n} ((\underset{1\leq i\leq n}{\bigwedge}RO(\dot{Q_i}))\rightarrow ST_{\mathrm{x}}(\varphi))$, 
\end{enumerate}
where $\dot{Q_1},\dots, \dot{Q_n}$ are the unary predicates occurring in $ST_{\mathrm{x}}(\varphi)$. 

Let $SOT(\varphi)=\forall \mathrm{x}SOT_{\mathrm{x}}(\varphi)$.\hfill$\triangleleft$
\end{definition}

From now on we will drop the dots over symbols when we trust there is no danger of confusion, and we will use abbreviations of the form `$\forall \mathrm{y}\sqsubseteq\mathrm{x}$' for restricted quantification.

\begin{proposition}[Correspondence on Frames]\label{CorFrameSOT} For any \textit{full} possibility frame $\mathcal{F}$, $s\in S$, and $\varphi\in\mathcal{L}(\sig,\ind)$:
\begin{enumerate}
\item $\mathcal{F},s\Vdash \varphi$ iff the formula $SOT_\mathrm{x}(\varphi)$ is satisfied by $s$ in $\mathcal{F}_\found$;
\item $\mathcal{F}\Vdash\varphi$ iff the sentence $SOT(\varphi)$ is true in $\mathcal{F}_\found$.
\end{enumerate}
\end{proposition}

\noindent Thus, a modal formula $\varphi$ having a first-order local/global correspondent $\psi$ over full possibility frames is the same as the second-order $SOT_\mathrm{x}(\varphi)$/$SOT(\varphi)$ being logically equivalent to the first-order $\psi$.

The first-order correspondence problem for full possibility frames is related to the first-order correspondence problem for full \textit{world} frames by the following fact, which is an immediate consequence of the fact that full world frames are a special case of full possibility frames.

\begin{fact}[Possibility Correspondence Implies World Correspondence]\label{PossToWorldCor}  Any modal formula that has a local/global first-order correspondent over all full possibility frames also has a local/global first-order correspondent over full world frames, viz., the same first-order formula, in which $\dot{\sqsubseteq}$ may be replaced by $=$.
\end{fact}
\noindent Thus, the modal formulas that lack first-order correspondents over full world frames, e.g., $\Box_i\Diamond_i p\rightarrow\Diamond_i \Box_i p$ \citep{Benthem1975,Goldblatt1975b}, also lack first-order correspondents over full possibility frames. This does not show, however, that such formulas lack a first-order correspondent over some restricted class of full possibility frames that does not include all full world frames, e.g., \textit{functional} full frames (\S~\ref{FuncFrames}). 

The converse question of whether every modal formula with a local/global first-order correspondent over full world frames also has one over full possibility frames is an open question (see \S~\ref{OpenProb}).

We will consider both local and global correspondence as in Definition \ref{RelCor}. By part \ref{LocalVsGlobal1} of the following obvious fact, having a local correspondent implies having a global one (part \ref{LocalVsGlobal2} will be useful later).
 
\begin{fact}[From Local to Global]\label{LocalVsGlobal} For any $\varphi,\psi\in\mathcal{L}(\sig,\ind)$, $\delta(\mathrm{x})\in\mathcal{L}^1(\ind)$ with $\mathrm{x}$ free, and class $\mathsf{F}$ of possibility frames:
\begin{enumerate}
\item\label{LocalVsGlobal1} if $\varphi$ locally corresponds to $\delta(\mathrm{x})$ relative to $\mathsf{F}$, then $\varphi$ globally corresponds to $\forall \mathrm{x}\,\delta(\mathrm{x})$ relative to $\mathsf{F}$;
\item\label{LocalVsGlobal2} if $\varphi\rightarrow\psi$ locally corresponds to $\forall \mathrm{x} \sqsubseteq \mathrm{x}_0 \,\delta(\mathrm{x})$ relative to $\mathsf{F}$, then $\varphi\rightarrow\psi$ globally corresponds to $\forall \mathrm{x}\,\delta(\mathrm{x})$ relative to $\mathsf{F}$.
\end{enumerate}
\end{fact}
As in possible world semantics, so too in possibility semantics, having a global correspondent does not imply having a local correspondent. The following example comes from \citealt[Theorem 7.1]{Benthem1983}.

\begin{fact}[Global without Local]\label{GlobalWithoutLocal} $\Box_i\Diamond_i\Box_i\Box_i p\rightarrow \Diamond_i\Diamond_i\Box_i\Diamond_i p$ globally corresponds to $\forall \mathrm{x}\,\exists\mathrm{y}\, \mathrm{x}R_i\mathrm{y}$ over all possibility frames, but does not locally correspond to any $\mathcal{L}^1(\ind)$ formula even over full possibility frames.
\end{fact}

\begin{proof} First, for the claim of global correspondence, observe that for any possibility model $\mathcal{M}$ and $x\in\mathcal{M}$, if there is no $y$ such that $xR_iy$, then every $\Box_i\varphi$ is true at $x$ and every $\Diamond_i\psi$ is false at $x$. Thus, if a possibility frame $\mathcal{F}$ does not satisfy $\forall \mathrm{x}\,\exists\mathrm{y}\, \mathrm{x}R_i\mathrm{y}$, so there is an $x$ with no $y$ such that $xR_i y$, then $\Box_i\Diamond_i\Box_i\Box_i p\rightarrow \Diamond_i\Diamond_i\Box_i\Diamond_i p$ is false at $x$ in any model based on $\mathcal{F}$. In the other direction, observe that if $\mathcal{F}$  satisfies $\forall \mathrm{x}\,\exists\mathrm{y}\, \mathrm{x}R_i\mathrm{y}$, then for any $\mathcal{M}$ based on $\mathcal{F}$ and $x\in\mathcal{M}$, $\mathcal{M},x\Vdash \Box_i\varphi$ implies $\mathcal{M},x\Vdash\Diamond_i\varphi$. For if $\mathcal{M},x\Vdash \Box_i\varphi$, then by Persistence, for all $x'\sqsubseteq x$ we have $\mathcal{M},x'\Vdash\Box_i\varphi$, which with $\forall \mathrm{x}\,\exists\mathrm{y}\, \mathrm{x}R_i\mathrm{y}$ implies that there is a $y'$ such that $x'R_iy'$ and $\mathcal{M},y'\Vdash\varphi$, which implies $\mathcal{M},x\Vdash \Diamond_i\varphi$.  Then since $\mathcal{M},x\Vdash\Box_i\varphi$ implies $\mathcal{M},x\Vdash\Diamond_i\varphi$, it clearly follows that $\mathcal{M},x\Vdash\Box_i\Diamond_i\Box_i\Box_i p$ implies $\mathcal{M},x\Vdash \Diamond_i\Diamond_i\Box_i\Diamond_i p$. Thus, $\Box_i\Diamond_i\Box_i\Box_i p\rightarrow \Diamond_i\Diamond_i\Box_i\Diamond_i p$ is valid over $\mathcal{F}$.

The claim of no local correspondent follows from the fact that $\Box_i\Diamond_i\Box_i\Box_i p\rightarrow \Diamond_i\Diamond_i\Box_i\Diamond_i p$ has no local correspondent over full world frames (\citealt[Theorem 7.1]{Benthem1983}) plus Fact \ref{PossToWorldCor}.
\end{proof}

The proof of Fact \ref{GlobalWithoutLocal} also shows that $\Box_i p\rightarrow\Diamond_i p$ globally corresponds to $\forall \mathrm{x}\,\exists \mathrm{y}\, \mathrm{x}R_i\mathrm{y}$. In contrast to Fact \ref{GlobalWithoutLocal}, however, it is easy to see that $\Box_i p\rightarrow\Diamond_i p$ locally corresponds to $\forall\mathrm{x}\sqsubseteq\mathrm{x}_0\,\exists \mathrm{y}\,\mathrm{x}R_i\mathrm{y}$. This is an atypical case of correspondence, insofar as we can show that if $\mathcal{F},x_0$ does not satisfy $\forall\mathrm{x}\sqsubseteq\mathrm{x}_0\,\exists \mathrm{y}\,\mathrm{x}R_i\mathrm{y}$, then for \textit{every} model $\mathcal{M}$ based on $\mathcal{F}$, $\mathcal{M},x_0\nVdash \Box_i p\rightarrow\Diamond_i p$. Typically we have to carefully choose a falsifying model.

Indeed, the tricky part of establishing correspondence is usually that of showing that if $\mathcal{F},x\Vdash\varphi$, then $\mathcal{F},x$ satisfies the putative local first-order correspondent $\psi(\mathrm{x})$ of $\varphi$, and similarly in the global case. The \textit{direct strategy} for local correspondence is to show that we can find \textit{minimal valuations} $\pi$ in such a way that from the fact that $\varphi$ is true at $x$ in the models $\langle\mathcal{F},\pi\rangle$ (by the assumption that $\mathcal{F},x\Vdash\varphi$), it follows that $\psi(\mathrm{x})$ is satisfied by $x$ in $\mathcal{F}$ as a first-order structure. (See \citealt[\S~9.4]{Benthem2010} for an introduction to this strategy as applied to full world frames.) The \textit{contrapositive strategy} for local or global correspondence is to show that if $\mathcal{F},x$ (resp.~$\mathcal{F}$) does not satisfy the putative first-order correspondent, then we can add a valuation to obtain a model on $\mathcal{F}$ witnessing $\mathcal{F},x\nVdash\varphi$ (resp.~$\mathcal{F}\nVdash\varphi$). For both strategies, the trick is to pick the right valuations. In the context of possibility semantics, we have the additional constraint that an admissible valuation $\pi$ on a frame must be such that $\pi(p)$ satisfies \textit{persistence} and \textit{refinability}. In this respect, correspondence theory for possibility semantics is similar to intuitionistic correspondence theory, as developed in \citealt{Rodenburg1986}, which shares the \textit{persistence} constraint on valuations.  

Thus, we need methods of constructing admissible valuations on full possibility frames. From Fact \ref{CofGenerated}, we know that if for each $p\in\sig$, we set $\pi(p)=\{x'\in S\mid x'\cof x\}$ for some $x\in S$, then $\pi$ is admissible. More generally, from Fact \ref{RefReg}.\ref{RefReg2.5}, we know that if for each $p\in\sig$ we set 
\[\pi(p)=\mathrm{int}(\mathrm{cl}(\mathord{\Downarrow}X))=\{x\in S\mid \forall x'\sqsubseteq x\,\exists x''\sqsubseteq x'\colon x''\in \mathord{\Downarrow}X\}\] 
for some $X\subseteq S$, then $\pi$ is admissible. The following gives another way of ensuring that $\pi$ is admissible.

\begin{fact}[$\incomp$-Generated Propositions]\label{ConVal} For any poset $\langle S,\sqsubseteq\rangle$ and $Y\subseteq S$, $\{x\in S\mid \forall y\in Y,\, x\incomp y\}$ is in $\mathrm{RO}(S,\sqsubseteq)$.
\end{fact}

\begin{proof} Since $x'\sqsubseteq x$ and $x\incomp y$ together imply $x'\incomp y$, $\{x\in S\mid \forall y\in Y,\, x\incomp y\}$ satisfies \textit{persistence}. 

For \textit{refinability}, suppose $x\not\in\{x\in S\mid \forall y\in Y,\, x\incomp y\}$, so there is a $y\in Y$ such that $x\comp y$. Then there is an $x'\sqsubseteq x$ such that $x'\sqsubseteq y$, which implies that for all $x''\sqsubseteq x'$, $x''\comp  y$ and hence $x''\not\in \{x\in S\mid \forall y\in Y,\, x\incomp y\}$. Thus, $\{x\in S\mid \forall y\in Y,\, x\incomp y\}$ satisfies \textit{refinability}.
\end{proof} 

In the rest of this section, we will state our results for correspondence relative to the class of full \textit{standard} possibility frames (Definition \ref{Standard}), i.e., full possibility frames satisfying the \Rdown{} condition that if $xR_iy$ and $y'\sqsubseteq y$, then $xR_iy'$ (recall Example \ref{IntMod} and \S~\ref{FullFrames}). The \Rdown{} condition significantly simplifies the first-order correspondents of modal formulas (recall Fact \ref{ItMod}). The first-order correspondents are further simplified by assuming the conditions of \textit{strong} possibility frames, especially \Rdense{}: if $\forall y'\sqsubseteq y$ $\exists y''\sqsubseteq y'$ $xR_iy''$, then $xR_iy$. Recall from Fact \ref{R(x)RO} that \Rdown{} and \Rdense{} are jointly equivalent to the condition that for each $x\in\mathcal{F}$, $R_i(x)\in\mathrm{RO}(\mathcal{F})$. Also recall from Proposition \ref{Representation} that any full possibility frame can be transformed into a modally equivalent full, strong and hence standard possibility frame. Our correspondence results for standard full frames will imply correspondence results for all full frames, by the following reasoning.

\begin{lemma}[Transferring Correspondence]\label{TransferCorr} Given $\psi\in\mathcal{L}^1(\ind)$, let $\psi_{\mathord{\downarrow}}$ be the result of replacing each subformula of $\psi$ of the form $\mathrm{x}\dot{R_i} \mathrm{y}$ with $\exists \mathrm{y}' (\mathrm{x}\dot{R_i}\mathrm{y}' \wedge \mathrm{y}\sqsubseteq \mathrm{y}')$ for a fresh $\mathrm{y}'$. If $\psi$ is a local (resp. global) first-order correspondent of $\varphi$ over \textit{standard} full possibility frames, then $\psi_{\mathord{\downarrow}}$ is a local (resp. global) first-order correspondent of $\varphi$ over \textit{all} full possibility frames.
\end{lemma}
\begin{proof} Given a full possibility frame $\mathcal{F}=\langle S,\sqsubseteq, \{R_i\}_{i\in\ind},\adm\rangle$, define $\mathcal{F}_{\mathord{\downarrow}}=\langle S,\sqsubseteq, \{R_{i\mathord{\downarrow}}\}_{i\in\ind}, \adm\rangle$ by $xR_{i\mathord{\downarrow}}y$ iff $\exists y'$: $xR_iy'$ and $y\sqsubseteq y'$. Then it is easy to check that $\mathcal{F}_{\mathord{\downarrow}}$ is a full possibility frame, and $\mathcal{F}_{\mathord{\downarrow}}$ satisfies \Rdown{} by construction, so it is a standard possibility frame. Moreover, the identity map on $S$ is a surjective robust possibility morphism from $\mathcal{F}$ to $\mathcal{F}_{\mathord{\downarrow}}$, so these frames validate the same formulas at each state by the proof of Proposition \ref{Preservation}. Since $\mathcal{F}$ satisfies $\exists \mathrm{y}' (\mathrm{x}\dot{R_i}\mathrm{y}' \wedge \mathrm{y}\sqsubseteq \mathrm{y}')$ with some variable assignment iff $\mathcal{F}_{\mathord{\downarrow}}$ satisfies $\mathrm{x}\dot{R_i} \mathrm{y}$ with the same variable assignment, $\mathcal{F}$ satisfies $\psi_{\mathord{\downarrow}}$ iff $\mathcal{F}_{\mathord{\downarrow}}$ satisfies $\psi$. Then assuming that  $\mathcal{F}_{\mathord{\downarrow}},x$ (resp.~$\mathcal{F}_{\mathord{\downarrow}}$) satisfies $\psi$ iff $\mathcal{F}_{\mathord{\downarrow}},x$ (resp.~$\mathcal{F}_{\mathord{\downarrow}}$) validates $\varphi$, it follows that $\mathcal{F},x$ (resp.~$\mathcal{F}$) satisfies $\psi_{\mathord{\downarrow}}$ iff $\mathcal{F},x$ (resp.~$\mathcal{F}$) validates~$\varphi$.\end{proof}

By the same kind of argument, we can show that if $\psi$ is a first-order correspondent of $\varphi$ over \textit{strong} full possibility frames, then a related formula $\psi'$ is a first-order correspondent of $\varphi$ over all full possibility frames (using the relation $R_i^*$ defined by $xR_i^*y$ iff $y\in\mathrm{int}(\mathrm{cl}(\mathord{\Downarrow}R_i(x)))$, which for full frames is $R_i^\tight$ from Proposition~\ref{Representation}). Note, however, that none of our arguments show, e.g., that if $\varphi$ has a first-order correspondent over \textit{functional} full possibility frames (\S~\ref{FuncFrames}), then it has a first-order correspondent over all full frames. For we cannot always turn a full possibility frame into a modally equivalent functional one simply by modifying the accessibility relations of the frame, let alone in a first-order definable way. It is an open question whether more modal formulas have first-order correspondents over functional full possibility frames.

Let us first consider the direct strategy for establishing correspondence with minimal valuations. As an example, we compute a first-order correspondent for $\Box_i q\rightarrow\Box_i\Box_i q$ (again~cf. \citealt[\S~9.4]{Benthem2010}).
 
\begin{example}[Correspondence by Minimal Valuations]\label{MinVal1}
Given a frame $\mathcal{F}$, what is the minimal way of making the antecedent of $\Box_i q\rightarrow \Box_i\Box_i q$ true at a state $x$? If $\mathcal{F}$ were a full world frame,  the minimal way of making $\Box_i q$ true at $x$ would be with a valuation $\pi$ such that $\pi(q)=R_i(x)$. However, if $\mathcal{F}$ is a standard full possibility frame, the minimal \textit{admissible} way of making $\Box_iq$ true at $x$ is with a valuation $\pi$ such that 
\begin{equation}\pi (q)=\mathrm{int}(\mathrm{cl}(R_i(x)))=\{z\in S\mid \forall w\sqsubseteq z\,\exists w'\sqsubseteq w\colon xR_iw'\}.\label{MinVal}\end{equation}
Since $R_i(x)$ is a downset by \Rdown{}, this $\pi(q)$ is the minimal regular open set including $R_i(x)$ by Fact~\ref{RefReg}.\ref{RefReg2.5}. Now consider the second-order translation with respect to $\mathrm{x}_0$ of $\Box_i q\rightarrow\Box_i\Box_iq$,
\[\forall Q(RO(Q)\rightarrow ST_{\mathrm{x}_0}(\Box_iq\rightarrow \Box_i\Box_iq)),\]
which is equivalent to
\begin{eqnarray}
&&\forall Q\Big(\forall \mathrm{v}(Q\wo{v}\leftrightarrow \forall \wo{v}'\sqsubseteq \wo{v}\,\exists\wo{v}''\sqsubseteq\wo{v}'\; Q\wo{v}'')\rightarrow \nonumber\\
&& \qquad\forall \mathrm{x}\sqsubseteq \mathrm{x}_0\big(\forall \mathrm{y} (\mathrm{x}R_i\mathrm{y}\rightarrow Q\mathrm{y})\rightarrow \forall \mathrm{z} (\mathrm{x}R_i^2\mathrm{z}\rightarrow Q\mathrm{z})\big)\Big),\label{SOTran}\end{eqnarray}
using $\mathrm{x}R_i^2\mathrm{z}$ as the obvious abbreviation. If we move $\forall \mathrm{x}\sqsubseteq\mathrm{x}_0$ in (\ref{SOTran}) before $\forall Q$ and then, for each $
  \mathrm{x}$, plug in for $Q$ the first-order description of the minimal valuation from (\ref{MinVal}), i.e., for each variable $\alpha$, replace $Q\alpha$ with $\forall \mathrm{w}\sqsubseteq\alpha\,\exists \mathrm{w}'\sqsubseteq \mathrm{w}\; \mathrm{x}R_i\mathrm{w}'$, we obtain:
\begin{eqnarray} 
&&\forall \mathrm{x}\sqsubseteq\mathrm{x}_0 \Big( \forall \mathrm{v}\big((\forall \mathrm{w}\sqsubseteq \mathrm{v}\,\exists \mathrm{w}'\sqsubseteq \mathrm{w}\; \mathrm{x}R_i\mathrm{w}') \leftrightarrow (\forall \mathrm{v}'\sqsubseteq \mathrm{v}\,\exists\mathrm{v}''\sqsubseteq\mathrm{v}'\, \forall \mathrm{w}\sqsubseteq \mathrm{v}''\,\exists \mathrm{w}'\sqsubseteq\mathrm{w} \;\mathrm{x}R_i \mathrm{w}')\big)\rightarrow \nonumber\\
&& \qquad\qquad\big(\forall \mathrm{y} (\mathrm{x}R_i\mathrm{y}\rightarrow \forall \mathrm{w}\sqsubseteq \mathrm{y}\,\exists \mathrm{w}'\sqsubseteq \mathrm{w}\, \mathrm{x}R_i\mathrm{w}')\rightarrow \forall \mathrm{z} (\mathrm{x}R_i^2\mathrm{z}\rightarrow \forall \mathrm{w}\sqsubseteq \mathrm{z}\,\exists \mathrm{w}'\sqsubseteq \mathrm{w}\; \mathrm{x}R_i\mathrm{w}')\big)\Big).\label{instantiate}\end{eqnarray}
Every frame satisfies the main antecedent of (\ref{instantiate}), so we can reduce (\ref{instantiate}) to 
\begin{eqnarray} 
&&\forall \mathrm{x}\sqsubseteq\mathrm{x}_0\big(\forall \mathrm{y} (\mathrm{x}R_i\mathrm{y}\rightarrow \forall \mathrm{w}\sqsubseteq \mathrm{y}\,\exists \mathrm{w}'\sqsubseteq \mathrm{w}\, \mathrm{x}R_i\mathrm{w}')\rightarrow \forall \mathrm{z} (\mathrm{x}R_i^2\mathrm{z}\rightarrow \forall \mathrm{w}\sqsubseteq \mathrm{z}\,\exists \mathrm{w}'\sqsubseteq \mathrm{w}\, \mathrm{x}R_i\mathrm{w}')\big).\label{instantiate2}\end{eqnarray}
Moreover, every frame satisfying \Rdown{} satisfies the main antecedent in (\ref{instantiate2}), so we can reduce (\ref{instantiate2}) to
\begin{equation}\forall\mathrm{x}\sqsubseteq\mathrm{x}_0\,\forall \mathrm{z} (\mathrm{x}R_i^2\mathrm{z}\rightarrow \forall \mathrm{w}\sqsubseteq \mathrm{z}\,\exists \mathrm{w}'\sqsubseteq \mathrm{w}\, \mathrm{x}R_i\mathrm{w}'),\label{PreCor}\end{equation}
which over frames satisfying \Rdown{} is equivalent to the simpler
\begin{equation}\forall\mathrm{x}\sqsubseteq\mathrm{x}_0\,\forall \mathrm{w} (\mathrm{x}R_i^2\mathrm{w}\rightarrow \exists \mathrm{w}'\sqsubseteq \mathrm{w}\, \mathrm{x}R_i\mathrm{w}').\label{LastCor}\end{equation}
Over frames satisfying both \Rdown{} and \Rdense{}, (\ref{PreCor}) is equivalent to the still simpler
\begin{equation}\forall\mathrm{x}\sqsubseteq\mathrm{x}_0\,\forall \mathrm{w} (\mathrm{x}R_i^2\mathrm{w}\rightarrow \mathrm{x}R_i\mathrm{w}).\label{LastCor2}\end{equation}
Now we claim that (\ref{SOTran}) and (\ref{PreCor})/(\ref{LastCor}) are equivalent over standard frames. We already have the implication from (\ref{SOTran}) to (\ref{PreCor}), since (\ref{PreCor}) is equivalent to (\ref{instantiate}), which is simply an \textit{instantiation} of (\ref{SOTran}). The implication from (\ref{PreCor}) to (\ref{SOTran}) relies on the following two key points:
\begin{itemize}
\item[(A)] If $\alpha:=\forall \mathrm{v}(Q\wo{v}\leftrightarrow \forall \wo{v}'\sqsubseteq \wo{v}\,\exists\wo{v}''\sqsubseteq\wo{v}'\; Q\wo{v}'')\wedge \mathrm{x}\sqsubseteq\mathrm{x}_0\wedge \forall \mathrm{y} (\mathrm{x}R_i\mathrm{y}\rightarrow Q\mathrm{y})$ from (\ref{SOTran}) is satisfied in $\mathcal{F}$ with a variable assignment $\nu$ such that $\nu(\mathrm{x})=x$, then $\nu(Q)$ is a regular open set that includes $R_i(x)$. Thus, $\nu(Q)$ is a superset of our $\mathrm{int}(\mathrm{cl}(R_i(x)))$ from (\ref{MinVal}), which we noted above is the \textit{minimal} regular open set that includes $R_i(x)$.
\item[(B)] Thanks to its \textit{positive} syntactic form,  $\beta:=\forall \mathrm{z} (\mathrm{x}R_i^2\mathrm{z}\rightarrow Q\mathrm{z})$ from (\ref{SOTran}) is semantically \textit{upward monotone} in $Q$: for any assignments $\nu$ and $\nu'$ that agree on $\mathrm{x}$, if $\beta$ is satisfied by $\mathcal{F}$ with $\nu'$, and $\nu'(Q)\subseteq\nu(Q)$, then $\beta$ is  satisfied by $\mathcal{F}$ with $\nu$. 
\end{itemize}
Let $\gamma$ be the formula that comes after $\forall \mathrm{x}\sqsubseteq\mathrm{x}_0$ in (\ref{PreCor}), and suppose $\gamma$ is satisfied in $\mathcal{F}$ with some assignment $\nu$. It follows that $\beta$ is satisfied in $\mathcal{F}$ with the assignment $\nu'$ that differs from $\nu$ only in that $\nu'(Q)=\mathrm{int}(\mathrm{cl}(R_i(\nu(\mathrm{x}))))=\{z\in S\mid \forall w\sqsubseteq z\,\exists w'\sqsubseteq w\colon \nu(\mathrm{x})R_iw'\}$. Finally, suppose that $\alpha$ is satisfied in $\mathcal{F}$ with $\nu$. Then by (A), $\nu'(Q)\subseteq\nu(Q)$, so by (B) and the previous point that $\beta$ is satisfied in $\mathcal{F}$ with $\nu'$, we have that $\beta$ is satisfied by $\mathcal{F}$ with $\nu$. This shows that (\ref{PreCor}) implies (\ref{SOTran}). 

Thus, (\ref{LastCor}) is a local first-order correspondent of $\Box_i q\rightarrow \Box_i\Box_i q$ over standard frames, and (\ref{LastCor2}) is a local first-order correspondent over strong frames. Then by Fact \ref{LocalVsGlobal}, its global first-order correspondents are
\[\forall\mathrm{x}\,\forall \mathrm{w} (\mathrm{x}R_i^2\mathrm{w}\rightarrow \exists \mathrm{w}'\sqsubseteq \mathrm{w}\, \mathrm{x}R_i\mathrm{w}')\mbox{ and } \forall\mathrm{x}\,\forall \mathrm{w} (\mathrm{x}R_i^2\mathrm{w}\rightarrow \mathrm{x}R_i\mathrm{w})\]
over standard and strong possibility frames, respectively. Thus,  $\Box_i q\rightarrow \Box_i\Box_i q$ corresponds to a kind of ``delayed transitivity'' over standard frames and to ordinary transitivity over strong frames.
  \hfill $\triangleleft$ \end{example}

Example \ref{MinVal1} can be turned into a general result, which we state as Proposition \ref{SafeSahl} below. First, let us state the limitation of the minimal valuation method in the context of possibility semantics: it does not work when a diamond is in the antecedent of an implication, as in, e.g., the 5 axiom $\Diamond_i q\rightarrow \Box_i\Diamond_i q$. Recall from \S~\ref{Acc&Poss} that the derived clause for $\Diamond_i\varphi:=\neg\Box_i\neg\varphi$ is (in its simplified form using \Rdown{}):
\begin{itemize}
\item $\mathcal{M},x\Vdash \Diamond_i \varphi$ iff $\forall x'\sqsubseteq x$ $\exists y'$: $x'R_i y'$ and $\mathcal{M},y'\Vdash \varphi$.
\end{itemize}
This $\forall\exists$ pattern in an antecedent blocks the standard minimal valuation method, just as the $\Box_i\Diamond_i$ pattern in the antecedent of $\Box_i\Diamond_i q\rightarrow \Diamond_i\Box_i q$ blocks the method in Kripke semantics.
 
However, this obstacle involving diamonds in antecedents can be overcome. For simplicity, suppose we have a modal formula $\varphi$ containing only the propositional variable $q$. If the minimal valuation method would work on $\varphi$, it would give us a formula $\sigma\in\mathcal{L}^1(\ind)$ with a free variable $\mathrm{u}$ such that if we replace, for each variable $\alpha$, each occurrence of $Q\alpha$ in $ST_\mathrm{x}(\varphi)$  by $\sigma[\alpha/\mathrm{u}]$, then we obtain a formula $ST_\mathrm{x}(\varphi)[\sigma/Q]\in\mathcal{L}^1(\ind)$ with exactly $\mathrm{x}$ free that is equivalent to the second-order translation $\forall Q(RO(Q)\rightarrow ST_\mathrm{x}(\varphi))$ of $\varphi$. But this is more than we need for local first-order correspondence. It suffices to find a formula $\sigma\in\mathcal{L}^1(\ind)$ with free variables including $\mathrm{u}$ and $\mathrm{u}_1,\dots,\mathrm{u}_n$ such that the formula $\forall\mathrm{u}_1\dots\forall\mathrm{u}_n ST_\mathrm{x}(\varphi)[\sigma/Q]\in\mathcal{L}^1(\ind)$ has exactly $\mathrm{x}$ free and is equivalent to the second-order translation $\forall Q(RO(Q)\rightarrow ST_\mathrm{x}(\varphi))$ of $\varphi$. Thus, we are using $\sigma[\alpha/\mathrm{u}]$ to pick out a \textit{family} of subsets by varying parameters. This is an example of the more general method of substitutions for first-order correspondence (see \citealt[pp.~108-9]{Benthem1983}) that Yamamoto \citeyearpar{Yamamoto2016} uses to prove that every Sahlqvist modal formula has a first-order correspondent over full possibility frames.

Before leaving the minimal valuation method, let us generalize Example~\ref{MinVal1}. For this we need a quick review of relevant notions (see \citealt[pp.~151-3]{Blackburn2001}). Recall that for a propositional variable $q\in\sig$, a modal formula $\varphi\in\mathcal{L}(\sig,\ind)$ is \textit{positive} (resp.~\textit{negative}) \textit{in $q$} iff every occurrence of $q$ in $\varphi$ is under the scope of an even (resp.~odd) number of negations (where we count $\varphi\rightarrow\psi$ as $\neg\varphi\vee\psi$). Then $\varphi$ is \textit{positive} (resp.~\textit{negative}) iff it is positive (resp.~negative) in every propositional variable that occurs in it. The same notions apply to a predicate variable $\dot{Q}$ and a second-order formula $\varphi\in\mathcal{L}^2(\sig,\ind)$ in the analogous ways. If $\varphi$ is positive (resp.~negative), then obviously $\neg\varphi$ is negative (resp.~positive). Also note that $\varphi\in\mathcal{L}(\sig,\ind)$ is positive (resp.~negative) iff $ST_\mathrm{x}(\varphi)$ is positive (resp.~negative). The key property of positive (resp.~negative) formulas is that they are semantically \textit{upward} (resp.~\textit{downward}) \textit{monotone} in their variables: if a positive (resp.~negative) formula is true in a structure given an interpretation of the variables, then it remains true under expanding (resp.~shrinking) the interpretation of the variables. 

We now define the class of Sahlqvist formulas without diamonds in the antecedents of implications or disjunctions outside of antecedents, which we call \textit{safe} Sahlqvist formulas (cf.~\citealt[Def.~3.51]{Blackburn2001}).

\begin{definition}[Safe Sahlqvist Formulas] Fix a set $\mathcal{L}(\sig,\ind)$ of modal formulas.  

A \textit{boxed atom} is a formula of the form $\Box_{i_1}\dots\Box_{i_n}p$ for $i_1,\dots,i_n\in \ind$, $p\in\sig$. If $n=0$, $\Box_{i_1}\dots\Box_{i_n}p$ is $p$.

A \textit{safe Sahlqvist antecedent} is a formula built from boxed atoms and negative formulas using $\wedge$ and $\vee$.

A \textit{safe Sahlqvist implication} is a formula $\varphi\rightarrow\psi$ where $\varphi$ is a safe Sahlqvist antecedent and $\psi$ is positive. 

A \textit{safe Sahlqvist formula} is a formula built from safe Sahlqvist implications by applying boxes and conjunctions.\hfill $\triangleleft$
\end{definition}

The correspondence problem for safe Sahlqvist formulas can be reduced to the problem for safe Sahlqvist implications using the following standard lemma \citep[Lem.~3.53]{Blackburn2001} that continues to hold in possibility semantics. 

\begin{lemma}[Inductive Local Correspondence]\label{IndCor} For any $\varphi,\psi\in\mathcal{L}(\sig,\ind)$ and $\delta,\gamma\in\mathcal{L}^1(\ind)$:
\begin{enumerate}
\item if $\varphi$ locally corresponds to $\delta(\mathrm{x})$, then $\Box_i\varphi$ locally corresponds to $\forall \mathrm{y} (\mathrm{x}R_i\mathrm{y}\rightarrow \delta[\mathrm{y}/\mathrm{x}])$;
\item\label{IndCor2} if $\varphi$ locally corresponds to $\delta$, and $\psi$ locally corresponds to $\gamma$, then $\varphi\wedge\psi$ locally corresponds to $\delta\wedge\gamma$.
\end{enumerate}
\end{lemma}

Now we have the following significant generalization of Example \ref{MinVal1}.

\begin{proposition}[Safe Sahlqvist Correspondence]\label{SafeSahl} Every safe Sahlqvist formula $\varphi$ locally corresponds over full standard possibility frames to a $\psi(\mathrm{x})\in\mathcal{L}^1(\ind)$ that can be effectively computed from $\varphi$.
\end{proposition}

\begin{proof} Given Lemma \ref{IndCor}, we need only explain what to do with each safe Sahlqvist implication $\psi\rightarrow\chi$. First, using propositional logic, rewrite $\psi$ as an equivalent disjunction  $\psi_1\vee\dots \vee\psi_n$ where each $\psi_k$ is a conjunction of boxed atoms and negative formulas. Then rewrite $\psi\rightarrow\chi$ as the equivalent conjunction $\underset{1\leq k\leq n}\bigwedge (\psi_k\rightarrow\chi)$, which is a safe Sahlqvist formula. Then by Lemma \ref{IndCor}.\ref{IndCor2}, it suffices to find a local correspondent for each implication $\psi_k\rightarrow\chi$. Using the equivalence of $(\alpha\wedge\beta)\rightarrow\gamma$ and $\alpha\rightarrow (\neg\beta\vee\gamma)$, pull out any negative conjuncts from $\psi_k$, negate them, and disjoin them with the consequent. Then the consequent is still positive, so the result is a safe Sahlqvist implication $\psi_k'\rightarrow\chi'$ where $\psi_k'$ is a conjunction of boxed atoms and $\chi'$ is positive. Now we determine the local first-order correspondent of $\psi_k'\rightarrow\chi'$ as in the standard Sahlqvist-van Benthem algorithm, as presented in \citealt{Blackburn2001}, Theorems 3.40, 3.42, and 3.49, except for three main differences evident from Example \ref{MinVal1}. First, our second-order translation has the $RO(\dot{Q})$ constraint on each predicate variable $\dot{Q}$; but these disappear as in Example \ref{MinVal1} when we substitute the description of the minimal valuation in for $\dot{Q}$. Second, our implications introduce universal quantification $\forall \mathrm{x}\sqsubseteq\mathrm{x}_0$ as in Example \ref{MinVal1}, but that poses no problem. Third, our minimal valuation is the first-order definable $\mathrm{int}(\mathrm{cl}(\mathord{\Downarrow}X))$ where $X$ is the minimal valuation used in Blackburn et al.'s proof.
\end{proof} 

As noted above, the restriction that diamonds not appear in the antecedents of implications is not necessary. We will illustrate this in a way that is useful for calculating first-order correspondents of many standard modal axioms. Where $\sigma$ is a sequence $i_1,\dots,i_n$ of modal operator indices from $\ind$, let $\Diamond_\sigma \varphi$ be $\Diamond_{i_1}\dots\Diamond_{i_n} \varphi$ and $\Box_\sigma\varphi$ be $\Box_{i_1}\dots\Box_{i_n}\varphi$. If $\sigma$ is the empty sequence, then $\Diamond_\sigma\varphi$ and $\Box_\sigma\varphi$ are simply $\varphi$. We will prove a possibility-semantic analogue of the correspondence result from Lemmon and Scott \citeyearpar[\S~4]{Lemmon1977} for formulas of the form $\Diamond_\alpha\Box_\beta p\rightarrow\Box_\delta\Diamond_\gamma p$, where $\alpha$, $\beta$, $\gamma$, and $\delta$ are sequences of indices from $\ind$.\footnote{The result in \citealt{Lemmon1977} is only stated for the unimodal language, but the polymodal version is a straightforward generalization.}

First, we recall the original result of Lemmon and Scott for Kripke frames. Given a Kripke frame $\mathfrak{F}=\langle \wo{W},\{\wo{R}_i\}_{i\in\ind}\rangle$ and a sequence $\alpha=\langle i_1,\dots,i_n\rangle$ of indices from $\ind$, let $x\mathrm{R}_\alpha y$ iff there are $x_0,\dots, x_{n}$ such $x_0=x$, $x_n=y$, and $x_0\mathrm{R}_{i_1}x_1$, $x_1\mathrm{R}_{i_2}x_2$, $\dots$, $x_{n-1}\mathrm{R}_{i_n}x_n$. If $\alpha$ is the empty sequence, then $x\mathrm{R}_\alpha y$ iff $x=y$.

In the following, we will blur the distinction between our metalanguage for talking about frame conditions and our formal first-order language $\mathcal{L}^1(\ind)$, trusting that no confusion will arise.

\begin{proposition}[Lemmon-Scott Kripke Frame Correspondence]\label{LS} Let $\mathfrak{F}$ be a Kripke frame. Then for any sequences $\alpha$, $\beta$, $\delta$, and $\gamma$ of indices from $\ind$, $\Diamond_\alpha\Box_\beta p\rightarrow \Box_\delta\Diamond_\gamma p$ is valid over $\mathfrak{F}$ iff $\mathfrak{F}$ satisfies:
\begin{equation}\forall {x}\forall {y}\forall {z} (({x}\mathrm{R}_\delta{y}\wedge {x}\mathrm{R}_\alpha {z})\rightarrow \exists {u} ( {y} \mathrm{R}_\gamma {u} \wedge {z}\mathrm{R}_\beta {u})).\label{LSeq}\end{equation}
For a \textit{local} correspondent of $\Diamond_\alpha\Box_\beta p\rightarrow \Box_\delta\Diamond_\gamma p$, simply delete the $\forall x$ from (\ref{LSeq}).
\end{proposition}

For possibility frames, we have a very similar result in Proposition \ref{Gklmn}. Given a possibility frame $\mathcal{F}=\langle S,\sqsubseteq,\{R_i\}_{i\in\ind},\adm\rangle$ and a sequence $\alpha=\langle i_1,\dots,i_n\rangle$ of indices from $\ind$, define $xR_\alpha y$ as above except that if $\alpha$ is the empty sequence, let $xR_\alpha y$ iff $y\sqsubseteq x$. Note that the \Rdown{} condition of standard frames implies:
\begin{itemize}
\item if $xR_\alpha y$ and $y'\sqsubseteq y$, then $xR_\alpha y'$.
\end{itemize}
For the empty sequence, this means that if $y\sqsubseteq x$ and $y'\sqsubseteq y$, then $y'\sqsubseteq x$, which is just the transitivity of $\sqsubseteq$. 

We are now ready to prove the analogue of Proposition \ref{LS} in possibility semantics. 

\begin{proposition}[Lemmon-Scott Possibility Frame Correspondence]\label{Gklmn} Let $\mathcal{F}$ be a full standard possibility frame. Then for any sequences $\alpha$, $\beta$, $\delta$, and $\gamma$ of indices from $\ind$, $\Diamond_\alpha  \Box_\beta  p\rightarrow \Box_\delta   \Diamond_\gamma  p$ is valid over $\mathcal{F}$ iff $\mathcal{F}$ satisfies
\begin{equation}\forall{x}\forall{y}\big({x}R_\delta {y}\rightarrow \exists {x'}\sqsubseteq{x} \;\forall {z}( {x'}R_\alpha   {z} \rightarrow\exists {u} ({y} R_\gamma  {u}\wedge {z} R_\beta  {u}))\big).\label{Gklmn-rel-con}\end{equation} 
For the case where $\alpha$ is empty, $\Box_\beta  p\rightarrow \Box_\delta   \Diamond_\gamma  p$ is valid over $\mathcal{F}$ iff $\mathcal{F}$ satisfies
\begin{equation}\forall{x}\forall{y}({x}R_\delta {y}\rightarrow \exists {u} ({y} R_\gamma  {u}\wedge x R_\beta  {u})).\label{Gklmn-rel-con3}\end{equation} 
For \textit{local} correspondents, replace $\forall \mathrm{x}$ in (\ref{Gklmn-rel-con}) and (\ref{Gklmn-rel-con3}) with $\forall \mathrm{x}\sqsubseteq\mathrm{x}_0$ so $\mathrm{x}_0$ is free.
\end{proposition}

\begin{proof} Suppose $\mathcal{F}$ satisfies (\ref{Gklmn-rel-con}). Further suppose that there is a model $\mathcal{M}$ based on $\mathcal{F}$ and an $x\in\mathcal{M}$ such that $\mathcal{M},x\nVdash \Box_\delta   \Diamond_\gamma  p$. Then as in Figure \ref{GklmnFig}, there is a $y'$ such that $xR_\delta y'$ and $\mathcal{M},y'\nVdash \Diamond_\gamma  p$, i.e., $\mathcal{M},y'\nVdash \neg \Box_\gamma  \neg p$, so there is a $y\sqsubseteq y'$ such that $\mathcal{M},y\Vdash \Box_\gamma  \neg p$. By \Rdown{}, $xR_\delta    y'$ implies $xR_\delta    y$. Now take an $x'$ as in (\ref{Gklmn-rel-con}) and consider any $z'$ such that $x'R_\alpha   z'$. We claim that $\mathcal{M},z'\Vdash \neg\Box_\beta  p $. So consider any $z\sqsubseteq z'$. By \Rdown{}, $x'R_\alpha   z'$ implies $x'R_\alpha   z$, so by (\ref{Gklmn-rel-con}) there is a $u$ with $y R_\gamma  u$ and $z R_\beta  u$. Then since $\mathcal{M},y\Vdash \Box_\gamma  \neg p$, $y R_\gamma  u$ implies $\mathcal{M},u\Vdash \neg p$, which with $z R_\beta  u$ implies $\mathcal{M},z\nVdash \Box_\beta p$. Since this holds for all $z\sqsubseteq z'$, we have $\mathcal{M},z'\Vdash \neg\Box_\beta  p $, as claimed. Then since $z'$ was an arbitrary $R_\alpha  $-successor of $x'$, we have $\mathcal{M},x'\Vdash \Box_\alpha   \neg\Box_\beta  p $, which with $x'\sqsubseteq x$ implies $\mathcal{M},x\nVdash \neg \Box_\alpha   \neg\Box_\beta  p$, i.e., $\mathcal{M},x\nVdash  \Diamond_\alpha   \Box_\beta  p$. Thus, $\Diamond_\alpha  \Box_\beta  p\rightarrow \Box_\delta   \Diamond_\gamma  p$ is valid over $\mathcal{F}$.

 \begin{figure}[h]
\begin{center}
\begin{tikzpicture}[->,>=stealth',shorten >=1pt,shorten <=1pt, auto,node
distance=2cm,thick,every loop/.style={<-,shorten <=1pt}]
\tikzstyle{every state}=[fill=gray!20,draw=none,text=black]
\node (x) at (0,0) {{$x$}};
\node (y') at (2,0) {{$y'$}};
\node (y'dash) at (2.75,-0.05) {{$\nVdash \Diamond_\gamma p$}};
\node (y) at (2,-2) {{$y$}};
\node (ydash) at (2.85,-2) {{$\Vdash \Box_\gamma \neg p$}};

\node (x') at (0,-4) {{$x'$}};
\node (z') at (2,-4) {{$z'$}};
\node (y'dash) at (2.85,-4.05) {{$\Vdash \neg\Box_\beta p$}};
\node (z) at (2,-6) {{$z$}};
\node (zdash) at (2.75,-6.05) {{$\nVdash \Box_\beta p$}};

\node (u) at (5,-4) {{$u$}};
\node (udash) at (5.65,-4) {{$\Vdash \neg p$}};

\path (x) edge[->] node {{}} (x');
\path (y') edge[->] node {{$\exists$}} (y);
\path (z') edge[->] node {{$\forall$}} (z);

\path (x) edge[->,dashed] node {{$\delta$}} (y');
\path (x) edge[->,dashed] node {{$\delta$}} (y);

\path (x') edge[->,dashed] node {{$\alpha$}} (z');
\path (x') edge[->,dashed] node {{$\alpha$}} (z);

\path (y) edge[->,dashed] node {{$\gamma$}} (u);
\path (u) edge[<-,dashed] node {{$\beta$}} (z);

\end{tikzpicture}
\end{center}
\caption{diagram for the proof of Proposition \ref{Gklmn}.}\label{GklmnFig}
\end{figure}
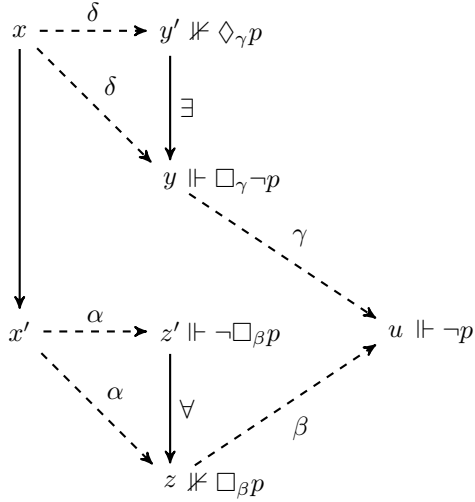

In the other direction, suppose $\mathcal{F}$ does not satisfy (\ref{Gklmn-rel-con}), so it satisfies
\begin{equation}
\exists{x}\exists{y}\big({x}R_\delta {y}\wedge \forall {x'}\sqsubseteq{x} \;\exists {z}( {x'}R_\alpha   {z} \wedge\forall {u} ({y} R_\gamma  {u}\rightarrow\neg {z} R_\beta  {u}))\big).\label{notGklmn}
\end{equation}
Define $\mathcal{M}=\langle \mathcal{F},\pi\rangle$ such that $s\in \pi(p)$ iff for all $t\in R_\gamma (y)$, $s\incomp t$, so by Fact \ref{ConVal}, $\mathcal{M}$ is an admissible model based on $\mathcal{F}$.  Now consider any $x'\sqsubseteq x$, so we have a $z$ as in (\ref{notGklmn}). We claim that $\mathcal{M},z\Vdash\Box_\beta  p$, so consider any $u'$ with $zR_\beta u'$. Suppose for reductio that $u'\not\in \pi(p)$, so there is a $t\in R_\gamma (y)$ with $u'\comp t$, so there is a $u\sqsubseteq u'$ with $u\sqsubseteq t$. By \Rdown{}, $zR_\beta u'$ and $u\sqsubseteq u'$ together imply $zR_\beta u$, and  $y R_\gamma  t$ and $u\sqsubseteq t$ together imply $y R_\gamma  u$. But $zR_\beta u$ and $y R_\gamma  u$ together contradict (\ref{notGklmn}). Thus, $u'\in \pi(p)$.  Since $u'$ was an arbitrary $R_\beta $-successor of $z$, $\mathcal{M},z\Vdash \Box_\beta p$, which with  $x'R_\alpha   z$ implies $\mathcal{M},x'\nVdash \Box_\alpha   \neg  \Box_\beta p$. Since this holds for all $x'\sqsubseteq x$, $\mathcal{M},x\Vdash \neg \Box_\alpha   \neg  \Box_\beta p$, i.e., $\mathcal{M},x\Vdash  \Diamond_\alpha     \Box_\beta p$. By definition of $\pi$ together with  \Rdown{}, we also have $\mathcal{M},y\Vdash \Box_\gamma \neg p$ and hence $\mathcal{M},y\nVdash \Diamond_\gamma  p$, which with $xR_\delta    y$ implies $\mathcal{M},x\nVdash \Box_\delta    \Diamond_\gamma  p$. Thus, $\Diamond_\alpha  \Box_\beta  p\rightarrow \Box_\delta   \Diamond_\gamma  p$ is not valid over $\mathcal{F}$.

For the claim about the empty $\alpha$ case, suppose $\mathcal{F}$ satisfies (\ref{Gklmn-rel-con3}) and there is a model $\mathcal{M}$ based on $\mathcal{F}$ and an $x\in\mathcal{M}$ such that $\mathcal{M},x\nVdash \Box_\delta   \Diamond_\gamma  p$. Thus, there is a $y'$ such that $xR_\delta    y'$ and $\mathcal{M},y'\nVdash \Diamond_\gamma  p$, i.e., $\mathcal{M},y'\nVdash \neg \Box_\gamma  \neg p$, so there is a $y\sqsubseteq y'$ such that $\mathcal{M},y\Vdash \Box_\gamma  \neg p$. By \Rdown{}, $xR_\delta    y'$ implies $xR_\delta    y$. Then by (\ref{Gklmn-rel-con3}) there is a $u$ with $y R_\gamma  u$ and $x R_\beta  u$. Then since $\mathcal{M},y\Vdash \Box_\gamma  \neg p$, $y R_\gamma  u$ implies $\mathcal{M},u\Vdash \neg p$, which with $x R_\beta  u$ implies $\mathcal{M},x\nVdash \Box_\beta p$. Thus, $\Box_\beta  p\rightarrow \Box_\delta   \Diamond_\gamma  p$ is valid over $\mathcal{F}$. 

In the other direction for the empty $\alpha$ case, suppose $\mathcal{F}$ does not satisfy (\ref{Gklmn-rel-con3}), so
\begin{equation}\exists{x}\exists{y}({x}R_\delta {y}\wedge \forall {u} ({y} R_\gamma  {u}\rightarrow \neg x R_\beta  {u})).\label{Gklmn-rel-con3NOT}\end{equation} 
Define $\mathcal{M}$ exactly as we did after (\ref{notGklmn}). Then by the same form of argument used to show that $\mathcal{M},z\Vdash \Box_\beta p$ in the proof from (\ref{notGklmn}), one can show that $\mathcal{M},x\Vdash\Box_\beta  p$ from (\ref{Gklmn-rel-con3NOT}). Moreover, the exact same argument used to show $\mathcal{M},x\nVdash \Box_\delta    \Diamond_\gamma  p$ above works here. Thus, $\Box_\beta  p\rightarrow \Box_\delta   \Diamond_\gamma  p$ is not valid over $\mathcal{F}$.

Finally, for the claim about local correspondence, where $\forall \mathrm{x}\,\psi(\mathrm{x})$ is (\ref{Gklmn-rel-con}), we claim that $\forall \mathrm{x}\sqsubseteq\mathrm{x}_0\, \psi(\mathrm{x})$ is a local correspondent of $\Diamond_\alpha\Box_\beta p\rightarrow \Box_\delta\Diamond_\gamma p$. If there is a model $\mathcal{M}$ based on $\mathcal{F}$ and an $x_0\in\mathcal{M}$ such that $\mathcal{M},x_0\nVdash \Diamond_\alpha\Box_\beta p\rightarrow \Box_\delta\Diamond_\gamma p$, then there is an $x\sqsubseteq x_0$ such that $\mathcal{M},x\Vdash \Diamond_\alpha\Box_\beta p$ but $\mathcal{M},x\nVdash \Box_\delta\Diamond_\gamma p$. Then by the same reasoning as in the first paragraph of this proof, it follows that $\psi(\mathrm{x})$ is falsified by $x$, so $\forall \mathrm{x}\sqsubseteq\mathrm{x}_0\,\psi(\mathrm{x})$ is falsified by $x_0$. In the other direction, if the local correspondent is not satisfied at a state $x_0$ in $\mathcal{F}$, then there is an $x\sqsubseteq x_0$ such that $x$ satisfies the formula obtained from (\ref{notGklmn}) by deleting the initial existential quantifier. Then the proof that $\mathcal{M},x\Vdash  \Diamond_\alpha\Box_\beta p$ and $\mathcal{M},x\nVdash \Box_\delta \Diamond_\gamma p$ proceeds as before.\end{proof}

Note that the proof that the Lemmon-Scott axiom is valid if $\mathcal{F}$ satisfies (\ref{Gklmn-rel-con})/(\ref{Gklmn-rel-con3}) does not rely on the assumption that $\mathcal{F}$ is \textit{full}; thus it holds for all standard possibility frames.

Before applying Proposition \ref{Gklmn} to some example axioms, it is worth seeing how the frame conditions in Propositions \ref{LS} and \ref{Gklmn} relate between a Kripke frame and its \textit{powerset possibilization} (Example \ref{PowerPoss}).

\begin{example}[Frame Properties and Powerset Possibilization]\label{FramePropsIII} For any Kripke frame $\mathfrak{F}=\langle \wo{W},\{\wo{R}_i\}_{i\in\ind }\rangle$ and sequences $\alpha$, $\beta$, $\delta$, and $\gamma$ of indices from $\ind$, the condition 
\begin{equation}\forall {x}\forall {y}\forall {z}(({x}\mathrm{R}_\delta  {y}\wedge {x}\mathrm{R}_\alpha {z})\rightarrow \exists {u} ( {y} \mathrm{R}_\gamma  {u} \wedge {z}\mathrm{R}_\beta  {u})\big)\label{Cond1}\end{equation}
holds of $\mathfrak{F}$ iff the condition 
\begin{equation}\forall X \forall Y (  X R_\delta  Y \rightarrow \exists  X' \sqsubseteq X  \;\forall  Z (  X' R_\alpha  Z  \rightarrow\exists  U  ( Y  R_\gamma   U \wedge  Z  R_\beta   U ))\big)\label{Cond2}\end{equation}
holds of the powerset possibilization $\mathfrak{F}^\pow=\langle S, \sqsubseteq, \{R_i\}_{i\in\ind},\adm\rangle$.\footnote{Note that despite the capitalization of `$X$', `$Y$', etc., we are using them as first-order variables ranging over the domain of $\mathfrak{F}^\pow$. We are using capital letters only as a reminder that the possibilities in $\mathfrak{F}^\pow$ are sets of worlds.} This follows from Propositions \ref{LS} and \ref{Gklmn} plus the fact that $\mathfrak{F}$ and $\mathfrak{F}^\pow$ validate the same formulas (Fact \ref{WtoP1}.\ref{WtoP1b}), but it is also easy to prove directly.

First, recall that for possibility frames, we defined $XR_i^0Y$ as $Y\sqsubseteq X$, which for $\mathfrak{F}^\pow$ means $Y\subseteq X$, which is equivalent to $Y\subseteq\mathrm{R}_i^0[X]$, since we defined $x\mathrm{R}_i^0y$ as $x=y$ for Kripke frames.

Now suppose $\mathfrak{F}$ satisfies (\ref{Cond1}), and consider ${X},{Y}\in \wp(\wo{W})\setminus\{\emptyset\}$ with $XR_\delta Y$, so ${Y}\subseteq\mathrm{R}_\delta [{X}]$. Let ${X'}=\{x\in {X}\mid \exists y\in{Y}\colon x\mathrm{R}_\delta y\}$, so $\emptyset\not={X'}\subseteq{X}$. Now consider any $Z\in \wp(\wo{W})\setminus\{\emptyset\}$ with $X'R_\alpha Z$, so ${Z}\subseteq \mathrm{R}_\alpha[{X'}]$, and pick an $x\in {X'}$ such that for some $z\in {Z}$, $x\mathrm{R}_\alpha z$. Since $x\in{X'}$, we also have a $y\in{Y}$ with $x\mathrm{R}_\delta y$. Then by (\ref{Cond1}), there is a $u$ such that $y\mathrm{R}_\gamma  u$ and $z\mathrm{R}_\beta u$. Setting ${U}=\{u\}$, we have ${U}\subseteq \mathrm{R}_\gamma [{Y}]$ and  ${U}\subseteq \mathrm{R}_\beta [{Z}]$, so ${Y} R_\gamma  {U}$ and ${Z} R_\beta  {U}$. Thus, $\mathfrak{F}^\pow$ satisfies (\ref{Cond2}).

Next, suppose $\mathfrak{F}^\pow$ satisfies (\ref{Cond2}), and consider $x,y,z\in \wo{W}$ with $x\mathrm{R}_\delta  y$ and $x\mathrm{R}_\alpha z$. Then setting ${X}=\{x\}$ and ${Y}=\{y\}$, there is a ${X'}\sqsubseteq{X}$ for which the rest of (\ref{Cond2}) holds. Now ${X'}\sqsubseteq{X}=\{x\}$ implies ${X'}=\{x\}$. Then since $x\mathrm{R}_\alpha z$ gives us $\{x\}R_\alpha \{z\}$, setting ${Z}=\{z\}$ in (\ref{Cond2}) gives us a ${U}\not=\emptyset$ such that $\{y\}R_\gamma {U}$ and $\{z\}R_\beta {U}$, which implies there is a $u\in{U}$ such that $y\mathrm{R}_\gamma  u$ and $z\mathrm{R}_\beta  u$. Thus, $\mathfrak{F}$ satisfies (\ref{Cond1}).\hfill $\triangleleft$\end{example}

Example \ref{FramePropsIII} enables a quick proof that where \textbf{L} is the least normal extension of \textbf{K} with a set of Lemmon-Scott axioms, \textbf{L} is complete with respect to the class of \textit{rich} possibility frames (\S~\ref{RichFrames}) satisfying the conditions corresponding to those axioms as in Proposition \ref{Gklmn}. For each such \textbf{L} is complete with respect to the class of Kripke frames satisfying the corresponding conditions in Proposition \ref{LS} \citep[\S~4]{Lemmon1977}; powerset possibilization preserves the truth of modal formulas (Fact \ref{WtoP1}.\ref{WtoP1a}); and powerset possibilizations are a special case of rich frames.

In Example \ref{AxiomExamples} below, we check how Proposition \ref{Gklmn} applies to some familiar axioms. While one typically visualizes the Kripke-frame conditions corresponding to modal axioms two-dimensionally, it may help to visualize the possibility-frame conditions corresponding to modal axioms three-dimensionally, e.g., with accessibility arrows parallel to the x-y plane and refinement arrows parallel to the x-z plane. The following fact is useful for simplifying some of the first-order frame conditions produced by Proposition \ref{Gklmn}.

\begin{fact}[\textbf{\textit{R}-common}]\label{CommonFact} Every full standard possibility frame satisfies \textbf{\textit{R}-common}: if $x'\sqsubseteq x$ and $x'R_iy'$, then $\exists z\sqsubseteq y'$: $x'R_iz$ and $xR_i z$.
\end{fact}
\begin{proof} If $x'\sqsubseteq x$ and $x'R_iy'$, then by the \Rrule{} property of full frames (\S~\ref{FullFrames}), $\exists y$: $xR_iy\comp y'$, so $\exists z$: $z\sqsubseteq y$ and $z\sqsubseteq y'$. Then since $xR_iy$ and $x'R_iy'$, we have $xR_iz$ and $x'R_iz$ by \Rdown{} for standard frames.
\end{proof}

Some of the first-order correspondents produced by Proposition \ref{Gklmn} have even simpler forms over \textit{strong} possibility frames (Definition \ref{StrongPoss}), which we will note in the following example.

\begin{example}[Familiar Axioms]\label{AxiomExamples} Consider the following special cases of the correspondence between $\Diamond_\alpha  \Box_\beta  p\rightarrow \Box_\delta   \Diamond_\gamma  p$ and $\forall{x}\forall{y}\big({x}R_\delta {y}\rightarrow \exists {x'}\sqsubseteq{x} \;\forall {z}( {x'}R_\alpha   {z} \rightarrow\exists {u} ({y} R_\gamma  {u}\wedge {z} R_\beta  {u}))\big)$, and between $\Box_\beta  p\rightarrow \Box_\delta   \Diamond_\gamma  p$ and $\forall{x}\forall{y}({x}R_\delta {y}\rightarrow \exists {u} ({y} R_\gamma  {u}\wedge x R_\beta  {u}))$, over full standard frames from Proposition \ref{Gklmn} (recall that if $\sigma$ is empty, then $xR_\sigma y$ means $y\sqsubseteq x$):
\begin{itemize}
\item $\Box_i p\rightarrow \Diamond_i p$ corresponds to $\forall{x}\forall{y}(y\sqsubseteq x\rightarrow \exists {u} ({y} R_i {u}\wedge x R_i {u}))$, which by Fact \ref{CommonFact} is equivalent to \textit{seriality}, $\forall x\,\exists u\, xR_iu$.

\item $\Diamond_i p\rightarrow\Box_i p$ corresponds to $\forall x\forall y (xR_iy\rightarrow \exists x'\sqsubseteq x\, \forall z (x'R_i z\rightarrow \exists u (u\sqsubseteq y\wedge u\sqsubseteq z)))$, i.e., $\forall x\forall y (xR_iy\rightarrow \exists x'\sqsubseteq x\, \forall z (x'R_i z\rightarrow y\comp z))$.

\item $\Box_ip\rightarrow p$ corresponds to $\forall{x}\forall y\sqsubseteq x\,\exists u\sqsubseteq y\,x R_i u$, which by Fact \ref{CommonFact} is equivalent to $\forall{x}\,\exists u\sqsubseteq x \, x R_i u$. Over \textit{strong} frames, it is equivalent to \textit{reflexivity}, $\forall x\, xR_ix$, by \Rdense{} and \Rdown{}.

\item  $p\rightarrow \Box_i p$ corresponds to $\forall x\forall y (xR_i y\rightarrow \exists u (u\sqsubseteq y\wedge u\sqsubseteq x))$, i.e., $\forall x\forall y(xR_i y\rightarrow x\comp y)$.

\item $\Box_i p\rightarrow \Box_i\Box_ip$ corresponds to $\forall{x}\forall{y}({x}R_i^2{y}\rightarrow \exists {u} (u\sqsubseteq y\wedge x R_i {u}))$. Over strong frames, this is equivalent to \textit{transitivity}, $\forall x\forall y (xR_i^2y\rightarrow xR_iy)$, by \Rdense{} and \Rdown{}.

\item $p\rightarrow \Box_i\Diamond_i p$ corresponds to $\forall{x}\forall{y}({x}R_i{y}\rightarrow \exists {u} ({y} R_i {u}\wedge u\sqsubseteq x))$.\footnote{Humberstone \citeyearpar[p.~329f]{Humberstone1981} identifies a different first-order condition that together with his  \RrefPlusPlus{} (recall Remark \ref{HumbFrame}) defines a class of frames relative to which the minimal normal modal logic including $p\rightarrow \Box_i\Diamond_i p$ is sound and complete---but with no claim of correspondence.}

\item $\Diamond_i p\rightarrow \Box_i\Diamond_i p$ corresponds to $\forall{x}\forall{y}\big({x}R_i{y}\rightarrow \exists {x'}\sqsubseteq{x} \;\forall {z}( {x'}R_i {z} \rightarrow\exists {u} ({y} R_i {u}\wedge u\sqsubseteq z))\big)$. Over strong frames, this becomes $\forall{x}\forall{y}\big({x}R_i{y}\rightarrow \exists {x'}\sqsubseteq{x} \;\forall {z}( {x'}R_i {z} \rightarrow yR_i z )\big)$ by \Rdense{} and \Rdown{}.
\item  $\Box_i p\rightarrow \Box_j p$ corresponds to $\forall x\forall y(xR_j y\rightarrow \exists u (u\sqsubseteq y\wedge xR_i u))$. Over strong frames, this is equivalent to $\forall x\forall y(xR_jy\rightarrow xR_iy)$ by \Rdense{} and \Rdown{}.
\end{itemize}
If we combine Proposition \ref{Gklmn} with Lemma \ref{IndCor} and Fact \ref{LocalVsGlobal}, we can treat other familiar axioms such as:
\begin{itemize}
\item $\Box_i(\Box_i p \rightarrow p)$ (globally) corresponds to $\forall x\forall y(xR_iy\rightarrow Loc_y(\Box_ip\rightarrow p))$, where $Loc_y(\Box_ip\rightarrow p)$ is the local correspondent with respect to $y$ given by Proposition \ref{Gklmn}. Given \Rdown{}, the result is equivalent to $\forall x\forall y(xR_iy\rightarrow \exists u\sqsubseteq y\; yR_iu)$. Over strong frames, this is equivalent to \textit{shift-reflexivity}, $\forall x\forall y (xR_iy\rightarrow yR_iy)$, by \Rdense{}, \Rdown{}, and \upR{}. \hfill $\triangleleft$
\end{itemize}
\end{example}

Note that when we apply the first-order conditions above to the special case of Kripke frames regarded as possibility frames as in Examples \ref{KripkeExample} and \ref{KripkeAgain}, so $\sqsubseteq$ is the identity relation, then these conditions reduce to the familiar conditions corresponding to the axioms over Kripke frames (recall Fact \ref{PossToWorldCor}). 

Also note that over strong possibility frames, the correspondents for the axioms in Example \ref{AxiomExamples} with boxes in the antecedent and without any diamonds are already the same as the familiar correspondents over Kripke frames.\footnote{On this point, we can add a follow-up to Example \ref{FramePropsIII}: for axioms of the form $\Box_\beta p\rightarrow \Box_\delta p$, a Kripke frame $\mathfrak{F}=\langle \wo{W},\{\wo{R}_i\}_{i\in\ind}\rangle$ satisfies the corresponding condition $\forall x \forall y (x\mathrm{R}_\delta  y\rightarrow x\mathrm{R}_\beta  y)$ iff its powerset possibilization $\mathfrak{F}^\pow=\langle S, \sqsubseteq, \{R_i\}_{i\in\ind},\adm\rangle$ satisfies $\forall X \forall Y (XR_\delta  Y\rightarrow XR_\beta  Y)$, i.e., $\forall X\forall Y(Y\subseteq \wo{R}_\delta [X]\rightarrow Y\subseteq \wo{R}_\beta [X])$.} This observation can be turned into a general result, but we will not go into it here. (See \citealt[Ch.~4]{Kojima2012} for a study of when a modal formula has the same first-order correspondent over \textit{intuitionistic} modal frames as over classical Kripke frames. We leave for future work a comparison of first-order correspondence over intuitionistic modal frames vs. classical possibility frames.)

Finally, let us note how the first-order correspondents look over \textit{functional} frames (\S~\ref{FuncFrames}). Using Lemma \ref{TransferCorr}, all of our correspondence results can be applied to functional full possibility frames. But over functional frames the correspondents can be further simplified, especially if we assume the frames are \textit{separative} as in \S~\ref{SepSec}. Recall that in a separative full frame, for any $x\in S$, $\mathord{\downarrow}x$ is an admissible proposition (Fact \ref{Sep&Princ}). Thus, in a full frame that is both separative and functional, a minimal valuation $\pi$ making $\Box_i p$ true at $x$ is such that $\pi(p)=\mathord{\downarrow}f_i(x)$, or $\pi(p)=\emptyset$ if $f_i(x)$ is undefined. For simplicity, let us consider frames in which each $f_i$ is a total function, which validate the D axiom $\Box_i p\rightarrow \Diamond_ip$. Note that in functional frames the $\Diamond_i$ clause (recall \S~\ref{Acc&Poss}) is that $\mathcal{M},x\Vdash \Diamond_i \varphi$ iff $\forall x'\sqsubseteq x$ $\exists y\sqsubseteq f_i(x')$: $\mathcal{M},y\Vdash \varphi$. Also note that if we have a functional frame $\mathcal{F}$ and obtain the modally equivalent standard frame $\mathcal{F}_{\mathord{\downarrow}}$ as in Lemma \ref{TransferCorr}, then $xR_{i\mathord{\downarrow}}y$ iff $y\sqsubseteq f_i(x)$. This equivalence is helpful when comparing the first-order correspondents in Examples \ref{AxiomExamples} and  \ref{FuncCorEx}.

\begin{example}[Functional Correspondence]\label{FuncCorEx} Over separative full possibility frames in which each accessibility relation is a total function $f_i$, we have the following global correspondences:
\begin{itemize}
\item $\Diamond_i p\rightarrow\Box_ip$ corresponds to $\forall x\forall y\sqsubseteq f_i(x)\, \exists x'\sqsubseteq x\,\forall z\sqsubseteq f_i(x')\; y\comp z$.
\item $\Box_i p\rightarrow p$ corresponds to $\forall x\, x\sqsubseteq f_i(x)$.
\item $p\rightarrow \Box_i p$ corresponds to $\forall x\, f_i(x)\sqsubseteq x$.
\item $\Box_i p\rightarrow \Box_i\Box_i p$ corresponds to $\forall x\, f_i(f_i(x))\sqsubseteq f_i(x)$ (cf.~Example \ref{BethFlowEx}).
\item $p\rightarrow\Box_i\Diamond_i p$ corresponds to $\forall x\forall y\sqsubseteq f_i(x)\, \exists x'\sqsubseteq x\; x'\sqsubseteq f_i(y)$.
\item $\Diamond_ip\rightarrow \Box_i\Diamond_i p$ corresponds to $\forall x\forall y\sqsubseteq f_i(x)\, \exists x'\sqsubseteq x\; f_i(x')\sqsubseteq f_i(y)$.
\item $\Box_i p\rightarrow \Box_j p$ corresponds to $\forall x \, f_j(x)\sqsubseteq f_i(x)$.
\item $\Box_i (\Box_i p\rightarrow p)$ corresponds to $\forall x\forall y\sqsubseteq f_i(x)\; y\sqsubseteq f_i(y)$.\hfill $\triangleleft$ 
\end{itemize}
\end{example}

\section{Beginnings of Completeness Theory}\label{CompletenessTheory}

In this section, we begin the study of classes of logics complete with respect to classes of possibility frames. Our starting point is Theorem \ref{Linked}, which follows immediately from results we have already established.

\
\begin{theorem}[From Algebraic to Possibility Completeness]\label{Linked} Let $\mathbf{L}$ be a consistent normal modal logic.
\begin{enumerate}
\item\label{Linked1} $\mathbf{L}$ is sound and complete with respect to a \textit{filter-descriptive} possibility frame;
\item\label{Linked2} if $\mathbf{L}$ is sound and complete with respect to a class of $\mathcal{V}$-BAOs, then $\mathbf{L}$ is sound and complete with respect to a class of \textit{principal} possibility frames (and vice versa);
\item\label{Linked3} if $\mathbf{L}$ is sound and complete with respect to a class of $\mathcal{T}$-BAOs, then $\mathbf{L}$ is sound and complete with respect to a class of \textit{functional principal} possibility frames (and vice versa);
\item\label{Linked4} if $\mathbf{L}$ is sound and complete with respect to a class of $\mathcal{CV}$-BAOs, then $\mathbf{L}$ is sound and complete with respect to a class of \textit{full} (indeed \textit{rich}) possibility frames (and vice versa).
\end{enumerate}
\end{theorem}

\begin{proof} For part \ref{Linked1}, $\mathbf{L}$ is sound and complete with respect to a BAO (Theorem \ref{AdAlg}), which can be turned into a modally equivalent possibility frame (Theorem \ref{SatGenFil}) that is filter-descriptive (Proposition \ref{GenAb}.\ref{GenAb1}).

Parts \ref{Linked2}-\ref{Linked4} follow from Theorem \ref{VtoPossFrames} and Proposition \ref{QtoF}.\ref{QtoF5} (for the move from quasi-functional to functional completeness), with the \textit{vice versa} reminders coming from Corollary \ref{TransferPossAlg}.\end{proof}

Using Theorem \ref{Linked}.\ref{Linked2}-\ref{Linked4}, we can transfer knowledge about the algebraic completeness notions of $\mathcal{V}$-completeness, $\mathcal{T}$-completeness, and $\mathcal{CV}$-completeness to knowledge about the associated possibility-semantic completeness notions. As shown in \citealt{Litak2005}, these three completeness notions are distinct and generalize Kripke-completeness; and as shown in \citealt{Litak2015}, the most general of these notions, $\mathcal{V}$-completeness, is still a nontrivial notion of completeness---not all normal modal logics are $\mathcal{V}$-complete. More detailed knowledge about $\mathcal{CV}$- and $\mathcal{V}$-completeness would be of great value for understanding possibility semantics. At present, only $\mathcal{T}$-completeness is fairly well understood, as we briefly review in \S~\ref{SyntacticProp}.

In the rest of this section, we report some salient results concerning completeness with respect to the following special classes of possibility frames: \textit{full} possibility frames (\S~\ref{CompFull}), \textit{principal} possibility frames (\S~\ref{SyntacticProp}), \textit{atomless} possibility frames (\S~\ref{AtomlessFull}), and finally, \textit{canonical} possibility frames (\S~\ref{Canonical}). 

\subsection{Completeness for Full Possibility Frames}\label{CompFull}

By the fact that every Kripke frame may be regarded as a modally equivalent full possibility frame, every Kripke-frame complete logic is also sound and complete with respect to a class of full possibility frames. As we have seen, the converse does not hold. We already showed in \S~\ref{NoKripke} that there are continuum many full possibility frames for a polymodal language whose logics are pairwise distinct and Kripke-frame \textit{incomplete}. By using our duality theory and the polymodal-to-unimodal reduction techniques of Thomason \citeyearpar{Thomason1974b,Thomason1975b,Thomason1975c} and Kracht and Wolter \citeyearpar{Kracht1999b}, we can now prove the analogous result in the \textit{unimodal} case.  

\UNCOUNT*
  
\begin{proof} Theorem~\ref{NoEquiv} gives us continuum many full possibility frames for a polymodal language whose logics are pairwise distinct and Kripke-frame incomplete. Thus, by Theorem \ref{PtoB}, there are continuum many $\mathcal{CV}$-BAOs for a polymodal language whose logics are pairwise distinct and Kripke-frame incomplete.  

Kracht and Wolter \citeyearpar{Kracht1999b} show that there is a function $sim$, the \textit{Thomason-Simulation}, sending each normal polymodal logic \textbf{L} to a normal \textit{unimodal} logic $\textbf{L}^{sim}$ such that: (i) for each $n\in\mathbb{N}$, $sim$ is an isomorphism from the lattice of normal modal logics with $n$ modal operators onto an interval in the lattice of normal unimodal logics, so for any normal modal logics $\textbf{L}_1$ and $\textbf{L}_2$ for $n$ modal operators, $\textbf{L}_1\subseteq\textbf{L}_2$ iff $\textbf{L}_1^{sim}\subseteq\textbf{L}_2^{sim}$; and (ii) if $\mathbf{L}$ is Kripke-frame incomplete, so is $\mathbf{L}^{sim}$. Kracht and Wolter (p.~116) also give an algebraic characterization of the Thomason-Simulation: given a polymodal BAO $\mathbb{A}$, they define a unimodal BAO $\mathbb{A}^{sim}$ such that if \textbf{L} is the logic of $\mathbb{A}$, then $\textbf{L}^{sim}$ is the logic of $\mathbb{A}^{sim}$. (They only give the definition of $\mathbb{A}^{sim}$ for the case where $\mathbb{A}$ has two operators, but one can work out the general definition from the remarks on p.~121 of their paper.) Now we add the following observation, a proof of which is sketched below: if $\mathbb{A}$ is a $\mathcal{CV}$-BAO, then so is $\mathbb{A}^{sim}$. Thus, since there are continuum many $\mathcal{CV}$-BAOs for a polymodal language whose logics are pairwise distinct and Kripke-frame incomplete, there are also continuum many $\mathcal{CV}$-BAOs for the \textit{unimodal} language whose logics are pairwise distinct (by (i) above) and Kripke-frame incomplete (by (ii) above). Then by Theorem \ref{VtoPossFrames}, there are continuum many full possibility frames for the unimodal language whose logics are pairwise distinct and Kripke-frame incomplete.

Let us now sketch the proof that if $\mathbb{A}$ is a $\mathcal{CV}$-BAO, then so is $\mathbb{A}^{sim}$.  For simplicity, consider the case where $\mathbb{A}$ has only two modal operators. Kracht and Wolter (p.~116) define $\mathbb{A}^{sim}$ as follows. Where $A$ is the Boolean reduct of $\mathbb{A}$ and $2$ is the two-element Boolean algebra, the Boolean reduct $A^{sim}$ of $\mathbb{A}^{sim}$ is $A\times A\times 2$, which is complete if $A$ is complete. Where $\blacksquare_1$ and $\blacksquare_2$ are the dual operators in $\mathbb{A}$, the dual operator $\boxminus$ in $\mathbb{A}^{sim}$ is defined as follows:
\[\boxminus \langle a,b,c\rangle := \begin{cases}\langle b\wedge \blacksquare_1 a,a\wedge\blacksquare_2 b,\{1\}\rangle & \mbox{if }c=\{1\} \\ \langle 0, a\wedge\blacksquare_2 b,\{1\}\rangle &\mbox{if }c=\emptyset \end{cases}.\]
We claim that if $\blacksquare_1$ and $\blacksquare_2$ are completely multiplicative, then so is $\boxminus$. Suppose that $\underset{i\in I}{\bigwedge} \langle a_i,b_i,c_i\rangle$ exists in $\mathbb{A}^{sim}$. Since meets are calculated pointwise, this means that $a:=\underset{i\in I}{\bigwedge} a_i$ and $b:=\underset{i\in I}{\bigwedge} b_i$ exist in $\mathbb{A}$, and $\underset{i\in I}{\bigwedge} \langle a_i,b_i,c_i\rangle = \langle a,b,c\rangle$, where $c:=\underset{i\in I}{\bigwedge} c_i$ in the algebra $2$. Now we claim that $\boxminus \langle a,b,c\rangle$ is the greatest lower bound of $\{\boxminus \langle a_i,b_i,c_i\rangle \mid i\in I\}$. That it is a lower bound is immediate from the monotonicity of $\boxminus$, so we need only show that it is the greatest. 

Case 1: $c=\{1\}$. Then by the definition of $\boxminus$, $a$, and $b$, and the complete multiplicativity of $\blacksquare_1$ and $\blacksquare_2$:
\begin{eqnarray}
\boxminus \langle a,b,c\rangle &=& \langle b\wedge \blacksquare_1 a,a\wedge\blacksquare_2 b,\{1\}\rangle \nonumber \\
 &=& \langle b\wedge \blacksquare_1 \underset{i\in I}{\bigwedge}a_i,a\wedge\blacksquare_2 \underset{i\in I}{\bigwedge}b_i,\{1\}\rangle \nonumber\\
&=& \langle b\wedge \underset{i\in I}{\bigwedge}\blacksquare_1 a_i, a\wedge \underset{i\in I}{\bigwedge}\blacksquare_2 b_i,\{1\}\rangle.\label{SimEq1}
\end{eqnarray}
Now suppose $\langle e,f,g\rangle$ is a lower bound of $\{\boxminus \langle a_i,b_i,c_i\rangle \mid i\in I\}$, so for each $i\in I$, we have $\langle e,f,g\rangle\leq \boxminus \langle a_i,b_i,c_i\rangle$, which implies $\langle e,f,g\rangle\leq \langle b\wedge\blacksquare_1 a_i,a\wedge\blacksquare_2 b_i,\{1\}\rangle$ or $\langle e,f,g\rangle\leq \langle 0,a\wedge\blacksquare_2 b_i,\{1\}\rangle$. In either case, we have $e\leq b\wedge \blacksquare_1 a_i$ and $f\leq a\wedge\blacksquare_2 b_i$. Since this holds for each $i\in I$, we have $e\leq b\wedge \underset{i\in I}{\bigwedge}\blacksquare_1 a_i$ and $f\leq a\wedge\underset{i\in I}{\bigwedge}\blacksquare_2 b_i$, which with (\ref{SimEq1}) implies $\langle e,f,g\rangle \leq \boxminus \langle a,b,c\rangle$.

Case 2: $c=\emptyset$. Then by the definition of $\boxminus$, $b$, and the complete multiplicativity of $\blacksquare_2$, we have:
\begin{eqnarray}\boxminus \langle a,b,c\rangle &=& \langle 0, a\wedge\blacksquare_2 b,\{1\}\rangle = \langle 0, a\wedge\underset{i\in I}{\bigwedge}\blacksquare_2 b_i,\{1\}\rangle.\label{SimEq2}\end{eqnarray}
Again suppose $\langle e,f,g\rangle$ is a lower bound of $\{\boxminus \langle a_i,b_i,c_i\rangle \mid i\in I\}$, so for each $i\in I$, we have $\langle e,f,g\rangle\leq \boxminus \langle a_i,b_i,c_i\rangle$, which implies $\langle e,f,g\rangle\leq \langle b\wedge\blacksquare_1 a_i,a\wedge\blacksquare_2 b_i,\{1\}\rangle$ or $\langle e,f,g\rangle\leq \langle 0,a\wedge\blacksquare_2 b_i,\{1\}\rangle$. Now since $c=\underset{i\in I}\bigwedge c_i=\emptyset$, there is some $i\in I$ for which $c_i=\emptyset$ and hence $\langle e,f,g\rangle\leq \langle 0,a\wedge\blacksquare_2 b_i,\{1\}\rangle$. Thus, $e\leq 0$, and for each $i\in I$, $f\leq a\wedge\blacksquare_2 b_i$, whence $f\leq a\wedge \underset{i\in I}{\bigwedge} \blacksquare_2 b_i$. It follows by (\ref{SimEq2}) that $\langle e,f,g\rangle \leq \boxminus \langle a,b,c\rangle$. This completes the proof that $\boxminus \langle a,b,c\rangle$ is the greatest lower bound.\end{proof}

Thus, Kripke-completeness is a sufficient but far from necessary condition for a unimodal logic to be sound and complete with respect to a class of full possibility frames. Whether we can find weaker sufficient conditions that are illuminating remains to be seen. For the case of completeness with respect to \textit{principal} possibility frames, we will see examples of such sufficient conditions in \S~\ref{SyntacticProp}.

Recall that the \textit{degree of Kripke-incompleteness} of a normal modal logic \textbf{L} for $\mathcal{L}(\sig,\ind)$ is the cardinality of the set of normal modal logics for  $\mathcal{L}(\sig,\ind)$ that are valid over exactly the same Kripke frames as $\mathbf{L}$ \citep{Fine1974b}. More generally, we can consider a notion of \textit{degree of} relative \textit{Kripke-incompleteness}, suggested by Tadeusz Litak (see \citealt{Litak2015}), which we can apply as follows. Fixing the language $\mathcal{L}(\sig,\ind)$, let $\mathrm{ML}(\mathsf{FP})$ be the set of normal modal logics that are sound and complete with respect to some class of full possibility frames. Then the degree of Kripke-incompleteness \textit{relative to }$\mathsf{FP}$ of a normal modal logic \textbf{L} is the cardinality of the set of $\mathbf{L}'\in\mathrm{ML}(\mathsf{FP})$ such that $\mathbf{L}$ and $\mathbf{L}'$ are valid over exactly the same Kripke frames. For example, for a polymodal language, Theorem \ref{NoEquiv} showed that the logic of our full possibility frame $\mathcal{F}$ in \S~\ref{NoKripke} has degree $2^{\aleph_0}$, because it showed that there are continuum many other logics in $\mathrm{ML}(\mathsf{FP})$ that are valid over the same set of Kripke frames, namely the \textit{empty} set. Of course, this special example involving the empty set of frames gives little hint of the general situation. While a great deal is known about degrees of unrelativized incompleteness \citep{Blok1980,Litak2008}, little is currently known about relativized degrees. Knowing more about degrees of Kripke-incompleteness relative to $\mathsf{FP}$ would give us a better measure of how much more fine-grained full possibility frames are than Kripke frames for distinguishing between logics. 

We will see more about completeness for full possibility frames in \S\S~\ref{AtomlessFull}-\ref{Canonical}, where we treat completeness for \textit{atomless} full possibility frames and \textit{canonical} full possibility frames.

\subsection{Completeness for Principal Possibility Frames}\label{SyntacticProp}

We turn now to classes of logics that are sound and complete with respect to \textit{principal} possibility frames (equivalently, $\mathcal{V}$-BAOs), starting with the special case of \textit{functional} principal frames (equivalently, $\mathcal{T}$-BAOs). 

Recall that a \textit{tense logic} is a normal bimodal logic containing the axioms $p\rightarrow \Box_<\Diamond_> p$ and $p\rightarrow \Box_>\Diamond_< p$. For any tense logic \textbf{L}, $\vdash_\mathbf{L} p\rightarrow \Box_< q$  iff $\vdash_\mathbf{L}\Diamond_>p\rightarrow q$, and similarly for $\Box_>$ and $\Diamond_<$. Thus, the Lindenbaum algebra (Definition \ref{LinAlg}) of any consistent tense logic is a $\mathcal{T}$-BAO. Then since every such logic is sound and complete with respect to its Lindenbaum algebra (Theorem \ref{AdAlg}), every such logic is sound and complete with respect to a $\mathcal{T}$-BAO. By contrast, it is well-known that there are tense logics that are not sound and complete with respect to any class of Kripke frames \citep[\S~4.4]{Thomason1972,Blackburn2001}.  

More generally, say that a normal modal logic \textbf{L} is a logic \textit{with converses} iff for every $a\in \ind$ there is a $b\in \ind$ such that $\vdash_\mathbf{L} p\rightarrow \Box_a\Diamond_b p$ and $\vdash_\mathbf{L} p\rightarrow \Box_b\Diamond_a p$ (this class includes what are called the  \textit{connected} logics in \citealt[p.~71]{Kracht1999}). By the same argument as for tense logics, every consistent logic with converses is sound and  complete with respect to a $\mathcal{T}$-BAO, namely its Lindenbaum algebra. Note that since we do not require that the $a$ and $b$ be distinct, any extension of  \textbf{KB} (in its unimodal or polymodal fusion version) is also a logic with converses. Just as there are tense logics that are not sound and complete with respect to any class of Kripke frames, there are also extensions of \textbf{KB}---in fact, continuum many extensions of \textbf{KB}---that are not sound and complete with respect to any class of Kripke frames \citep{Miyazaki2007}.

From the previous two paragraphs, Theorem \ref{VtoPossFrames}, and Proposition \ref{QtoF}.\ref{QtoF5}, we obtain the following.

\begin{corollary}[General Completeness for Logics with Converses]\label{GenTense} 
Every consistent normal modal logic with converses (and hence every consistent tense logic and every consistent extension of \textbf{KB}) is sound and complete with respect to a functional principal possibility frame.
\end{corollary} 

Corollary \ref{GenTense} gives us a simple syntactic characterization of a class of logics that are $\mathcal{T}$-complete. We can give a somewhat more complex syntactic characterization of a strictly larger class of $\mathcal{T}$-complete logics using the following notion studied in \citealt{Holliday2014}.

\begin{restatable}[Internal Adjointness]{definition}{INTAD}\label{intad} A modal logic \textbf{L} has \textit{internal adjointness} iff for every $i\in  \ind $ and $\varphi\in\mathcal{L}(\sig, \ind)$, there is a $\mathsf{f}^\mathbf{L}_i(\varphi)\in\mathcal{L}(\sig, \ind)$ such that for all $\psi\in\mathcal{L}(\sig, \ind)$: \begin{equation}\vdash_\mathbf{L}\varphi\rightarrow\Box_i\psi\mbox{ iff }\vdash_\mathbf{L}\mathsf{f}_i^\mathbf{L}(\varphi)\rightarrow\psi.\tag*{$\triangleleft$}\end{equation} 
\end{restatable}

It is easy to see that internal adjointness is related to $\mathcal{T}$-BAOs as follows.

\begin{lemma}[Internal Adjointness and $\mathcal{T}$-BAOs]\label{AdjClose0} A consistent normal modal logic \textbf{L} has internal adjointness iff the Lindenbaum algebra of $\mathbf{L}$ is a $\mathcal{T}$-BAO.
\end{lemma} 

Then we can strengthen the general completeness result given in Corollary \ref{GenTense} as follows.

\begin{corollary}[General Completeness for Logics with Internal Adjointness]\label{GenAdj0} Every consistent normal modal logic with internal adjointness is sound and complete with respect to a functional principal possibility frame.
\end{corollary}

Corollary \ref{GenAdj0} is stronger than Corollary \ref{GenTense} because every logic with converses has internal adjointness, but there are logics with internal adjointness that are not logics with converses.  \citealt{Holliday2014} determines whether several well-known Kripke-frame complete modal logics have internal adjointness: in addition to all normal modal logics with converses, the logics \textbf{K} (cf.~\citealt[Thm.~6.3]{Ghilardi1995}), \textbf{KD}, \textbf{T}, \textbf{KD4}, and \textbf{S4} have internal adjointness, whereas \textbf{K5} and its extensions with the D and/or 4 axioms lack internal adjointness. It would be interesting to know more about internal adjointness for Kripke-frame incomplete modal logics.

So far we have shown that various logics are $\mathcal{T}$-complete by observing that their Lindenbaum algebras are $\mathcal{T}$-BAOs. In other cases, we may know that a logic is sound and complete with respect to a $\mathcal{T}$-BAO that is not its Lindenbaum algebra. For example, consider the much-used \textit{veiled recession frame} \citep{Makinson1969,Thomason1974,Benthem1978,Blok1979}, which is the world frame $\mathfrak{g}=\langle \wo{W},\wo{R},\wo{A}\rangle$ where $\wo{W}=\mathbb{N}$, $n\wo{R}m$ iff $n-1\leq m $, and $\wo{A}$ is the set of all finite and co-finite subsets of $\mathbb{N}$. It is easy to check that $\mathfrak{g}$ is a general world frame. In addition, $\mathrm{A}$ has the property that if $X\in\wo{A}$, then $\wo{R}[X]\in \wo{A}$, because for any nonempty $X\subseteq\mathbb{N}$, $\wo{R}[X]$ is co-finite, since $\wo{W}\setminus \wo{R}[X]=\{k\in\mathbb{N} \mid k< min(X)-1\}$. From here it is easy to see that the BAO underlying $\mathfrak{g}$ is a $\mathcal{T}$-BAO. But the logic of $\mathfrak{g}$ is not a logic with converses, since $\varphi\rightarrow\Box\Diamond\varphi$ is not valid over $\mathfrak{g}$. Like the logic of any world frame, the logic of the veiled recession frame is a normal modal logic. Moreover, it is finitely axiomatizable \citep{Blok1979}. This logic is \textit{Kripke-frame incomplete} \citep{Blok1979,Blok1980}, but since the BAO underlying $\mathfrak{g}$ is a $\mathcal{T}$-BAO, we have the following result as a corollary of Theorem~\ref{VtoPossFrames}.

\begin{corollary}[Logic of the Veiled Recession Frame] The logic of the veiled recession frame $\mathfrak{g}$ is sound and complete with respect to a quasi-functional principal possibility frame, viz., the principal frame $(\mathfrak{g}^\under)_\rela$ of its underlying BAO $\mathfrak{g}^\under$, and hence with respect to a functional principal possibility frame (by Proposition \ref{QtoF}).
\end{corollary} 

The syntactic conditions mentioned above---being a logic with converses, having internal adjointness---are sufficient for $\mathcal{T}$-completeness, but far from necessary. For a necessary and sufficient syntactic condition, let us say that for a normal unimodal logic \textbf{L}, the \textit{minimal tense extension} $\textbf{L}.t$ of \textbf{L} is the minimal tense logic containing \textbf{L} \citep[\S~3.3]{Kracht1997}. More generally, given a normal modal logic \textbf{L} for the language $\mathcal{L}(\sig,\ind)$, let the minimal tense extension $\textbf{L}.t$ of \textbf{L} be the smallest normal modal logic in the language $\mathcal{L}(\sig,\ind \cup\{i^{-1}\mid i\in \ind\})$ such that $\textbf{L}.t$ extends \textbf{L} and for every $i\in\ind$, $\textbf{L}.t$ contains the axioms $p\rightarrow \Box_{i}\Diamond_{i^{-1}}p$ and $p\rightarrow \Box_{i^{-1}}\Diamond_i p$. The following result is straightforward to prove (see \citealt[Corollary 3.25]{Litak2005b}).

\begin{lemma}[Minimal Tense Extensions and $\mathcal{T}$-Completeness]\label{TensEx} For any normal modal logic $\textbf{L}$, the following are equivalent:
\begin{enumerate}
\item\label{TensEx1} \textbf{L} is sound and complete with respect to a class of $\mathcal{T}$-BAOs;
\item\label{TensEx2} the minimal tense extension $\textbf{L}.t$ of \textbf{L} is a conservative extension of \textbf{L}.
\end{enumerate}
\end{lemma}

Characterizations of $\mathcal{V}$-completeness and $\mathcal{CV}$-completeness in terms of conservativity of extensions are open questions. The only new result we will report here is the analogue of Lemma \ref{AdjClose0} for $\mathcal{V}$-BAOs. Recall that the standard ``existence lemma'' for a normal modal logic \textbf{L} states that for any $i\in\ind$, $\Sigma\subseteq\mathcal{L}(\sig,\ind)$, and $\psi\in\mathcal{L}(\sig,\ind)$, if $\Sigma\nvdash_\mathbf{L}\Box_i\psi$, then there is a (maximally \textbf{L}-consistent) set $\Delta\subseteq\mathcal{L}(\sig,\ind)$ such that: 1. $\Delta\nvdash_\mathbf{L}\psi$, and 2. for all $\alpha\in\mathcal{L}(\sig,\ind)$, $\Sigma\vdash_\mathbf{L} \Box_i \alpha$ implies $\Delta\vdash_\mathbf{L}\alpha$. Consider the following finitary version of this property.

\begin{restatable}[Finite Existence Lemma]{definition}{FinEx}
A modal logic \textbf{L} satisfies the \textit{finite existence lemma} iff for every $i\in \ind$ and $\varphi, \psi\in\mathcal{L}(\sig,\ind)$, if $\nvdash_\mathbf{L} \varphi\rightarrow\Box_i\psi$, then there is a $\mathrm{g}_i^\mathbf{L}(\varphi,\psi)\in\mathcal{L}(\sig,\ind)$ such that:
\begin{enumerate}
\item $\nvdash_\mathbf{L}\mathrm{g}_i^\mathbf{L}(\varphi,\psi)\rightarrow\psi$;
\item for all $\alpha\in\mathcal{L}(\sig,\ind)$, $\vdash_\mathbf{L}\varphi\rightarrow \Box_i\alpha$ implies $\vdash_\mathbf{L}\mathrm{g}_i^\mathbf{L}(\varphi,\psi)\rightarrow\alpha$. \hfill $\triangleleft$
\end{enumerate}
\end{restatable}

Clearly \textbf{L} satisfies the finite existence lemma if it has internal adjointness. Conversely, arguments showing a lack of internal adjointness often extend to show a lack of the finite existence lemma. For example, the arguments in \citealt{Holliday2014} showing that the logic \textbf{K5} and its extensions with the D and/or 4 axioms lack internal adjointness can be easily adapted to show that these logics lack the finite existence lemma.\footnote{For negative results on internal adjointness and the finite existence lemma in the context of \textit{quantified} modal logic, see \citealt{HT2016a}.}

The following is the analogue of Lemma \ref{AdjClose0} for complete additivity.

\begin{lemma}[Finite Existence Lemma and Complete Additivity]\label{FinExLem} A consistent normal modal logic \textbf{L} satisfies the finite existence lemma iff the Lindenbaum algebra of \textbf{L} is a $\mathcal{V}$-BAO.
\end{lemma}

\begin{proof} To show that the Lindenbaum algebra $\mathbb{A}^\mathbf{L}$ of $\mathbf{L}$ is a $\mathcal{V}$-BAO, by Lemma \ref{VandH}.\ref{HtoV} it suffices to show that it is an $\mathcal{R}$-BAO as in Definition \ref{plentiful}, i.e., that for all $x,y\in A$, if $x\meet \blacklozenge_i y\not=\bot$, then there is a $y'\in A$ such that $\bot\not=y'\leq y$ and $xR_iy'$, where $xR_iy'$ iff for all $y''\in A$, if $\bot\not=y''\leq y'$, then $x\meet\blacklozenge_i y''\not=\bot$. So suppose that for some $x,y\in A$, $x\meet \blacklozenge_i y \not = \bot$. By definition of $\mathbb{A}^\mathbf{L}$ (Definition \ref{LinAlg}), there are $\varphi,\psi\in\mathcal{L}(\sig,\ind)$ such that $x=[\varphi]_\mathbf{L}$ and $y=[\psi]_\mathbf{L}$. (We will henceforth omit the subscripts for the equivalence classes.) Thus, $x\meet \blacklozenge_i y \not = \bot$ implies $\nvdash_\mathbf{L} \varphi\rightarrow \Box_i\neg\psi$.  Assuming that \textbf{L} satisfies the finite existence lemma, it follows that there is a $\chi\in\mathcal{L}(\sig,\ind)$ such that (i) $\nvdash_\mathbf{L} \chi\rightarrow \neg\psi$ and (ii) for all $\alpha\in\mathcal{L}(\sig,\ind)$, if $\vdash_\mathbf{L} \varphi\rightarrow \Box_i\alpha$, then $\vdash_\mathbf{L} \chi\rightarrow\alpha$. By (i), $\psi ' =\psi\wedge\chi$ is \textbf{L}-consistent, so  $\bot\not= [\psi']\leq [\psi]$. Now we claim that (ii) implies $[\varphi]R_i [\psi']$. Suppose not, so there is a $[\psi'']\in A$ such that $\bot\not= [\psi'']\leq [\psi']$ and $[\varphi]\meet\blacklozenge_i [\psi'']=[\varphi]\meet [\Diamond_i\psi'']=\bot$. It follows that $\vdash_\mathbf{L}\varphi\rightarrow\Box_i\neg\psi''$, which with (ii) implies $\vdash_\mathbf{L}\chi\rightarrow \neg\psi''$. But then since $\psi'=\psi\wedge\chi$, $\vdash_\mathbf{L}\psi'\rightarrow \neg\psi''$, which contradicts $\bot\not=[\psi'']\leq [\psi']$. Thus, $[\varphi] R_i [\psi']$, which shows that $\mathbb{A}^\mathbf{L}$ is an $\mathcal{R}$-BAO.

Now suppose \textbf{L} does not satisfy the finite existence lemma, so there are $\varphi,\psi\in\mathcal{L}(\sig,\ind)$ such that (a) $\nvdash_\mathbf{L}\varphi\rightarrow\Box_i\psi$ but (b) for every $\chi\in\mathcal{L}(\sig,\ind)$ such that $\nvdash_\mathbf{L}\chi\rightarrow\psi$, there is an $\alpha\in\mathcal{L}(\sig,\ind)$ such that $\vdash_\mathbf{L}\varphi\rightarrow\Box_i\alpha$ but $\nvdash_\mathbf{L} \chi\rightarrow\alpha$. By (a),  $[\varphi]\meet \blacklozenge_i [\neg\psi]=[\varphi]\meet [\Diamond_i\neg\psi]\not=\bot$. Now consider any $\chi\in\mathcal{L}(\sig,\ind)$ such that $\bot\not= [\chi]\leq[\neg\psi]$, so $\nvdash_\mathbf{L} \chi\rightarrow\psi$. Then by (b), there is an $\alpha\in\mathcal{L}(\sig,\ind)$ such that $\vdash_\mathbf{L}\varphi\rightarrow\Box_i\alpha$ but $\nvdash_\mathbf{L} \chi\rightarrow\alpha$. Where $\chi'=\chi\wedge\neg\alpha$, $\nvdash_\mathbf{L} \chi\rightarrow\alpha$ implies $\bot\not=[\chi']\leq [\chi]$, and  $\vdash_\mathbf{L}\varphi\rightarrow\Box_i\alpha$ implies $[\varphi]\meet\blacklozenge_i [\chi']=\bot$. Thus, $\mathbb{A}^\mathbf{L}$ is not an $\mathcal{R}$-BAO, so by Lemma \ref{VandH}.\ref{VtoH}, it is not a $\mathcal{V}$-BAO.\end{proof}

From Lemma \ref{FinExLem} and Theorem \ref{VtoPossFrames} we obtain the analogue of Corollary \ref{GenAdj0} for the finite existence lemma in place of internal adjointness.

\begin{corollary}[General Completeness for Logics with the Finite Existence Lemma] Every consistent normal modal logic satisfying the finite existence lemma is sound and complete with respect to a principal possibility frame.
\end{corollary}

\subsection{Completeness for Atomless Possibility Frames}\label{AtomlessFull}

As noted in \S~\ref{intro}, a topic emphasized in previous work on possibility semantics \citep{Humberstone1981,Holliday2014} is the completeness of some modal logics with respect to classes of \textit{atomless} possibility frames, which we define as possibility frames whose poset $\langle S,\sqsubseteq\rangle$ is such that for every $x\in S$, there is an $x'\in S$ such that $x' \sqsubset x$, i.e., $x'\sqsubseteq x$ but $x\not\sqsubseteq x'$. Compare this with the notion of an atomless BAO, which is a BAO such that for every $x\in A\setminus\{\bot\}$, there is an $x'\in A\setminus\{\bot\}$ such that $x'<x$, i.e., $x'\leq x$ but $x\not\leq x'$. Note that if the set $\sig$ of propositional variables is infinite, then the Lindenbaum algebra $\mathbb{A}^\mathbf{L}$ of any normal modal logic $\mathbf{L}\subseteq\mathcal{L}(\sig,\ind)$ (Definition \ref{LinAlg}) is atomless. We assume in this section that $\sig$ is infinite.

Using the fact that for any atomless BAO $\mathbb{A}$, its principal frame $\mathbb{A}_\rela$ as in Definition \ref{Hposs} is atomless, we obtain the following result concerning completeness with respect to atomless \textit{principal} possibility frames.

\begin{theorem}[Completeness for Atomless Principal Possibility Frames]\label{CompPrincAtomless} For any consistent normal modal logic $\mathbf{L}$: 
\begin{enumerate}
\item\label{CompPrincAtomless1} if the Lindenbaum algebra $\mathbb{A}^\mathbf{L}$ of $\mathbf{L}$ is a $\mathcal{V}$-BAO (equivalently, if $\mathbf{L}$ satisfies the finite existence lemma), then $\mathbf{L}$ is sound and complete with respect to an atomless principal possibility frame, viz., $(\mathbb{A}^\mathbf{L})_\rela$; 
\item\label{CompPrincAtomless2} if \textbf{L} is sound and complete with respect to a class of $\mathcal{T}$-BAOs (equivalently, if the minimal tense extension $\mathbf{L}.t$ of $\mathbf{L}$ is a conservative extension), then \textbf{L} is sound and complete with respect to an atomless and functional principal possibility frame.
\end{enumerate}
\end{theorem}

\begin{proof} Part \ref{CompPrincAtomless1} follows from Theorem \ref{VtoPossFrames} and the fact that the Lindenbaum algebra is atomless.

For part \ref{CompPrincAtomless2}, it suffices to show that \textbf{L} is sound and complete with respect to a principal possibility frame $\mathcal{F}$ that is atomless and \textit{quasi}-functional, for then \textbf{L} is also sound and complete with respect to the \textit{functionalization} of $\mathcal{F}$ (Proposition \ref{QtoF}), which is principal and atomless if $\mathcal{F}$ is. Where $\mathbb{A}^{\mathbf{L}.t}$ is the Lindenbaum algebra of $\mathbf{L}.t$, which is a $\mathcal{T}$-BAO, let $\mathbb{B}$ be the reduct of $\mathbb{A}^{\mathbf{L}.t}$ to the signature of $\mathbf{L}$, i.e., the result of dropping from $\mathbb{A}^{\mathbf{L}.t}$ the converse operators $\blacksquare_{i^{-1}}$.  $\mathbb{B}$ is still a $\mathcal{T}$-BAO, because the functions that must exist for a BAO to be a $\mathcal{T}$-BAO need not be included among the operations of the BAO; and since $\mathbb{A}^{\mathbf{L}.t}$ is an atomless BAO, so is $\mathbb{B}$. Now for any $\varphi\in\mathcal{L}(\sig,\ind)$, we have the following equivalences: $\varphi$ is consistent in $\mathbf{L}$ iff $\varphi$ is consistent in $\mathbf{L}.t$ (because $\mathbf{L}.t$ is a conservative extension of $\mathbf{L}$) iff $\varphi$ is satisfiable in  $\mathbb{A}^{\mathbf{L}.t}$ iff $\varphi$ is satisfiable in $\mathbb{B}$ (because $\varphi\in\mathcal{L}(\sig,\ind)$ and $\mathbb{B}$ is the reduct of $\mathbb{A}^{\mathbf{L}.t}$ to this signature) iff  $\varphi$ is satisfiable in $\mathbb{B}_\rela$ (by Theorem \ref{VtoPossFrames}.\ref{VtoPoss5}). Thus, $\mathbf{L}$ is sound and complete with respect to $\mathbb{B}_\rela$. Then since $\mathbb{B}$ is atomless, so is $\mathbb{B}_\rela$, and since $\mathbb{B}$ is a $\mathcal{T}$-BAO, $\mathbb{B}_\rela$ is a quasi-functional principal frame by Theorem \ref{VtoPossFrames}.\ref{VtoPoss2}, so we are done.
\end{proof}

In addition to proving completeness with respect to atomless \textit{principal} possibility frames, we would like to do so with respect to atomless \textit{full} possibility frames.  We will first prove that every normal modal logic obtained from \textbf{K} by adding a set of axioms of the Lemmon-Scott form $\Diamond_\alpha\Box_\beta  p\rightarrow \Box_\delta \Diamond_\gamma  p$ from \S~\ref{LemmScottCorr} is sound and complete with respect to an \textit{atomless full} possibility frame. To do so, we first show that if a $\mathcal{T}$-BAO $\mathbb{A}$ validates a Lemmon-Scott axiom, then its full frame $\mathbb{A}_\fullV$ as in Definition \ref{Hposs} has the property corresponding to the axiom according to Proposition \ref{Gklmn}. 

\begin{lemma}[Full Frames of $\mathcal{T}$-BAOs and Lemmon-Scott Correspondence]\label{T-LemScott} For any $\mathcal{T}$-BAO $\mathbb{A}$ and sequences $\alpha$, $\beta$, $\delta$, and $\gamma$ of indices from $\ind$, if for all $x\in\mathbb{A}$, $\blacklozenge_\alpha\blacksquare_\beta  x\leq \blacksquare_\delta \blacklozenge_\gamma  x$, then $\mathbb{A}_\fullV$ satisfies:
\[\forall{x}\forall{y}\big({x}R_\delta {y}\rightarrow \exists {x'}\sqsubseteq{x} \;\forall {z}( {x'}R_\alpha {z} \rightarrow\exists {u} ({y} R_\gamma  {u}\wedge {z} R_\beta  {u}))\big).\]
\end{lemma}
 
\begin{proof} Let $f_i$ be the left adjoint of $\blacksquare_i$ that is given since $\mathbb{A}$ is a $\mathcal{T}$-BAO. For a sequence $\sigma=\langle i_1,\dots,i_n\rangle$ of indices from $\ind$, let $f_\sigma(x)=f_{i_n}(f_{i_{n-1}}(\dots f_{i_1}(x)\dots))$, and if $\sigma$ is empty, let $f_\sigma(x)=x$. Recall that by Lemma \ref{RelT}, $\mathbb{A}_\fullV$ is such that $xR_iy$ iff $y\sqsubseteq f_i(x)$. Also note that since $f_i$ is monotone with respect to the order $\leq$ of $\mathbb{A}$ (for since $f_i(x)\leq f_i(x)$, we have $x\leq \blacksquare_i f_i(x)$ by adjointness, so $x'\leq x$ implies $x'\leq \blacksquare_i f_i(x)$ and hence $f_i(x')\leq f_i(x)$ by adjointness), in  $\mathbb{A}_\fullV$ we have that if $x'\sqsubseteq x$ and $R_i(x')\not=\emptyset$, then $f_i(x')\sqsubseteq f_i(x)$. It follows that $xR_\sigma y$ as in \S~\ref{LemmScottCorr} iff $y\sqsubseteq f_\sigma(x)$, so the property in the statement of the lemma is equivalent to:
\[\forall x\forall y\sqsubseteq f_\delta  ({x})\,\exists {x'}\sqsubseteq{x} \,\forall {z}\sqsubseteq f_\alpha({x'}) \,\exists u( u\sqsubseteq f_\gamma ({y})\wedge u\sqsubseteq f_\beta ({z})).\] 
Recall that the poset $\langle S,\sqsubseteq\rangle$ of $\mathbb{A}_\fullV$ is the Boolean lattice $\langle A,\leq\rangle$ of $\mathbb{A}$ restricted to $A\setminus\{\bot\}$. So for $y\in\mathbb{A}_\fullV$, $y\not=\bot$. Now suppose ${y}\sqsubseteq f_\delta  ({x})$. Then ${x}\meet \blacklozenge_\delta  {y}\not=\inc$, for otherwise ${x}\sqsubseteq \blacksquare_\delta \mathord{-}{y}$, which implies $f_\delta  ({x})\sqsubseteq -{y}$, which contradicts ${y}\sqsubseteq f_\delta  ({x})$ given $y\not=\bot$. Let ${x'}={x}\meet \blacklozenge_\delta  {y}$. Now consider some $ {z}\sqsubseteq f_\alpha({x'})$, so ${z}\not=\inc$. For reductio, suppose $f_\gamma ({y})\meet f_\beta ({z})=\inc$, so $f_\gamma ({y})\sqsubseteq - f_\beta ({z})$, which implies ${y}\sqsubseteq \blacksquare_\gamma  \mathord{-} f_\beta ({z})$, which in turn implies  $\blacklozenge_\delta {y}\sqsubseteq \blacklozenge_\delta \blacksquare_\gamma  \mathord{-} f_\beta ({z})$. By our choice of ${x'}$, ${x'}\sqsubseteq \blacklozenge_\delta  {y}$, and by the assumption of the lemma, $\blacklozenge_\delta \blacksquare_\gamma  \mathord{-} f_\beta ({z})\sqsubseteq \blacksquare_\alpha \blacklozenge_\beta  \mathord{-} f_\beta ({z})$, so ${x'}\sqsubseteq \blacksquare_\alpha \blacklozenge_\beta  \mathord{-} f_\beta ({z})$, which implies $f_\alpha({x'})\sqsubseteq \blacklozenge_\beta  \mathord{-} f_\beta ({z})$. Then by our choice of ${z}$, ${z}\sqsubseteq \blacklozenge_\beta  \mathord{-} f_\beta ({z})$, so ${z}\sqsubseteq \mathord{-} \blacksquare_\beta  f_\beta ({z})$. But then since ${z}\sqsubseteq \blacksquare_\beta  f_\beta ({z})$, we have ${z}=\inc$, which is inconsistent with our choice of ${z}$. Hence $f_\gamma ({y})\meet f_\beta ({z})\not=\inc$, so taking $u=f_\gamma ({y})\meet f_\beta ({z})$ completes the proof.\end{proof}

We can now prove the promised completeness result for atomless full possibility frames. 

\begin{theorem}[Lemmon-Scott Completeness for Atomless Full Possibility Frames]\label{CompAtomless} Where \textbf{L} is the least normal modal logic extending \textbf{K} with a set of axioms of the Lemmon-Scott form $\Diamond_\alpha\Box_\beta  p\rightarrow \Box_\delta \Diamond_\gamma  p$, \textbf{L} is sound and complete with respect to an atomless and functional full possibility frame.
\end{theorem}

\begin{proof} As in the proof of Theorem \ref{CompPrincAtomless}, it suffices to show that \textbf{L} is sound and complete with respect to an atomless and \textit{quasi}-functional full possibility frame. Any \textbf{L} as in the statement of the current theorem is Kripke-frame complete \citep{Lemmon1977}, so it is complete with respect to a $\mathcal{CAV}$-BAO and hence a $\mathcal{T}$-BAO, so the minimal tense extension $\mathbf{L}.t$ of \textbf{L} is a conservative extension of $\mathbf{L}$ by Lemma \ref{TensEx}. Where $\mathbb{A}^{\mathbf{L}.t}$ is the Lindenbaum algebra of $\mathbf{L}.t$, which is a $\mathcal{T}$-BAO, let $\mathbb{B}$ be the reduct of $\mathbb{A}^{\mathbf{L}.t}$ to the signature of $\mathbf{L}$, as in the proof of Theorem \ref{CompPrincAtomless}, so $\mathbb{B}$ is still a $\mathcal{T}$-BAO and atomless since $\mathbb{A}^{\mathbf{L}.t}$ is. Then as in the proof of Theorem \ref{CompPrincAtomless}, for any $\varphi\in\mathcal{L}(\sig,\ind)$, $\varphi$ is consistent in $\mathbf{L}$ iff  $\varphi$ is satisfiable in $\mathbb{B}_\rela$. Now if $\varphi$ is satisfiable in the principal frame $\mathbb{B}_\rela$, then $\varphi$ is satisfiable in the full frame $\mathbb{B}_\fullV$, so we have that $\mathbf{L}$ is complete with respect to $\mathbb{B}_\fullV$. It only remains to show that $\mathbf{L}$ is \textit{sound} with respect to $\mathbb{B}_\fullV$. For each axiom $\Diamond_\alpha\Box_\beta  p\rightarrow \Box_\delta \Diamond_\gamma  p$ added to \textbf{K} to obtain \textbf{L}, $\mathbb{A}^{\mathbf{L}.t}$ and hence $\mathbb{B}$ is such that for all $x\in\mathbb{B}$, $\blacklozenge_\alpha\blacksquare_\beta  x\leq \blacksquare_\delta \blacklozenge_\gamma  x$. Then by Lemma \ref{T-LemScott}, $\mathbb{B}_\fullV$ satisfies the properties that correspond to the axioms of \textbf{L} in Proposition \ref{Gklmn}. Thus, $\mathbf{L}$ is sound with respect to $\mathbb{B}_\fullV$. Since $\mathbb{B}_\fullV$ is atomless, quasi-functional (Theorem \ref{VtoPossFrames}.\ref{VtoPoss2}), and full, we are done. \end{proof} 

In fact, we have the resources to generalize Theorem \ref{CompAtomless} to any consistent \textit{Sahlqvist logic}, i.e., any least normal extension of $\mathbf{K}$ with some set of Sahlqvist axioms (see \citealt[\S~3.6]{Blackburn2001}). The key resources are Theorem \ref{MonkComp}, which showed that the underlying BAO $(\mathbb{A}_\fullV)^\under$ of the full frame $\mathbb{A}_\fullV$ of a $\mathcal{V}$-BAO $\mathbb{A}$ is isomorphic to the Monk completion of $\mathbb{A}$, plus Givant and Venema's \citeyearpar[Cor.~34(ii)]{Givant1999} result that all Sahlqvist equations are preserved in going from a $\mathcal{T}$-BAO to its Monk completion.

\begin{theorem}[Sahlqvist Completeness for Atomless Full Possibility Frames]\label{SahlqComp} Every consistent Sahlqvist logic is sound and complete with respect to an atomless and functional full possibility frame.
\end{theorem}

\begin{proof} The proof uses the same argument as in the proof of Theorem \ref{CompAtomless}, only now we must show that every Sahlqvist formula valid over $\mathbb{B}$ is also valid over $\mathbb{B}_\fullV$. Since $\mathbb{B}$ is a $\mathcal{T}$-BAO, the cited result of Givant and Venema \citeyearpar{Givant1999} gives us that every Sahlqvist formula valid over $\mathbb{B}$ is valid over its Monk completion. By Theorem \ref{MonkComp}, the Monk completion of $\mathbb{B}$ is isomorphic to $(\mathbb{B}_\fullV)^\under$, which has the same logic as $\mathbb{B}_\fullV$ by Theorem \ref{PtoB}.\ref{PtoB6}. Thus, the validity of Sahlqvist formulas is indeed preserved from $\mathbb{B}$ to $\mathbb{B}_\fullV$.\end{proof}

Theorem \ref{SahlqComp} represents a quite general realization of the goal from \citealt{Humberstone1981} of proving completeness of modal logics with respect to full possibility frames that contain no worlds.

\subsection{Completeness for Canonical Possibility Frames}\label{Canonical}

Let us finally treat the topic of \textit{canonical} possibility frames and models for normal modal logics. In Kripke semantics, the set of worlds in the canonical frame for a logic \textbf{L} is the set of all \textit{maximally} \textbf{L}-consistent sets of formulas. From an algebraic perspective, these worlds arise from \textit{ultra}filters in the Lindenbaum algebra for \textbf{L}. By contrast, in possibility semantics, we can take the set of possibilities in the canonical frame for \textbf{L} to be the set of all \textbf{L}-consistent and \textbf{L}-deductively closed sets of formulas. From an algebraic perspective, these possibilities arise from proper filters in the Lindenbaum algebra for \textbf{L}. We will define the canonical frame for \textbf{L} in this algebraic way, based on the notions of the filter frame $\mathbb{A}_\ff$ and general filter frame $\mathbb{A}_\gff$ of a BAO from Definition \ref{FiltEx}. The more syntactic definition is then straightforward to work out.

\begin{definition}[Canonical Frames, Models, Formulas, and Logics]\label{CanFrame&Mod} Given a consistent normal modal logic \textbf{L}, where $\mathbb{A}^\mathbf{L}$ is the Lindenbaum algebra for \textbf{L}, and $\mathbb{M}^\mathbf{L}=\langle\mathbb{A}^\mathbf{L},\theta\rangle$ is the algebraic model such that for all $p\in \sig$, $\theta(p)=[p]_\mathbf{L}$ (see Definition \ref{LinAlg}):
\begin{enumerate}
\item the \textit{canonical \textnormal{(}general\textnormal{)} possibility frame} for \textbf{L} is the frame $\mathcal{G}^\mathbf{L}=(\mathbb{A}^\mathbf{L})_{\gff}$;
\item the \textit{canonical full possibility frame} for \textbf{L} is the frame $\mathcal{F}^\mathbf{L}=(\mathbb{A}^\mathbf{L})_{\ff}$;
\item the \textit{canonical possibility model} for \textbf{L} is the model $\mathcal{M}^\mathbf{L}=(\mathbb{M}^\mathbf{L})_{\gff}$. 
\end{enumerate}
A formula $\varphi\in\mathcal{L}(\sig,\ind)$ is \textit{filter-canonical} iff for any consistent normal modal logic $\mathbf{L}\subseteq \mathcal{L}(\sig,\ind)$, $\vdash_\mathbf{L}\varphi$ implies $\mathcal{F}^\mathbf{L}\Vdash\varphi$.

A consistent normal modal logic \textbf{L} is \textit{filter-canonical} iff it is sound with respect to $\mathcal{F}^\mathbf{L}$.\footnote{In the context of Kripke semantics, the term `canonical' has been used to mean several things when applied to a normal modal logic \textbf{L} (with a countably infinite set of  propositional variables): (a) that \textbf{L} is sound with respect to its canonical Kripke frame \citep{Goldblatt1976b}; (b) that for any infinite cardinal $\kappa$, \textbf{L} is sound with respect to the canonical Kripke frame for $\textbf{L}_\kappa$, the conservative extension of \textbf{L} to the modal language with $\kappa$-many propositional variables \citep{Fine1975c}; and (c) that if \textbf{L} is sound with respect to a descriptive frame, then \textbf{L} is sound with respect to the associated Kripke frame \citep{Benthem1979b} (called `d-persistent' in \citealt{Goldblatt1974}). While (b) and (c) are equivalent \citep[p. 294]{Sambin1988}, and (b)/(c) obviously implies (a), it is unknown whether (a) implies (b)/(c). Note that our notion of filter-canonical is analogous to (a).} \hfill$\triangleleft$
\end{definition}

From Theorem \ref{SatGenFil} and the fact that every normal modal logic is sound and complete with respect to its Lindenbaum algebra (Theorem \ref{AdAlg}), we have the following result for $\mathcal{G}^\mathbf{L}$.
 
\begin{corollary}[General Soundness and Completeness]\label{GenCompCan} For any consistent normal modal logic \textbf{L}, \textbf{L} is sound and complete with respect to its canonical (general) possibility frame $\mathcal{G}^\mathbf{L}$, and \textbf{L} is complete with respect to its canonical full possibility frame $\mathcal{F}^\mathbf{L}$.
\end{corollary}
\noindent Recall that a general filter frame such as $\mathcal{G}^\mathbf{L}$ is \textit{filter-descriptive} (Proposition \ref{GenAb}), so Corollary \ref{GenCompCan} is another expression of general soundness and completeness with respect to filter-descriptive possibility frames.

Of course, we cannot have such a general soundness result for the canonical \textit{full} possibility frame $\mathcal{F}^\mathbf{L}$---not every consistent normal modal logic is \textit{filter-canonical}---because we know from Corollary \ref{TransferPossAlg} that full possibility frames are no more general than $\mathcal{CV}$-BAOs. There are several routes to proving that certain logics are filter-canonical. First, if we assume the ultrafilter axiom, then the logic of $\mathcal{F}^\mathbf{L}$ is exactly the logic of the canonical \textit{Kripke} frame for $\mathbf{L}$, or equivalently, the \textit{ultra}filter frame of $\mathbb{A}^\mathbf{L}$. For these frames have the same logics as their underlying BAOs (Theorem \ref{PtoB}.\ref{PtoB6}), and as we observed in \S~\ref{MacNeille}, the underlying BAOs of the filter and ultrafilter frames are isomorphic under the assumption of the ultrafilter axiom.  Thus, assuming that axiom, a logic is filter-canonical iff it is canonical in the traditional sense of Kripke semantics.

A second, more direct route to proving filter-canonicity, which does not require the ultrafilter axiom, takes advantage of correspondence results for possibility semantics. As an example, we will prove that any normal modal logic obtained from \textbf{K} by the addition of Lemmon-Scott axioms as in \S~\ref{LemmScottCorr} is filter-canonical.

\begin{lemma}[Filter Frames and Lemmon-Scott Correspondence]\label{FFFO} For any BAO $\mathbb{A}$ and sequences $\alpha$, $\beta$, $\delta$, and $\gamma$ of indices from $\ind$, if for all $x\in\mathbb{A}$, $\blacklozenge_\alpha\blacksquare_\beta  x\leq \blacksquare_\delta \blacklozenge_\gamma  x$, then $\mathbb{A}_{\gff}$ and $\mathbb{A}_{\ff}$ satisfy:
\[\forall{X}\,\forall{Y}\big({X}R_\delta {Y}\Rightarrow \exists {X'}\sqsubseteq{X} \;\forall {Z}( {X'}R_\alpha {Z} \Rightarrow\exists {U} ({Y} R_\gamma  {U}\wedge {Z} R_\beta  {U}))\big).\]
\end{lemma}

\begin{proof} It is helpful to reformulate the assumption as: $\blacklozenge_\delta \blacksquare_\gamma  x\leq \blacksquare_\alpha\blacklozenge_\beta  x$ for all $x\in\mathbb{A}$. Suppose ${X}R_\delta  {Y}$, so by Definition \ref{FiltEx}.\ref{FiltEx3}, we have that for all $x\in \mathbb{A}$, $\blacksquare_\delta  x\in{X}$ implies $x\in{Y}$. Where
\begin{equation}\mathrm{X}'={X}\cup \{\blacklozenge_\delta  y\mid y\in{Y}\},\label{X'def}\end{equation}
by the same form of argument as in the proof of Theorem \ref{BAOtoGenPos}.\ref{IsRelGen}, ${X'}=[\mathrm{X}')$ is a proper filter, and $X'\supseteq X$ gives us $X'\sqsubseteq X$. Now consider a proper filter ${Z}$ such that ${X'}R_\alpha{Z}$. Where
\begin{equation}\mathrm{U}=\{y\mid \blacksquare_\gamma  y\in{Y}\}\cup\{z\mid \blacksquare_\beta  z\in{Z}\},\label{Udef2}\end{equation}
suppose for reductio that $[\mathrm{U})$ is not a proper filter, i.e., $\inc\in[\mathrm{U})$. Then by (\ref{Udef2}) and Fact \ref{FiltGen}, there are $\blacksquare_\gamma y\in{Y}$ and $\blacksquare_\beta  z\in{Z}$ such that 
$y\meet z\leq \inc$, so $y\leq \mathord{-}z$, which implies $\blacksquare_\gamma  y\leq \blacksquare_\gamma  \mathord{-}z$ by the monotonicity of $\blacksquare_\gamma$. Then since ${Y}$ is a filter, $\blacksquare_\gamma  y\in{Y}$ implies $\blacksquare_\gamma  \mathord{-} z\in{Y}$. It follows by (\ref{X'def}) that $\blacklozenge_\delta \blacksquare_\gamma  \mathord{-}z\in {X'}$, which with our initial assumption and the fact that ${X'}$ is a filter implies $\blacksquare_\alpha\blacklozenge_\beta  \mathord{-} z\in{X'}$. Then since ${X'}R_\alpha{Z}$, we have $\blacklozenge_\beta  \mathord{-} z\in{Z}$, which contradicts the fact that $\blacksquare_\beta  z\in{Z}$. Thus, ${U}=[\mathrm{U})$ is a proper filter, and by (\ref{Udef2}), we have ${Y}R_\gamma  {U}$ and ${Z}R_\beta {U}$. 
\end{proof} 

From Lemma \ref{FFFO}, we immediately obtain the following.

\begin{corollary}[Lemmon-Scott Filter-Canonicity]\label{LSCan} For any consistent normal modal logic \textbf{L} and sequences $\alpha$, $\beta$, $\delta$, and $\gamma$ of indices from $\ind$, if $\vdash_\mathbf{L}\Diamond_\alpha\Box_\beta  p\rightarrow \Box_\delta \Diamond_\gamma  p$, then the canonical full possibility frame $\mathcal{F}^\mathbf{L}$ satisfies the first-order correspondent of $\Diamond_\alpha\Box_\beta  p\rightarrow \Box_\delta \Diamond_\gamma  p$ as in Proposition \ref{Gklmn}, so $\mathcal{F}^\mathbf{L}\Vdash\Diamond_\alpha\Box_\beta  p\rightarrow \Box_\delta \Diamond_\gamma  p$.
\end{corollary}

Together with Theorem \ref{BAOtoGenPos}, Corollary \ref{LSCan} gives us the following general completeness result. 

\begin{theorem}[Lemmon-Scott Soundness and Completeness for $\mathcal{F}^\mathbf{L}$]\label{LSagain} Where \textbf{L} is the least normal modal logic extending \textbf{K} with a set of axioms of the Lemmon-Scott form $\Diamond_\alpha\Box_\beta  p\rightarrow \Box_\delta \Diamond_\gamma  p$, \textbf{L} is sound and complete with respect to its canonical full possibility frame $\mathcal{F}^\mathbf{L}$.
\end{theorem} 

A third route to proving filter-canonicity, which also does not require the ultrafilter axiom, uses Theorem \ref{CanExtProp}. It follows from Theorem \ref{CanExtProp} that $(\mathcal{F}^\mathbf{L})^\under=((\mathbb{A}^\mathbf{L})_\ff)^\under$ is a \textit{constructive} canonical extension of the Lindenbaum algebra $\mathbb{A}^\mathbf{L}$ of $\mathbf{L}$. Conradie and Palmigiano \citeyearpar{Conradie2016} show in a constructive meta-theory that all \textit{inductive formulas}, which include Sahlqvist formulas as a special case, are preserved under constructive canonical extensions (cf.~\citealt{Ghilardi1997}). Thus, if $\mathbb{A}^\mathbf{L}$ validates some inductive formulas, then so does $(\mathcal{F}^\mathbf{L})^\under$ and hence $\mathcal{F}^\mathbf{L}$, which yields the following.

\begin{theorem}[Inductive Filter-Canonicity]\label{InductiveFilterCan} All inductive formulas are filter-canonical. Hence every consistent normal modal logic $\mathbf{L}$ axiomatized by inductive formulas is sound and complete with respect to its canonical full possibility frame $\mathcal{F}^\mathbf{L}$.
\end{theorem}

The results of this section, \S~\ref{AtomlessFull}, and \S~\ref{CompFull} concerning completeness with respect to \textit{full} possibility frames lead to questions about the ``persistence'' of modal formulas (see, e.g., \citealt[\S5.6]{Blackburn2001}), not in the sense of \textit{persistence} that we have used since Remark~\ref{Persp1}, but in the following sense: a formula $\varphi$ is persistent with respect to a class $\mathsf{F}$ of possibility frames that have associated full frames as in Definition \ref{Foundation} iff whenever $\varphi$ is valid over a frame $\mathcal{F}\in\mathsf{F}$, then $\varphi$ is also valid over the associated full frame $(\mathcal{F}_\found)^\full$. For example, we can consider persistence of $\varphi$ with respect to \textit{filter-descriptive} possibility frames as in \S~\ref{Fdes}, which implies that $\varphi$ is filter-canonical as in Definition \ref{CanFrame&Mod}, since $\mathcal{G}^\mathbf{L}$ is filter-descriptive and its associated full frame $((\mathcal{G}^\mathbf{L})_\found)^\full$ is $\mathcal{F}^\mathbf{L}$. Or we can consider persistence with respect to \textit{principal} possibility frames, which with Theorem \ref{MonkComp} and our duality theory implies that $\varphi$ is preserved under Monk completions of $\mathcal{V}$-BAOs. We leave it to future work to study various notions of persistence in the context of possibility semantics.

\section{Conclusion}\label{Conclusion}

It is remarkable that as a semantics for modal logic originally motivated mainly on philosophical and conceptual grounds \citep{Humberstone1981}, possibility semantics also turns out to be natural from a mathematical perspective, e.g., with its topological and set-theoretical connections, and so fruitful for pure modal logic, e.g., with its new categorical connections with classes of modal algebras. Something similar could be argued for standard possible world semantics, of which possibility semantics is a generalization. So it is remarkable that the confluence of the philosophical-conceptual and the mathematical continues. 

We will conclude our study of possibility semantics in this paper by looking to the past and present of related work in \S~\ref{RelatedWork} and then looking to the future of open problems in \S~\ref{OpenProb}.

\subsection{Related Work}\label{RelatedWork}

Below we organize our pointers to related work by topic, indicated by the italicized headings. We will not repeat here discussion of the related work in modal logic and set theory to which we have already referred. \\

\noindent\textit{Possibility Semantics for Non-modal Logic}. An early source of what is essentially possibility semantics for classical first-order logic is \citealt{Fine1975}. Here we focus on the propositional version of Fine's semantics, for comparison with our framework. One can think of Fine's propositional models as tuples $M=\langle S,\sqsubseteq,V\rangle$ where $\langle S,\sqsubseteq\rangle$ is a poset and $V\colon \sig \times S\pfun \{0,1\}$ is a partial function satisfying the following conditions (see \citealt{Fine1975}, p.~272 and pp.~278-9):
\begin{itemize}
\item \textit{stability}: if $V(p,x)$ is defined and $x'\sqsubseteq x$, then $V(p,x')$ is defined and $V(p,x')=V(p,x)$;
\item \textit{resolution}: if $V(p,x)$ is undefined, then there are $y\sqsubseteq x$ and $z\sqsubseteq x$ such that $V(p,y)=1$ and $V(p,z)=0$.
\end{itemize}
Note that \textit{stability} is the analogue of our \textit{persistence}, and \textit{resolution} is the analogue of our \textit{refinability}. 

 Fine (p.~279) defines truth and falsity relations $\vDash$ and $\rotatebox[origin=c]{180}{$\vDash$}$ as follows: 

\begin{itemize}
\item $M,x\vDash p$ iff $V(p,x)=1$;
\item $M,x\;\rotatebox[origin=c]{180}{$\vDash$}\; p$ iff $V(p,x)=0$;
\item $M,x\vDash \neg\varphi$ iff $M,x\;\rotatebox[origin=c]{180}{$\vDash$}\;\varphi$;
\item $M,x\;\rotatebox[origin=c]{180}{$\vDash$}\;\neg\varphi$ iff $M,x\vDash \varphi$;
\item $M,x\vDash \varphi\wedge\psi$ iff $M,x\vDash\varphi$ and $M,x\vDash\psi$;
\item $M,x\;\rotatebox[origin=c]{180}{$\vDash$}\; \varphi\wedge\psi$ iff $\forall x'\sqsubseteq x$ $\exists x''\sqsubseteq x'$: $M,x''\;\rotatebox[origin=c]{180}{$\vDash$}\;\varphi$ or $M,x''\;\rotatebox[origin=c]{180}{$\vDash$}\;\psi$.
\end{itemize}
The connectives $\vee$ and $\rightarrow$ are defined classically in terms of $\neg$ and $\wedge$. 

It is easy to see from the falsity clause for $\varphi\wedge\psi$, which matches our forcing clause for $\neg\varphi\vee\neg\psi$ (recall Fact \ref{ForcingDis}), that this setup is equivalent to  possibility semantics for propositional logic. Let us say that a \textit{propositional} possibility model is a tuple $\mathcal{M}=\langle S,\sqsubseteq,\pi\rangle$ where $\pi\colon \sig\to \mathrm{RO}(S,\sqsubseteq)$. Recall from \S~\ref{ClassicalFrames} that $\mathrm{RO}(S,\sqsubseteq)$ is the set of all $X\subseteq S$ that satisfy \textit{persistence} and \textit{refinability}, or equivalently, that are regular open sets in the downset topology on $\langle S,\sqsubseteq\rangle$. We can make the claimed equivalence precise as in Fact \ref{FineEq}, the proof of which is straightforward and left as an exercise.

\begin{fact}[Equivalence of Fine's Semantics and Possibility Semantics]\label{FineEq}
Given any Fine model $M=\langle S,\sqsubseteq,V\rangle$, define $M^p=\langle S,\sqsubseteq,\pi_V\rangle$ by $\pi_V(p)=\{x\in S\mid V(p,x)=1\}$. Then:
\begin{enumerate}[label=\arabic*.,ref=\arabic*]
\item $M^p$ is a propositional possibility model;
\item for any $x\in S$ and $\varphi\in\mathcal{L}(\sig,\emptyset)$, $M,x\vDash\varphi$ iff $M^p,x\Vdash \varphi$, 
\end{enumerate}
where $\Vdash$ is our forcing relation from Definition \ref{pmtruth1}. 

Conversely, given any propositional possibility model $\mathcal{M}=\langle S,\sqsubseteq,\pi\rangle$, define $\mathcal{M}^f=\langle S,\sqsubseteq, V_\pi\rangle$ by: $V_\pi(p,x)=1$ if $x\in \pi(p)$; $V_\pi(p,x)=0$ if $\forall x'\sqsubseteq x$ $x'\not\in \pi(p)$; and otherwise $V_\pi(p,x)$ is undefined. Then:
\begin{enumerate}[label=\arabic*.,ref=\arabic*,resume]
\item $\mathcal{M}^f$ is a Fine model;
\item for any $x\in S$ and $\varphi\in\mathcal{L}(\sig,\emptyset)$, $\mathcal{M},x\Vdash\varphi$ iff $\mathcal{M}^f,x\vDash \varphi$.
\end{enumerate}
\end{fact}
In addition to the conditions on Fine models given above, Fine proposes that every element in the poset $\langle S,\sqsubseteq\rangle$ should be refined by an endpoint, as in our Definition \ref{AtDef} for \textit{atomic} frames (see his condition of \textit{completeability} on p.~272). However, he also observes (p.~280) on the basis of a Cohen-style argument as in our Lemma \ref{classicality} that the assumption about endpoints is not necessary to obtain classical logic.

Another source of possibility semantics for classical first-order logic is  \citealt{Benthem1981}. Starting with Kripke models for intuitionistic first-order logic and then imposing the same kind of \textit{refinability} condition on valuations as in this paper (but there called \textit{cofinality}---recall footnote \ref{CofinalFootnote}), van Benthem observes that one obtains a semantics for classical first-order logic by retaining the intuitionistic semantic clauses for $\neg$, $\rightarrow$, $\wedge$, and $\forall$ and defining $(\varphi\vee\psi)$ as $\neg (\neg\varphi\wedge\neg\psi)$ and $\exists$ as $\neg\forall \neg$. Van Benthem calls this `possible world semantics' for classical logic, but we prefer `possibility semantics' for reasons that should now be clear. The main model-theoretic result of van Benthem's paper is a characterization of the classes of first-order possibility structures that are definable by a set of first-order sentences, viz., those classes that are closed under formation of generated submodels, disjoint unions, zig-zag images, filter products, and filter bases.

Van Benthem \citeyearpar{Benthem1981,Benthem1986,Benthem1988} also observes that the universe of all classically consistent sets of formulas, ordered by inclusion, is a canonical first-order possibility model. (\citealt{Benthem1981} considers taking only the finitely axiomatized consistent sets.) As he remarks: ``There is something inelegant to an ordinary Henkin argument. One has a consistent set of sentences $S$, perhaps quite small, that one would like to see satisfied semantically. Now, some arbitrary \textit{maximal} extensions $S^+$ of $S$ is to be taken to obtain a model (for $S^+$, and hence for $S$)---but the added part $S^+-S$ plays no role subsequently. We started out with something partial, but the method forces us to be total'' \citeyearpar[p.~78]{Benthem1988}. This ``problem of the `irrelevant extension'\,'' \citep[1]{Benthem1981} need not arise in completeness proofs for possibility semantics, as we saw with our use of filters instead of ultrafilters in \S~\ref{Canonical} and with our atomless possibility frames in \S~\ref{AtomlessFull}.

A more recent study of possibility semantics is \citealt{Garson2013}. Garson's (\S~8.8) notions of possibility models and forcing for propositional logic are essentially the same as ours, following Humberstone (see Remark \ref{ThreeVals}). Garson's main result about propositional possibility semantics can be seen as motivating the semantics \textit{proof-theoretically}, starting from the natural deduction rules for classical propositional logic. To state the result, we need some definitions. Let a (nontrivial) \textit{valuation} for $\mathcal{L}(\sig,\emptyset)$ be a $v\colon \mathcal{L}(\sig,\emptyset)\to\{0,1\}$ such that for some $\varphi\in\mathcal{L}(\sig,\emptyset)$, $v(\varphi)=0$. Given a set $\mathbb{V}$ of valuations for $\mathcal{L}(\sig,\emptyset)$, define $\mathcal{M}_\mathbb{V}=\langle \mathbb{V},\sqsubseteq_\mathbb{V},\pi_\mathbb{V}\rangle$ by:
\begin{itemize}
\item $\sqsubseteq_\mathbb{V}$ is the binary relation on $\mathbb{V}$ such that $v'\sqsubseteq_\mathbb{V} v$ iff for all $\varphi\in\mathcal{L}(\sig,\emptyset)$, if $v(\varphi)=1$, then $v'(\varphi)=1$;
\item $\pi_\mathbb{V}\colon \sig\to \wp (\mathbb{V})$ is such that $\pi_\mathbb{V}(p)=\{v\in\mathbb{V}\mid v(p)=1\}$.
\end{itemize}
Given a set $\mathbb{V}$ of valuations for $\mathcal{L}(\sig,\emptyset)$ and an \textit{argument} $\Gamma/\varphi$, which is a pair of $\Gamma\subseteq\mathcal{L}(\sig,\emptyset)$ and $\varphi\in\mathcal{L}(\sig,\emptyset)$, say that $\Gamma/\varphi$ is $\mathbb{V}$-\textit{valid}  iff for all $v\in\mathbb{V}$, if $v(\gamma)=1$ for each $\gamma\in\Gamma$, then $v(\varphi)=1$; and say that a natural deduction rule (an introduction or elimination rule for $\bot$, $\neg$, $\wedge$, $\vee$, $\rightarrow$, or $\leftrightarrow$)  \textit{preserves $\mathbb{V}$-validity} iff whenever the premise arguments of the rule are all $\mathbb{V}$-valid, the conclusion argument of the rule is also $\mathbb{V}$-valid.

\begin{theorem}[\citealt{Garson2013}, \S~8.8]\label{GarsonThm} For any nonempty set $\mathbb{V}$ of valuations for $\mathcal{L}(\sig,\emptyset)$, the following are equivalent:
\begin{enumerate} 
\item the natural deduction rules for classical propositional logic preserve $\mathbb{V}$-validity;
\item $\mathcal{M}_\mathbb{V}=\langle \mathbb{V},\sqsubseteq_\mathbb{V},\pi_\mathbb{V}\rangle$ is a propositional possibility model such that for all $v\in \mathbb{V}$ and $\varphi\in \mathcal{L}(\sig,\emptyset)$, $v(\varphi)=1$ iff $\mathcal{M}_\mathbb{V},v\Vdash \varphi$.
\end{enumerate}
\end{theorem}

\noindent Thus, one can ``read off'' propositional possibility semantics just from the assumption that the natural deduction rules for classical propositional logic preserve validity. As Garson \citeyearpar[\S~4.4]{Garson2013} shows, the same cannot be said for classical truth-table semantics. (For assumptions about the connection between natural deduction and propositional semantics sufficient to fix the classical truth tables, see \citealt{Bonay2015}.)\\    

\noindent\textit{Beth Semantics for Intuitionistic and Classical Logic}, \textit{Closure Operators on Heyting Algebras}, and \textit{Closure Frames for Substructural Logics}. As shown in \S~\ref{FromPartToPoss}, possibility semantics for classical (modal) logic is closely related to semantics for intuitionistic (modal) logic. However, the treatments of intuitionistic and classical logic using partial-state frames in \S~\ref{FromPartToPoss} were not as unified as they could be, due to the handling of disjunction on the intuitionistic side. Recall that we did not include disjunction as a primitive symbol of the classical language (Definition \ref{language}). Instead, we defined $\varphi\vee\psi$ as an abbreviation in terms of $\neg$ and $\wedge$. Of course, for the (full) intuitionistic language we need a primitive disjunction symbol, for which we used $\primvee$ (Definition~\ref{IntLanguage}). The disunity arose because we used the forcing clause for $\primvee$ from intuitionistic Kripke semantics (Example \ref{IntMod}), so $\vee$ and $\primvee$ had very different semantics:
\begin{itemize}
\item $\mathcal{M},x\Vdash \varphi\vee\psi$ iff $\forall x'\sqsubseteq x$ $\exists x''\sqsubseteq x'$: $\mathcal{M},x''\Vdash\varphi$ or $\mathcal{M},x''\Vdash\psi$;
\item $\mathcal{M},x\Vdash \varphi\primvee\psi$ iff $\mathcal{M},x\Vdash\varphi$ or $\mathcal{M},x\Vdash\psi$.
\end{itemize} 
In addition, this clause for $\primvee$ requires that the set $\adm$ of admissible sets in a frame be closed under \textit{unions}---in order for the logic of the frame to be closed under uniform substitution---which clashes with the requirement on full possibility frames that $\adm$ be the set of all regular open sets, which need not be closed under unions.

For a more unified treatment, we can move from Kripke semantics for intuitionistic logic to \textit{Beth} semantics \citep{Beth1956}. Dragalin \citeyearpar[p.~72ff]{Dragalin1988} presents a version of Beth semantics in which the difference between intuitionistic and classical logic emerges at the level of different frame classes, rather than different forcing clauses.\footnote{Thanks to Guram Bezhanishvili for pointing out connections between possibility semantics and Dragalin's version of Beth semantics, as well as connections with nuclei on Heyting algebras discussed below.} We will present a modified version of Dragalin's frames, which he calls \textit{Beth-Kripke frames}. (One of the differences will be that we include a set $\adm$ of admissible propositions.) For the case of propositional logic, instead of starting with partial-state frames $\mathcal{F}=\langle S,\sqsubseteq, \adm\rangle$ as in Definition \ref{PosetMod}, we start with richer frames $\mathcal{F}=\langle S,\sqsubseteq, \adm, Q\rangle$ where $\langle S,\sqsubseteq, \adm\rangle$ is a partial-state frame as before and $Q$ is a function assigning to each state $s\in S$ a set $Q(s)\subseteq \wp(S)$. An $X\in Q(s)$ is called a \textit{path starting from $s$}, so $Q(s)$ is the set of all paths starting from $s$. For $X,X'\subseteq S$, define $X'\sqsubseteq X$ iff $\forall x\in X$ $\exists x'\in X'$: $x'\sqsubseteq x$. Then the function $Q$ is required to satisfy at least the following natural conditions (for simplicity, we state the stronger version of Dragalin's second condition): first, $\emptyset\not\in Q(s)$; second, $X\in Q(s)$ implies $X\subseteq \mathord{\downarrow}s$; third, $s'\sqsubseteq s$ implies that $\forall X'\in Q(s')$ $\exists X\in Q(s)$ such that $X'\sqsubseteq X$; fourth, $s'\in X\in Q(s)$ implies $\exists X'\in Q(s')$ such that $X\sqsubseteq X'$; and finally, $Q(s)\not=\emptyset$ (for Dragalin's ``normal'' frames). As for $P$, it is required to be closed not only under intersection and the operation $\supset$ defined by $X\supset Y=\{s\in S\mid \forall s'\sqsubseteq s,\, s'\in X\Rightarrow s'\in Y\}\in \adm$ as in Definition \ref{PosetMod}, but also under the operation $+$ defined by $X+Y= \{s\in S\mid \forall Z\in Q(s) \;\exists z\in Z\colon z\in X\mbox{ or }z\in Y\}$. In other words, $s$ is in $X+Y$ iff every path starting from $s$ hits a state in $X$ or in $Y$. In addition, to obtain at least intuitionistic logic, we require that every $X\in P$ satisfies \textit{persistence} and
\begin{itemize}
\item \textit{barring} -- if $\forall Z\in Q(x)$ $\exists z\in Z$: $z\in X$, then $x\in X$. 
\end{itemize}
In other words, if every path starting from $x$ hits a state in $X$, then $x\in X$. Naturally, we define a \textit{full} frame as one in which $P$ is the set of all $X\subseteq S$ satisfying \textit{persistence} and \textit{barring}.

For the semantics, we keep the forcing clauses for $p$, $\neg$, $\wedge$, and $\rightarrow$  as in Definition \ref{pmtruth1}, throw away $\primvee$, and add a new forcing clause for a new primitive disjunction $\Bethvee$:
\begin{itemize}
\item $\mathcal{M},x\Vdash \varphi\Bethvee\psi$ iff $\forall Z\in Q(x)$ $\exists z\in Z$: $\mathcal{M},z\Vdash \varphi$ or $\mathcal{M},z\Vdash\psi$,
\end{itemize}
so $\llbracket \varphi\Bethvee\psi\rrbracket^\mathcal{M}=\llbracket \varphi\rrbracket^\mathcal{M}+\llbracket \psi\rrbracket^\mathcal{M}$ for the $+$ operation defined above. It is with this treatment of disjunction that we can achieve a more unified treatment of intuitionistic and classical logic. 

Intuitionistic propositional logic is sound and complete with respect to the class of \textit{all} frames $\mathcal{F}=\langle S,\sqsubseteq,\adm, Q\rangle$ satisfying the conditions above, according to the forcing relation just described. For soundness, one can check that the operations $\cap$, $\supset$, and $+$ above give rise to a Heyting algebra on $P$. For completeness, in the case where $Q(s)=\{\{s\}\}$ for each $s\in S$, the requirements on $Q$ hold and the forcing clause for $\Bethvee$ becomes the same as for $\primvee$ in Kripke semantics, so Kripke completeness implies Beth completeness. 

Now an observation of Dragalin (p.~74, Ex.~3) shows that we can use the same forcing relation and obtain soundness and completeness for \textit{classical} propositional logic as well. To do so, we do not consider all frames $\mathcal{F}=\langle S,\sqsubseteq,\adm, Q\rangle$ as above, but only those in which for every $x\in S$, $Q(x)=\{\mathord{\downarrow}y\mid y\sqsubseteq x\}$, i.e., the paths starting from $x$ are the principal downsets of refinements of $x$. Observe that given this definition of $Q$, the \textit{barring} condition above is equivalent to our \textit{refinability} condition, which in contrapositive form says that if $\forall x'\sqsubseteq x$ $\exists x''\sqsubseteq x'$ $x''\in X$, then $x\in X$. Also observe that given this definition of $Q$, the forcing clause for $\Bethvee$ becomes equivalent to the forcing clause for our classical $\vee$ above. From these observations it is a short step to the soundness and completeness of classical propositional logic. Furthermore, one can extend this analysis to intuitionistic and classical \textit{modal} logic based on Beth-style semantics.
   
There is a deeper perspective on this connection between Beth semantics and possibility semantics, explored in \citealt{BH2016,BH2017}. Given a frame $\mathcal{F}=\langle S,\sqsubseteq, P,Q\rangle$ as above, consider the Heyting algebra $\mathbb{H}(S,\sqsubseteq)$ of all downsets in $\langle S,\sqsubseteq\rangle$. Then consider the function $j$ on $\mathbb{H}(S,\sqsubseteq)$ such that for any downset $O$, $j(O)=\{x\in S\mid \forall Z\in Q(x) \,\exists z\in Z\colon z\in O\}$. Dragalin shows (pp.~72-3, using `$\boldsymbol{D}$' instead of `$j$') that this $j$ is a \textit{closure operator} on $\mathbb{H}(S,\sqsubseteq)$ in the sense of lattice theory \citep[\S~7.1]{Davey2002}, i.e., where $\leq$ is the natural order on the algebra, we have that for all elements $a$ and $b$ of the algebra: $a\leq b$ implies $j(a)\leq j(b)$; $a\leq j(a)$; and $j(j(a))=j(a)$. In fact, Dragalin shows that $j$ is a \textit{nucleus} on the Heyting algebra in the sense of point-free topology, i.e., a closure operator that also satisfies $j(a)\meet j(b)\leq j(a\meet b)$. There are two other important facts about the $j$ just defined using $Q$. First, together \textit{persistence} and \textit{barring} above imply that the sets $X\in\adm$ are \textit{fixed points} of $j$ in $\mathbb{H}(S,\sqsubseteq)$, i.e., $j(X)=X$; and if the frame is full, then $\adm$ is the set of \textit{all} fixed point of $j$. Second, the definition $Q(x)=\{\mathord{\downarrow}y\mid y\sqsubseteq x\}$ for classical frames implies that $j$ is the operation of \textit{double negation}. The first fact is important because for \textit{any} Heyting algebra (resp.~complete Heyting algebra) $\mathbb{H}$ and nucleus $j$, one obtains a new Heyting algebra (resp.~complete Heyting algebra) $\mathbb{H}_j$ by taking the fixed points of $j$ in $\mathbb{H}$, with the same meet and implication as in $\mathbb{H}$, but with a new join defined by $A+B =j(A\sqcup B)$, where $\sqcup$ is the join in $\mathbb{H}$. This is in essence what Beth semantics does. If $P$ is full, then $P$ together with $\cap$, $\supset$, and $+$ form the complete Heyting algebra $\mathbb{H}(S,\sqsubseteq)_j$. Even if $P$ is not full, our requirement that $P$ be closed under $\cap$, $\supset$, and $+$ guarantees that $P$ gives rise to a subalgebra of $\mathbb{H}(S,\sqsubseteq)_j$ and therefore a Heyting algebra. Finally, the second fact is important because if we form $\mathbb{H}_j$ with $j$ as double negation, then $\mathbb{H}_j$ is a Boolean algebra, which is complete if $\mathbb{H}$ is complete. We thereby obtain exactly the \textit{regular open algebra} from Remark \ref{Persp2} (which Dragalin calls the \textit{MacNeille algebra}). As we have seen, this is what possibility semantics does.

 The foregoing points lead to the idea of replacing $P$ and $Q$ in full frames $\mathcal{F}=\langle S,\sqsubseteq, P, Q\rangle$ by a nucleus $j$ on $\mathbb{H}(S,\sqsubseteq)$, or rather by something (such as a special binary relation or subframe) that determines a nucleus $j$ on $\mathbb{H}(S,\sqsubseteq)$, of which $j$ as double negation is but one example. Frames of this kind are studied in \citealt{Goldblatt1981}, \citealt{Fairtlough1997}, and \citealt{BH2016,BH2017}.

Dragalin's approach in terms of closure operators also appears in the semantics for substructural logics, including substructural modal logics, using \textit{closure frames} in \citealt[\S~12.2]{Restall2000}. In the non-modal case, a (full) closure frame is a tuple $\mathcal{F}=\langle S,\sqsubseteq, \Gamma\rangle$ where $\langle S,\sqsubseteq\rangle$ is a poset and $\Gamma\colon \wp(S)\to\wp(S)$ is a closure operator as above, but on the full powerset algebra rather than just the Heyting algebra of downsets, and by not requiring that $\Gamma$ is a nucleus, one can give semantics for \textit{non-distributive} logics. Otherwise the idea is as above: propositional variables must be interpreted as fixed points of $\Gamma$, and the semantic clause for disjunction is $\mathcal{M},x\Vdash \varphi\vee\psi$ iff $x\in\Gamma(\llbracket \varphi\rrbracket^\mathcal{M}\cup\llbracket \psi\rrbracket^\mathcal{M})$. (There is also a more general semantics for substructural negation.) Possibility semantics emerges in the case where $\Gamma (X)=\mathrm{int}(\mathrm{cl}(\mathord{\Downarrow}X))$ as in Fact \ref{RefReg}.\ref{RefReg2.5}.
 
Finally, the closure operator approach also appears in the \textit{truth-ground semantics} for intuitionistic propositional logic in \citealt{Rumfitt2012,Rumfitt2015}. Rumfitt considers nontrivial lower semilattices $\langle S,\sqsubseteq\rangle$ with a bottom element\footnote{Rumfitt goes up rather than down for refinements (writing `$\bullet$' for join), but we maintain our convention from Remark \ref{Flipped}.} and picks a specific closure operator on $\wp(S)$, namely the operator $(\cdot)^{ul}$ used for the MacNeille completion in \S~\ref{MacNeille}. As above, admissible propositions are taken to be the fixed points of the closure operator, and since the fixed points of $(\cdot)^{ul}$ are downsets, they are elements of the Heyting algebra $\mathbb{H}(S,\sqsubseteq)$. Under an assumption Rumfitt calls \textit{stability}, $(\cdot)^{ul}$ becomes a nucleus $j$ on $\mathbb{H}(S,\sqsubseteq)$. Since $\langle S,\sqsubseteq\rangle$ is a lower semilattice, the relative pseudocomplement in $\mathbb{H}(S,\sqsubseteq)$ may be defined by $X\rightarrowtriangle Y=\{s\in S\mid \forall x\in X,\, s\meet x\in Y\}$, which is how Rumfitt defines implication. As above, the conjunction of propositions is their intersection and the disjunction is the closure of their union. Thus, also as above, from the complete Heyting algebra $\mathbb{H}(S,\sqsubseteq)$ we obtain the complete Heyting algebra $\mathbb{H}(S,\sqsubseteq)_j$. Possibility semantics would instead take $j$ to be $\mathrm{int}(\mathrm{cl}(\mathord{\Downarrow}(\cdot)))$ with respect to the poset obtained from the bounded lower semilattice by deleting $\bot$. 

Since both Humberstone \citeyearpar{Humberstone1981} and Rumfitt \citeyearpar{Rumfitt2012,Rumfitt2015} speak of `possibilities', one could reasonably use the term `possibility semantics' for the general approach using some closure operator or other, reserving the term `classical possibility semantics' for the specific choice of $j$ as $\mathrm{int}(\mathrm{cl}(\mathord{\Downarrow}(\cdot)))$. \hfill $\triangleleft$\\

\noindent\textit{Possibility Semantics for Modal Logic}. The origin of possibility semantics for propositional modal logic is \citealt{Humberstone1981}. The important differences between Humberstone's frames and our possibility frames are discussed in \S~\ref{FullFrames} and Appendix \S~\ref{Strengths}. There is also a more superficial difference, namely that Humberstone's \citeyearpar{Humberstone1981} valuations were partial functions $V\colon \sig \times S\pfun \{0,1\}$ satisfying the condition of \textit{stability} and \textit{resolution} from above, but which Humberstone called `persistence' and `refinability'. Since the different approaches to valuations in the literature may be confusing, we provide the following guide.

\begin{remark}[Three Approaches to Valuations in Possibility Semantics]\label{ThreeVals} There are three approaches to valuation functions in the literature on possibility semantics:
\begin{enumerate}
\item\label{ThreeVals1} The approach in, e.g., \citealt{Humberstone1981} and \citealt{Holliday2014}: a valuation is a partial function $V\colon \sig \times S\pfun \{0,1\}$ satisfying \textit{stability} and \textit{resolution} above (but called `persistence' and `refinability' in the cited papers); $V(p,x)=1$ means that $x$ determines that $p$ is true; $V(p,x)=0$ means that $x$ determines that $p$ is false; $V(p,x)$ being undefined means that $x$ does not determine the truth or falsity of $p$.
\item\label{ThreeVals2} The approach in, e.g., \citealt{Garson2013} (\S~8.8): a valuation is a total function $u\colon \sig\times S\to \{0,1\}$ such that $\{x\in S\mid u(p,x)=1\}$ satisfies \textit{persistence} and \textit{refinability} in the sense of this paper; $u(p,x)=1$ means that $x$ determines that $p$ is true; $u(p,x)=0$ means that $x$ does not determine that $p$ is true, i.e., either $x$ determines that $p$ is false or $x$ does not determine the truth or falsity of $p$---which explains why this approach, unlike approach \ref{ThreeVals1}, does not require that if $u(p,x)=0$ and $x'\sqsubseteq x$, then $u(p,x')=0$.\footnote{Cf.~Kripke \citeyearpar[98]{Kripke1965} on his intuitionistic semantics: ``if $\phi(A, \mathbf{H}) = \mathbf{T}$ we can say that $A$ has been verified at the point $\mathbf{H}$ \dots; if $\phi(A,\mathbf{H}) = \mathbf{F}$, then $A$ has not been verified at $\mathbf{H}$. Notice, then, that $\mathbf{T}$ and $\mathbf{F}$ do not denote intuitionistic truth and falsity; if $\phi(A, \mathbf{H}) = \mathbf{T}$, then $A$ has been verified to be true at the time $\mathbf{H}$; but $\phi(A, \mathbf{H}) = \mathbf{F}$ does not mean that $A$ has been proved false at $\mathbf{H}$. It simply is not (yet) proved at $\mathbf{H}$; but may be established later.''}
\item\label{ThreeVals3} The approach in, e.g., \citealt{Humberstone2011} (\S~6.44) and the present paper: a valuation is a total function $\pi\colon \sig\to\wp(S)$ such that $\pi(p)$ satisfies \textit{persistence} and \textit{refinability}; $x\in\pi(p)$ means that $x$ determines that $p$ is true; $x\not\in\pi(p)$ means that $x$ does not determine that $p$ is true, i.e., that either $x$ determines that $p$ is false or $x$ does not determine the truth or falsity of $p$.
\end{enumerate}
The three approaches are all mathematically equivalent. We have already seen in Fact \ref{FineEq} how to go back and forth between valuations as in approaches \ref{ThreeVals1} and \ref{ThreeVals3}, and it is obvious how to go back and forth between valuations as in approaches \ref{ThreeVals2} and \ref{ThreeVals3}. An advantage of approach \ref{ThreeVals3} is that we can conveniently restate the constraint on admissible valuations as the constraint that $\pi\colon \sig\to\mathrm{RO}(S,\sqsubseteq)$ as in \S~\ref{ClassicalFrames}. \hfill $\triangleleft$\end{remark}

A direct follow-up to \citealt{Humberstone1981} is \citealt{Holliday2014}, which focuses on: functional possibility semantics as in \S~\ref{FuncFrames}; the construction of atomless canonical possibility models in which each possibility is given by a single finite formula (cf.~Appendix \S~\ref{Strengths}); and the closely related issue of internal adjointness mentioned in \S~\ref{SyntacticProp}. The results of \S~\ref{AtomlessFull} here can be seen as generalizing the results on atomless models in \citealt{Holliday2014}. 

The idea of giving a semantic clause for $\Box_i$ of the form $\mathcal{M},X\Vdash\Box_i\varphi$ iff $\mathcal{M},f_i(X)\Vdash \varphi$ appears earlier in Fine's \citeyearpar[p.~359]{Fine1974} discussion of relevance logic, in \citealt{Humberstone1988} (p.~418), and in \citealt{Humberstone2011} [p.~899]. The idea also appears in \citealt{Puncochar2014} (cf.~\citealt{Puncochar2017}), which presents a semantics for modal logic (using ``regular information models'') that is essentially equivalent to possibility semantics over functional and \textit{principal} possibility models as in \S~\ref{PrincFrames}. Although Pun\v{c}och\'{a}\v{r} uses an apparently different semantic clause for disjunction, namely the \textit{split disjunction} discussed below, we show in Lemma \ref{SplitDis} below that split disjunction is equivalent to our $\forall\exists$ forcing clause for disjunction over principal possibility models.

A more indirect follow-up to \citealt{Humberstone1981} is \citealt[Ch.~16]{Garson2013}, which discusses the extent to which something like Theorem \ref{GarsonThm} above extends to quantified modal logic. Explaining Garson's results for modal logic is beyond the scope of this overview, so we refer the reader to \citealt{Garson2013}.

Finally, we will briefly summarize the recent work on possibility semantics for modal logic in \citealt{Benthem2015} and \citealt{HT2016a,HT2016b}. The starting point of \citealt{Benthem2015} is the following \textit{bimodal} perspective on possibility semantics, which parallels previous bimodal perspectives on intuitionistic modal semantics \citep{Wolter1999}. A possibility frame $\mathcal{F}=\langle S,\sqsubseteq, R,\adm\rangle$ for a unimodal language gives rise to a frame $\langle S,\sqsubseteq, R\rangle$ for a bimodal language with modalities $[\sqsubseteq]$ and $[R]$ with the following semantics: $\mathcal{M},x\vDash [\sqsubseteq]\varphi$ iff for all $x'\sqsubseteq x$, $\mathcal{M},x'\vDash\varphi$; $\mathcal{M},x\vDash [R]\varphi$ iff for all $y\in R(x)$, $\mathcal{M},y\vDash\varphi$; and $\vDash$ treats the Boolean connectives as in ordinary possible world semantics. This bimodal perspective sheds light on both of our titular topics: possibility frames and possibility forcing. First, the frame conditions relating $\sqsubseteq$ and $R$ in \S~\ref{FullFrames} (recall Figure \ref{InterplayTable}) can be analyzed in terms of corresponding---in the precise sense of correspondence theory---bimodal interaction axioms relating $[\sqsubseteq]$ and $[R]$. Metatheoretic facts about the relations between frame conditions as in \S~\ref{FullFrames} can then be established by formal derivations in bimodal logic, as shown by van Benthem et al. Second, the possibility forcing relation $\Vdash$ and the requirement that admissible propositions be regular open sets suggests a \textit{translation} of the unimodal language with $\Box$ into the bimodal language with $[\sqsubseteq]$ and $[R]$. As shown by van Benthem et al., this translation embeds unimodal logics into a range of bimodal logics, including \textit{dynamic topological logics} as in \citealt{Kremer2005}.
 
The topic of \citealt{HT2016a} is possibility semantics for quantified modal logic. In addition to working out the basic setup of possibility semantics for quantified modal logic, Harrison-Trainor investigates the extent to which it is possible to do in the quantified modal case what \citealt{Humberstone1981} (p.~326) suggests and \citealt{Holliday2014} does in the propositional modal case: prove the completeness of standard modal logics using a simple canonical possibility model construction in which each possibility is identified with a \textit{finite} set of formulas (cf. \S~\ref{Strengths}). The construction in \citealt{Holliday2014} relies on the property of internal adjointness of certain propositional modal logics, discussed in \S~\ref{SyntacticProp}.\footnote{Note that the proof of Theorem \ref{CompAtomless} in \S~\ref{AtomlessFull} shows how to obtain completeness of a normal modal logic \textbf{L} with respect to an atomless and functional full possibility frame without assuming that \textbf{L} itself has internal adjointness, by taking a detour through the minimal tense extension of \textbf{L}, which always has internal adjointness.}  Harrison-Trainor shows that the quantified versions of those propositional modal logics lose the property of internal adjointness, as well as the weaker finite existence lemma from \S~\ref{SyntacticProp}. Thus, the direct analogue of the construction from \citealt{Holliday2014} does not work in these quantified modal cases. However, Harrison-Trainor also observes that one can prove completeness with a canonical possibility model construction in which each possibility is identified with a \textit{computable} set of formulas, which retains the spirit of the previous completeness results using finitary possibilities.

As reflected in the title of the present paper, our focus has been on possibility \textit{frames}---duality, definability, and completeness for  frames. This parallels the state of possible world semantics in its early decades, which focused on the theory of world frames (see, e.g., \citealt{Benthem1983}). The notion of validity over a class of frames always gives rise to a normal modal logic, whereas validity over a class of \textit{models} does not, since the set of formulas valid over a class of models may fail to be closed under Uniform Substitution. This is one reason for the focus on frames (cf.~\citealt[p.~159]{Hughes1996}). But just as the theory of possible world semantics expanded to include the study of models per se (see, e.g., \citealt{Blackburn2001}), we can expand the theory of possibility semantics in this direction. \citealt{HT2016b} adopts the model perspective and investigates how possibility models may be turned into and arise from Kripke models that validate the same formulas. (Recall from \S~\ref{NoKripke} that there can be no such transformation in general from full possibility frames to Kripke frames, so when Harrison-Trainor turns a possibility model into a Kripke model, their associated full \textit{frames} may validate different formulas; e.g., this will happen whenever starting from a possibility model based on the possibility frames in  \S~\ref{NoKripke}.) Harrison-Trainor shows how the method of generic chains used for Lemma \ref{classicality} can be extended to turn any countable possibility model for a countable language into a Kripke model---a \textit{worldization} of the possibility model---that has the same modal theory and bears a natural structural relation to the possibility model.\footnote{A less direct way of turning a possibility model $\mathcal{M}=\langle \mathcal{F},\pi\rangle$ (of any cardinality) based on $\mathcal{F}=\langle S,\sqsubseteq, \{R_i\}_{i\in\ind},\adm\rangle$ into a Kripke model that has the same modal theory is to first turn $\mathcal{M}$ into an algebraic model $\mathcal{M}^\under=\langle \mathcal{F}^\under, \pi\rangle$ as in Theorem \ref{PtoB}.\ref{PtoB5} and then turn $\mathcal{M}^\under$ into the Kripke model $(\mathcal{M}^\under)_+=\langle (\mathcal{F}^\under)_+, \pi_+\rangle$ based on the ultrafilter frame $(\mathcal{F}^\under)_+$ of $\mathcal{F}^\under$ as described in \S~\ref{AlgSem}. Whereas this construction builds the worlds of the Kripke model out of ultrafilters in $\langle \adm,\subseteq\rangle$, Harrison-Trainor's construction builds the worlds of the Kripke model out of ultrafilters in $\langle S,\sqsubseteq\rangle$.} In addition, Harrison-Trainor gives a general definition of a \textit{possibilization} of a Kripke model, generalizing the powerset possibilization of \S~\ref{PSFramesSem}, and shows that if $\mathcal{M}$ is a countable possibility model for a countable language and is separative as in \S~\ref{SepSec} and strong as in \S~\ref{FullFrames}, then $\mathcal{M}$ is isomorphic to a possibilization of a worldization of $\mathcal{M}$ (cf.~Propositions \ref{AtomProp} and \ref{PowChar}). Thus, we can see every such possibility model as arising from a Kripke model via possibilization. Crucially, this construction does not build in that \textit{every} set of worlds in the Kripke model becomes a possibility in the possibility model. Thus, unlike the views to be discussed below, Harrison-Trainor's construction does not imply a picture according to which the space $\langle S,\sqsubseteq\rangle$ of possibilities must be isomorphic to the powerset of a set of worlds (minus $\emptyset$) ordered by $\subseteq$.\\
 
\noindent\textit{Possibilities as Sets of Worlds}. Before Example \ref{PowerPoss}, we discussed the view of possibilities as arbitrary sets of worlds, leading to the definition of \textit{powerset possibilization}. Semantics for classical and intuitionistic propositional logic that evaluate formulas at sets of worlds appear in \citealt{Cresswell2004}. Semantics for modal logics that evaluate formulas at sets of worlds appear in recent work on \textit{inquisitive epistemic logic} \citep{Ciardelli2014b,Ciardelli2014}, where sets of worlds are taken to be ``information states,'' and in recent work on \textit{modal dependence logic} \citep{Vaananen2008,Hella2014} and \textit{modal independence logic} \citep{Kontinen2014}, where sets of worlds are called ``teams'' (cf.~the  discussion of inquisitive and team semantics in \citealt{Humberstone2019}, as well as the generalizations of inquisitive semantics in \citealt{Puncochar2016,Puncochar2017}). We will briefly discuss each of these frameworks.

For team semantics, we will not discuss the dependence and independence formulas that are the main point of modal dependence and independence logic, respectively. We mention only team semantics for the basic propositional modal language. For comparison with possibility semantics, recall from \S~\ref{ExtFrames} that the \textit{extended powerset possibilization} of a Kripke model $\mathfrak{M}=\langle \wo{W},\{\wo{R}_i\}_{i\in\ind},\wo{V}\rangle$ is the extended possibility model $\mathfrak{M}^\pow_\bot=\langle \wp(\wo{W}),\subseteq,\emptyset,\{R_i^\pow\}_{i\in\ind},\pi \rangle$ where $XR_i^\pow Y$ iff $Y\subseteq \wo{R}_i[X]$ and $\pi(p)=\{X\subseteq\wo{W}\mid X\subseteq\wo{V}(p)\}$. This $\mathfrak{M}^\pow_\bot$ is a \textit{quasi-functional} possibility model in the sense of \S~\ref{FuncFrames}, i.e., $R_i^\pow(X)$ has a maximum element, namely $\wo{R}_i[X]$, so the modal clause is equivalent to: $\mathfrak{M}^\pow_\bot,X\Vdash \Box_i\varphi$ iff $\mathfrak{M}^\pow_\bot,\wo{R}_i[X]\Vdash \varphi$. Similarly, \textit{team semantics} evaluates a formula at a set $T$ of worlds from a Kripke model $\mathfrak{M}=\langle \wo{W},\{\wo{R}_i\}_{i\in\ind},\wo{V}\rangle$, using the functional clause for the box modality: $\mathfrak{M},T\Vdash \Box_i\varphi$ iff $\mathfrak{M},\wo{R}_i[T]\Vdash\varphi$. What we wish to highlight is the clause for disjunction used in team semantics, which is called \textit{split disjunction}: $\mathfrak{M},T\Vdash \varphi\vee\psi$ iff $\exists T_1,T_2\subseteq\wo{W}$: $T=T_1\cup T_2$, $\mathfrak{M},T_1\Vdash \varphi$, and $\mathfrak{M},T_2\Vdash\psi$. It is easy to check that this is equivalent to our $\forall\exists$ forcing clause for $\vee$ from Fact \ref{ForcingDis} (holding fixed the other clauses from Definition \ref{pmtruth1}) when the space of possibilities/teams is the whole of $\wp(\wo{W})$, as in the extended powerset possibilization; but the clauses can differ over extended possibility models where $S$ is a proper subset of $\wp(\wo{W})$. For extended possibility models $\mathcal{M}=\langle S,\sqsubseteq,\bot, \{R_i\}_{i\in\ind},\pi \rangle$ as in \S~\ref{ExtFrames} where the states in $S$ are not assumed to be sets of worlds, we may consider the general split disjunction clause: $\mathcal{M},X\Vdash \varphi\vee\psi$ iff $\exists X_1,X_2$: $X=X_1\vee X_2$, $\mathcal{M},X_1\Vdash \varphi$, and $\mathcal{M},X_2\Vdash\psi$. Here $X_1\vee X_2$ is the least upper bound of $\{X_1,X_2\}$, if there is one, in the poset $\langle S,\sqsubseteq\rangle$. Clearly this clause is not equivalent to our $\forall\exists$ clause for $\vee$ over arbitrary possibility frames; but they are equivalent over \textit{principal} possibility frames as in \S~\ref{PrincFrames}.

\begin{lemma}[Split Disjunction in Principal Models]\label{SplitDis} For any extended principal possibility model $\mathcal{M}$, $X\in\mathcal{M}$, and $\varphi,\psi\in\mathcal{L}(\sig,\ind)$: $\mathcal{M},X\Vdash \neg(\neg\varphi\wedge\neg\psi)$ iff $\exists X_1,X_2$: $X=X_1\vee X_2$, $\mathcal{M},X_1\Vdash \varphi$, and $\mathcal{M},X_2\Vdash\psi$.
\end{lemma}

\begin{proof} By Fact \ref{propositions1} and the fact that the poset $\langle S,\sqsubseteq\rangle$ in an extended principal  $\mathcal{M}$ is a Boolean lattice (Fact~\ref{PrincEquiv}) with complement $-$, meet $\wedge$, and join $\vee$ we have $\mathcal{M},X\Vdash \neg(\neg\varphi\wedge\neg\psi)$ iff \[X\sqsubseteq \lVert \neg(\neg\varphi\wedge\neg\psi)\lVert^\mathcal{M}= -(-\lVert \varphi\rVert^\mathcal{M}\meet -\lVert \psi\rVert^\mathcal{M})=\lVert\varphi\rVert^\mathcal{M}\join \lVert\psi\rVert^\mathcal{M},\] i.e., $X=X\meet (\lVert\varphi\rVert^\mathcal{M}\join \lVert\psi\rVert^\mathcal{M})$, which is equivalent to $X= (X\meet \lVert\varphi\rVert^\mathcal{M})\join (X \meet\lVert\psi\rVert^\mathcal{M})$ by the Boolean laws. Then where $X_1=X\meet \lVert\varphi\rVert^\mathcal{M}$ and $X_2=X \meet\lVert\psi\rVert^\mathcal{M}$, the right side holds. Conversely, suppose $X=X_1\join X_2$, $X_1\sqsubseteq\lVert \varphi\rVert^\mathcal{M}$, and  ${X_2\sqsubseteq\lVert \psi\rVert^\mathcal{M}}$. Then we have $X= (X_1\meet \lVert \varphi\rVert^\mathcal{M})\join (X_2\meet \lVert \psi\rVert^\mathcal{M}) $, and since $X_1\sqsubseteq X$ and $X_2\sqsubseteq X$, we have that ${(X_1\meet \lVert \varphi\rVert^\mathcal{M})\join (X_2\meet \lVert \psi\rVert^\mathcal{M})}\sqsubseteq (X\meet \lVert \varphi\rVert^\mathcal{M})\join (X\meet \lVert \psi\rVert^\mathcal{M})\sqsubseteq X$. Thus, $X= (X\meet \lVert\varphi\rVert^\mathcal{M})\join (X \meet\lVert\psi\rVert^\mathcal{M})$, which we saw is equivalent to the left side.\end{proof}

Let us now return to \citealt{Cresswell2004}. The following observation is close to the main idea of that paper: while split disjunction behaves classically when the underlying poset is a Boolean lattice, it behaves intuitionistically when the underlying poset is such that all elements are \textit{join irreducible}, i.e., for every $X,Y,Z\in S$, if $X=Y\join Z$, then $X=Y$ or $X=Z$. In the join irreducible case, split disjunction requires either $\mathcal{M},X\Vdash\varphi$ or $\mathcal{M},X\Vdash\psi$. Roughly, Cresswell's idea is  to use the intuitionistic forcing clauses for $\neg$, $\wedge$, and $\rightarrow$, and the split clause for $\vee$, for \textit{both} classical and intuitionistic logic, while locating the difference between classical and intuitionistic logic in the assumptions about the underlying poset. (This is not exactly right, since Cresswell only evaluates formulas at nonempty sets of worlds, in which case split disjunction will not behave classically. Instead, he uses the following clause: $\mathfrak{M},X\Vdash \varphi\vee\psi$ iff $\exists X_1,X_2\in\wp(\wo{W})\setminus\{\emptyset\}$: $X\subseteq X_1\cup X_2$, [$\mathfrak{M},X_1\Vdash \varphi$ or $\mathfrak{M},X_1\Vdash\psi$], and [$\mathfrak{M},X_2\Vdash \varphi$ or $\mathfrak{M},X_2\Vdash\psi$]. See \citealt{Cresswell2004}, p.~22, for comparison of his disjunction and split disjunction, which he associates with Beth-Kripke-Joyal semantics for local set theory.)
 
Finally, let us draw some connections with the semantics for \textit{inquisitive epistemic logic} from \citealt{Ciardelli2014b} (Defs.~3 and 5) and \citealt{Ciardelli2014} (Defs.~2.2 and 2.4). For a direct comparison with the present paper, we will consider only the case of inquisitive modal logic without inquisitive modalities, as presented in Ch.~6 of \citealt{Ciardelli2016} (Def.~6.1.3). Inquisitive logic has important conceptual motivations, but we will not go into them here (see \citealt{Ciardelli2011,Ciardelli2013,Ciardelli2013b,Roelofsen2013}). We take the language of inquisitive modal logic to be the language $\mathcal{L}'(\sig,\ind)$ from Definition \ref{IntLanguage} that we used for intuitionistic modal logic.  Formulas of the form $\varphi \primvee \psi$ (written in the cited papers as `$?\{\varphi,\psi\}$' or `$\varphi\,\rotatebox[origin=c]{-90}{$\geqslant$}\,\psi$') are no longer thought of as \textit{declarative} disjunctions, but rather as \textit{interrogatives}. From the point of view of this paper, the semantics for inquisitive modal logic is equivalent to the following, as one can check by comparing the cited definitions from \citealt{Ciardelli2014b} and \citealt{Ciardelli2014,Ciardelli2016}.

\begin{definition}[Inquisitive Semantics for $\mathcal{L}'(\sig,\ind)$] $\,$ \\ An \textit{inquisitive frame} is a tuple $\mathfrak{F}^q=\langle S,\sqsubseteq,\bot, \{R_i^q\}_{i\in\ind},\adm \rangle$ that arises from a Kripke frame $\mathfrak{F}=\langle \wo{W},\{\wo{R}_i\}_{i\in\ind}\rangle$ as follows: $S=\wp (\wo{W})$; $X\sqsubseteq Y$ iff $X\subseteq Y$; $\bot=\emptyset$; $XR_i^qY$ iff $Y=\emptyset$ or for some $x\in X$, $\wo{R}_i(x)=Y$; and $\adm=\{\mathord{\downarrow}X\mid X\in S\}$, where as always, $\mathord{\downarrow}X=\{Y\in S\mid Y\sqsubseteq X\}$. 

An \textit{inquisitive model} is a tuple $\mathfrak{M}^q=\langle\mathfrak{F}^q,\pi \rangle$ that arises from a Kripke model $\mathfrak{M}=\langle \mathfrak{F},\wo{V}\rangle$ by setting $X\in \pi(p)$ iff $X\subseteq\wo{V}(p)$.  The inquisitive \textit{support} relation between pointed inquisitive models $\mathfrak{M}^q,X$ and formulas of $\mathcal{L}'(\sig,\ind)$ is the same as our forcing relation $\Vdash$ for $\mathcal{L}(\sig,\ind)$ from Definition \ref{ExtendedSemantics} extended to $\mathcal{L}'(\sig,\ind)$ as in Example \ref{IntMod}, so that $\mathfrak{M}^q,X\Vdash \varphi\primvee\psi$ iff $\mathfrak{M}^q,X\Vdash\varphi$ or $\mathfrak{M}^q,X\Vdash\psi$.\hfill $\triangleleft$
\end{definition}

Note that the inquisitive $\mathfrak{F}^q$/$\mathfrak{M}^q$ differs from the \textit{extended powerset possibilization} $\mathfrak{F}^\pow_\bot$/$\mathfrak{M}^\pow_\bot$ of $\mathfrak{F}$/$\mathfrak{M}$ (Example \ref{ExtPowerPoss}) only in the definition of the accessibility relation $R_i^q$.  Recall that $XR_i^\pow Y$ iff $Y\subseteq \wo{R}_i[X]$. Thus, $R_i^q$ is a subrelation of $R_i^\pow$, and it may be a proper subrelation.

Fact \ref{InqPoss} shows that inquisitive frames and models are a special case of extended possibility frames and models.  It also shows that for a Kripke frame $\mathfrak{F}$, the inquisitive frame $\mathfrak{F}^q$ is equivalent to the extended powerset possibilization $\mathfrak{F}^\pow_\bot$ with respect to $\mathcal{L}(\sig,\ind)$, though not necessarily $\mathcal{L}'(\sig,\ind)$ (see Example \ref{InqDiff}).

\begin{fact}[Inquisitive Frames as Possibility Frames]\label{InqPoss} For any Kripke frame $\mathfrak{F}$ and Kripke model $\mathfrak{M}$:
\begin{enumerate}
\item\label{InqPoss1} $\mathfrak{F}^q$ is an extended full possibility frame;
\item\label{InqPoss2} for all $X\in \mathfrak{M}^q$ and $\varphi\in\mathcal{L}(\sig,\ind)$, $\mathfrak{M}^q,X\Vdash\varphi$ iff $\mathfrak{M}^\pow_\bot,X\Vdash \varphi$.
\end{enumerate}
\end{fact}
\begin{proof} For part \ref{InqPoss1}, by Definition \ref{ExtendedFrames}, to say that $\mathfrak{F}^q$ is an extended full possibility frame is to say that the frame $\mathfrak{F}^q_{-}$ that results from deleting the bottom element $\emptyset$ from $\mathfrak{F}^q$ is a full possibility frame in the ordinary sense of Definition \ref{PosFrames}. The fullness requirement that $P_-=\mathrm{RO}(\mathfrak{F}^q_-)$ holds by the same reasoning as in the proof of Fact \ref{FullTFAE}. Then the key observation is that $R_i^q$ satisfies \upR{} and \Rwinweak{} from \S~\ref{FullFrames}, so by Proposition \ref{ROtoRO}, $P_-$ is closed under $\blacksquare_i$ as required for a partial-state frame.

For an inductive proof of part \ref{InqPoss2}, since $\mathfrak{M}^q$ and $\mathfrak{M}^\pow_\bot$ differ only with respect to their accessibility relations, the only case we need to check is the $\Box_i\varphi$ case. Since $R_i^q$ is a subrelation of $R_i^\pow$, $\mathfrak{M}^q,X\nVdash\Box_i\varphi$ implies $\mathfrak{M}^\pow_\bot,X\nVdash \Box_i\varphi$. Conversely, suppose $\mathfrak{M}^\pow_\bot,X\nVdash \Box_i\varphi$, so there is a $Y$ such that $XR_i^\pow Y$ and $\mathfrak{M}^\pow_\bot,Y\nVdash\varphi$. Since $\varphi\in\mathcal{L}(\sig,\ind)$, by Fact \ref{WtoP1}.\ref{WtoP1a} (which clearly also holds for the \textit{extended} powerset possibilization) $\mathfrak{M}^\pow_\bot,Y\nVdash\varphi$ implies that there is a $y\in Y$ such that $\mathfrak{M}^\pow_\bot,\{y\}\nVdash \varphi$, so $\mathfrak{M}^q,\{y\}\nVdash \varphi$ by the inductive hypothesis. Since $XR_i^\pow Y$, $Y\subseteq \wo{R}_i[X]$, so there is an $x\in X$ such that $x\mathrm{R}_iy$, so $\{y\}\subseteq \mathrm{R}_i(x)$ and hence $\{y\}\sqsubseteq \mathrm{R}_i(x)$. Since $\llbracket \varphi\rrbracket^{\mathfrak{M}^q}\in \adm$ satisfies \textit{persistence} in $\langle S,\sqsubseteq\rangle$, from $\mathfrak{M}^q,\{y\}\nVdash \varphi$ and $\{y\}\sqsubseteq \mathrm{R}_i(x)$ we have $\mathfrak{M}^q,\mathrm{R}_i(x)\nVdash \varphi$. Then since $x\in X$, we have $X R_i^q\mathrm{R}_i(x)$, so $\mathfrak{M}^q,\mathrm{R}_i(x)\nVdash \varphi$ implies $\mathfrak{M}^q,X\nVdash \Box_i\varphi$.\end{proof}

For part \ref{InqPoss1}, note that although the set $\adm$ of admissible propositions in $\mathfrak{F}^q$ is closed under the semantic operations corresponding to the operators $\neg$, $\wedge$, $\rightarrow$, and $\Box_i$ of $\mathcal{L}(\sig,\ind)$, it is not necessarily closed under the semantic operation corresponding to the operator $\primvee$ of $\mathcal{L}'(\sig,\ind)$, namely union. That is, given $\mathord{\downarrow}X\in \adm$ and $\mathord{\downarrow}Y\in\adm$, it does not follow that there is a $Z\in \adm$ such that $Z=\mathord{\downarrow}X\cup\mathord{\downarrow}Y$. This is the source of the fact that \textit{inquisitive logic}---the set of $\mathcal{L}'(\sig,\emptyset)$ formulas valid over all inquisitive frames---is not closed under uniform substitution (see \citealt{Ciardelli2009}). For example, although $\neg\neg p\rightarrow p$ is valid, $\neg\neg (p\primvee\neg p)\rightarrow (p\primvee\neg p)$ is not.

For part \ref{InqPoss2}, the semantic equivalence of $\mathfrak{M}^q$ and $\mathfrak{M}^\pow_\bot$ does not necessarily extend to $\mathcal{L}'(\sig,\ind)$. This is demonstrated by the following example, for which it is relevant that the intended meaning of $\Box_i (p\primvee\neg p)$ in inquisitive epistemic logic is that \textit{agent $i$ knows whether or not} $p$.

\begin{example}[Distinguishing $\mathfrak{M}^q$ and $\mathfrak{M}^\pow_\bot$ with $\primvee$]\label{InqDiff} Consider a Kripke model $\mathfrak{M}=\langle \wo{W},\{\wo{R}_i\}_{i\in\ind}, \wo{V}\rangle$ in which $\wo{W}=\{x_1, x_2, y_1, y_2\}$, $\wo{R}_i(x_1)=\{y_1\}$, $\wo{R}_i(x_2)=\{y_2\}$, and $\wo{V}(p)=\{y_1\}$. In $\mathfrak{M}^q$, we have $R_i^q(\{x_1,x_2\})=\{\emptyset,\{y_1\},\{y_2\}\}$, so $\mathfrak{M}^q,\{x_1,x_2\}\Vdash \Box_i (p\,\primvee\,\neg p)$; but in $\mathfrak{M}^\pow_\bot$, we have $R_i^\pow (\{x_1,x_2\})=\{\emptyset,\{y_1\},\{y_2\},\{y_1,y_2\}\}$, so $\mathfrak{M}^\pow_\bot,\{x_1,x_2\}\nVdash \Box_i (p \primvee\neg p)$ because $\{x_1,x_2\}R^\pow_i \{y_1,y_2\}$ and $\mathfrak{M}^\pow_\bot,\{y_1,y_2\}\nVdash p\primvee\neg p$. \hfill $\triangleleft$
\end{example}

In \S~\ref{OpenProb}, we will state as a problem for future work to investigate possibility semantics for languages extending the basic polymodal language $\mathcal{L}(\sig,\ind)$.  But in light of the above observations, we can see that the program of inquisitive semantics has already been doing this for $\mathcal{L}'(\sig,\ind)$. Something similar could be said about the team semantics mentioned above, but for an extended language with (in)dependence formulas.

\subsection{Open Problems}\label{OpenProb} 

We finish by listing a few of the open problems immediately suggested by our results.

As in \S~\ref{intro}, for a class $\mathcal{X}$ of BAOs or frames, let $\mathrm{ML}(\mathcal{X})$ be the set of modal logics \textbf{L} such that \textbf{L} is the logic of some subclass of $\mathcal{X}$. Where \textsf{K} is the class of \textit{Kripke} frames, \textsf{FP} is the class of \textit{full} possibility frames, \textsf{PR} is the class of \textit{principal} possibility frames, \textsf{f-PR} is the class of \textit{functional} principal possibility frames, and \textsf{P} is the class of all possibility frames---or we could take just filter-descriptive frames---we now know:
\[\begin{array}{ccccccccc}
\mathrm{ML}(\mathcal{CAV})& & \mathrm{ML}(\mathcal{CV})&& \mathrm{ML}(\mathcal{T})&& \mathrm{ML}(\mathcal{V})&&\mathrm{ML}(\mathcal{ALG})\\
 \rotatebox{90}{=} && \rotatebox{90}{=} && \rotatebox{90}{=} && \rotatebox{90}{=} &&\rotatebox{90}{=} \\
\mathrm{ML}(\mathsf{K})&\subsetneq& \mathrm{ML}(\mathsf{FP})&\subsetneq& \mathrm{ML}(\mbox{\textsf{f-PR}})&\subsetneq& \mathrm{ML}(\mathsf{PR})&\subsetneq&\mathrm{ML}(\mathsf{P}).\end{array}\]

 We know from the first strict inclusion that full possibility frames are more general than Kripke frames for characterizing normal modal logics. Indeed, we showed in \S~\ref{CompFull} that $\mathrm{ML}(\mathsf{FP})\setminus \mathrm{ML}(\mathsf{K})$ contains continuum many logics in the unimodal case alone. What more can we say about $\mathrm{ML}(\mathsf{FP})\setminus \mathrm{ML}(\mathsf{K})$?

\begin{problem}\label{Prob1} Find additional examples of logics in $\mathrm{ML}(\mathsf{FP})\setminus \mathrm{ML}(\mathsf{K})$, or equivalently, $\mathrm{ML}(\mathcal{CV})\setminus \mathrm{ML}(\mathcal{CAV})$, based on a different idea than the (\ref{Splitting}) formula of \S~\ref{NoKripke}.
\end{problem}

\begin{problem}\label{Prob2} Investigate degrees of \textsf{K}-incompleteness (Kripke-incompleteness) relative to $\textsf{FP}$ as in \S~\ref{CompFull}.
\end{problem}

\begin{problem}\label{Prob3} Give a syntactic or alternative semantic characterization of the logics in $\mathrm{ML}(\mathsf{FP})$, or of the logics in the difference $\mathrm{ML}(\mathsf{FP})\setminus \mathrm{ML}(\mathsf{K})$.
\end{problem}

A direction that may lead to a better understanding of $\mathsf{FP}$, mentioned at the end of \S~\ref{Canonical}, is the following.

\begin{problem} Investigate notions of persistence of modal formulas from possibility frames to their associated full possibility frames.
\end{problem}

The same kinds of questions as in Problems \ref{Prob1}-\ref{Prob3} arise from the other strict inclusions above. For example, we know from the third strict inclusion that arbitrary relations are more general than \textit{functions} for principal possibility frames, but what more can we say about this difference?

\begin{problem} Investigate analogues of Problems \ref{Prob1}-\ref{Prob3} for the other classes of frames/BAOs.
\end{problem} 

We discussed a number of other special classes of possibility frames in \S~\ref{FullFrames} and \S~\ref{SpecialClasses}, showing in several cases that any (full) possibility frame is semantically equivalent to one in a special class. One special class for which we did not prove such results is the class of Humberstone's \citeyearpar{Humberstone1981} original frames for possibility semantics (recall Remark \ref{HumbFrame} and see \S~\ref{Strengths}). Since every Kripke frame is a Humberstone frame and every Humberstone frame is a full possibility frame, we have $\mathrm{ML}(\mathsf{K})\subseteq \mathrm{ML}(\mathsf{H})\subseteq\mathrm{ML}(\mathsf{FP})$, where \textsf{H} is the class of Humberstone frames. We know that at least one of the inclusions is strict, but which?

\begin{problem}\label{HumProb} Which of the inclusions $\mathrm{ML}(\mathsf{K})\subseteq \mathrm{ML}(\mathsf{H})$ and $\mathrm{ML}(\mathsf{H})\subseteq \mathrm{ML}(\mathsf{FP})$ is strict?
\end{problem}

Turning from completeness to correspondence, the next problem is a follow-up to the discussion of \S~\ref{LemmScottCorr}.

\begin{problem}\label{Prob6} Does every modal formula that has a first-order correspondent over Kripke frames also have a first-order correspondent over full possibility frames?
\end{problem} 

Finally, although here we have restricted attention to the basic polymodal language for a direct comparison of possibility semantics and world semantics, a natural step is to see how the richer refinement structure of possibility frames could be exploited for the semantics of extended modal languages.

\begin{problem}\label{ExtProb} Investigate possibility semantics for extended modal languages.
\end{problem}

\section*{Acknowledgments}

For helpful comments or discussion, I wish to give special thanks to Johan van Benthem, Guram Bezhanishvili, Nick Bezhanishvili, Matthew Harrison-Trainor, Tadeusz Litak, and an anonymous referee for \textit{The Australasian Journal of Logic}. I also wish to thank Ivano Ciardelli, Josh Dever, David Gabelaia, Davide Grossi, Lloyd Humberstone, Thomas Icard, Mamuka Jibladze, Hans Kamp, Larry Moss, Lawrence Valby, Yanjing Wang, and Dag Westerst\aa hl, as well as the participants in my Fall 2014 or Spring 2015 graduate seminars at UC Berkeley: Russell Buehler, Sophia Dandelet, Matthew Harrison-Trainor, Alex Kocurek, Alex Kruckman, James Moody, James Walsh, and Kentaro Yamamoto. I am also thankful for feedback I received when presenting some of this material at the following venues: the Modal Logic Workshop on Consistency and Structure at Carnegie Mellon University in April 2014; my course at the 3rd East-Asian School on Logic, Language and Computation (EASSLLC) at Tsinghua University in July 2014; the Advances in Modal Logic conference at the University of Groningen in August 2014 (see \citealt{Holliday2014}); the Workshop on the Future of Logic in honor of Johan van Benthem in Amsterdam in September 2014; the Hans Kamp Seminar in Logic and Language at the University of Texas at Austin in April 2015; the 4th CSLI Workshop on Logic, Rationality and Intelligent Interaction at Stanford University in May 2015; and the New Mexico State University Mathematics Colloquium in November 2015. Finally, I wish to gratefully acknowledge an HRF grant from UC Berkeley that allowed me to complete most of this work in Fall 2015 (see \citealt{Holliday2015}).

\section*{Appendices}

\appendix

\section{Review of Standard Semantics}\label{Review}

\subsection{Kripke Semantics}\label{WRM}
 
To fix terminology and notation, we review the standard definitions for Kripke semantics here.
 
\begin{definition}[Kripke Frames and Models]\label{RelWorld} A \textit{Kripke frame} is a tuple $\mathfrak{F}=\langle \wo{W},\{\wo{R}_i\}_{i\in\ind }\rangle$ where $\wo{W}$ is a nonempty set (the set of \textit{worlds}) and $\wo{R}_i$ is a binary relation on $\wo{W}$ (the \textit{$i$-accessibility relation}). A \textit{Kripke model} based on $\mathfrak{F}$ is a tuple $\mathfrak{M}=\langle \mathfrak{F},\wo{V}\rangle$ where $\wo{V}\colon \sig\to \wp(\wo{W})$ (a \textit{valuation}). \hfill $\triangleleft$
\end{definition}

We use `$\vDash$' for the standard satisfaction relation in Kripke models, in contrast with `$\Vdash$' in Definition \ref{pmtruth1}.

\begin{definition}[Kripke Semantics]\label{pwmtruth1} Given a Kripke model $\mathfrak{M}=\langle \wo{W},\{\wo{R}_i\}_{i\in\ind },\wo{V}\rangle$ with $w\in \wo{W}$ and $\varphi\in\mathcal{L}(\sig,\ind)$,  define $\mathfrak{M},w\vDash\varphi$ (``$\varphi$ is true at $w$ in $\mathfrak{M}$'') recursively as follows:
\begin{enumerate}
\item $\mathfrak{M},w\vDash p$ iff $w\in\wo{V}(p)$;
\item $\mathfrak{M},w\vDash\neg\varphi$ iff $\mathfrak{M},w\nvDash\varphi$;
\item $\mathfrak{M},w\vDash \varphi\wedge\psi$ iff $\mathfrak{M},w\vDash\varphi$ and $\mathfrak{M},w\vDash\psi$;
\item $\mathfrak{M},w\vDash \varphi\rightarrow\psi$ iff $\mathfrak{M},w\nvDash\varphi$ or $\mathfrak{M},w\vDash\psi$;
\item $\mathfrak{M},w\vDash\Box_i\varphi$ iff $\forall v\in\mathrm{R}_i(w)$: $\mathfrak{M},v\vDash\varphi$,
\end{enumerate}
where as in Notation \ref{notation}, $\mathrm{R}_i(w)=\{v\in\wo{W}\mid w\mathrm{R}_iv\}$.

Validity, satisfiability, soundness, and completeness are defined as in Definition \ref{pmtruth1}.\hfill $\triangleleft$
\end{definition}

Recall the standard Henkin-style canonical model construction for Kripke semantics.

\begin{definition}[Canonical Kripke Frames and Models]\label{CanKrip} The \textit{canonical Kripke frame} for a consistent normal modal logic $\mathbf{L}$ is the structure $\mathfrak{F}^\mathbf{L}=\langle \wo{W}^\mathbf{L},\{\wo{R}_i^\mathbf{L}\}_{i\in\ind}\rangle$ where:
\begin{enumerate}[label=\arabic*.,ref=\arabic*]
\item $\wo{W}^\mathbf{L}$ is the set of all maximally \textbf{L}-consistent sets of formulas;
\item for $\Gamma,\Delta\in\wo{W}^\mathbf{L}$, $\Gamma \wo{R}_i^\mathbf{L}\Delta$ iff for all $\varphi\in\mathcal{L}(\sig,\ind)$, $\Box_i\varphi\in \Gamma$ implies $\varphi\in\Delta$.
\end{enumerate}
The \textit{canonical Kripke model} for $\mathbf{L}$ is the structure $\mathfrak{M}^\mathbf{L}=\langle \mathfrak{F}^\mathbf{L},\wo{V}^\mathbf{L}\rangle$ where $\mathfrak{F}^\mathbf{L}$ is the canonical Kripke frame for $\mathbf{L}$ and for all $p\in\sig$ and $\Gamma\in\wo{W}^\mathbf{L}$: $\Gamma\in\wo{V}^\mathbf{L}(p)$ iff $p\in\Gamma$. \hfill $\triangleleft$
\end{definition}

As shown in textbooks on modal logic, every consistent normal modal logic \textbf{L} is sound and complete with respect to its canonical Kripke model $\mathfrak{M}^\mathbf{L}$, so every such \textbf{L} is complete with respect to its canonical Kripke frame $\mathfrak{F}^\mathbf{L}$; but not  every such \textbf{L} is \textit{sound} with respect to $\mathfrak{F}^\mathbf{L}$---not every such \textbf{L} is \textit{canonical}. In fact, there are many normal modal logics that are not sound and complete with respect to \textit{any} class of Kripke frames (see Footnote \ref{IncompRefs}), which leads to the following distinction.

\begin{definition}[Kripke Completeness] A normal modal logic \textbf{L} is \textit{Kripke-frame complete} iff there is a class $\mathsf{F}$ of Kripke frames for which \textbf{L} is sound and complete. Otherwise it is \textit{Kripke-frame incomplete}. \hfill $\triangleleft$
\end{definition}

The existence of Kripke-frame incomplete normal modal logics is one of the motivations for the semantics of the following sections.

\subsection{General Frame Semantics}\label{GFS} 

The following more general semantics for normal modal logics is due to \citealt{Thomason1972} (cf.~\citealt{Makinson1970}, and see \citealt[\S~1.4, \S~5.5]{Blackburn2001} for a textbook treatment). 

\begin{definition}[General Frame Semantics]\label{GenSem} A \textit{general frame}---or in the terminology of this paper, a \textit{world frame}---is a tuple $\mathfrak{g}=\langle \wo{W},\{\wo{R}_i\}_{i\in\ind },\wo{A}\rangle$ where $\langle \wo{W},\{\wo{R}_i\}_{i\in\ind } \rangle$ is a Kripke frame and $\wo{A}\subseteq \wp(\wo{W})$ is a nonempty set (of \textit{admissible propositions}) such that if $X,Y\in \wo{A}$, then $\wo{W}\setminus X\in \wo{A}$, $X\cap Y\in\wo{A}$, and $\{w\in \wo{W}\mid \wo{R}_i(w)\subseteq X\}\in \wo{A}$ for each $i\in \ind$. 

A Kripke model $\mathfrak{M}=\langle \wo{W},\{\wo{R}_i\}_{i\in\ind },\wo{V}\rangle$ is \textit{based on} the general frame $\mathfrak{g}$ iff for every $p\in \sig$, $\wo{V}(p)\in\wo{A}$, in which case $\wo{V}$ is an \textit{admissible valuation} for $\mathfrak{g}$. We may regard such an $\mathfrak{M}$ as the pair $\langle \mathfrak{g},\wo{V}\rangle$. 

The notions of validity, soundness, and completeness with respect to a class $\mathsf{G}$ of general frames are defined in terms of these notions with respect to the class of Kripke models based on general frames in $\mathsf{G}$. \hfill $\triangleleft$
\end{definition} 

Note that in a Kripke model $\mathfrak{M}$, we have $\wo{W}\setminus\llbracket \varphi\rrbracket^\mathfrak{M}=\llbracket \neg\varphi\rrbracket^\mathfrak{M}$, $\llbracket \varphi\rrbracket^\mathfrak{M}\cap\llbracket\psi\rrbracket^\mathfrak{M}=\llbracket \varphi\wedge\psi\rrbracket^\mathfrak{M}$, and $\{{w\in \wo{W}}\mid \wo{R}_i(w)\subseteq\llbracket \varphi\rrbracket^\mathfrak{M}\}=\llbracket \Box_i\varphi\rrbracket^\mathfrak{M}$, so the closure conditions on $\wo{A}$ in $\mathfrak{g}$ ensure that for any $\mathfrak{M}$ based on $\mathfrak{g}$ and $\varphi\in\mathcal{L}(\sig,\ind)$, $\llbracket\varphi\rrbracket^\mathfrak{M}\in\wo{A}$. This guarantees that the set of formulas valid over a class of general frames is closed under Uniform Substitution (Definition \ref{NML}), as is guaranteed for Kripke frames but not for Kripke models.

Also note that any Kripke frame $\mathfrak{F}=\langle \wo{W},\{\wo{R}_i\}_{i\in\ind}\rangle$ can be equivalently viewed as the general frame $\mathfrak{F}^\sharp =\langle\wo{W},\{\wo{R}_i\}_{i\in\ind},\wp(\wo{W})\rangle$, called the \textit{full general frame of} $\mathfrak{F}$, in the sense that the class of Kripke models based on $\mathfrak{F}$ is the same as the class of Kripke models based on $\mathfrak{F}^\sharp$, so $\mathfrak{F}$ and $\mathfrak{F}^\sharp$ determine the same logic.

Note that while in this section we use lower-case Fraktur letters for general frames and upper-case Fraktur letters for Kripke frames (as in \citealt{Blackburn2001}), for simplicity in the main text we use upper-case Fraktur letters for all ``world frames,'' i.e., general frames and Kripke frames regarded as full general frames. 

\begin{definition}[Canonical General Frame]\label{CanGenFrame} The \textit{canonical general frame} for a consistent normal modal logic $\mathbf{L}$  \citep[p.~306]{Blackburn2001} is the structure $\mathfrak{g}^\mathbf{L}=\langle \mathfrak{F}^\mathbf{L},\wo{A}^\mathbf{L}\rangle$ where $\mathfrak{F}^\mathbf{L}$ is the canonical Kripke frame for \textbf{L} as in Definition \ref{CanKrip} and $\wo{A}^\mathbf{L}=\{X\subseteq \wo{W}^\mathbf{L}\mid \exists \varphi\in\mathcal{L}(\sig,\ind)\colon X= \{\Gamma\in \wo{W}^\mathbf{L}\mid \varphi\in\Gamma\}\}$. By the Truth Lemma \citep[\S~4.2]{Blackburn2001}, $\wo{A}^\mathbf{L}=\{\llbracket \varphi\rrbracket^{\mathfrak{M}^\mathbf{L}}\mid \varphi\in\mathcal{L}(\sig,\ind)\}$, so the admissible propositions of $\mathfrak{g}^\mathbf{L}$ are the sets of worlds definable in the canonical model $\mathfrak{M}^\mathbf{L}$ by formulas of $\mathcal{L}(\sig,\ind)$. \hfill $\triangleleft$
\end{definition}

The following result (see, e.g., \citealt[Thm.~5.64]{Blackburn2001}) shows that general frames, unlike Kripke frames, can characterize any normal modal logic.

\begin{theorem}[Adequacy of General Frame Semantics]\label{GenAdeq} Every consistent normal modal logic \textbf{L} is sound and complete with respect to its canonical general frame $\mathfrak{g}^\mathbf{L}$.\end{theorem}

\subsection{Algebraic Semantics}\label{AlgSem}

In Kripke and general frame semantics, we first define the truth of a formula at a world in a model $\mathfrak{M}$ and then derivatively obtain a mapping $\varphi\mapsto \llbracket\varphi\rrbracket^\mathfrak{M}$ of formulas to ``propositions,'' understood as sets of worlds. In the algebraic semantics for modal logic (see \citealt[\S~5.2]{Blackburn2001} for a textbook treatment), we cut out the worlds and directly map formulas to propositions, taken as primitive objects. In the following definition, due to \citealt{Jonsson1952a,Jonsson1952b}, one may think of elements of $A$ as propositions.

\begin{definition}[Boolean Algebra with Operators]\label{BAOs} A \textit{Boolean algebra with \textnormal{(}unary, dual\textnormal{)} operators} (BAO) is a tuple $\mathbb{A}=\langle A, \meet, -, \top, \{\blacksquare_i\}_{i\in \ind}\rangle$ where $\langle A, \wedge,-,\top\rangle$ is a Boolean algebra with $\meet$ as meet, $-$ as complement, and $\top$ as the top element, and each $\blacksquare_i\colon A\to A$ satisfies:
\begin{enumerate}[label=\arabic*.,ref=\arabic*]
\item $\blacksquare_i \top =\top$;
\item for all $x,y\in A$, $\blacksquare_i(x\wedge y)= \blacksquare_i x\wedge \blacksquare_i y$. 
\end{enumerate}
Equivalently, where $\join$ is join, $\bot$ is the bottom element, and for $x\in A$, $\blacklozenge_i x:= -\blacksquare_i - x$:
\begin{enumerate}[label=\arabic*.,ref=\arabic*,resume]
\item $\blacklozenge_i \bot=\bot$;
\item for all $x,y\in A$, $\blacklozenge_i(x\vee y)=\blacklozenge_i x\join \blacklozenge_i y$. 
\end{enumerate}
A BAO $\mathbb{A}=\langle A, \meet, -, \top, \{\blacksquare_i\}_{i\in \ind}\rangle$ is \textit{trivial} if $|A|=1$ and \textit{nontrivial} otherwise.\hfill $\triangleleft$
\end{definition}
BAOs are often defined with the additive \textit{operators} $\blacklozenge_i$ as primitive, rather than the multiplicative \textit{dual operators} $\blacksquare_i$. But for our purposes, it will be more convenient to take the dual operators as primitive. We also trust that no confusion will arise by using the same symbol `$\wedge$' for conjunction in our formal language and the meet operation in our Boolean algebras, and similarly for `$\vee$' with disjunction and join. 

Turning to the semantics, instead of mapping each $p\in\sig$ to a set of worlds where it holds, as in a Kripke model, an algebraic model based on $\mathbb{A}$ maps each $p\in\sig$ to an element of $A$, thought of as a proposition, and this mapping extends to a mapping of each formula of $\mathcal{L}(\sig,\ind)$ to an element of $A$.

\begin{definition}[Algebraic Semantics]\label{AlgSemantics} An \textit{algebraic model} is a pair $\mathbb{M}=\langle \mathbb{A},\theta\rangle$ where $\mathbb{A}$ is a BAO and $\theta\colon \sig\to A$. We extend $\theta$ to a meaning function $\tilde{\theta}\colon \mathcal{L}(\sig,\ind)\to A$ defined by: $\tilde{\theta}(p)=\theta(p)$;  $\tilde{\theta}(\neg\varphi)= - \tilde{\theta}(\varphi)$;   $\tilde{\theta}(\varphi\wedge\psi)=\tilde{\theta}(\varphi)\wedge \tilde{\theta}(\psi)$; $\tilde{\theta}(\varphi\to\psi)=-\tilde{\theta}(\varphi)\vee \tilde{\theta}(\psi)$; and $\tilde{\theta}(\Box_i\varphi)=\blacksquare_i\tilde{\theta}(\varphi)$.

A formula $\varphi\in\mathcal{L}(\sig,\ind)$ is \textit{valid} over a class $\mathsf{C}$ of BAOs iff for every algebraic model $\mathbb{M}=\langle \mathbb{A},\theta\rangle$ with $\mathbb{A}\in\mathsf{C}$, $\tilde{\theta}(\varphi)=\top$; and $\varphi$ is \textit{satisfiable} over $\mathsf{C}$ iff there is an algebraic model $\mathbb{M}=\langle \mathbb{A},\theta\rangle$ with $\mathbb{A}\in\mathsf{C}$ and $\tilde{\theta}(\varphi)\not=\bot$. Soundness and completeness of a modal logic \textbf{L} with respect to a class $\mathsf{C}$ of BAOs are defined in terms of validity over $\mathsf{C}$, as usual.\hfill $\triangleleft$
\end{definition} 

For any general frame $\mathfrak{g}=\langle \wo{W},\{\wo{R}_i\}_{i\in\ind},\wo{A}\rangle$, the structure $\mathfrak{g}^*=\langle \wo{A},\cap, - ,\wo{W}, \{\blacksquare_i\}_{i\in\ind}\rangle$ with $-\colon \wo{A}\to \wo{A}$ and $\blacksquare_i\colon \wo{A}\to \wo{A}$ defined by $-X=\wo{W}\setminus X$ and $\blacksquare_i X=\{w\in \wo{W}\mid \wo{R}_i(w)\subseteq X \}$ is a BAO, called the \textit{underlying BAO of $\mathfrak{g}$}. For a Kripke frame $\mathfrak{F}=\langle \wo{W},\{\wo{R}_i\}_{i\in\ind}\rangle$, the structure $\mathfrak{F}^+=\langle \wp(\wo{W}),\cap, - ,\wo{W}, \{\blacksquare_i\}_{i\in\ind} \rangle$ with $-$ and $\blacksquare_i$ as just defined is called the \textit{full complex algebra of $\mathfrak{F}$}, which is the underlying BAO of the full general frame of $\mathfrak{F}$ (recall \S~\ref{GFS}), i.e., $\mathfrak{F}^+=(\mathfrak{F}^\sharp)^*$. One can check that for any general frame $\mathfrak{g}$ and $\varphi\in\mathcal{L}(\sig,\ind)$, $\varphi$ is valid over $\mathfrak{g}$ according to Definition \ref{GenSem} iff $\varphi$ is valid over $\mathfrak{g}^*$ according to Definition \ref{AlgSemantics}. Thus, any general frame or Kripke frame can be turned into a BAO that validates the same formulas. 

In the other direction, for any BAO $\mathbb{A}=\langle A, \meet, -, \top, \{\blacksquare_i\}_{i\in \ind}\rangle$, the structure $\mathbb{A}_+=\langle Uf\mathbb{A},\{\wo{R}_i\}_{i\in\ind}\rangle$ where $Uf\mathbb{A}$ is the set of ultrafilters in $\mathbb{A}$, and $u\mathrm{R}_iu'$ iff $\forall x\in A$, $\blacksquare_i x\in u$ implies $x\in u'$, is a Kripke frame, called the \textit{ultrafilter frame} of $\mathbb{A}$. (Here one goes beyond ZF set theory and assumes the ultrafilter axiom.) The structure $\mathbb{A}_*=\langle \mathbb{A}_+, \{\hat{a}\mid a\in A\}\rangle$ with $\hat{a}=\{u\in Uf\mathbb{A} \mid a\in u\}$ is a general frame, called the \textit{general ultrafilter frame} of $\mathbb{A}$. One can check that for any BAO $\mathbb{A}$ and $\varphi\in\mathcal{L}(\sig,\ind)$, $\varphi$ is valid over $\mathbb{A}$ according to Definition \ref{AlgSemantics} iff $\varphi$ is valid over $\mathbb{A}_*$ according to Definition \ref{GenSem}. But that `iff' may fail for $\mathbb{A}_+$: although for any algebraic model $\mathbb{M}=\langle \mathbb{A},\theta\rangle$ based on $\mathbb{A}$, the Kripke model $\mathbb{M}_+=\langle \mathbb{A}_+,\theta_+\rangle$ with $\theta_+(p)=\{u\in Uf\mathbb{A}\mid \theta(p)\in u\}$ is modally equivalent to $\mathbb{M}$, there may be Kripke models based on $\mathbb{A}_+$ that are not modally equivalent to any algebraic model based on $\mathbb{A}$, because $\mathbb{A}_+$ imposes no constraints on admissible valuations.

Given a normal modal logic \textbf{L} and $\varphi,\psi\in\mathcal{L}(\sig,\ind)$, let $\varphi\sim_\mathbf{L} \psi$ iff $\vdash_\mathbf{L} \varphi\leftrightarrow\psi$, and $[\varphi]_\mathbf{L}=\{\psi\in\mathcal{L}(\sig,\ind)\mid \varphi\sim_\mathbf{L}\psi\}$. One can check that $\sim_\mathbf{L}$ is a congruence relation with respect to the structure $\langle \mathcal{L}(\sig,\ind), O_\wedge, O_\neg,\top,\{O_i\}_{i\in\ind}\rangle$ where $O_\wedge (\varphi,\psi)=(\varphi\wedge\psi)$, $O_\neg (\varphi)=\neg\varphi$, $\top = (p\vee\neg p)$, and $O_i(\varphi)=\Box_i\varphi$. Thus, we can take the quotient of this structure with respect to $\sim_\mathbf{L}$, obtaining the following.

\begin{definition}[Lindenbaum Algebra]\label{LinAlg}  The \textit{Lindenbaum algebra} for a normal modal logic \textbf{L} is the structure  $\mathbb{A}^\mathbf{L}=\langle A, \wedge, -, \top, \{\blacksquare_i\}_{i\in\ind} \rangle$ where: $A=\{[\varphi]_\mathbf{L}\mid \varphi\in\mathcal{L}(\sig,\ind)\}$; $[\varphi]_\mathbf{L}\meet [\psi]_\mathbf{L}=[(\varphi\wedge\psi)]_\mathbf{L}$;  $-[\varphi]_\mathbf{L}=[\neg\varphi]_\mathbf{L}$; $\top=[(p\vee\neg p)]_\mathbf{L}$; and  $\blacksquare_i [\varphi]_\mathbf{L}=[\Box_i\varphi]_\mathbf{L}$. \hfill$\triangleleft$
\end{definition}

\noindent One can check that $\mathbb{A}^\mathbf{L}$ is indeed a BAO. If $\mathbf{L}$ is consistent, then $\mathbb{A}^\mathbf{L}$ is isomorphic to the underlying BAO of the canonical general frame $\mathfrak{g}^\mathbf{L}$ from \S~\ref{GFS}, and $\mathfrak{g}^\mathbf{L}$ is isomorphic to the general ultrafilter frame of $\mathbb{A}^\mathbf{L}$. 

The following result (see \citealt[\S~5.2]{Blackburn2001}) is an algebraic analogue of Theorem \ref{GenAdeq}.

\begin{theorem}[Adequacy of Algebraic Semantics]\label{AdAlg} Every normal modal logic \textbf{L} is sound and complete with respect to its Lindenbaum algebra $\mathbb{A}^\mathbf{L}$.
\end{theorem}

\section{Deferred Topics}\label{Deferred}

\subsection{Stronger Refinability}\label{Strengths}

Recall from \S~\ref{FullFrames} that Humberstone's \citeyearpar{Humberstone1981} original frames for possibility semantics were, in the terminology of this paper, full possibility frames satisfying \upR{}, \Rdown{}, and \RrefPlusPlus{}. In this section, we compare the following refinability conditions, from weaker to stronger:
 
\begin{itemize}
\item[] \Rref{} -- if ${x}R_i{y}$, then $\exists{x'}\sqsubseteq {x}$ $\forall {x''}\sqsubseteq {x'}$ \underline{$\exists{y'}\sqsubseteq {y}$}: ${x''}R_i{y'}$;
\item[] \RrefPlus{} -- if ${x}R_i{y}$, then \underline{$\exists{y'}\sqsubseteq {y}$} $\exists{x'}\sqsubseteq {x}$ $\forall {x''}\sqsubseteq {x'}$\; ${x''}R_i{y'}$;
\item[] \RrefPlusPlus{} -- if ${x}R_i{y}$, then $\exists{x'}\sqsubseteq {x}$ $\forall {x''}\sqsubseteq {x'}$\; ${x''}R_i{y}$.
\end{itemize}

In \S~\ref{FullFrames}, we mentioned the following fact.

\begin{fact}[Powerset Possibilizations and \RrefPlusPlus{}] There are Kripke frames $\mathfrak{F}$ whose powerset possibilizations $\mathfrak{F}^\pow$ do not satisfy \RrefPlusPlus{}.
\end{fact}

\begin{proof} For $\mathfrak{F}^\pow$, \RrefPlusPlus{} requires the following:
\begin{itemize}
\item if $Y\subseteq \mathrm{R}_a[X]$, then there is some nonempty $X'\subseteq X$ such that for all nonempty $X''\subseteq X'$, $Y\subseteq \mathrm{R}_a[X'']$,
\end{itemize}
which implies:
\begin{itemize}
\item if every world in $Y$ can be ``seen'' by some world or other in $X$, then there is some single world in $X$ that can ``see'' every world in $Y$. 
\end{itemize}
This obviously does not hold for all Kripke frames $\mathfrak{F}$, so \RrefPlusPlus{} does not hold for all $\mathfrak{F}^\pow$.
\end{proof}

\noindent For similar reasons, we cannot always transform a $\mathcal{V}$-BAO into a Humberstone frame as in \S~\ref{GFPF}.

In \citealt{Humberstone1981} (p.~326), it is stated that one can prove the completeness of some standard modal logics with respect to classes of \textit{atomless} Humberstone frames (recall \S~\ref{AtomlessFull}) using a canonical model construction in which each possibility is the set of syntactic consequences of a consistent \textit{finite} set of formulas, the refinement relation is the subset relation between sets of formulas, and the accessibility relations and valuation function are defined as for the canonical Kripke frame (Definition \ref{CanKrip}). However, the \RrefPlusPlus{} condition is too strong for such a construction to work, as the following fact shows. 

\begin{fact}[Infinitary \RrefPlusPlus{}] If $\sig$ is infinite and $\mathbf{L}\subseteq\mathcal{L}(\sig,\ind)$ is a normal modal logic, then there is no partial-state model $\mathcal{M}$ with a nonempty relation $R_i$ satisfying \RrefPlusPlus{} such that for every $x\in\mathcal{M}$, there is a \textit{finite} $\Gamma_x\subseteq\mathcal{L}(\sig,\ind)$ such that $\{\varphi\in\mathcal{L}(\sig,\ind)\mid \mathcal{M},x\Vdash\varphi\}=\{\varphi\in\mathcal{L}(\sig,\ind)\mid \Gamma_x\vdash_\mathbf{L} \varphi\}$.
\end{fact}
\begin{proof} For reductio, suppose there is such a model $\mathcal{M}$. Since $R_i$ is nonempty, there are $x,y\in\mathcal{M}$ such that $xR_iy$. Since $\{\varphi\in\mathcal{L}(\sig,\ind)\mid \mathcal{M},y\Vdash\varphi\}=\{\varphi\in\mathcal{L}(\sig,\ind)\mid \Gamma_y\vdash_\mathbf{L} \varphi\}$, the finite set $\Gamma_y$ is \textbf{L}-consistent. It follows that  $\sig(y)=\{p\in\sig\mid \mathcal{M},y\Vdash p\}$ is finite, because no finite \textbf{L}-consistent set $\Gamma_y$ entails infinitely many $p\in\sig$ (For if $\vdash_\mathbf{L}\bigwedge\Gamma_y \rightarrow p$ for a $p$ not occurring in $\Gamma_y$, then $\vdash_\mathbf{L}\bigwedge\Gamma_y\rightarrow\bot$ by Uniform Substitution.)  Thus, $\sig\setminus \sig(y)$ is infinite. Since $xR_i y$, by \RrefPlusPlus{}, there is an $x'\sqsubseteq x$ such that for all $x''\sqsubseteq x'$, $x''R_i y$, which implies $\mathcal{M},x''\nVdash\Box_i p$ for every $p \in \sig\setminus \sig(y)$. Since this holds for all $x''\sqsubseteq x'$, we have $\mathcal{M},x'\Vdash \neg\Box_i p$ for every $p\in \sig\setminus \sig(y)$. Thus, there are infinitely many $p\in\sig$ such that $\mathcal{M},x'\Vdash \neg\Box_i p$. But then there is no finite $\Gamma_{x'}$ such that $\{\varphi\in\mathcal{L}(\sig,\ind)\mid \mathcal{M},x'\Vdash\varphi\}=\{\varphi\in\mathcal{L}(\sig,\ind)\mid \Gamma_{x'}\vdash_\mathbf{L} \varphi\}$, because no finite \textbf{L}-consistent set $\Gamma_{x'}$ entails $\neg\Box_i p$ for infinitely many $p\in\sig$ . For if $\vdash_\mathbf{L} \bigwedge \Gamma_{x'}\rightarrow \neg\Box_ip$ for a $p$ not occurring in $\Gamma_{x'}$, then $\vdash_\mathbf{L} \bigwedge \Gamma_{x'}\rightarrow \neg\Box_i\top$ by Uniform Substitution, which with $\vdash_\mathbf{L}\Box_i\top$ gives us $\vdash_\mathbf{L} \neg\bigwedge \Gamma_{x'}$.\end{proof} 

By contrast, we can prove the completeness of various modal logics with respect to classes of atomless possibility frames satisfying \Rref{} using a model construction in which each possibility is (an equivalence class of) a single finite formula, as in \S~\ref{AtomlessFull} (also see \citealt{Holliday2014}). 

In \S~\ref{DualEquiv}, we showed that full possibility frames and $\mathcal{CV}$-BAOs can be turned into semantically equivalent \textit{rich} possibility frames, which satisfy, among other conditions, \Rref{} and \Rmax{}, the latter being definitive of \textit{quasi-functional} possibility frames as in \S~\ref{FuncFrames}. This kind of result cannot be proved with \RrefPlusPlus{} in place of \Rref{}, as shown by Fact \ref{DisMod}.\ref{DisMod4}.
 
\begin{fact}[\RrefPlusPlus{} and \Rmax{}]\label{DisMod} For any possibility model $\mathcal{M}=\langle S,\sqsubseteq,\{R_i\}_{i\in\ind},\pi\rangle$ with $R_i$ satisfying \RrefPlusPlus{} and $x,y\in S$:
\begin{enumerate}
\item\label{DisMod1} if $\mathcal{M},x\Vdash \Box_i p \vee\Box_i q$ and $xR_iy$, then $\mathcal{M},y\Vdash p$ or $\mathcal{M},y\Vdash q$;
\item\label{DisMod2} if $R_i$ satisfies \Rmax{} and $\mathcal{M},x\Vdash\Box_i p\vee\Box_i q$, then $\mathcal{M},x\Vdash \Box_i p$ or $\mathcal{M},x\Vdash\Box_i q$;
\item\label{DisMod3} if $R_i$ and $R_j$ satisfy \Rmax{}, then $\Box_j (\Box_i p\vee\Box_i q)\rightarrow (\Box_j \Box_i p\vee\Box_j\Box_iq)$ is valid over the underlying frame of $\mathcal{M}$;
\item\label{DisMod4} \textbf{K} is not complete with respect to the class of quasi-functional models satisfying \RrefPlusPlus{}.
\end{enumerate}
\end{fact}
\begin{proof} For part \ref{DisMod1}, suppose $\mathcal{M},y\nVdash p$ and $\mathcal{M},y\nVdash q$. By \RrefPlusPlus{}, $xR_iy$ implies that $\exists x'\sqsubseteq x$ $\forall x''\sqsubseteq x'$, $x''R_i y$, so $\mathcal{M},x''\nVdash \Box_ip$ and $\mathcal{M},x''\nVdash \Box_i q$. Since this holds for all $x''\sqsubseteq x'$, we have $\mathcal{M},x'\Vdash \neg\Box_i p$ and $\mathcal{M},x'\Vdash \neg\Box_i q$, which with $x'\sqsubseteq x$ implies $\mathcal{M},x\nVdash \Box_i p \vee\Box_i q$ by Fact \ref{ForcingDis}.

For part \ref{DisMod2}, if $R_i(x)$ is empty, then $\mathcal{M},x\Vdash \Box_i \varphi$ for every $\varphi\in\mathcal{L}(\sig,\ind)$, so we are done. If $R_i(x)$ is nonempty, then by \Rmax{} it has a maximum, $f_i(x)$, in which case $\mathcal{M},x\Vdash\Box_i p\vee\Box_i q$ and part \ref{DisMod1} together imply that $\mathcal{M},f_i(x)\Vdash p$ or $\mathcal{M},f_i(x)\Vdash q$. If $\mathcal{M},f_i(x)\Vdash p$, then since $f_i(x)$ is the maximum of $R_i(x)$, it follows by Persistence that $\mathcal{M},y\Vdash p$ for all $y\in R_i(x)$, so $\mathcal{M},x\Vdash\Box_i p$. By the same reasoning, if $\mathcal{M},f_i(x)\Vdash q$, then $\mathcal{M},x\Vdash\Box_i q$. Thus, $\mathcal{M},x\Vdash\Box_i p$ or $\mathcal{M},x\Vdash\Box_i q$.

For part \ref{DisMod3}, suppose $\mathcal{M},x\Vdash \Box_j (\Box_i p\vee\Box_i q)$. If $R_j(x)$ is empty, then $\mathcal{M},x\Vdash \Box_j \varphi$ for every $\varphi\in\mathcal{L}(\sig,\ind)$, so we are done. If $R_j(x)$ is nonempty, then by \Rmax{} it has a maximum $f_j(x)$, in which case $\mathcal{M},x\Vdash \Box_j (\Box_i p\vee\Box_i q)$ implies $\mathcal{M},f_j(x)\Vdash \Box_i p\vee\Box_i q$. Then by part \ref{DisMod2}, either $\mathcal{M},f_j(x)\Vdash \Box_i p$ or $\mathcal{M},f_j(x)\Vdash\Box_i q$, so either $\mathcal{M},x\Vdash \Box_j\Box_ip$ or $\mathcal{M},x\Vdash \Box_j\Box_iq$.

Part \ref{DisMod4} follows from part \ref{DisMod3} and the fact that $\Box_j (\Box_i p\vee\Box_i q)\rightarrow (\Box_j \Box_i p\vee\Box_j\Box_iq)$ is not a theorem of \textbf{K}.
\end{proof}

We conclude with an observation about the intermediate \RrefPlus{}. If one goes beyond ZF set theory and assumes the ultrafilter axiom, one can prove that filter frames as in \S~\ref{GFPF} satisfy \RrefPlus{}.

\begin{proposition}[\RrefPlus{} for Filter Frames]\label{can+} For any BAO $\mathbb{A}$, $\mathbb{A}_{\gff}$ and $\mathbb{A}_{\ff}$ satisfy \RrefPlus{}.
\end{proposition}

\begin{proof} For proper filters ${X},{Y}$ in $\mathbb{A}$, suppose ${X}R_i{Y}$. By the ultrafilter axiom, there is an ultrafilter ${Y'}\supseteq{Y}$, so ${Y'}\sqsubseteq{Y}$. Where
\begin{equation}\mathrm{X}'={X}\cup\{-\blacksquare_i y\mid y\not\in {Y'}\},\label{AnotherXDef}\end{equation}
suppose for reductio that $[\mathrm{X}')$ is not a proper filter, i.e., $\inc\in [\mathrm{X}')$. Then by (\ref{AnotherXDef}) and Fact \ref{FiltGen}, there are $x\in{X}$ and $y_1,\dots,y_k\not\in {Y'}$ such that $x\meet -\blacksquare_i y_1\meet\dots\meet -\blacksquare_i y_k\leq \inc$, which implies
$x\leq \blacksquare_i (y_1\join \dots\join y_k)$
by the properties of $\blacksquare_i$. Then since ${X}$ is a filter and $x\in {X}$,  $ \blacksquare_i (y_1\join \dots\join y_k)\in{X}$, which with ${X}R_i{Y}$ implies $y_1\join \dots\join y_k\in{Y}$, which in turn implies $y_1\join \dots\join y_k\in{Y'}$. But since ${Y'}$ is an ultrafilter, $y_1,\dots,y_k\not\in {Y'}$ implies $-y_1,\dots,-y_k\in {Y'}$, which contradicts $y_1\join \dots\join y_k\in{Y'}$. Thus, ${X'}=[\mathrm{X}')$ is a proper filter. 

Now consider any proper filter ${X''}\sqsubseteq{X'}$, i.e., ${X''}\supseteq{X'}$. If $y\not\in {Y'}$, then by (\ref{AnotherXDef}), $-\blacksquare_i y\in{X'}$, so $\blacksquare_i y\not\in{X''}$ by the fact that ${X''}$ is a proper filter. Thus, ${X''}R_i{Y'}$, which establishes \RrefPlus{}.
\end{proof}

\subsection{Separative Quotients}\label{SepProof}

Finally, we prove the following proposition from \S~\ref{SepSec}.

\SepQuo*

For a possibility frame $\mathcal{F}=\langle S,\sqsubseteq, \{R_i\}_{i\in\ind},\adm\rangle$, recall the equivalence relation $\simeq_\cofsub$ on $S$ from Definition~\ref{CoRef}, defined by: $x\simeq_\cofsub y$ iff $x\cof y$ and $y\cof x$. For $x\in S$, let $[x]_\simeq =\{x'\in S\mid x\simeq_\cofsub x'\}$ be the equivalence class of $x$. Then the frame $\mathcal{F}^\simeq$ for Proposition \ref{sep-quo} is defined as follows.

\begin{definition}[Separative Quotient] Given a possibility frame $\mathcal{F}=\langle S,\sqsubseteq, \{R_i\}_{i\in\ind},\adm\rangle$, define the \textit{separative quotient} $\mathcal{F}^\simeq=\langle S^\simeq,\sqsubseteq^\simeq,\{R_i^\simeq\}_{i\in\ind },\adm^\simeq\rangle$ of $\mathcal{F}$ by:
\begin{enumerate}
\item $S^\simeq = \{[x]_{\simeq}\mid x\in S\}$;
\item $[x]_{\simeq}\sqsubseteq^\simeq [y]_{\simeq}$ iff $x\cof y$;
\item $[x]_{\simeq}R_i^\simeq [y]_{\simeq}$ iff $\exists x'\in [x]_{\simeq}$ $\exists y'\in [y]_{\simeq}$: $x'R_iy'$;
\item $\adm^\simeq = \{X^\simeq \mid  X\in \adm\}$, where $X^\simeq = \{[x]_{\simeq} \mid x\in X\}$. \hfill $\triangleleft$
\end{enumerate}
\end{definition}

\noindent We will now prove the proposition, defining the morphism $h$ from $\mathcal{F}$ to $\mathcal{F}^\simeq$ by $h(x)=[x]_{\simeq}$ (so $X^\simeq=h[X]$).

\begin{proof}[Proof of Proposition \ref{sep-quo}] First, we show that $\langle S^\simeq,\sqsubseteq^\simeq\rangle$ is separative. Consider the relation $\cof^\simeq$ defined in terms of $\sqsubseteq^\simeq$ in $\mathcal{F}^\simeq$. Since $[x]_{\simeq}\sqsubseteq^\simeq [y]_{\simeq}$ implies $[x]_{\simeq}\cof^\simeq [y]_{\simeq}$, we need only prove the converse. If $[x]_{\simeq}\cof^\simeq [y]_{\simeq}$, then  for all $ [x']_{\simeq}\sqsubseteq^\simeq [x]_{\simeq}$ there is a $[x'']_{\simeq}\sqsubseteq^\simeq [x']_{\simeq}$ with $[x'']_{\simeq}\sqsubseteq^\simeq [y]_{\simeq}$. This means that (a) $\forall x'\cof x$ $\exists x''\cof x'$: $x''\cof y$.  To show $x\cof y$, consider any $x'\sqsubseteq x$. Then $x'\cof x$, so by (a), there is an ${x''}\cof{x'}$ with ${x''}\cof y$.  Since ${x''}\cof {x'}$, there is an ${x'''}\sqsubseteq {x''}$ with ${x'''}\sqsubseteq {x'}$. Then together ${x''}\cof y$ and ${x'''}\sqsubseteq {x''}$ imply that there is an ${x''''}\sqsubseteq {x'''}$ with ${x''''}\sqsubseteq y$. By the transitivity of $\sqsubseteq$, ${x''''}\sqsubseteq {x'''}\sqsubseteq {x'}$ implies ${x''''}\sqsubseteq {x'}$. Thus, for any ${x'}\sqsubseteq x$, there is an ${x''''}\sqsubseteq {x'}$ with ${x''''}\sqsubseteq y$. Hence $x\cof y$, so $[x]_{\simeq}\sqsubseteq^\simeq [y]_{\simeq}$.

Second, we show that $\adm^\simeq\subseteq\mathrm{RO}(\mathcal{F}^\simeq)$. Consider an $X^\simeq\in \adm^\simeq$. To show that $X^\simeq$ satisfies \textit{persistence}, suppose $[x]_{\simeq}\in X^\simeq$ and $[x']_{\simeq}\sqsubseteq^\simeq [x]_{\simeq}$, so $x\in X$ and $x'\cof x$. Then since $X\in\mathrm{RO}(\mathcal{F})$, it follows that $x'\in X$ by Fact \ref{CofClose}, so $[x']_{\simeq}\in X^\simeq$. To show that $X^\simeq$ satisfies \textit{refinability}, suppose $[x]_{\simeq}\not\in X^\simeq$, so $x\not\in X$. Then since $X\in\mathrm{RO}(\mathcal{F})$, $x\not\in X$ implies that there is an $x'\sqsubseteq x$ such that (b) for all $x''\sqsubseteq x'$, $x''\not\in X$. Now consider any $y'\cof x'$, so there is a $y''\sqsubseteq y'$ such that $y''\sqsubseteq x'$, which with (b) implies $y''\not\in X$, which with $y''\sqsubseteq y'$ and \textit{persistence} for $X$ implies $y'\not\in X$. Thus, for all $y'\cof x'$, $y'\not\in X$, so for all $[y']_{\simeq}\sqsubseteq^\simeq [x']_{\simeq}$, $[y']_{\simeq}\not\in X^\simeq$. Finally, from $x'\sqsubseteq x$ we have $[x']_{\simeq}\sqsubseteq^\simeq [x]_{\simeq}$, so $[x']_{\simeq}$ is the witness we need for \textit{refinability}.

Next, we must prove that $\mathcal{F}^\simeq$ is a partial-state frame. To do so, we first show:
\begin{itemize}
\item[] (i) $(X\cap Y)^\simeq=X^\simeq\cap Y^\simeq$; (ii) $(X\supset Y)^\simeq=X^\simeq\supset^\simeq Y^\simeq$; and (iii) $(\blacksquare_i X)^\simeq = \blacksquare_i^\simeq X^\simeq$.
\end{itemize}
Checking  (i) is straightforward. For (ii), from right to left, if $[z]_{\simeq}\not\in (X\supset Y)^\simeq$, so  $z\not\in X\supset Y$, then there is a $z'\sqsubseteq z$ such that $z'\in X$ but $z'\not\in Y$, which implies $[z']_{\simeq}\sqsubseteq^\simeq [z]_{\simeq}$ and $[z']_{\simeq}\in X^\simeq$ but $[z']_{\simeq}\not\in Y^\simeq$, so $[z]_{\simeq}\not\in X^\simeq \supset^\simeq Y^\simeq$. From left to right, if $[z]_{\simeq}\not\in X^\simeq \supset^\simeq Y^\simeq$, then there is a $[z']_{\simeq}\sqsubseteq^\simeq [z]_{\simeq}$ such that $[z']_{\simeq}\in X^\simeq$ but $[z']_{\simeq}\not\in Y^\simeq$, so $z'\in X$ but $z'\not\in Y$, so $z'\not\in X\supset Y$. Since $[z']_{\simeq}\sqsubseteq^\simeq [z]_{\simeq}$, we have $z'\cof z$. Finally, since $X\supset Y\in\mathrm{RO}(\mathcal{F})$, $z'\not\in X\supset Y$ and $z'\cof z$ together imply $z\not\in X\supset Y$ by Fact \ref{CofClose}, so $[z]_{\simeq}\not\in (X\supset Y)^\simeq$. For (iii), from right to left, if $[x]_{\simeq}\not\in (\blacksquare_iX)^\simeq$, so $x\not\in\blacksquare_i X$, then there is a $y\in S$ such that $xR_iy$ and $y\not\in X$, which implies $[x]_{\simeq}R_i^\simeq [y]_{\simeq}$ and $[y]_{\simeq}\not\in X^\simeq$, so $[x]_{\simeq}\not\in\blacksquare_i^\simeq X^\simeq$. From left to right, if $[x]_{\simeq}\not\in \blacksquare_i^\simeq X^\simeq $, then there is a $[y]_{\simeq}\in S^\simeq$ such that $[x]_{\simeq}R_i^\simeq [y]_{\simeq}$ and $[y]_{\simeq}\not\in X^\simeq$. From $[x]_{\simeq}R_i^\simeq [y]_{\simeq}$, we have $\exists x'\in [x]_{\simeq}$ $\exists y'\in [y]_{\simeq}$: $x'R_iy'$.  From $[y]_{\simeq}\not\in X^\simeq$, we have $y\not\in X$. Since $X\in \mathrm{RO}(\mathcal{F})$, together $y\not\in X$ and $y\simeq_\cofsub y'$ imply $y'\not\in X$ by Fact \ref{CofClose}, which with $x'R_iy'$ implies $x'\not\in\blacksquare_i X$. Then since $\blacksquare_iX\in\mathrm{RO}(\mathcal{F})$, together $x'\not\in\blacksquare_i X$ and $x\simeq_\cofsub x'$ imply $x\not\in\blacksquare_i X$ by Fact \ref{CofClose}, so $[x]_{\simeq}\not\in (\blacksquare_iX)^\simeq$.

Since $\mathcal{F}$ is a partial-state frame, $\adm$ is closed under $\cap$, $\supset$, and $\blacksquare_i$. It follows from (i)-(iii) and the definition of $\adm^\simeq$ that $\adm^\simeq$ is also closed under $\cap$, $\supset^\simeq$, and $\blacksquare_i^\simeq$, so $\mathcal{F}^\simeq$ is a partial-state frame. Since we also showed above that $\adm^\simeq\subseteq \mathrm{RO}(\mathcal{F}^\simeq)$, we conclude that $\mathcal{F}^\simeq$ is a possibility frame. To see that $\mathcal{F}^\simeq$ is \textit{full} if $\mathcal{F}$ is, one can easily show that if $\mathcal{X}\in \mathrm{RO}(\mathcal{F}^\simeq)$, then $h^{-1}[\mathcal{X}]\in \mathrm{RO}(\mathcal{F})=\adm$, so $\mathcal{X}=h[h^{-1}[\mathcal{X}]]=(h^{-1}[\mathcal{X}])^\simeq\in \adm^\simeq$. 

The function $h\colon \mathcal{F}\to\mathcal{F}^\simeq$ defined by $h(x)=[x]_{\simeq}$ is surjective (which gave us $\mathcal{X}=h[h^{-1}[\mathcal{X}]]$ above). Given surjectivity, the condition for a \textit{robust} morphism is that for all $X\in \adm$, $X=h^{-1}[h[X]]$ and $h[X]\in \adm^\simeq$. By the definition of $\adm^\simeq$, $X\in \adm$ implies $h[X]=X^\simeq\in \adm^\simeq$. Also observe that $h^{-1}[X^\simeq]=X$, so $h^{-1}[h[X]]=X$. Next, it follows from (ii) above, taking $Y=\emptyset$, that $h$ satisfies the \SqMatch{} clause of possibility morphisms; it follows from (iii) above that $h$ satisfies \RMatch{}; and it follows from the definition of $\adm^\simeq$ and $h^{-1}[X^\simeq]=X$ that $h$ satisfies \PullBack{}. Thus, $h$ is a surjective robust possibility morphism.
\end{proof} 
                
\bibliographystyle{plainnat}
\bibliography{possibility-frames}
    
\end{document}